\documentclass[reqno,10pt,centertags]{amsart}
\usepackage{amsmath,amsthm,amscd,amssymb,latexsym,esint,upref,%stmaryrd,
enumerate,color,verbatim,%yfonts,
mathrsfs,mathtools,accents}
\usepackage{hyperref}

\newcommand*{\mailto}[1]{\href{mailto:#1}{\nolinkurl{#1}}}

\newcommand{\bbC}{{\mathbb{C}}}

\newcommand{\bbN}{{\mathbb{N}}}

\newcommand{\bbR}{{\mathbb{R}}}
\newcommand{\bbS}{{\mathbb{S}}}

\newcommand{\cB}{{\mathcal B}}

\newcommand{\cD}{{\mathcal D}}

\newcommand{\cH}{{\mathcal H}}

\newcommand{\cN}{{\mathcal N}}

\newcommand{\cS}{{\mathcal S}}

\newcommand{\cX}{{\mathcal X}}

\newcommand{\ga}{{\mathfrak{a}}}

\newcommand{\gs}{{\mathfrak{s}}}
\newcommand{\gt}{{\mathfrak{t}}}

\DeclareMathOperator{\supp}{supp}

\DeclareMathOperator{\ran}{ran}
\DeclareMathOperator{\dom}{dom}
\DeclareMathOperator{\essinf}{essinf}

\DeclareMathOperator*{\slim}{s-lim}

\renewcommand{\Im}{\text{\rm Im}}
\renewcommand{\ln}{\text{\rm ln}}

\newcommand{\norm}[1]{\lVert#1\rVert}

\newcommand{\no}{\notag}

\newcommand{\ol}{\overline}

\newcommand{\bi}{\bibitem}
\newcommand{\wti}{\widetilde}

\newcommand{\lb}{\label}
\newcommand{\dott}{\,\cdot\,}

\newcommand{\Om}{\Omega}
\newcommand{\dOm}{{\partial\Omega}}

\let\geq\geqslant
\let\leq\leqslant

\renewcommand{\dot}{\overset{\textbf{\Large.}}}

\renewcommand{\dotplus}{\overset{\textbf{\Large.}} +}

\makeatletter
\def\theequation{\@arabic\c@equation}
\allowdisplaybreaks 
\numberwithin{equation}{section}

\newtheorem{theorem}{Theorem}[section]
\newtheorem{proposition}[theorem]{Proposition}
\newtheorem{lemma}[theorem]{Lemma}
\newtheorem{corollary}[theorem]{Corollary}
\newtheorem{definition}[theorem]{Definition}
\newtheorem{hypothesis}[theorem]{Hypothesis}

\theoremstyle{remark}
\newtheorem{remark}[theorem]{Remark}
\begin{document}

\numberwithin{equation}{section}
\allowdisplaybreaks

\title[Sharp Boundary Trace Theory and Schr\"{o}dinger Operators]{Sharp Boundary Trace Theory and Schr\"{o}dinger Operators on Bounded Lipschitz Domains}

\author[J.\ Behrndt]{Jussi Behrndt}  
\address{Institut f\"{u}r Angewandte Mathematik, Technische Universit\"at Graz, Steyrergasse 30, 
8010 Graz, Austria and Department of Mathematics, Stanford University, 450 Jane Stanford Way,\
Stanford CA 94305-2125, USA} 
\email{\mailto{behrndt@tugraz.at}} 
%\email{behrndt@tugraz.at}  
\urladdr{\url{http://www.math.tugraz.at/~behrndt/}}
%\urladdr{www.math.tugraz.at/~behrndt/}
  
\author[F.\ Gesztesy]{Fritz Gesztesy}
\address{Department of Mathematics,
Baylor University, Sid Richardson Bldg., 1410 S.\,4th Street,
Waco, TX 76706, USA}
\email{\mailto{Fritz\_Gesztesy@baylor.edu}}
%\email{Fritz\_Gesztesy@baylor.edu}
\urladdr{\url{http://www.baylor.edu/math/index.php?id=935340}}
%\urladdr{http://www.baylor.edu/math/index.php?id=935340}

\author[M.\ Mitrea]{Marius Mitrea}
\address{Department of Mathematics,
Baylor University, Sid Richardson Bldg., 1410 S.\,4th Street,
Waco, TX 76706, USA}
\email{\mailto{Marius\_Mitrea@baylor.edu}}
%\email{Marius\_Mitrea@baylor.edu}
\urladdr{\url{https://www.baylor.edu/math/index.php?id=962939}}
%\urladdr{https://www.baylor.edu/math/index.php?id=962939}

%\dedicatory{}
\thanks{J.~B.~gratefully acknowledges support for the Distinguished Visiting Austrian Chair at Stanford University by the Europe Center and the Freeman Spogli Institute for International Studies. Work of M.~M.~was partially supported by the Simons Foundation Grant $\#$\,637481.} 

\thanks{To appear in Memoirs AMS}

\date{\today}

%%%%%%%%%%%%
%%%%%%%%%%%%
\begin{abstract} 
We develop a sharp boundary trace theory in arbitrary bounded Lipschitz domains 
which, in contrast to classical results, allows ``forbidden'' endpoints and permits the 
consideration of functions exhibiting very limited regularity. This is done at the (necessary) expense of 
stipulating an additional regularity condition involving the action of the Laplacian on the functions in question 
which, nonetheless, works perfectly with the Dirichlet and Neumann realizations of the Schr\"odinger 
differential expression 
$-\Delta+V$. In turn, this boundary trace theory serves as a platform for developing a spectral theory for Schr\"odinger 
operators on bounded Lipschitz domains, along with their associated Weyl--Titchmarsh operators. Overall, this pushes the present 
state of knowledge a significant step further. For example, we succeed in extending the Dirichlet and Neumann trace 
operators in such a way that all self-adjoint extensions of a Schr\"odinger operator on a bounded Lipschitz domain 
may be described with explicit boundary conditions, thus providing a final answer to a problem that has been investigated 
for more than 60 years in the mathematical literature. Along the way, a number of other open problems are solved. 
The most general geometric and analytic setting in which the theory developed here yields satisfactory results is that 
of Lipschitz subdomains of Riemannian manifolds and for the corresponding Laplace--Beltrami operator (in place of 
the standard flat-space Laplacian). In particular, such an extension yields results for variable coefficient 
Schr\"odinger operators on bounded Lipschitz domains.
\end{abstract}
%%%%%%%%%%%%
%%%%%%%%%%%%

\maketitle 

%\vspace*{-3mm} 
{\scriptsize{\tableofcontents}}
%\normalsize

%%%%%%%%%%%%%%%%%%%%%%%%%%%%%%%%%%%%%%%%%%
%%%%%%%%%%%%%%%%%%%%%%%%%%%%%%%%%%%%%%%%%%
\section{Introduction}  
\label{s1}
%%%%%%%%%%%%%%%%%%%%%%%%%%%%%%%%%%%%%%%%%%
%%%%%%%%%%%%%%%%%%%%%%%%%%%%%%%%%%%%%%%%%%

Given an open set $\Omega\subseteq{\mathbb{R}}^n$, let $H^s(\Omega)$ denote the $L^2$-based 
Sobolev space of (fractional) order $s\in{\mathbb{R}}$ in $\Omega$. When 
$\Omega={\mathbb{R}}^{n-1}\oplus{\mathbb{R}}_{+}$, the upper half-space, starting from the realization that 
$C^\infty(\overline{\Omega})\cap H^1(\Omega)$ is dense in $H^1(\Omega)$ (as may be seen via translation and 
a standard mollifying argument) and the restriction-to-the-boundary map 
$C^\infty(\overline{\Omega})\ni u\mapsto f:=u(\cdot,0)\in C^\infty({\mathbb{R}}^{n-1})$ satisfies
$\|f\|_{L^2({\mathbb{R}}^{n-1})}\leq C_n\|u\|_{H^1(\Omega)}$ for each 
$u\in C^\infty(\overline{\Omega})\cap H^1(\Omega)$, one concludes that 
the assignment $u\mapsto f$ extends uniquely to a linear and bounded mapping, henceforth referred to 
as the Dirichlet boundary trace operator $\gamma_D$, from $H^1(\Omega)$ into $L^2({\mathbb{R}}^{n-1})$.
This trace operator is not surjective, since N.\ Aronszajn \cite{Ar55} (see also \cite{Sl66}) has noted 
that its image may be described as 
\begin{equation}\label{8uygg-Y}
\gamma_D(H^1(\Omega))=\bigg\{f\in L^2({\mathbb{R}}^{n-1}) \, \bigg| \
\int_{{\mathbb{R}}^{n-1}}|\xi||\widehat{f}(\xi)|^2\,d^{n-1}\xi<\infty\bigg\}
\end{equation}
where ``hat'' stands for the Fourier transform in ${\mathbb{R}}^{n-1}$. This result has been 
subsequently extended by E.\ Gagliardo, whose work in \cite{Ga57} marks the beginning of a flurry 
of activities concerning trace theory which, in turn, has firmly established this topic in the 
present day mathematical landscape. 

For example, we now know that if $\Omega\subseteq{\mathbb{R}}^n$ is a bounded Lipschitz domain then the 
restriction-to-the-boundary map $C^\infty(\overline{\Omega})\ni u\mapsto f:=u\big|_{\partial\Omega}\in C^0(\partial\Omega)$ 
extends uniquely to a linear and continuous operator 
\begin{equation}\label{7yfFc-P}
\gamma_D:H^s(\Omega)\to H^{s-(1/2)}(\partial\Omega)\,\text{ whenever }\, 1/2 < s < 3/2.
\end{equation}
Furthermore, the Dirichlet trace operator $\gamma_D$ is surjective in the above context and, in fact, 
admits a continuous linear right-inverse. 

The study of trace operators like \eqref{7yfFc-P} interfaces tightly with the issue of extending functions 
from Sobolev spaces (and other smoothness scales) defined intrinsically in $\Omega$ to the entire Euclidean space 
${\mathbb{R}}^n$ with preservation of class. More generally, given a set $F\subseteq{\mathbb{R}}^n$ which is $d$-dimensional 
in a certain sense for some $d\in(0,n]$, the question arises whether it is possible to extend any function 
$f$ belonging to a Besov space $B^{p,p}_\beta(F)$, suitably defined on $F$, to a function in 
$B^{p,p}_\alpha({\mathbb{R}}^n)$ where the smoothness exponents $\alpha,\beta$ satisfy 
$\alpha=\beta+[(n-d)/p]$. 
As far as traces are concerned, in place of \eqref{7yfFc-P} one may ask if the trace on $F$ of any function from 
$B^{p,p}_\alpha({\mathbb{R}}^n)$ lies in $B^{p,p}_\beta(F)$. For example, such an extension/restriction problem has an affirmative solution if $F$ is a $d$-dimensional plane in ${\mathbb{R}}^n$, say $F:={\mathbb{R}}^d\times\{0\}^{n-d}$, for any $d\in\{1,\dots,n\}$. 

The extension/restriction problems leading to this and other related results 
have been studied by many authors. Early contributors include N.\ Aronszajn, F.\ Mulla, and P.\ Szeptycki \cite{AMS63}, 
O.\ V.\ Besov \cite{Be64}, \cite{Be68}, V.\ I.\ Burenkov \cite{Bu65}, 
A.\ P.\ Calder\'on \cite{Ca61}, E.\ Gagliardo \cite{Ga57}, J.\ L.\ Lions and E.\ Magenes \cite{LM60}-\cite{LM72}, 
P.\ I.\ Lizorkin \cite{Li63}, J.\ Ne\v{c}as \cite{Ne67}, S.\ M.\ Nikol'ski\u{\i} \cite{Ni61}, \cite{Ni75}, 
E.\ M.\ Stein \cite{St61}, \cite{St70}, and M.\ H.\ Taibleson \cite{Ta64}, and S.\ V.\ Uspenski\u{\i} \cite{Us70}, among others. 
Let us also note that the case when $F$ is a surface in ${\mathbb{R}}^n$ satisfying a local Lipschitz condition has been studied by 
O.\ V.\ Besov in \cite{Be72}, \cite{Be72b}, \cite{Be72c}, while extension and 
restriction problems for $F$ an arbitrary $d$-dimensional closed subset of ${\mathbb{R}}^n$ (see \eqref{Fq-A.4} below)
have been investigated by D.\ R.\ Adams \cite{Ad71},  A.\ Jonsson \cite{Jo73}, J.\ Petree \cite{Pe75}, T.\ Sj\"odin \cite{Sj75}, and H.\ Wallin \cite{Wa63}.

In \cite{JW78} A.\ Jonsson and H.\ Wallin have initiated a breakthrough, proving a very general extension/restriction 
theorem on the Besov scale for $d$-sets. We recall that a closed set $F\subseteq{\mathbb{R}}^n$ is said to be a
$d$-set for some $d\in(0,n]$, provided there exists some finite constant $C\geq 1$ with the property that
\begin{eqnarray}\label{Fq-A.4}
C^{-1}r^{d}\leq {\mathscr{H}}^{d}(B(x,r)\cap F)\leq Cr^{d},\,\,\,\forall\,x\in F,\,\,\,0<r\leq{\rm diam}\,(F),
\end{eqnarray}
where ${\mathscr{H}}^d$ is the $d$-dimensional Hausdorff measure in ${\mathbb{R}}^n$. 
(For example, the closure $\overline{\Omega}$ of a Lipschitz domain $\Omega\subseteq{\mathbb{R}}^n$ is an $n$-set, 
while its topological boundary $\partial\Omega$ is an $(n-1)$-set; parenthetically, we also note that the boundary 
of Koch's snowflake in ${\mathbb{R}}^2$ is a $d$-set with $d:= \ln (4)/\ln (3)$.) In this context, a brand of 
Besov spaces has been introduced by A.\ Jonsson and H.\ Wallin in \cite{JW78} as follows. Given $p\in[1,\infty)$ and
$s\in(0,\infty)\backslash {\mathbb{N}}$, define the Besov space $B^{p,p}_{s}(F)$ as 
the collection of families $\dot{f}:=\{f_\alpha\}_{|\alpha|\leq [s]}$
(where $[s]$ denotes the integer part of $s$), whose components are ${\mathscr{H}}^{d}$-measurable functions 
on $F$ with the property that if for each multi-index $\alpha\in{\mathbb{N}}_0^n$ with $|\alpha|\leq[s]$ one introduces
\begin{equation}\label{Fq-A.5PL1}
R_\alpha(x,y):=f_{\alpha}(x)-\sum_{|\beta|\leq[s]-|\alpha|}\frac{(x-y)^\beta}{\beta!}
f_{\alpha+\beta}(y)\,\,\mbox{ for ${\mathscr{H}}^{d}$-a.e. }\,x,y\in F,
\end{equation}
then $\big\|\dot{f}\big\|_{B^{p,p}_{s}(F)}<+\infty$, where
\begin{eqnarray}\label{Fq-A.5PL2}
&&\hskip -0.20in
\big\|\dot{f}\big\|_{B^{p,p}_{s}(F)}:=\sum_{|\alpha|\leq[s]}\Big(\int_{F}|f_\alpha|^p\,d{\mathscr{H}}^{d}\Big)^{1/p}
\\[6pt]
&&\hskip 0.20in
+\left(\sum_{j=0}^\infty \sum_{|\alpha|\leq[s]}2^{j(s-|\alpha|)p+jd}
\iint\limits_{\substack{x,y\in F\\ |x-y|<2^{-j}}}|R_\alpha(x,y)|^p
\,d{\mathscr{H}}^{d}(x)\,d{\mathscr{H}}^{d}(y)\right)^{1/p}.
\nonumber
\end{eqnarray}

The following fundamental result regarding traces and extensions on (and from) arbitrary $d$-sets in ${\mathbb{R}}^n$ 
has been proved by A.\ Jonsson and H.\ Wallin in \cite[Main Theorem, p.\,146]{JW78}.

\begin{theorem}[Jonsson--Wallin Trace/Extension Theory]\label{GCC-67}
Assume $F\subseteq{\mathbb{R}}^n$ is a given $d$-set for some $d\in(0,n]$. 
Fix a number $k\in{\mathbb{N}}_0$ along with some integrability exponent $p\in[1,\infty)$. 
Also, select two smoothness exponents $\alpha,\beta$ satisfying 
\begin{equation}\label{GGG-89J5g}
\beta\in(k,k+1)\,\,\text{ and }\,\,\alpha=\beta+[(n-d)/p].
\end{equation}
Finally, it is agreed that a barred integral sign denotes an integral average.

Then, for every scalar function $u\in B^{p,p}_\alpha({\mathbb{R}}^n)$, the vector-valued limit 
\begin{equation}\label{THna}
\big({\mathscr{R}}_F^{(k)}u\big)(x):=\left\{\lim\limits_{r\to 0^{+}}\fint_{B(x,r)}
(\partial^\alpha u)(y)\,d^ny\right\}_{|\alpha|\leq k-1}
\,\,\,\text{ exists at ${\mathscr{H}}^d$-a.e. }\,x\in F
\end{equation}
and, defined as such, this higher-order trace operator on $F$ induces a well defined, linear, and bounded mapping  
\begin{equation}\label{GGG-89}
{\mathscr{R}}_F^{(k)}:B^{p,p}_\alpha({\mathbb{R}}^n)\longrightarrow B^{p,p}_{\beta}(F).
\end{equation}

In the converse direction, there exists a linear and bounded operator 
\begin{equation}\label{GGG-90}
{\mathscr{E}}^{(k)}_F:B^{p,p}_{\beta}(F)\longrightarrow B^{p,p}_\alpha({\mathbb{R}}^n)
\end{equation}
with the property that
\begin{equation}\label{GGG-91}
{\mathscr{R}}_F^{(k)}\circ{\mathscr{E}}^{(k)}_F=I,
\,\,\,\mbox{ the identity on }\,\,\,B^{p,p}_{\beta}(F).
\end{equation}
\end{theorem}

Subsequently, the program initiated in \cite{JW78} has been amply expanded by A.\ Jonsson and H.\ Wallin in their monograph \cite{JW84}. The body of work described so far in the introduction is of immense practical value and various refinements (allowing two integrability exponents $p\not=q$, other scales of spaces measuring smoothness, alternative proofs, etc.) have since come to light. See, for instance, \cite{AF03}, \cite{BMMM13}, \cite{Ch84}, \cite{Ch92}, \cite{Co88}, \cite{Gr85}, \cite{JK95}, \cite{Jo81}, \cite{KMM07}, \cite{Ma93}, \cite{MM04}, \cite{Ma81}, \cite{Ma11}, \cite{MMS10}, \cite{Mc00}, \cite{Mi18}, \cite{MMM22a}, \cite{Mi93}, \cite{Ro06}, \cite{RS96}, \cite{Ry99}, \cite{Sa95}, \cite{Sa96}, \cite{Se89}, \cite{Tr92a}, \cite{Tr01}, \cite{Tr02}, \cite{Tr06}, \cite{Tr08}, and the references therein. This is but an indicative sample of a large body of works on the subject of traces and extensions, which remains an active topic of research to date. 

Note that the well definiteness, boundedness, and surjectivity of the trace operator $\gamma_D$ in \eqref{7yfFc-P}
is a very special case of Theorem~\ref{GCC-67} when $s\not=1$ (corresponding to $p=2$ and $d=n-1$).
Indeed, if $\Omega$ is a Lipschitz domain then any function in $H^s(\Omega)$ with $s\in\big(\tfrac{1}{2},\tfrac{3}{2}\big) \big\backslash \{1\}$ 
may be extended to $H^s({\mathbb{R}}^n)=B^{2,2}_s({\mathbb{R}}^n)$ and \eqref{THna}--\eqref{GGG-89} apply to this extension, 
bearing in mind that $F:=\partial\Omega$ is an $(n-1)$-set. The case when $\Omega$ is a Lipschitz domain and $s=1$ may be reduced to 
the situation when $\Omega={\mathbb{R}}^{n-1}\oplus{\mathbb{R}}_{+}$, the upper half-space, via a simple localization and a 
bi-Lipschitz change of variables flattening the boundary.

For a bounded Lipschitz domain $\Omega\subseteq{\mathbb{R}}^n$, the end-point cases $s=1/2$ 
and $s=3/2$ in \eqref{7yfFc-P} are problematic. As regards the limiting value $s=1/2$, 
it has been pointed out in the last paragraph of \cite[p.~180]{JK95} that $C^\infty_0(\Omega)$ is dense 
in $H^{1/2}(\Omega)$. Consequently, the restriction-to-the-boundary map 
\begin{equation}\label{8uygg-Y.1111}
C^\infty(\overline{\Omega})\ni u\mapsto u\big|_{\partial\Omega}\in C^0(\partial\Omega)
\end{equation}
vanishes identically on a dense subspace of $H^{1/2}(\Omega)$, so its unique extension to $H^{1/2}(\Omega)$ 
is the trivial map $\gamma_D(u)=0$ for all $u\in H^{1/2}(\Omega)$. The space $\gamma_D(H^{3/2}(\Omega))$, 
identified in \cite{JW84}, has a rather technical description. Even in the case of a bounded $C^1$-domain $\Omega$ 
this space looks very different from the natural candidate in the smooth case (when $\Omega$ is a bounded 
$C^\infty$-domain, the Dirichlet boundary trace maps $H^{3/2}(\Omega)$ continuously onto $H^{1}(\partial\Omega)$).
Hence, in sharp contrast with the $C^\infty$ case, there is a substantial change in the character of the trace 
operator on the boundary of a bounded $C^1$-domain corresponding to $s=3/2$. 
In fact, in \cite[Proposition~3.2, p.~176]{JK95} a bounded $C^1$-domain $\Omega\subseteq{\mathbb{R}}^2$ 
and a function $u\in H^{3/2}(\Omega)$ are constructed with the property that $\gamma_Du\notin H^1(\partial\Omega)$. 
This goes to show that, corresponding to the limiting value $s=3/2$, the range of the 
Dirichlet trace operator in \eqref{7yfFc-P} is {\it strictly larger} than $H^{1}(\partial\Omega)$.

In the present work we succeed in incorporating the end-points $\{\tfrac{1}{2},\tfrac{3}{2}\}$ 
in the range of indices for which the Dirichlet trace operator behaves naturally, in a given bounded
Lipschitz domain $\Omega\subseteq{\mathbb{R}}^n$. As is apparent from our earlier discussion, 
for this to happen we need to restrict $\gamma_D$ to a smaller domain than $H^s(\Omega)$ with 
$s\in\big[\tfrac{1}{2},\tfrac{3}{2}\big]$, that is, demand that $\gamma_D$ acts from a subspace 
of $H^s(\Omega)$ consisting of functions satisfying further regularity assumptions. The novel 
idea is that, starting with $u\in H^{s}(\Omega)$ for some $s\in\big[\tfrac{1}{2},\tfrac{3}{2}\big]$,  
if $\Delta u$ is slightly more regular than typical action of the Laplacian on functions
from $H^{s}(\Omega)$, that is, more regular than $H^{s-2}(\Omega)$, then we may meaningfully define 
its Dirichlet boundary trace $\gamma_D u$ for the full range $s\in\big[\tfrac{1}{2},\tfrac{3}{2}\big]$.

Simply put, if the function $u\in H^{s}(\Omega)$ with $s\in\big[\tfrac{1}{2},\tfrac{3}{2}\big]$
has a ``better-than-expected'' Laplacian (in the sense of membership to a certain smoothness scale)
then $\gamma_Du$ is well defined in $H^{s-(1/2)}(\partial\Omega)$. An embodiment of this principle 
on the scale of Sobolev spaces is Theorem~\ref{YTfdf-T} which states that if $\Omega\subset\bbR^n$ 
is a bounded Lipschitz domain and $\varepsilon>0$ is arbitrary, then the restriction-to-the-boundary 
operator \eqref{8uygg-Y.1111} induces a unique, well defined, linear, surjective, continuous map
\begin{equation}\label{eqn:gammaDs.2aux.INTRO}
\gamma_D:\big\{u\in H^s(\Omega)\,\big|\,\Delta u\in H^{s-2+\varepsilon}(\Omega)\big\}
\rightarrow H^{s-(1/2)}(\partial\Omega),\quad\forall\,s\in\big[\tfrac{1}{2},\tfrac{3}{2}\big],    
\end{equation}
if the space on the left is equipped with the norm 
$u\mapsto\|u\|_{H^{s}(\Omega)}+\|\Delta u\|_{H^{s-2+\varepsilon}(\Omega)}$. 
For example, this implies that for each $\varepsilon>0$,
\begin{equation}\label{eqn:gammaDs.2aux.INTRO.WACO}
\big\{u\in H^{3/2}(\Omega)\,\big|\,\Delta u\in H^{-(1/2)+\varepsilon}(\Omega)\big\}
\ni u\mapsto\gamma_D(\nabla u)\in[L^2(\partial\Omega)]^n
\end{equation}
is a well defined, linear, and bounded operator. 
In this context, it is also significant to observe that the domain of the Dirichlet trace operator in 
\eqref{eqn:gammaDs.2aux.INTRO} embeds (strictly) in certain Triebel--Lizorkin spaces. Specifically,  
as noted in \eqref{eqn:gammaDs.2aux.WACO}, we have the continuous strict embeddings 
\begin{equation}\label{eqn:gammaDs.2aux.WACO.INTRO}
\begin{array}{c}
\big\{u\in H^s(\Omega)\,\big|\,\Delta u\in H^{s-2+\varepsilon}(\Omega)\big\}
\hookrightarrow F^{2,q}_s(\Omega)\hookrightarrow H^s(\Omega)
\\[2pt]
\quad\text{for any }\,s\in\big[\tfrac{1}{2},\tfrac{3}{2}\big],\,\text{ any }\,\varepsilon>0,\,\text{ and any }\,q\in(0,2).
\end{array}
\end{equation}
Thus, the demand that $\Delta u\in H^{s-2+\varepsilon}(\Omega)$ improves the regularity 
of $u\in H^s(\Omega)$, albeit in a subtle fashion. 

Employing Besov spaces allows us to express in an even more precise fashion the amount of regularity one needs 
to impose on $\Delta u$ in order to be able to allow the end-points $s\in\big\{\tfrac{1}{2},\tfrac{3}{2}\big\}$
in \eqref{7yfFc-P}. Concretely, given any bounded Lipschitz domain $\Omega\subset\bbR^n$, the 
restriction-to-the-boundary operator \eqref{8uygg-Y.1111} induces a unique, well defined, linear, 
surjective, continuous map
\begin{equation}\label{eqn:gammaDs.2aux-Mi.INTRO}
\gamma^{\#}_D:\big\{u\in H^s(\Omega)\,\big|\,\Delta u\in B^{2,1}_{s-2}(\Omega)\big\}
\rightarrow H^{s-(1/2)}(\partial\Omega),\quad\forall\,s\in\big[\tfrac{1}{2},\tfrac{3}{2}\big],    
\end{equation}
where, this time, the space on the left-hand side of \eqref{eqn:gammaDs.2aux-Mi.INTRO} 
is equipped with the norm $u\mapsto\|u\|_{H^{s}(\Omega)}+\|\Delta u\|_{B^{2,1}_{s-2}(\Omega)}$. 
Reassuringly, the sharp Dirichlet trace $\gamma^{\#}_D$ from \eqref{eqn:gammaDs.2aux-Mi.INTRO} 
is compatible with $\gamma_D$ in \eqref{eqn:gammaDs.2aux.INTRO}.
Also, from \eqref{eqn:gammaDs.2aux-Mi.INTRO} (with $s=1/2$) we see 
that for each $\varepsilon>0$ we have a well defined, linear, and bounded operator 
\begin{equation}\label{eqn:gammaDs.2aux.INTRO.WACO.TEXAS}
\big\{u\in H^{3/2}(\Omega)\,\big|\,\Delta u\in B^{2,1}_{-1/2}(\Omega)\big\}
\ni u\mapsto\gamma^{\#}_D(\nabla u)\in [L^2(\partial\Omega)]^n.
\end{equation}
See Theorem~\ref{YTfdf-T-Mi} for a more expansive and nuanced result of this flavor.  
In particular, it has been noted in \eqref{eqn:gammaDs.2aux.WACO.BBBB} that the domain of the sharp Dirichlet 
trace operator in \eqref{eqn:gammaDs.2aux-Mi.INTRO} embeds (strictly) in certain Triebel--Lizorkin spaces. 
Specifically, we have the continuous strict embeddings 
\begin{equation}\label{eqn:gammaDs.2aux.WACO.BBBB.INTRO}
\big\{u\in H^s(\Omega)\,\big|\,\Delta u\in B^{2,1}_{s-2}(\Omega)\big\}
\hookrightarrow F^{2,1}_s(\Omega)\hookrightarrow H^s(\Omega),\quad s\in\big[\tfrac{1}{2},\tfrac{3}{2}\big].
\end{equation}

In addition to the results for the Dirichlet boundary trace operator, we develop in Section~\ref{s5} 
a similar theory for the Neumann boundary trace operator $\gamma_N$ in the context of Sobolev 
spaces in a given bounded Lipschitz domain $\Omega\subset\bbR^n$. More specifically, the Neumann trace 
map originally defined as $u\mapsto\nu\cdot(\nabla u)|_{\partial\Omega}$ for functions $u\in C^\infty(\overline{\Omega})$, 
where $\nu$ denotes the outward unit normal to $\Omega$, extends uniquely to linear, continuous, surjective operators  
\begin{equation}\label{eqn:gammaN-pp.INTRO}
\gamma_N:\big\{u\in H^s(\Omega)\,\big|\,\Delta u\in L^2(\Omega)\big\} 
\rightarrow H^{s-(3/2)}(\partial\Omega),\quad s\in\big[\tfrac{1}{2},\tfrac{3}{2}\big],
\end{equation}
that are compatible with one another, when the space on the left-hand side of \eqref{eqn:gammaN-pp.INTRO} 
is equipped with the natural graph norm $u\mapsto\|u\|_{H^{s}(\Omega)}+\|\Delta u\|_{L^2(\Omega)}$.
See Theorem~\ref{YTfdf.NNN.2-Main} and Corollary~\ref{YTfdf.NNN.2} in this regard. Here we only wish to mention 
that, with $\nu$ denoting the outward unit normal vector to $\Omega$, 
\begin{align}\label{eq:Nnan7yg-P.CCC.INTRO}
\begin{split}
& \text{if $u\in H^{3/2}(\Omega)$ has $\Delta u\in L^2(\Omega)$ then 
$\gamma_Nu=\nu\cdot\gamma_D(\nabla u)\in L^2(\partial\Omega)$} 
\\[2pt]
& \quad\text{with the Dirichlet trace taken in the sense of \eqref{eqn:gammaDs.2aux.INTRO}.}
\end{split}
\end{align}
It is remarkable that $\gamma_N$ in \eqref{eqn:gammaN-pp.INTRO} acts on a class of functions $u$  
for which the notion of the ``classical'' Neumann trace of $\nu\cdot\gamma_D(\nabla u)$ is utterly ill defined.
To illustrate this via an example, take $\Omega:=B(0,1)$ the unit ball in ${\mathbb{R}}^n$ and for each $\alpha\in(0,1)$ 
consider $u_\alpha(x):=(1-|x|^2)^\alpha$ for each $x\in\Omega$. Then $u_\alpha\in H^s(\Omega)$ for each 
$s<\alpha+(1/2)$, yet $\nabla u_\alpha$ blows up (in the limit) at each boundary point $x\in\partial\Omega$.

Compared to earlier work, the crucial new ingredient here is the use of well-posedness results for the 
$L^2$ Dirichlet, Neumann, and Regularity boundary value problems in bounded Lipschitz domains in which 
the size of the solution is measured using the nontangential maximal operator and boundary traces are 
taken in a nontangential pointwise sense. In this regard, we heavily rely on the basic work in \cite{Da77}, 
\cite{JK81}, \cite{MT99}, \cite{MT00a}, \cite{MT01}, \cite{Ve84}. We also make essential use 
of solvability results and estimates for the corresponding inhomogeneous problems from 
\cite{FMM98}, \cite{JK95}, \cite{MT00}.

One of the primary motivations for developing a sharp boundary trace theory in bounded Lipschitz domains
(which now includes the traditionally forbidden end-points $1/2$ and $3/2$) is because 
this provides a platform for the study of Schr\"{o}dinger operators in this class of domains. 
The format of our brand of trace theorems (cf.~\eqref{eqn:gammaDs.2aux.INTRO}) is perfectly suited 
for such a study, which we take up in Section~\ref{s6}. There, among a variety of topics, we discuss the self-adjoint 
Friedrichs extension and the self-adjoint Dirichlet and Neumann realizations of $-\Delta+V$ where the 
potential $V$ is a real-valued essentially bounded function. We then proceed to introduce $z$-dependent 
Dirichlet-to-Neumann maps, otherwise known as Weyl--Titchmarsh operators, for Schr\"odinger operators on 
bounded Lipschitz domains in Section~\ref{s7}. In turn, these results are used in Section~\ref{s8} to construct 
what we call maximal extensions of the Dirichlet and Neumann trace operators 
on arbitrary bounded Lipschitz domains in ${\mathbb{R}}^n$. 

More specifically, the goal in Section~\ref{s8} is to further extend the Dirichlet trace operator
\eqref{eqn:gammaDs.2aux.INTRO}, and its Neumann counterpart $\gamma_N$, by continuity onto the 
domain of $A_{max,\Omega}$, the maximal realization of $-\Delta+V$ defined as 
(with all derivatives taken in the sense of distributions)
\begin{equation}\label{We-Q.2.INTR}
A_{max,\Omega}:=-\Delta+V,\quad
\dom(A_{max,\Omega}):=\big\{f\in L^2(\Omega)\,\big|\,\Delta f\in L^2(\Omega)\big\}.
\end{equation}
To describe the said extensions of the Dirichlet and Neumann traces, we bring to the forefront the 
spaces
\begin{equation}\label{eqn:G0G1.INTR}
\mathscr{G}_D(\partial\Omega):=\ran\big(\gamma_D\big|_{\dom(A_{N,\Omega})}\big)
\,\text{ and }\,\mathscr{G}_N(\partial\Omega):=\ran\big(\gamma_N\big|_{\dom(A_{D,\Omega})}\big),
\end{equation} 
where 
\begin{equation}\label{We-Q.10EE.WACO}
\begin{split}
& A_{D,\Omega}=-\Delta+V, 
\\[2pt]
& \dom(A_{D,\Omega})=\big\{f\in H^1(\Omega)\,\big|\,\Delta f\in L^2(\Omega)\,\text{ and }\,\gamma_D f=0\big\},
\end{split}
\end{equation}
and
\begin{equation}\label{We-Q.10EE-jussiNN.WACO}
\begin{split}
& A_{N,\Omega}=-\Delta+V, 
\\[2pt]
&\dom(A_{N,\Omega})=\big\{f\in H^1(\Omega)\,\big|\,\Delta f\in L^2(\Omega)\,\text{ and }\,\gamma_N f=0\big\},
\end{split}
\end{equation}
are, respectively, the Dirichlet and Neumann self-adjoint realizations of the differential expression 
$-\Delta+V$ in the Lipschitz domain $\Omega$ (studied in Section~\ref{s6}). In the rough setting considered here, 
the spaces $\mathscr{G}_D(\partial\Omega)$, $\mathscr{G}_N(\partial\Omega)$ turn out to be the correct 
substitutes for $H^{3/2}(\partial\Omega)$ and, respectively, $H^{1/2}(\partial\Omega)$, to which they reduce 
if $\Omega$ were to be a bounded $C^\infty$-domain. Indeed, work in \cite{GM11} shows that 
\begin{align}\label{Tan-C5.WACO}
\begin{split} 
& \mathscr{G}_D(\partial\Omega)=H^{3/2}(\partial\Omega) \, \text{ and } \, 
\mathscr{G}_N(\partial\Omega)=H^{1/2}(\partial\Omega)    \\
& \quad \text{when $\Omega$ is a bounded $C^{1,r}$-domain with $r>1/2$}
\end{split}
\end{align} 
(where the parameter $r$ refers to the H\"older regularity of the first order derivatives of the functions whose graphs 
locally describe $\partial\Omega$). In fact, (cf. \cite{GM11}), 
\begin{equation}\label{Tan-C5.WACO.2}
\parbox{9.00cm}{whenever $\Omega$ is some bounded open convex set,  
or some bounded $C^{1,r}$-domain for some $r>1/2$, 
it follows that $\dom(A_{D,\Omega}),\dom(A_{N,\Omega})\subset H^2(\Omega)$ and the Dirichlet trace operator 
$\gamma_D:H^2(\Omega)\to H^{3/2}(\partial\Omega)$ as well as the Neumann trace operator 
$\gamma_N:H^2(\Omega)\to H^{1/2}(\partial\Omega)$ are both well defined, bounded, and onto.}
\end{equation} 
As such, the duals 
$\mathscr{G}_D(\partial\Omega)^\ast$, $\mathscr{G}_N(\partial\Omega)^\ast$ should be thought of as natural 
substitutes for $H^{-3/2}(\partial\Omega)$ and, respectively, $H^{-1/2}(\partial\Omega)$, in the rough setting considered here. 
See also \cite{BM13} in this regard.

The following theorem, which presents the most complete 
result along the lines of work in \cite{BM13}, \cite{GM08}, \cite{GM11}, is one of the central results in this work. 

%%%%%%%%%%%%
\begin{theorem}\label{t5.5.INTR}
Assume that $\Omega\subset\bbR^n$ is a bounded Lipschitz domain, 
and that the potential $V\in L^\infty(\Omega)$ is a real-valued function. 
Then the following statements hold: 

\vskip 0.08in 
\noindent $(i)$ The spaces $\mathscr{G}_D(\partial\Omega),\mathscr{G}_N(\partial\Omega)$ carry a natural Hilbert 
space structure $($see item {\it $(vi)$} below for equivalent norms\,$)$ and the Dirichlet trace operator 
$\gamma_D$ {\rm (}from \eqref{eqn:gammaDs.2aux.INTRO}{\rm )} along with its counterpart, the Neumann trace operator 
$\gamma_N$ {\rm (}from \eqref{eqn:gammaN-pp.INTRO}{\rm )}, extend by continuity to continuous surjective mappings 
\begin{equation}\label{4.55.INTR}
\begin{split}
\widetilde{\gamma}_D:\dom(A_{max,\Omega})&\to\mathscr{G}_N(\partial\Omega)^*,
\\[2pt]  
\widetilde{\gamma}_N:\dom(A_{max,\Omega})&\to\mathscr{G}_D(\partial\Omega)^*,  
\end{split}
\end{equation}
where $\dom(A_{max,\Omega})$ is endowed with the graph norm of $A_{max,\Omega}$, 
and $\mathscr{G}_D(\partial\Omega)^*$, $\mathscr{G}_N(\partial\Omega)^*$ are, respectively, 
the adjoint $($conjugate dual\,$)$ spaces of $\mathscr{G}_D(\partial\Omega)$,
$\mathscr{G}_N(\partial\Omega)$ carrying the natural topology induced by the said Hilbert space structure. 

\vskip 0.08in  
\noindent $(ii)$ These extensions satisfy 
\begin{equation}\label{4.30.INTR}
\ker(\widetilde{\gamma}_D)=\dom(A_{D,\Omega})\,\text{ and }\,\ker(\widetilde{\gamma}_N)=\dom(A_{N,\Omega}).
\end{equation} 
Also, for each $s\in[0,1]$ there exists a constant $C\in(0,\infty)$ with the property that
\begin{align}\label{4.55.REE.D.INTR}
\begin{split}
& f\in\dom(A_{max,\Omega})\,\text{ and }\,\widetilde{\gamma}_D f\in H^s(\partial\Omega)
\,\text{ imply }\,f\in H^{s+(1/2)}(\Omega)
\\[2pt]
& \quad\text{and }\,\|f\|_{H^{s+(1/2)}(\Omega)}\leq C\big(\|\Delta f\|_{L^2(\Omega)}
+\|\widetilde{\gamma}_D f\|_{H^s(\partial\Omega)}\big),
\end{split}
\end{align}
and
\begin{align}\label{4.55.REE.N.INTR}
\begin{split}
& f\in\dom(A_{max,\Omega})\,\text{ and }\,\widetilde{\gamma}_N f\in H^{-s}(\partial\Omega)
\,\text{ imply }\,f\in H^{-s+(3/2)}(\Omega)
\\[2pt]
& \quad\text{and }\,\|f\|_{H^{-s+(3/2)}(\Omega)}\leq C\big(\|f\|_{L^2(\Omega)}+\|\Delta f\|_{L^2(\Omega)}
+\|\widetilde{\gamma}_N f\|_{H^{-s}(\partial\Omega)}\big).
\end{split}
\end{align} 

\noindent $(iii)$ 
With $\accentset{\circ}{H}^2(\Omega)$ denoting the closure of $C_0^{\infty}(\Omega)$ in $H^2(\Omega)$ and 
with $\widetilde{\gamma}_D,\widetilde{\gamma}_N$ as in \eqref{4.55.INTR}, one has
\begin{align}\label{ut444.INTR}
\accentset{\circ}{H}^2(\Omega)=\big\{f\in\dom(A_{max,\Omega})\,\big|\,
& \wti{\gamma}_D f=0\,\text{ in }\,\mathscr{G}_N(\partial\Omega)^*
\nonumber\\[2pt]
&\text{and }\,\wti{\gamma}_N f=0\,\text{ in }\,\mathscr{G}_D(\partial\Omega)^*\big\}.  
\end{align}

\noindent $(iv)$ The manner in which the mapping $\widetilde{\gamma}_D$ in \eqref{4.55.INTR} 
operates is as follows: Given $f\in\dom(A_{max,\Omega})$ and some arbitrary 
$\phi\in\mathscr{G}_N(\partial\Omega)$, there exists 
$g\in H^{3/2}(\Omega)\cap\dom(A_{max,\Om})$ such that $\gamma_D g=0$ 
and $\gamma_N g=\phi$, and the functional $\widetilde{\gamma}_D f\in\mathscr{G}_N(\partial\Omega)^*$ 
acts {\rm (}in a coherent fashion{\rm )} on the given $\phi$ according to 
\begin{equation}\label{eq:33f4iU.INTR}
{}_{\mathscr{G}_N(\partial\Omega)^\ast}\big\langle\widetilde{\gamma}_D f,
\phi\big\rangle_{\mathscr{G}_N(\partial\Omega)}=(f,\Delta g)_{L^2(\Omega)}-(\Delta f,g)_{L^2(\Omega)}.
\end{equation}
As a consequence, the following Green's formula holds: 
\begin{equation}\label{eq:33VaV.INTR}
{}_{\mathscr{G}_N(\partial\Omega)^\ast}\big\langle\widetilde{\gamma}_D f,
\gamma_N g\big\rangle_{\mathscr{G}_N(\partial\Omega)}=(f,\Delta g)_{L^2(\Omega)}-(\Delta f,g)_{L^2(\Omega)},
\end{equation}
for each $f\in\dom(A_{max,\Omega})$ and each $g\in\dom(A_{D,\Om})$. 

\vskip 0.08in  
\noindent $(v)$ The mapping $\widetilde{\gamma}_N$ in \eqref{4.55.INTR} operates in the 
following fashion: Given a function $f\in\dom(A_{max,\Omega})$ along with some arbitrary 
$\psi\in\mathscr{G}_D(\partial\Omega)$, there exists $g\in H^{3/2}(\Omega)\cap\dom(A_{max,\Om})$ 
such that $\gamma_N g=0$ and $\gamma_D g=\psi$, and the functional 
$\widetilde{\gamma}_N f\in\mathscr{G}_D(\partial\Omega)^*$ acts {\rm (}in a coherent fashion{\rm )} 
on the given $\psi$ according to 
\begin{equation}\label{eq:33f4iU.NN.INTR}
{}_{\mathscr{G}_D(\partial\Omega)^\ast}\big\langle\widetilde{\gamma}_N f,\psi\big\rangle_{\mathscr{G}_D(\partial\Omega)}
=-(f,\Delta g)_{L^2(\Omega)}+(\Delta f,g)_{L^2(\Omega)}.
\end{equation}
In particular, the following Green's formula holds: 
\begin{equation}\label{eq:33VaV.NN.INTR}
{}_{\mathscr{G}_D(\partial\Omega)^\ast}\big\langle\widetilde{\gamma}_N f,
\gamma_D g\big\rangle_{\mathscr{G}_D(\partial\Omega)}=-(f,\Delta g)_{L^2(\Omega)}+(\Delta f,g)_{L^2(\Omega)},
\end{equation}
for each $f\in\dom(A_{max,\Omega})$ and each $g\in\dom(A_{N,\Om})$. 

\vskip 0.08in 
\noindent $(vi)$ The operators 
\begin{align}\label{4.5ASdf.1.INTR}
&\gamma_D:\dom(A_{N,\Om})=H^{3/2}(\Omega)\cap\dom(A_{max,\Omega})\cap\ker(\gamma_N)
\to\mathscr{G}_D(\partial\Omega),
\\[2pt] 
&\gamma_N:\dom(A_{D,\Om})=H^{3/2}(\Omega)\cap\dom(A_{max,\Omega})\cap\ker(\gamma_D)
\to\mathscr{G}_N(\partial\Omega),    
\label{4.5ASdf.2.INTR}
\end{align}
are well defined, linear, surjective, and continuous if for some $s\in[0,3/2]$ both 
spaces on the left-hand sides of \eqref{4.5ASdf.1.INTR}, \eqref{4.5ASdf.2.INTR} are equipped 
with the norm $f\mapsto\|f\|_{H^s(\Omega)}+\|\Delta f\|_{L^2(\Omega)}$ 
{\rm (}which are all equivalent{\rm )}. In addition, 
\begin{equation}\label{eq:Ajf.INTR}
\text{the kernel of $\gamma_D$ and $\gamma_N$ in \eqref{4.5ASdf.1.INTR}--\eqref{4.5ASdf.2.INTR} 
is $\accentset{\circ}{H}^2(\Omega)$}.
\end{equation}
Moreover, 
\begin{align}\label{i7g5r.INTR}
\|\phi\|_{\mathscr{G}_D(\partial\Omega)} & \approx\inf_{\substack{f\in
H^{3/2}(\Omega)\cap\dom(A_{max,\Omega})\\ \gamma_N f=0,\,\,\gamma_D f=\phi}}
\big(\|f\|_{H^{3/2}(\Omega)}+\|\Delta f\|_{L^2(\Omega)}\big)
\nonumber\\[2pt]
& \approx\inf_{\substack{f\in H^{3/2}(\Omega)\cap\dom(A_{max,\Omega})\\ \gamma_N f=0,\,\,\gamma_D f=\phi}}
\big(\|f\|_{L^2(\Omega)}+\|\Delta f\|_{L^2(\Omega)}\big)
\nonumber\\[2pt]
& \approx\inf_{\substack{f\in\dom(A_{max,\Omega})\\ \widetilde{\gamma}_N f=0,\,\,\widetilde{\gamma}_D f=\phi}}
\big(\|f\|_{L^2(\Omega)}+\|\Delta f\|_{L^2(\Omega)}\big),
\end{align}
uniformly for $\phi\in\mathscr{G}_D(\partial\Omega)$, and 
\begin{align}\label{i7g5r.2N.INTR}
\|\psi\|_{\mathscr{G}_N(\partial\Omega)} & \approx\inf_{\substack{g\in
H^{3/2}(\Omega)\cap\dom(A_{max,\Omega})\\ \gamma_D g=0,\,\,\gamma_N g=\psi}}
\big(\|g\|_{H^{3/2}(\Omega)}+\|\Delta g\|_{L^2(\Omega)}\big)
\nonumber\\[2pt]
& \approx\inf_{\substack{g\in H^{3/2}(\Omega)\cap\dom(A_{max,\Omega})\\ \gamma_D g=0,\,\,\gamma_N g=\psi}}
\big(\|g\|_{L^2(\Omega)}+\|\Delta g\|_{L^2(\Omega)}\big)
\nonumber\\[2pt]
& \approx\inf_{\substack{g\in\dom(A_{max,\Omega})\\ \widetilde\gamma_D g=0,\,\,\widetilde\gamma_N g=\psi}}
\big(\|g\|_{L^2(\Omega)}+\|\Delta g\|_{L^2(\Omega)}\big)
\nonumber\\[2pt]
& \approx\inf_{\substack{g\in\dom(A_{max,\Omega})\\ \widetilde{\gamma}_D g=0,\,\,\widetilde{\gamma}_N g=\psi}}
\|\Delta g\|_{L^2(\Omega)},
\end{align}
uniformly for $\psi\in\mathscr{G}_N(\partial\Omega)$. As a consequence,
\begin{align}\label{ytrr555e.INTR}
\begin{split}
& \mathscr{G}_D(\partial\Omega)\hookrightarrow H^1(\partial\Omega)
\hookrightarrow L^2(\partial\Omega)\hookrightarrow H^{-1}(\partial\Omega)
\hookrightarrow\mathscr{G}_D(\partial\Omega)^*,  
\\[2pt] 
& \qquad\qquad\quad
\mathscr{G}_N(\partial\Omega)\hookrightarrow L^2(\partial\Omega)
\hookrightarrow\mathscr{G}_N(\partial\Omega)^*, 
\end{split}
\end{align} 
with all embeddings linear, continuous, and with dense range. Moreover, the duality 
pairings between $\mathscr{G}_D(\partial\Omega)$ and $\mathscr{G}_D(\partial\Omega)^\ast$,
as well as between $\mathscr{G}_N(\partial\Omega)$ and $\mathscr{G}_N(\partial\Omega)^\ast$, 
are both compatible with the inner product in $L^2(\partial\Omega)$. \\[1mm] 
\noindent $(vii)$ For each $z\in\rho(A_{D,\Omega})$, the boundary value problem 
\begin{equation}\label{eqn:bvp.2.INTR}
\begin{cases} 
(-\Delta+V-z)f=0\,\text{ in $\Omega,\quad f\in\dom(A_{max,\Omega})$,}   
\\[2pt]  
\widetilde{\gamma}_D f=\varphi
\,\text{ in $\mathscr{G}_N(\partial\Omega)^*,\quad\varphi\in\mathscr{G}_N(\partial\Omega)^*$,}
\end{cases} 
\end{equation}
is well posed. In particular, for each $z\in\rho(A_{D,\Omega})$ there exists a constant $C\in(0,\infty)$, which depends 
only on $\Omega$, $n$, $z$, and $V$, with the property that
\begin{align}\label{eqn:bvp.2.INTR.WACO} 
\begin{split} 
& \|f\|_{L^2(\Omega)}\leq C\|\widetilde{\gamma}_D f\|_{\mathscr{G}_N(\partial\Omega)^*}
\,\text{ for each }\,f\in\dom(A_{max,\Omega})      
\\[1pt]
& \quad\text{with }\,(-\Delta+V-z)f=0\,\text{ in }\,\Omega.
\end{split} 
\end{align}
Likewise, for each $z\in\rho(A_{N,\Omega})$, the boundary value problem 
\begin{equation}\label{eqn:bvpNeumann.INTR}
\begin{cases} 
(-\Delta+V-z)f=0\,\text{ in $\Omega,\quad f\in\dom(A_{max,\Omega})$,}   
\\[2pt]  
-\widetilde{\gamma}_N f=\varphi 
\,\text{ in $\mathscr{G}_D(\partial\Omega)^*,\quad\varphi\in\mathscr{G}_D(\partial\Omega)^*$,}
\end{cases} 
\end{equation}
is well posed. In particular, for each $z\in\rho(A_{N,\Omega})$ there exists a constant $C\in(0,\infty)$, which depends 
only on $\Omega$, $n$, $z$, and $V$, with the property that
\begin{align}\label{eqn:bvpNeumann.INTR.WACO}
\begin{split}
&\|f\|_{L^2(\Omega)}\leq C\|\widetilde{\gamma}_N f\|_{\mathscr{G}_D(\partial\Omega)^*}\,\text{ for each }\,f\in\dom(A_{max,\Omega})     
\\[1pt]
& \quad\text{with }\,(-\Delta+V-z)f=0\,\text{ in }\,\Omega.
\end{split} 
\end{align}
\end{theorem}
%%%%%%%%%%

The powerful machinery developed in Theorem~\ref{t5.5.INTR} allows us to settle a number of outstanding issues. First of all, 
this allows us to address the following question posed (to the current last-named author) by G.\ Uhlmann in 2004 (\cite{Uh04}):
\begin{quote}\label{Uhlmann}
{\it ``If $\Omega$ is a bounded Lipschitz domain in ${\mathbb{R}}^n$ and $f$ is in $H^{1/2}(\partial\Omega)$, 
there exists a unique harmonic function $u$ in $\Omega$ with {\rm [}Dirichlet{\rm ]} trace $f$, and $u$ satisfies 
$\|u\|_{H^1(\Omega)}\leq C\|f\|_{H^{1/2}(\partial\Omega)}$. Is it also true that 
$\|u\|_{L^2(\Omega)}\leq C\|f\|_{H^{-1/2}(\partial\Omega)}$? This holds for smooth domains.''}
\end{quote}
Specifically, since in the case $V=0$ we have $0\in\rho(A_{D,\Omega})$, one concludes from \eqref{eqn:bvp.2.INTR.WACO} that 
\begin{equation}\label{eqn:bvp.2.INTR.WACO.bis}
\|u\|_{L^2(\Omega)}\leq C\|\widetilde{\gamma}_D u\|_{\mathscr{G}_N(\partial\Omega)^*}
\,\text{ for each harmonic function }\,u\in L^2(\Omega).
\end{equation}
In fact, given the boundedness of $\widetilde{\gamma}_D$ in the context of 
\eqref{4.55.INTR}, the opposite inequality in \eqref{eqn:bvp.2.INTR.WACO.bis} also holds so that, ultimately, 
\begin{equation}\label{eqn:bvp.2.INTR.WACO.bis.567}
\|u\|_{L^2(\Omega)}\approx\|\widetilde{\gamma}_D u\|_{\mathscr{G}_N(\partial\Omega)^*}
\,\text{ uniformly in $u\in L^2(\Omega)$ a harmonic function}.
\end{equation}
In view of the fact that $\widetilde{\gamma}_D$ from \eqref{4.55.INTR} is an extension of the ordinary 
Dirichlet trace operator $\gamma_D$ {\rm (}from \eqref{eqn:gammaDs.2aux.INTRO}{\rm )}, we therefore have
\begin{equation}\label{eqn:bvp.2.INTR.WACO.bis.123} 
\|u\|_{L^2(\Omega)}\leq C\|{\gamma}_D u\|_{\mathscr{G}_N(\partial\Omega)^*}\,\text{ for each harmonic function }\,u\in H^1(\Omega). 
\end{equation}
Moreover, combining \eqref{Tan-C5.WACO} with \eqref{eqn:bvp.2.INTR.WACO.bis} yields that 
\begin{align}\label{Tan-C5.WACO.TX}
\begin{split} 
& \text{whenever $\Omega$ is a bounded $C^{1,r}$-domain with $r>1/2$, one has}    
\\
& \quad \text{$\|u\|_{L^2(\Omega)}\leq C\|\widetilde{\gamma}_D u\|_{H^{-1/2}(\partial\Omega)}$ 
for each harmonic function $u\in L^2(\Omega)$.}
\end{split} 
\end{align} 
More generally, in the case when the potential $V$ satisfies $L^\infty(\Omega)\ni V\geq 0$ at a.e.~point in the 
bounded Lipschitz domain $\Omega\subset{\mathbb{R}}^n$,  we continue to have $0\in\rho(A_{D,\Omega})$ 
so \eqref{eqn:bvp.2.INTR.WACO} yields 
\begin{equation}\label{eqn:bvp.2.INTR.WACO.bis.2}
\|u\|_{L^2(\Omega)}\leq C\|\widetilde{\gamma}_D u\|_{\mathscr{G}_N(\partial\Omega)^*} 
\,\text{ for each }\,u\in L^2(\Omega)\,\text{ with }\,(-\Delta+V)u=0\,\text{ in }\,\Omega.
\end{equation}
Upon recalling that $\widetilde{\gamma}_D$ is compatible with the ordinary Dirichlet trace $\gamma_D$ from 
\eqref{7yfFc-P} and keeping in mind the identifications in \eqref{Tan-C5.WACO}, these considerations provide a satisfactory answer to G.\ Uhlmann's question formulated above. 
The subtle aspect in this context  is that while measuring the size of the 
Dirichlet trace in the space $H^{-1/2}(\partial\Omega)$ is inadequate within the class of Lipschitz domains, the correct substitute which does the job is precisely our space ${\mathscr{G}}_N(\partial\Omega)^*$. 

In addition, similar results are valid for our generalized Neumann trace operator $\widetilde{\gamma}_N$ (cf.~\eqref{4.55.INTR}, 
\eqref{eqn:bvpNeumann.INTR.WACO}), namely, whenever $L^\infty(\Omega)\ni V\geq 0$ at a.e.~point in the bounded Lipschitz domain $\Omega\subset{\mathbb{R}}^n$, one has 
\begin{equation}\label{eqn:bvp.2.INTR.WACO.NEU} 
\|u\|_{L^2(\Omega)}\leq C\|\widetilde{\gamma}_N u\|_{\mathscr{G}_D(\partial\Omega)^*} 
\,\text{ for each }\,u\in L^2(\Omega)\,\text{ with }\,(-\Delta+V)u=0\,\text{ in }\,\Omega. 
\end{equation}
In particular, 
\begin{equation}\label{eqn:bvp.2.INTR.WACO.NEU.2} 
\|u\|_{L^2(\Omega)}\leq C\|\widetilde{\gamma}_N u\|_{\mathscr{G}_D(\partial\Omega)^*} 
\,\text{ for each harmonic function }\,u\in L^2(\Omega), 
\end{equation}
which may be regarded as the analogue of G.\ Uhlmann's question for the Neumann trace operator. 

Moreover, in Section~\ref{s9} we rely on the power of Theorem~\ref{t5.5.INTR} to describe the 
Krein--von Neumann extensions of Schr\"odinger operators on bounded Lipschitz domains. 
Our main result in this regard is Theorem~\ref{t5.6}, stating that if $\Omega\subset\bbR^n$ is a bounded Lipschitz domain, and if the potential $V\in L^\infty(\Omega)$ is real-valued a.e., then the  Krein--von Neumann 
extension $A_{K,\Omega}$ of $A_{min,\Omega}$ (the minimal realization of $-\Delta+V$, defined 
as the closure in $L^2(\Omega)$ of $-\Delta+V$ acting from $C^\infty_0(\Omega)$) is given by 
\begin{align}\label{eqn:A_K.INTR}
\begin{split} 
& A_{K,\Omega}=-\Delta+V,    
\\[2pt] 
& \dom(A_{K,\Omega})=\big\{f\in\dom(A_{max,\Omega})\,\big|\, 
\widetilde{\gamma}_N f+\widetilde M_{\Omega}(0)\widetilde{\gamma}_D f=0\big\},
\end{split} 
\end{align}
where $\widetilde{\gamma}_D,\widetilde{\gamma}_N$ are the maximal extensions of the Dirichlet and Neumann 
trace operators defined as in \eqref{4.55.INTR}, and where $\widetilde{M}_{\Omega}(\cdot)$ is (up to a sign) 
a spectral parameter dependent extended Dirichlet-to-Neumann map, or Weyl--Titchmarsh operator for the 
for Schr\"{o}dinger operator (cf. the discussion in Section~\ref{s7}). 

The concrete description of $\dom(A_{K,\Omega})$ in \eqref{eqn:A_K.INTR} has the distinct advantage of 
making explicit the underlying boundary condition. Nonetheless, as opposed to the classical Dirichlet and 
Neumann boundary condition, this boundary condition is {\it nonlocal} in nature, as it involves
$\widetilde{M}_{\Omega}(\cdot)$. When $\Om$ is smooth and $V=0$, $\widetilde{M}_{\Omega}(\cdot)$
is a boundary pseudodifferential operator of order $1$, and \eqref{eqn:A_K.INTR} becomes 
the appropriate rigorous interpretation in a very general geometric setting of the informal philosophy, outlined 
by A.\ Alonso and B.\ Simon in \cite{AS80}, asserting that the Krein Laplacian is realization of the Laplacian with 
the non-local boundary condition 
\begin{equation}\label{bcK}
\partial_\nu f=\partial_\nu H(f)\,\text{ on }\,\partial\Omega,     
\end{equation}  
where $\partial_\nu=\nu\cdot\nabla$, with $\nu$ denoting the outward unit normal to $\Omega$, 
is the normal directional derivative and, given a sufficiently nice function $f$ in $\Om$, 
the symbol $H(f)$ denotes the harmonic extension to $\Omega$ of the trace of $f$ on ${\partial\Omega}$.
Near the end of their paper \cite{AS80}, A.\ Alonso and B.\ Simon also raise the following issue:
\begin{quote}
{\it ``It seems to us that the Krein extension of $-\Delta$, that is, $-\Delta$ with the 
boundary condition $\eqref{bcK}$, is a natural object and therefore worthy of further study. 
For example: Are the asymptotics of its nonzero eigenvalues given by Weyl's formula?''}
\end{quote}
In the case where $\Omega$ is bounded and $C^{\infty}$-smooth, and 
$V\in C^{\infty}(\overline{\Omega})$, this has been shown to be the case three years 
later by G.\ Grubb \cite{Gr83}. More specifically, in \cite{Gr83} Grubb has proved that if 
$N(\lambda,A_{K,\Omega})$ denotes the number of nonzero eigenvalues $\lambda_j$ of 
$A_{K,\Omega}$ not exceeding $\lambda\in\bbR$,
\begin{equation}\label{1.6}
N(\lambda,A_{K,\Omega}):=\#\{j\in\bbN\,|\,0<\lambda_j\leq\lambda\},\quad\forall\,\lambda\in\bbR, 
\end{equation}
then
\begin{align}\label{Df-H2}
\begin{split} 
& \Omega\in C^\infty\,\text{ and }\,V\in C^{\infty}(\overline{\Omega})
\,\text{ imply }    
\\[2pt]
& \quad N(\lambda,A_{K,\Omega})\underset{\lambda\to\infty}{=}(2\pi)^{-n}v_n|\Omega|\,\lambda^{n/2} 
+O\big(\lambda^{(n-\theta)/2}\big), 
\end{split} 
\end{align} 
where 
\begin{equation}\label{Df-H4}
\theta:=\max\,\big\{\tfrac{1}{2}-\varepsilon\,,\,\tfrac{2}{n+1}\big\},
\,\text{ with $\varepsilon>0$ arbitrary}.
\end{equation} 
In fact, Grubb considers the case of strongly elliptic differential operators 
of order $2m$, $m\in\bbN$, strictly positive, with smooth coefficients, though 
we here restrict our discussion to the case $m=1$. The methods used by Grubb rely 
on pseudo-differential operator techniques (which are not applicable to the minimally 
smooth case we are aiming at in this work). See also \cite{Mi94}, \cite{Mi06}, and 
most recently, \cite{Gr12}, where the authors  derive a sharpening of the remainder 
in \eqref{Df-H2} to any $\theta<1$. 

To prove \eqref{Df-H2}--\eqref{Df-H4}, Grubb showed that the eigenvalue problem
\begin{equation}\label{Df-H7}
(-\Delta+V)f=\lambda\,f,\quad f\in\dom(A_{K,\Omega}),\,\,\lambda>0,
\end{equation} 
is spectrally equivalent to the following fourth-order pencil eigenvalue problem
\begin{align}\label{Df-H8} 
\begin{split}
& (-\Delta+V)^2 w=\lambda\,(-\Delta+V)w\,\text{ in }\,\Om,    
\\[2pt]
& \quad w\in\dom\big((-\Delta_{max,\Omega})(-\Delta_{min,\Omega})\big),\,\,\lambda>0. 
\end{split}
\end{align}
This is closely related to the so-called problem of the {\it buckling of a clamped plate}, 
\begin{equation}\label{Df-H8F}
-\Delta^2 w=\lambda\,\Delta w\,\text{ in }\,\Om,
\quad w\in\dom\big((-\Delta_{max,\Omega})(-\Delta_{min,\Omega})\big),\,\,\lambda>0,
\end{equation}
to which \eqref{Df-H8} reduces when $V\equiv 0$. In particular, this permits one to 
allude to the theory of generalized eigenvalue problems, that is, operator pencil problems 
of the form $Tu=\lambda\,Su$, where $T$ and $S$ are linear operators in a Hilbert space. 
However, given the present low regularity assumptions (cf. \eqref{1.7Mba}--\eqref{1.7} below) 
we find it more convenient to appeal to a version of this pencil problem which 
emphasizes the role of the following symmetric forms in $L^2(\Omega)$, 
\begin{align}\label{1.11}
\mathfrak{a}_{K,\Omega}(f,g) &:=\big((-\Delta +V)f,(-\Delta+V)g\big)_{L^2(\Omega)}, 
\quad\forall\,f,g\in\accentset{\circ}{H}^2(\Omega),      
\\[2pt] 
\mathfrak{b}_{K,\Omega}(f,g) &:=(\nabla f,\nabla g)_{L^2(\Omega)^n}
+\big(V^{1/2}f,V^{1/2}g\big)_{L^2(\Omega)},\quad\forall\,f,g\in\accentset{\circ}{H}^2(\Omega),    
\label{1.12}
\end{align}
and hence focus on the problem of finding $f\in\accentset{\circ}{H}^2(\Omega)$ satisfying 
\begin{equation}\label{Xkko-14}
\mathfrak{a}_{K,\Omega}(f,g)=\lambda\,\mathfrak{b}_{K,\Omega}(f,g), 
\quad\forall\,g\in\accentset{\circ}{H}^2(\Omega).
\end{equation}
This type of eigenvalue problem, in the language of bilinear forms associated
with differential operators, has been studied by V.\ A.\ Kozlov in a series of papers 
\cite{Ko79}, \cite{Ko83}, \cite{Ko84}. In particular, in \cite{Ko84}, Kozlov 
has obtained Weyl asymptotic formulas for \eqref{Xkko-14} in the case where the underlying
domain $\Omega$ is merely Lipschitz and $V\in L^{\infty}(\Omega)$. 

For rough domains $\Omega$, matters are much more delicate as the nature of the
boundary trace operators and the standard elliptic regularity theory are both 
fundamentally affected. Following work in \cite{GM11}, the class of 
{\it quasi-convex domains} was considered in great detail in \cite{AGMT10}. The latter is
a subclass of bounded, Lipschitz domains in $\bbR^n$ where only singularities pointing 
in the outward direction are permitted. For example, the class of 
of quasi-convex domains includes all bounded (geometrically) convex domains,  
all bounded Lipschitz domains satisfying a uniform exterior ball condition 
(which, informally speaking, means that a ball of fixed radius can be ``rolled'' 
along the boundary), and all bounded domains of class $C^{1,r}$ for some $r>1/2$. 
One of the key features of this class of quasi-convex domains is the fact that the classical
elliptic regularity property
\begin{equation}\label{Df-H1}
\dom(A_{D,\Om})\subset H^2(\Om),\quad\dom(A_{N,\Om})\subset H^2(\Om), 
\end{equation} 
holds (this property, however, is known to fail for general bounded Lipschitz domains;
for example, work in \cite{Da79} imply the existence of a bounded Lipschitz domain $\Omega$ 
and $f\in\dom(A_{D,\Om})$ with second-order derivatives not in $L^p(\Omega)$ for any $p>1$).
It was recognized in \cite{AGMT10} that Kozlov's analysis can be applied to the spectral 
asymptotics of perturbed Krein Laplacians. The main result proved in \cite{AGMT10} then 
established the Weyl-type spectral asymptotics 
\begin{equation}\label{1.16}
N(\lambda,A_{K,\Omega})\underset{\lambda\to\infty}{=}(2\pi)^{-n}v_n|\Omega|\,\lambda^{n/2} 
+O\big(\lambda^{(n-(1/2))/2}\big)      
\end{equation}
for the Krein--von Neumann extension, denoted by $A_{K,\Omega}$, of the perturbed Laplacian 
$(-\Delta+V)|_{ C^\infty_0(\Omega)}$ in the case where $0\leq V\in L^\infty(\Omega)$ 
and $\Omega\subset\bbR^n$ is a quasi-convex domain. 

Another principal goal of the current work is to take the final step in this development and 
prove the Weyl-type spectral asymptotics \eqref{1.16} for $A_{K,\Omega}$ in the case where again 
\begin{equation}\label{1.7Mba}
0\leq V\in L^\infty(\Omega), 
\end{equation} 
and 
\begin{equation}\label{1.7}
\text{$\Omega\subset\bbR^n$ is a bounded Lipschitz domain.}   
\end{equation} 
We emphasize that the potential coefficient $V$ is permitted to be nonsmooth and that
the underlying domain $\Omega$ is allowed to have irregularities of a more general 
nature than the class of quasi-convex domains discussed above. The methods employed in 
this work rely on the spectral equivalence to the underlying buckling problem (see \cite{AGMST10} 
for an abstract approach), 
on the use of spectral parameter dependent Dirichlet-to-Neumann map (the Weyl--Titchmarsh operator), 
and on appropriate Gelfand triples defined in terms of the Dirichlet and Neumann boundary trace maps. 
What underpins this entire approach is a sharp boundary trace theory, that continues to 
be effective outside of the traditional settings. 

Indeed, one of the challenges in the nonsmooth setting considered here pertains to the lack 
of $H^2(\Omega)$-regularity \eqref{Df-H1}, which will be replaced by $H^{3/2}(\Omega)$-regularity. 
It has long been understood that this regularity issue is intimately linked to the analytic and geometric 
properties of the underlying domain $\Omega$. To illustrate this point, we briefly consider the case 
when $\Omega\subset{\mathbb{R}}^2$ is a polygonal domain with at least one re-entrant corner. 
In this scenario, let $\omega_1,\dots,\omega_N$ be the internal angles of $\Omega$ satisfying 
$\pi<\omega_j<2\pi$, $1\leq j\leq N$, and denote by $P_1,\dots,P_N$ the corresponding vertices. 
Then (cf., e.g., \cite{KMR97}) the structure of a generic function $u$ belonging to $\dom(-\Delta_{D,\Om})$ is 
\begin{equation}\label{bf-G2}
u=\sum_{j=1}^N\lambda_jv_j+w,\text{ for some }\,\lambda_j\in{\mathbb{R}},\;1\leq j\leq N,
\end{equation}
where $w\in H^{2}(\Omega)\cap\accentset{\circ}{H}^1(\Omega)$ and, for each $j\in\{1,\dots,N\}$, 
the function $v_j$ exhibits a singular behavior at the vertex $P_j$ of the following nature. 
Given $j\in\{1,\dots,N\}$, choose polar coordinates $(r_j,\theta_j)$ taking $P_j$ 
as the origin and so that the internal angle is spanned by the half-lines $\theta_j=0$ 
and $\theta_j=\omega_j$. Then 
\begin{equation}\label{bf-G3}
v_j(r_j,\theta_j)=\phi_j(r_j,\theta_j)\,r_j^{\pi/\omega_j}\sin(\pi\theta_j/\omega_j),
\quad 1\leq j\leq N,
\end{equation}
where $\phi_j$ is a $C^\infty$ cut-off function of small support, which is identically one 
near $P_j$. In this scenario, $v_j\in H^s(\Omega)$ for every $s<1+(\pi/\omega_j)$, though 
$v_j\notin H^{1+(\pi/\omega_j)}(\Omega)$ (see Proposition~\ref{CC.WACO} in this regard). 
This analysis implies that the best regularity statement regarding a generic function 
$u\in\dom(A_{D,\Om})$ is 
\begin{equation}\label{bf-G4}
u\in H^s(\Omega)\,\text{ for every }\,s<1+\frac{\pi}{\max\,\{\omega_1,\dots,\omega_N\}}
\end{equation}
and this membership fails for the above critical value of $s$. We note that 
\begin{equation}\label{bf-G4.TT}
1+\frac{\pi}{\max\,\{\omega_1,\dots,\omega_N\}}\in\big(3/2,2\big)
\end{equation}
and, in particular, this provides a geometrically quantifiable way of measuring the failure of the 
inclusion $\dom(A_{D,\Om})\subset H^2(\Om)$ in \eqref{Df-H1} even for piecewise $C^\infty$-domains 
exhibiting inwardly directed irregularities. This being said, from \eqref{bf-G4}--\eqref{bf-G4.TT} 
(and a similar type of analysis corresponding to Neumann boundary conditions) we do have 
\begin{equation}\label{eq:4f4f}
\dom(A_{D,\Om})\subset H^{3/2}(\Om),\quad\dom(A_{N,\Om})\subset H^{3/2}(\Om)
\end{equation}
for this type of domains, and the exponent $3/2$ is sharp. We shall see later that this 
sharp regularity result holds in the more general class of arbitrary bounded Lipschitz domains.
The fact that \eqref{Df-H1} downgrades, in the said class of domains, to just \eqref{eq:4f4f} 
creates significant difficulties as, for example, the Dirichlet boundary trace operator fails 
to map $H^{3/2}(\Om)$ into $H^1(\partial\Omega)$. One of the key ingredients in dealing with 
\eqref{eq:4f4f} in lieu of \eqref{Df-H1} is devising a boundary trace theory which, in addition 
to making optimal use of the regularity (measured on the scale of Sobolev spaces) exhibited by 
functions belonging to $\dom(A_{D,\Om})$ and $\dom(A_{N,\Om})$, also employs the 
PDE aspect inherent to a such membership. See Theorem~\ref{YTfdf-T}, Theorem~\ref{YTfdf.NNN.2-Main}, 
and Theorem~\ref{t5.5} in this regard, which rely heavily on the theory of boundary value 
problem for the Laplacian in Lipschitz domains developed in \cite{FMM98}, \cite{JK95}, 
\cite{MT99}--\cite{MT01}.

Yet another fundamental application of Theorem~\ref{t5.5.INTR} is the classification of all 
self-adjoint extensions of the minimal Schr\"odinger operator in an arbitrary bounded Lipschitz 
domain $\Omega\subset{\mathbb{R}}^n$. The aforementioned family is parametrized in terms of closed 
subspaces $\mathscr{X}\subset\mathscr{G}_N(\partial\Omega)^*$ and self-adjoint operators 
$T:\mathscr{X}\supset\dom(T)\rightarrow\mathscr{X}^*$ in the manner described in Theorem~\ref{jussisaext}.
Specifically, for every closed subspace $\mathscr{X}\subset\mathscr{G}_N(\partial\Omega)^*$ and every 
self-adjoint operator $T:\mathscr{X}\supset\dom(T)\rightarrow\mathscr{X}^*$ the operator
\begin{equation}\label{atjussi.INTRO}
\begin{split}
& A_{T,\Omega}=-\Delta+V,
\\[2pt]
& \dom(A_{T,\Omega})=\big\{f\in\dom(A_{max,\Omega})\,\big|\,
T\,\widetilde\gamma_D f=P_{\mathscr{X}^*}\gamma_N f_D\big\}
\end{split}
\end{equation}
is a self-adjoint extension of $A_{min,\Omega}$ in $L^2(\Omega)$, where $P_{\mathscr{X}^*}$ denotes 
the orthogonal projection in $\mathscr{G}_N(\partial\Omega)$ onto $\mathscr{X}^*$ (cf.~\eqref{jussipq}) 
and, for some fixed $\mu\in\rho(A_{D,\Omega})\cap\mathbb R$, we have decomposed (see \eqref{fdeco0}) 
each $f\in\dom(A_{max,\Omega})$ as
\begin{equation}\label{fdeco.INTRO}
f=f_D+f_\mu\,\text{ with }\,f\in\dom(A_{D,\Omega})\,\text{ and }\,f_\mu\in\ker(A_{max,\Omega}-\mu). 
\end{equation}
Conversely, for every self-adjoint extension $A$ of $A_{min,\Omega}$ in $L^2(\Omega)$ there exists a 
closed subspace $\mathscr{X}\subset\mathscr{G}_N(\partial\Omega)^*$ and a self-adjoint operator 
$T:\mathscr{X}\supset\dom(T)\rightarrow\mathscr{X}^*$ such that $A=A_{T,\Omega}$, that is,
\begin{equation}\label{rafc.INTRO.2}
\begin{split}
& A=-\Delta+V,
\\[2pt]
& \dom(A)=\big\{f\in\dom(A_{max,\Omega})\,\big|\,T\,\widetilde\gamma_D f=P_{\mathscr{X}^*}\gamma_N f_D\big\}.
\end{split}
\end{equation}
A key feature of this result is the fact that all said extensions are characterized via explicit boundary conditions. 
Of course, the Dirichlet and Neumann self-adjoint realizations of $-\Delta+V$ are among these, but the said family 
also includes self-adjoint realizations of the Schr\"odinger operator with exotic boundary conditions of a non-local nature, 
as in the case of the Krein--von Neumann extension $A_{K,\Omega}$ of $A_{min,\Omega}$ described in \eqref{eqn:A_K.INTR}.
This provides a most satisfactory answer to a problem that has been investigated for more than 60 years in the mathematical 
literature (starting with the pioneering works of M.~I.~Vi{\u s}ik and G.~Grubb). In addition, this extends and unifies 
fundamental results going back to J.\ L.\ Lions and E.\ Magenes, as well as D.\ Jerison and C.\ Kenig.

Finally, in Section~\ref{s11} we initiate a treatment of variable coefficient second-order elliptic operators 
(in place of the ordinary Laplacian). While this topic is worth pursuing further, here we lay the foundations 
by demonstrating how the bulk of the material in Sections~\ref{s2}--\ref{s10} extends to the Laplace--Beltrami 
operator (perturbed by a scalar potential $V$) on a compact boundaryless Riemannian manifold. 

Our principal result in Section~\ref{s11} is the version of Theorem~\ref{t5.5.INTR} in the aforementioned geometric setting
(see Theorem~\ref{t5.5.MAN}). In Subsection \ref{ss.MM.EEE} we also indicate how to recast such results in the language of 
ordinary (Euclidean) elliptic differential operators with variable coefficients, of class $C^{1,1}$, on the closure 
a bounded Lipschitz domain $\Omega\subset{\mathbb{R}}^n$. A benefit of developing the aforementioned machinery 
for the Laplace--Beltrami operator on Riemannian manifolds is that we may painlessly reformulate results proved earlier 
in Subsections~\ref{ss.MM.1}--\ref{ss.MM.3} in the language of variable-coefficient differential operators. 
Given their intrinsic importance, we close Section~\ref{s11} by elaborating on the variable-coefficient versions 
of our earlier Euclidean trace results (from Theorem~\ref{YTfdf-T}, Theorem~\ref{GNT}, and Theorem~\ref{YTfdf.NNN.2-Main}) 
in Theorem~\ref{YTfdf-T-MMM.Euclid} and Corollary~\ref{YTfdf-T.NNN-MMM.Euclid} for the Dirichlet trace, 
and in Theorem~\ref{GNT.Euclid} and Corollary~\ref{YTfdf.NNN.2-MMM.CCC} for the Neumann trace.

\vspace{10pt} 

The layout of the manuscript is as follows. Section~\ref{s2} is devoted to Sobolev and Besov spaces on Lipschitz domains. 
After a thorough review of Lipschitz domains $\Omega\subset\bbR^n$, and nontangential maximal functions we turn 
to fractional Sobolev and Besov spaces on $\Omega$ and $\partial\Omega$. In Section~\ref{s3} we take up the task 
of developing, in a systematic manner, a sharp Dirichlet boundary trace theory in bounded Lipschitz domains in 
$\bbR^n$ involving Sobolev and Besov spaces that is particularly well-suited for the goals we have in mind in 
this work. Our main results there are Theorems~\ref{YTfdf-T}, \ref{YTfdf-T-Mi} with a brand of Dirichlet boundary 
trace operators which continue to remain meaningful in limiting cases when their ordinary versions fail to apply. 
Section~\ref{s4} employs the Dirichlet boundary trace operator introduced in Section~\ref{s3} to derive far-reaching 
divergence theorems culminating in Theorem~\ref{Ygav-75-BIS-BBB}. Given Sections~\ref{s3} and \ref{s4} we are in 
position to develop a sharp Neumann boundary trace theory on bounded Lipschitz domains in $\bbR^n$ involving Sobolev 
spaces, the principal result on the weak boundary trace map being recorded in Theorem~\ref{YTfdf.NNN.2-Main}. 
Section~\ref{s6} discusses Schr\"odinger operators and their Dirichlet and Neumann realizations (also, the Friedrichs 
extension of an appropriate minimal Schr\"odinger operator realization) in arbitrary nonempty open sets $\Omega\subseteq\bbR^n$ 
as well as on bounded Lipschitz domains. Section~\ref{s7} is devoted to a study of Weyl--Titchmarsh operators $M_{\Omega}(\cdot)$, 
that is, spectral parameter dependent Dirichlet-to-Neumann maps, associated with Schr\"odinger operators on bounded 
Lipschitz domains. The principal objective of Section~\ref{s8} is to extend the Dirichlet and Neumann traces by 
continuity onto the domain of the underlying maximal Schr\"odinger operator on bounded Lipschitz domains. 
The Krein--von Neumann extension of Schr\"odinger operators on bounded Lipschitz domains is the principal 
object of Section~\ref{s9}. We identify the nonlocal boundary condition characterizing the perturbed 
Krein Laplacian in terms of an appropriate extension of $M_{\Omega}(0)$, and invoking the spectral equivalence 
between the buckling problem (with potential $V$) and the perturbed Krein Laplacian, and, with the help of 
Kozlov's analysis of Weyl asymptotics for the buckling problem on Lipschitz domains, we derive the Weyl spectral 
asymptotics for the perturbed Krein Laplacian in bounded Lipschitz domains in Theorem~\ref{t6.2}. A description of all 
self-adjoint extensions of the minimal Schr\"odinger operator and Krein-type resolvent formulas in connection with 
bounded Lipschitz domains are the subject of Section~\ref{s10}. Our final Section~\ref{s11} offers a glimpse at the 
case of variable coefficient operators and treats Laplace--Beltrami operators perturbed by scalar potentials on 
boundaryless Riemannian manifolds. This section (a substantial one), initiates such a treatment and points the way 
to future research in this direction. In Section~\ref{s11} we also present variable-coefficient versions 
of our earlier Euclidean trace results. 

\vspace{10pt} 

We conclude this introduction by summarizing the notation used in this work. Throughout, the symbol 
$\cH$ is reserved to denote a separable complex Hilbert space with $(\dott,\dott)_{\cH}$ the scalar product in 
$\cH$ (linear in the second argument), and $I_{\cH}$ the identity operator in $\cH$. Next, let $T$ be a linear 
operator mapping (a subspace of) a Banach space into another, with $\dom(T)$ and $\ran(T)$ denoting the domain 
and range of $T$. The closure of a closable operator $S$ is denoted by $\ol S$. The kernel (null space) of $T$ 
is denoted by $\ker(T)$. The spectrum, point spectrum (i.e., the set of eigenvalues), discrete spectrum, 
essential spectrum, and resolvent set of a closed linear operator in $\cH$ will be denoted by $\sigma(\cdot)$, 
$\sigma_{p}(\cdot)$, $\sigma_{d}(\cdot)$, $\sigma_{ess}(\cdot)$, and $\rho(\cdot)$, respectively. The symbol 
$\slim$ abbreviates the limit in the strong (i.e., pointwise) operator topology. 

The Banach space of bounded linear operators on $\cH$ is denoted by $\cB(\cH)$. 
The analogous notation $\cB(\cX_1,\cX_2)$ will be used for bounded operators between two Banach spaces 
$\cX_1$ and $\cX_2$. Moreover, $\cX_1\hookrightarrow\cX_2$ denotes the continuous embedding of the Banach space 
$\cX_1$ into the Banach space $\cX_2$. In addition, $U_1\dotplus U_2$ denotes the direct sum of the subspaces 
$U_1$ and $U_2$ of a Banach space $\cX$; and $V_1\oplus V_2$ represents the 
orthogonal direct sum of the subspaces $V_1$ and $V_2$ of a Hilbert space $\cH$. 

Given a Banach space $X$, we let $X^*$ denote the {\it adjoint space} of continuous 
conjugate linear functionals on $X$, that is, the {\it conjugate dual space} of $X$ (rather than the usual 
dual space of continuous linear functionals on $X$). This avoids the well-known awkward distinction between 
adjoint operators in Banach and Hilbert spaces (cf., e.g., the pertinent discussion in \cite[pp.~3--4]{EE18}). 

The symbol $L^2(\Omega)$, with $\Omega\subseteq\bbR^n$ open, $n\in\bbN\backslash\{1\}$, is a shortcut for 
$L^2(\Omega,d^n x)$, whenever the $n$-dimensional Lebesgue measure is understood. (For simplicity we exclude 
the one-dimensional case $n=1$ in this work as the case $\Omega=(a,b)\subset\bbR$ has been treated in detail 
in \cite[Section~10.1]{AGMT10}.) Moreover, if $\Omega$ is a Lipschitz domain in $\bbR^n$, $L^2(\partial\Omega)$ 
represents the Lebesgue space of square integrable functions with respect to the canonical surface measure 
on $\partial\Omega$. For brevity, the identity operator in $L^2(\Omega)$ and $L^2(\partial\Omega)$ will typically 
be denoted by $I$ if no confusion can arise. The symbol $\cD(\Omega)$ is reserved for the set of test functions 
$C_0^{\infty}(\Omega)$ on $\Omega$, equipped with the standard inductive limit topology, and $\cD'(\Omega)$ represents 
its dual space, the set of distributions in $\Omega$. In addition, $\bbC_+$ (resp., $\bbC_-$) denotes the open complex 
upper (resp., lower) half-plane, while $\#(M)$ abbreviates the cardinality of the set $M$. We agree to define 
$\bbN_0:=\bbN\cup\{0\}$, so that $\bbN_0^n$ becomes the collection of all multi-indices with $n$ components. 
As is customary, for each $\alpha=(\alpha_1,\dots,\alpha_n)\in{\mathbb{N}}_0^n$ we denote by 
$|\alpha|:=\alpha_1+\ldots+\alpha_n$ the length of $\alpha$. Also, we shall let 
$\bbS^{n-1}:=\{x\in{\mathbb{R}}^n|\,|x|=1\}$ denote the unit sphere in ${\mathbb{R}}^n$ centered at the origin. 

We shall often use the common convention of denoting by the same letter $C$ possibly different multiplicative 
constants in various inequalities throughout the monograph. Moreover, writing ``$A(x)\approx B(x)$ uniformly in $x$''
signifies the existence of some number $C\in(1,\infty)$ which is independent of $x$ with the property that 
$A(x)/C\leq B(x)\leq CA(x)$ for every $x$.

Finally, a notational comment: For obvious reasons, which have their roots in quantum mechanics, we will, 
with a slight abuse of notation, dub the expression $-\Delta=-\sum_{j=1}^n\partial_j^2$ (rather than $\Delta$) 
as the ``Laplacian'' in this work. When acting on vector-valued functions (or distributions), the Laplacian 
is considered componentwise.

%%%%%%%%%%%%%%%%%%%%%%%%%%%%%%
%%%%%%%%%%%%%%%%%%%%%%%%%%%%%%
\section{Sobolev and Besov Spaces on Lipschitz Domains}  
\label{s2} 
%%%%%%%%%%%%%%%%%%%%%%%%%%%%%%
%%%%%%%%%%%%%%%%%%%%%%%%%%%%%%

In this section we recall a variety of background material including, a thorough review of Lipschitz domains 
in $\bbR^n$, nontangential maximal functions, fractional Sobolev and Besov spaces on arbitrary open sets and 
on bounded Lipschitz domains in $\bbR^n$, as well as on the boundaries of bounded Lipschitz domains, and 
Sobolev regularity in terms of nontangential maximal functions. 

%%%%%%%%%%%%
\subsection{The class of Lipschitz domains}\label{ss2.1}
%%%%%%%%%%%%

The reader is reminded that a function (acting between two metric spaces) is called 
Lipschitz if it does not distort distances by more than a fixed multiplicative constant. 
We begin by giving the formal definition of the category of Lipschitz domains 
(cf., e.g., \cite{MMM22}, for more on this topic). 

%%%%%%%
\begin{definition}\label{Def-LiP}
Let $\Omega$ be a nonempty, proper, open subset of ${\mathbb{R}}^n$. \\[1mm] 
$(i)$ Call $\Omega$ a {\tt Lipschitz domain near} 
$x_0\in\partial\Omega$ if there exist $r,\tau\in(0,\infty)$ with the following significance. 
For some choice of an $(n-1)$-dimensional plane $H\subseteq{\mathbb{R}}^n$ passing through 
the point $x_0$, some choice of a unit normal vector $N$ to $H$, the cylinder 
${\mathcal{C}}_{r,\tau}(x_0,N):=\{x'+tN\,|\,x'\in H,\,|x'-x_0|<r,\,\,|t|<\tau\}$ 
{\rm (}called coordinate cylinder near $x_0${\rm )} has the property that
\begin{align}\label{Def-Lip}
{\mathcal{C}}_{r,\tau}(x_0,N)\cap\Omega &={\mathcal{C}}_{r,\tau}(x_0,N)
\cap\{x'+tN\,|\,x'\in H,\,\,t>\varphi(x')\}
\nonumber\\[2pt]
&=\{x'+tN\,|\,x'\in H,\,\,|x'-x_0|<r,\,\,t\in(\varphi(x'),\tau)\},
\end{align}
for some Lipschitz function $\varphi:H\to{\mathbb{R}}$ {\rm (}called the defining function 
for $\partial\Omega$ near $x_0${\rm )}, satisfying 
\begin{equation}\label{Def-Lip3}
\varphi(x_0)=0\,\text{ and }\,|\varphi(x')|<\tau\,\text{ if }\,|x'-x_0|\leq r. 
\end{equation}
Collectively, the pair $\big\{{\mathcal{C}}_{r,\tau}(x_0,N),\varphi\big\}$ will 
be referred to as a local chart near $x_0$, whose geometrical characteristics 
consist of $r$, $\tau$, and the Lipschitz constant of $\varphi$. \\[1mm] 
$(ii)$ Call $\Omega$ a {\tt locally Lipschitz domain} 
if it is a Lipschitz domain near every point $x\in\partial\Omega$. \\[1mm] 
$(iii)$ Call $\Omega$ a {\tt Lipschitz domain} if $\Omega$ is a locally Lipschitz domain 
and at each boundary point there exists a local chart whose geometrical characteristics 
are independent of the point in question {\rm (}collectively, the said geometrical characteristics
are going to be referred to in the future as the {\tt Lipschitz character} of $\Omega${\rm )}. 
\\[1mm] 
$(iv)$ The category of $C^{k}$-domains with $k\in{\mathbb{N}}\cup\{\infty\}$ 
is defined analogously, requiring that the defining functions $\varphi$ 
are of class $C^k$. 
\end{definition} 
%%%%%%%

We emphasize that no topological conditions are placed on the class of bounded Lipschitz domains 
considered here; in particular, the boundaries of the domains in question are allowed to be disconnected. 

A few useful observations related to the property of an open set $\Omega\subseteq{\mathbb{R}}^n$ 
of being a Lipschitz domain near a point $x_0\in\partial\Omega$ are collected in the lemma below
(proved in \cite[Proposition~2.8]{ABMMZ11}). The reader is reminded that the complement of a set 
$E\subseteq{\mathbb{R}}^n$, relative to ${\mathbb{R}}^n$, is denoted by $E^c:={\mathbb{R}}^n\backslash E$. 
In addition, by $E^\circ$ and $\overline{E}$ we shall denote the interior and closure of $E$ in the standard 
topology of ${\mathbb{R}}^n$, respectively. 

%%%%%%%
\begin{lemma}\label{JGbv.2}
Assume that $\Omega$ is a nonempty, proper, open subset of ${\mathbb{R}}^n$, 
and fix some point $x_0\in\partial\Omega$. \\[1mm] 
$(i)$
If $\Omega$ is a Lipschitz domain near $x_0$ and if $\big\{{\mathcal{C}}_{r,\tau}(x_0,N),\varphi\big\}$ 
is a local chart near $x_0$ {\rm (}in the sense of Definition~\ref{Def-LiP}{\rm )} then, in addition to 
\eqref{Def-Lip}, one also has
\begin{align}\label{Def-Lip1}
& {\mathcal{C}}_{r,\tau}(x_0,N)\cap\partial\Omega={\mathcal{C}}_{r,\tau}(x_0,N)
\cap\{x'+tN\,|\,x'\in H,\,\,t=\varphi(x')\}, 
\\[2pt]
& {\mathcal{C}}_{r,\tau}(x_0,N)\cap(\overline{\Omega})^c
={\mathcal{C}}_{r,\tau}(x_0,N)\cap\{x'+tN\,|\,x'\in H,\,\,t<\varphi(x')\}.
\label{Def-Lip2}
\end{align}
Furthermore, 
\begin{align}\label{DR-a1}
& {\mathcal{C}}_{r,\tau}(x_0,N)\cap\overline{\Omega}={\mathcal{C}}_{r,\tau}(x_0,N)
\cap\{x'+tN\,|\,x'\in H,\,\,t\geq\varphi(x')\}, 
\\[2pt]
& {\mathcal{C}}_{r,\tau}(x_0,N)\cap(\overline{\Omega})^\circ
={\mathcal{C}}_{r,\tau}(x_0,N)\cap\{x'+tN\,|\,x'\in H,\,\,t>\varphi(x')\}.
\label{DR-a2}
\end{align}
$(ii)$ Suppose there exist an $(n-1)$-dimensional plane $H\subseteq{\mathbb{R}}^n$ 
passing through the point $x_0$, a choice of a unit normal vector $N$ to $H$, an open 
cylinder ${\mathcal{C}}_{r,\tau}(x_0,N)=\{x'+tN\,|\,x'\in H,\,|x'-x_0|<r,\,\,|t|<\tau\}$ 
and a Lipschitz function $\varphi:H\to{\mathbb{R}}$ satisfying \eqref{Def-Lip3}
such that \eqref{Def-Lip1} holds. Then, assuming $x_0\notin(\overline{\Omega})^{\circ}$, 
it follows that $\Omega$ is a Lipschitz domain near $x_0$. 
\end{lemma}
%%%%%%%

Definition~\ref{Def-LiP} and item $(i)$ in Lemma~\ref{JGbv.2} show that if $\Omega\subseteq{\mathbb{R}}^n$ 
is a Lipschitz domain near a boundary point $x_0$ then, in a neighborhood of $x_0$, the topological boundary 
$\partial\Omega$ agrees with the graph of a Lipschitz function $\varphi:{\mathbb{R}}^{n-1}\to{\mathbb{R}}$, 
considered in a suitably chosen system of coordinates (which is isometric with the original one). 
Then the outward unit normal $\nu=(\nu_1,\nu_2,\dots,\nu_n)$ to $\Omega$ has an explicit formula 
in terms of $\nabla'\varphi$, the $(n-1)$-dimensional gradient of $\varphi$. Specifically, if 
${\mathscr{H}}^{n-1}$ stands for the $(n-1)$-dimensional Hausdorff measure in $\bbR^n$, 
then in the new system of coordinates we have
\begin{align}\label{Nu-C2}
\nu\big(x',\varphi(x')\big) 
&=\frac{\big((\partial_1\varphi)(x'),\ldots,(\partial_{n-1}\varphi)(x'),-1\big)}{\sqrt{1+|(\nabla'\varphi)(x')|^2}}
\nonumber\\[2pt]
&=\frac{\big((\nabla'\varphi)(x'),-1\big)}
{\sqrt{1+|(\nabla'\varphi)(x')|^2}}\,\text{ for ${\mathscr{H}}^{n-1}$-a.e.~$x'$ near $x_0'$},  
\end{align}
where $(\nabla'\varphi)(x'):=\big((\partial_1\varphi)(x'),\ldots,(\partial_{n-1}\varphi)(x')\big)$ 
exists for ${\mathscr{H}}^{n-1}$-a.e.~$x'\in{\mathbb{R}}^{n-1}$ thanks to the classical Rademacher theorem
(in this vein, see, e.g., \cite{EG92}). 

For a Lipschitz domain $\Omega$ in $\bbR^n$ the surface measure on $\partial\Omega$ is defined via the formula 
\begin{equation}\label{j6444}
\sigma:={\mathscr{H}}^{n-1}\lfloor{\partial\Omega}.
\end{equation}
As a consequence, the outward unit normal $\nu$ to $\Omega$ exists $\sigma$-a.e.~on $\partial\Omega$.
We also note here that locally, near any boundary point $x_0\in\partial\Omega$, identifying 
$\partial\Omega$ with the graph of a Lipschitz function $\varphi:{\mathbb{R}}^{n-1}\to{\mathbb{R}}$
(in a suitable system of coordinates, isometric with the original one) permits us to express 
the surface measure in this new system of coordinates as
\begin{equation}\label{ytrre5}
d^{n-1}\sigma(x)=\sqrt{1+|(\nabla'\varphi)(x')|^2}\,d^{n-1}x'\,\text{ for $x=(x',\varphi(x'))$ near $x_0$}. 
\end{equation}

The theorem below, established in \cite[Theorem~2.10]{ABMMZ11}, formalizes the idea that a connected, 
proper, open subset of ${\mathbb{R}}^n$ whose boundary is a compact Lipschitz surface 
is a Lipschitz domain. Before stating this fact, we note that the connectivity assumption 
is necessary since, otherwise, $\Omega:=\{x\in{\mathbb{R}}^n\,|\,|x|<2\,\text{ and }\,|x|\not=1\}$ 
would serve as a counterexample. 

%%%%%%%
\begin{theorem}\label{J.2BBB}
Let $\Omega$ be a nonempty, connected, proper, open subset of ${\mathbb{R}}^n$, with $\partial\Omega$ 
bounded. In addition, suppose that for each $x_0\in\partial\Omega$ there exist an $(n-1)$-dimensional
plane $H\subseteq{\mathbb{R}}^n$ passing through $x_0$, a choice $N$ of the unit normal to $H$, an open 
cylinder ${\mathcal{C}}_{r,\tau}(x_0,N)=\{x'+tN\,|\,x'\in H,\,|x'-x_0|<r,\,\,|t|<\tau\}$ and a 
Lipschitz function $\varphi:H\to{\mathbb{R}}$ satisfying \eqref{Def-Lip3} such that \eqref{Def-Lip1} holds. 
Then $\Omega$ is a Lipschitz domain.
\end{theorem}
%%%%%%%

The proof of the above result relies on Lemma~\ref{JGbv.2} and the generalization of the Jordan-Brouwer separation
theorem for arbitrary compact topological hypersurfaces in ${\mathbb{R}}^n$ noted in \cite[Theorem~1, p.~284]{Al78}.
To proceed, we make the following definition. 

%%%%%%%
\begin{definition}\label{Def1A}
$(i)$ A nonempty set $\Omega\subseteq{\mathbb{R}}^n$ is called {\tt starlike with respect to} 
$x_0\in\Omega$ if $\mathcal{I}(x,x_0)\subseteq\Omega$ for every $x\in\Omega$, where $\mathcal{I}(x,x_0)$ 
denotes the open line segment in ${\mathbb{R}}^n$ with endpoints $x$ and $x_0$. \\[1mm] 
$(ii)$ A nonempty set $\Omega\subseteq{\mathbb{R}}^n$ is called 
{\tt starlike with respect to a ball} if there exists a ball $B\subseteq\Omega$ with 
the property that $\mathcal{I}(x,y)\subseteq\Omega$ for every $x\in\Omega$ and every $y\in B$ 
{\rm (}that is, $\Omega$ is starlike with respect to any point in $B${\rm )}.
\end{definition}
%%%%%%%

It turns out that local Lipschitzianity may be characterized in terms of local starlikeness 
(with respect to balls), in the precise sense described in the theorem below, proved in 
\cite[Theorem~3.9]{ABMMZ11}.

%%%%%%%
\begin{theorem}\label{Starlike}
Let $\Omega$ be an open, proper, nonempty subset of ${\mathbb{R}}^n$.
Then $\Omega$ is a locally Lipschitz domain if and only if every point
$x_0\in\partial\Omega$ has an open neighborhood ${\mathcal{O}}\subseteq{\mathbb{R}}^n$ 
with the property that $\Omega\cap{\mathcal{O}}$ is starlike with respect to a ball. 

Moreover, any nonempty bounded convex open set is a Lipschitz domain. 
\end{theorem}
%%%%%%%

Next, we discuss various types of cone properties possessed by locally Lipschitz domains. 
By an open, truncated, one-component circular cone in ${\mathbb{R}}^n$ we shall understand a set of the form 
\begin{equation}\label{Cg-77}
{\mathscr{U}}_{\theta,h}(x_0,v):=\big\{x\in{\mathbb{R}}^n\,\big|\,
\cos(\theta/2)\,|x-x_0|<(x-x_0)\cdot v<h\big\},
\end{equation}
where $x_0\in{\mathbb{R}}^n$ is the vertex of the cone, $v\in \bbS^{n-1}$ is the direction of 
the axis, $\theta\in(0,\pi)$ is the (full) aperture of the cone, $h\in(0,\infty)$ is the height 
of the cone, and where ``dot" denotes the standard inner product in ${\mathbb{R}}^n$.

Here is a characterization of local Lipschitzianity in terms of a two-sided cone condition
from \cite[Proposition~3.7]{ABMMZ11}. 

%%%%%%%
\begin{theorem}\label{Hv56.3}
Assume that $\Omega\subseteq{\mathbb{R}}^n$ is a nonempty, proper, open set 
and fix a point $x_0\in\partial\Omega$. Then $\Omega$ is a Lipschitz domain near $x_0$ 
if and only if there exist a height $h\in(0,\infty)$, an angle $\theta\in(0,\pi)$, along with 
a radius $r\in(0,\infty)$ and a function $v:B(x_0,r)\cap\partial\Omega\to \bbS^{n-1}$ which is continuous at $x_0$ 
and with the property that  
\begin{equation}\label{Hv56.2}
{\mathscr{U}}_{\theta,h}(x,v(x))\subseteq\Omega\,\text{ and }\,
{\mathscr{U}}_{\theta,h}(x,-v(x))\subseteq{\mathbb{R}}^n\backslash\Omega,
\quad\forall\,x\in B(x_0,r)\cap\partial\Omega. 
\end{equation}
\end{theorem}
%%%%%%%

The global two-sided cone property for bounded Lipschitz domains recorded below
is a direct consequence of Theorem~\ref{Hv56.3}.

%%%%%%%
\begin{corollary}\label{Hv56.3-CC}
Let $\Omega\subset{\mathbb{R}}^n$ be a bounded Lipschitz domain. 
Then there exist a height $h\in(0,\infty)$, an angle $\theta\in(0,\pi)$, and a continuous 
function $v:\partial\Omega\to \bbS^{n-1}$ such that 
\begin{equation}\label{Hv56.2-CC}
{\mathscr{U}}_{\theta,h}(x,v(x))\subseteq\Omega\,\text{ and }\,
{\mathscr{U}}_{\theta,h}(x,-v(x))\subseteq{\mathbb{R}}^n\backslash\Omega,
\quad\forall\,x\in\partial\Omega. 
\end{equation}
\end{corollary}
%%%%%%%

In fact, it is possible to characterize local Lipschitzianity in terms of one-sided
cone conditions. The case of an exterior cone condition is described in the next theorem, 
proved in \cite[Proposition~3.5]{ABMMZ11}.

%%%%%%%
\begin{theorem}\label{Pr-2Ciii}
Let $\Omega$ be a proper, nonempty open subset of ${\mathbb{R}}^n$ 
and fix $x_0\in\partial\Omega$. Then the set $\Omega$ is a Lipschitz domain 
near $x_{0}$ if and only if there exist two numbers $r,h\in(0,\infty)$, 
an angle $\theta\in(0,\pi)$, along with a function $v:B(x_{0},r)\cap\partial\Omega\to \bbS^{n-1}$ 
which is continuous at $x_{0}$ and such that 
\begin{equation}\label{Co-P1Bx}
{\mathscr{U}}_{\theta,h}(x,v(x))\subseteq{\mathbb{R}}^n\backslash\Omega,
\quad\forall\,x\in B(x_0,r)\cap\partial\Omega. 
\end{equation}
\end{theorem}
%%%%%%%

Finally, a characterization of local Lipschitzianity in terms of an interior cone condition
is contained in the theorem below (taken from \cite[Proposition~3.6]{ABMMZ11}). 

%%%%%%%
\begin{theorem}\label{Pr-2Cc.32}
Assume that $\Omega\subseteq{\mathbb{R}}^n$ is an open set and suppose $x_0\in\partial\Omega$. 
Then $\Omega$ is a Lipschitz domain near $x_{0}$ if and only if there 
exist two numbers $r,h\in(0,\infty)$, an angle $\theta\in(0,\pi)$, and a function
$v:B(x_{0},r)\cap\partial\Omega\to \bbS^{n-1}$ which is continuous at $x_{0}$ 
and such that 
\begin{align}\label{Co-P1Bx.TT} 
\begin{split} 
& B(x_0,r)\cap\partial\Omega=B(x_0,r)\cap\partial(\overline{\Omega})\,\text{ and }
\\[2pt] 
& \quad{\mathscr{U}}_{\theta,h}(x,v(x))\subseteq\Omega,\quad\forall\,x\in B(x_0,r)\cap\partial\Omega. 
\end{split} 
\end{align}
\end{theorem}
%%%%%%%

Next, we recall several basic definitions. Given a bounded Lipschitz domain $\Omega$ in 
${\mathbb{R}}^n$ and some fixed $\kappa\in(0,\infty)$, for each $x\in\partial\Omega$
we first define the nontangential approach region with vertex at $x$ and aperture parameter $\kappa$ by setting 
\begin{equation}\label{eq:MM1}
\Gamma_\kappa(x):=\big\{y\in\Omega\,\big|\,|x-y|<(1+\kappa)\,{\rm dist}\,(y,\partial\Omega)\big\}.
\end{equation}
Results in \cite{MMM22} prove that  
\begin{equation}\label{kj64d5-y543}
x\in\overline{\Gamma_\kappa(x)}\,\text{ for $\sigma$-a.e.~}\,x\in\partial\Omega.
\end{equation}
Second, given an arbitrary $u:\Omega\to{\mathbb{C}}$, we define its nontangential maximal function and
its pointwise nontangential boundary trace at $x\in\partial\Omega$, respectively, as 
\begin{equation}\label{eq:MM2}
\big({\mathcal{N}}_\kappa u\big)(x):=\sup\big\{|u(y)|\,\big|\,y\in\Gamma_\kappa(x)\big\}
\in[0,\infty],
\end{equation}
and
\begin{equation}\label{eq:MM2-BIS}
\Big(u\big|^{\kappa-{\rm n.t.}}_{\partial\Omega}\Big)(x):=\lim_{\Gamma_\kappa(x)\ni y\to x}u(y),
\end{equation}
whenever the above limit exists. In this connection remark that by 
\eqref{kj64d5-y543} the nontangential convergence $\Gamma_\kappa(x)\ni y\to x$ 
in \eqref{eq:MM2-BIS} makes sense for $\sigma$-a.e.~$x\in\partial\Omega$.

These definitions readily adapt to vector-valued functions, in a natural fashion 
(interpreting $|u(y)|$ as norm in \eqref{eq:MM2}, and considering
$u\big|^{\kappa-{\rm n.t.}}_{\partial\Omega}$ componentwise). In the sequel, we shall make 
no notation distinction between the scalar-valued and the vector-valued case.   
Clearly, 
\begin{equation}\label{eq:HByrr}
\big|u\big|^{\kappa-{\rm n.t.}}_{\partial\Omega}\big|\leq{\mathcal{N}}_\kappa u\,\text{ pointwise on }\,\partial\Omega. 
\end{equation}
It turns out that ${\mathcal{N}}_\kappa u$ is a lower semi-continuous 
function on $\partial\Omega$, hence Lebesgue measurable. In addition, the parameter $\kappa$ plays 
a somewhat secondary role, since for any $\kappa_1,\kappa_2\in(0,\infty)$ and $p\in(0,\infty)$ 
there exists $C=C(\kappa_1,\kappa_2,p)\in(1,\infty)$ with the property that, for each 
$u:\Omega\to{\mathbb{C}}$,  
\begin{equation}\label{eq:HBab}
C^{-1}\big\|{\mathcal{N}}_{\kappa_1}u\big\|_{L^p(\partial\Omega)}
\leq\big\|{\mathcal{N}}_{\kappa_2}u\big\|_{L^p(\partial\Omega)}
\leq C\big\|{\mathcal{N}}_{\kappa_1}u\big\|_{L^p(\partial\Omega)}.
\end{equation}
Also, whenever $u:\Omega\to{\mathbb{C}}$ is such that ${\mathcal{N}}_{\kappa}u\in L^p(\partial\Omega)$ for some 
$\kappa>0$ and $p\in(0,\infty)$ for any aperture parameters $\kappa_1,\kappa_2\in(0,\infty)$ it follows that
\begin{equation}\label{eq:FXGF.NEW}
\begin{array}{c}
u\big|^{\kappa_1-{\rm n.t.}}_{\partial\Omega}\,\,\text{ exists $\sigma$-a.e. on }\,\,\partial\Omega 
\,\,\text{ if and only if}
\\[2pt]
u\big|^{\kappa_2-{\rm n.t.}}_{\partial\Omega}\,\,\text{ exists $\sigma$-a.e. on }\,\,\partial\Omega.
\end{array}
\end{equation}

We shall need two additional properties of the nontangential maximal operator
{\rm (}i.e., the mapping $u\mapsto{\mathcal{N}}_{\kappa}u${\rm )}.
First, as proved in \cite[Proposition~2.3]{MM13}, for any $p\in(0,\infty)$ there exists $C_p\in(0,\infty)$ 
with the property that for every measurable function $u:\Omega\to{\mathbb{C}}$ one has
\begin{align}\label{eq:HBab.2}
\begin{split}
& {\mathcal{N}}_\kappa u\in L^p(\partial\Omega)\,\text{ implies }\,u\in L^{np/(n-1)}(\Omega) 
\\[2pt]
& \quad\text{and }\,\|u\|_{L^{np/(n-1)}(\Omega)}\leq C_p\|{\mathcal{N}}_\kappa u\|_{L^p(\partial\Omega)}.
\end{split}
\end{align}
The second property alluded to above is contained in the lemma below. 

%%%%%%%
\begin{lemma}\label{Hv56.3-LL}
For any bounded Lipschitz domain $\Omega\subset{\mathbb{R}}^n$ there exists a 
compact set $K\subset\Omega$ with the property that for each $p\in(0,\infty)$ and $\kappa>0$
one can find a constant $C\in(0,\infty)$ such that  
\begin{equation}\label{eqn.hAAbvvb-1}
\big\|{\mathcal{N}}_\kappa u\big\|_{L^p(\partial\Omega)}
\leq C\big(\big\|{\mathcal{N}}_\kappa(\nabla u)\big\|_{L^p(\partial\Omega)}+\sup_{x\in K}|u(x)|\big),
\end{equation}
for every function $u\in C^1(\Omega)$. 
\end{lemma}
%%%%%%%
\begin{proof}
Given a bounded Lipschitz domain $\Omega\subset{\mathbb{R}}^n$, 
pick the parameters $h\in(0,\infty)$, $\theta\in(0,\pi)$, and the continuous function 
$v:\partial\Omega\to \bbS^{n-1}$ as in Corollary~\ref{Hv56.3-CC}. Then, for a suitably small $r>0$,
define $K:=\{x\in\Omega\,|\,{\rm dist}\,(x,\partial\Omega)\geq r\}$. Specifically, we select 
$r>0$ such that for every $x\in\partial\Omega$ the entire flat portion of the boundary of the 
truncated circular cone ${\mathscr{U}}_{\theta,h}(x,v(x))$ is contained in $K$. 

Next, we pick an arbitrary point $x\in\partial\Omega$ along with some $y\in{\mathscr{U}}_{\theta,h}(x,v(x))$,
and consider 
\begin{equation}\label{k754rFF}
\text{$z:=x+t(y-x)$, where $t:=\frac{h}{(y-x)\cdot v(x)}$}.
\end{equation}
Then the fact that $(z-x)\cdot v(x)=h$ places the point $z$ on the flat portion 
of the boundary of ${\mathscr{U}}_{\theta,h}(x,v(x))$. In particular, $z\in K$. 
Keeping this in mind, it follows that for every function $u\in C^1(\Omega)$ we may estimate 
(using the Mean-Value Theorem and the fact that ${\mathscr{U}}_{\theta,h}(x,v(x))$ is a convex subset of $\Omega$)
\begin{align}\label{k754rFF-2}
& |u(y)|\leq|u(y)-u(z)|+|u(z)|\leq|y-z|\sup_{\xi\in[y,z]}|(\nabla u)(\xi)|+\sup_{\zeta\in K}|u(\zeta)|
\nonumber\\[2pt]
& \quad\leq C_{\theta,h}\sup\{|(\nabla u)(\xi)|\,|\,\xi\in{\mathscr{U}}_{\theta,h}(x,v(x))\}
+\sup_{\zeta\in K}|u(\zeta)|,
\end{align}
for some constant $C_{\theta,h}\in(0,\infty)$. These considerations suggest introducing 
the following version of the nontangential maximal operator 
\begin{equation}\label{eq:MM2-NEW}
\big(\widetilde{\mathcal{N}}_{\theta,h}w\big)(x):=\sup\{|w(y)|\,|\,y\in{\mathscr{U}}_{\theta,h}(x,v(x))\},
\quad\forall\,x\in\partial\Omega,
\end{equation}
where $w$ is an arbitrary (possibly vector-valued) continuous function defined in $\Omega$. 
In this notation, \eqref{k754rFF-2} yields 
\begin{equation}\label{eqn.hAAb-FFCV}
\big(\widetilde{\mathcal{N}}_{\theta,h}u\big)(x)\leq
C_{\theta,h}\big(\widetilde{\mathcal{N}}_{\theta,h}(\nabla u)\big)(x)+\sup_{\zeta\in K}|u(\zeta)|,
\quad\forall\,x\in\partial\Omega,
\end{equation}
hence, further, 
\begin{equation}\label{eqn.hAAb-FFCV-bbb}
\big\|\widetilde{\mathcal{N}}_{\theta,h}u\big\|_{L^p(\partial\Omega)}\leq
C\big(\big\|\widetilde{\mathcal{N}}_{\theta,h}(\nabla u)\big\|_{L^p(\partial\Omega)}
+\sup_{\zeta\in K}|u(\zeta)|\big),
\end{equation}
for every function $u\in C^1(\Omega)$. Having established \eqref{eqn.hAAb-FFCV-bbb}, 
Proposition~2.2 in \cite{MM13} and the remark following its proof (where the two brands of nontangential 
maximal operators, $\widetilde{\mathcal{N}}_{\theta,h}$ and ${\mathcal{N}}_\kappa$, are compared) 
then allow us to conclude that \eqref{eqn.hAAbvvb-1} holds for every function $u\in C^1(\Omega)$.
\end{proof}
%%%%%%%

Both the notion of nontangential maximal function and the notion of nontangential boundary trace 
are pivotal in the formulation of the following version of the divergence theorem recorded below.   
This is a particular case of a result established in \cite{MMM22} (see also \cite{MMM20}) for a more 
general category of sets than the class of Lipschitz domains.

%%%%%%%%%%%
\begin{theorem}\label{banff-3}
Let $\Omega\subset{\mathbb{R}}^n$ be a bounded Lipschitz domain and denote by $\nu$ the outward 
unit normal to $\Omega$, which is well defined $\sigma$-a.e.~on $\partial\Omega$, where 
$\sigma$ is the canonical surface measure defined as in \eqref{j6444}. Also, fix some 
aperture parameter $\kappa>0$. Then for every vector field satisfying
\begin{align}\label{banff-4.LL}
\begin{split}
& \text{$\vec{F}\in\big[L^1_{\rm loc}(\Omega)\big]^n$, the nontangential trace 
$\vec{F}\big|^{\kappa-{\rm n.t.}}_{\partial\Omega}$ exists $\sigma$-a.e.~on $\partial\Omega$},    
\\[2pt]
& \quad{\mathcal{N}}_\kappa(\vec{F}\,)\,\,\text{ belongs to }\,\,L^1(\partial\Omega),\,\,\text{ and }\,\,
{\rm div}\vec{F}\,\,\text{ belongs to }\,\,L^1(\Omega)
\end{split}
\end{align}
{\rm (}with the divergence taken in the sense of distributions in $\Omega${\rm )}, one has  
\begin{equation}\label{banff-5}
\int_{\Omega}{\rm div}\vec{F}\,d^n x
=\int_{\partial\Omega}\nu\cdot\Big(\vec{F}\,\big|^{\kappa-{\rm n.t.}}_{\partial\Omega}\Big)\,d^{n-1}\sigma,
\end{equation}
where, as before, ``dot" denotes the standard inner product in ${\mathbb{R}}^n$. As a corollary of this and \eqref{kj64d5-y543},
\begin{align}\label{banff-5bFG}
& \displaystyle\int_{\Omega}{\rm div}\vec{F}\,d^n x
=\int_{\partial\Omega}\nu\cdot\big(\vec{F}\big|_{\partial\Omega}\big)\,d^{n-1}\sigma\,\text{ for every}
\nonumber\\[2pt]
& \quad\text{vector field $\vec{F}\in\big[C^0(\overline{\Omega})\big]^n$ with ${\rm div}\vec{F}\in L^1(\Omega)$}
\\[2pt]
& \quad\text{{\rm\big(}hence, in particular, for every $\vec{F}\in\big[C^1(\overline{\Omega})\big]^n${\rm\big)}}.
\nonumber
\end{align}
\end{theorem}
%%%%%%%%%%%

In the next lemma we record an approximation procedure developed in \cite{Ca85}, \cite{MMM13}, \cite{MMY10}, \cite{Ve84}. 

%%%%%%%
\begin{lemma}\label{OM-OM}
Given a bounded Lipschitz domain $\Omega\subset{\mathbb{R}}^n$, there exists a family 
$\{\Omega_\ell\}_{\ell\in{\mathbb{N}}}$ of domains in ${\mathbb{R}}^n$ satisfying 
the following properties: \\[1mm] 
$(i)$ Each $\Omega_\ell$ is a bounded Lipschitz domain, with Lipschitz character 
bounded uniformly in $\ell\in{\mathbb{N}}$. \\[1mm] 
$(ii)$ For every $\ell\in{\mathbb{N}}$ one has 
$\overline{\Omega_\ell}\subset\Omega_{\ell+1}\subset\Omega$, and 
$\Omega=\bigcup_{\ell\in{\mathbb{N}}}\Omega_\ell$. \\[1mm] 
$(iii)$ There exist $\kappa\in(0,\infty)$ and bi-Lipschitz homeomorphisms 
$\Lambda_\ell:\partial\Omega\to\partial\Omega_\ell$, $\ell\in{\mathbb{N}}$, 
such that for every $x\in\partial\Omega$ one has $\Lambda_\ell(x)\to x$ as $\ell\to\infty$, 
and $\Lambda_\ell(x)\in\Gamma_\kappa(x)$ for each $\ell\in{\mathbb{N}}$. \\[1mm] 
$(iv)$ If for each $\ell\in{\mathbb{N}}$ we let $\nu^\ell$ be the outward unit 
normal to $\Omega_\ell$ and if $\nu$ denotes the outward unit normal to $\Omega$, then 
$\nu^\ell\circ\Lambda_\ell\to\nu$ as $\ell\to\infty$ both pointwise $\sigma$-a.e.~and in 
$\big[L^2(\partial\Omega)\big]^n$. \\[1mm] 
$(v)$ There exist non-negative, measurable functions $\omega_\ell$ on $\partial\Omega$ which are bounded
away from zero and infinity uniformly in $\ell\in{\mathbb{N}}$, converge pointwise $\sigma$-a.e.~to $1$ 
as $\ell\to\infty$, and which have the property that for each integrable function 
$g:\partial\Omega_\ell\to{\mathbb{R}}$ the following change of variable formula holds 
\begin{equation}\label{eQQ-14}
\int_{\partial\Omega_\ell}g\,d^{n-1}\sigma_\ell
=\int_{\partial\Omega}g\circ\Lambda_\ell\,\omega_\ell\,d^{n-1}\sigma,
\end{equation}
where $\sigma_\ell$ is the canonical surface measure on $\partial\Omega_\ell$. 
\end{lemma}
%%%%%%%

\noindent We shall use the notation $\Omega_\ell\nearrow\Omega$ as $\ell\to\infty$ to indicate that 
the family $\{\Omega_\ell\}_{\ell\in{\mathbb{N}}}$ approximates $\Omega$ 
in the manner described in Lemma~\ref{OM-OM} above.

%%%%%%%%%%%%%%%%%%%%%%%%%%%%%%%%%%%
\subsection{Fractional Sobolev, Besov, and Triebel--Lizorkin spaces in arbitrary open sets}\label{ss2.2}
%%%%%%%%%%%%%%%%%%%%%%%%%%%%%%%%%%%

Given a nonempty open set $\Omega\subseteq\bbR^n$, we denote by $H^s(\Omega)$ the scale of 
$L^2$-based Sobolev spaces of (fractional) order $s\in\mathbb{R}$ in $\Omega$. 
More specifically, with $\cS'(\bbR^n)$ and $\mathscr{F}$ denoting, respectively, the 
space of tempered distributions and the Fourier transform in $\bbR^n$, for each $s\in\bbR$ set 
\begin{equation}\label{u54ff}
H^s(\bbR^n):=\big\{f\in\cS'(\bbR^n)\,\big|\,(1+|\xi|^2)^{s/2}\mathscr{F}f\in L^2(\bbR^n)\big\},
\end{equation}
equipped with the natural norm 
\begin{align}\label{uyrtgoii}
\|f\|_{H^s(\bbR^n)} &:=\big\|(1+|\cdot|^2)^{s/2}(\mathscr{F}f)(\cdot)\big\|_{L^2(\bbR^n)}
\nonumber\\[2pt]
&=\Big(\int_{\bbR^n}(1+|\xi|^2)^{s}|(\mathscr{F}f)(\xi)|^2\,d^n\xi\Big)^{1/2}.
\end{align} 
Then define 
\begin{equation}\label{tgBBn}
H^s(\Omega):=\big\{f\in\cD'(\Omega)\,\big|\,\text{there exists $g\in H^s(\bbR^n)$ 
such that $f=g|_{\Omega}$}\big\},
\end{equation}
where $g|_{\Omega}\in\cD'(\Omega)$ stands for the restriction of the distribution 
$g\in\cD'({\mathbb{R}}^n)$ to the open set $\Omega$, and endow the space \eqref{tgBBn} with the norm 
\begin{equation}\label{HGaYga.3}
\|f\|_{H^s(\Omega)}:=\inf_{\substack{g\in H^s(\bbR^n)\\ f=g|_\Omega}}
\|g\|_{H^s(\bbR^n)},\quad\forall\,f\in H^s(\Omega).
\end{equation}
The above definition allows for more or less directly transferring a number of properties 
of the scale of fractional Sobolev spaces in ${\mathbb{R}}^n$ to the corresponding version 
of that scale considered in an arbitrary open subset $\Omega$ of the Euclidean space. 
For example, we have
\begin{align}\label{eq:DDEEn.2}
H^{s_1}(\Omega)\hookrightarrow H^{s_2}(\Omega)
\,\text{ continuously, if }\,s_1,s_2\in{\mathbb{R}},\,\,s_1\geq s_2,
\end{align}
and 
\begin{equation}\label{eq:DDEEn.3}
\partial^\alpha:H^s(\Omega)\rightarrow H^{s-|\alpha|}(\Omega)
\,\text{ continuously, for each $\alpha\in{\mathbb{N}}_0^n$ and $s\in{\mathbb{R}}$}.
\end{equation}
Furthermore, if we set 
\begin{equation}\label{uam-mi77}
C^\infty(\overline{\Omega}):=\big\{\psi|_{\Omega}\,\big|\,\psi\in C^\infty_0(\bbR^n)\big\}
\end{equation}
then 
\begin{equation}\label{uam-mi78}
C^\infty(\overline{\Omega})\hookrightarrow H^s(\Omega)\,\text{ densely, for every }\,s\in\bbR,
\end{equation}
and 
\begin{align}\label{eq:DDEj6g5}
\begin{split} 
& \text{for every $\psi\in C^\infty(\overline{\Omega})$ and every $s\in{\mathbb{R}}$, the assignment} 
\\[2pt]  
& \quad\text{$H^s(\Omega)\ni u\mapsto\psi u\in H^s(\Omega)$ is well defined, linear, and bounded.}
\end{split} 
\end{align}

Given an open set $\Omega\subseteq{\mathbb{R}}^n$ and some $p\in(0,\infty)$, 
we use $L^p_{\rm loc}(\Omega)$ to denote the space of functions which are locally 
$p$-th power integrable in $\Omega$. We shall also occasionally work with the 
local version of the scale \eqref{tgBBn}, defined for $s\in{\mathbb{R}}$ as
\begin{equation}\label{tgBBn.344f}
H^s_{\rm loc}(\Omega):=\big\{f\in\cD'(\Omega)\,\big|\,\zeta f\in H^s(\Omega)
\text{ for every }\zeta\in C^\infty_0(\Omega)\big\}.
\end{equation}
In addition, for each $s\in\bbR$, by $\accentset{\circ}{H}^s(\Omega)$ we shall denote the closure of 
$ C^\infty_0(\Omega)$ in $H^s(\Omega)$, that is, 
\begin{equation}\label{Rdac}
\accentset{\circ}{H}^s(\Omega):=\overline{C^\infty_0(\Omega)}^{H^s(\Omega)},\quad\forall\,s\in\bbR. 
\end{equation} 

Finally, we consider $L^2$-based Sobolev spaces of integer order, that is, $W^{k}(\Omega)$ with 
$k\in{\mathbb{N}}_0$, intrinsically defined in $\Omega$ as
\begin{equation}\label{eq:WWW-sp}
W^{k}(\Omega):=\big\{u\in L^1_{\rm loc}(\Omega)\,\big|\,\partial^\alpha u\in L^2(\Omega)
\text{ for each $\alpha\in{\mathbb{N}}_0^n$ with $|\alpha|\leq k$}\big\},
\end{equation}
and equipped with the natural norm 
\begin{equation}\label{eq:WWW-nrm}
\|u\|_{W^{k}(\Omega)}:=\sum_{|\alpha|\leq k}\|\partial^\alpha u\|_{L^2(\Omega)},
\quad\forall\,u\in W^{k}(\Omega).
\end{equation}
Furthermore, given $k\in{\mathbb{N}}_0$ set 
\begin{equation}\label{Rdac.WWW}
\accentset{\circ}{W}^k(\Omega):=\ol{ C^\infty_0(\Omega)}^{W^k(\Omega)}. 
\end{equation}

While for arbitrary open sets $\Omega\subset{\mathbb{R}}^n$ one only has
$H^k(\Omega)\subset W^k(\Omega)$ for each $k\in{\mathbb{N}}_0$, equality actually holds
in the class of bounded Lipschitz domains (to be discussed later; cf. \eqref{eq:WH}). 

\medskip 

Fix a family of Schwartz functions $\{\zeta_j\}_{j=0}^\infty\subset\cS(\bbR^n)$ 
possessing the following properties: \\[1mm] 
(a) there exist constants $a,b,c\in(0,\infty)$ such that
\begin{equation}\label{besov-jussi3}
\begin{cases}
\supp\,(\zeta_0)\subset\{x\in\mathbb{R}^n\,|\,\vert x\vert\leq a\},
\\[2pt]
\supp\,(\zeta_j)\subset\{x\in\mathbb{R}^n\,|\,b\,2^{j-1}\leq\vert x\vert\leq c\,2^{j+1}\}\,\text{ for each }\,j\in{\mathbb{N}};
\end{cases}
\end{equation}
(b) for every multi-index $\alpha\in{\mathbb{N}}_0^n$ there exists a number $C_\alpha\in(0,\infty)$ such that
\begin{equation}\label{besov-jussi4}
\sup_{x\in\mathbb{R}^n}\sup_{j\in\mathbb{N}}\,2^{j\vert\alpha\vert}\vert\partial^\alpha\zeta_j(x)\vert\leq C_\alpha;
\end{equation}
(c) for every $x\in\mathbb{R}^n$ one has 
\begin{equation}\label{besov-jussi5}
\sum_{j=0}^\infty\zeta_j(x)=1.
\end{equation}
Then the standard Besov scale in ${\mathbb{R}}^n$ consists of spaces $B^{p,q}_{s}({\mathbb{R}}^n)$ 
defined for each $p,q\in(0,\infty]$ and $s\in{\mathbb{R}}$ as
\begin{equation}\label{besov-jussi}
B^{p,q}_{s}({\mathbb{R}}^n):=\bigg\{f\in\cS'(\bbR^n)\,\bigg|\,
\sum_{j=0}^\infty\Vert 2^{sj}\mathscr{F}^{-1}(\zeta_j\mathscr{F}f)\Vert^q_{L^p(\mathbb{R}^n)}<\infty\bigg\}.
\end{equation}
Each such space is equipped with the natural quasi-norm
\begin{equation}\label{besov-jussi2}
B^{p,q}_{s}({\mathbb{R}}^n)\ni f\mapsto\Vert f\Vert_{B^{p,q}_{s}({\mathbb{R}}^n)}:=
\bigg(\sum_{j=0}^\infty\Vert 2^{sj}\mathscr{F}^{-1}(\zeta_j\mathscr{F}f)\Vert^q_{L^p(\mathbb{R}^n)}\bigg)^\frac{1}{q},
\end{equation}
rendering $B^{p,q}_{s}({\mathbb{R}}^n)$ a quasi-Banach space (which is actually a genuine Banach 
space in the range $1\leq p,q\leq\infty$). We mention that a different choice of a family of functions
$\{\zeta_j\}_{j=0}^\infty\subset\cS(\bbR^n)$ satisfying (a)--(c) in \eqref{besov-jussi}--\eqref{besov-jussi2} 
yields the same vector space, which is now equipped with an equivalent quasi-norm. We note also that for 
$0<p,q<\infty$ and $s\in{\mathbb{R}}$ the class of Schwartz functions in $\mathbb{R}^n$ is dense in 
$B^{p,q}_{s}({\mathbb{R}}^n)$. There is a wealth of material pertaining to Besov spaces in the 
Euclidean setting and the interested reader is referred to the monographs \cite{BL76} by J.\ Bergh 
and J.\ L\"{o}fstr\"{o}m, \cite{RS96} by T.\ Runst and W.\ Sickel, and \cite{Tr83} by H.\ Triebel. 

Moving on, having fixed an arbitrary open set $\Omega\subseteq{\mathbb{R}}^n$, 
whenever $0<p,q\leq\infty$ and $s\in{\mathbb{R}}$ it is meaningful to define
\begin{align}\label{restr}
\begin{split}
& B^{p,q}_s(\Omega):=\big\{f\in{\mathcal{D}}'(\Omega)\,\big|\,\text{there exists }\,
g\in B^{p,q}_s({\mathbb{R}}^n)\,\text{ such that }\,g|_{\Omega}=f\big\},
\\[2pt]
& \quad\|f\|_{B^{p,q}_s(\Omega)}:=\inf\big\{\|g\|_{B^{p,q}_s({\mathbb{R}}^n)}\,\big|\,
g\in B^{p,q}_s({\mathbb{R}}^n),\,\,g|_{\Omega}=f\big\},\quad\forall f\in B^{p,q}_s(\Omega).
\end{split}
\end{align}
This definition permits transferring with ease a number of properties shared by Besov spaces 
in the Euclidean setting (cf., e.g., the discussion in \cite[Section~2.2]{RS96}) to arbitrary 
open subsets of ${\mathbb{R}}^n$, such as
\begin{equation}\label{Inc-P3}
B^{2,2}_s(\Omega)=H^s(\Omega)\,\text{ for each }\,s\in{\mathbb{R}},
\end{equation}
(identical vector spaces with equivalent norms) and, with continuous inclusions, 
\begin{align}\label{q-TWU-iii}
& B^{p,\infty}_{s_0}(\Omega)\hookrightarrow B^{p,q}_{s_1}(\Omega)\,\text{ if }\,s_0>s_1,\quad 0<p,q\leq\infty, 
\\[2pt]
& B^{p,q_0}_{s}(\Omega)\hookrightarrow B^{p,q_1}_{s}(\Omega)
\,\text{ if }\,0<q_0\leq q_1\leq\infty,\,\,\,0<p\leq\infty,\,\,s\in{\mathbb{R}}.
\label{q-TWC}
\end{align}
Moreover, we note that \eqref{q-TWC} (used with $q_1:=\infty$ and $s:=s_0$) together with 
\eqref{q-TWU-iii} (used with $q:=q_1$) imply 
\begin{equation}\label{q-TWU}
B^{p,q_0}_{s_0}(\Omega)\hookrightarrow B^{p,q_1}_{s_1}(\Omega)
\,\text{ if }\,s_0>s_1\,\text{ and }\,0<p,q_0,q_1\leq\infty.
\end{equation}

In addition, for each multi-index $\alpha\in{\mathbb{N}}_0^n$, the partial derivative operator 
\begin{align}\label{q-TWU.ted}
\begin{split}
& \partial^\alpha:B^{p,q}_{s}(\Omega)\rightarrow B^{p,q}_{s-|\alpha|}(\Omega)
\,\text{ is well defined and bounded}
\\[2pt]
& \quad\text{whenever }\,0<p,q\leq\infty\,\text{ and }\,s\in{\mathbb{R}}.
\end{split}
\end{align}
In particular, from \eqref{Inc-P3} and \eqref{q-TWU} one concludes that for 
any open set $\Omega\subseteq{\mathbb{R}}^n$ one has the continuous inclusion (to be relevant shortly) 
\begin{equation}\label{q-TWU-yt}
H^{s_0}(\Omega)=B^{2,2}_{s_0}(\Omega)\hookrightarrow B^{2,1}_{s_1}(\Omega)\,\text{ whenever }\,s_0>s_1. 
\end{equation}
Finally, for each $\varphi\in C^\infty_0({\mathbb{R}}^n)$, the operator of multiplication by $\varphi$ 
(in the sense of distributions) 
\begin{align}\label{q-TWU.ted-DD}
\begin{split}
& B^{p,q}_{s}(\Omega)\ni u\longmapsto\varphi u\in B^{p,q}_{s}(\Omega)\,\text{ is well defined and bounded}
\\[2pt]
& \quad\text{whenever }\,0<p,q<\infty\,\text{ and }\,s\in{\mathbb{R}}.
\end{split}
\end{align}

The scale of Triebel--Lizorkin spaces in ${\mathbb{R}}^n$ may be introduced in a similar fashion
(using the same approach based on Littlewood--Paley theory). Specifically, having fixed a family 
$\{\zeta_j\}_{j=0}^{\infty}$ satisfying properties (a)--(c) listed in \eqref{besov-jussi3}--\eqref{besov-jussi5}, 
for each $s\in {\mathbb{R}}$ and $0<q<\infty$, define the Triebel--Lizorkin space $F_{s}^{p,q}({\mathbb{R}}^n)$ as
\begin{equation}\label{eqd1.8}
F_{s}^{p,q}({\mathbb{R}}^n):=\Bigg\{f\in{\mathcal{S}}'({\mathbb{R}}^n)\,\Bigg|\,
\bigg(\sum_{j=0}^{\infty}|2^{sj}{\mathcal{F}}^{-1}(\zeta_j{\mathscr{F}}f)|^q\bigg)^{1/q}\in L^p({\mathbb{R}}^n)\Bigg\}
\end{equation}
and equip it with the semi-norm 
\begin{equation}\label{eqd1.8.SN}
F_{s}^{p,q}({\mathbb{R}}^n)\ni f\mapsto\|f\|_{F^{p,q}_s({\mathbb{R}}^n)}
:=\Bigg\|\bigg(\sum_{j=0}^{\infty}|2^{sj}{\mathcal{F}}^{-1}(\zeta_j{\mathscr{F}}f)|^q\bigg)^{1/q}\Bigg\|_{L^p({\mathbb{R}}^n)}.
\end{equation}
See \cite{FJ3} for a precise definition of $F^{\infty,q}_s({\mathbb{R}}^n)$ (cf. also \cite{RS96}).
Then, as is well-known, $F_s^{p,q}({\mathbb{R}}^n)$ is a quasi-Banach space whenever $s\in{\mathbb{R}}$, $0<p<\infty$, 
and $0<q\leq\infty$, which is actually a Banach space if $1\leq p<\infty$ and $1\leq q\leq\infty$. In all cases, 
\begin{equation}\label{SD2}
{\mathcal{S}}({\mathbb{R}}^n)\hookrightarrow F_s^{p,q}({\mathbb{R}}^n)\hookrightarrow{\mathcal{S}}'({\mathbb{R}}^n).
\end{equation}
Also, given $s\in{\mathbb{R}}$ along with $0<p<\infty$, 
\begin{equation}\label{SD2-XV}
{\mathcal{S}}({\mathbb{R}}^n)\hookrightarrow F_s^{p,q}({\mathbb{R}}^n)\,\text{ densely, if and only if }\,q<\infty.
\end{equation}
For further reference we also point out that, for each $0<p\leq\infty$ and $s\in{\mathbb{R}}$, 
one has (cf., e.g., \cite{RS96}):
\begin{equation}\label{F-qqq}
F^{p,q_0}_s({\mathbb{R}}^n)\hookrightarrow F^{p,q_1}_s({\mathbb{R}}^n)\,\text{ whenever }\,0<q_0\leq q_1\leq\infty.
\end{equation}
Also, for each $0<p,q\leq\infty$, $s\in{\mathbb{R}}$, and $m\in{\mathbb{N}}$, 
\begin{align}\label{lift}
F^{p,q}_s({\mathbb{R}}^n) &=\big\{f\in{\mathcal{S}}'({\mathbb{R}}^n)\,|\,\partial^\alpha f\in F^{p,q}_{s-m}({\mathbb{R}}^n)
\,\text{ for all $\alpha\in{\mathbb{N}}_0^n$ with }\,|\alpha|\leq m\big\}
\nonumber\\[2pt]
&=\big\{f\in F^{p,q}_{s-m}({\mathbb{R}}^n)\,|\,\partial^\alpha f\in F^{p,q}_{s-m}({\mathbb{R}}^n)
\,\text{ for all $\alpha\in{\mathbb{N}}_0^n$ with }\,|\alpha|=m\big\},
\end{align}
and 
\begin{align}\label{liftF2}
\|f\|_{F^{p,q}_s({\mathbb{R}}^n)} &\approx\sum_{|\alpha|\leq m}\|\partial^\alpha f\|_{F^{p,q}_{s-m}({\mathbb{R}}^n)}
\nonumber\\[2pt]
&\approx\|f\|_{F^{p,q}_{s-m}({\mathbb{R}}^n)}+\sum_{|\alpha|=m}\|\partial^\alpha f\|_{F^{p,q}_{s-m}({\mathbb{R}}^n)},
\end{align}
uniformly in $f\in F^{p,q}_s({\mathbb{R}}^n)$. In particular, for each multi-index $\alpha\in{\mathbb{N}}_0^n$, 
one has the well defined, linear, and bounded operator
\begin{equation}\label{deriv}
\partial^\alpha:F^{p,q}_s({\mathbb{R}}^n)\longrightarrow F^{p,q}_{s-|\alpha|}({\mathbb{R}}^n).
\end{equation} 
Furthermore, one has continuous embeddings (cf., e.g., \cite[p.~30]{RS96})
\begin{equation}\label{LLa-45}
B^{p,\min\{p,q\}}_{s}({\mathbb{R}}^n)\hookrightarrow F^{p,q}_{s}({\mathbb{R}}^n) 
\hookrightarrow B^{p,\max\{p,q\}}_{s}({\mathbb{R}}^n)\,\text{ for }\,0<p,q\leq\infty,\,\,s\in{\mathbb{R}}. 
\end{equation}
In particular, 
\begin{equation}\label{LLa-45.222}
F^{p,p}_{s}({\mathbb{R}}^n)=B^{p,p}_{s}({\mathbb{R}}^n)\,\text{ for }\,0<p\leq\infty,\,\,s\in{\mathbb{R}}
\end{equation}
(identical vector spaces with equivalent quasi-norms).

As in the case of Besov spaces, given an arbitrary open set $\Omega\subseteq{\mathbb{R}}^n$, 
whenever $0<p,q\leq\infty$ and $s\in{\mathbb{R}}$, we define
\begin{align}\label{restr.WACO}
\begin{split}
& F^{p,q}_s(\Omega):=\big\{f\in{\mathcal{D}}'(\Omega)\,\big|\,\text{there exists }\,
g\in F^{p,q}_s({\mathbb{R}}^n)\,\text{ such that }\,g|_{\Omega}=f\big\},
\\[2pt]
& \|f\|_{F^{p,q}_s(\Omega)}:=\inf\big\{\|g\|_{F^{p,q}_s({\mathbb{R}}^n)}\,\big|\,
g\in F^{p,q}_s({\mathbb{R}}^n),\,\,g|_{\Omega}=f\big\},\quad\forall f\in F^{p,q}_s(\Omega).
\end{split}
\end{align}
As in the past, this allows us to readily transfer various properties enjoyed by 
Triebel--Lizorkin spaces in the Euclidean setting (cf., e.g., \cite[Section~2.2]{RS96}) to arbitrary 
open subsets of ${\mathbb{R}}^n$. For instance, for each $\varphi\in C^\infty_0({\mathbb{R}}^n)$, 
the operator of multiplication by $\varphi$ (in the sense of distributions) 
\begin{align}\label{q-TWU.ted-DD.WACO}
\begin{split}
& F^{p,q}_{s}(\Omega)\ni u\longmapsto\varphi u\in F^{p,q}_{s}(\Omega)\,\text{ is well defined and bounded}
\\[2pt]
& \quad\text{whenever }\,0<p,q<\infty\,\text{ and }\,s\in{\mathbb{R}},
\end{split}
\end{align}
and one has the continuous inclusions
\begin{align}\label{q-TWU-iii.WACO}
& F^{p,\infty}_{s_0}(\Omega)\hookrightarrow F^{p,q}_{s_1}(\Omega)\,\text{ if }\,s_0>s_1,\quad 0<p,q\leq\infty,
\\[2pt]
& F^{p,q_0}_{s}(\Omega)\hookrightarrow F^{p,q_1}_{s}(\Omega)
\,\text{ if }\,0<q_0\leq q_1\leq\infty,\,\,\,0<p\leq\infty,\,\,s\in{\mathbb{R}}.
\label{q-TWC.WACO}
\end{align}
Moreover, 
\begin{equation}\label{eq:414151.WACO}
\parbox{5.80cm}{the inclusion in \eqref{q-TWU-iii.WACO} is strict, and so is the inclusion in \eqref{q-TWC.WACO} if $q_0<q_1$.}
\end{equation}
In particular, \eqref{q-TWC.WACO} (used with $q_1:=\infty$ and $s:=s_0$) 
together with \eqref{q-TWU-iii.WACO} (used with $q:=q_1$) imply 
\begin{equation}\label{q-TWU.WACO}
F^{p,q_0}_{s_0}(\Omega)\hookrightarrow F^{p,q_1}_{s_1}(\Omega)\,\text{ if }\,s_0>s_1\,\text{ and }\,0<p,q_0,q_1\leq\infty.
\end{equation}
Finally, \eqref{LLa-45} implies
\begin{equation}\label{LLa-45.WACO}
B^{p,\min\{p,q\}}_{s}(\Omega)\hookrightarrow F^{p,q}_{s}(\Omega)\hookrightarrow B^{p,\max\{p,q\}}_{s}(\Omega)
\,\text{ for }\,0<p,q\leq\infty,\,\,s\in{\mathbb{R}}.
\end{equation}
As a consequence, 
\begin{equation}\label{LLa-45.222.WACO}
F^{p,p}_{s}(\Omega)=B^{p,p}_{s}(\Omega)\,\text{ for }\,0<p\leq\infty,\,\,s\in{\mathbb{R}}
\end{equation}
(identical vector spaces with equivalent quasi-norms).

%%%%%%%%%%%%%%%%%%%%%%%%%%%%
\subsection{Fractional Sobolev and Besov spaces in Lipschitz domains}\label{ss2.3}
%%%%%%%%%%%%%%%%%%%%%%%%%%%%

Hence forth, unless otherwise mentioned, $\Omega\subset{\mathbb{R}}^n$ is a bounded Lipschitz domain. 
In such a setting, one has  
\begin{equation}\label{eq:WH}
H^s(\Omega)=W^s(\Omega)\,\,(\text{hence also }\,\accentset{\circ}{H}^s(\Omega)
=\accentset{\circ}{W}^s(\Omega)),\,\text{ for each }\,s\in\bbN_0,
\end{equation}
in the sense that $H^s(\Omega)$ and $W^s(\Omega)$ coincide as vector spaces, 
and the norm on $H^s(\Omega)$ (from \eqref{HGaYga.3}) is equivalent with 
\begin{equation}\label{uagaIIuY}
f\mapsto\sum_{|\alpha|\leq s}\|\partial^\alpha f\|_{L^2(\Omega)},\quad\forall\,f\in H^s(\Omega).
\end{equation}
Continue to assume that $\Omega\subset\mathbb{R}^n$ is a bounded Lipschitz domain and, 
for each $s\in\bbR$, define
\begin{align}\label{u54ruy}
\begin{split}
& H_0^s(\Omega):=\big\{f\in H^s(\mathbb{R}^n)\,\big|\,
{\rm supp}\,f\subseteq\overline{\Omega}\,\big\}
\\[2pt] 
& \quad\text{viewed as a closed subspace of $H^s(\mathbb{R}^n)$.}
\end{split}
\end{align}
Then \eqref{eq:DDEEn.2} (used with $\Omega:={\mathbb{R}}^n$) implies
\begin{align}\label{eq:DDEEn.2.WACO}
H^{s_1}_0(\Omega)\hookrightarrow H^{s_2}_0(\Omega)
\,\text{ continuously, if }\,s_1,s_2\in{\mathbb{R}},\,\,s_1\geq s_2.
\end{align}
In addition, if $\widetilde{C_0^\infty(\Om)}$ denotes the set of functions from 
$C_0^\infty(\Om)$ extended to all of $\bbR^n$ by zero outside their supports, then 
(cf. \cite[Remark~2.7, p.~170]{JK95})
\begin{equation}\label{eq:12ddC}
\widetilde{C_0^\infty(\Om)}\hookrightarrow H^s_0(\Om)
\,\text{ densely for each }\,s\in\bbR.
\end{equation}

For each $s\in{\mathbb{R}}$, it is of interest to also introduce 
\begin{equation}\label{restr-0}
\begin{array}{l}
H^s_{z}(\Omega):=\big\{u\in{\mathcal{D}}'(\Omega)\,\big|\,\text{there exists }\,f\in H^s_0(\Omega)
\,\text{ with }\,f\big|_{\Omega}=u\big\},
\\[6pt]
\|u\|_{H^s_{z}(\Omega)}:=\inf\,\big\{\|f\|_{H^s({\mathbb{R}}^n)}\,\big|\,
f\in H^s_0(\Omega),\,f\big|_{\Omega}=u\big\},\,\,\forall\,u\in H^s_{z}(\Omega).
\end{array}
\end{equation}
In particular, 
\begin{align}\label{restr-0.EMB}
\begin{split} 
& H^s_{z}(\Omega)=\big\{f\big|_{\Omega}\,\big|\,f\in H^s_0(\Omega)\big\}
\subseteq H^s(\Omega)\,\text{ and the inclusion}
\\[2pt]  
& \quad H^s_{z}(\Omega)\hookrightarrow H^s(\Omega)\,\text{ is continuous for each }\,s\in{\mathbb{R}}.
\end{split}
\end{align}
As is apparent from definitions, the operator of restriction (in the sense of distributions) 
$H^s({\mathbb{R}}^n)\ni f\mapsto f|_{\Omega}\in H^s(\Omega)$ 
maps $H^s_0(\Omega)$ continuously onto $H^s_{z}(\Omega)$ for each $s\in{\mathbb{R}}$. 
Together with \eqref{eq:12ddC} this implies that 
\begin{eqnarray}\label{F-dense-3}
C^\infty_0(\Omega)\hookrightarrow H^s_{z}(\Omega)\,\text{ densely, for each }\,s\in{\mathbb{R}}.
\end{eqnarray} 

We also record the identification (cf. the discussion in \cite{JK95}, \cite{MM13})
\begin{equation}\label{eq:Redxax.1}
\big(H^s(\Om)\big)^*=H^{-s}_0(\Om),\quad\forall\,s\in\bbR, 
\end{equation}
where each $V\in H^{-s}_0(\Om)$ is identified with the functional 
\begin{equation}\label{jussi:pairing}
\big(H^s(\Om)\big)^*\ni u\mapsto V(u):={}_{(H^s(\Om))^*}\big\langle\overline V,u\big\rangle_{H^s(\Om)}
\end{equation}
acting on an arbitrary $u\in H^s(\Om)$ according to the (unambiguous, due to \eqref{eq:12ddC}) recipe: 
\begin{align}\label{eq:Redxax.2}
\begin{split} 
& {}_{(H^s(\Om))^*}\big\langle\overline V,u\big\rangle_{H^s(\Om)}:=
{}_{H^{-s}({\mathbb{R}}^n)}\big\langle\overline V,U\big\rangle_{H^s({\mathbb{R}}^n)},     
\,\text{where $U$ is} 
\\[2pt]
& \quad\text{any distribution in $H^s({\mathbb{R}}^n)$ such that $U\big|_{\Omega}=u$},
\end{split}
\end{align} 
where ${}_{H^{-s}({\mathbb{R}}^n)}\big\langle\cdot,\cdot\big\rangle_{H^s({\mathbb{R}}^n)}$ is the canonical duality 
pairing between distributions in $H^{-s}({\mathbb{R}}^n)$ and, respectively, $H^s({\mathbb{R}}^n)=(H^{-s}({\mathbb{R}}^n))^\ast$. 
Moreover, if $\psi\in C^\infty(\overline{\Omega})$ and $u\in H^s_0(\Omega)$ for some 
$s\in{\mathbb{R}}$, then $\psi u:=\Psi u$ (considered in the sense of distributions), where
$\Psi\in C^\infty({\mathbb{R}}^n)$ is any smooth extension of $\psi$, is unambiguously defined
(due to \eqref{eq:12ddC}), belongs to $H^{s}_0(\Omega)$, and for every $v\in H^{-s}(\Omega)$ 
one has
\begin{equation}\label{tvdee}
{}_{H^{s}_0(\Omega)}\big\langle\psi u,v\big\rangle_{H^{-s}(\Omega)}
={}_{H^{s}_0(\Omega)}\big\langle u,\overline{\psi}v\big\rangle_{H^{-s}(\Omega)}.
\end{equation}

Since $H^s(\Om)$ is a reflexive Banach space for each $s\in\bbR$ (again, see the discussion in 
\cite{JK95}, \cite{MM13}), from \eqref{eq:Redxax.1} we also conclude that 
\begin{equation}\label{eq:Redxax.1REF}
\big(H^{s}_0(\Om)\big)^*=H^{-s}(\Om),\quad\forall\,s\in\bbR.
\end{equation}

For future use we note the identification
\begin{equation}\label{u5iiui} 
\big(H^{s}(\Omega)\big)^*=H^{-s}(\Omega),\,\text{ whenever }\,-\tfrac{1}{2}<s<\tfrac{1}{2},
\end{equation} 
in the sense that
\begin{align}\label{u5iiui-bbb} 
\begin{split}
& \big(H^{s}(\Omega)\big)^*=H^{-s}_0(\Omega)\ni u\mapsto u\big|_{\Omega}\in H^{-s}(\Omega)
\\[2pt]
& \quad\text{is an isomorphism whenever }\,-\tfrac{1}{2}<s<\tfrac{1}{2}.
\end{split}
\end{align} 
Furthermore, if tilde denotes the extension of a function 
to the entire Euclidean space by zero outside its original domain, then
\begin{equation}\label{eq:12d-yw5}
\parbox{10.70cm}
{for each $s\in\big(-\tfrac12,\tfrac12\big)$, the inclusion 
$C_0^\infty(\Omega)\hookrightarrow H^s(\Omega)$ has dense range, 
and the assignment $C^\infty(\overline{\Omega})\ni\varphi\mapsto\widetilde{\varphi}\in H^s_0(\Omega)$
extends by density to a linear and bounded isomorphism from $H^s(\Omega)$ onto $H^s_0(\Omega)$, which is 
the inverse of the restriction map $H^s_0(\Omega)\ni u\mapsto u\big|_{\Omega}\in H^s(\Omega)$
(cf.~\eqref{u5iiui-bbb}).}
\end{equation}
It has been shown in \cite{MM13} that 
\begin{equation}\label{equal.HH}
\accentset{\circ}{H}^s(\Omega)=H^s_{z}(\Omega)\,\text{ whenever }\,s>-\tfrac{1}{2}\,\text{ and }\,s-\tfrac{1}{2}\notin{\mathbb{N}}_0.
\end{equation}
As a consequence of this and \eqref{restr-0}, one therefore has 
\begin{equation}\label{u54r5t}
\accentset{\circ}{H}^{s}(\Omega)=\big\{f|_{\Omega}\,\big|\,f\in H^s_0(\Omega)\big\}
\,\text{ if }\,s>-\tfrac{1}{2}\,\text{ and }\,s-\tfrac{1}{2}\notin\mathbb{N}_0.
\end{equation}
From Lemma~2.2 in \cite{MiMo08} it follows that if $s\in\big(-\tfrac12,\tfrac12\big)$ 
then for each $u\in H^{s}(\Omega)$ and $v\in H^{-s}({\mathbb{R}}^n)$ one has 
\begin{align}\label{trans-U}
{}_{H^{s}({\mathbb{R}}^n)}\langle\widetilde{u},v\rangle_{H^{-s}({\mathbb{R}}^n)}
={}_{H^{s}(\Omega)}\big\langle u,v\big|_{\Omega}\big\rangle_{H^{-s}(\Omega)}.
\end{align}
Together, \eqref{Rdac}, the first line in \eqref{eq:12d-yw5}, and \eqref{u54r5t} also imply that 
\begin{equation}\label{Rdac.WACO}
H^s(\Omega)=\accentset{\circ}{H}^s(\Omega)=H^s_{z}(\Omega)\,\text{ for each }\,s\in\big(-\tfrac12,\tfrac12\big).
\end{equation} 
Later on, we shall use the fact that 
\begin{equation}\label{SDErr99-1}
\big\{u\in H^s({\mathbb{R}}^n)\,\big|\,{\rm supp}\,u\subseteq\partial\Omega\big\}=\{0\}\,\text{ whenever }\,s>-\tfrac{1}{2}.
\end{equation}
In addition, we shall need the following lifting result, valid for each $s>0$:
\begin{align}\label{eq:TRfav}
\begin{split}
& u\in H^{1+s}(\Omega)\,\text{ if and only if }\,u\in L^2(\Omega)
\,\text{ and }\,\nabla u\in\big[H^{s}(\Omega)\big]^n,
\\[2pt]
& \quad\text{and }\,\|u\|_{H^{1+s}(\Omega)}\approx\|u\|_{L^2(\Omega)}
+\big\|\nabla u\big\|_{[H^{s}(\Omega)]^n},\,\text{ uniformly in $u$}. 
\end{split}
\end{align}
See, for instance, \cite{JK95,Mc00,MM13,Tr02} for these and other related properties. 
We also note that when $\Omega$ is a bounded Lipschitz domain in $\bbR^n$ and $s\in(0,1)$,  
then there exists $C\in(1,\infty)$ such that for every $f\in H^s(\Omega)$ there holds 
\begin{equation}\label{FSA-2X}
C^{-1}\|f\|_{H^s(\Omega)}\leq\|f\|_{L^2(\Omega)}
+\bigg(\int_{\Omega}\int_{\Omega}\frac{|f(x)-f(y)|^2}{|x-y|^{n+2s}}\,d^n x d^n y\bigg)^{1/2}
\leq C\|f\|_{H^s(\Omega)}.
\end{equation}
See \cite{Di03}, \cite[Proposition~2.28, pp.~51--52]{MM13} for more general results of this nature. 

We continue by discussing a very useful density result, refining work in 
\cite{CD98}, \cite[Lemma~1.5.3.9, p.~59]{Gr85}, \cite[Theorem~6.4, Chapter~2]{LM72}. 

%%%%%%%%%%%
\begin{lemma}\label{Dense-LLLe}
Let $\Omega\subset\bbR^n$ be a bounded Lipschitz domain, and fix two arbitrary 
numbers $s_1,s_2\in{\mathbb{R}}$. Define 
\begin{equation}\label{gag-8gb.222}
H^{s_1,s_2}_{\Delta}(\Omega):=\big\{u\in H^{s_1}(\Omega)\,\big|\,\Delta u\in H^{s_2}(\Omega)\big\}
\end{equation}
and equip this space with the natural graph norm 
$u\mapsto\|u\|_{H^{s_1}(\Omega)}+\|\Delta u\|_{H^{s_2}(\Omega)}$. Then $H^{s_1,s_2}_{\Delta}(\Omega)$ becomes a Banach space and
\begin{equation}\label{gag-8gb}
 C^\infty(\overline{\Omega})\subset H^{s_1,s_2}_{\Delta}(\Omega)\,\text{ densely}.
\end{equation}
\end{lemma}
%%%%%%%%%%%
\begin{proof}
To check that $H^{s_1,s_2}_{\Delta}(\Omega)$ is complete, assume $\{u_j\}_{j\in{\mathbb{N}}}\subseteq H^{s_1,s_2}_{\Delta}(\Omega)$ 
is a Cauchy sequence (with respect to the graph norm described in the statement). Then  
$\{u_j\}_{j\in{\mathbb{N}}}$ is a Cauchy sequence in $H^{s_1}(\Omega)$ and
$\{\Delta u_j\}_{j\in{\mathbb{N}}}$ is a Cauchy sequence in $H^{s_2}(\Omega)$. Since 
both $H^{s_1}(\Omega)$ and $H^{s_2}(\Omega)$ are complete, this implies that there exist 
$u\in H^{s_1}(\Omega)$ and $w\in H^{s_2}(\Omega)$ such that 
$\{u_j\}_{j\in{\mathbb{N}}}$ converges to $u$ in $H^{s_1}(\Omega)$ and
$\{\Delta u_j\}_{j\in{\mathbb{N}}}$ converges to $w$ in $H^{s_2}(\Omega)$. 
Given that both $H^{s_1}(\Omega)$ and $H^{s_2}(\Omega)$ embed continuously into ${\mathcal{D}}'(\Omega)$
(itself, a Hausdorff topological vector space), and that $\Delta:{\mathcal{D}}'(\Omega)\to{\mathcal{D}}'(\Omega)$ is continuous, 
we may then conclude that $\Delta u=w$ in ${\mathcal{D}}'(\Omega)$. Hence, $u\in H^{s_1,s_2}_{\Delta}(\Omega)$ and 
$\{u_j\}_{j\in{\mathbb{N}}}$ converges to $u$ in $H^{s_1,s_2}_{\Delta}(\Omega)$. This proves that 
$H^{s_1,s_2}_{\Delta}(\Omega)$ is indeed a Banach space when equipped with the graph norm. 

Moving on, in the case when $s_1-s_2\geq 2$ it follows from \eqref{eq:DDEEn.2}--\eqref{eq:DDEEn.3} that
$H^{s_1,s_2}_{\Delta}(\Omega)$ and $H^{s_1}(\Omega)$ coincide as vector spaces and have equivalent norms. 
Hence, in this scenario, the claim in \eqref{gag-8gb} is a direct consequence of \eqref{uam-mi78}.

Next, consider the situation where $s_1,s_2\in{\mathbb{R}}$ satisfy 
$s_1-s_2<2$. To proceed, define the isometric embedding 
\begin{align}\label{eq:IS-EM.1}
\iota: 
\begin{cases}
H^{s_1,s_2}_{\Delta}(\Omega)\rightarrow H^{s_1}(\Omega)\oplus H^{s_2}(\Omega),    
\\[2pt] 
u\mapsto\iota(u):=(u,\Delta u), 
\end{cases}
\end{align}
and note that its image, $\ran(\iota)$, is a closed subspace of $H^{s_1}(\Omega)\oplus H^{s_2}(\Omega)$.
In particular, $\iota:H^{s_1,s_2}_{\Delta}(\Omega)\rightarrow\ran(\iota)$ is a continuous isomorphism, and we 
denote by $\iota^{-1}:\ran(\iota)\rightarrow H^{s_1,s_2}_{\Delta}(\Omega)$ its inverse. Let now 
$\Lambda:H^{s_1,s_2}_{\Delta}(\Omega)\to{\mathbb{C}}$ be an arbitrary continuous functional. 
Then $\Lambda\circ\iota^{-1}$ is a continuous functional on the closed subspace $\ran(\iota)$
of the Banach space $H^{s_1}(\Omega)\oplus H^{s_2}(\Omega)$. As such, the Hahn--Banach theorem 
ensures that this extends to a functional (cf.~\eqref{eq:Redxax.1})
\begin{align}\label{eq:fccaRF}
\widehat{\Lambda}\in\big(H^{s_1}(\Omega)\oplus H^{s_2}(\Omega)\big)^\ast
&\equiv\big(H^{s_1}(\Omega)\big)^\ast\oplus\big(H^{s_2}(\Omega)\big)^\ast
\nonumber\\[2pt]
&=H^{-s_1}_0(\Omega)\oplus H^{-s_2}_0(\Omega).
\end{align}
Together with \eqref{eq:Redxax.1}--\eqref{eq:Redxax.2}, this implies that there exist 
\begin{equation}\label{eq:ACD}
h_1\in H^{-s_1}_0(\Omega)\,\text{ and }\,h_2\in H^{-s_2}_0(\Omega)
\end{equation}
with the property that for each $u\in H^{s_1,s_2}_{\Delta}(\Omega)$ one has 
\begin{equation}\label{HGav-jt6t}
\Lambda(u)={}_{H^{-s_1}({\mathbb{R}}^n)}\big\langle\overline h_1,F\big\rangle_{H^{s_1}({\mathbb{R}}^n)}
+{}_{H^{-s_2}({\mathbb{R}}^n)}\big\langle\overline h_2,G\big\rangle_{H^{s_2}({\mathbb{R}}^n)}, 
\end{equation}
whenever 
\begin{equation}\label{e:jFa}
F\in H^{s_1}({\mathbb{R}}^n)\,\text{ and }\,G\in H^{s_2}({\mathbb{R}}^n)
\,\text{ satisfy }\,F\big|_{\Omega}=u,\quad G\big|_{\Omega}=\Delta u.
\end{equation}
To proceed, we consider an arbitrary $\varphi\in C^\infty_0({\mathbb{R}}^n)$ and note that 
\eqref{HGav-jt6t}--\eqref{e:jFa} applied with $u:=\varphi\big|_{\Omega}\in H^{s_1,s_2}_{\Delta}(\Omega)$, 
$F:=\varphi\in H^{s_1}({\mathbb{R}}^n)$, and $G:=\Delta\varphi\in H^{s_2}({\mathbb{R}}^n)$, yields 
\begin{align}\label{HGav-jt6t.32}
\Lambda\big(\varphi\big|_{\Omega}\big) &={}_{H^{-s_1}({\mathbb{R}}^n)}
\big\langle\overline h_1,\varphi\big\rangle_{H^{s_1}({\mathbb{R}}^n)}
+{}_{H^{-s_2}({\mathbb{R}}^n)}\big\langle\overline h_2,\Delta\varphi\big\rangle_{H^{s_2}({\mathbb{R}}^n)}
\nonumber\\[2pt]
&={}_{\cD'({\mathbb{R}}^n)}\big\langle h_1,\varphi\big\rangle_{\cD({\mathbb{R}}^n)}
+{}_{\cD'({\mathbb{R}}^n)}\big\langle h_2,\Delta\varphi\big\rangle_{\cD({\mathbb{R}}^n)}
\end{align} 
which ultimately leads to the conclusion that 
\begin{equation}\label{eq:Nvav.1}
\Lambda\big(\varphi\big|_{\Omega}\big)
={}_{\cD'({\mathbb{R}}^n)}\big\langle h_1+\Delta h_2,\varphi\big\rangle_{\cD({\mathbb{R}}^n)},
\quad\forall\,\varphi\in C^\infty_0({\mathbb{R}}^n).
\end{equation}

Next, make the assumption that 
\begin{equation}\label{eq:GfapkII.1}
\Lambda(v)=0\,\text{ for each }\,v\in C^\infty(\overline{\Omega}),
\end{equation}
and note that, by virtue of \eqref{eq:Nvav.1}, this forces 
\begin{equation}\label{eq:GfapkII.2}
h_1+\Delta h_2=0\,\text{ in }\,\mathcal{D}'({\mathbb{R}}^n). 
\end{equation}
Hence, $\Delta h_2=-h_1\in H^{-s_1}({\mathbb{R}}^n)$ so $h_2\in H^{2-s_1}({\mathbb{R}}^n)$ by 
elliptic regularity. Moreover, since $h_2\in H^{-s_2}_0(\Omega)$ entails 
${\rm supp}\,(h_2)\subseteq\overline{\Omega}$, one actually has $h_2\in H^{2-s_1}_0(\Omega)$.
This fact and \eqref{eq:12ddC} imply the existence of a sequence 
$\{\phi_j\}_{j\in{\mathbb{N}}}\subset C^\infty_0(\Omega)$ with the property that
\begin{equation}\label{eq:Ascv-1}
\widetilde{\phi_j}\to h_2\,\text{ in }\, H^{2-s_1}({\mathbb{R}}^n)
\,\text{ as }\, j\to\infty,
\end{equation}
where tilde denotes the extension by zero outside the support to the entire ${\mathbb{R}}^n$.
In turn, from \eqref{eq:Ascv-1}, \eqref{eq:DDEEn.3}, and \eqref{eq:GfapkII.2} one deduces that
\begin{equation}\label{eq:Ascv-2}
\widetilde{\Delta\phi_j}=
\Delta(\widetilde{\phi_j})\to\Delta h_2=-h_1\,\text{ in }\, H^{-s_1}({\mathbb{R}}^n)
\,\text{ as }\, j\to\infty.
\end{equation}
In addition, our assumption $s_1-s_2<2$ implies 
$H^{2-s_1}({\mathbb{R}}^n)\hookrightarrow H^{-s_2}({\mathbb{R}}^n)$
(cf.~\eqref{eq:DDEEn.2}) which, together with \eqref{eq:Ascv-1}, also implies 
\begin{equation}\label{eq:Ascv-3}
\widetilde{\phi_j}\to h_2\,\text{ in }\,H^{-s_2}({\mathbb{R}}^n)
\,\text{ as }\,j\to\infty.
\end{equation}
Pick now an arbitrary $u\in H^{s_1,s_2}_{\Delta}(\Omega)$ and let $F,G$ be as in \eqref{e:jFa}. 
Then based on \eqref{HGav-jt6t}, \eqref{eq:Ascv-2}--\eqref{eq:Ascv-3}, and \eqref{e:jFa}, one can write 
\begin{align}\label{HGav-jt6t.77}
\Lambda(u) &={}_{H^{-s_1}({\mathbb{R}}^n)}\big\langle\overline h_1,F\big\rangle_{H^{s_1}({\mathbb{R}}^n)}
+{}_{H^{-s_2}({\mathbb{R}}^n)}\big\langle\overline h_2,G\big\rangle_{H^{s_2}({\mathbb{R}}^n)}
\nonumber\\[2pt]
&=\lim_{j\to\infty}\Big\{
{}_{H^{-s_1}({\mathbb{R}}^n)}\Big\langle -\overline{\widetilde{\Delta\phi_j}},F\Big\rangle_{H^{s_1}({\mathbb{R}}^n)}
+{}_{H^{-s_2}({\mathbb{R}}^n)}\Big\langle\overline{\widetilde{\phi_j}},G\Big\rangle_{H^{s_2}({\mathbb{R}}^n)}\Big\}
\nonumber\\[2pt]
&=\lim_{j\to\infty}\big\{
{}_{\cD({\mathbb{R}}^n)}\big\langle -\widetilde{\Delta\phi_j},F\big\rangle_{\cD'({\mathbb{R}}^n)}
+{}_{\cD({\mathbb{R}}^n)}\big\langle\widetilde{\phi_j},G\big\rangle_{\cD'({\mathbb{R}}^n)}\big\}
\nonumber\\[2pt]
&=\lim_{j\to\infty}\big\{
{}_{\cD(\Omega)}\big\langle -\Delta\phi_j,F\big|_{\Omega}\big\rangle_{\cD'(\Omega)}
+{}_{\cD(\Omega)}\big\langle\phi_j,G|_{\Omega}\big\rangle_{\cD'(\Omega)}\big\}
\nonumber\\[2pt] 
&=\lim_{j\to\infty}\big\{
{}_{\cD(\Omega)}\big\langle-\Delta\phi_j,u\big\rangle_{\cD'(\Omega)}
+{}_{\cD(\Omega)}\big\langle\phi_j,\Delta u\big\rangle_{\cD'(\Omega)}\big\}
\nonumber\\[2pt] 
&=0.
\end{align} 
This shows that any linear and continuous functional $\Lambda$ on $H^{s_1,s_2}_{\Delta}(\Omega)$ 
satisfying \eqref{eq:GfapkII.1} ultimately vanishes identically, from which the claim 
in \eqref{gag-8gb} readily follows. This finishes the proof of Lemma~\ref{Dense-LLLe}.
\end{proof}
%%%%%%%%%%

For later purposes a variant of Lemma~\ref{Dense-LLLe} with the Sobolev space $H^{s_2}(\Omega)$ replaced by 
a suitable Besov space will be useful.

%%%%%%%%%%%
\begin{lemma}\label{Dense-LLLe-BBB}
Let $\Omega\subset\bbR^n$ be a bounded Lipschitz domain, and fix an arbitrary 
number $s\in{\mathbb{R}}$. Define the hybrid space
\begin{equation}\label{gag-8gb.222-BBB}
H\!B^{s}_{\Delta}(\Omega):=\big\{u\in H^{s}(\Omega)\,\big|\,\Delta u\in B^{2,1}_{s-2}(\Omega)\big\}
\end{equation}
and equip it with the natural graph norm 
$u\mapsto\|u\|_{H^{s}(\Omega)}+\|\Delta u\|_{B^{2,1}_{s-2}(\Omega)}$. Then
\begin{equation}\label{gag-8gb-Mi}
C^\infty(\overline{\Omega})\subset H\!B^{s}_{\Delta}(\Omega)\,\text{ densely}.
\end{equation}
\end{lemma}
%%%%%%%%%%%
\begin{proof}
Pick an arbitrary function $u\in H\!B^{s}_{\Delta}(\Omega)$ and 
extend $\Delta u\in B^{2,1}_{s-2}(\Omega)$ to a compactly supported distribution 
$U\in B^{2,1}_{s-2}(\mathbb{R}^n)$. 
Let $E_0$ denote 
the standard fundamental solution for the Laplacian in ${\mathbb{R}}^n$, that is,
\begin{equation}\label{EEE-jussi}
E_0(x):=\begin{cases} 
\displaystyle{\frac{1}{\omega_{n-1}(2-n)}|x|^{2-n}},\,\text{  if }\,n\geq 3,
\\[12pt] 
\displaystyle{\frac{1}{2\pi}}\,\ln|x|,\,\text{ if }\,n=2,
\end{cases}
\quad\forall\,x\in{\mathbb{R}}^{n}\backslash\{0\},
\end{equation}
\noindent where $\omega_{n-1}$ is the surface measure of the unit sphere
$\bbS^{n-1}$ in ${\mathbb{R}}^n$. Classical Calder\'on--Zygmund theory gives that the 
operator of convolution with $E_0$ is (locally) smoothing of order two on the fractional 
Besov scale (see, e.g., the discussion in \cite[Section~4]{KMM07}). 
Hence, considering, $\eta:=(E_0\ast U)|_{\Omega}$ then 
\begin{equation}\label{aaa-Mi-BB}
\eta\in B^{2,1}_{s}(\Omega)\subseteq B^{2,2}_{s}(\Omega)=H^s(\Omega)
\end{equation}
and $\Delta\eta=(\Delta E_0\ast U)|_{\Omega}=U|_{\Omega}=\Delta u$ in $\Omega$. Considering  
$v:=u-\eta$, then $v\in H^s(\Omega)$ and $\Delta v=0$ in $\Omega$. In the notation introduced 
in \eqref{gag-8gb.222}, this implies
\begin{equation}\label{gag-8gb.222-BBB.22}
v\in H^{s,s_\ast}_{\Delta}(\Omega)\,\text{ for each }\,s_\ast\in{\mathbb{R}}. 
\end{equation}
To proceed, fix a real number $s_\ast$ satisfying 
\begin{equation}\label{gag-8gb.222-BBB.33}
s_\ast>s-2
\end{equation}
and invoke Lemma~\ref{Dense-LLLe} to produce a sequence 
$\{v_j\}_{j\in{\mathbb{N}}}\subset C^\infty(\overline{\Omega})$ with the property that
\begin{equation}\label{ajtrtv5-Mi}
v_j\to v\,\text{ in }\,H^s(\Omega)\,\text{ and }\,\Delta v_j\to 0
\,\text{ in }\,H^{s_\ast}(\Omega),\,\text{ as }\,j\to\infty.
\end{equation}
In light of \eqref{gag-8gb.222-BBB.33} and \eqref{q-TWU-yt}, the last convergence above also implies
\begin{equation}\label{ajtrtv5-Mi.2}
\Delta v_j\to 0\,\text{ in }\,B^{2,1}_{s-2}(\Omega),\,\text{ as }\,j\to\infty.
\end{equation}
On the other hand, from \eqref{q-TWU.ted-DD}, the fact that  
$U$ is a compactly supported distribution in ${\mathbb{R}}^n$, and 
${\mathcal{S}}({\mathbb{R}}^n)\subset B^{2,1}_{s-2}({\mathbb{R}}^n)$ densely,
one deduces that there exists a sequence $\{\phi_j\}_{j\in{\mathbb{N}}}\subset C^\infty_0({\mathbb{R}}^n)$ 
with supports contained in a common compact subset of ${\mathbb{R}}^n$ and such that
\begin{equation}\label{ajtrtv5-Mi.4}
\phi_j\to U\,\text{ in }\,B^{2,1}_{s-2}({\mathbb{R}}^n)\,\text{ as }\,j\to\infty.
\end{equation}
If for each $j\in{\mathbb{N}}$ we now define $\eta_j:=(E_0\ast\phi_j)|_{\Omega}$, 
then $\eta_j\in C^\infty(\overline{\Omega})$ and 
\begin{align}\label{aaa-Mi-BB.ee}
\begin{split}
& \eta_j\to(E_0\ast\phi)|_{\Omega}=\eta\,\text{ in }\,B^{2,1}_{s}(\Omega),
\,\text{ hence}
\\[2pt]
& \quad\text{also in }\,B^{2,2}_{s}(\Omega)=H^{s}(\Omega),\,\text{ as }\,j\to\infty.
\end{split}
\end{align}
In addition, 
\begin{align}\label{aaa-Mi-BB.ff}
\Delta\eta_j &=(\Delta E_0\ast\phi_j)|_{\Omega}
\nonumber\\[2pt]
&=\phi_j|_{\Omega}\to U|_{\Omega}=\Delta u\,\text{ in }\,B^{2,1}_{s-2}(\Omega)\,\text{ as }\,j\to\infty.
\end{align}
Next, consider $\psi_j:=v_j+\eta_j\in C^\infty(\overline{\Omega})$ for each $j\in{\mathbb{N}}$.
Then from \eqref{ajtrtv5-Mi} and \eqref{aaa-Mi-BB.ee} one concludes that  
\begin{equation}\label{ajtrtv5-Mi.6y}
\psi_j\to v+\eta=u\,\text{ in }\,H^s(\Omega)\,\text{ as }\,j\to\infty,
\end{equation}
while from \eqref{ajtrtv5-Mi.2} and \eqref{aaa-Mi-BB.ff} one infers that 
\begin{equation}\label{ajtrtv5-Mi.7y}
\Delta\psi_j=\Delta v_j+\Delta\eta_j\to\Delta u\,\text{ in }\,B^{2,1}_{s-2}(\Omega)\,\text{ as }\,j\to\infty.
\end{equation}
In view of the nature of the norm on the space $H\!B^{s}_{\Delta}(\Omega)$, this finishes the 
proof of \eqref{gag-8gb-Mi}.
\end{proof}
%%%%%%%

Loosely speaking, the result in the proposition below may be interpreted as saying that, 
for a function $u$ belonging to a Triebel--Lizorkin space in a bounded Lipschitz domain, 
having a ``better-than-expected'' Laplacian $\Delta u$ (again, measured on the Triebel--Lizorkin scale) 
translates into better regularity for the function $u$ than originally assumed. 

%%%%%%%%%%%
\begin{proposition}\label{Dense-LLLe-BBB.WACO}
Let $\Omega\subset\bbR^n$ be a bounded Lipschitz domain and fix some integrability exponents $0<p,q_0,q_1<\infty$ 
along with a smoothness index $s\in{\mathbb{R}}$. Suppose $u\in F^{p,q_0}_s(\Omega)$ is a function with the property 
that $\Delta u\in F^{p,q_1}_{s-2}(\Omega)$. Then $u$ belongs to $F^{p,q_1}_s(\Omega)$ and there exists a constant 
$C\in(0,\infty)$ which is independent of $u$ such that 
\begin{equation}\label{eq:H.WACO.1}
\|u\|_{F^{p,q_1}_s(\Omega)}\leq C\big(\|u\|_{F^{p,q_0}_s(\Omega)}+\|\Delta u\|_{F^{p,q_1}_{s-2}(\Omega)}\big).
\end{equation}
\end{proposition}
%%%%%%%%%%%
\begin{proof}
Since for $q_0\leq q_1$ the desired conclusions follow directly from \eqref{q-TWC.WACO}, 
it remains to treat the case $q_1<q_0$. In view of \eqref{restr.WACO}, \eqref{q-TWU.ted-DD.WACO}, it is possible to extend 
$\Delta u\in F^{p,q_1}_{s-2}(\Omega)$ to a compactly supported distribution $U\in F^{p,q_1}_{s-2}(\mathbb{R}^n)$ satisfying 
$\|U\|_{F^{p,q_1}_{s-2}(\mathbb{R}^n)}\leq C\|\Delta u\|_{F^{p,q_1}_{s-2}(\Omega)}$ for some $C\in(0,\infty)$ independent of $u$. 
Let $E_0$ denote the standard fundamental solution for the Laplacian in ${\mathbb{R}}^n$ (cf.~\eqref{EEE-jussi}).
Then the operator of convolution with $E_0$ is (locally) smoothing of order two on the Triebel--Lizorkin scale 
(cf., e.g., \cite[Section~4]{KMM07}). As such, if we consider $w:=(E_0\ast U)|_{\Omega}$ then 
\begin{equation}\label{aaa-Mi-BB.WACO}
w\in F^{p,q_1}_{s}(\Omega)\hookrightarrow F^{p,q_0}_{s}(\Omega),
\end{equation}
with the continuous inclusion provided by \eqref{q-TWC.WACO}, and 
\begin{equation}\label{aaa-Mi-BB.WACO.22}
\|w\|_{F^{p,q_0}_{s}(\Omega)}\leq C\|w\|_{F^{p,q_1}_{s}(\Omega)}\leq C\|U\|_{F^{p,q_1}_{s-2}(\mathbb{R}^n)}
\leq C\|\Delta u\|_{F^{p,q_1}_{s-2}(\Omega)}.
\end{equation}
In addition, one has $\Delta w=(\Delta E_0\ast U)|_{\Omega}=U|_{\Omega}=\Delta u$ in $\Omega$. 
Consequently, if we introduce $v:=u-w$, then $v\in F^{p,q_0}_s(\Omega)$ satisfies $\Delta v=0$ in $\Omega$ and 
\begin{align}\label{eq:H.WACO.2}
\|v\|_{F^{p,q_0}_s(\Omega)} &\leq C\big(\|u\|_{F^{p,q_0}_s(\Omega)}+\|w\|_{F^{p,q_0}_s(\Omega)}\big)
\nonumber\\[2pt]
&\leq C\big(\|u\|_{F^{p,q_0}_s(\Omega)}+\|\Delta u\|_{F^{p,q_1}_{s-2}(\Omega)}\big).
\end{align}
Next, we recall from \cite[Theorem~1.6]{KMM07} that 
\begin{equation}\label{Ind-X} 
\parbox{9.00cm}{the space of harmonic functions in $F^{p,q}_s(\Omega)$ is actually independent of the index $q\in(0,\infty)$
and all quasi-norms $\|\cdot\|_{F^{p,q}_s(\Omega)}$ with $q\in(0,\infty)$ are equivalent when considered 
on the space of harmonic functions in $\Omega$.}
\end{equation}
We then conclude that $v$ belongs to $F^{p,q_1}_s(\Omega)$ and satisfies
\begin{equation}\label{eq:H.WACO.3}
\|v\|_{F^{p,q_1}_s(\Omega)}\leq C\|v\|_{F^{p,q_0}_s(\Omega)}.
\end{equation}
Hence, $u=v+w\in F^{p,q_1}_s(\Omega)$ and \eqref{eq:H.WACO.2}, \eqref{eq:H.WACO.3}, \eqref{aaa-Mi-BB.WACO.22} 
prove that \eqref{eq:H.WACO.1} holds. 
\end{proof}
%%%%%%%

Our next lemma brings to light the compatibility of the Sobolev space pairing with the 
ordinary integral pairing, when both turn out to be meaningful. Given an open set 
$\Omega\subseteq{\mathbb{R}}^n$, we denote by $L^\infty_{\rm loc}(\Omega)$ the space 
of measurable functions defined in $\Omega$ which become essentially bounded when 
restricted to compact subsets of $\Omega$. In addition, for each $p\in(0,\infty]$,
we let $L^p_{\rm comp}(\Omega)$ stand for the subspace of $L^p(\Omega)$ consisting of 
functions with compact support in $\Omega$.

%%%%%%%
\begin{lemma}\label{utrCOMPA}
Assume that $\Omega$ is a bounded Lipschitz domain in ${\mathbb{R}}^n$ and fix some 
$s\in\big(-\tfrac12,\tfrac12\big)$. Then 
\begin{align}\label{khyg-886-AAA}
{}_{H^{s}(\Omega)}\big\langle u,v\big\rangle_{H^{-s}(\Omega)}
=\int_{\Omega}\overline{u(x)}v(x)\,d^nx
\end{align}
provided either 
\begin{equation}\label{ii76g4r-AAA}
u\in H^s(\Omega)\cap L^1_{\rm loc}(\Omega)\,\text{ and }\,v\in H^{-s}(\Omega)\cap L^\infty_{\rm comp}(\Omega),
\end{equation}
or 
\begin{equation}\label{ii76g4r-AAA.2}
u\in H^s(\Omega)\cap L^1_{\rm comp}(\Omega)\,\text{ and }\,v\in H^{-s}(\Omega)\cap L^\infty_{\rm loc}(\Omega). 
\end{equation}
\end{lemma}
%%%%%%%
\begin{proof}  
Let $\eta\in C^\infty_0({\mathbb{R}}^n)$ be a real-valued, even function,
satisfying $\eta=1$ on $B(0,1)$, $\eta=0$ outside $B(0,2)$,
$\int_{{\mathbb{R}}^n}\eta(x)\,d^n x=1$. In addition, for each $t>0$, set 
$\eta_t(x):=t^{-n}\eta(x/t)$ 
for each $x\in{\mathbb{R}}^n$. For each $t>0$, consider the operator 
\begin{equation}\label{trew54.d}
I_t:{\mathcal{D}}'({\mathbb{R}}^n)\to C^\infty({\mathbb{R}}^n),\quad
I_tu:=u\ast\eta_t,\quad\forall\,u\in{\mathcal{D}}'({\mathbb{R}}^n).
\end{equation}
Then $I_t$ is bounded on $L^2({\mathbb{R}}^n)$ for each $t>0$ with operator norm 
controlled independently of $t$, and for each $u\in L^2({\mathbb{R}}^n)$ one has  
$I_tu\to u$ as $t\to 0_{+}$ in $L^2({\mathbb{R}}^n)$. Moreover, given any $k\in{\mathbb{N}}$, 
if $u\in H^k({\mathbb{R}}^n)$ and $\alpha$ is a multi-index of length at most $k$, then 
\begin{equation}\label{eq:JBa6gg-BIS}
\partial^\alpha(I_tu)=(\partial^\alpha u)\ast\eta_t
\to\partial^\alpha u\,\text{ as }\,t\to 0_{+}\,\text{ in }\,L^2({\mathbb{R}}^n).
\end{equation}
As a consequence, it follows that $I_t$ is bounded on $H^k({\mathbb{R}}^n)$ for each $t>0$ 
with operator norm controlled independently on $t$. Hence, by interpolation, for each $t>0$
the operator $I_t$ is bounded on any $H^s({\mathbb{R}}^n)$ with $s\geq 0$,  
with operator norm controlled independently of $t$.

Next, consider an arbitrary number $s>0$ and pick $k\in{\mathbb{N}}$, $k>s$, and $\theta\in(0,1)$ 
such that $s=\theta k$. Then for every $u\in C^\infty_0({\mathbb{R}}^n)$, the interpolation 
inequality 
\begin{equation}\label{Jvav-utrf-BIS}
\|I_tu-u\|_{H^s({\mathbb{R}}^n)}\leq\|I_tu-u\|_{H^k({\mathbb{R}}^n)}^\theta
\|I_tu-u\|_{L^2({\mathbb{R}}^n)}^{1-\theta}
\end{equation}
ultimately proves (in light of the density of $C^\infty_0({\mathbb{R}}^n)$ in 
$H^s({\mathbb{R}}^n)$) that if $s\geq 0$ then
\begin{equation}\label{eq:Hac-BIS}
I_tu\to u\,\text{ as }\,t\to 0_{+}\,\text{ in }\,H^s({\mathbb{R}}^n),\quad\forall\,u\in H^s({\mathbb{R}}^n).
\end{equation}
Moreover, for each $u,v\in C^\infty_0({\mathbb{R}}^n)$ one has 
$I_t u,\,I_tv\in C^\infty_0({\mathbb{R}}^n)$ and, given any $s\geq 0$, one can write (since $\eta_t$ is even)
\begin{align}\label{khyg-y4}
{}_{H^{-s}({\mathbb{R}}^n)}\big\langle I_tu,v\big\rangle_{H^s({\mathbb{R}}^n)}
&={}_{{\mathcal{D}}'({\mathbb{R}}^n)}\big\langle\,\overline{I_tu},v\big\rangle_{{\mathcal{D}}({\mathbb{R}}^n)}
\nonumber\\[2pt]
&=\int_{{\mathbb{R}}^n}\overline{(u\ast\eta_t)(x)}v(x)\,d^nx
=\int_{{\mathbb{R}}^n}\overline{u(x)}(v\ast\eta_t)(x)\,d^nx
\nonumber\\[2pt]
&={}_{{\mathcal{D}}'({\mathbb{R}}^n)}\big\langle\overline{u},I_tv\big\rangle_{{\mathcal{D}}({\mathbb{R}}^n)}
={}_{H^{-s}({\mathbb{R}}^n)}\big\langle u,I_tv\big\rangle_{H^s({\mathbb{R}}^n)}.
\end{align}
From \eqref{khyg-y4} one then concludes that 
\begin{align}\label{uttfvccc}
& \text{for any $s\in{\mathbb{R}}$ the operator $I_t$ induces a linear and bounded} 
\nonumber\\[2pt]
& \quad\text{mapping on $H^s({\mathbb{R}}^n)$ for each $t>0$, with operator norm} 
\\[2pt]
& \quad\text{controlled independently of $t$.}  
\nonumber 
\end{align}
One notes that \eqref{khyg-y4} also implies that for each $u,v\in C^\infty_0({\mathbb{R}}^n)$ and each 
$s\geq 0$ one has 
\begin{align}\label{khyg-y4-bis}
{}_{H^{-s}({\mathbb{R}}^n)}\big\langle u-I_tu,v\big\rangle_{H^s({\mathbb{R}}^n)}
={}_{H^{-s}({\mathbb{R}}^n)}\big\langle u,v-I_tv\big\rangle_{H^s({\mathbb{R}}^n)}.
\end{align}
On account of \eqref{khyg-y4-bis}, \eqref{uttfvccc}, and the density of $C^\infty_0({\mathbb{R}}^n)$ in 
$H^s({\mathbb{R}}^n)$, it follows that 
\begin{equation}\label{eq:WACO.1}
\text{\eqref{eq:Hac-BIS} actually holds for any $s\in{\mathbb{R}}$}.
\end{equation}

Next, let $\eta$ be as in the first part of the proof. 
For each fixed $t>0$, we now introduce the operator $J_t$ assigning to each 
$\varphi\in L^1(\Omega)$ the function 
\begin{equation}\label{trew54.a}
J_t\varphi:=\big(\widetilde{\varphi}\ast\eta_t\big)\big|_{\Omega}\in C^\infty(\overline{\Omega})
\subset L^\infty(\Omega)
\end{equation}
where, as usual, tilde denotes the extension to ${\mathbb{R}}^n$ by zero outside
$\Omega$. Then 
\begin{equation}\label{trew54.b}
J_t\varphi\to\varphi\,\text{ in }\,L^1(\Omega)\,\text{ as }\,t\to 0_{+},
\quad\forall\,\varphi\in L^1(\Omega),
\end{equation}
and one can easily check that, for each $t>0$, the operator $J_t$ satisfies
\begin{equation}\label{trew54.c}
\int_{\Omega}(J_t\varphi)(x)\psi(x)\,d^nx=\int_{\Omega}\varphi(x)(J_t\psi)(x)\,d^nx,
\quad\forall\,\varphi,\psi\in L^1(\Omega).
\end{equation}
In addition, by \eqref{eq:12d-yw5} and \eqref{eq:WACO.1}, 
\begin{equation}\label{eq:Hac-BIS.22}
J_tu\to u\,\text{ as }\,t\to 0_{+}\,\text{ in }\,H^s(\Omega),
\quad\forall\,u\in H^s(\Omega)\,\text{ with }\,s\in\big(-\tfrac12,\tfrac12\big).
\end{equation}
Assume that $s\in\big(-\tfrac12,\tfrac12\big)$ and fix $u,v$ as in \eqref{ii76g4r-AAA}.
Pick a real-valued function $\zeta\in C^\infty_0(\Omega)$ such that $\zeta=1$ in a neighborhood 
of ${\rm supp}\, (v)$. Given that $J_t v\in C^\infty_0(\Omega)$ for $t>0$ sufficiently 
small, and $\zeta u,v\in L^1(\Omega)$, one can write (here \eqref{tvdee} is relevant) 
\begin{align}\label{khyg-886}
{}_{H^{s}(\Omega)}\big\langle u,v\big\rangle_{H^{-s}(\Omega)}
&={}_{H^{s}(\Omega)}\big\langle u,\zeta v\big\rangle_{H^{-s}(\Omega)}
={}_{H^{s}(\Omega)}\big\langle\zeta u,v\big\rangle_{H^{-s}(\Omega)}
\nonumber\\[2pt]
&=\lim_{t\to 0_{+}}{}_{H^{s}(\Omega)}\big\langle\zeta u,J_tv\big\rangle_{H^{-s}(\Omega)}
=\lim_{t\to 0_{+}}{}_{{\mathcal{D}}'(\Omega)}\big\langle\,\overline{\zeta u},J_tv\big\rangle_{{\mathcal{D}}(\Omega)}
\nonumber\\[2pt]
&=\lim_{t\to 0_{+}}\int_{\Omega}\overline{(\zeta u)(x)}(J_tv)(x)\,d^nx
=\lim_{t\to 0_{+}}\int_{\Omega}\overline{J_t(\zeta u)(x)}v(x)\,d^nx
\nonumber\\[2pt]
&=\int_{\Omega}\overline{(\zeta u)(x)}v(x)\,d^nx=\int_{\Omega}\overline{u(x)}v(x)\,d^nx, 
\end{align}
as wanted. In the case when $u,v$ are as indicated in \eqref{ii76g4r-AAA.2}, pick some real-valued function 
$\zeta\in C^\infty_0(\Omega)$ with $\zeta=1$ in a neighborhood of ${\rm supp}\, (u)$.
Observing that $J_t(\zeta v)\in C^\infty_0(\Omega)$ for $t>0$ sufficiently small, 
and $u,\zeta v\in L^1(\Omega)$, then permits us to write
\begin{align}\label{khyg-886-B}
{}_{H^{s}(\Omega)}\big\langle u,v\big\rangle_{H^{-s}(\Omega)}
&={}_{H^{s}(\Omega)}\big\langle\zeta u, v\big\rangle_{H^{-s}(\Omega)}
={}_{H^{s}(\Omega)}\big\langle u,\zeta v\big\rangle_{H^{-s}(\Omega)}
\nonumber\\[2pt]
&=\lim_{t\to 0_{+}}{}_{H^{s}(\Omega)}\big\langle u,J_t(\zeta v)\big\rangle_{H^{-s}(\Omega)}
=\lim_{t\to 0_{+}}{}_{{\mathcal{D}}'(\Omega)}
\big\langle\,\overline{u},J_t(\zeta v)\big\rangle_{{\mathcal{D}}(\Omega)}
\nonumber\\[2pt]
&=\lim_{t\to 0_{+}}\int_{\Omega}\overline{u(x)}J_t(\zeta v)(x)\,d^nx
=\lim_{t\to 0_{+}}\int_{\Omega}\overline{(J_t u)(x)}(\zeta v)(x)\,d^nx
\nonumber\\[2pt]
&=\int_{\Omega}\overline{u(x)}(\zeta v)(x)\,d^nx=\int_{\Omega}\overline{u(x)}v(x)\,d^nx,
\end{align}
once again as desired. 
\end{proof}
%%%%%%%

We continue with a result complementing \eqref{eq:DDEj6g5}. 
To state it, let $\operatorname{Lip}(\Omega)$ stand for the space of Lipschitz functions in $\Omega$.

%%%%%%%
\begin{lemma}\label{utfr}
If $\Omega$ is a bounded Lipschitz domain, then for every $s\in[-1,1]$ it follows that
multiplication with a function from $\operatorname{Lip}(\Omega)$ induces a well defined, 
linear, and bounded operator from $H^s(\Omega)$ into itself. 
\end{lemma}
%%%%%%%
\begin{proof}
The case when $s\in[0,1]$ is seen via interpolation between $s=0$ and $s=1$. 
Furthermore, since pointwise multiplication with a function does not 
increase the support, pointwise multiplication by a Lipschitz function 
also preserves $H^s_0(\Omega)$, for each $s\in[0,1]$. Hence, by duality, this 
also preserves $\big(H^s_0(\Omega)\big)^*=H^{-s}(\Omega)$ for every $s\in[0,1]$. 
\end{proof}
%%%%%%%

We conclude this subsection with a discussion aimed at identifying the amount of smoothness, 
measured on the scales of fractional Sobolev spaces, possessed by certain functions defined in 
bounded Lipschitz domains. Here is our first concrete result in this regard.

%%%%%%%
\begin{lemma}\label{LL.TEXAS}
Fix $\beta\in(\tfrac{1}{2},1)$ and consider the planar open set 
\begin{equation}\label{eq:Texas.1}
\Omega_\beta:=\{z\in{\mathbb{C}}\,|\,0<|z|<1\,\text{ and }\,0<{\rm arg}\,z<\pi/\beta\}.
\end{equation}
Suppose $w\in C^1(\Omega_\beta)$ is a function with the property that there exists some constant 
$C\in(0,\infty)$ such that 
\begin{equation}\label{eq:Texas.2}
|w(x)|\leq C|x|^{\beta-1}\,\text{ and }\,
|(\nabla w)(x)|\leq C|x|^{\beta-2}\,\text{ for each }\,x\in\Omega_\beta.
\end{equation}
Then $w$ belongs to the Sobolev space $H^s(\Omega_\beta)$ whenever $s<\beta$.
\end{lemma}
%%%%%%%
\begin{proof}
First, one observes that the first inequality in \eqref{eq:Texas.2} implies
\begin{align}\label{eq:Texas.0}
\|w\|^2_{L^2(\Omega_\beta)} &\leq C\int_{\{x\in{\mathbb{R}}^2|\,|x|<1\}}|x|^{2\beta-2}\,d^2x
\nonumber\\[2pt]
&=C\int_0^1\rho^{2\beta-1}\,d\rho<\infty.
\end{align}
Next, elementary geometry shows that 
\begin{equation}\label{eq:Texas.4} 
\text{$B(x,r)\cap\Omega_\beta$ is a convex set}\,\text{ for each }\,x\in\Omega_\beta\,\text{ and }\,r\in(0,|x|).
\end{equation}
To proceed, given any $x\in{\mathbb{R}}^2$ and $h\in{\mathbb{R}}^2\backslash\{0\}$, define the first-order difference 
\begin{equation}\label{eq:Texas.5}
(\Delta_h w)(x):=\left\{
\begin{array}{ll}
w(x+h)-w(x) & \text{ if $x\in\Omega_\beta$ and $x+h\in\Omega_\beta$},
\\[2pt]
0 & \text{ if either $x\notin\Omega_\beta$ or $x+h\notin\Omega_\beta$}.
\end{array}
\right.
\end{equation}

Suppose $x\in{\mathbb{R}}^2$ and $h\in{\mathbb{R}}^2\backslash \{0\}$ are such that 
$x\in\Omega_\beta$, $x+h\in\Omega_\beta$, and $|x|>2|h|$. Denote by $(x,x+h)$ the open line 
segment with endpoints $x$ and $x+h$ and pick an arbitrary point $y$ belonging to $(x,x+h)$. 
Then $|x-y|<|h|$ which, in turn, permits us to estimate 
\begin{equation}\label{eq:Texas.6}
|x|\leq|x-y|+|y|<|h|+|y|<2^{-1}|x|+|y|.
\end{equation}
This ultimately implies 
\begin{equation}\label{eq:Texas.7}
2^{-1}|x|<|y|\,\text{ for each }\,y\in(x,x+h). 
\end{equation}
One also observes that since both $x$ and $x+h$ belong to $B(x,|h|)\cap\Omega_\beta$, the property recorded in 
\eqref{eq:Texas.4} ensures that $(x,x+h)\subseteq B(x,|h|)\cap\Omega_\beta\subseteq\Omega_\beta$. Granted these facts, one invokes the Mean Value Theorem which, in view of \eqref{eq:Texas.2} and \eqref{eq:Texas.7}, permits one to estimate 
\begin{align}\label{eq:Texas.8}
|(\Delta_h w)(x)| &=|w(x+h)-w(x)|\leq |h|\sup_{y\in(x,x+h)}|(\nabla w)(y)|
\nonumber\\[2pt]
&\leq C|h|\sup_{y\in(x,x+h)}|y|^{\beta-2}\leq C|h||x|^{\beta-2},
\end{align}
for some constant $C\in(0,\infty)$ which depends only on $w$ and $\beta$. Consequently, 
for each given $h\in{\mathbb{R}}^2\backslash \{0\}$, we may rely on \eqref{eq:Texas.8} to write 
(keeping in mind that $2\beta-3<-1$)
\begin{align}\label{eq:Texas.9}
\int_{\{x\in\Omega_\beta|\,|x|>2|h|\}}|(\Delta_h w)(x)|^2\,d^2x &
\leq C|h|^2\int_{\{x\in{\mathbb{R}}^2|\,|x|>2|h|\}}|x|^{2\beta-4}\,d^2x
\nonumber\\[2pt]
&=C|h|^2\int_{2|h|}^{\infty}\rho^{2\beta-3}\,d\rho=C|h|^{2\beta},
\end{align}
for some constant $C\in(0,\infty)$ independent of $h$.

Next, assume that $x\in{\mathbb{R}}^2$ and $h\in{\mathbb{R}}^2\backslash\{0\}$ are such that 
$x\in\Omega_\beta$, $x+h\in\Omega_\beta$, and $|x|\leq 2|h|$. From \eqref{eq:Texas.2} we know that
\begin{align}\label{eq:Texas.11}
|(\Delta_h w)(x)| &\leq|w(x)|+|w(x+h)|
\nonumber\\[2pt]
&\leq C|x|^{\beta-1}+C|x+h|^{\beta-1}.
\end{align}
As such, 
\begin{align}\label{eq:Texas.12}
\int_{\{x\in\Omega_\beta|\,|x|\leq 2|h|\}}|(\Delta_h w)(x)|^2\,d^2x\leq{\rm I}+{\rm II} 
\end{align}
where, for some constant $C\in(0,\infty)$ independent of $h$, 
\begin{align}\label{eq:Texas.12.WACO.a}
{\rm I}:=C\int_{\{x\in{\mathbb{R}}^2|\,|x|\leq 2|h|\}}|x|^{2\beta-2}\,d^2x
=C\int_0^{2|h|}\rho^{2\beta-1}\,d\rho=C|h|^{2\beta},
\end{align}
and
\begin{align}\label{eq:Texas.12.WACO.b}
{\rm II} &:=C\int_{\{x\in{\mathbb{R}}^2|\,|x|\leq 2|h|\}}|x+h|^{2\beta-2}\,d^2x
\nonumber\\[6pt]
&\,\leq C\int_{\{x\in{\mathbb{R}}^2|\,|x+h|\leq 3|h|\}}|x+h|^{2\beta-2}\,d^2x
\nonumber\\[6pt]
&\,=C\int_{\{y\in{\mathbb{R}}^2|\,|y|\leq 3|h|\}}|y|^{2\beta-2}\,d^2y
\nonumber\\[6pt]
&\,=C\int_0^{3|h|}\rho^{2\beta-1}\,d\rho=C|h|^{2\beta}.
\end{align}
Collectively, the estimates established in \eqref{eq:Texas.9} and \eqref{eq:Texas.12}-\eqref{eq:Texas.12.WACO.b} imply that 
there exists some constant $C\in(0,\infty)$ with the property that
\begin{align}\label{eq:Texas.13}
\int_{\Omega_\beta}|(\Delta_h w)(x)|^2\,d^2x\leq C|h|^{2\beta}\,\text{ for each }\,h\in{\mathbb{R}}^2.
\end{align}
In turn, this allows us to conclude that 
\begin{align}\label{eq:Texas.14}
\sup_{|h|\leq t}\|\Delta_h w\|^2_{L^2(\Omega_\beta)}\leq Ct^{2\beta}\,\text{ for each }\,t\in(0,\infty),
\end{align}
hence, further, 
\begin{align}\label{eq:Texas.15}
\Big(\int_0^1 t^{-2s}\sup_{|h|\leq t}\|\Delta_h w\|^2_{L^2(\Omega_\beta)}\frac{dt}{t}\Big)^{1/2}<\infty
\,\text{ for each }\,s\in(0,\beta).
\end{align}
Since from \cite[Theorem~3.18, p.~30]{Di03} we know that for each $s\in(0,1)$ 
the norm of $w$ in $B^{2,2}_s(\Omega_\beta)=H^s(\Omega_\beta)$ is equivalent to  
\begin{align}\label{eq:Texas.16}
\|w\|_{L^2(\Omega_\beta)}+\Big(\int_0^1 t^{-2s}\sup_{|h|\leq t}\|\Delta_h w\|^2_{L^2(\Omega_\beta)}\frac{dt}{t}\Big)^{1/2},
\end{align}
one finally concludes from \eqref{eq:Texas.0} and \eqref{eq:Texas.15} that $w$ belongs 
to the Sobolev space $H^s(\Omega_\beta)$ whenever $s<\beta$.
\end{proof}
%%%%%%%

In turn, Lemma~\ref{LL.TEXAS} is an ingredient in the proof of the following regularity result
(answering a question which arose in discussions with Volodymyr Derkach). 

%%%%%%%
\begin{proposition}\label{CC.WACO}
For some fixed $\beta\in(\tfrac{1}{2},1)$, consider the planar open set $\Omega_\beta$ as in 
\eqref{eq:Texas.1} and define the function $u_\beta(z):={\rm Im}\,(z^\beta)$ for each $z\in\Omega_\beta$ or, in polar 
coordinates, 
\begin{equation}\label{eq:Texas.WACO.1}
u_\beta(\rho,\theta):=\rho^\beta\sin(\beta\theta)\,\text{ for each }\,z=\rho e^{i\theta}\in\Omega_\beta.
\end{equation}
Then the function $u_\beta$ belongs to the Sobolev space $H^s(\Omega_\beta)$ whenever $s<1+\beta$, 
however, $u_\beta\notin H^{1+\beta}(\Omega_\beta)$. Consequently, for each cutoff function 
$\phi\in C^\infty_0({\mathbb{R}}^2)$ with ${\rm supp}\,(\phi)\subseteq B(0,1/2)$ and 
$\phi=1$ near the origin one has
\begin{equation}\label{eq:MOND.waco.1}
\phi\,u_\beta\in\accentset{\circ}{H}^1(\Omega_\beta)\cap H^s(\Omega_\beta)\,\text{ for each }\,s<1+\beta,
\,\text{ but }\,\phi\,u_\beta\notin H^{1+\beta}(\Omega_\beta).
\end{equation}
\end{proposition}
%%%%%%%
\begin{proof}
First we show that $u_\beta\in H^s(\Omega_\beta)$ whenever $s<1+\beta$. 
In view of the monotonicity property of the fractional Sobolev scale (cf.~\eqref{eq:DDEEn.2}) 
it suffices to consider the case $1<s<1+\beta$. Since clearly $u_\beta\in L^2(\Omega_\beta)$, 
from the lifting result recorded in \eqref{eq:TRfav} (presently used with $s$ replaced by $s-1$) 
we know that $u_\beta$ belongs to the Sobolev space $H^s(\Omega_\beta)$ if and only if
\begin{align}\label{eq:TRfav.WACO}
w_j:=\partial_j u_\beta\,\text{ belongs to }\,H^{s-1}(\Omega_\beta)\,\text{ for each }\,j\in\{1,2\}.
\end{align}
Given that fact, as seen from \eqref{eq:Texas.WACO.1}, for each $j\in\{1,2\}$ one has  
$w_j\in C^1(\Omega_\beta)$ and there exists some constant $C\in(0,\infty)$ such that 
\begin{equation}\label{eq:Texas.2.WACO}
|w_j(x)|\leq C|x|^{\beta-1}\,\text{ and }\,
|(\nabla w_j)(x)|\leq C|x|^{\beta-2}\,\text{ for all }\,x\in\Omega_\beta,
\end{equation}
Lemma~\ref{LL.TEXAS} applies and leads to the conclusion that $w_j\in H^t(\Omega_\beta)$ whenever 
$t<\beta$ for each $j\in\{1,2\}$, which then establishes \eqref{eq:TRfav.WACO}. In turn, this completes 
the proof of the fact that $u_\beta\in H^s(\Omega_\beta)$ whenever $s<1+\beta$, as claimed.

Next, we turn our attention to the task of showing that $u_\beta\notin H^{1+\beta}(\Omega_\beta)$. 
Because of \eqref{eq:TRfav}, this boils down to proving that we cannot have
\begin{align}\label{eq:TRfav.WACO.II}
w_j:=\partial_j u_\beta\in H^{\beta}(\Omega_\beta)\,\text{ for each }\,j\in\{1,2\}.
\end{align}
Arguing by contradiction, we assume that \eqref{eq:TRfav.WACO.II} holds. Then, 
with the first-order difference operator defined as in \eqref{eq:Texas.5}, we may invoke 
\cite[Theorem~3.18, p.~30]{Di03} to conclude that
\begin{align}\label{eq:Texas.16.II}
\sum_{j=1}^2\int_0^1 t^{-2\beta}\sup_{|h|\leq t}\|\Delta_h w_j\|^2_{L^2(\Omega_\beta)}\frac{dt}{t}<\infty.
\end{align}
Since in polar coordinates one has  
\begin{align}\label{eq:Texas.16.III.a}
w_1(\rho,\theta) &=\cos\theta\frac{\partial u(\rho,\theta)}{\partial\rho}
-\frac{1}{\rho}\sin\theta\frac{\partial u(\rho,\theta)}{\partial\theta}
\nonumber\\[2pt]
&=\beta\rho^{\beta-1}\big(\sin(\beta\theta)\cos\theta-\cos(\beta\theta)\sin\theta\big)
\nonumber\\[2pt]
&=\beta\rho^{\beta-1}\sin(\beta\theta-\theta),
\end{align}
and
\begin{align}\label{eq:Texas.16.III.b}
w_2(\rho,\theta) &=\sin\theta\frac{\partial u(\rho,\theta)}{\partial\rho}
+\frac{1}{\rho}\cos\theta\frac{\partial u(\rho,\theta)}{\partial\theta}
\nonumber\\[2pt]
&=\beta\rho^{\beta-1}\big(\sin(\beta\theta)\sin\theta+\cos(\beta\theta)\cos\theta\big)
\nonumber\\[2pt]
&=\beta\rho^{\beta-1}\cos(\beta\theta-\theta),
\end{align}
one concludes that 
\begin{align}\label{eq:Texas.16.III.c}
\sum_{j=1}^2|w_j(x)|^2=\beta^2|x|^{2\beta-2}\,\text{ for each }\,x\in\Omega_\beta. 
\end{align}
In addition, let us observe that if $x\in\Omega_\beta$ has $|x|< 1/2$ then also $2x\in\Omega_\beta$ and one obtains 
$w_j(2x)=2^{\beta-1}w_j(x)$ for each $j\in\{1,2\}$. Consequently, for each $t\in(0,1/2)$, one estimates  
\begin{align}\label{eq:Texas.16.IV}
\sum_{j=1}^2\sup_{|h|\leq t}\|\Delta_h w_j\|^2_{L^2(\Omega_\beta)}
&\geq\sum_{j=1}^2\sup_{|h|\leq t}\int_{\{x\in\Omega_\beta|\,|x|<t\}}|(\Delta_h w_j)(x)|^2\,d^2x
\nonumber\\[2pt]
&\geq\sum_{j=1}^2\int_{\{x\in\Omega_\beta|\,|x|<t\}}|w_j(2x)-w_j(x)|^2\,d^2x
\nonumber\\[2pt]
&=(1-2^{\beta-1})^2\int_{\{x\in\Omega_\beta|\,|x|<t\}}\sum_{j=1}^2|w_j(x)|^2\,d^2x
\nonumber\\[2pt]
&=\beta^2(1-2^{\beta-1})^2\int_{\{x\in\Omega_\beta|\,|x|<t\}}|x|^{2\beta-2}\,d^2x
\nonumber\\[2pt]
&=\pi\beta(1-2^{\beta-1})^2\int_0^t\rho^{2\beta-1}\,d\rho
\nonumber\\[2pt]
&=\frac{\pi}{2}(1-2^{\beta-1})^2 t^{2\beta}.
\end{align}
However, this implies
\begin{align}\label{eq:Texas.16.V}
\sum_{j=1}^2\int_0^1 t^{-2\beta}\sup_{|h|\leq t}\|\Delta_h w_j\|^2_{L^2(\Omega_\beta)}\frac{dt}{t}
&\geq\frac{\pi}{2}(1-2^{\beta-1})^2\int_0^{1/2}t^{-2\beta}t^{2\beta}\frac{dt}{t}
\nonumber\\[2pt]
&=\frac{\pi}{2}(1-2^{\beta-1})^2\int_0^{1/2}\frac{dt}{t}=\infty,
\end{align}
contradicting \eqref{eq:Texas.16.II}. In turn, this contradiction shows that 
$u_\beta\notin H^{1+\beta}(\Omega_\beta)$.

It remains to justify the claims made in \eqref{eq:MOND.waco.1}. To this end, fix a cutoff function 
$\phi\in C^\infty_0({\mathbb{R}}^2)$ with ${\rm supp}\,(\phi)\subseteq B(0,1/2)$ and $\phi=1$ near the origin.  
From \eqref{eq:DDEj6g5} and what we have proved already, we see that $\phi\,u_\beta\in H^s(\Omega_\beta)$ 
for each $s<1+\beta$. One notes that $u_\beta$ in \eqref{eq:Texas.WACO.1} is designed so that it extends 
continuously to the closure of $\Omega_\beta$ and this extension vanishes on 
$\{z\in\partial\Omega_\beta\,|\,{\rm arg}\,z\in\{0,\pi/\beta\}\}$. Granted these properties, 
\eqref{eq:VCca.1-tt5} and Lemma~\ref{Hgg-uh} then imply that 
$\phi\,u_\beta\in\accentset{\circ}{H}^1(\Omega_\beta)$. Next, observe that $(1-\phi)u_\beta$ is of class 
$C^\infty$ and has bounded derivatives of any order in $\Omega_\beta$. As such, $(1-\phi)u_\beta$ belongs
to any Sobolev space in $\Omega_\beta$. In particular, $(1-\phi)u_\beta\in H^{1+\beta}(\Omega_\beta)$.
Since we already know that $u_\beta\notin H^{1+\beta}(\Omega_\beta)$, this ultimately implies that 
actually $\phi\,u_\beta$ does not belong to $H^{1+\beta}(\Omega_\beta)$.
\end{proof}

%%%%%%%%%%%%%%%%%%%%%%%%%%%%%%%%%%%
\subsection{Fractional Sobolev spaces on the boundaries of Lipschitz domains}\label{ss2.4}
%%%%%%%%%%%%%%%%%%%%%%%%%%%%%%%%%%%

In a first stage, assume that $\Omega\subset\bbR^n$ is the domain lying above
the graph of a Lipschitz function $\varphi\colon\bbR^{n-1}\to\bbR$, and let $0\leq s\leq 1$. 
Then the Sobolev space $H^s(\partial\Omega)$ consists of functions 
$f\in L^2(\partial\Omega)$ such that $f(x',\varphi(x'))$,
as a function of $x'\in\bbR^{n-1}$, belongs to $H^s(\bbR^{n-1})$. 
To define $H^{-s}(\partial\Omega)$, let ${\rm Lip}_{\rm comp}(\partial\Omega)$ be 
the space of compactly supported Lipschitz functions on $\partial\Omega$ 
(equipped with the usual inductive limit topology). Then a functional 
$f\in\big({\rm Lip}_{\rm comp}(\partial\Omega)\big)^*$ is said to belong to $H^{-s}(\partial\Omega)$ 
provided $\sqrt{1+|(\nabla'\varphi)(\dott)|^2}f(\dott,\varphi(\dott))\in H^{-s}(\bbR^{n-1})$. 
Here, $\sqrt{1+|(\nabla'\varphi)(\dott)|^2}f(\dott,\varphi(\dott))$ is understood
as the distribution in ${\mathbb{R}}^{n-1}$ acting according to 
\begin{align}\label{i6f444}
& C^\infty_0({\mathbb{R}}^{n-1})\ni\psi\mapsto{}_{{\rm Lip}_{\rm comp}(\partial\Omega)}
\big\langle\widetilde\psi,f\big\rangle_{({\rm Lip}_{\rm comp}(\partial\Omega))^*}\,\text{  where},
\nonumber\\[2pt]
& \quad\text{given any $\psi\in C^\infty_0({\mathbb{R}}^{n-1})$, the function 
$\widetilde\psi\in{\rm Lip}_{\rm comp}(\partial\Omega)$} 
\\[2pt]
& \quad\text{is given by $\widetilde\psi(x):=\psi(x')$ for each $x=(x',\varphi(x'))\in\partial\Omega$}.  
\nonumber
\end{align}

Next, to define $H^s(\dOm)$ for $-1\leq s\leq 1$, when $\Om$ is a Lipschitz domain with compact 
boundary, we use a smooth partition of unity to reduce matters to the graph case just discussed. 
More precisely, if $0\leq s\leq 1$ then $f\in H^s(\partial\Omega)$ if and only if 
$f\in L^2(\partial\Omega)$ and the assignment ${\mathbb{R}}^{n-1}\ni x'\mapsto(\zeta f)(x',\varphi(x'))$ 
is in $H^s({\mathbb{R}}^{n-1})$ whenever $\zeta\in C^\infty_0({\mathbb{R}}^n)$ and 
$\varphi\colon {\mathbb{R}}^{n-1}\to{\mathbb{R}}$ is a Lipschitz function with the property 
that if $\Sigma$ is an appropriate rotation and translation of the graph 
$\{(x',\varphi(x'))\in\bbR^n\,|\,x'\in{\mathbb{R}}^{n-1}\}$, then 
$(\supp\,(\zeta)\cap\partial\Omega)\subset\Sigma$. Then Sobolev spaces 
with a negative amount of smoothness are defined in an analogous fashion. 

From the above characterization of $H^s(\partial\Omega)$ it follows that 
properties of Sobolev spaces $H^s(\bbR^{n-1})$ with $s\in[-1,1]$ which are invariant 
under multiplication by smooth, compactly supported functions, as well as composition 
by bi-Lipschitz maps, readily extend to the setting of $H^s(\partial\Omega)$ 
(via localization and pullback). In particular, 
\begin{equation}\label{y6rf8HB}
\big(H^s(\partial\Omega)\big)^*=H^{-s}(\partial\Omega),\,\text{ whenever }\,-1\leq s\leq 1,
\end{equation}
and one has a continuous (in fact compact) and dense embedding 
\begin{equation}\label{y6rf8HB.ttr}
H^{s_2}(\partial\Omega)\hookrightarrow H^{s_1}(\partial\Omega),\,\text{ whenever }\,-1\leq s_1<s_2\leq 1.
\end{equation}
In addition, if $\Omega$ is a bounded Lipschitz domain in ${\mathbb{R}}^n$ then
\begin{equation}\label{lip-d1}
C^\infty({\mathbb{R}}^n)\big|_{\partial\Omega}
\hookrightarrow H^s(\partial\Omega)\,\text{ densely, }\,\forall\,s\in[-1,1].
\end{equation}
See, for instance, \cite{JK95}, \cite[Chapter~3]{Mc00}, \cite{MM13}.

Later on, we shall employ the following characterization of the Sobolev space 
of order one on the boundary of a Lipschitz domain; 
see \cite[Propositions~2.8--2.9, p.~33]{MM13} for a proof.

%%%%%%%
\begin{lemma}\label{Dida}
Let $\Omega\subset{\mathbb{R}}^n$ be a bounded Lipschitz domain, with outward unit normal 
$\nu=(\nu_1,\dots,\nu_n)$ and surface measure $\sigma$. Then $H^1(\partial\Omega)$ is the 
collection of functions $\varphi\in L^2(\partial\Omega)$ with the property that there exists 
a constant $C\in(0,\infty)$ such that
\begin{equation}\label{eq:Mvc}
\sum_{j,k=1}^n\Big|\int_{\partial\Omega}\varphi\,\partial_{\tau_{jk}}\psi\,d^{n-1}\sigma\Big|
\leq C\big\|\psi\big|_{\partial\Omega}\big\|_{L^2(\partial\Omega)},\quad
\forall\,\psi\in C^\infty_0({\mathbb{R}}^n),
\end{equation}
where 
\begin{equation}\label{eq:Mvcuyg}
\partial_{\tau_{jk}}\psi:=\nu_j\big(\partial_k\psi\big)\big|_{\partial\Omega}
-\nu_k\big(\partial_j\psi\big)\big|_{\partial\Omega}\,\text{ for each }\,j,k\in\{1,\dots,n\}.
\end{equation}
In addition, 
\begin{align}\label{trRFj6t}
\|\varphi\|_{H^1(\partial\Omega)}\approx
\|\varphi\|_{L^2(\partial\Omega)}
+\sup_{\substack{\psi\in C^\infty_0({\mathbb{R}}^n)\\ \|\psi|_{\partial\Omega}\|_{L^2(\partial\Omega)}\leq 1}}
\left\{\sum_{j,k=1}^n\Big|\int_{\partial\Omega}\varphi\,\partial_{\tau_{jk}}\psi\,d^{n-1}\sigma\Big|\right\},
\end{align}
uniformly for $\varphi\in H^1(\partial\Omega)$. 
\end{lemma}
%%%%%%%

In closing, we note that if $\Omega\subset\bbR^n$ is a bounded
Lipschitz domain then for any $j,k\in\{1,\dots,n\}$ the first-order tangential differential 
operator $\partial_{\tau_{jk}}$ extends to a well defined, linear and bounded mapping in the context 
\begin{equation}\label{Pf-2}
\partial_{\tau_{jk}}:H^s(\partial\Omega)\rightarrow H^{s-1}(\partial\Omega),
\quad 0\leq s\leq 1, 
\end{equation}
This is proved by interpolating the case $s=1$ and its dual version. In fact, the 
following more general result (extending Lemma~\ref{Dida}) is true. Specifically, 
assuming that $\Omega\subset\bbR^n$ is a bounded Lipschitz domain, for every 
$s\in[0,1]$ one has  
\begin{equation}\label{Pf-2.3}
H^s(\partial\Omega)=\big\{f\in L^2(\partial\Omega)\,\big|\,
\partial_{\tau_{jk}}f\in H^{s-1}(\partial\Omega),\,1\leq j,k\leq n\big\}
\end{equation}
and 
\begin{equation}\label{Pf-2.4}
\|f\|_{H^s(\partial\Omega)}\approx\|f\|_{L^2(\partial\Omega)}
+\sum_{j,k=1}^n\|\partial_{\tau_{jk}}f\|_{H^{s-1}(\partial\Omega)},
\end{equation}
uniformly for $f\in H^s(\partial\Omega)$ (see the discussion in \cite{GM08}). 

%%%%%%%%%%%%%%%%%%%%%%%%%%%%
\subsection{Sobolev regularity in terms of the nontangential maximal function}\label{ss2.5}
%%%%%%%%%%%%%%%%%%%%%%%%%%%%

We begin by recalling a standard elliptic regularity result to the effect that 
\begin{align}\label{eq:MM3.ELL.WACO}
\left.
\begin{array}{r}
\text{$\Omega\subseteq{\mathbb{R}}^n$ open, $V\in L^p_{\rm loc}(\Omega)$ with $p>n/2$} 
\\[2pt]
\text{$u\in L^2_{\rm loc}(\Omega)$ with $(-\Delta+V)u=0$ in $\Omega$}
\end{array}
\right\}\Longrightarrow u\in C^1(\Omega).
\end{align}
See, e.g., \cite{MMT01}, \cite[Proposition~3.1]{MT00}, \cite{Ta96}, \cite{Ta11} in this regard. 

In the class of functions that are null-solutions of zeroth-order perturbation 
of the Laplacian in a bounded Lipschitz domain $\Omega\subset{\mathbb{R}}^n$, the relationship 
between membership to Sobolev spaces in $\Omega$, on the one hand, and the membership of the 
nontangential maximal function to Lebesgue spaces on the boundary $\partial\Omega$, 
on the other hand, becomes rather precise. First, one has the following characterization:
\begin{align}\label{eq:MM3}
\begin{split} 
& \text{if $0\leq V\in L^p(\Omega)$ with $p>n$ and $u\in C^1(\Omega)$ with $(-\Delta+V)u=0$ in $\Omega$,}
\\[4pt]
& \quad\text{then }\,{\mathcal{N}}_\kappa u\in L^2(\partial\Omega)\,\text{ if and only if }\,u\in H^{1/2}(\Omega),
\end{split}
\end{align}
with naturally accompanying estimates, namely, 
\begin{align}\label{eq:MM3-AC}
\|{\mathcal{N}}_\kappa u\|_{L^2(\partial\Omega)}\approx\|u\|_{H^{1/2}(\Omega)},
\end{align}
uniformly for $u$ as in \eqref{eq:MM3}.
See \cite{Fa88}, \cite{JK95}, \cite{MMP97}, for $V=0$, and \cite{MT00} for the general case.
From \eqref{eq:MM3}, \eqref{eq:TRfav}, and iterations, one deduces that 
\begin{align}\label{eq:MM3BIS}
& \text{if $V\in[0,\infty)$ is constant, $k\in{\mathbb{N}}_0$, and $u\in C^\infty(\Omega)$ with $(-\Delta+V)u=0$ in $\Omega$, then}    
\nonumber\\ 
& \quad{\mathcal{N}}_\kappa(\partial^\alpha u)\in L^2(\partial\Omega)
\,\text{ for }\,|\alpha|\leq k\,\text{ if and only if }\,u\in H^{k+(1/2)}(\Omega), 
\end{align}
with the naturally accompanying estimates
\begin{align}\label{eq:MM3BIS-AC}
\sum_{|\alpha|\leq k}\|{\mathcal{N}}_\kappa(\partial^\alpha u)\|_{L^2(\partial\Omega)}
\approx\|u\|_{H^{k+(1/2)}(\Omega)},
\end{align}
uniformly for $u$ as in \eqref{eq:MM3BIS}. In this regard, we also record the following Fatou-type result 
(cf. \cite[Proposition~3.1]{MT01}, \cite[Proposition~4.7, Proposition~5.6]{MT00a}):
\begin{align}\label{eq:MM4} 
& \text{if $u\in C^1(\Omega)$ with $(-\Delta+V)u=0$ in $\Omega$ for some $0\leq V\in L^\infty(\Omega)$, then}
\nonumber\\[2pt]
& \quad{\mathcal{N}}_\kappa u\in L^2(\partial\Omega)\,\text{ implies }\,u\big|^{\kappa-{\rm n.t.}}_{\partial\Omega}
\,\text{ exists $\sigma$-a.e.~and $u\big|^{\kappa-{\rm n.t.}}_{\partial\Omega}\in L^2(\partial\Omega)$},      
\\[2pt]
& \quad\text{while }\,{\mathcal{N}}_\kappa(\nabla u)\in L^2(\partial\Omega)\,\text{ implies }\, 
(\nabla u)\big|^{\kappa-{\rm n.t.}}_{\partial\Omega}\,\text{ exists $\sigma$-a.e.~and in 
$\big[L^2(\partial\Omega)\big]^n$}.  
\nonumber
\end{align}

%%%%%%%%%%%%%%%%%%%%%%%%%%%%%%%%%%%%%%%%%%
%%%%%%%%%%%%%%%%%%%%%%%%%%%%%%%%%%%%%%%%%%
\section{A Sharp Dirichlet Trace Involving Sobolev and Besov Spaces}\label{s3}
%%%%%%%%%%%%%%%%%%%%%%%%%%%%%%%%%%%%%%%%%%
%%%%%%%%%%%%%%%%%%%%%%%%%%%%%%%%%%%%%%%%%%

The prime object in this section will be a detailed treatment of the Dirichlet trace operator 
$\gamma_D\colon H^s(\Omega)\rightarrow H^{s-(1/2)}(\partial\Omega)$ associated with bounded 
Lipschitz domains $\Omega\subset\bbR^n$, at first studied for all $s\in(1/2, 3/2)$. Upon 
noticing the difficulties extending the Dirichlet trace to the endpoints $s=1/2$ and $s=3/2$, 
we employ additional regularity of the Laplacian in Section~\ref{ss3.2} and \ref{ss3.3} to 
arrive at sharp Dirichlet trace results in the Sobolev and Besov space context for the full 
scale $s\in[1/2,3/2]$.

%%%%%%%%%%%%%%%%
\subsection{A first look at the Dirichlet trace}\label{ss3.1}
%%%%%%%%%%%%%%%%

Let $\Omega$ be a bounded Lipschitz domain in ${\mathbb{R}}^n$. In this context, 
the Dirichlet boundary trace map $f\mapsto f\big|_{\partial\Omega}$, originally considered for 
functions $f\in C^\infty(\overline\Omega)$, extends to operators (compatible with one another) 
\begin{equation}\label{eqn:gammaDs.1}
\gamma_D:H^s(\Omega)\rightarrow H^{s-(1/2)}(\partial\Omega),\quad\forall\, 
s\in\big(\tfrac{1}{2},\tfrac{3}{2}\big) 
\end{equation}
(see also \cite[Lemma~3.6]{Co88}), that are linear, continuous, surjective, and whose operator norm depend on the underlying Lipschitz domain 
only via the Lipschitz character of the latter. (We agree that for vector-valued functions the Dirichlet 
trace is applied componentwise.) In fact, there exist linear and bounded operators 
\begin{align}\label{eqn:gammaDs.1-INV1}
\vartheta_D:H^{s-(1/2)}(\partial\Omega)\rightarrow H^s(\Omega),\quad\forall\, 
s\in\big(\tfrac{1}{2},\tfrac{3}{2}\big),
\end{align}
which are right-inverses for those in \eqref{eqn:gammaDs.1}, that is,  
\begin{align}\label{eqn:gammaDs.1-INV2}
\gamma_D(\vartheta_D f)=f,\quad\forall\,f\in H^{s-(1/2)}(\partial\Omega),\quad
\forall\,s\in\big(\tfrac{1}{2},\tfrac{3}{2}\big).
\end{align}
As a consequence, 
\begin{align}\label{eqn:gammaDs.1-INV3}
& \text{given any $s\in\big(\tfrac{1}{2},\tfrac{3}{2}\big)$ 
there exists a constant $C_s\in(0,\infty)$ with the} 
\nonumber\\[2pt]
& \quad\text{property that for every $f\in H^{s-(1/2)}(\partial\Omega)$ there exists $u\in H^s(\Omega)$} 
\\[2pt] 
& \quad\text{satisfying $\gamma_Du=f$ on $\partial\Omega$ and 
$\|u\|_{H^s(\Omega)}\leq C_s\|f\|_{H^{s-(1/2)}(\partial\Omega)}$.}   
\nonumber 
\end{align}
Moreover, 
\begin{align}\label{eqn:gammaDs.1-TR}
\gamma_D(\Phi u)=\big(\Phi\big|_{\partial\Omega}\big)\gamma_D u,\quad
\forall\,u\in H^s(\Omega)\,\text{ with }\,s\in\big(\tfrac{1}{2},\tfrac{3}{2}\big),
\quad\forall\,\Phi\in C^\infty(\overline\Omega).
\end{align}
While the Dirichlet trace operator fails to be bounded in the context of 
\eqref{eqn:gammaDs.1} in the limiting case $s=1/2$, one still obtains that  
\begin{equation}\label{eq:VCca.1-tt5}
\gamma_D:H^{(1/2)+\varepsilon}(\Omega)\rightarrow L^2(\partial\Omega)
\,\text{ is well defined, linear, and bounded},
\end{equation}
for every $\varepsilon>0$. For future reference we also note that for any bounded Lipschitz 
domain $\Omega\subset{\mathbb{R}}^n$ one has (see \eqref{equal.HH} for the first equality) 
\begin{equation}\label{incl-Yb.EE}
H^s_z(\Omega)=\accentset{\circ}{H}^{s}(\Omega)=\big\{u\in H^s(\Omega)\,\big|\,\gamma_Du=0\big\},\quad 
\forall\,s\in\big(\tfrac{1}{2},\tfrac{3}{2}\big).
\end{equation}
See \cite{JK95}, \cite{MM04}, \cite{MM13} for general results of this type; 
cf. also the discussion in \cite{GM08}.

It turns out that the Dirichlet trace operator $\gamma_D$ from \eqref{eqn:gammaDs.1} and the 
pointwise nontangential boundary trace from \eqref{eq:MM2} are compatible, in the sense 
that they agree a.e., whenever they both exist: 

%%%%%%%%%%
\begin{lemma}\label{Hgg-uh}
Let $\Omega\subset\bbR^n$ be a bounded Lipschitz domain, and fix some aperture parameter $\kappa>0$. Then
\begin{align}\label{eq:MM3bis}
\begin{split}
& \text{whenever $u\in H^s(\Omega)$ for some $s\in\big(\tfrac{1}{2},\tfrac{3}{2}\big)$ 
and $u\big|^{\kappa-{\rm n.t.}}_{\partial\Omega}$ exists}    
\\[2pt]
& \quad\text{at $\sigma$-a.e.~point on $\partial\Omega$, then 
$u\big|^{\kappa-{\rm n.t.}}_{\partial\Omega}=\gamma_D u\in H^{s-(1/2)}(\partial\Omega)$}. 
\end{split}
\end{align}
\end{lemma}
%%%%%%%%%%
\begin{proof}
From \cite[Theorem~8.7(iii)]{BMMM13} (cf. also \cite[Corollary~5.7]{BMMM13})
one knows that if $u\in H^s(\Omega)$ for some $s\in\big(\tfrac{1}{2},\tfrac{3}{2}\big)$
then its trace $\gamma_D u\in H^{s-(1/2)}(\partial\Omega)$ has the property that 
\begin{equation}\label{eq:246t4}
(\gamma_D u)(x)=\lim_{r\to 0_+}\bigg\{-\hskip -0.15in
\int_{\Gamma_\kappa(x)\cap B(x,r)}u(y)\,dy\bigg\} 
\,\text{ at $\sigma$-a.e.~$x\in\partial\Omega$},
\end{equation}
where the barred integral, $-\hskip -0.12in\int$\,, indicates the mean average. 
Finally, whenever $\big(u\big|^{\kappa-{\rm n.t.}}_{\partial\Omega}\big)(x)$ exists at 
some point $x\in\partial\Omega$ it is given by the limit in the right-hand side of 
\eqref{eq:246t4}, hence the desired conclusion follows.
\end{proof}
%%%%%%%%%%%

The end-point $s=\tfrac{1}{2}$ is naturally excluded in \eqref{eqn:gammaDs.1} 
since it turns out that $ C^\infty_0(\Omega)$ is dense in $H^{1/2}(\Omega)$ 
(cf. the discussion at the bottom of p.~180 in \cite{JK95}). The Dirichlet trace 
operator \eqref{eqn:gammaDs.1} also fails to be well defined corresponding to the 
end-point case $s=\tfrac{3}{2}$ although, of course, \eqref{eqn:gammaDs.1} 
implies that for each $\varepsilon\in(0,1)$ the operator
\begin{equation}\label{eq:VCHBb}
\gamma_D:H^{3/2}(\Omega)\rightarrow H^{1-\varepsilon}(\partial\Omega)
\,\text{ is well defined, linear, and bounded},
\end{equation}
(though, \eqref{eq:VCHBb} does not hold with $\varepsilon=0$). Indeed, in 
\cite[Proposition~3.2, p.~176]{JK95}
an example of a bounded $C^1$-domain (hence, also Lipschitz) in ${\mathbb{R}}^2$ 
and of a function $u\in H^{3/2}(\Omega)$ are given with the property 
that $\gamma_D u\notin H^1(\partial\Omega)$. 
Hence, what goes wrong when $s=\tfrac{3}{2}$ is that in the class of bounded $C^1$ and 
Lipschitz domains $\Omega$, the Dirichlet boundary trace operator $\gamma_D$, when applied to 
$H^{3/2}(\Omega)$, has a larger range than the usual range $H^1(\partial\Omega)$. 
Nonetheless, the Dirichlet traces of smoother functions in $\Omega$ do belong to 
$H^1(\partial\Omega)$ as our next result shows.

%%%%%%%%%%
\begin{lemma}\label{tafc-644}
Let $\Omega\subset{\mathbb{R}}^n$ be a bounded Lipschitz domain. 
Then, for each $\varepsilon>0$, the Dirichlet trace operator 
\begin{equation}\label{eq:VCca.1}
\gamma_D:H^{(3/2)+\varepsilon}(\Omega)\rightarrow H^1(\partial\Omega)
\,\text{ is well defined, linear, and bounded}.
\end{equation}
\end{lemma}
%%%%%%%%%%
\begin{proof}
In the justification of \eqref{eq:VCca.1} we shall employ the characterization 
of $H^1(\partial\Omega)$ from Lemma~\ref{Dida}. Regarding the tangential derivatives 
$\partial_{\tau_{jk}}$ defined in \eqref{eq:Mvcuyg} one notes that  
for every function $\Phi\in C^\infty(\overline{\Omega})$, the divergence theorem 
(see the last part of Theorem~\ref{banff-3}) yields 
\begin{align}\label{u5fre4}
\int_{\partial\Omega}\partial_{\tau_{jk}}\Phi\,d^{n-1}\sigma
&=\int_{\partial\Omega}\big\{\nu_j(\partial_k\Phi)\big|_{\partial\Omega}
-\nu_k(\partial_j\Phi)\big|_{\partial\Omega}\big\}\,d^{n-1}\sigma
\nonumber\\[2pt]
&=\int_{\Omega}\big\{\partial_j\partial_k\Phi-\partial_k\partial_j\Phi\big\}\,d^n x
\nonumber\\[2pt]
&=0.
\end{align}
Suppose now that some $\zeta,\xi\in C^\infty(\overline{\Omega})$ have been given, and
use the product rule to expand 
\begin{align}\label{u5fre4.2}
\partial_{\tau_{jk}}(\zeta\xi)=\big(\xi\big|_{\partial\Omega}\big)\partial_{\tau_{jk}}\zeta
+\big(\zeta\big|_{\partial\Omega}\big)\partial_{\tau_{jk}}\xi.
\end{align}
Combining \eqref{u5fre4} (written for $\Phi:=\zeta\xi$) and \eqref{u5fre4.2} one therefore arrives at the identity 
\begin{equation}\label{eq:ytrd335}
\int_{\partial\Omega}\big(\xi\big|_{\partial\Omega}\big)\partial_{\tau_{jk}}\zeta\,d^{n-1}\sigma
=-\int_{\partial\Omega}\big(\zeta\big|_{\partial\Omega}\big)\partial_{\tau_{jk}}\xi\,d^{n-1}\sigma, 
\quad\forall\,\zeta,\xi\in C^\infty(\overline{\Omega}).
\end{equation}

Considering an arbitrary function $\eta\in H^{(3/2)+\varepsilon}(\Omega)$, 
then there exists a sequence $\{\eta_m\}_{m\in{\mathbb{N}}}\subset C^\infty(\overline{\Omega})$ 
such that $\eta_m\to\eta$ in $H^{(3/2)+\varepsilon}(\Omega)$ as $m\to\infty$.
In particular, $\nabla\eta_m\to\nabla\eta$ in $\big[H^{(1/2)+\varepsilon}(\Omega)\big]^n$ 
as $m\to\infty$ which, together with the continuity of \eqref{eqn:gammaDs.1}, 
further implies $\nabla\eta_m\big|_{\partial\Omega}\to\gamma_D(\nabla\eta)$ 
in $\big[H^{\varepsilon}(\partial\Omega)\big]^n$, hence also in 
$\big[L^2(\partial\Omega)\big]^n$, as $m\to\infty$.
We also note that $\eta_m|_{\partial\Omega}\rightarrow\gamma_D\eta$ in 
$L^2(\partial\Omega)$ as $m\to\infty$. Based on these facts and the identity 
in \eqref{eq:ytrd335}, given any $\psi\in C^\infty_0({\mathbb{R}}^n)$, 
for each $m\in{\mathbb{N}}$ and $j,k\in\{1,\dots,n\}$, one estimates  
\begin{align}\label{eq:MVcAA}
& \Big|\int_{\partial\Omega}(\gamma_D\eta)\,\partial_{\tau_{jk}}\psi\,d^{n-1}\sigma\Big|
=\lim_{m\to\infty}\Big|\int_{\partial\Omega}\big(\eta_m\big|_{\partial\Omega}\big)
\partial_{\tau_{jk}}\psi\,d^{n-1}\sigma\Big|
\nonumber\\[2pt]
& \quad=\lim_{m\to\infty}\Big|\int_{\partial\Omega}\Big[\nu_k\big(\partial_j\eta_m\big)\big|_{\partial\Omega}
-\nu_j\big(\partial_k\eta_m\big)\big|_{\partial\Omega}\Big]\,(\psi\big|_{\partial\Omega})\,
d^{n-1}\sigma\Big|
\nonumber\\[2pt] 
& \quad\leq C\lim_{m\to\infty}\big\|(\nabla\eta_m)\big|_{\partial\Omega}\big\|_{[L^2(\partial\Omega)]^n}
\big\|\psi\big|_{\partial\Omega}\big\|_{L^2(\partial\Omega)}
\nonumber\\[2pt]
& \quad=C\|\gamma_D(\nabla\eta)\|_{[L^2(\partial\Omega)]^n}\big\|\psi\big|_{\partial\Omega}\big\|_{L^2(\partial\Omega)}
\nonumber\\[2pt]
& \quad\leq C\|\gamma_D(\nabla\eta)\|_{[H^{\delta}(\partial\Omega)]^n}\big\|\psi\big|_{\partial\Omega}\big\|_{L^2(\partial\Omega)}
\nonumber\\[2pt]
& \quad\leq C\|\nabla\eta\|_{[H^{(1/2)+\delta}(\Omega)]^n}\big\|\psi\big|_{\partial\Omega}\big\|_{L^2(\partial\Omega)}
\nonumber\\[2pt]
& \quad\leq C\|\eta\|_{H^{(3/2)+\delta}(\Omega)}\big\|\psi\big|_{\partial\Omega}\big\|_{L^2(\partial\Omega)}
\nonumber\\[2pt]
& \quad\leq C\|\eta\|_{H^{(3/2)+\varepsilon}(\Omega)}\big\|\psi\big|_{\partial\Omega}\big\|_{L^2(\partial\Omega)},
\end{align}
where $0<\delta<\min\{\varepsilon,1\}$. 
In light of Lemma~\ref{Dida}, this proves that $\gamma_D\eta\in H^1(\partial\Omega)$.
Moreover, \eqref{trRFj6t} and \eqref{eq:VCca.1-tt5} imply that there exists a constant 
$C\in(0,\infty)$, independent of $\eta$, with the property that
\begin{align}\label{iu5fed}
\|\gamma_D\eta\|_{H^1(\partial\Omega)} &\leq C\big(\|\gamma_D\eta\|_{L^2(\partial\Omega)} 
+\|\eta\|_{H^{(3/2)+\varepsilon}(\Omega)}\big)
\nonumber\\[2pt]
&\leq C\big(\|\eta\|_{H^{(1/2)+\varepsilon}(\Omega)}+\|\eta\|_{H^{(3/2)+\varepsilon}(\Omega)}\big)
\nonumber\\[2pt]
&\leq C\|\eta\|_{H^{(3/2)+\varepsilon}(\Omega)}.
\end{align}
The proof of \eqref{eq:VCca.1} is therefore complete.
\end{proof}
%%%%%%%%%%

A useful consequence of \eqref{eqn:gammaDs.1} and Lemma~\ref{tafc-644} is recorded below. 

%%%%%%%%%%
\begin{corollary}\label{tafc-644-TRE}
Assume that $\Omega\subset{\mathbb{R}}^n$ is a bounded Lipschitz domain. 
Then, for each $s\in\big[\tfrac12,\tfrac32\big]$ and $\varepsilon>0$ with 
$\varepsilon\not=\tfrac32-s$, the Dirichlet trace operator 
\begin{align}\label{eq:VCca.1htrr}
\begin{split}
\gamma_D:H^{s+\varepsilon}(\Omega)\rightarrow H^{\min\{1,s+\varepsilon-(1/2)\}}(\partial\Omega)
\\[2pt]
\text{is well defined, linear, and bounded}.
\end{split}
\end{align}
\end{corollary}
%%%%%%%%%%

The following technical lemma is going to play a role in the proof of the version of 
the divergence theorem discussed later, in Theorem~\ref{Ygav-75}.

%%%%%%%
\begin{lemma}\label{Tgav9jy}
Assume $\Omega\subset\bbR^n$ is a bounded Lipschitz domain, and suppose that 
$\Omega_\ell\nearrow\Omega$ as $\ell\to\infty$, in the sense described in Lemma~\ref{OM-OM}.
For each $\ell\in{\mathbb{N}}$, denote by $\gamma_{\ell,D}$ the Dirichlet boundary 
trace operator \eqref{eqn:gammaDs.1} associated with the domain $\Omega_\ell$. 
Then for any $u\in\accentset{\circ}{H}^s(\Omega)$, with $s\in\big(\tfrac{1}{2},1\big)$, 
it follows that $u\big|_{\Omega_\ell}\in H^s(\Omega_\ell)$ for each $\ell\in{\mathbb{N}}$ 
and 
\begin{equation}\label{eq:Rfbbba}
\lim\limits_{\ell\to\infty}
\big\|\gamma_{\ell,D}\big(u\big|_{\Omega_\ell}\big)\big\|_{H^{s-(1/2)}(\partial\Omega_\ell)}=0.
\end{equation}
\end{lemma}
%%%%%%%
\begin{proof}
Fix some function $u\in\accentset{\circ}{H}^s(\Omega)$, with $s\in\big(\tfrac{1}{2},1\big)$. 
That $u\big|_{\Omega_\ell}\in H^s(\Omega_\ell)$ for each $\ell\in{\mathbb{N}}$ follows 
directly from definitions. Next, fix an arbitrary function $v\in C^\infty_0(\Omega)$. 
Making use of the fact that dependence on the underlying Lipschitz domain of the operator 
norm of the Dirichlet boundary trace operator manifests itself only via its Lipschitz character one obtains
\begin{align}\label{eq:MNVC.1}
\limsup\limits_{\ell\to\infty}
\big\|\gamma_{\ell,D}\big(u\big|_{\Omega_\ell}\big)\big\|_{H^{s-(1/2)}(\partial\Omega_\ell)}
&=\limsup\limits_{\ell\to\infty}
\big\|\gamma_{\ell,D}\big((u-v)\big|_{\Omega_\ell}\big)\big\|_{H^{s-(1/2)}(\partial\Omega_\ell)}
\nonumber\\[2pt]
&\leq C\limsup\limits_{\ell\to\infty}
\big\|(u-v)\big|_{\Omega_\ell}\big\|_{H^{s}(\Omega_\ell)}
\nonumber\\[2pt]
&\leq C\|u-v\|_{H^s(\Omega)},
\end{align}
where the last inequality is a consequence of \eqref{FSA-2X}.
With \eqref{eq:MNVC.1} in hand, the desired conclusion follows from \eqref{Rdac}.
\end{proof}
%%%%%%%

Admitting the full scale of Besov spaces instead of Sobolev spaces  permits the 
consideration of the Dirichlet boundary trace operator in a more general context than before. 
Specifically, one has the following result which, in contrast to $\gamma_D$ in \eqref{eqn:gammaDs.1}, 
allows including the end-points of the interval $\big(\tfrac{1}{2},\tfrac{3}{2}\big)$. 

%%%%%%%
\begin{proposition}\label{P-New-Trace}
Let $\Omega\subset{\mathbb{R}}^n$ be a bounded Lipschitz domain, and fix an aperture parameter $\kappa>0$. Then the mapping
\begin{equation}\label{Sta-V6}
\big({\rm Tr}\,u\big)(x):=\lim_{r\to 0_{+}}{-\hskip -0.15in}\int_{B(x,r)\cap\Omega}u(y)\,d^ny,
\,\text{ for $\sigma$-a.e.~$x\in\partial\Omega$},
\end{equation}
induces a well defined, linear, and bounded operator in the context
\begin{equation}\label{fancy-Tr2}
{\rm Tr}:B^{2,1}_{s}(\Omega)\rightarrow H^{s-(1/2)}(\partial\Omega),
\quad\forall\,s\in\big[\tfrac{1}{2},\tfrac{3}{2}\big],
\end{equation}
which is compatible both with the Dirichlet trace operator $\gamma_D$ considered in \eqref{eqn:gammaDs.1}
and with the nontangential boundary trace $u\mapsto u\big|^{\kappa-{\rm n.t.}}_{\partial\Omega}$ 
whenever the latter exists. 
\end{proposition}
%%%%%%%
\begin{proof}
This follows directly from \cite[Proposition~2.61 on p.~107, Remarks~(i)--(ii) on p.~90--91]{MM13} 
specialized to the case $p=2$. 
\end{proof}
%%%%%%%

In relation to \eqref{fancy-Tr2} it is worth pointing out that, as seen from 
\eqref{Inc-P3}--\eqref{q-TWC}, one has  
\begin{equation}\label{fancy-Tr2-ytf}
B^{2,1}_{s}(\Omega)\subsetneq B^{2,2}_{s}(\Omega)=H^{s}(\Omega),
\quad\forall\,s\in\big[\tfrac{1}{2},\tfrac{3}{2}\big].
\end{equation}
Thus, compared to \eqref{eqn:gammaDs.1}, in \eqref{fancy-Tr2} we are now permitted to include 
the end-points of the interval $\big(\tfrac{1}{2},\tfrac{3}{2}\big)$, the price to be paid is 
the consideration of the strictly smaller Besov space $B^{2,1}_{s}(\Omega)$ in place of 
the Sobolev space $H^{s}(\Omega)$ as the domain on which the trace operator now acts.

%%%%%%%%%%%%%%%%%%%%%%%%%%%%
\subsection{A sharp Dirichlet trace involving Sobolev spaces} 
\label{ss3.2}
%%%%%%%%%%%%%%%%%%%%%%%%%%%%

Let $\Omega$ be a bounded Lipschitz domain in ${\mathbb{R}}^n$. 
As already mentioned in the context of Sobolev spaces, the Dirichlet trace operator \eqref{eqn:gammaDs.1} 
fails to be well defined for the end-point cases $s\in\big\{\tfrac{1}{2},\tfrac{3}{2}\big\}$. A remedy that 
allows the inclusion of these prohibitive limiting values is to restrict $\gamma_D$ to a suitably smaller space.

Specifically, starting with $u\in H^{s}(\Omega)$ for some $s\in\big[\tfrac{1}{2},\tfrac{3}{2}\big]$,  
if $\Delta u$ is slightly more regular than the typical action of the Laplacian on functions
from $H^{s}(\Omega)$, that is, more regular than $H^{s-2}(\Omega)$, then one can  
meaningfully define its trace $\gamma_D u$ for the full range $s\in\big[\tfrac{1}{2},\tfrac{3}{2}\big]$.

Here is the theorem about this extended trace result for functions with a better-than-expected 
Laplacian (in the sense of membership to the Sobolev scale). The reader is alerted to the fact that 
having a better-than-expected Laplacian forces the function to be more regular than originally assumed, 
in the manner indicated in \eqref{eqn:gammaDs.2aux.WACO} below. 

%%%%%%%%%%
\begin{theorem}\label{YTfdf-T}
Assume that $\Omega\subset\bbR^n$ is a bounded Lipschitz domain and fix an arbitrary 
$\varepsilon>0$. Then the restriction of the boundary trace operator \eqref{eqn:gammaDs.1} 
to the space $\big\{u\in H^s(\Omega)\,\big|\,\Delta u\in H^{s-2+\varepsilon}(\Omega)\big\}$, 
originally considered for $s\in\big(\tfrac{1}{2},\tfrac{3}{2}\big)$, induces 
a well defined, linear, continuous operator 
\begin{equation}\label{eqn:gammaDs.2aux}
\gamma_D:\big\{u\in H^s(\Omega)\,\big|\,\Delta u\in H^{s-2+\varepsilon}(\Omega)\big\}
\rightarrow H^{s-(1/2)}(\partial\Omega),\quad\forall\,s\in\big[\tfrac{1}{2},\tfrac{3}{2}\big],    
\end{equation}
{\rm (}throughout, the space on the left-hand side of \eqref{eqn:gammaDs.2aux} is equipped with the natural graph norm 
$u\mapsto\|u\|_{H^{s}(\Omega)}+\|\Delta u\|_{H^{s-2+\varepsilon}(\Omega)}${\rm )}, which continues 
to be compatible with \eqref{eqn:gammaDs.1} when $s\in\big(\tfrac{1}{2},\tfrac{3}{2}\big)$.
Thus defined, the Dirichlet trace operator possesses the following additional properties: \\[1mm] 
$(i)$ The Dirichlet boundary trace operator in \eqref{eqn:gammaDs.2aux} is surjective.
In fact, there exist linear and bounded operators
\begin{equation}\label{2.88X-NN-ii-RRDD}
\Upsilon_D:H^{s-(1/2)}(\partial\Omega)\to
\big\{u\in H^s(\Omega)\,\big|\,\Delta u\in L^2(\Omega)\big\},\quad s\in\big[\tfrac12,\tfrac32\big],
\end{equation}
which are compatible with one another and serve as right-inverses for the Dirichlet trace, that is, 
\begin{equation}\label{2.88X-NN2-ii-RRDD}
\gamma_D(\Upsilon_D\psi)=\psi,
\quad\forall\,\psi\in H^{s-(1/2)}(\partial\Omega)\,\text{ with }\,s\in\big[\tfrac12,\tfrac32\big].
\end{equation}
In fact, matters may be arranged so that each function in the range of $\Upsilon_D$ is harmonic, that is, 
\begin{equation}\label{2.88X-NN2-ii-RRDD.bis}
\Delta(\Upsilon_D\psi)=0,\quad\forall\,\psi\in H^{s-(1/2)}(\partial\Omega)
\,\text{ with }\,s\in\big[\tfrac12,\tfrac32\big].
\end{equation}
$(ii)$ The Dirichlet boundary trace operator \eqref{eqn:gammaDs.2aux} is compatible
with the pointwise nontangential trace in the sense that, given any aperture parameter $\kappa>0$, 
\begin{align}\label{eqn:gammaDs.2auxBBB}
\begin{split}
& \text{if $u\in H^s(\Omega)$ has $\Delta u\in H^{s-2+\varepsilon}(\Omega)$ for some 
$s\in\big[\tfrac{1}{2},\tfrac{3}{2}\big]$,}    
\\[2pt]
& \quad\text{and if $u\big|^{\kappa-{\rm n.t.}}_{\partial\Omega}$ exists $\sigma$-a.e.~on $\partial\Omega$,  
then $u\big|^{\kappa-{\rm n.t.}}_{\partial\Omega}=\gamma_D u\in H^{s-(1/2)}(\partial\Omega)$.}
\end{split}
\end{align}
$(iii)$ The Dirichlet boundary trace operator $\gamma_D$ in 
\eqref{eqn:gammaDs.2aux} is the unique extension by continuity and density of the mapping 
$C^\infty(\overline{\Omega})\ni f\mapsto f\big|_{\partial\Omega}$. \\[1mm] 
$(iv)$ For each $s\in\big[\tfrac{1}{2},\tfrac{3}{2}\big]$, the Dirichlet 
boundary trace operator satisfies
\begin{align}\label{eqn:gammaDs.1-TR.2}
\begin{split}
& \gamma_D(\Phi u)=\big(\Phi\big|_{\partial\Omega}\big)\gamma_D u
\,\text{ at $\sigma$-a.e.~point on $\partial\Omega$, for all}
\\[2pt] 
& \quad u\in H^s(\Omega)\,\text{ with }\,\Delta u\in H^{s-2+\varepsilon}(\Omega)
\,\text{ and all }\,\Phi\in C^\infty(\overline\Omega).
\end{split}
\end{align}
$(v)$ For each $s\in\big[\tfrac{1}{2},\tfrac{3}{2}\big]$, and each $\varepsilon>0$ 
such that $\varepsilon\not=\tfrac{3}{2}-s$, the null space   of the Dirichlet boundary 
trace operator \eqref{eqn:gammaDs.2aux} satisfies
\begin{equation}\label{eq:EFFa}
{\rm ker}(\gamma_D)\subseteq H^{\,\min\!\{s+\varepsilon,3/2\}}(\Omega).
\end{equation}
In fact, the inclusion in \eqref{eq:EFFa} is quantitative in the sense that,
whenever $s\in\big[\tfrac{1}{2},\tfrac{3}{2}\big]$ and $\varepsilon>0$ is such that $\varepsilon\not=\tfrac{3}{2}-s$, 
then there exists a constant $C\in(0,\infty)$ with the property that
\begin{align}\label{gafvv.655} 
& \text{if $u\in H^s(\Omega)$ satisfies $\Delta u\in H^{s-2+\varepsilon}(\Omega)$ and $\gamma_D u=0$}   
\nonumber\\[2pt]
& \quad\text{then the function $u$ belongs to $H^{\,\min\!\{s+\varepsilon,3/2\}}(\Omega)$ and}    
\\[2pt]
& \quad\text{$\|u\|_{H^{\,\min\!\{s+\varepsilon,3/2\}}(\Omega)}\leq C\big(\|u\|_{H^s(\Omega)}
+\|\Delta u\|_{H^{s-2+\varepsilon}(\Omega)}\big)$.}   
\nonumber 
\end{align}
$(vi)$ For each $s\in\big[\tfrac{1}{2},\tfrac{3}{2}\big]$, the space on the left-hand side of 
\eqref{eqn:gammaDs.2aux} {\rm (}equipped with the natural graph norm\,{\rm )} embeds continuously into 
the Triebel--Lizorkin space $F^{2,q}_s(\Omega)$ for any $q\in(0,\infty)$. In particular, one has the continuous strict embeddings 
\begin{align}\label{eqn:gammaDs.2aux.WACO}
\begin{split} 
& \big\{u\in H^s(\Omega)\,\big|\,\Delta u\in H^{s-2+\varepsilon}(\Omega)\big\}
\hookrightarrow F^{2,q}_s(\Omega)\hookrightarrow H^s(\Omega)
\\[1pt]
& \quad\text{for any }\,s\in\big[\tfrac{1}{2},\tfrac{3}{2}\big]\,\text{ and any }\,q\in(0,2).
\end{split}
\end{align}
$(vii)$ The operator
\begin{equation}\label{eqn:gammaDs.2aux.INTRO.WACO.b}
\big\{u\in H^{3/2}(\Omega)\,\big|\,\Delta u\in H^{-(1/2)+\varepsilon}(\Omega)\big\}
\ni u\mapsto\gamma_D(\nabla u)\in [L^2(\partial\Omega)]^n 
\end{equation}
{\rm (}with the Dirichlet trace acting componentwise, in the sense of \eqref{eqn:gammaDs.2aux} 
with $s:= 1/2${\rm ),} is well defined, linear, and bounded. 
\end{theorem}
%%%%%%%%%%
\begin{proof}
We split the proof of the claims in the opening part of the statement of the theorem 
into the following three cases: \\[1mm] 
\noindent{\bf Case~1:} {\it Assume $s\in\big(\tfrac{1}{2},\tfrac{3}{2}\big)$}. 
Since $\big\{u\in H^s(\Omega)\,\big|\,\Delta u\in H^{s-2+\varepsilon}(\Omega)\big\}\subset H^s(\Omega)$, 
we let $\gamma_D$ in \eqref{eqn:gammaDs.2aux} act in the same manner as the trace operator 
from \eqref{eqn:gammaDs.1}. This, by design, ensures that $\gamma_D$ is well defined, linear, 
continuous, and compatible with its restrictions defined previously. 
\\[1mm] 
\noindent{\bf Case~2:} {\it Assume $s=\tfrac{3}{2}$}. Given that 
$\big\{u\in H^{3/2}(\Omega)\,\big|\,\Delta u\in H^{-(1/2)+\varepsilon}(\Omega)\big\}\subset H^1(\Omega)$, 
we once again let $\gamma_D$ in \eqref{eqn:gammaDs.2aux} act in the same fashion as the trace operator 
from \eqref{eqn:gammaDs.1} (when $s=1$). Of course, this choice ensures linearity and compatibility. 
We claim that there exists a constant $C\in(0,\infty)$ with the property that
\begin{align}\label{eq:DDa}
\begin{split}
& \text{if $u\in H^{3/2}(\Omega)$ has $\Delta u\in H^{-(1/2)+\varepsilon}(\Omega)$ 
for some $\varepsilon>0$, then actually}  
\\[2pt]
& \quad\text{$\gamma_D u\in H^1(\partial\Omega)$ with 
$\|\gamma_D u\|_{H^1(\partial\Omega)}\leq C\big(\|u\|_{H^{3/2}(\Omega)} 
+\|\Delta u\|_{H^{-(1/2)+\varepsilon}(\Omega)}\big)$.}
\end{split}
\end{align}
To justify this claim, let $u$ be as in the first line of \eqref{eq:DDa} and solve
\begin{equation}\label{eqn:bvp-DIR.22}
\begin{cases}
\Delta v=\Delta u\,\text{ in $\Omega,\quad v\in H^{3/2}(\Omega)$,}     
\\[2pt]
\gamma_D v=0\,\text{ on $\partial\Omega$,}
\end{cases}    
\end{equation}
by proceeding as follows. First, it is possible to extend $\Delta u\in H^{-(1/2)+\varepsilon}(\Omega)$ 
to a compactly supported distribution $U\in H^{-(1/2)+\varepsilon}(\mathbb{R}^n)$ such that,
for some constant $C\in(0,\infty)$, independent of $u$, one has $\|U\|_{H^{-(1/2)+\varepsilon}(\mathbb{R}^n)}
\leq C\|\Delta u\|_{H^{-(1/2)+\varepsilon}(\Omega)}$ (cf.~\eqref{tgBBn}). As in \eqref{EEE-jussi} let $E_0$ denote 
the standard fundamental solution for the Laplacian in ${\mathbb{R}}^n$, that is,
\begin{equation}\label{EEE}
E_0(x)=\begin{cases} 
\displaystyle{\frac{1}{\omega_{n-1}(2-n)}|x|^{2-n}},\,\text{ if }\,n\geq 3,
\\[12pt] 
\displaystyle{\frac{1}{2\pi}}\,\ln|x|,\,\text{ if }\,n=2,
\end{cases}
\quad\forall\,x\in{\mathbb{R}}^{n}\backslash\{0\},
\end{equation}
where $\omega_{n-1}$ is the surface measure of the unit sphere $\bbS^{n-1}$ in ${\mathbb{R}}^n$. 
Classical Calder\'on--Zygmund theory gives that the operator of convolution with $E_0$ is (locally) 
smoothing of order two on the fractional Sobolev scale. Hence, considering $\eta:=(E_0\ast U)|_{\Omega}$, 
then $\eta\in H^{(3/2)+\varepsilon}(\Omega)$ and $\|\eta\|_{H^{(3/2)+\varepsilon}(\Omega)}
\leq C\|U\|_{H^{-(1/2)+\varepsilon}(\mathbb{R}^n)}$. Moreover, $\Delta\eta=(\Delta E_0\ast U)|_{\Omega}=U|_{\Omega}=\Delta u$ 
in $\Omega$. In addition, by \eqref{eq:VCca.1}, one has $\gamma_D\eta\in H^1(\partial\Omega)$ and 
$\|\gamma_D\eta\|_{H^1(\partial\Omega)}\leq C\|\eta\|_{H^{(3/2)+\varepsilon}(\Omega)}$.  
Second, from \cite{JK81}, \cite{Ve84}, one knows that for each aperture parameter $\kappa>0$ the boundary value problem 
\begin{equation}\label{eqn:bvp-REG.222}
\begin{cases}
\Delta h=0\,\text{ in $\Omega,\quad{\mathcal{N}}_\kappa h,
{\mathcal{N}}_\kappa(\nabla h)\in L^2(\partial\Omega)$,}     
\\[2pt]  
h\big|^{\kappa-{\rm n.t.}}_{\partial\Omega}=\gamma_D\eta\,\text{ $\sigma$-a.e.~on $\partial\Omega$,}
\end{cases}    
\end{equation}
has a unique solution, satisfying the naturally accompanying estimate
\begin{equation}\label{eqn:bvp-REG.222.aa}
\big\|{\mathcal{N}}_\kappa h\big\|_{L^2(\partial\Omega)}
+\big\|{\mathcal{N}}_\kappa(\nabla h)\big\|_{L^2(\partial\Omega)}
\leq C\|\gamma_D\eta\|_{H^1(\partial\Omega)},
\end{equation}
for some constant $C\in(0,\infty)$ independent of $\eta$. 
Due to \eqref{eq:MM3BIS}--\eqref{eq:MM3BIS-AC} (with $k=1$) one concludes that
$h\in H^{3/2}(\Omega)$, and from \eqref{eq:MM3BIS-AC} and \eqref{eqn:bvp-REG.222.aa} 
one obtains the estimate $\|h\|_{H^{3/2}(\Omega)}\leq C\|\gamma_D\eta\|_{H^1(\partial\Omega)}$. 
Keeping in mind \eqref{eq:MM3bis}, one then deduces that the function $v:=(\eta-h)\in H^{3/2}(\Omega)$ 
solves \eqref{eqn:bvp-DIR.22}. For later reference we note that 
\begin{align}\label{eqHmn}
\|v\|_{H^{3/2}(\Omega)} & \leq\|\eta\|_{H^{3/2}(\Omega)}+\|h\|_{H^{3/2}(\Omega)}
\nonumber\\[2pt]
&\leq C\big(\|U\|_{H^{-(1/2)+\varepsilon}(\mathbb{R}^n)}+\|\gamma_D\eta\|_{H^1(\partial\Omega)}\big)
\nonumber\\[2pt] 
&\leq C\big(\|\Delta u\|_{H^{-(1/2)+\varepsilon}(\Omega)}
+\|\eta\|_{H^{(3/2)+\varepsilon}(\Omega)}\big)
\nonumber\\[2pt]
&\leq C\|\Delta u\|_{H^{-(1/2)+\varepsilon}(\Omega)}.
\end{align}

Next, with $v$ as in \eqref{eqn:bvp-DIR.22}, consider $w:=u-v\in H^{3/2}(\Omega)$
and note that $\Delta w=\Delta u-\Delta v=0$ in $\Omega$. In particular, $w\in C^\infty(\Omega)$ by elliptic regularity. 
Given these facts, it follows from \eqref{eq:MM3BIS} and \eqref{eq:MM4} that 
\begin{align}\label{eq:MM4.bb}
\begin{split}
& {\mathcal{N}}_\kappa w,{\mathcal{N}}_\kappa(\nabla w)\in L^2(\partial\Omega)
\,\text{ and both }\,w\big|^{\kappa-{\rm n.t.}}_{\partial\Omega}\,\text{ and }\, 
\nabla w\big|^{\kappa-{\rm n.t.}}_{\partial\Omega}\,\text{ exist at}    
\\[2pt]
& \quad\text{$\sigma$-a.e.~point on $\partial\Omega$ and lie in $L^2(\partial\Omega)$
and $\big[L^2(\partial\Omega)\big]^n$, respectively}.
\end{split}
\end{align}
Moreover, \eqref{eqHmn} and the definition of $w$ imply 
\begin{align}\label{Yrr6433}
\|w\|_{H^{3/2}(\Omega)} &\leq\|u\|_{H^{3/2}(\Omega)}+\|v\|_{H^{3/2}(\Omega)}
\nonumber\\[2pt]
& \leq C\big(\|u\|_{H^{3/2}(\Omega)}+\|\Delta u\|_{H^{-(1/2)+\varepsilon}(\Omega)}\big).
\end{align}
Next, fix $j,k\in\{1,\dots,n\}$ along with some arbitrary $\psi\in C^\infty_0({\mathbb{R}}^n)$,
and consider the vector fields defined in $\Omega$ as
\begin{equation}\label{eq:DccB.1}
\vec{F}:=w\,\partial_k\psi\,e_j-w\,\partial_j\psi\,e_k,\quad
\vec{G}:=\psi\,\partial_j w\,e_k-\psi\,\partial_k w\,e_j,
\end{equation}
where $\{e_m\}_{1\leq m\leq n}$ is the standard orthonormal basis in ${\mathbb{R}}^n$.
From \eqref{eq:MM4.bb}, \eqref{eq:DccB.1}, and \eqref{eq:HBab.2}, one deduces that
\begin{align}\label{banff-4}
\begin{split}
& \vec{F},\vec{G}\in\big[L^1_{\rm loc}(\Omega)\big]^n\,\text{ and }\, 
{\mathcal{N}}_\kappa(\vec{F}),{\mathcal{N}}_\kappa(\vec{G})
\in L^2(\partial\Omega)\subset L^1(\partial\Omega), 
\\[2pt]
& \vec{F}\big|^{\kappa-{\rm n.t.}}_{\partial\Omega},\,\vec{G}\big|^{\kappa-{\rm n.t.}}_{\partial\Omega}\,\text{ exist }\,\sigma
\text{-a.e.~on $\partial\Omega$ and lie in $\big[L^2(\partial\Omega)\big]^n\subset\big[L^1(\partial\Omega)\big]^n$},
\\[2pt]
& {\rm div}\,\vec{F},{\rm div}\,\vec{G}\in L^{2n/(n-1)}(\Omega)\subset L^1(\Omega)
\,\text{ and }\,{\rm div}\,\vec{F}={\rm div}\,\vec{G}\,\text{ in }\,\Omega.
\end{split}
\end{align}
Based on these facts and Theorem~\ref{banff-3}, one computes  
\begin{align}\label{eq:MVcAA.4R}
& \Big|\int_{\partial\Omega}(\gamma_D u)\,\partial_{\tau_{jk}}\psi\,d^{n-1}\sigma\Big|
=\bigg|\int_{\partial\Omega}(\gamma_D w)\,\big(\nu_j\big(\partial_k\psi\big)\big|_{\partial\Omega}
-\nu_k\big(\partial_j\psi\big)\big|_{\partial\Omega}\big)\,d^{n-1}\sigma\bigg|
\nonumber\\[2pt]
& \quad=\bigg|\int_{\partial\Omega}\nu\cdot\Big(\vec{F}\,\big|^{\kappa-{\rm n.t.}}_{\partial\Omega}\Big)\,d^{n-1}\sigma\Big|
=\Big|\int_{\Omega}{\rm div}\vec{F}\,d^n x\bigg|
\nonumber\\[2pt]
& \quad=\bigg|\int_{\Omega}{\rm div}\,\vec{G}\,d^n x\Big|
=\Big|\int_{\partial\Omega}\nu\cdot\Big(\vec{G}\,\big|^{\kappa-{\rm n.t.}}_{\partial\Omega}\Big)
\,d^{n-1}\sigma\bigg|
\nonumber\\[2pt]
& \quad=\bigg|\int_{\partial\Omega}\Big[\nu_k\Big((\partial_j w)\big|^{\kappa-{\rm n.t.}}_{\partial\Omega}\Big)
-\nu_j\Big((\partial_k w)\big|^{\kappa-{\rm n.t.}}_{\partial\Omega}\Big)\Big]
\big(\psi\big|_{\partial\Omega}\big)\,d^{n-1}\sigma\bigg|
\nonumber\\[2pt]
& \quad\leq C\Big\|\nabla w\big|^{\kappa-{\rm n.t.}}_{\partial\Omega}\Big\|_{[L^2(\partial\Omega)]^n}
\big\|\psi\big|_{\partial\Omega}\big\|_{L^2(\partial\Omega)}
\nonumber\\[2pt] 
& \quad\leq C\big\|{\mathcal{N}}_\kappa(\nabla w)\big\|_{L^2(\partial\Omega)}
\big\|\psi\big|_{\partial\Omega}\big\|_{L^2(\partial\Omega)}
\nonumber\\[2pt]
& \quad\leq C\|w\|_{H^{3/2}(\Omega)}\big\|\psi\big|_{\partial\Omega}\big\|_{L^2(\partial\Omega)}
\nonumber\\[2pt]
& \quad\leq C\big(\|u\|_{H^{3/2}(\Omega)}
+\|\Delta u\|_{H^{-(1/2)+\varepsilon}(\Omega)}\big)\big\|\psi\big|_{\partial\Omega}\big\|_{L^2(\partial\Omega)},
\end{align}
where the second inequality comes from \eqref{eq:HByrr}, and the penultimate inequality uses \eqref{eq:MM3BIS-AC}. 
In light of the characterization of $H^1(\partial\Omega)$ proved in Lemma~\ref{Dida} 
(cf.~\eqref{eq:Mvc}) and \eqref{eq:VCca.1-tt5}, estimate \eqref{eq:MVcAA.4R} shows that 
the claim in \eqref{eq:DDa} holds. In turn, this implies that the operator $\gamma_D$ 
in \eqref{eqn:gammaDs.2aux} is well defined and continuous when $s=\tfrac{3}{2}$. \\[1mm] 
\noindent{\bf Case~3:} {\it Assume $s=\tfrac{1}{2}$}. 
In this scenario, $\big\{u\in H^{1/2}(\Omega)\,\big|\,\Delta u\in H^{-(3/2)+\varepsilon}(\Omega)\big\}$ 
is not included in $\bigcup_{\frac{1}{2}<s<\frac{3}{2}}H^s(\Omega)$, so we start by assigning 
a meaning to the action of the Dirichlet trace $\gamma_D$ in \eqref{eqn:gammaDs.2aux} when $s=\tfrac{1}{2}$.
Specifically, assume that $u\in H^{1/2}(\Omega)$ satisfies 
$\Delta u\in H^{-(3/2)+\varepsilon}(\Omega)$ for some $\varepsilon\in(0,1)$ 
(which suffices for our purposes). Invoke \cite[Theorem~0.5(b), pp.~164--165]{JK95} to solve 
\begin{equation}\label{eqn:bvp-DIR.2}
\begin{cases}
\Delta v=\Delta u\in H^{-(3/2)+\varepsilon}(\Omega)
\,\text{ in $\Omega,\quad v\in H^{(1/2)+\varepsilon}(\Omega)$,}     
\\[2pt]
\gamma_D v=0\,\text{ on $\partial\Omega$,}
\end{cases}  
\end{equation}
with the Dirichlet trace understood in the sense of \eqref{eqn:gammaDs.1}. 
The solution $v$ is unique and satisfies a naturally accompanying estimate, namely
\begin{equation}\label{eqn:bvp-DIR.2ggd}
\|v\|_{H^{(1/2)+\varepsilon}(\Omega)}\leq C\|\Delta u\|_{H^{-(3/2)+\varepsilon}(\Omega)}
\end{equation}
for some $C\in(0,\infty)$ independent of $u,v$. To proceed, consider 
\begin{equation}\label{eq:ASSD.eee}
w:=u-v\,\text{ in }\,\Omega. 
\end{equation}
Then, by design, $w\in H^{1/2}(\Omega)$ and $\Delta w=0$ in $\Omega$, (hence also $w\in C^\infty(\Omega)$, by elliptic regularity). 
Given these facts, \eqref{eq:MM3} implies that ${\mathcal{N}}_\kappa w\in L^2(\partial\Omega)$.
Together with the Fatou-type result recorded in \eqref{eq:MM4} this ensures that
\begin{equation}\label{eq:ASSD}
w\big|^{\kappa-{\rm n.t.}}_{\partial\Omega}\,\text{ exists at $\sigma$-a.e.~point on 
$\partial\Omega$ and $w\big|^{\kappa-{\rm n.t.}}_{\partial\Omega}\in L^2(\partial\Omega)$}.
\end{equation}
Then we define the action of the Dirichlet trace operator $\gamma_D$ 
from \eqref{eqn:gammaDs.2aux} when $s=\tfrac{1}{2}$ on the function $u$ to be precisely 
the nontangential pointwise trace of $w$, that is, 
\begin{equation}\label{eq:ee3ee}
\gamma_D u:=w\big|^{\kappa-{\rm n.t.}}_{\partial\Omega}.
\end{equation}
The operator just introduced is well defined, linear, and continuous since, thanks 
to \eqref{eq:ee3ee}, \eqref{eq:HByrr}, \eqref{eq:MM3-AC}, \eqref{eq:ASSD.eee}, \eqref{eq:DDEEn.2}, 
and \eqref{eqn:bvp-DIR.2ggd}, we have
\begin{align}\label{eq:ejgrFb}
\|\gamma_D u\|_{L^2(\partial\Omega)} &=\big\|w\big|^{\kappa-{\rm n.t.}}_{\partial\Omega}\big\|_{L^2(\partial\Omega)}
\leq\big\|{\mathcal{N}}_\kappa w\big\|_{L^2(\partial\Omega)}
\nonumber\\[2pt]
&\leq C\|w\|_{H^{1/2}(\Omega)}\leq C\|u\|_{H^{1/2}(\Omega)}+C\|v\|_{H^{1/2}(\Omega)}
\nonumber\\[2pt]
&\leq C\|u\|_{H^{1/2}(\Omega)}+C\|v\|_{H^{(1/2)+\varepsilon}(\Omega)}
\nonumber\\[2pt]
&\leq C\big(\|u\|_{H^{1/2}(\Omega)}+\|\Delta u\|_{H^{-(3/2)+\varepsilon}(\Omega)}\big),
\end{align}
for some $C\in(0,\infty)$ independent of $u$.
To show that this operator is compatible with the Dirichlet trace from \eqref{eqn:gammaDs.1},
assume that $u\in H^{s}(\Omega)$ for some $s\in\big(\tfrac{1}{2},\tfrac{3}{2}\big)$
satisfies $\Delta u\in H^{-(3/2)+\varepsilon}(\Omega)$ for some $\varepsilon>0$. 
Then, following the same procedure as above that has led to the definition in \eqref{eq:ee3ee}, 
one observes that the function $w$ now exhibits better regularity on the Sobolev scale, 
namely $w\in H^{(1/2)+\delta}(\Omega)$, where $\delta:=\min\{\varepsilon,s-(1/2)\}>0$.
Granted this fact, and employing \eqref{eq:ASSD}, one can invoke \eqref{eq:MM3bis} for $w$ 
in order to conclude that 
\begin{equation}\label{eq:ee3ee.2}
\gamma_D w=w\big|^{\kappa-{\rm n.t.}}_{\partial\Omega}.
\end{equation}
Since by design $u=w+v$ in $\Omega$ and $\gamma_D v=0$, it follows from \eqref{eq:ee3ee.2} that 
$\gamma_D u$ considered in the sense of \eqref{eqn:gammaDs.1} is consistent with
our definition in \eqref{eq:ee3ee}. \\[1mm] 
We now address the claims made in itemized portion of the statement of the theorem: \\[1mm] 
{\it Proof of $(i)$.} 
Given any $s\in\big[\tfrac12,\tfrac32\big]$, consider the operator 
\begin{equation}\label{2.88X-NN-44aaa}
\Upsilon_D:H^{s-(1/2)}(\partial\Omega)\rightarrow
\big\{u\in H^s(\Omega)\,\big|\,\Delta u=0\,\text{ in }\,\Omega\big\}
\end{equation}
given by $\Upsilon_D\varphi:=u$, where $u$ is, respectively, the unique solution of
\begin{equation}\label{2.88X-NN-44aaa.1}
\begin{cases}
\Delta u=0\,\text{ in $\Omega,\quad u\in H^{s}(\Omega)$,}     
\\[2pt]  
\gamma_D u=\varphi\,\text{ on $\partial\Omega,\quad\varphi\in H^{s-(1/2)}(\partial\Omega)$,}
\end{cases}   
\end{equation}
if $s\in\big(\tfrac12,\tfrac32\big)$, of 
\begin{equation}\label{eqn:bvp-DDD}
\begin{cases}
\Delta u=0\,\text{ in $\Omega,\quad{\mathcal{N}}_\kappa u\in L^2(\partial\Omega)$,}     
\\[2pt]   
u\big|^{\kappa-{\rm n.t.}}_{\partial\Omega}=\varphi\,\text{ $\sigma$-a.e.~on $\partial\Omega$,}
\end{cases}    
\end{equation}
if $s=\tfrac12$, and of 
\begin{equation}\label{eqn:bvp-REG}
\begin{cases}
\Delta u=0\,\text{ in $\Omega,\quad{\mathcal{N}}_\kappa u,
{\mathcal{N}}_\kappa(\nabla u)\in L^2(\partial\Omega)$,}     
\\[2pt]  
u\big|^{\kappa-{\rm n.t.}}_{\partial\Omega}=\varphi\,\text{ $\sigma$-a.e.~on $\partial\Omega$,}
\end{cases}    
\end{equation}
if $s=\tfrac32$. That the above Dirichlet boundary value problems
are indeed well posed has been proved in \cite[Theorem~10.1]{FMM98} (for $\tfrac12<s<\tfrac32$) and 
\cite{Ve84} (for $s\in\big\{\tfrac12,\tfrac32\big\}$, utilizing \eqref{eq:MM3} and \eqref{eq:MM3BIS})
via boundary layer potential methods. As such, $\Upsilon_D$ is well defined, linear, and bounded. 
In addition, when considered as a family indexed by the parameter $s\in\big[\tfrac12,\tfrac32\big]$, 
the operators $\Upsilon_D$ act in a coherent fashion. Then from \eqref{eq:MM3bis} and \eqref{eq:ee3ee} 
one deduces that 
\begin{align}\label{2.9NEjht-rre}
\gamma_D(\Upsilon_D\varphi)=\varphi,\quad
\forall\,\varphi\in H^{s-(1/2)}(\partial\Omega)\,\text{ with }\,s\in\big[\tfrac12,\tfrac32\big],
\end{align}
proving \eqref{2.88X-NN2-ii-RRDD}. Of course, this also shows that the Dirichlet 
boundary trace operator $\gamma_D$ is surjective in the context of \eqref{eqn:gammaDs.2aux}. \\[1mm] 
{\it Proof of $(ii)$.} 
We start by considering the case where the function $u\in H^{1/2}(\Omega)$ 
satisfies $\Delta u\in H^{-(3/2)+\varepsilon}(\Omega)$,
and assume that $u\big|^{\kappa-{\rm n.t.}}_{\partial\Omega}$ exists at $\sigma$-a.e.~point on $\partial\Omega$. 
In addition, we recall the function $v$ from \eqref{eqn:bvp-DIR.2} and the function $w$ from \eqref{eq:ASSD.eee}.  
In particular, it follows from \eqref{eq:ASSD} and the current assumptions on $u$ that 
$v\big|^{\kappa-{\rm n.t.}}_{\partial\Omega}$ exists at $\sigma$-a.e.~point on $\partial\Omega$.
Since $v\in H^{(1/2)+\varepsilon}(\Omega)$, this further implies (by  Lemma~\ref{Hgg-uh}) that     
$v\big|^{\kappa-{\rm n.t.}}_{\partial\Omega}=\gamma_D v=0$ at $\sigma$-a.e.~point on $\partial\Omega$. 
Granted this fact, one writes (upon recalling the definition of 
$\gamma_D$ from \eqref{eqn:gammaDs.2aux} in the case $s=\tfrac{1}{2}$\,; cf.~\eqref{eq:ee3ee})
\begin{equation}\label{eq:ee3ee.2bgfd}
u\big|^{\kappa-{\rm n.t.}}_{\partial\Omega}=w\big|^{\kappa-{\rm n.t.}}_{\partial\Omega}
+v\big|^{\kappa-{\rm n.t.}}_{\partial\Omega}=w\big|^{\kappa-{\rm n.t.}}_{\partial\Omega}=\gamma_D u,
\end{equation}
as wanted. To complete the proof of \eqref{eqn:gammaDs.2auxBBB} there remains to observe that when 
$s\in\big(\tfrac12,\tfrac32\big]$ the desired compatibility property follows from the 
manner in which the Dirichlet trace has been defined in \eqref{eqn:gammaDs.2aux} and Lemma~\ref{Hgg-uh}. \\[1mm] 
{\it Proof of $(iii)$.} 
That $\gamma_D$ in \eqref{eqn:gammaDs.2aux} is the unique extension by continuity and density
of the mapping $C^\infty(\overline{\Omega})\ni f\mapsto f\big|_{\partial\Omega}$
follows from Lemma~\ref{Dense-LLLe} and \eqref{eqn:gammaDs.2auxBBB}.
\\[1mm] 
{\it Proof of $(iv)$.}
Pick $u\in H^s(\Omega)$ satisfying $\Delta u\in H^{s-2+\varepsilon}(\Omega)$ 
for some $s\in\big[\tfrac{1}{2},\tfrac{3}{2}\big]$, along 
with some $\Phi\in C^\infty(\overline\Omega)$. By the density result proved in 
Lemma~\ref{Dense-LLLe} there exists a sequence 
$\{u_j\}_{j\in{\mathbb{N}}}\subset C^\infty(\overline{\Omega})$ with the property that
\begin{equation}\label{ajtrtv5}
u_j\to u\,\text{ in }\,H^s(\Omega)\,\text{ and }\,\Delta u_j\to\Delta u
\,\text{ in }\,H^{s-2+\varepsilon}(\Omega),\,\text{ as }\,j\to\infty.
\end{equation}
In particular, $\Phi u_j\to\Phi u$ in $H^s(\Omega)$ and $\Delta(\Phi u_j)\to\Delta(\Phi u)$ 
in $H^{s-2+\varepsilon}(\Omega)$ as $j\to\infty$. On account of the continuity of the Dirichlet 
trace operator, this permits us to write, in the sense of $H^{s-(1/2)}(\partial\Omega)$, 
\begin{align}\label{ajtrtv7}
\gamma_D(\Phi u)&=\lim_{j\to\infty}\gamma_D(\Phi u_j)
=\lim_{j\to\infty}(\Phi u_j)\big|_{\partial\Omega}
\nonumber\\[2pt]
&=\lim_{j\to\infty}\big(\Phi\big|_{\partial\Omega}\big)\gamma_Du_j
=\big(\Phi\big|_{\partial\Omega}\big)\gamma_Du,
\end{align}
as wanted. \\[1mm] 
{\it Proof of $(v)$.}
Fix $s\in\big[\tfrac{1}{2},\tfrac{3}{2}\big]$ such that $\varepsilon\not=\tfrac{3}{2}-s$, and choose some 
\begin{equation}\label{eq:GFav}
u\in H^{s}(\Omega)\,\text{ with }\,\Delta u\in H^{s-2+\varepsilon}(\Omega)
\,\text{ and }\,\gamma_D u=0.
\end{equation}
Next, consider a compactly supported distribution 
$U\in H^{s-2+\varepsilon}(\mathbb{R}^n)$ 
with the property that $U\big|_{\Omega}=\Delta u$ and such that 
$\|U\|_{H^{s-2+\varepsilon}(\mathbb{R}^n)}\leq C\|\Delta u\|_{H^{s-2+\varepsilon}(\Omega)}$
where $C\in(0,\infty)$ is a constant independent of $u$. Then, with $E_0$ as in \eqref{EEE}, 
define $v:=(E_0\ast U)|_{\Omega}\in H^{s+\varepsilon}(\Omega)$ and note that this entails 
$\Delta v=\Delta u$ as well as $\|v\|_{H^{s+\varepsilon}(\Omega)}\leq C\|U\|_{H^{s-2+\varepsilon}(\mathbb{R}^n)}
\leq C\|\Delta u\|_{H^{s-2+\varepsilon}(\Omega)}$ again with $C\in(0,\infty)$ independent of $u$.
Hence, introducing $h:=v-u$, it follows from \eqref{eq:GFav}, \eqref{eqn:gammaDs.1}, 
\eqref{eq:VCca.1}, and Corollary~\ref{tafc-644-TRE} that, for some constant $C\in(0,\infty)$, independent of $u$,
\begin{align}\label{eq:GFav.2}
\begin{split}
& h\in H^{s}(\Omega)\,\text{ satisfies }\,
\|h\|_{H^{s}(\Omega)}\leq C\big(\|u\|_{H^{s}(\Omega)}+\|\Delta u\|_{H^{s-2+\varepsilon}(\Omega)}\big),
\\[2pt]
& \quad\text{has }\,\Delta h=0\,\text{ in }\,\Omega,\,\text{ and }\,  
\gamma_D h=\gamma_D v\in H^{\,\min\{s+\varepsilon-(1/2),1\}}(\partial\Omega).
\end{split}
\end{align}
The regularity results for the Dirichlet problem for the Laplacian from 
\cite{FMM98}, \cite{JK95}, and \cite{Ve84} (cf. also \eqref{eq:MM3BIS}) then imply 
\begin{align}\label{eq:EFFa.2323}
\begin{split}
& h\in H^{\,\min\!\{s+\varepsilon, 3/2\}}(\Omega)\,\text{ and, for $C\in(0,\infty)$, independent of $h$,}
\\[2pt]
& \quad\|h\|_{H^{\,\min\!\{s+\varepsilon, 3/2\}}(\Omega)}\leq C\big(\|h\|_{H^{s}(\Omega)}
+\|\gamma_D h\|_{H^{\,\min\{s+\varepsilon-(1/2),1\}}(\partial\Omega)}\big).
\end{split}
\end{align}
In turn, this forces $u=v-h\in H^{\,\min\!\{s+\varepsilon,3/2\}}(\Omega)$ and, making use
of \eqref{eq:VCca.1htrr} as well as \eqref{eq:GFav.2}, one estimates  
\begin{align}\label{eq:EFFa.annb}
&\|u\|_{H^{\,\min\!\{s+\varepsilon,3/2\}}(\Omega)} 
\leq\|v\|_{H^{\,\min\!\{s+\varepsilon,3/2\}}(\Omega)}
+\|h\|_{H^{\,\min\!\{s+\varepsilon,3/2\}}(\Omega)}
\nonumber\\[2pt]
&\quad 
\leq C\|\Delta u\|_{H^{s-2+\varepsilon}(\Omega)}
+C\big(\|h\|_{H^{s}(\Omega)}
+\|\gamma_D h\|_{H^{\,\min\{s+\varepsilon-(1/2),1\}}(\partial\Omega)}\big)
\nonumber\\[2pt]
&\quad 
=C\|\Delta u\|_{H^{s-2+\varepsilon}(\Omega)}
+C\big(\|h\|_{H^{s}(\Omega)}
+\|\gamma_D v\|_{H^{\,\min\{s+\varepsilon-(1/2),1\}}(\partial\Omega)}\big)
\nonumber\\[2pt]
&\quad 
\leq C\|\Delta u\|_{H^{s-2+\varepsilon}(\Omega)}
+C\big(\|h\|_{H^{s}(\Omega)}+\|v\|_{H^{s+\varepsilon}(\Omega)}\big)
\nonumber\\[2pt]
&\quad 
\leq C\|\Delta u\|_{H^{s-2+\varepsilon}(\Omega)}+C\|h\|_{H^{s}(\Omega)}
\nonumber\\[2pt]
&\quad 
\leq C\big(\|u\|_{H^{s}(\Omega)}+\|\Delta u\|_{H^{s-2+\varepsilon}(\Omega)}\big),
\end{align}
for some constant $C\in(0,\infty)$, independent of $u$.
This justifies \eqref{eq:EFFa}, as well as the claim in \eqref{gafvv.655}. \\[1mm] 
{\it Proof of $(vi)$.}
Fix $q\in(0,\infty)$, assume $s\in\big[\tfrac{1}{2},\tfrac{3}{2}\big]$ and $u\in H^s(\Omega)$ is such that 
$\Delta u\in H^{s-2+\varepsilon}(\Omega)$. Since $H^s(\Omega)=B^{2,2}_s(\Omega)=F^{2,2}_s(\Omega)$ 
(with equivalent norms) and 
\begin{align}\label{eq:EFFa.annb.WACO}
H^{s-2+\varepsilon}(\Omega)=B^{2,2}_{s-2+\varepsilon}(\Omega)=F^{2,2}_{s-2+\varepsilon}(\Omega)
\hookrightarrow F^{2,\infty}_{s-2+\varepsilon}(\Omega)\hookrightarrow F^{2,q}_{s-2}(\Omega),
\end{align}
(cf.~\eqref{Inc-P3}, \eqref{LLa-45.222.WACO}, and \eqref{q-TWU-iii.WACO}), one concludes that 
$u\in F^{2,2}_s(\Omega)$, $\Delta u\in F^{2,q}_{s-2}(\Omega)$, and there exists $C\in(0,\infty)$,  
independent of $u$,  such that 
\begin{equation}\label{eq:H.WACO.1.ii}
\|u\|_{F^{2,2}_s(\Omega)}\leq C\|u\|_{H^s(\Omega)},\quad
\|\Delta u\|_{F^{2,q}_{s-2}(\Omega)}\leq C\|\Delta u\|_{H^{s-2+\varepsilon}(\Omega)}.
\end{equation}
Granted these facts, Proposition~\ref{Dense-LLLe-BBB.WACO} applies and yields that 
$u$ belongs to $F^{2,q}_s(\Omega)$ and 
\begin{equation}\label{eq:H.WACO.1.aaa}
\|u\|_{F^{2,q}_s(\Omega)}\leq C\big(\|u\|_{H^s(\Omega)}+\|\Delta u\|_{H^{s-2+\varepsilon}(\Omega)}\big).
\end{equation}
This proves that the space on the left-hand side of \eqref{eqn:gammaDs.2aux}, equipped with the natural graph norm, 
embeds continuously into $F^{2,q}_s(\Omega)$ (from which \eqref{eqn:gammaDs.2aux.WACO} also follows). 

That the first embedding in \eqref{eqn:gammaDs.2aux.WACO} is strict whenever $q\in(0,2)$ is a consequence 
of the fact that there exist functions $u\in F^{2,q}_s(\Omega)$ with $\Delta u\notin H^{s-2+\varepsilon}(\Omega)$.
For example, one may start with $w\in F^{2,q}_{s-2}({\mathbb{R}}^n)\backslash F^{2,2}_{s-2+\varepsilon}({\mathbb{R}}^n)$ 
which has compact support (which may be always arranged via a suitable truncation), then take 
$u:=(E_0\ast w)\big|_{\Omega}$, with $E_0$ as in \eqref{EEE}. Finally, the fact that the second embedding 
in \eqref{eqn:gammaDs.2aux.WACO} is strict whenever $q\in(0,2)$ is seen from \eqref{eq:414151.WACO}, 
\eqref{LLa-45.222.WACO}, and \eqref{Inc-P3}.\\[1mm] 
{\it Proof of $(vii)$.}
Pick some function $u\in H^{3/2}(\Omega)$ satisfying $\Delta u\in H^{-(1/2)+\varepsilon}(\Omega)$
and fix an arbitrary index $j\in\{1,\dots,n\}$. Then \eqref{eq:DDEEn.3} implies that 
$\partial_j u\in H^{1/2}(\Omega)$ and $\|\partial_j u\|_{H^{1/2}(\Omega)}\leq C\|u\|_{H^{3/2}(\Omega)}$
for some constant $C\in(0,\infty)$ independent of $u$. From \eqref{eq:DDEEn.3} and the assumptions made one also infers that
\begin{align}\label{eq:H.WACO.1.aaa.TX.1}
\begin{split}
& \Delta(\partial_j u)=\partial_j(\Delta u)\in H^{-(3/2)+\varepsilon}(\Omega)\,\text{ and}
\\[2pt]
& \quad\|\Delta(\partial_j u)\|_{H^{-(3/2)+\varepsilon}(\Omega)}\leq C\|\Delta u\|_{H^{-(1/2)+\varepsilon}(\Omega)}
\end{split}
\end{align}
again, with $C\in(0,\infty)$ independent of $u$. Together with the fact that \eqref{eqn:gammaDs.2aux} is 
well defined and bounded when $s=\tfrac{1}{2}$, these properties then imply that $\gamma_D(\partial_j u)$ 
belongs to $L^2(\partial\Omega)$ and 
\begin{align}\label{eq:H.WACO.1.aaa.TX.2}
\|\gamma_D(\partial_j u)\|_{L^2(\partial\Omega)} &\leq C\big(\|\partial_j u\|_{H^{1/2}(\Omega)}
+\|\Delta(\partial_j u)\|_{H^{-(3/2)+\varepsilon}(\Omega)}\big)
\nonumber\\[2pt]
&\leq C\big(\|u\|_{H^{3/2}(\Omega)}+\|\Delta u\|_{H^{-(1/2)+\varepsilon}(\Omega)}\big).
\end{align}
All together, this shows that the operator \eqref{eqn:gammaDs.2aux.INTRO.WACO.b} is indeed 
well defined, linear, and bounded. 
\end{proof}
%%%%%%%%%%

For simplicity of notation, we will use the same symbol $\gamma_D$ in connection with 
either \eqref{eqn:gammaDs.1} or \eqref{eqn:gammaDs.2aux}, as the setting in which this 
is used will always be clear from the context. Furthermore, we will continue to employ 
the symbol $\gamma_D$ for vector-valued functions (in which case the Dirichlet trace is applied 
componentwise).

The following special case of Theorem~\ref{YTfdf-T} is particularly useful in applications. 

%%%%%%%%%%
\begin{corollary}\label{YTfdf-T.NNN}
Suppose $\Omega\subset\bbR^n$ is a bounded Lipschitz domain. Then the restriction of the operator 
\eqref{eqn:gammaDs.1} to $\big\{u\in H^s(\Omega)\,\big|\,\Delta u\in L^2(\Omega)\big\}$, originally 
considered for $s\in\big(\tfrac{1}{2},\tfrac{3}{2}\big)$, induces a well defined, linear, continuous 
operator 
\begin{align}\label{eqn:gammaDs.2}
\gamma_D:\big\{u\in H^s(\Omega)\,\big|\,\Delta u\in L^2(\Omega)\big\}
\rightarrow H^{s-(1/2)}(\partial\Omega),\quad\forall\,s\in\big[\tfrac{1}{2},\tfrac{3}{2}\big],    
\end{align}
{\rm (}throughout, the space on the left-hand side of \eqref{eqn:gammaDs.2} being equipped with 
the natural graph norm $u\mapsto\|u\|_{H^{s}(\Omega)}+\|\Delta u\|_{L^{2}(\Omega)}${\rm )}, which 
continues to be compatible with \eqref{eqn:gammaDs.1} when $s\in\big(\tfrac{1}{2},\tfrac{3}{2}\big)$,
and also with the pointwise nontangential trace, whenever the latter exists.

Moreover, the following additional properties are true:

\begin{enumerate}
\item[(i)] The Dirichlet boundary trace operator in \eqref{eqn:gammaDs.2} is surjective
and, in fact, there exist linear and bounded operators
\begin{equation}\label{2.88X-NN-ii-RRDD-2}
\Upsilon_D:H^{s-(1/2)}(\partial\Omega)\to
\big\{u\in H^s(\Omega)\,\big|\,\Delta u\in L^2(\Omega)\big\},\quad s\in\big[\tfrac12,\tfrac32\big],
\end{equation}
which are compatible with one another and serve as right-inverses for the Dirichlet trace, that is, 
\begin{equation}\label{2.88X-NN2-ii-RRDD-2}
\gamma_D(\Upsilon_D\psi)=\psi,\quad\forall\,\psi\in H^{s-(1/2)}(\partial\Omega)
\,\text{ with }\,s\in\big[\tfrac12,\tfrac32\big].
\end{equation}
Actually, matters may be arranged so that each function in the range of $\Upsilon_D$ is harmonic, that is, 
\begin{equation}\label{2.88X-NN2-ii-RRDD.bis.2}
\Delta(\Upsilon_D\psi)=0,\quad\forall\,\psi\in H^{s-(1/2)}(\partial\Omega)
\,\text{ with }\,s\in\big[\tfrac12,\tfrac32\big].
\end{equation}

\item[(ii)] For each $s\in\big[\tfrac{1}{2},\tfrac{3}{2}\big]$, the null space   
of the Dirichlet boundary trace operator \eqref{eqn:gammaDs.2} satisfies
\begin{equation}\label{eq:EFFa.NNN}
{\rm ker}(\gamma_D)\subseteq H^{3/2}(\Omega).
\end{equation}
In fact, the inclusion in \eqref{eq:EFFa.NNN} is quantitative in the sense that
there exists a constant $C\in(0,\infty)$ with the property that
\begin{align}\label{gafvv.6577}
\begin{split}
&\text{whenever $u\in H^{1/2}(\Omega)$ with $\Delta u\in L^2(\Omega)$ satisfies $\gamma_D u=0$, then}
\\[2pt]
& \quad\text{$u\in H^{3/2}(\Omega)$ and 
$\|u\|_{H^{3/2}(\Omega)}\leq C\big(\|u\|_{L^2(\Omega)}+\|\Delta u\|_{L^2(\Omega)}\big)$.}
\end{split}
\end{align}

\item[(iii)] Regarding the domain of the Dirichlet trace operator in \eqref{eqn:gammaDs.2}, 
one has the continuous strict embedding 
\begin{equation}\label{eqn:gammaDs.2aux.WACO.2} 
\begin{array}{c}
\big\{u\in H^s(\Omega)\,\big|\,\Delta u\in L^2(\Omega)\big\}\hookrightarrow F^{2,q}_s(\Omega)
\\[4pt]
\text{for any }\,s\in\big[\tfrac{1}{2},\tfrac{3}{2}\big]\,\text{ and any }\,q\in(0,\infty). 
\end{array}
\end{equation}

\item[(iv)] The operator
\begin{equation}\label{eqn:gammaDs.2aux.INTRO.WACO.c}
\big\{u\in H^{3/2}(\Omega)\,\big|\,\Delta u\in L^2(\Omega)\big\}
\ni u\mapsto\gamma_D(\nabla u)\in [L^2(\partial\Omega)]^n  
\end{equation}
{\rm (}with the Dirichlet trace considered in the sense of \eqref{eqn:gammaDs.2} with $s:= 1/2${\rm ),}
is well defined, linear, and bounded. 
\end{enumerate}
\end{corollary}
%%%%%%%
\begin{proof}
All claims, up to (and including) \eqref{eq:EFFa.NNN}, 
as well as \eqref{eqn:gammaDs.2aux.WACO.2} and \eqref{eqn:gammaDs.2aux.INTRO.WACO.c}, 
are particular cases of the corresponding statement in Theorem~\ref{YTfdf-T}, choosing 
$\varepsilon=2-s$. To prove \eqref{gafvv.6577}, assume that $u\in H^{1/2}(\Omega)$ satisfies 
$\Delta u\in L^2(\Omega)$ and $\gamma_D u=0$. From \eqref{gafvv.655} with $s=\tfrac12$ 
and $\varepsilon=\tfrac32$ it follows that
\begin{align}\label{gafvv.6577-bb.1}
\text{$u\in H^{3/2}(\Omega)$ and 
$\|u\|_{H^{3/2}(\Omega)}\leq C\big(\|u\|_{H^{1/2}(\Omega)}+\|\Delta u\|_{L^2(\Omega)}\big)$}
\end{align}
for some constant $C\in(0,\infty)$, independent of $u$. 
In view of \eqref{incl-Yb.EE} one therefore has
\begin{align}\label{gafvv.6577-bb}
\text{$u\in\accentset{\circ}{H}^1(\Omega)\cap H^{3/2}(\Omega)$ and 
$\|u\|_{H^{3/2}(\Omega)}\leq C\big(\|u\|_{H^1(\Omega)}+\|\Delta u\|_{L^2(\Omega)}\big)$.}
\end{align}
From \eqref{Rdac} one knows that there exists a sequence 
$\{\varphi_j\}_{j\in{\mathbb{N}}}\subset C^\infty_0(\Omega)$ 
with the property that
\begin{equation}\label{eq:fre322}
\varphi_j\to u\,\text{ in }\,H^1(\Omega)\,\text{ as }\,j\to\infty.
\end{equation}
Thus, one can write
\begin{align}\label{gafvv.6577-cc}
(\Delta u,u)_{L^2(\Omega)} &=\lim_{j\to\infty}(\Delta u,\varphi_j)_{L^2(\Omega)}
=\lim_{j\to\infty}{}_{{\mathcal{D}}'(\Omega)}\big\langle\,\overline{\Delta u},\varphi_j
\big\rangle_{{\mathcal{D}}(\Omega)}
\nonumber\\[2pt]
&=-\lim_{j\to\infty}\sum_{k=1}^n{}_{{\mathcal{D}}'(\Omega)}
\big\langle\,\overline{\partial_k u},\partial_k\varphi_j
\big\rangle_{{\mathcal{D}}(\Omega)}
\nonumber\\[2pt]
&=-\lim_{j\to\infty}\sum_{k=1}^n\big(\partial_k u,\partial_k\varphi_j\big)_{L^2(\Omega)}
\nonumber\\[2pt]
&=-\|\nabla u\|^2_{[L^2(\Omega)]^n}.
\end{align}
This fact and the Cauchy--Schwartz inequality imply 
\begin{align}\label{gafvv.6577-dd}
\|\nabla u\|^2_{[L^2(\Omega)]^n} &\leq\big|(\Delta u,u)_{L^2(\Omega)}\big|
\leq\|\Delta u\|_{L^2(\Omega)}\|u\|_{L^2(\Omega)}
\nonumber\\[2pt]
&\leq\big(\|\Delta u\|_{L^2(\Omega)}+\|u\|_{L^2(\Omega)}\big)^2,
\end{align}
and hence, for some dimensional constant $C\in(0,\infty)$, 
\begin{align}\label{gafvv.6577-ee}
\|u\|_{H^1(\Omega)}\leq C\big(\|u\|_{L^2(\Omega)}+\|\Delta u\|_{L^2(\Omega)}\big).
\end{align}
When used back in \eqref{gafvv.6577-bb}, this yields the estimate in \eqref{gafvv.6577}.
\end{proof}
%%%%%%%

Once again, we will continue to employ the same symbol $\gamma_D$ as before in connection with 
the operator in \eqref{eqn:gammaDs.2}.

%%%%%%%%%%%%%%%%%%%%%
\subsection{A sharp Dirichlet trace involving Besov spaces} 
\label{ss3.3}
%%%%%%%%%%%%%%%%%%%%%

We are now ready to study the Dirichlet boundary 
trace operator in the Besov space context. 
In a nutshell, the next theorem asserts that 
given any bounded Lipschitz domain $\Omega$ in ${\mathbb{R}}^n$, the Dirichlet boundary 
trace operator $\gamma_D$ considered in Theorem~\ref{YTfdf-T} extends to a linear and bounded
mapping on the hybrid space $H\!B^{s}_{\Delta}(\Omega)$ defined in Lemma~\ref{Dense-LLLe-BBB}
for each $s\in\big[\tfrac{1}{2},\tfrac{3}{2}\big]$, while at the same time retaining all 
the nice features shared by $\gamma_D$ in the previous smaller setting.  Indeed, for each $s\in\big[\tfrac{1}{2},\tfrac{3}{2}\big]$
and $\varepsilon>0$ one obtains 
\begin{equation}\label{eqn:gammaDs.2aux-YRE}
\big\{u\in H^s(\Omega)\,|\,\Delta u\in H^{s-2+\varepsilon}(\Omega)\big\}\subsetneq H\!B^{s}_{\Delta}(\Omega)
\end{equation}
since by \eqref{q-TWU-yt} one has
\begin{equation}\label{eqn:gammaDs.2aux-YRE.2}
H^{s-2+\varepsilon}(\Omega)\subsetneq B^{2,1}_{s-2}(\Omega).
\end{equation}
So, while Theorem~\ref{YTfdf-T} pertaining to the nature of $\gamma_D$ is optimal as far as the 
Sobolev scale is concerned, the consideration of the hybrid scale $H\!B^{s}_{\Delta}(\Omega)$, 
involving Besov spaces, opens the door for pushing this theory to its natural limit. Specifically, 
we have the following result about what we shall refer to as the sharp Dirichlet boundary trace 
operator $\gamma^{\#}_D$.

%%%%%%%
\begin{theorem}\label{YTfdf-T-Mi}
Let $\Omega\subset\bbR^n$ be a bounded Lipschitz domain. Then the boundary trace operator \eqref{eqn:gammaDs.2aux} 
extends to a well defined, linear, continuous mapping 
\begin{equation}\label{eqn:gammaDs.2aux-Mi}
\gamma^{\#}_D:\big\{u\in H^s(\Omega)\,\big|\,\Delta u\in B^{2,1}_{s-2}(\Omega)\big\}
\rightarrow H^{s-(1/2)}(\partial\Omega),\quad\forall\,s\in\big[\tfrac{1}{2},\tfrac{3}{2}\big],    
\end{equation}
when the space on the left-hand side of \eqref{eqn:gammaDs.2aux-Mi} is equipped with the natural graph 
norm $u\mapsto\|u\|_{H^{s}(\Omega)}+\|\Delta u\|_{B^{2,1}_{s-2}(\Omega)}$.
Defined as such, this sharp Dirichlet trace operator is compatible with \eqref{eqn:gammaDs.2aux} 
for each $\varepsilon>0$ {\rm (}hence also with \eqref{eqn:gammaDs.1} when 
$s\in\big(\tfrac{1}{2},\tfrac{3}{2}\big)${\rm )}, and possesses the following additional properties: \\[1mm] 
$(i)$ The sharp Dirichlet boundary trace operator \eqref{eqn:gammaDs.2aux-Mi} is surjective.
In fact, there exist linear and bounded operators
\begin{equation}\label{2.88X-NN-ii-RRDD-Mi}
\Upsilon_D:H^{s-(1/2)}(\partial\Omega)\rightarrow 
\big\{u\in H^s(\Omega)\,\big|\,\Delta u=0\big\},\quad s\in\big[\tfrac12,\tfrac32\big],
\end{equation}
which are compatible with one another and serve as right-inverses for the Dirichlet trace, that is, 
\begin{equation}\label{2.88X-NN2-ii-RRDD-Mi}
\gamma^{\#}_D(\Upsilon_D\psi)=\psi,\quad\forall\,\psi\in H^{s-(1/2)}(\partial\Omega)
\,\text{ with }\,s\in\big[\tfrac12,\tfrac32\big].
\end{equation}

\noindent $(ii)$ The sharp Dirichlet boundary trace operator \eqref{eqn:gammaDs.2aux-Mi} is compatible
with the pointwise nontangential trace in the sense that:
\begin{align}\label{eqn:gammaDs.2auxBBB-Mi}
\begin{split}
& \text{if $u\in H^s(\Omega)$ has $\Delta u\in B^{2,1}_{s-2}(\Omega)$ for some 
$s\in\big[\tfrac{1}{2},\tfrac{3}{2}\big]$ and if}    
\\[2pt]
& \quad\text{$u\big|^{\kappa-{\rm n.t.}}_{\partial\Omega}$ exists $\sigma$-a.e.~on $\partial\Omega$,
then $u\big|^{\kappa-{\rm n.t.}}_{\partial\Omega}=\gamma^{\#}_D u\in H^{s-(1/2)}(\partial\Omega)$.}
\end{split}
\end{align}

\noindent $(iii)$ The sharp Dirichlet boundary trace operator $\gamma^{\#}_D$ in 
\eqref{eqn:gammaDs.2aux-Mi} is the unique extension by continuity and density 
of the mapping $C^\infty(\overline{\Omega})\ni f\mapsto f\big|_{\partial\Omega}$. 
\\[1mm] 
$(iv)$ For each $s\in\big[\tfrac{1}{2},\tfrac{3}{2}\big]$, the sharp Dirichlet 
boundary trace operator satisfies
\begin{align}\label{eqn:gammaDs.1-TR.2-Mi}
\begin{split}
& \gamma^{\#}_D(\Phi u)=\big(\Phi\big|_{\partial\Omega}\big)\gamma_D u
\,\text{ at $\sigma$-a.e.~point on $\partial\Omega$, for all}
\\[2pt]
& \quad u\in H^s(\Omega)\,\text{ with }\,\Delta u\in B^{2,1}_{s-2}(\Omega)
\,\text{ and all }\,\Phi\in C^\infty(\overline\Omega).
\end{split}
\end{align}
$(v)$ For each $s\in\big[\tfrac{1}{2},\tfrac{3}{2}\big]$, the space on the left-hand side of 
\eqref{eqn:gammaDs.2aux-Mi} {\rm (}equipped with the natural graph norm\,{\rm )} embeds continuously into 
the Triebel--Lizorkin space $F^{2,1}_s(\Omega)$. In particular, one has the continuous strict embeddings 
\begin{equation}\label{eqn:gammaDs.2aux.WACO.BBBB}
\big\{u\in H^s(\Omega)\,\big|\,\Delta u\in B^{2,1}_{s-2}(\Omega)\big\}
\hookrightarrow F^{2,1}_s(\Omega)\hookrightarrow H^s(\Omega),\quad s\in\big[\tfrac{1}{2},\tfrac{3}{2}\big].
\end{equation}
$(vi)$ The operator
\begin{equation}\label{eqn:gammaDs.2aux.INTRO.WACO.d}
\big\{u\in H^{3/2}(\Omega)\,\big|\,\Delta u\in B^{2,1}_{-1/2}(\Omega)\big\}
\ni u\mapsto\gamma^{\#}_D(\nabla u)\in [L^2(\partial\Omega)]^n  
\end{equation}
{\rm (}with the sharp Dirichlet trace acting componentwise, in the sense of \eqref{eqn:gammaDs.2aux-Mi} 
with $s:= 1/2${\rm ),} is well defined, linear, and bounded. 
\end{theorem}
%%%%%%%
\begin{proof}
We split the proof of the claims in the opening part of the statement of the theorem 
into three cases, starting with: \\[1mm] 
\noindent{\bf Case~1:} {\it Assume $s\in\big(\tfrac{1}{2},\tfrac{3}{2}\big)$}. 
Since $\big\{u\in H^s(\Omega)\,\big|\,\Delta u\in B^{2,1}_{s-2}(\Omega)\big\}\subset H^s(\Omega)$, 
we let $\gamma^{\#}_D$ in \eqref{eqn:gammaDs.2aux-Mi} act in the same manner as the trace operator 
from \eqref{eqn:gammaDs.1}. This, by design, ensures that $\gamma^{\#}_D$ is well defined, linear, 
continuous, and compatible with its restrictions defined previously. 
\\[1mm] 
\noindent{\bf Case~2:} {\it Assume $s=\tfrac{3}{2}$}. Given that 
$\big\{u\in H^{3/2}(\Omega)\,\big|\,\Delta u\in B^{2,1}_{-1/2}(\Omega)\big\}\subset H^1(\Omega)$, 
we once again let $\gamma^{\#}_D$ in \eqref{eqn:gammaDs.2aux-Mi} act in the same fashion as the trace operator 
from \eqref{eqn:gammaDs.1} (when $s=1$). Of course, this choice ensures linearity and compatibility. 
We claim that there exists a constant $C\in(0,\infty)$ with the property that
\begin{align}\label{eq:DDa-Mi}
\begin{split}
& \text{if $u\in H^{3/2}(\Omega)$ has the property that $\Delta u\in B^{2,1}_{-1/2}(\Omega)$ then actually}  
\\[2pt]
& \quad\text{$\gamma^{\#}_D u\in H^1(\partial\Omega)$ with 
$\|\gamma^{\#}_D u\|_{H^1(\partial\Omega)}\leq C\big(\|u\|_{H^{3/2}(\Omega)} 
+\|\Delta u\|_{B^{2,1}_{-1/2}(\Omega)}\big)$.}
\end{split}
\end{align}
To justify this claim, fix a function $u\in H^{3/2}(\Omega)$ with 
$\Delta u\in B^{2,1}_{-1/2}(\Omega)$ and solve
\begin{equation}\label{eqn:bvp-DIR.22-Mi}
\begin{cases}
\Delta v=\Delta u\,\text{ in $\Omega,\quad v\in H^{3/2}(\Omega)$,}     
\\[2pt]
\gamma_D v=0\,\text{ on $\partial\Omega$,}
\end{cases}    
\end{equation}
by proceeding as follows. First, it is possible to extend $\Delta u\in B^{2,1}_{-1/2}(\Omega)$ 
to a compactly supported distribution $U\in B^{2,1}_{-1/2}(\mathbb{R}^n)$ such that,
for some constant $C\in(0,\infty)$, independent of $u$, one has $\|U\|_{B^{2,1}_{-1/2}(\mathbb{R}^n)}
\leq C\|\Delta u\|_{B^{2,1}_{-1/2}(\Omega)}$ (cf.~\eqref{restr}). We recall that $E_0$ denotes 
the standard fundamental solution for the Laplacian in ${\mathbb{R}}^n$ defined in \eqref{EEE}.
Calder\'on--Zygmund theory then gives that the operator of convolution with $E_0$ is locally 
smoothing of order two on the Besov scale (see, e.g., \cite{KMM07}). Hence, considering $\eta:=(E_0\ast U)|_{\Omega}$, then 
\begin{equation}\label{aaa-Mi}
\eta\in B^{2,1}_{3/2}(\Omega)\subseteq B^{2,2}_{3/2}(\Omega)=H^{3/2}(\Omega),
\end{equation}
and $\|\eta\|_{B^{2,1}_{3/2}(\Omega)}\leq C\|U\|_{B^{2,1}_{-1/2}(\mathbb{R}^n)}$.
Moreover, $\Delta\eta=(\Delta E_0\ast U)|_{\Omega}=U|_{\Omega}=\Delta u$ in $\Omega$.
In addition, Proposition~\ref{P-New-Trace} (cf.~\eqref{Sta-V6}, \eqref{fancy-Tr2}) used with 
$s=\tfrac{3}{2}$ ensures that ${\rm Tr}\,\eta\in H^1(\partial\Omega)$ and 
$\|{\rm Tr}\,\eta\|_{H^1(\partial\Omega)}\leq C\|\eta\|_{B^{2,1}_{3/2}(\Omega)}$.  
Second, from \cite{JK81}, \cite{Ve84}, one obtains the existence of some constant $C\in(0,\infty)$ 
with the property that the boundary value problem 
\begin{equation}\label{eqn:bvp-REG.222-Mi}
\begin{cases}
\Delta h=0\,\text{ in $\Omega,\quad{\mathcal{N}}_\kappa h,
{\mathcal{N}}_\kappa(\nabla h)\in L^2(\partial\Omega)$,}     
\\[2pt]  
h\big|^{\kappa-{\rm n.t.}}_{\partial\Omega}={\rm Tr}\,\eta\,\text{ $\sigma$-a.e.~on $\partial\Omega$,}
\end{cases}    
\end{equation}
has a unique solution, satisfying the naturally accompanying estimate
\begin{equation}\label{eqn:bvp-REG.222.aa-Mi}
\big\|{\mathcal{N}}_\kappa h\big\|_{L^2(\partial\Omega)}
+\big\|{\mathcal{N}}_\kappa(\nabla h)\big\|_{L^2(\partial\Omega)}
\leq C\|{\rm Tr}\,\eta\|_{H^1(\partial\Omega)}.
\end{equation}
Due to \eqref{eq:MM3BIS}--\eqref{eq:MM3BIS-AC} (with $k=1$) one concludes that
$h\in H^{3/2}(\Omega)$, and from \eqref{eq:MM3BIS-AC} and \eqref{eqn:bvp-REG.222.aa-Mi} 
one obtains the estimate $\|h\|_{H^{3/2}(\Omega)}\leq C\|{\rm Tr}\,\eta\|_{H^1(\partial\Omega)}$. 
Keeping in mind \eqref{eq:MM3bis}, the compatibility properties of ${\rm Tr}$ recorded in 
Proposition~\ref{P-New-Trace}, and \eqref{aaa-Mi}, one then deduces that the function 
$v:=\eta-h\in H^{3/2}(\Omega)$ solves \eqref{eqn:bvp-DIR.22-Mi}. For later reference we note that 
\begin{align}\label{eqHmn-Mi}
\|v\|_{H^{3/2}(\Omega)} & \leq\|\eta\|_{H^{3/2}(\Omega)}+\|h\|_{H^{3/2}(\Omega)}
\nonumber\\[2pt]
& \leq C\|\eta\|_{B^{2,1}_{3/2}(\Omega)}+\|h\|_{H^{3/2}(\Omega)}
\nonumber\\[2pt]
&\leq C\big(\|U\|_{B^{2,1}_{-1/2}(\mathbb{R}^n)}+\|{\rm Tr}\,\eta\|_{H^1(\partial\Omega)}\big)
\nonumber\\[2pt] 
&\leq C\big(\|\Delta u\|_{B^{2,1}_{-1/2}(\Omega)}+\|\eta\|_{B^{2,1}_{3/2}(\Omega)}\big)
\nonumber\\[2pt]
&\leq C\|\Delta u\|_{B^{2,1}_{-1/2}(\Omega)}.
\end{align}
Next, with $v$ as in \eqref{eqn:bvp-DIR.22-Mi}, one considers $w:=(u-v)\in H^{3/2}(\Omega)$
and note that 
\begin{equation}\label{bbb-Mi}
\text{$\Delta w=\Delta u-\Delta v=0$ in $\Omega$, and $\gamma_D w=\gamma_D u=\gamma^{\#}_D u$},
\end{equation}
where the last equality is a consequence of the manner in which $\gamma^{\#}_D u$ has been 
defined in the current case. Moreover, \eqref{eqHmn-Mi} and the definition of $w$ imply 
\begin{align}\label{Yrr6433-Mi}
\|w\|_{H^{3/2}(\Omega)} &\leq\|u\|_{H^{3/2}(\Omega)}+\|v\|_{H^{3/2}(\Omega)}
\nonumber\\[2pt]
& \leq C\big(\|u\|_{H^{3/2}(\Omega)}+\|\Delta u\|_{B^{2,1}_{-1/2}(\Omega)}\big).
\end{align}
Applying \eqref{bbb-Mi} and \eqref{eq:DDa} in connection with the function $w$ one then concludes that 
\begin{align}\label{eq:DDa-Mi.2eeD}
\gamma^{\#}_D u=\gamma_D w\in H^1(\partial\Omega)
\end{align}
and 
\begin{align}\label{eq:DDa-Mi.2}
\|\gamma^{\#}_D u\|_{H^1(\partial\Omega)} &=\|\gamma_D w\|_{H^1(\partial\Omega)}
\leq C\|w\|_{H^{3/2}(\Omega)}
\nonumber\\[2pt]
&\leq C\big(\|u\|_{H^{3/2}(\Omega)}+\|\Delta u\|_{B^{2,1}_{-1/2}(\Omega)}\big),
\end{align}
for some constant $C\in(0,\infty)$ independent of $u$. This finishes the proof of 
the claim in \eqref{eq:DDa-Mi}. In turn, \eqref{eq:DDa-Mi} implies that the operator $\gamma^{\#}_D$ 
in \eqref{eqn:gammaDs.2aux-Mi} is well defined and continuous when $s=\tfrac{3}{2}$. \\[1mm] 
\noindent{\bf Case~3:} {\it Assume $s=\tfrac{1}{2}$}. 
In this scenario, $\big\{u\in H^{1/2}(\Omega)\,\big|\,\Delta u\in B^{2,1}_{-3/2}(\Omega)\big\}$ 
is not included in $\bigcup_{\frac{1}{2}<s<\frac{3}{2}}H^s(\Omega)$, so we start by assigning meaning 
to the action of the sharp Dirichlet trace $\gamma^{\#}_D$ in \eqref{eqn:gammaDs.2aux-Mi} 
when $s=\tfrac{1}{2}$. Specifically, assuming that $u\in H^{1/2}(\Omega)$ satisfies 
$\Delta u\in B^{2,1}_{-3/2}(\Omega)$, we extend the latter distribution in $\Omega$ 
to a compactly supported distribution $\widetilde{U}\in B^{2,1}_{-3/2}(\mathbb{R}^n)$ such that,
for some constant $C\in(0,\infty)$, independent of $u$, one has $\|\widetilde{U}\|_{B^{2,1}_{-3/2}(\mathbb{R}^n)}
\leq C\|\Delta u\|_{B^{2,1}_{-3/2}(\Omega)}$ (cf.~\eqref{restr}). Considering 
$\widetilde{\eta}:=(E_0\ast\widetilde{U})|_{\Omega}$, then 
\begin{equation}\label{aaa-Mi.TL}
\widetilde{\eta}\in B^{2,1}_{1/2}(\Omega)\subseteq B^{2,2}_{1/2}(\Omega)=H^{1/2}(\Omega),
\end{equation}
and $\|\widetilde{\eta}\|_{B^{2,1}_{1/2}(\Omega)}\leq C\|\widetilde{U}\|_{B^{2,1}_{-3/2}(\mathbb{R}^n)}$.
Also, 
\begin{equation}\label{aaa-Mi.TL.222}
\Delta\widetilde{\eta}=(\Delta E_0\ast\widetilde{U})|_{\Omega}
=\widetilde{U}|_{\Omega}=\Delta u\,\text{ in }\,\Omega. 
\end{equation}
Moreover, Proposition~\ref{P-New-Trace} used with $s=\tfrac{1}{2}$ ensures that 
${\rm Tr}\,\widetilde{\eta}$ belongs to $L^2(\partial\Omega)$ and 
$\|{\rm Tr}\,\widetilde{\eta}\|_{L^2(\partial\Omega)}\leq C\|\widetilde{\eta}\|_{B^{2,1}_{1/2}(\Omega)}$.  
Second, from \cite{JK81}, \cite{Ve84}, one knows that there exists some constant $C\in(0,\infty)$ 
with the property that the boundary value problem 
\begin{equation}\label{eqn:bvp-REG.222-Mi.TL}
\begin{cases}
\Delta\widetilde{h}=0\,\text{ in $\Omega,\quad{\mathcal{N}}_\kappa\widetilde{h}\in L^2(\partial\Omega)$,}     
\\[4pt]  
\widetilde{h}\big|^{\kappa-{\rm n.t.}}_{\partial\Omega}={\rm Tr}\,\widetilde{\eta}
\,\text{ $\sigma$-a.e.~on $\partial\Omega$,}
\end{cases}    
\end{equation}
has a unique solution, satisfying the naturally accompanying estimate
\begin{equation}\label{eqn:bvp-REG.222.aa-Mi.TL}
\big\|{\mathcal{N}}_\kappa\widetilde{h}\big\|_{L^2(\partial\Omega)}
\leq C\|{\rm Tr}\,\widetilde{\eta}\|_{L^2(\partial\Omega)}.
\end{equation}
Due to \eqref{eq:MM3BIS}--\eqref{eq:MM3BIS-AC} (with $k=0$) one concludes that
$\widetilde{h}\in H^{1/2}(\Omega)$, and from \eqref{eq:MM3BIS-AC} and \eqref{eqn:bvp-REG.222.aa-Mi.TL} 
one obtains the estimate $\|\widetilde{h}\|_{H^{1/2}(\Omega)}\leq C\|{\rm Tr}\,\widetilde{\eta}\|_{L^2(\partial\Omega)}$. 
In turn, from this estimate, \eqref{aaa-Mi.TL}, \eqref{aaa-Mi.TL.222}, our earlier estimates for $\widetilde{\eta}$, 
$\widetilde{U}$, and the boundedness of ${\rm Tr}$ corresponding to $s=\tfrac{1}{2}$ in Proposition~\ref{P-New-Trace}, 
one then deduces that the function 
\begin{equation}\label{eHg-Mi.1}
\widetilde{v}:=\big(\widetilde{\eta}-\widetilde{h}\big)\in H^{1/2}(\Omega)
\end{equation}
satisfies the estimate 
\begin{align}\label{eqHmn-Mi.TL}
\|\widetilde{v}\|_{H^{1/2}(\Omega)} & \leq\|\widetilde{\eta}\|_{H^{1/2}(\Omega)}
+\|\widetilde{h}\|_{H^{1/2}(\Omega)}
\nonumber\\[2pt]
& \leq C\|\widetilde{\eta}\|_{B^{2,1}_{1/2}(\Omega)}+\|\widetilde{h}\|_{H^{1/2}(\Omega)}
\nonumber\\[2pt]
&\leq C\big(\|\widetilde{U}\|_{B^{2,1}_{-3/2}(\mathbb{R}^n)}
+\|{\rm Tr}\,\widetilde{\eta}\|_{L^2(\partial\Omega)}\big)
\nonumber\\[2pt] 
&\leq C\big(\|\Delta u\|_{B^{2,1}_{-3/2}(\Omega)}
+\|\widetilde{\eta}\|_{B^{2,1}_{1/2}(\Omega)}\big)
\nonumber\\[2pt]
&\leq C\|\Delta u\|_{B^{2,1}_{-3/2}(\Omega)}, 
\end{align}
for some constant $C\in(0,\infty)$, independent of $u$, and satisfies
\begin{equation}\label{eHg-Mi}
\Delta\widetilde{v}=\Delta u\,\text{ in }\,\Omega.
\end{equation}
To proceed, one considers  
\begin{equation}\label{eHg-Mi.2}
\widetilde{w}:=u-\widetilde{v}\,\text{ in }\,\Omega. 
\end{equation}
Then, by design, $\widetilde{w}\in H^{1/2}(\Omega)$ and $\Delta\widetilde{w}=0$ in $\Omega$. 
Given these facts, \eqref{eq:MM3} implies that ${\mathcal{N}}_\kappa\widetilde{w}\in L^2(\partial\Omega)$.
Together with the Fatou-type result recorded in \eqref{eq:MM4} this ensures that
\begin{align}\label{eq:ASSD-Mi}
& \text{the nontangential trace $\widetilde{w}\big|^{\kappa-{\rm n.t.}}_{\partial\Omega}$ 
exists at $\sigma$-a.e.~point on $\partial\Omega$,}   
\nonumber\\[2pt] 
& \quad\text{the function $\widetilde{w}\big|^{\kappa-{\rm n.t.}}_{\partial\Omega}$ belongs to 
$L^2(\partial\Omega)$, and one has} 
\\[2pt] 
& \quad\text{$\big\|\widetilde{w}\big|^{\kappa-{\rm n.t.}}_{\partial\Omega}\big\|_{L^2(\partial\Omega)}
\leq C\|{\mathcal{N}}_\kappa\widetilde{w}\|_{L^2(\partial\Omega)}\leq C\|\widetilde{w}\|_{H^{1/2}(\Omega)}$.}   
\nonumber
\end{align}
Then we define the action of the sharp Dirichlet trace operator $\gamma^{\#}_D$ 
from \eqref{eqn:gammaDs.2aux-Mi} when $s=\tfrac{1}{2}$ on the function $u$ to be precisely 
the nontangential pointwise trace of $\widetilde{w}$, that is, 
\begin{equation}\label{eq:ee3ee-Mi}
\gamma^{\#}_D u:=\widetilde{w}\big|^{\kappa-{\rm n.t.}}_{\partial\Omega}.
\end{equation}
The operator just introduced is well defined, linear, and continuous since 
there exists some $C\in(0,\infty)$ independent of $u$ for which one can write 
\begin{align}\label{eq:ejgrFb-Mi}
\|\gamma^{\#}_D u\|_{L^2(\partial\Omega)} 
&=\big\|\widetilde{w}\big|^{\kappa-{\rm n.t.}}_{\partial\Omega}\big\|_{L^2(\partial\Omega)}
\leq C\|\widetilde{w}\|_{H^{1/2}(\Omega)}
\nonumber\\[2pt]
&\leq C\|u\|_{H^{1/2}(\Omega)}+C\|\widetilde{v}\|_{H^{1/2}(\Omega)}
\nonumber\\[2pt]
&\leq C\big(\|u\|_{H^{1/2}(\Omega)}+\|\Delta u\|_{B^{2,1}_{-3/2}(\Omega)}\big).
\end{align}
To show that the operator defined in \eqref{eq:ee3ee-Mi} is compatible with the Dirichlet 
trace from \eqref{eqn:gammaDs.1}, assume that $u\in H^{s}(\Omega)$ with 
$s\in\big(\tfrac{1}{2},\tfrac{3}{2}\big)$. Then $\Delta u\in H^{-(3/2)+\varepsilon}(\Omega)$ 
for some sufficiently small $\varepsilon>0$. Without loss of generality one can assume that $\varepsilon\in(0,1)$. 
Then the functions $\widetilde{\eta}$, $\widetilde{h}$, $\widetilde{v}$, and $\widetilde{w}$ now 
exhibit better regularity on the Sobolev scale than in the previous case. Specifically, in place 
of \eqref{aaa-Mi.TL} one now has $\widetilde{\eta}\in H^{(1/2)+\varepsilon}(\Omega)$, which further 
translates into ${\rm Tr}\,\widetilde{\eta}\in H^{\varepsilon}(\partial\Omega)$. When the latter 
function is regarded as the boundary datum in the Dirichlet problem \eqref{eqn:bvp-REG.222-Mi.TL}, 
this extra regularity forces the solution $\widetilde{h}$ to be correspondingly more regular. 
Indeed, since the solution of that Dirichlet problem is constructed 
via boundary layer potentials, the mapping properties of these integral operators on 
fractional Sobolev spaces established in \cite{FMM98}, \cite{MT00} then imply that 
$\widetilde{h}\in H^{(1/2)+\varepsilon}(\Omega)$. Ultimately, this guarantees that the function
$\widetilde{v}:=\widetilde{\eta}-\widetilde{h}$ belongs to $H^{(1/2)+\varepsilon}(\Omega)$ and 
\begin{equation}\label{uhgG6ggG32}
\gamma_D\widetilde{v}=\gamma_D\widetilde{\eta}-\gamma_D\widetilde{h}
=\gamma_D\widetilde{\eta}-\widetilde{h}\big|^{\kappa-{\rm n.t.}}_{\partial\Omega}
=\gamma_D\widetilde{\eta}-{\rm Tr}\,\widetilde{\eta}=0,
\end{equation}
by Lemma~\ref{Hgg-uh} (applied to $\widetilde{h}$), the boundary condition in 
\eqref{eqn:bvp-REG.222-Mi.TL}, and the compatibility of ${\rm Tr}$ with $\gamma_D$ described 
in Proposition~\ref{P-New-Trace}. Next, following the same procedure as above that has led 
to the definition in \eqref{eq:ee3ee-Mi}, one observes that the function $\widetilde{w}$ now 
exhibits better regularity on the Sobolev scale, namely $\widetilde{w}\in H^{(1/2)+\delta}(\Omega)$, 
where $\delta:=\min\{\varepsilon,s-(1/2)\}>0$. Granted this fact and  \eqref{eq:ASSD-Mi}, 
one then invokes \eqref{eq:MM3bis} for $\widetilde{w}$ to conclude that 
\begin{equation}\label{eq:ee3ee.2-Mi}
\gamma_D\widetilde{w}=\widetilde{w}\big|^{\kappa-{\rm n.t.}}_{\partial\Omega}.
\end{equation}
Since by design $u=\widetilde{w}+\widetilde{v}$ in $\Omega$, it follows from \eqref{uhgG6ggG32}
and \eqref{eq:ee3ee.2-Mi} that $\gamma_D u$ considered in the sense of \eqref{eqn:gammaDs.1} 
is consistent with our definition in \eqref{eq:ee3ee-Mi}. \\[1mm] 
We now address the claims made in the itemized portion of the statement of the theorem. \\[1mm] 
{\it Proof of $(i)$.} 
Fix $s\in\big[\tfrac12,\tfrac32\big]$. Since, obviously, $\big\{u\in H^s(\Omega)\,\big|\,\Delta u=0\big\}$ 
is a subspace of $\big\{u\in H^s(\Omega)\,\big|\,\Delta u\in B^{2,1}_{s-2}(\Omega)\big\}$,
the same operator $\Upsilon_D$ as in \eqref{2.88X-NN-44aaa}--\eqref{eqn:bvp-REG} may be employed  
as a right-inverse for $\gamma^{\#}_D$ (since the compatibility of the 
present sharp trace operator $\gamma^{\#}_D$ with $\gamma_D$ from \eqref{eqn:gammaDs.2aux}, has 
already been established). As a corollary, this also proves that the sharp Dirichlet boundary trace 
operator $\gamma^{\#}_D$ is surjective in the context of \eqref{eqn:gammaDs.2aux-Mi}. \\[1mm] 
{\it Proof of $(ii)$.}
Fix a function $u\in H^{1/2}(\Omega)$ satisfying $\Delta u\in B^{2,1}_{-3/2}(\Omega)$ 
and such that $u\big|^{\kappa-{\rm n.t.}}_{\partial\Omega}$ exists at $\sigma$-a.e.~point on $\partial\Omega$.
Since by \eqref{eq:ASSD-Mi} one knows that $\widetilde{w}\big|^{\kappa-{\rm n.t.}}_{\partial\Omega}$ 
also exists at $\sigma$-a.e.~point on $\partial\Omega$, one concludes from \eqref{eHg-Mi.2}
that $\widetilde{v}\big|^{\kappa-{\rm n.t.}}_{\partial\Omega}$ exists at $\sigma$-a.e.~point on $\partial\Omega$.
Together with \eqref{eHg-Mi.1} and the fact that, by design, $\widetilde{h}\big|^{\kappa-{\rm n.t.}}_{\partial\Omega}$ 
does exist at $\sigma$-a.e.~point on $\partial\Omega$, this implies that 
$\widetilde{\eta}\big|^{\kappa-{\rm n.t.}}_{\partial\Omega}$ exists at $\sigma$-a.e.~point on $\partial\Omega$.
Having established this fact, the compatibility of ${\rm Tr}$ with the nontangential boundary trace  
guaranteed by Proposition~\ref{P-New-Trace} then forces
\begin{equation}\label{eqn:bvp-REG.222-Mi.TL.Ba}
\widetilde{\eta}\big|^{\kappa-{\rm n.t.}}_{\partial\Omega}={\rm Tr}\,\widetilde{\eta}
\,\text{ $\sigma$-a.e.~on $\partial\Omega$}.
\end{equation}
Consequently, on account of \eqref{eqn:bvp-REG.222-Mi.TL.Ba} and the boundary 
condition in \eqref{eqn:bvp-REG.222-Mi.TL}, one can write 
\begin{align}\label{eq:ee3ee.2bgfd-Mi}
u\big|^{\kappa-{\rm n.t.}}_{\partial\Omega} &=\widetilde{w}\big|^{\kappa-{\rm n.t.}}_{\partial\Omega}
+\widetilde{v}\big|^{\kappa-{\rm n.t.}}_{\partial\Omega}
\nonumber\\[2pt]
&=\widetilde{w}\big|^{\kappa-{\rm n.t.}}_{\partial\Omega}
+\widetilde{\eta}\big|^{\kappa-{\rm n.t.}}_{\partial\Omega}-\widetilde{h}\big|^{\kappa-{\rm n.t.}}_{\partial\Omega}
\nonumber\\[2pt]
&=\widetilde{w}\big|^{\kappa-{\rm n.t.}}_{\partial\Omega}
+{\rm Tr}\,\widetilde{\eta}-\widetilde{h}\big|^{\kappa-{\rm n.t.}}_{\partial\Omega}
\nonumber\\[2pt]
&=\gamma^{\#}_D u,
\end{align}
as wanted. To complete the proof of \eqref{eqn:gammaDs.2auxBBB-Mi} there remains to 
observe that when $s\in\big(\tfrac12,\tfrac32\big]$ the desired compatibility property 
follows from the manner in which the sharp Dirichlet trace has been defined in 
\eqref{eqn:gammaDs.2aux-Mi} and Lemma~\ref{Hgg-uh}. \\[1mm] 
{\it Proof of $(iii)$.}
That $\gamma^{\#}_D$ in \eqref{eqn:gammaDs.2aux-Mi} is the unique extension by continuity and density
of the mapping $C^\infty(\overline{\Omega})\ni f\mapsto f\big|_{\partial\Omega}$
follows from Lemma~\ref{Dense-LLLe-BBB} and \eqref{eqn:gammaDs.2auxBBB-Mi}. \\[1mm] 
{\it Proof of $(iv)$.}
Pick $u\in H^s(\Omega)$ satisfying $\Delta u\in B^{2,1}_{s-2}(\Omega)$ for some 
$s\in\big[\tfrac{1}{2},\tfrac{3}{2}\big]$, along with some $\Phi\in C^\infty(\overline\Omega)$. 
By the density result proved in Lemma~\ref{Dense-LLLe-BBB} there exists a sequence 
$\{u_j\}_{j\in{\mathbb{N}}}\subset C^\infty(\overline{\Omega})$ with the property that
\begin{equation}\label{ajtrtv5-Mi-678}
u_j\to u\,\text{ in }\,H^s(\Omega)\,\text{ and }\,\Delta u_j\to\Delta u
\,\text{ in }\,B^{2,1}_{s-2}(\Omega),\,\text{ as }\,j\to\infty.
\end{equation}
In particular, $\Phi u_j\to\Phi u$ in $H^s(\Omega)$ and $\Delta(\Phi u_j)\to\Delta(\Phi u)$ in $B^{2,1}_{s-2}(\Omega)$ as $j\to\infty$.
On account of the continuity of the sharp Dirichlet trace operator, this permits us to write,
in the sense of $H^{s-(1/2)}(\partial\Omega)$, 
\begin{align}\label{ajtrtv7-Mi}
\gamma^{\#}_D(\Phi u) &=\lim_{j\to\infty}\gamma^{\#}_D(\Phi u_j)
=\lim_{j\to\infty}(\Phi u_j)\big|_{\partial\Omega}
\nonumber\\[2pt]
&=\lim_{j\to\infty}\big(\Phi\big|_{\partial\Omega}\big)\gamma^{\#}_Du_j
=\big(\Phi\big|_{\partial\Omega}\big)\gamma^{\#}_Du,
\end{align}
as wanted. 
\\[1mm] 
{\it Proof of $(v)$.}
Suppose that $s\in\big[\tfrac{1}{2},\tfrac{3}{2}\big]$ and $u\in H^s(\Omega)$ is such that 
$\Delta u\in B^{2,1}_{s-2}(\Omega)$. Since $H^s(\Omega)=B^{2,2}_s(\Omega)=F^{2,2}_s(\Omega)$ 
(with equivalent norms) and 
\begin{align}\label{eq:EFFa.annb.WACO.BB}
B^{2,1}_{s-2}(\Omega)\hookrightarrow F^{2,1}_{s-2}(\Omega) 
\end{align}
(cf.~\eqref{Inc-P3}, \eqref{LLa-45.222.WACO}, and \eqref{LLa-45}), it follows that 
$u\in F^{2,2}_s(\Omega)$, $\Delta u\in F^{2,1}_{s-2}(\Omega)$, and there exists some constant $C\in(0,\infty)$, independent of $u$, such that 
\begin{equation}\label{eq:H.WACO.1.ii.BB}
\|u\|_{F^{2,2}_s(\Omega)}\leq C\|u\|_{H^s(\Omega)},\quad
\|\Delta u\|_{F^{2,1}_{s-2}(\Omega)}\leq C\|\Delta u\|_{B^{2,1}_{s-2}(\Omega)}.
\end{equation}
With these in hand, one can invoke Proposition~\ref{Dense-LLLe-BBB.WACO} to conclude that  
$u$ belongs to $F^{2,1}_s(\Omega)$ and that 
\begin{equation}\label{eq:H.WACO.1.aaa.BB}
\|u\|_{F^{2,1}_s(\Omega)}\leq C\big(\|u\|_{H^s(\Omega)}+\|\Delta u\|_{B^{2,1}_{s-2}(\Omega)}\big).
\end{equation}
Hence, the space on the left-hand side of \eqref{eqn:gammaDs.2aux-Mi}, equipped with the natural graph norm, 
embeds continuously into $F^{2,1}_s(\Omega)$. This implies that the embeddings in \eqref{eqn:gammaDs.2aux.WACO.BBBB} are 
well defined mappings. The fact that said embeddings are strict is then justified much as in the case of \eqref{eqn:gammaDs.2aux.WACO}.
\\[1mm] 
{\it Proof of $(vi)$.}
Consider a function $u\in H^{3/2}(\Omega)$ with $\Delta u\in B^{2,1}_{-1/2}(\Omega)$
and fix some arbitrary index $j\in\{1,\dots,n\}$. Based on the assumptions made and \eqref{eq:DDEEn.3} one concludes that $\partial_j u\in H^{1/2}(\Omega)$ and $\|\partial_j u\|_{H^{1/2}(\Omega)}\leq C\|u\|_{H^{3/2}(\Omega)}$ for some constant $C\in(0,\infty)$ independent of $u$. Due to \eqref{q-TWU.ted} and the assumptions made, one also has 
\begin{equation}\label{eq:H.WACO.1.aaa.TX.1.BE}
\begin{array}{c}
\Delta(\partial_j u)=\partial_j(\Delta u)\in B^{2,1}_{-3/2}(\Omega)\,\text{ and}
\\[6pt]
\|\Delta(\partial_j u)\|_{B^{2,1}_{-3/2}(\Omega)}\leq C\|\Delta u\|_{B^{2,1}_{-1/2}(\Omega)}, 
\end{array}
\end{equation}
with $C\in(0,\infty)$ independent of $u$. Upon recalling that \eqref{eqn:gammaDs.2aux-Mi} is well defined 
and bounded when $s=\tfrac{1}{2}$, these properties guarantee that $\gamma^{\#}_D(\partial_j u)$ belongs to $L^2(\partial\Omega)$ and 
\begin{align}\label{eq:H.WACO.1.aaa.TX.2.BE}
\|\gamma^{\#}_D(\partial_j u)\|_{L^2(\partial\Omega)} &\leq C\big(\|\partial_j u\|_{H^{1/2}(\Omega)}
+\|\Delta(\partial_j u)\|_{B^{2,1}_{-3/2}(\Omega)}\big)
\nonumber\\[2pt]
&\leq C\big(\|u\|_{H^{3/2}(\Omega)}+\|\Delta u\|_{B^{2,1}_{-1/2}(\Omega)}\big).
\end{align}
Hence, the operator \eqref{eqn:gammaDs.2aux.INTRO.WACO.d} is well defined, linear, and bounded.
The proof of Theorem~\ref{YTfdf-T-Mi} is therefore complete. 
\end{proof}
%%%%%%%%

%%%%%%%%%%%%%%%%%%%%%%%%%%%%%%%
%%%%%%%%%%%%%%%%%%%%%%%%%%%%%%%
\section{Divergence Theorems with Sobolev Traces} 
\label{s4}
%%%%%%%%%%%%%%%%%%%%%%%%%%%%%%%
%%%%%%%%%%%%%%%%%%%%%%%%%%%%%%%

The goal in this section is to test the versatility of the brand of the Dirichlet boundary 
trace developed in Theorem~\ref{YTfdf-T} in the context of the divergence theorem.  

A first result of this nature is presented in Theorem~\ref{Ygav-75}.
As a preamble, we first deal with the weaker result below. 

%%%%%%
\begin{lemma}\label{YgLLam}
Let $\Omega\subset\bbR^n$ be a bounded Lipschitz domain, and fix some open neighborhood 
${\mathcal{O}}$ of $\overline{\Omega}$, along with some number $\varepsilon>0$. In addition, 
assume that the vector field $\vec{G}\in\big[H^{(1/2)+\varepsilon}_{\rm loc}({\mathcal{O}})\big]^n$
satisfies ${\rm div}\,\vec{G}\in L^1_{\rm loc}({\mathcal{O}})$. Then, if $\nu$ and $\sigma$ are,
respectively, the outward unit normal and surface measure to $\partial\Omega$, it follows that
\begin{equation}\label{eq:GCxa.22}
\int_{\Omega}{\rm div}\,\vec{G}\,d^n x
=\int_{\partial\Omega}\nu\cdot\gamma_D\vec{G}\,d^{n-1}\sigma,
\end{equation}
where the Dirichlet boundary trace operator acts componentwise.
\end{lemma}
%%%%%%%
\begin{proof}
Consider a function $\eta\in C^\infty_0({\mathbb{R}}^n)$ 
such that $\eta=1$ on $B(0,1)$, $\eta=0$ outside $B(0,2)$,
$\int_{{\mathbb{R}}^n}\eta(x)\,d^n x=1$ and, for each $t>0$, set $\eta_t(x):=t^{-n}\eta(x/t)$ 
for $x\in{\mathbb{R}}^n$. Next, fix a cutoff function $\zeta\in C^\infty_0({\mathcal{O}})$ 
with the property that $\zeta=1$ near $\overline{\Omega}$ and, for each $t>0$, consider the operator 
\begin{equation}\label{eq:Jba.1}
T_tu:=\big[\eta_t\ast(\zeta u)\big]\big|_{\Omega}\in C^\infty(\overline{\Omega})
\,\text{ for }\,u\in L^1_{\rm loc}({\mathcal{O}}).
\end{equation}
Then for each $u\in L^2_{\rm loc}({\mathcal{O}})$ one has  
$T_tu\to u\big|_{\Omega}$ as $t\to 0_{+}$ in $L^2(\Omega)$.
Moreover, if $u\in H^k_{\rm loc}({\mathcal{O}})$ for some $k\in{\mathbb{N}}$ and if $\alpha$ 
is a multi-index of length at most $k$, then 
\begin{equation}\label{eq:JBa6gg}
\partial^\alpha(T_tu)=\big[\eta_t\ast\big(\partial^\alpha(\zeta u)\big)\big]\big|_{\Omega}
\to\partial^\alpha u\big|_{\Omega}\,\text{ as }\,t\to 0_{+}\,\text{ in }\,L^2(\Omega).
\end{equation}
Next, consider an arbitrary number $s>0$ and pick $k\in{\mathbb{N}}$, $k>s$, and $\theta\in(0,1)$ 
such that $s=\theta k$. Then for every $u\in C^\infty({\mathcal{O}})$, the interpolation 
inequality 
\begin{equation}\label{Jvav-utrf}
\big\|T_tu-u\big|_{\Omega}\big\|_{H^s(\Omega)}
\leq\big\|T_tu-u\big|_{\Omega}\big\|_{H^k(\Omega)}^\theta
\big\|T_tu-u\big|_{\Omega}\big\|_{L^2(\Omega)}^{1-\theta}
\end{equation}
proves that 
\begin{equation}\label{eq:Hac}
T_tu\to u\big|_{\Omega}\,\text{ as }\,t\to 0_{+}\,\text{ in }\,H^s(\Omega),
\quad\forall\,u\in C^\infty({\mathcal{O}}).
\end{equation}
To proceed, select a bounded Lipschitz domain $\widetilde{\Omega}$ whose closure is contained 
in ${\mathcal{O}}$ and such that ${\rm supp}\,(\zeta)\subset\widetilde{\Omega}$. 
When viewed as an operator acting from $H^k(\widetilde{\Omega})$, $k\in{\mathbb{N}}\cup\{0\}$, 
via the same recipe as in \eqref{eq:Jba.1}, the same type of argument as in \eqref{eq:JBa6gg}
shows that $T_t$ is bounded into $H^k(\Omega)$, uniformly in $t>0$. Hence, by interpolation, 
$T_t$ is bounded from $H^s\big(\widetilde{\Omega}\big)$ into $H^s(\Omega)$ for each $s>0$, 
uniformly in $t>0$. 

At this point, consider an arbitrary $u\in H^s_{\rm loc}({\mathcal{O}})$ and pick some arbitrary
$\delta>0$. Then there exists $v\in C^\infty({\mathcal{O}})$ such that 
$\big\|u|_{\widetilde{\Omega}}-v|_{\widetilde{\Omega}}\big\|_{H^s(\widetilde{\Omega})}<\delta$. Then 
\begin{align}\label{eNba}
\big\|T_t u-u\big|_{\Omega}\big\|_{H^s(\Omega)}
& \leq\|T_t(u-v)\|_{H^s(\Omega)}
+\big\|T_t v-v\big|_{\Omega}\big\|_{H^s(\Omega)}+\big\|u|_{\Omega}-v|_{\Omega}\big\|_{H^s(\Omega)}
\nonumber\\[2pt]
& \leq C\big\|u|_{\widetilde{\Omega}}-v|_{\widetilde{\Omega}}\big\|_{H^s(\widetilde{\Omega})}
+\big\|T_t v-v\big|_{\Omega}\big\|_{H^s(\Omega)}+\delta
\nonumber\\[2pt]
& \leq C\delta+\big\|T_t v-v\big|_{\Omega}\big\|_{H^s(\Omega)}.
\end{align}
Together with \eqref{eq:Hac} this ultimately proves that 
\begin{equation}\label{eq:Hac.2}
T_tu\to u\big|_{\Omega}\,\text{ as }\,t\to 0_{+}\,\text{ in }\, H^s(\Omega),
\,\text{ for every }\,u\in H^s_{\rm loc}({\mathcal{O}}).
\end{equation}

Next, we extend the definition of $T_t$ by allowing it to act componentwise (as in 
\eqref{eq:Jba.1}) on vector fields. In this regard, we note that if 
$\vec{F}\in\big[L^1_{\rm loc}({\mathcal{O}})\big]^n$ is such that 
${\rm div}\vec{F}\in L^1_{\rm loc}({\mathcal{O}})$ then 
\begin{equation}\label{eq:Jba.7}
{\rm div}(T_t\vec{F})=\big[\eta_t\ast\big({\rm div}(\zeta\vec{F})\big)\big]\big|_{\Omega}
=\big[\eta_t\ast\big(\zeta{\rm div}\vec{F}\big)\big]\big|_{\Omega}
+\big[\eta_t\ast\big(\nabla\zeta\cdot\vec{F}\big)\big]\big|_{\Omega}
\end{equation}
hence, in this case, 
\begin{equation}\label{eq:BBa27}
{\rm div}(T_t\vec{F})\to\big({\rm div}\vec{F}\big)\big|_{\Omega}\,\text{ in }\, L^1(\Omega)
\,\text{ as }\, t\to 0_{+}.
\end{equation}

Given a vector field $\vec{G}\in\big[H^{(1/2)+\varepsilon}_{\rm loc}({\mathcal{O}})\big]^n$
with ${\rm div}\,\vec{G}\in L^1_{\rm loc}({\mathcal{O}})$, one can write 
\begin{align}\label{eq:GCxa.25yy}
\int_{\Omega}{\rm div}\,\vec{G}\,d^n x
&=\lim_{t\to 0_{+}}\int_{\Omega}{\rm div}\,(T_t\vec{G})\,d^n x
=\lim_{t\to 0_{+}}\int_{\partial\Omega}\nu\cdot\gamma_D\big(T_t\vec{G}\big)\,d^{n-1}\sigma
\nonumber\\[2pt]
&=\int_{\partial\Omega}\nu\cdot\gamma_D\vec{G}\,d^{n-1}\sigma.
\end{align}
Above, we used \eqref{eq:BBa27} in the first equality. The second 
equality is based on the divergence theorem for the vector field 
$T_t\vec{G}\in\big[ C^\infty(\overline{\Omega})\big]^n$. The final equality 
relies on the fact that \eqref{eq:Hac.2} implies
\begin{equation}\label{eq:Hac.2Bn}
T_t\vec{G}\to\vec{G}\big|_{\Omega}\,\text{ as }\, 
t\to 0_{+}\,\text{ in }\,\big[H^{(1/2)+\varepsilon}(\Omega)\big]^n, 
\end{equation}
hence, by the continuity of the Dirichlet trace,  
\begin{equation}\label{eq:Hac.2Bn.AA}
\gamma_D\big(T_t\vec{G}\big)\to\gamma_D\vec{G}\,\text{ as }\,t\to 0_{+}
\,\text{ in }\,\big[H^{\varepsilon}(\partial\Omega)\big]^n
\hookrightarrow\big[L^1(\partial\Omega)\big]^n.
\end{equation}
This finishes the proof of \eqref{eq:GCxa.22}.
\end{proof}
%%%%%%%

We are now ready to discuss a version of the divergence theorem which makes use
of the brand of Dirichlet boundary trace from Theorem~\ref{YTfdf-T} (when $s= 1/2$).
In turn, results of this type are going to be instrumental in the proof of Theorem~\ref{YTfdf.NNN.2-Main}, 
dealing with the Neumann boundary trace operator. 

%%%%%%%
\begin{theorem}\label{Ygav-75}
Let $\Omega\subset\bbR^n$ be a bounded Lipschitz domain, with surface 
measure $\sigma$ and outward unit normal $\nu$. Then for every vector field 
$\vec{F}\in\big[H^{1/2}(\Omega)\big]^n$ with ${\rm div}\vec{F}\in L^1(\Omega)$ 
and satisfying $\Delta\vec{F}\in\big[H^{-(3/2)+\varepsilon}(\Omega)\big]^n$ 
for some $\varepsilon>0$ one has
\begin{equation}\label{eq:GCxa}
\int_{\Omega}{\rm div}\vec{F}\,d^n x
=\int_{\partial\Omega}\nu\cdot\gamma_D\vec{F}\,d^{n-1}\sigma,
\end{equation}
where the action of $\gamma_D$ on $\vec{F}$ is considered componentwise, in the sense
of \eqref{eqn:gammaDs.2aux} with $s=1/2$ {\rm (}which places
$\gamma_D\vec{F}$ in $\big[L^2(\partial\Omega)\big]^n${\rm )}.

As a corollary, \eqref{eq:GCxa} holds for every vector field 
$\vec{F}\in\big[H^{(1/2)+\varepsilon}(\Omega)\big]^n$ for some $\varepsilon>0$
with the property that ${\rm div}\vec{F}\in L^1(\Omega)$ {\rm (}hence, in particular, 
for every vector field $\vec{F}\in\big[H^1(\Omega)\big]^n${\rm )}.
\end{theorem}
%%%%%%%
\begin{proof}
To get started, one invokes \cite[Theorem~0.5(b), pp.~164--165]{JK95} in order to solve the boundary value problem 
\begin{equation}\label{eqn:bNN.1}
\begin{cases}
\Delta\vec{G}=\Delta\vec{F}\,\text{ in $\Omega,\quad 
\vec{G}\in\big[H^{(1/2)+\varepsilon}(\Omega)\big]^n$,}     
\\[2pt]  
\gamma_D\vec{G}=0\,\text{ on $\partial\Omega$.}
\end{cases}   
\end{equation}
Next, consider $\vec{h}:=\vec{F}-\vec{G}$ in $\Omega$. It follows that $\Delta\vec{h}=0$ 
in $\Omega$, thus $\vec{h}\in\big[ C^\infty(\Omega)\big]^n$. In particular, 
${\rm div}\,\vec{G}={\rm div}\vec{F}-{\rm div}\vec{h}\in L^1_{\rm loc}(\Omega)$.
Moreover, $\vec{h}\in\big[H^{1/2}(\Omega)\big]^n$, hence by \eqref{eq:MM3} and 
\eqref{eq:MM4} one concludes that ${\mathcal{N}}_\kappa\vec{h}\in L^2(\partial\Omega)$ 
and $\vec{h}\big|^{\kappa-{\rm n.t.}}_{\partial\Omega}$ exists $\sigma$-a.e.~on $\partial\Omega$,
and belongs to $\big[L^2(\partial\Omega)\big]^n$. By the last condition in \eqref{eqn:bNN.1}, this forces 
$\gamma_D\vec{F}=\gamma_D\vec{h}=\vec{h}\big|^{\kappa-{\rm n.t.}}_{\partial\Omega}$, where the last 
equality is a consequence of item $(ii)$ in Theorem~\ref{YTfdf-T} (cf.~\eqref{eqn:gammaDs.2auxBBB}).

To proceed, we consider an approximating family $\Omega_\ell\nearrow\Omega$ as $\ell\to\infty$ 
of the sort described in Lemma~\ref{OM-OM} and recall that 
$\nu^\ell\circ\Lambda_\ell\to\nu$ as $\ell\to\infty$ both pointwise $\sigma$-a.e.~on 
$\partial\Omega$ and in $\big[L^2(\partial\Omega)\big]^n$. Moreover, the properties
of the homeomorphisms $\Lambda_\ell$ allow us to conclude that 
$\big(\vec{h}\big|_{\partial\Omega_\ell}\big)\circ\Lambda_\ell\to\vec{h}\big|^{\kappa-{\rm n.t.}}_{\partial\Omega}$ 
as $\ell\to\infty$ both pointwise and in $\big[L^2(\partial\Omega)\big]^n$, by Lebesgue's dominated 
convergence theorem (with uniform domination provided by 
${\mathcal{N}}_\kappa\vec{h}\in L^2(\partial\Omega)$). 
Finally, one recalls that the $\omega_\ell$'s appearing in the change of variable formula 
\eqref{eQQ-14} are uniformly bounded and converge to $1$ as $\ell\to\infty$ pointwise 
$\sigma$-a.e.~on $\partial\Omega$. Given these facts and keeping in mind that $\vec{h}\in\big[C^\infty(\Omega)\big]^n$, one computes  
\begin{align}\label{Utfv}
& \lim_{\ell\to\infty}\int_{\partial\Omega_\ell}\nu^\ell\cdot 
\Big(\vec{h}\big|_{\partial\Omega_\ell}\Big)\,d^{n-1}\sigma_\ell
\nonumber\\[2pt]
& \quad=\lim_{\ell\to\infty}\int_{\partial\Omega}(\nu^\ell\circ\Lambda_\ell)\cdot
\Big(\vec{h}\big|_{\partial\Omega_\ell}\Big)\circ\Lambda_\ell\,\omega_\ell\,d^{n-1}\sigma
\nonumber\\[2pt]
& \quad=\int_{\partial\Omega}\nu\cdot\Big(\vec{h}\big|^{\kappa-{\rm n.t.}}_{\partial\Omega}\Big)\,d^{n-1}\sigma
=\int_{\partial\Omega}\nu\cdot\gamma_D\vec{F}\,d^{n-1}\sigma.
\end{align}
On the other hand, applying the divergence theorem in each Lipschitz domain $\Omega_\ell$ 
for the vector field $\vec{h}\big|_{\Omega_\ell}\in\big[C^\infty(\overline{\Omega_\ell})\big]^n$ 
(cf.\  Theorem~\ref{banff-3}), relying on Lebesgue's dominated convergence theorem, and invoking Lemma~\ref{YgLLam}, yields
\begin{align}\label{Utfv.2}
& \lim_{\ell\to\infty}\int_{\partial\Omega_\ell}\nu^\ell\cdot  
\Big(\vec{h}\big|_{\partial\Omega_\ell}\Big)\,d^{n-1}\sigma_\ell
\nonumber\\[2pt]
& \quad=\lim_{\ell\to\infty}\int_{\Omega_\ell}{\rm div}\,\vec{h}\,d^n x
\nonumber\\[2pt]
& \quad=\lim_{\ell\to\infty}\int_{\Omega_\ell}{\rm div}\vec{F}\,d^n x
-\lim_{\ell\to\infty}\int_{\Omega_\ell}{\rm div}\,\vec{G}\,d^n x
\nonumber\\[2pt]
& \quad=\int_{\Omega}{\rm div}\vec{F}\,d^n x-\lim_{\ell\to\infty}\int_{\partial\Omega_\ell}
\nu^\ell\cdot\gamma_{\ell,D}\big(\vec{G}\big|_{\Omega_\ell}\big)\,d^{n-1}\sigma_\ell, 
\end{align}
where, for each $\ell\in{\mathbb{N}}$, we denoted by $\gamma_{\ell,D}$ the Dirichlet boundary 
trace operator associated with the Lipschitz domain $\Omega_\ell$. The next step in the proof is to 
pick a number $\delta\in\big(0,\min\{\tfrac{1}{2},\varepsilon\}\big)$ then estimate 
\begin{align}\label{jab0uhb}
\begin{split} 
\bigg|\int_{\partial\Omega_\ell}\nu^\ell\cdot\gamma_{\ell,D}
\big(\vec{G}\big|_{\Omega_\ell}\big)\,d^{n-1}\sigma_\ell\bigg|
&\leq\big\|\gamma_{\ell,D}\big(\vec{G}\big|_{\Omega_\ell}\big)\big\|_{[L^1(\partial\Omega_\ell)]^n}    
\\[2pt] 
&\leq C\big\|\gamma_{\ell,D}\big(\vec{G}\big|_{\Omega_\ell}\big)\big\|_{[H^\delta(\partial\Omega_\ell)]^n}
\end{split} 
\end{align}
for some constant $C\in(0,\infty)$, independent of $\ell\in{\mathbb{N}}$. Since by \eqref{incl-Yb.EE} 
and \eqref{eqn:bNN.1} one has $\vec{G}\in\big[\accentset{\circ}{H}^{(1/2)+\delta}(\Omega)\big]^n$, 
it follows from Lemma~\ref{Tgav9jy} (used with $s=\tfrac{1}{2}+\delta\in(\tfrac{1}{2},1)$) that 
\begin{align}\label{jab0uhb.2}
\lim_{\ell\to\infty}\big\|\gamma_{\ell,D}
\big(\vec{G}\big|_{\Omega_\ell}\big)\big\|_{[H^\delta(\partial\Omega_\ell)]^n}=0.
\end{align}
At this stage, \eqref{eq:GCxa} follows from \eqref{Utfv}--\eqref{jab0uhb.2}.
\end{proof}
%%%%%%%

The technical result contained in our next lemma is going to be useful shortly, in the 
proof of Theorem~\ref{Ygav-75-BIS} below.

%%%%%%%
\begin{lemma}\label{Lapp-L1}
Let $\Omega\subset\bbR^n$ be a bounded Lipschitz domain, and consider an approximating 
family $\Omega_\ell\nearrow\Omega$ as $\ell\to\infty$ as described in Lemma~\ref{OM-OM}. Assume that 
$f\in L^1_{\rm loc}(\Omega)\cap H^{-(1/2)+\varepsilon}(\Omega)$ for some $\varepsilon\in(0,1)$. Then 
\begin{equation}\label{Grdd-1}
\lim_{\ell\to\infty}\int_{\Omega_\ell}f(x)\,d^nx=
{}_{H^{(1/2)-\varepsilon}(\Omega)}\big\langle{\bf 1},f\big\rangle_{H^{-(1/2)+\varepsilon}(\Omega)},
\end{equation}
where ${\bf 1}$ denotes the constant function, identically equal to $1$, in $\Omega$. 
\end{lemma}
%%%%%%%
\begin{proof}
For each $\ell\in{\mathbb{N}}$ denote by $\chi_{\Omega_\ell}$ the characteristic function 
of $\Omega_\ell$. That is, $\chi_{\Omega_\ell}:{\mathbb{R}}^n\to{\mathbb{R}}$ given by 
$\chi_{\Omega_\ell}(x)=1$ if $x\in\Omega_\ell$, and $\chi_{\Omega_\ell}(x)=0$ if 
$x\in{\mathbb{R}}^n\backslash\Omega_\ell$. By \cite[Lemma~4, p.~52]{RS96} and item 
(4) in the proposition from \cite[pp.~29--30]{RS96}, for every $\ell\in{\mathbb{N}}$ 
one has (with $B^{p,q}_s({\mathbb{R}}^n)$ denoting the standard scale of Besov spaces 
in ${\mathbb{R}}^n$ defined in \eqref{besov-jussi}--\eqref{besov-jussi2})
\begin{equation}\label{CHI-iuh.1}
\chi_{\Omega_\ell}\in B^{2,\infty}_{1/2}({\mathbb{R}}^n)\hookrightarrow 
B^{2,2}_{(1/2)-\varepsilon}({\mathbb{R}}^n)=H^{(1/2)-\varepsilon}({\mathbb{R}}^n)
\end{equation}
and, in fact, 
\begin{equation}\label{CHI-iuh.2}
\sup_{\ell\in{\mathbb{N}}}\|\chi_{\Omega_\ell}\|_{H^{(1/2)-\varepsilon}({\mathbb{R}}^n)}<\infty.
\end{equation}
Consequently, if one considers 
${\bf 1}_\ell:=\chi_{\Omega_\ell}\big|_{\Omega}$ 
for each $\ell\in{\mathbb{N}}$, it follows that
\begin{equation}\label{CHI-iuh.3}
{\bf 1}_\ell\in H^{(1/2)-\varepsilon}(\Omega)\,\text{ for every }\,\ell\in{\mathbb{N}},
\,\text{ and }\,\sup_{\ell\in{\mathbb{N}}}\|{\bf 1}_\ell\|_{H^{(1/2)-\varepsilon}(\Omega)}<\infty.
\end{equation}
We claim that actually 
\begin{equation}\label{CHI-iuh.4}
{\bf 1}_\ell\to{\bf 1}\,\text{ in }\,H^{(1/2)-\varepsilon}(\Omega)\,\text{ as }\,\ell\to\infty.
\end{equation}
Indeed, since 
\begin{equation}\label{CHI-iuh.6}
C^\infty_0(\Omega)\,\text{ is dense in }\,H^{-(1/2)+\varepsilon}(\Omega),\quad\forall\,\varepsilon\in(0,1),
\end{equation}
the claim in \eqref{CHI-iuh.4} follows with the help of \eqref{CHI-iuh.3}, upon noting that 
for each function $\varphi\in C^\infty_0(\Omega)$ one has
\begin{align}\label{CHI-iuh.7}
& \lim_{\ell\to\infty}{}_{H^{(1/2)-\varepsilon}(\Omega)}\big\langle{\bf 1}_\ell,\varphi
\big\rangle_{H^{-(1/2)+\varepsilon}(\Omega)}
=\lim_{\ell\to\infty}{}_{{\mathcal{D}}'(\Omega)}\big\langle{\bf 1}_\ell,\varphi
\big\rangle_{{\mathcal{D}}(\Omega)}
\nonumber\\[2pt]
& \quad=\lim_{\ell\to\infty}\int_{\Omega_\ell}\varphi(x)\,d^n x=\int_{\Omega}\varphi(x)\,d^n x
\nonumber\\[2pt]
& \quad={}_{H^{(1/2)-\varepsilon}(\Omega)}\big\langle{\bf 1},\varphi\big\rangle_{H^{-(1/2)+\varepsilon}(\Omega)}.
\end{align}
Having established this fact, for every 
$f\in L^1_{\rm loc}(\Omega)\cap H^{-(1/2)+\varepsilon}(\Omega)$ 
with $\varepsilon\in(0,1)$ one then computes 
\begin{align}\label{CHI-iuh.8}
\lim_{\ell\to\infty}\int_{\Omega_\ell}f(x)\,d^nx &=
{}_{H^{(1/2)-\varepsilon}(\Omega)}\big\langle{\bf 1}_\ell,f
\big\rangle_{H^{-(1/2)+\varepsilon}(\Omega)}
\nonumber\\[2pt]
&={}_{H^{(1/2)-\varepsilon}(\Omega)}\big\langle{\bf 1},f
\big\rangle_{H^{-(1/2)+\varepsilon}(\Omega)},
\end{align}
where the first equality is a consequence of Lemma~\ref{utrCOMPA}, while the second 
one uses \eqref{CHI-iuh.4}. The desired conclusion follows.
\end{proof}
%%%%%%%

Here is a version of the divergence theorem for vector fields whose divergence 
is not necessarily an absolutely integrable function. 

%%%%%%%
\begin{theorem}\label{Ygav-75-BIS}
Suppose $\Omega\subset\bbR^n$ is a bounded Lipschitz domain, with surface measure $\sigma$ 
and outward unit normal $\nu$. Let $\vec{F}\in\big[H^{1/2}(\Omega)\big]^n$ be a vector 
field with the property that $\Delta\vec{F}\in\big[H^{-(3/2)+\varepsilon}(\Omega)\big]^n$
and ${\rm div}\vec{F}\in H^{-(1/2)+\varepsilon}(\Omega)$ for some $\varepsilon\in(0,1)$. Then 
\begin{equation}\label{eq:GCxa-BIS}
{}_{H^{(1/2)-\varepsilon}(\Omega)}\big\langle{\bf 1},{\rm div}\vec{F}
\big\rangle_{H^{-(1/2)+\varepsilon}(\Omega)}
=\int_{\partial\Omega}\nu\cdot\gamma_D\vec{F}\,d^{n-1}\sigma,
\end{equation}
where ${\bf 1}$ denotes the constant function identically to $1$ in $\Omega$, and the action 
of $\gamma_D$ on $\vec{F}$ is considered componentwise, in the sense of \eqref{eqn:gammaDs.2aux} 
with $s=1/2$ {\rm (}which places $\gamma_D\vec{F}$ in $\big[L^2(\partial\Omega)\big]^n${\rm )}.
\end{theorem}
%%%%%%%
\begin{proof}
We shall reuse part of the proof of Theorem~\ref{Ygav-75}. In particular, we 
let $\vec{G}$ solve \eqref{eqn:bNN.1} and set $\vec{h}:=\vec{F}-\vec{G}$ in $\Omega$. 
As before, this satisfies 
\begin{align}\label{NEWDVT.1}
& \vec{h}\in\big[C^\infty(\Omega)\cap H^{1/2}(\Omega)\big]^n,
\\[2pt]
& \Delta\vec{h}=0\,\text{ in }\,\Omega,\quad{\mathcal{N}}_\kappa\vec{h}\in L^2(\partial\Omega),
\label{NEWDVT.2}
\\[2pt]
& \gamma_D\vec{F}=\gamma_D\vec{h}=\vec{h}\big|^{\kappa-{\rm n.t.}}_{\partial\Omega}\in\big[L^2(\partial\Omega)\big]^n.
\label{NEWDVT.3}
\end{align}
Granted the current hypotheses, one also has  
\begin{align}\label{NEWDVT.4}
{\rm div}\,\vec{h}={\rm div}\vec{F}-{\rm div}\,\vec{G}\in L^1_{\rm loc}(\Omega)
\cap H^{-(1/2)+\varepsilon}(\Omega).
\end{align}
Since $\vec{G}\in\big[\accentset{\circ}{H}^{(1/2)+\varepsilon}(\Omega)\big]^n$, by 
\eqref{eqn:bNN.1} and \eqref{incl-Yb.EE}, it follows that there exists a sequence 
$\{\vec{G}_j\}_{j\in{\mathbb{N}}}\subset\big[C^\infty_0(\Omega)\big]^n$ with the property that
\begin{align}\label{NEWDVT.5}
\vec{G}_j\to\vec{G}\,\text{ in }\,H^{(1/2)+\varepsilon}(\Omega)\,\text{ as }\,j\to\infty.
\end{align}
As a consequence, 
\begin{align}\label{NEWDVT.6}
{\rm div}\,\vec{G}_j\to{\rm div}\,\vec{G}\,\text{ in }\,
H^{-(1/2)+\varepsilon}(\Omega)\,\text{ as }\,j\to\infty,
\end{align}
hence
\begin{align}\label{NEWDVT.7}
{}_{H^{(1/2)-\varepsilon}(\Omega)}\big\langle{\bf 1},{\rm div}\,\vec{G}
\big\rangle_{H^{-(1/2)+\varepsilon}(\Omega)}
&=\lim_{j\to\infty}{}_{H^{(1/2)-\varepsilon}(\Omega)}\big\langle{\bf 1},{\rm div}\,\vec{G}_j
\big\rangle_{H^{-(1/2)+\varepsilon}(\Omega)}
\nonumber\\[2pt]
&=\lim_{j\to\infty}\int_{\Omega}\big({\rm div}\,\vec{G}_j\big)(x)\,d^nx
\nonumber\\[2pt]
&=\lim_{j\rightarrow\infty}\int_{\partial\Omega}\nu
\cdot\big(\vec{G}_j\big|_{\partial\Omega}\big)\,d^{n-1}\sigma=0,
\end{align}
given that $\vec{G}_j\in\big[C^\infty_0(\Omega)\big]^n$ for every $j\in{\mathbb{N}}$. 
This fact and \eqref{NEWDVT.4} then imply 
\begin{align}\label{NEWDVT.8}
{}_{H^{(1/2)-\varepsilon}(\Omega)}\big\langle{\bf 1},{\rm div}\vec{F}
\big\rangle_{H^{-(1/2)+\varepsilon}(\Omega)}
={}_{H^{(1/2)-\varepsilon}(\Omega)}\big\langle{\bf 1},{\rm div}\,\vec{h}
\big\rangle_{H^{-(1/2)+\varepsilon}(\Omega)}.
\end{align}
As in the past, consider an approximating family $\Omega_\ell\nearrow\Omega$ as $\ell\to\infty$ 
(described in Lemma~\ref{OM-OM}). Then one writes  
\begin{align}\label{Utfv.2BIS}
{}_{H^{(1/2)-\varepsilon}(\Omega)}\big\langle{\bf 1},{\rm div}\,\vec{h}
\big\rangle_{H^{-(1/2)+\varepsilon}(\Omega)}
&=\lim_{\ell\to\infty}\int_{\Omega_\ell}{\rm div}\,\vec{h}\,d^n x
\nonumber\\[2pt]
&=\lim_{\ell\to\infty}\int_{\partial\Omega_\ell}\nu^\ell\cdot  
\Big(\vec{h}\big|_{\partial\Omega_\ell}\Big)\,d^{n-1}\sigma_\ell
\nonumber\\[2pt]
&=\int_{\partial\Omega}\nu\cdot\gamma_{D}\vec{F}\,d^{n-1}\sigma, 
\end{align}
where the first equality is implied by Lemma~\ref{Lapp-L1} and \eqref{NEWDVT.4}, 
the second equality is a consequence of \eqref{NEWDVT.1} and the 
divergence theorem in the Lipschitz domain $\Omega_\ell$ for the vector field 
$\vec{h}\big|_{\Omega_\ell}\in\big[C^\infty(\overline{\Omega_\ell})\big]^n$ 
(Theorem~\ref{banff-3} is more than adequate in this context), while the third 
equality is seen from \eqref{Utfv}. Formula \eqref{eq:GCxa-BIS} now follows 
by combining \eqref{NEWDVT.8} and \eqref{Utfv.2BIS}.
\end{proof}
%%%%%%%

It turns out that Theorem~\ref{Ygav-75-BIS} self-improves in the manner described below. 

%%%%%%%
\begin{corollary}\label{Ygav-75-BIS-CCC}
Assume that $\Omega\subset\bbR^n$ is a bounded Lipschitz domain with outward unit normal $\nu$, 
and fix some $\varepsilon\in(0,1)$. Let $\vec{F}\in\big[H^{1/2}(\Omega)\big]^n$ be a vector 
field with the property that $\Delta\vec{F}\in\big[H^{-(3/2)+\varepsilon}(\Omega)\big]^n$ 
and ${\rm div}\vec{F}\in H^{-(1/2)+\varepsilon}(\Omega)$. In addition, consider a scalar function
$u\in H^{(1/2)+\varepsilon}(\Omega)$. Then 
\begin{align}\label{eq:GCxa-BIS-CCC}
& \big(\gamma_Du\,,\,\nu\cdot\gamma_D\vec{F}\big)_{L^2(\partial\Omega)}
\nonumber\\[2pt]
& \quad={}_{H^{(1/2)-\varepsilon}(\Omega)}\big\langle u,{\rm div}\vec{F}
\big\rangle_{H^{-(1/2)+\varepsilon}(\Omega)}
\nonumber\\[2pt]
& \qquad +{}_{[H^{-(1/2)+\varepsilon}(\Omega)]^n}\big\langle\nabla u,\vec{F}
\big\rangle_{[H^{(1/2)-\varepsilon}(\Omega)]^n}.
\end{align}
\end{corollary}
%%%%%%%
\begin{proof}
From \eqref{uam-mi78} one infers that there exists a sequence 
$\{\Phi_j\}_{j\in{\mathbb{N}}}\subset C^\infty(\overline{\Omega})$ 
with the property that 
\begin{equation}\label{m7g4d5}
\Phi_j\to u\,\text{ in }\,H^{(1/2)+\varepsilon}(\Omega)\,\text{ as }\,j\to\infty.
\end{equation}
By virtue of \eqref{eqn:gammaDs.1} and \eqref{eq:DDEEn.3}, this implies
\begin{align}\label{m7g4d5.1}
\begin{split}
& \gamma_D\Phi_j\to\gamma_Du\,\text{ in }\,H^{\varepsilon}(\partial\Omega)
\hookrightarrow L^2(\partial\Omega)\,\text{ as }\,j\to\infty,
\\[2pt]
& \nabla\Phi_j\to\nabla u\,\text{ in }\,\big[H^{-(1/2)+\varepsilon}(\Omega)\big]^n\,\text{ as }\,j\to\infty.
\end{split}
\end{align}
In addition, by \eqref{eq:DDEj6g5}, for each $j\in{\mathbb{N}}$, the vector field 
$\overline{\Phi}_j\vec{F}$ satisfies the same properties as the original $\vec{F}$.
As such, with $\sigma$ denoting the surface measure on $\partial\Omega$, one can write, 
\begin{align}\label{eq:GCxa-Bvf}
& \big(\gamma_Du\,,\,\nu\cdot\gamma_D\vec{F}\big)_{L^2(\partial\Omega)}
=\lim_{j\to\infty}\big(\gamma_D\Phi_j\,,\,\nu\cdot\gamma_D\vec{F}\big)_{L^2(\partial\Omega)}
\nonumber\\[2pt]
& \quad=\lim_{j\to\infty}\int_{\partial\Omega}\overline{\Phi}_j\,\nu\cdot\gamma_D\vec{F}\,d^{n-1}\sigma
=\lim_{j\to\infty}\int_{\partial\Omega}\nu\cdot\gamma_D\big(\,\overline{\Phi}_j\vec{F}\,\big)\,d^{n-1}\sigma
\nonumber\\[2pt]
& \quad=\lim_{j\to\infty}
{}_{H^{(1/2)-\varepsilon}(\Omega)}\big\langle{\bf 1},{\rm div}(\overline{\Phi}_j\vec{F})
\big\rangle_{H^{-(1/2)+\varepsilon}(\Omega)}
\nonumber\\[2pt]
& \quad=\lim_{j\to\infty}
{}_{H^{(1/2)-\varepsilon}(\Omega)}\big\langle{\bf 1},\nabla\overline{\Phi}_j\cdot\vec{F}
\big\rangle_{H^{-(1/2)+\varepsilon}(\Omega)}
\nonumber\\[2pt]
& \qquad+\lim_{j\to\infty}{}_{H^{(1/2)-\varepsilon}(\Omega)}\big\langle{\bf 1},\overline{\Phi}_j{\rm div}\vec{F}
\big\rangle_{H^{-(1/2)+\varepsilon}(\Omega)}
\nonumber\\[2pt]
& \quad=\lim_{j\to\infty}{}_{[H^{-(1/2)+\varepsilon}(\Omega)]^n}\big\langle\nabla\Phi_j,\vec{F}
\big\rangle_{[H^{(1/2)-\varepsilon}(\Omega)]^n}
\nonumber\\[2pt]
& \qquad+\lim_{j\to\infty}{}_{H^{(1/2)-\varepsilon}(\Omega)}\big\langle\Phi_j,{\rm div}\vec{F}
\big\rangle_{H^{-(1/2)+\varepsilon}(\Omega)}
\nonumber\\[2pt]
& \quad={}_{[H^{-(1/2)+\varepsilon}(\Omega)]^n}\big\langle\nabla u,\vec{F}
\big\rangle_{[H^{(1/2)-\varepsilon}(\Omega)]^n}
\nonumber\\[2pt]
& \qquad+{}_{H^{(1/2)-\varepsilon}(\Omega)}\big\langle u,{\rm div}\vec{F}
\big\rangle_{H^{-(1/2)+\varepsilon}(\Omega)},
\end{align}
on account of Theorem~\ref{Ygav-75-BIS} together with \eqref{m7g4d5}, \eqref{m7g4d5.1}, as well as
\eqref{eqn:gammaDs.1-TR.2} and \eqref{tvdee}. This establishes \eqref{eq:GCxa-BIS-CCC}.
\end{proof}
%%%%%%%

It turns out that there is a more general result encompassing both Theorem~\ref{Ygav-75} 
and Theorem~\ref{Ygav-75-BIS}. Stating this requires a piece of notation, clarified below. 
Given a nonempty open set $\Omega\subseteq{\mathbb{R}}^n$ and some $s\in{\mathbb{R}}$, 
both $H^s(\Omega)$ and $L^1(\Omega)$ may be regarded as subspaces of ${\mathcal{D}}'(\Omega)$. 
In this context, it makes sense to consider their algebraic sum
\begin{align}\label{utrr-kb6}
H^s(\Omega)+L^1(\Omega):=\big\{u\in{\mathcal{D}}'(\Omega)\,\big|\,\,&\text{there exist }\,v\in H^s(\Omega)
\,\text{ and }\,w\in L^1(\Omega)
\nonumber\\[2pt]
&\text{with }\,u=v+w\,\text{ in }\,{\mathcal{D}}'(\Omega)\big\}.
\end{align}
Equipping this with the norm associating to each $u\in H^s(\Omega)+L^1(\Omega)$ the number
\begin{align}\label{utrr-kb7}
\|u\|_{H^s(\Omega)+L^1(\Omega)}:=\inf_{\substack{u=v+w\,\text{ in }\,{\mathcal{D}}'(\Omega)\\ 
v\in H^s(\Omega),\,w\in L^1(\Omega)}}\big(\|v\|_{H^s(\Omega)}+\|w\|_{L^1(\Omega)}\big),
\end{align}
turns $H^s(\Omega)+L^1(\Omega)$ into a Banach space, for which the natural inclusions 
\begin{align}\label{utrr-kb8}
\begin{split}
& H^s(\Omega)\hookrightarrow H^s(\Omega)+L^1(\Omega)\hookrightarrow{\mathcal{D}}'(\Omega),
\\[2pt]
& L^1(\Omega)\hookrightarrow H^s(\Omega)+L^1(\Omega)\hookrightarrow{\mathcal{D}}'(\Omega),
\end{split}
\end{align}
are continuous. Moreover, assuming that $\Omega$ is a bounded Lipschitz domain, it follows that
\begin{equation}\label{utrr-kb9}
C^\infty_0(\Omega)\hookrightarrow H^s(\Omega)+L^1(\Omega)\,\text{ densely, provided }\, 
s\in\big(-\tfrac12,\tfrac12\big).
\end{equation}
After this preamble, here is the general result alluded to earlier. 

%%%%%%%
\begin{theorem}\label{Ygav-75-BIS-BBB}
Let $\Omega\subset\bbR^n$ be a bounded Lipschitz domain, and suppose that 
$\vec{F}\in\big[H^{1/2}(\Omega)\big]^n$ is a vector field with the property that there 
exists $\varepsilon\in(0,1)$ such that $\Delta\vec{F}\in\big[H^{-(3/2)+\varepsilon}(\Omega)\big]^n$
and ${\rm div}\vec{F}\in H^{-(1/2)+\varepsilon}(\Omega)+L^1(\Omega)$. Then 
\begin{equation}\label{eq:GCxa-BIS-BBB}
{}_{(H^{-(1/2)+\varepsilon}(\Omega)+L^1(\Omega))^*}\big\langle{\bf 1},{\rm div}\vec{F}
\big\rangle_{H^{-(1/2)+\varepsilon}(\Omega)+L^1(\Omega)}
=\int_{\partial\Omega}\nu\cdot\gamma_D\vec{F}\,d^{n-1}\sigma,
\end{equation}
where ${\bf 1}$ denotes the constant function identically to $1$ in $\Omega$, 
and the action of $\gamma_D$ on $\vec{F}$ is considered componentwise, in the sense
of \eqref{eqn:gammaDs.2aux} with $s=1/2$ {\rm (}which places
$\gamma_D\vec{F}$ in $\big[L^2(\partial\Omega)\big]^n${\rm )}.
\end{theorem}
%%%%%%%
\begin{proof}
We shall follow the general outline of the proof of Theorem~\ref{Ygav-75-BIS}. To get started, let 
$\vec{G}$ solve \eqref{eqn:bNN.1} and set $\vec{h}:=\vec{F}-\vec{G}$ in $\Omega$. Once again, this 
satisfies \eqref{NEWDVT.1}--\eqref{NEWDVT.3}. In the present setting, in place of \eqref{NEWDVT.4} 
one has  
\begin{align}\label{NEWDVT.4bbb}
{\rm div}\,\vec{h}={\rm div}\vec{F}-{\rm div}\,\vec{G}\in L^1_{\rm loc}(\Omega)\cap 
\big(H^{-(1/2)+\varepsilon}(\Omega)+L^1(\Omega)\big).
\end{align}
Arguing as in \eqref{NEWDVT.5}--\eqref{NEWDVT.7} gives
\begin{align}\label{NEWDVT.7bbb}
{}_{(H^{-(1/2)+\varepsilon}(\Omega)+L^1(\Omega))^*}
\big\langle{\bf 1},{\rm div}\,\vec{G}\,\big\rangle_{H^{-(1/2)+\varepsilon}(\Omega)+L^1(\Omega)}
=0
\end{align}
which, in light of \eqref{NEWDVT.4bbb}, forces
\begin{align}\label{NEWDVT.8bbb}
& {}_{(H^{-(1/2)+\varepsilon}(\Omega)+L^1(\Omega))^*}\big\langle{\bf 1},{\rm div}\vec{F}\,
\big\rangle_{H^{-(1/2)+\varepsilon}(\Omega)+L^1(\Omega)}
\\[2pt] 
& \quad={}_{(H^{-(1/2)+\varepsilon}(\Omega)+L^1(\Omega))^*}\big\langle{\bf 1},{\rm div}\,\vec{h}\,
\big\rangle_{H^{-(1/2)+\varepsilon}(\Omega)+L^1(\Omega)}.
\nonumber
\end{align}
At this stage we recall the approximating family of domains, 
$\Omega_j\nearrow\Omega$ as $j\to\infty$ (cf.\ Lemma~\ref{OM-OM}). 
An inspection of the proof of Lemma~\ref{Lapp-L1} reveals that this easily extends to imply 
\begin{align}\label{Grdd-1bbb}
\begin{split}
& \displaystyle\lim_{j\to\infty}\int_{\Omega_j}f(x)\,d^n x
={}_{(H^{-(1/2)+\varepsilon}(\Omega)+L^1(\Omega))^*}\big\langle{\bf 1},f
\big\rangle_{H^{-(1/2)+\varepsilon}(\Omega)+L^1(\Omega)}
\\[2pt]
& \quad\text{for every function }\,
f\in L^1_{\rm loc}(\Omega)\cap\big(H^{-(1/2)+\varepsilon}(\Omega)+L^1(\Omega)\big).
\end{split}
\end{align}
Indeed, the key ingredients in the justification of \eqref{Grdd-1bbb} are: the density result 
recorded in \eqref{utrr-kb9}, along with the fact that if $s\in\big(-\tfrac12,\tfrac12\big)$ then, 
with $J_t$ as in \eqref{trew54.a} (cf. also \eqref{eq:12d-yw5}), 
\begin{align}\label{yt5tr.1}
& J_tu\to u\,\text{ in }\,H^s(\Omega)+L^1(\Omega)\,\text{ as }\,t\to 0_{+},
\quad\forall\,u\in H^s(\Omega)+L^1(\Omega),
\\[2pt]
& \quad \text{and }\,{\bf 1}_j\to 1\,\text{ in }\,\big(H^s(\Omega)+L^1(\Omega)\big)^*\,\text{ as }\,j\to\infty.
\label{yt5tr.2}
\end{align}

Continuing, using \eqref{Grdd-1bbb} for $f:={\rm div}\,\vec{h}$ (cf.~\eqref{NEWDVT.4bbb}) 
and then reasoning as in \eqref{Utfv.2BIS}, one arrives at the conclusion that 
\begin{align}\label{Utfv.2BISbbb}
\begin{split} 
& {}_{(H^{-(1/2)+\varepsilon}(\Omega)+L^1(\Omega))^*}\big\langle {\bf 1},{\rm div}\,\vec{h}
\big\rangle_{H^{-(1/2)+\varepsilon}(\Omega)+L^1(\Omega)}
\\[2pt]
& \quad=\lim_{j\to\infty}\int_{\Omega_j}{\rm div}\,\vec{h}\,d^n x
=\int_{\partial\Omega}\nu\cdot\gamma_{D}\vec{F}\,d^{n-1}\sigma. 
\end{split} 
\end{align}
Now \eqref{NEWDVT.8bbb} and \eqref{Utfv.2BISbbb} establish \eqref{eq:GCxa-BIS-BBB}, 
finishing the proof of the theorem.
\end{proof}
%%%%%%%

%%%%%%%%%%%%%%%%%%%%%%%%%%%%%%%%%
%%%%%%%%%%%%%%%%%%%%%%%%%%%%%%%%%
\section{A Sharp Neumann Trace Involving Sobolev Spaces} 
\label{s5}
%%%%%%%%%%%%%%%%%%%%%%%%%%%%%%%%%
%%%%%%%%%%%%%%%%%%%%%%%%%%%%%%%%%

Having dealt with the Dirichlet trace $\gamma_D$ in Section~\ref{s3}, we now turn 
our attention to the task of defining the Neumann boundary trace operator $\gamma_N$
in the class of bounded Lipschitz domains. In a first stage, we shall introduce a weak version 
$\widetilde\gamma_N$ of the aforementioned Neumann boundary trace operator,
whose definition is inspired by the ``half" Green's formula for the Laplacian.
Specifically, we make the following definition. 

%%%%%%%
\begin{definition}\label{h6r4d5UU}
Let $\Omega\subset{\mathbb{R}}^n$ be a bounded Lipschitz domain. For some fixed smoothness exponent 
$s\in\big(\tfrac12,\tfrac32\big)$, the {\tt weak Neumann trace operator} 
is considered acting in the context  
\begin{equation}\label{2.88X}
\widetilde\gamma_N:\big\{(f,F)\in H^s(\Omega)\times H^{s-2}_0(\Omega)
\,\big|\,
\Delta f=F|_{\Omega}\text{ in }\mathcal{D}'(\Omega)\big\}\rightarrow H^{s-(3/2)}(\partial\Omega).
\end{equation}
Specifically, suppose a function $f\in H^s(\Omega)$ along with a distribution 
$F\in H^{s-2}_0(\Omega)\subset H^{s-2}({\mathbb{R}}^n)$ 
satisfying $\Delta f=F|_{\Omega}$ in $\mathcal{D}'(\Omega)$ have been given. In particular, 
\eqref{eq:DDEEn.3} and \eqref{u5iiui} entail 
\begin{equation}\label{u543d}
\partial_j f\in H^{s-1}(\Omega)=\big(H^{1-s}(\Omega)\big)^*,\quad\forall\,j\in\{1,\dots,n\}.
\end{equation}
Then the manner in which $\widetilde\gamma_N(f,F)$ is now defined as a functional in the space
$H^{s-(3/2)}(\partial\Omega)=\big(H^{(3/2)-s}(\partial\Omega)\big)^*$ is as follows: 
Given $\phi\in H^{(3/2)-s}(\partial\Omega)$, then for any $\Phi\in H^{2-s}(\Omega)$ 
such that $\gamma_D\Phi=\phi$ {\rm (}whose existence is ensured by the surjectivity of 
\eqref{eqn:gammaDs.1}{\rm )}, set
\begin{align}\label{2.9NEW}
{}_{H^{(3/2)-s}(\partial\Omega)}\big\langle\phi,
\widetilde\gamma_N(f,F)\big\rangle_{(H^{(3/2)-s}(\partial\Omega))^*}
&:=\sum_{j=1}^n {}_{H^{1-s}(\Omega)}\big\langle\partial_j\Phi,
\partial_j f\big\rangle_{(H^{1-s}(\Omega))^*}
\nonumber\\[2pt]
&\quad\,+{}_{H^{2-s}(\Omega)}\big\langle\Phi,F\big\rangle_{(H^{2-s}(\Omega))^*}.
\end{align}
\end{definition}
%%%%%%%

Regarding Definition~\ref{h6r4d5UU}, one observes that, in the context described there, 
$\partial_j\Phi\in H^{1-s}(\Omega)$ for each $j\in\{1,\dots,n\}$, by \eqref{eq:DDEEn.3}. 
By \eqref{u543d}, this shows that the pairings under the summation symbol in the right-hand 
side of \eqref{2.9NEW} are meaningful. In addition, one can canonically identify the distribution 
$F$, originally belonging to $H^{s-2}_0(\Omega)$, with a functional in $(H^{2-s}(\Omega))^*$ 
(cf. the discussion pertaining to \eqref{eq:Redxax.2} and \eqref{eq:Redxax.1REF}), so the last pairing 
in \eqref{2.9NEW} is also meaningfully defined as
\begin{align}\label{uyree-y5}
\begin{split}
& {}_{H^{2-s}(\Omega)}\big\langle\Phi,F\big\rangle_{(H^{2-s}(\Omega))^*}
={}_{H^{2-s}({\mathbb{R}}^n)}\big\langle\Theta,F\big\rangle_{H^{s-2}({\mathbb{R}}^n)}
\\[2pt]
& \quad\text{for any $\Theta\in H^{2-s}({\mathbb{R}}^n)$ 
satisfying $\Theta\big|_{\Omega}=\Phi$ in ${\mathcal{D}}'(\Omega)$}.
\end{split}
\end{align}

Our next theorem elaborates on the main properties of the weak Neumann trace operator defined above. 

%%%%%%%
\begin{theorem}\label{GNT}
Let $\Omega\subset{\mathbb{R}}^n$ be a bounded Lipschitz domain, and fix 
$s\in\big(\tfrac12,\tfrac32\big)$. Then the weak Neumann trace mapping 
\begin{equation}\label{2.88X-jussi}
\widetilde\gamma_N:\big\{(f,F)\in H^s(\Omega)\times H^{s-2}_0(\Omega)\,\big|\,
\Delta f=F|_{\Omega}\text{ in }\mathcal{D}'(\Omega)\big\}\rightarrow H^{s-(3/2)}(\partial\Omega)
\end{equation}
from Definition~\ref{h6r4d5UU} yields an operator which is unambiguously defined, 
linear, and bounded {\rm (}assuming the space on the left-hand side of \eqref{2.88X-jussi} 
is equipped with the natural norm $(f,F)\mapsto\|f\|_{H^{s}(\Omega)}+\|F\|_{H^{s-2}({\mathbb{R}}^n)}${\rm )}.
The weak Neumann boundary trace map possesses the following properties: \\[1mm] 
$(i)$ The weak Neumann trace operators corresponding to various values of the 
parameter $s\in\big(\tfrac12,\tfrac32\big)$ are compatible with one another and each 
of them is surjective. In fact, there exist linear and bounded operators
\begin{equation}\label{2.88X-NN}
\Upsilon_N:H^{s-(3/2)}(\partial\Omega)\rightarrow 
\big\{u\in H^s(\Omega)\,\big|\,\Delta u\in L^2(\Omega)\big\},\quad s\in\big(\tfrac12,\tfrac32\big),
\end{equation}
which are compatible with one another and satisfy {\rm (}with tilde denoting the extension by zero 
outside $\Omega${\rm )} 
\begin{equation}\label{2.88X-NN2}
\widetilde\gamma_N\big(\Upsilon_N\psi,\widetilde{\Delta(\Upsilon_N\psi)}\,\big)=\psi,
\quad\forall\,\psi\in H^{s-(3/2)}(\partial\Omega)\,\text{ with }\,s\in\big(\tfrac12,\tfrac32\big).
\end{equation}
$(ii)$ Given any two pairs, 
\begin{align}\label{ju7653es}
\begin{split}
& \text{$(f,F)\in H^s(\Omega)\times H^{s-2}_0(\Omega)$ 
such that $\Delta f=F|_{\Omega}$ in $\mathcal{D}'(\Omega)$,}
\\[2pt] 
& \quad\text{and $(g,G)\in H^{2-s}(\Omega)\times H^{-s}_0(\Omega)$ such that $\Delta g=G|_{\Omega}$ 
in $\mathcal{D}'(\Omega)$},
\end{split}
\end{align}
the following Green's formula holds:
\begin{align}\label{GGGRRR-prim}
& {}_{H^{(3/2)-s}(\partial\Omega)}\big\langle\gamma_D g,\widetilde\gamma_N(f,F)
\big\rangle_{(H^{(3/2)-s}(\partial\Omega))^*}
\nonumber\\[2pt]
& \qquad -{}_{(H^{s-(1/2)}(\partial\Omega))^*}\big\langle\widetilde\gamma_N(g,G),\gamma_D f
\big\rangle_{H^{s-(1/2)}(\partial\Omega)}     
\nonumber\\[2pt]
&\quad={}_{H^{2-s}(\Omega)}\big\langle g,F\big\rangle_{(H^{2-s}(\Omega))^*}
-{}_{(H^s(\Omega))^*}\big\langle G,f\big\rangle_{H^s(\Omega)}.
\end{align}
\end{theorem}
%%%%%%%
\begin{proof}
We start by presenting the proof of the opening statement of the theorem. 
Pick a pair $(f,F)$ belonging to the domain of $\widetilde\gamma_N$ in \eqref{2.88X}. 
We note that the right-hand side of \eqref{2.9NEW} is independent of the particular 
extension $\Phi$ of $\phi$, as may be seen with the help of \eqref{incl-Yb.EE} and 
\eqref{Rdac}. Hence, $\widetilde\gamma_N(f,F)$ is well defined as a functional 
in $\big(H^{(3/2)-s}(\partial\Omega)\big)^*$ and satisfies the natural estimate
\begin{equation}\label{eq:thgcd-NEW}
\|\widetilde\gamma_N(f,F)\|_{H^{s-(3/2)}(\partial\Omega)}\leq C\big(
\|f\|_{H^{s}(\Omega)}+\|F\|_{H^{s-2}({\mathbb{R}}^n)}\big),
\end{equation}
for some constant $C\in(0,\infty)$ independent of $(f,F)$. Indeed, 
\begin{align}\label{eq:thgcd-NEW.rE}
\begin{split} 
& \|\widetilde\gamma_N(f,F)\|_{H^{s-(3/2)}(\partial\Omega)}
=\|\widetilde\gamma_N(f,F)\|_{(H^{(3/2)-s}(\partial\Omega))^*}
\\[2pt]
& \quad=\sup_{\substack{\phi\in H^{(3/2)-s}(\partial\Omega)\\ \|\phi\|_{H^{(3/2)-s}(\partial\Omega)}\leq 1}}
\Big|{}_{H^{(3/2)-s}(\partial\Omega)}\big\langle\phi,
\widetilde\gamma_N(f,F)\big\rangle_{(H^{(3/2)-s}(\partial\Omega))^*}\Big|. 
\end{split} 
\end{align}
Moreover, for every $\phi\in H^{(3/2)-s}(\partial\Omega)$ with $\|\phi\|_{H^{(3/2)-s}(\partial\Omega)}\leq 1$, 
if $\vartheta_D$ is the extension operator described in \eqref{eqn:gammaDs.1-INV1}--\eqref{eqn:gammaDs.1-INV2} 
one estimates  
\begin{align}\label{eq:thgcd-NEW.rE.2}
& \Big|{}_{H^{(3/2)-s}(\partial\Omega)}\big\langle\phi, 
\widetilde\gamma_N(f,F)\big\rangle_{(H^{(3/2)-s}(\partial\Omega))^*}\Big|
\nonumber\\[2pt]
& \quad\leq\sum_{j=1}^n\Big|{}_{H^{1-s}(\Omega)}\big\langle\partial_j(\vartheta_D\phi),
\partial_j f\big\rangle_{(H^{1-s}(\Omega))^*}\Big|
\nonumber\\[2pt]
&\qquad\,+\Big|{}_{H^{2-s}(\Omega)}\big\langle\vartheta_D\phi,F\big\rangle_{(H^{2-s}(\Omega))^*}\Big|
\nonumber\\[2pt]
& \quad\leq\sum_{j=1}^n\|\partial_j(\vartheta_D\phi)\|_{H^{1-s}(\Omega)}
\|\partial_j f\|_{(H^{1-s}(\Omega))^*}
\nonumber\\[2pt]
&\qquad\,+\|\vartheta_D\phi\|_{H^{2-s}(\Omega)}\|F\big\|_{(H^{2-s}(\Omega))^*}
\nonumber\\[2pt]
& \quad\leq C\|\vartheta_D\phi\|_{H^{2-s}(\Omega)}\big(\|f\|_{H^{s}(\Omega)}+\|F\|_{H^{s-2}({\mathbb{R}}^n)}\big)
\nonumber\\[2pt]
& \quad\leq C\big(\|f\|_{H^{s}(\Omega)}+\|F\|_{H^{s-2}({\mathbb{R}}^n)}\big),
\end{align}
using \eqref{u5iiui}, \eqref{eq:Redxax.1}, \eqref{u54ruy}, \eqref{eq:DDEEn.3}, and the fact that 
$\|\vartheta_D\phi\|_{H^{2-s}(\Omega)}\leq C$, for some constant $C\in(0,\infty)$ independent of $f$. 
This proves \eqref{eq:thgcd-NEW}. \\[1mm] 

We now address the claims made in itemized portion of the statement of the theorem. \\[1mm] 
{\it Proof of $(i)$.}
That the weak Neumann trace operators corresponding to various values of the 
parameter $s\in\big(\tfrac12,\tfrac32\big)$ are compatible with one another
is implied by the compatibility of the duality pairings intervening in \eqref{2.9NEW}.

Next, given any $s\in\big(\tfrac12,\tfrac32\big)$, consider the operator 
\begin{equation}\label{2.88X-NN-44}
\Upsilon_N:
\begin{cases} 
H^{s-(3/2)}(\partial\Omega)\rightarrow\big\{u\in H^s(\Omega)\,\big|\,\Delta u\in L^2(\Omega)\big\},  
\\[2pt]
\psi\mapsto\Upsilon_N\psi:=u,
\end{cases}
\end{equation}
where $u$ is the unique solution of 
\begin{equation} 
\begin{cases}
(-\Delta+1)u=0\,\text{ in $\Omega,\quad u\in H^{s}(\Omega)$,}     
\\[2pt]   
\widetilde\gamma_N(u,\widetilde{u})=\psi\in H^{s-(3/2)}(\partial\Omega).
\end{cases}   
\end{equation}
In this regard, it is worth noting that, since $s\in\big(\tfrac12,\tfrac32\big)$, 
picking some $r\in\big(-\tfrac12,\tfrac12\big)$ (e.g., $r=0$ will do) allows us 
to write, on account of \eqref{eq:DDEEn.2} and \eqref{eq:12d-yw5},
\begin{equation}\label{ug5d443}
u\in H^s(\Omega)\subset H^r(\Omega)\Longrightarrow
\widetilde u\in H^r_0(\Omega)\subset H^{s-2}_0(\Omega).
\end{equation}
Hence, $\widetilde u\in H_0^{s-2}(\Omega)$ and, in addition, 
$\widetilde u\big|_{\Omega}=u=\Delta u$ in $\Omega$. This ensures that the weak Neumann boundary 
trace $\widetilde\gamma_N(u,\widetilde{u})$ has meaning (cf. Definition~\ref{h6r4d5UU}). 
That the Neumann boundary value problem for the Helmholtz operator formulated in \eqref{2.88X-NN-44}
is well posed is a consequence of work in \cite{FMM98}, \cite{Mi96}, \cite{MT00}. This implies that
$\Upsilon_N$ is well defined, linear, and bounded. Moreover, when viewed as a family 
indexed by the parameter $s\in\big(\tfrac12,\tfrac32\big)$, the operators $\Upsilon_N$ 
act in a compatible fashion. Then for each $\psi\in H^{s-(3/2)}(\partial\Omega)$ 
with $s\in\big(\tfrac12,\tfrac32\big)$ one has 
\begin{align}\label{2.9NEjht}
\widetilde\gamma_N\big(\Upsilon_N\psi,\widetilde{\Delta(\Upsilon_N\psi)}\,\big)
=\widetilde\gamma_N\big(u,\widetilde{u}\,\big)=\psi,
\end{align}
proving \eqref{2.88X-NN2}. Of course, this also shows that each weak Neumann trace operator
$\widetilde\gamma_N$ is surjective in the context of \eqref{2.88X-jussi}. \\[1mm] 
{\it Proof of $(ii)$.}
Green's formula \eqref{GGGRRR-prim} readily follows by a two-fold application of \eqref{2.9NEW}.
\end{proof}
%%%%%%%

We shall build in the direction of including the end-point cases 
$s=\tfrac{1}{2}$ and $s=\tfrac{3}{2}$ in \eqref{2.88X}. As a preamble, we first define 
a Neumann trace operator acting from spaces of null-solutions of the Helmholtz 
operator $-\Delta+1$ from $H^{1/2}(\Omega)$ and $H^{3/2}(\Omega)$. The underlying reason 
why we prefer to work with a Helmholtz operator in place of the Laplacian is that 
we employ layer potentials, and the layer potentials associated with the Laplacian are, 
as opposed to those associated with the Helmholtz operator, sensitive to the 
topology of the underlying domain (cf. \cite{Mi97} in this regard).

%%%%%%%
\begin{lemma}\label{YTfBBVa}
Assume that $\Omega\subset\bbR^n$ is a bounded Lipschitz domain with outward 
unit normal $\nu$. Fix $\kappa>0$ and introduce
\begin{align}\label{eq:Dcaav.1}
{\mathscr{V}}(\Omega) & :=\big\{v\in H^{1/2}(\Omega)\,\big|\,(-\Delta+1)v=0\,\text{ in }\,\Omega\big\},
\\[2pt]
{\mathscr{W}}(\Omega) & :=\big\{w\in H^{3/2}(\Omega)\,\big|\,(-\Delta+1)w=0\,\text{ in }\,\Omega\big\}.
\label{eq:Dcaav.2}
\end{align}
Then ${\mathscr{V}}(\Omega) $ and ${\mathscr{W}}(\Omega)$ 
are closed subspaces of $H^{1/2}(\Omega)$ and $H^{3/2}(\Omega)$, respectively. Moreover,
\begin{align}\label{eq:Dcaav.1a}
\begin{split}
& {\mathscr{V}}(\Omega)=\big\{v\in C^\infty(\Omega)\,\big|\,(-\Delta+1)v=0\,\text{ in }\,\Omega
\,\text{ and }\,{\mathcal{N}}_\kappa v\in L^2(\partial\Omega)\big\}, 
\\[2pt]
& {\mathscr{W}}(\Omega)=\big\{w\in C^\infty(\Omega)\,\big|\,(-\Delta+1)w=0\,\text{ in }\Omega,
\,\,{\mathcal{N}}_\kappa w,{\mathcal{N}}_\kappa(\nabla w)\in L^2(\partial\Omega)\big\},
\end{split}
\end{align}
and the Dirichlet trace induces continuous isomorphisms in the following contexts:
\begin{equation}\label{eq:JbaiGV}
\gamma_D:{\mathscr{V}}(\Omega)\rightarrow L^2(\partial\Omega),\quad
\gamma_D:{\mathscr{W}}(\Omega)\rightarrow H^1(\partial\Omega).
\end{equation}
In addition, considering 
\begin{equation}\label{eq:Gva444rt}
\gamma^{\mathscr{V}}_N:{\mathscr{V}}(\Omega)\rightarrow H^{-1}(\partial\Omega)
=\big(H^{1}(\partial\Omega)\big)^\ast,
\end{equation}
defined by setting for each $v\in{\mathscr{V}}(\Omega)$ and each $\phi\in H^{1}(\partial\Omega)$,
\begin{equation}\label{eq:Jba-hr4ee}
{}_{H^{-1}(\partial\Omega)}\big\langle\gamma^{\mathscr{V}}_N v,\phi\big\rangle_{H^{1}(\partial\Omega)}
:=\big(\gamma_D v,\nu\cdot\gamma_D(\nabla w)\big)_{L^2(\partial\Omega)}, 
\end{equation}
where $w$ is the unique function in ${\mathscr{W}}(\Omega)$ such that $\gamma_D w=\phi$, 
then the operator $\gamma^{\mathscr{V}}_N$ in \eqref{eq:Gva444rt}--\eqref{eq:Jba-hr4ee} is a continuous isomorphism.

Finally, the assignment 
\begin{equation}\label{eq:ERfv274}
{\mathscr{W}}(\Omega)\ni w\mapsto
\nu\cdot\Big((\nabla w)\big|^{\kappa-{\rm n.t.}}_{\partial\Omega}\Big)\in L^2(\partial\Omega)
\end{equation}
is also a continuous isomorphism.
\end{lemma}
%%%%%%%
\begin{proof}
That ${\mathscr{V}}(\Omega) $ and ${\mathscr{W}}(\Omega)$ 
are closed subspaces of $H^{1/2}(\Omega)$ and $H^{3/2}(\Omega)$, respectively, is clear from definitions. 
The fact that the spaces ${\mathscr{V}}(\Omega),{\mathscr{W}}(\Omega)$, originally defined 
as in \eqref{eq:Dcaav.1}--\eqref{eq:Dcaav.2} may be alternatively described as in 
\eqref{eq:Dcaav.1a} is a direct consequence of \eqref{eq:MM3BIS}. Next, let $E_1(\cdot)$ 
denote the standard fundamental solution for the Helmholtz operator 
$-\Delta+1$ in ${\mathbb{R}}^n$, $n\geq 2$, that is,
\begin{equation}\label{trree.P}
E_1(x):=(i/4)\big(-2\pi i|x|\big)^{(2-n)/2} H^{(1)}_{(n-2)/2}\big(i|x|\big),
\quad\forall\,x\in\bbR^n\backslash\{0\}, 
\end{equation}
where $H^{(1)}_{\lambda}(\cdot)$ denotes the Hankel function of the first kind 
with index $\lambda\geq 0$ (cf. \cite[Section~9.1]{AS72}). In addition, given 
$f\in L^2(\partial\Omega)$, consider the integral operators 
\begin{align}\label{eq:GCa.1}
{\mathscr{S}}f(x) &:=\int_{\partial\Omega}E_1(x-y)f(y)\,d^{n-1}\sigma(y),\quad\forall\,x\in\Omega,
\\[2pt]
Sf(x) &:=\int_{\partial\Omega}E_1(x-y)f(y)\,d^{n-1}\sigma(y),\quad\forall\,x\in\partial\Omega,
\label{eq:GCa.2}
\\[2pt]
Kf(x) &:=\lim_{\varepsilon\to 0_{+}}\int_{\partial\Omega\backslash B(x,\varepsilon)}
\nu(y)\cdot(\nabla E_1)(x-y)f(y)\,d^{n-1}\sigma(y),\quad\forall\,x\in\partial\Omega.
\label{eq:GCa.2K}
\end{align}
Then from the work in \cite{Mi96}, \cite{MT00}, \cite{MT00a}, \cite{MT01}, it is known that 
for each $f\in L^2(\partial\Omega)$ the principal value defining $Kf(x)$ exists for $\sigma$-a.e.
$x\in\partial\Omega$, and $K$ is a well defined and bounded operator both on $L^2(\partial\Omega)$ 
and on $H^1(\partial\Omega)$. In addition, for each $f\in L^2(\partial\Omega)$ one has  
\begin{equation}\label{eKJa77}
{\mathscr{S}}f\big|^{\kappa-{\rm n.t.}}_{\partial\Omega}(x)=Sf(x)\,\text{ for $\sigma$-a.e.~$x\in\partial\Omega$},
\end{equation}
and
\begin{equation}\label{eq:DabV6h}
\nu(x)\cdot\Big(\nabla{\mathscr{S}}f\big|^{\kappa-{\rm n.t.}}_{\partial\Omega}\Big)(x) 
=\big(-\tfrac{1}{2}I+K^\ast\big)f(x)\,\text{ for $\sigma$-a.e.~$x\in\partial\Omega$},
\end{equation}
where $K^\ast$ is the adjoint of $K$ acting on $L^2(\partial\Omega)$. 
In addition, these operators induce continuous isomorphisms in the following contexts:
\begin{align}\label{eq:GCa.3}
& {\mathscr{S}}:L^2(\partial\Omega)\rightarrow{\mathscr{W}}(\Omega),\quad
{\mathscr{S}}:H^{-1}(\partial\Omega)\rightarrow{\mathscr{V}}(\Omega),   
\\[2pt] 
& S:H^{-1}(\partial\Omega)\rightarrow L^2(\partial\Omega),\quad
S:L^2(\partial\Omega)\rightarrow H^{1}(\partial\Omega),
\label{eq:GCa.6}  
\\[2pt] 
& \pm\tfrac{1}{2}I+K:L^2(\partial\Omega)\rightarrow L^2(\partial\Omega),\quad
\pm\tfrac{1}{2}I+K:H^1(\partial\Omega)\rightarrow H^{1}(\partial\Omega).
\label{eq:GCa.7}
\end{align}
In fact, the operators in \eqref{eq:GCa.6} are adjoints to one another.
In addition, the two Dirichlet boundary traces from \eqref{eq:JbaiGV} coincide with 
the operator $S\circ{\mathscr{S}}^{-1}$, acting from ${\mathscr{V}}(\Omega)$ onto 
$L^2(\partial\Omega)$, and from ${\mathscr{W}}(\Omega)$ onto $H^1(\partial\Omega)$, 
respectively. Hence, they induce continuous isomorphisms in the context of 
\eqref{eq:JbaiGV}. Consequently, given any $\phi\in H^{1}(\partial\Omega)$, 
if $w$ is the unique function in ${\mathscr{W}}(\Omega)$ such that $\gamma_D w=\phi$, 
then necessarily $w={\mathscr{S}}(S^{-1}\phi)$ in $\Omega$. Based on this, \eqref{eq:DabV6h}, 
and \eqref{eqn:gammaDs.2auxBBB}, for each function $v\in{\mathscr{V}}(\Omega)$ one can then write 
\begin{align}\label{eq:Jba-hfV}
\big(\gamma_D v,\nu\cdot\gamma_D(\nabla w)\big)_{L^2(\partial\Omega)}
&=\big(\gamma_D v\,,\,\big(-\tfrac{1}{2}I+K^\ast\big)(S^{-1}\phi)\big)_{L^2(\partial\Omega)}
\nonumber\\[2pt]
&={}_{H^{-1}(\partial\Omega)}\big\langle S^{-1}\big(-\tfrac{1}{2}I+K\big)(\gamma_D v)\,,\,
\phi\big\rangle_{H^{1}(\partial\Omega)}.
\end{align}
In light of \eqref{eq:Jba-hr4ee}, this proves that
\begin{equation}\label{eq:Gva4kpp}
\gamma^{\mathscr{V}}_N v=S^{-1}\big(-\tfrac{1}{2}I+K\big)(\gamma_D v)
\,\text{ for each $v\in{\mathscr{V}}(\Omega)$}.
\end{equation}
From \eqref{eq:GCa.6}--\eqref{eq:GCa.7}, the fact that 
$\gamma_D:{\mathscr{V}}(\Omega)\rightarrow L^2(\partial\Omega)$ 
is a continuous isomorphism, and \eqref{eq:Gva4kpp} one concludes that the 
operator $\gamma^{\mathscr{V}}_N$ in \eqref{eq:Gva444rt}--\eqref{eq:Jba-hr4ee} is a continuous isomorphism.

Finally, regarding \eqref{eq:ERfv274}, starting from the fact that any function 
$w\in{\mathscr{W}}(\Omega)$ may be represented as $w={\mathscr{S}}(S^{-1}\gamma_D w)$ 
in $\Omega$, one deduces from \eqref{eq:DabV6h} 
\begin{equation}\label{eq:ERfv2ii}
\nu\cdot\Big((\nabla w)\big|^{\kappa-{\rm n.t.}}_{\partial\Omega}\Big)
=\big(-\tfrac{1}{2}I+K^\ast\big)(S^{-1}\gamma_D w),\quad\forall\,w\in{\mathscr{W}}(\Omega).
\end{equation}
Then the claim about \eqref{eq:ERfv274} becomes a consequence of this and 
the fact that the mappings in \eqref{eq:JbaiGV} and \eqref{eq:GCa.3}--\eqref{eq:GCa.7}
are continuous isomorphisms.
\end{proof}
%%%%%%%

Our main result pertaining to the Neumann boundary trace operator is contained in the 
theorem below. As in the case of the Dirichlet trace, by restricting ourselves to functions
with a better-than-expected Laplacian (in the sense of membership within the Sobolev
scale) we are able to include the end-point cases $s=\tfrac12$ and $s=\tfrac32$ in \eqref{2.88X}.
Expanding the action of the weak Neumann boundary trace map in this fashion is 
going to be crucially important in our future endeavors. 

%%%%%%%
\begin{theorem} \label{YTfdf.NNN.2-Main}
Assume that $\Omega\subset\mathbb R^n$ is a bounded Lipschitz domain. Then for each $\varepsilon>0$ 
the weak Neumann boundary trace map, originally introduced in Definition~\ref{h6r4d5UU}, induces 
linear and continuous operators in the context 
\begin{equation}\label{eqn:gammaN}
\begin{array}{c}
\widetilde\gamma_N:\big\{(f,F)\in H^s(\Omega)\times H^{s-2+\varepsilon}_0(\Omega)\,\big|\,
\Delta f=F\big|_{\Omega}\,\text{ in }\,{\mathcal{D}}'(\Omega)\big\}\rightarrow H^{s-(3/2)}(\partial\Omega)
\\[6pt]
\text{with }\,\,s\in\big[\tfrac{1}{2},\tfrac{3}{2}\big]
\end{array}
\end{equation}
{\rm (}throughout, the space on the left-hand side of \eqref{eqn:gammaN} equipped with the natural norm 
$(f,F)\mapsto\|f\|_{H^{s}(\Omega)}+\|F\|_{H^{s-2+\varepsilon}({\mathbb{R}}^n)}${\rm )}
which are compatible with those in Definition~\ref{h6r4d5UU} when $s\in\big(\tfrac{1}{2},\tfrac{3}{2}\big)$.
Thus defined, the weak Neumann boundary trace map possesses the following additional properties: 
\\[1mm] 
$(i)$ Each weak Neumann boundary trace map in \eqref{eqn:gammaN} is surjective.
In fact, there exist linear and bounded operators
\begin{equation}\label{2.88X-NN-ii}
\Upsilon_N:H^{s-(3/2)}(\partial\Omega)\rightarrow 
\big\{u\in H^s(\Omega)\,\big|\,\Delta u\in L^2(\Omega)\big\},\quad s\in\big[\tfrac12,\tfrac32\big],
\end{equation}
which are compatible with one another and satisfy {\rm (}with tilde denoting the extension by zero 
outside $\Omega${\rm )} 
\begin{equation}\label{2.88X-NN2-ii}
\widetilde\gamma_N\big(\Upsilon_N\psi,\widetilde{\Delta(\Upsilon_N\psi)}\,\big)=\psi,
\quad\forall\,\psi\in H^{s-(3/2)}(\partial\Omega)\,\text{ with }\,s\in\big[\tfrac12,\tfrac32\big].
\end{equation}
$(ii)$ If $\varepsilon\in(0,1)$ and $s\in\big[\tfrac{1}{2},\tfrac{3}{2}\big]$ then for any two pairs
\begin{align}\label{n7b66}
\begin{split}
& \text{$(f,F)\in H^s(\Omega)\times H^{s-2+\varepsilon}_0(\Omega)$ 
such that $\Delta f=F|_{\Omega}$ in $\mathcal{D}'(\Omega)$,}
\\[2pt]
& \quad\text{and $(g,G)\in H^{2-s}(\Omega)\times H^{-s+\varepsilon}_0(\Omega)$ 
such that $\Delta g=G|_{\Omega}$ in $\mathcal{D}'(\Omega)$}, 
\end{split}
\end{align}
the following Green's formula holds:
\begin{align}\label{GGGRRR-prim222}
& {}_{H^{(3/2)-s}(\partial\Omega)}\big\langle\gamma_D g,\widetilde\gamma_N(f,F)
\big\rangle_{(H^{(3/2)-s}(\partial\Omega))^*}
\nonumber\\[2pt]
& \qquad -{}_{(H^{s-(1/2)}(\partial\Omega))^*}\big\langle\widetilde\gamma_N(g,G),\gamma_D f
\big\rangle_{H^{s-(1/2)}(\partial\Omega)}     
\nonumber\\[2pt]
& \quad={}_{H^{2-s}(\Omega)}\big\langle g,F\big\rangle_{(H^{2-s}(\Omega))^*}
-{}_{(H^s(\Omega))^*}\big\langle G,f\big\rangle_{H^s(\Omega)}.
\end{align}
$(iii)$ There exists a constant $C\in(0,\infty)$ with the property that
\begin{align}\label{gafvv.6588-P}
\begin{split}
& \text{if $f\in H^{1/2}(\Omega)$ and $F\in H^{-(3/2)+\varepsilon}_0(\Omega)$ with 
$0<\varepsilon\leq 1$ satisfy} 
\\[2pt]
& \quad\text{$\Delta f=F\big|_{\Omega}$ in ${\mathcal{D}}'(\Omega)$ and 
$\widetilde\gamma_N(f,F)=0$, then $f\in H^{(1/2)+\varepsilon}(\Omega)$ }
\\[2pt]
& \quad\text{and $\|f\|_{H^{(1/2)+\varepsilon}(\Omega)}\leq 
C\big(\|f\|_{L^2(\Omega)}+\|F\|_{H^{-(3/2)+\varepsilon}(\mathbb{R}^n)}\big)$ holds.}
\end{split}
\end{align} 
$(iv)$ Denote by $\nu$ the outward unit normal vector to $\Omega$. Then
\begin{align}\label{eq:Nnan7yg-P}
\begin{split}
& \text{if $f\in H^{3/2}(\Omega)$ and $F\in H^{-(1/2)+\varepsilon}_0(\Omega)$ 
for some $\varepsilon\in(0,1)$ satisfy} 
\\[2pt]
& \quad\text{$\Delta f=F\big|_{\Omega}$ in ${\mathcal{D}}'(\Omega)$
then, actually, $\widetilde\gamma_N(f,F)\in L^2(\partial\Omega)$ and, in fact,} 
\\[2pt]
& \quad\text{$\widetilde\gamma_N(f,F)=\nu\cdot\gamma_D(\nabla f)$
with the Dirichlet trace taken as in \eqref{eqn:gammaDs.2aux}.}
\end{split}
\end{align}
Moreover, there exists a constant $C\in(0,\infty)$ with the property that in the context of \eqref{eq:Nnan7yg-P} one has
\begin{align}\label{ut543i-Ai.WACO}
\big\|\widetilde\gamma_N(f,F)\big\|_{L^2(\partial\Omega)}
\leq C\big(\|f\|_{H^{3/2}(\Omega)}+\|F\|_{H^{-(1/2)+\varepsilon}(\bbR^n)}\big).
\end{align}
$(v)$ Recall \eqref{restr-0}. Under the assumption that 
\begin{equation}\label{ihg65f+AAA}
\varepsilon>0,\quad s\in\big[\tfrac{1}{2},\tfrac{3}{2}\big],\,\,\text{ and }\,\,\varepsilon>\tfrac{3}{2}-s,
\end{equation}
it follows that the mapping 
\begin{align}\label{eqn:gammaN-III.a}
& {\mathcal{I}}:\big\{(f,F)\in H^s(\Omega)\times H^{s-2+\varepsilon}_0(\Omega)\,\big|\,
\Delta f=F\big|_{\Omega}\,\text{ in }\,{\mathcal{D}}'(\Omega)\big\}
\\[2pt]
& \quad\longrightarrow\big\{f\in H^s(\Omega):\,\Delta f\in H^{s-2+\varepsilon}_{z}(\Omega)\big\}
\nonumber
\end{align}
given by 
\begin{align}\label{eqn:gammaN-III.b}
\begin{split}
& {\mathcal{I}}(f,F):=f\,\text{ for each pair }\,(f,F)\in H^s(\Omega)\times H^{s-2+\varepsilon}_0(\Omega) 
\\[2pt]
& \quad\text{with }\,\Delta f=F\big|_{\Omega}\,\text{ in }\,{\mathcal{D}}'(\Omega),
\end{split}
\end{align}
is actually a continuous linear isomorphism. As a consequence of this and \eqref{Rdac.WACO}, 
under the assumption made in \eqref{ihg65f+AAA} it follows that the mapping 
\begin{align}\label{eqn:gammaN-III.c}
\dot\gamma_N:=\widetilde\gamma_N\circ{\mathcal{I}}^{\,-1}
\end{align}
is well defined, linear, and continuous in the context
\begin{align}\label{eqn:gammaN-III.d}
\dot\gamma_N:\big\{f\in H^s(\Omega)\,\big|\,\Delta f\in H^{s-2+\varepsilon}_z(\Omega)\big\}
\rightarrow H^{s-(3/2)}(\partial\Omega).
\end{align}
In view of the fact that \eqref{ihg65f+AAA} is satisfied if $\varepsilon>0$ and $s=\tfrac{3}{2}$, this together with 
\eqref{Rdac.WACO} further imply that the mapping in \eqref{eqn:gammaN-III.d} yields the following 
{\rm (}well defined, linear, continuous, surjective{\rm )} brand of Neumann trace operator:
\begin{align}\label{eqn:gammaN-III.e}
\begin{split}
&\dot\gamma_N:\big\{f\in H^{3/2}(\Omega)\,\big|\,\Delta f\in H^{-(1/2)+\varepsilon}(\Omega)\big\}
\longrightarrow L^2(\partial\Omega),\quad\dot\gamma_N(f):=\nu\cdot\gamma_D(\nabla f), 
\\[6pt]
& \quad\text{for each $\varepsilon\in(0,1)$, with the Dirichlet trace taken as in \eqref{eqn:gammaDs.2aux}.}
\end{split}
\end{align} 
\end{theorem}
%%%%%%%
\begin{proof}
We start by considering the claims made in the opening part and in item $(i)$ in the statement of the theorem. 
It is convenient to analyze three distinct cases, depending on the nature of the smoothness parameter 
$s\in\big[\tfrac12,\tfrac32\big]$. For the goals we have in mind, there is no loss of generality in 
assuming that $\varepsilon\in(0,1)$.
\\[1mm] 
\noindent{\bf Case~1:} {\it Assume} $s\in\big(\tfrac{1}{2},\tfrac{3}{2}\big)$. 
In this scenario, all desired conclusions follow from Theorem~\ref{GNT} (as well as its proof) 
simply by observing that $\big\{(f,F)\in H^s(\Omega)\times H^{s-2+\varepsilon}_0(\Omega)\,\big|\,
\Delta f=F\big|_{\Omega}\,\text{ in }\,{\mathcal{D}}'(\Omega)\big\}$, the domain of
$\widetilde\gamma_N$ in \eqref{eqn:gammaN}, is a subspace of 
$\big\{(f,F)\in H^s(\Omega)\times H^{s-2}_0(\Omega)\,\big|\,
\Delta f=F\big|_{\Omega}\,\text{ in }\,{\mathcal{D}}'(\Omega)\big\}$, the domain of
$\widetilde\gamma_N$ in \eqref{2.88X}. In addition, the same operators $\Upsilon_N$ 
from \eqref{2.88X-NN} will work in the current context. 
\\[1mm] 
\noindent{\bf Case~2:} {\it Assume} $s=\tfrac{3}{2}$.
Suppose now that some $f\in H^{3/2}(\Omega)$ along with some $F\in H^{-(1/2)+\varepsilon}_0(\Omega)$ 
satisfying $\Delta f=F\big|_{\Omega}$ in ${\mathcal{D}}'(\Omega)$ have been given. 
In particular, 
\begin{equation}\label{j6f432}
\Delta f\in H^{-(1/2)+\varepsilon}(\Omega)
\end{equation}
and, for each $j\in\{1,\dots,n\}$, the function $\partial_j f\in H^{1/2}(\Omega)$ satisfies 
\begin{equation}\label{j6f433}
\Delta(\partial_j f)=\partial_j(\Delta f)=\partial_j\big(F\big|_{\Omega}\big)
=(\partial_j F)\big|_{\Omega}\in H^{-(3/2)+\varepsilon}(\Omega).
\end{equation}
Hence, by \eqref{eqn:gammaDs.2aux} (used with $s=1/2$),
\begin{equation}\label{eq:Gav7gt5}
\text{$\gamma_D(\partial_j f)$ exists in $L^2(\partial\Omega)$ for each $j\in\{1,\dots,n\}$}.
\end{equation}
Pick now an arbitrary $\Phi\in C^\infty(\overline{\Omega})$ and set 
$\phi:=\Phi\big|_{\partial\Omega}$. In addition, consider the vector field 
\begin{equation}\label{elhgtf}
\vec{F}:=\overline{\Phi}\nabla f\,\text{ in }\,\Omega. 
\end{equation}
Then \eqref{eq:DDEj6g5} implies that $\vec{F}\in\big[H^{1/2}(\Omega)\big]^n$ and 
\begin{equation}\label{eq:NBb65a}
\Delta\vec{F}=(\overline{\Delta\Phi})\nabla f
+\overline{\Phi}\nabla(\Delta f)+2\big(\overline{\nabla\Phi}\cdot\nabla\partial_j f\big)_{1\leq j\leq n}
\in\big[H^{-(3/2)+\varepsilon}(\Omega)\big]^n,
\end{equation}
as well as 
\begin{equation}\label{lk8ht5}
{\rm div}\vec{F}=\overline{\nabla\Phi}\cdot\nabla f+\overline{\Phi}\Delta f
\in H^{-(1/2)+\varepsilon}(\Omega).
\end{equation}
As such, Theorem~\ref{Ygav-75-BIS} applies and yields, with the Dirichlet 
trace $\gamma_D(\nabla f)$ understood in the sense of \eqref{eq:Gav7gt5} 
(cf.~\eqref{tvdee} and \eqref{eqn:gammaDs.1-TR.2}),
\begin{align}\label{2.9MMa-ii}
& \big(\phi\,,\,\nu\cdot\gamma_D(\nabla f)\big)_{L^2(\partial\Omega)}
=\int_{\partial\Omega}\nu\cdot\gamma_D\vec{F}\,d^{n-1}\sigma
={}_{H^{(1/2)-\varepsilon}(\Omega)}\big\langle{\bf 1},{\rm div}\vec{F}
\big\rangle_{H^{-(1/2)+\varepsilon}(\Omega)}
\nonumber\\[2pt]
& \quad={}_{H^{(1/2)-\varepsilon}(\Omega)}\big\langle{\bf 1},\overline{\nabla\Phi}\cdot\nabla f
\big\rangle_{H^{-(1/2)+\varepsilon}(\Omega)}
+{}_{H^{(1/2)-\varepsilon}(\Omega)}\big\langle{\bf 1},\overline{\Phi}\Delta f
\big\rangle_{H^{-(1/2)+\varepsilon}(\Omega)}
\nonumber\\[2pt]
& \quad=\sum_{j=1}^n{}_{H^{(1/2)-\varepsilon}(\Omega)}\big\langle
\partial_j\Phi,\partial_j f\big\rangle_{H^{-(1/2)+\varepsilon}(\Omega)}
+{}_{H^{(1/2)-\varepsilon}(\Omega)}\big\langle\Phi,\Delta f
\big\rangle_{H^{-(1/2)+\varepsilon}(\Omega)}
\nonumber\\[2pt]
& \quad=\sum_{j=1}^n\big(\partial_j\Phi,\partial_j f\big)_{L^2(\Omega)}
+{}_{H^{(1/2)-\varepsilon}(\Omega)}\big\langle\Phi,F
\big\rangle_{(H^{(1/2)-\varepsilon}(\Omega))^*}
\end{align}
with $\nu$ and $\sigma$ denoting, respectively, the outward unit normal and 
surface measure on $\partial\Omega$. Above, the last step relies on the manner in which 
$(H^{(1/2)-\varepsilon}(\Omega))^*$ is identified with $H^{-(1/2)+\varepsilon}(\Omega)$
(see \eqref{u5iiui}--\eqref{u5iiui-bbb}). 

Of course, the fact that $f\in H^{3/2}(\Omega)$ entails 
$f\in H^{s}(\Omega)$ for any $s\in\big(\tfrac{1}{2},\tfrac{3}{2}\big)$ 
and, as such, a direct comparison of \eqref{2.9MMa-ii} and \eqref{2.9NEW} reveals that 
\begin{align}\label{2.9MMabn-P}
\begin{split}
& {}_{H^{(3/2)-s}(\partial\Omega)}\big\langle\phi,
\widetilde\gamma_N(f,F)\big\rangle_{(H^{(3/2)-s}(\partial\Omega))^*}
=\big(\phi\,,\,\nu\cdot\gamma_D(\nabla f)\big)_{L^2(\partial\Omega)} 
\\[2pt]
& \quad\text{for every $s\in\big(\tfrac{1}{2},\tfrac{3}{2}\big)$ and every 
$\phi\in\big\{\Phi\big|_{\partial\Omega}\,\big|\,\,\Phi\in C^\infty(\overline{\Omega})\big\}$.}
\end{split}
\end{align}
Since the latter space is dense in $L^2(\partial\Omega)$, this ultimately proves \eqref{eq:Nnan7yg-P}.
Moreover, based on \eqref{eqn:gammaDs.2aux} with $s=\tfrac{1}{2}$, \eqref{eq:DDEEn.3}, the fact that 
$\Delta f=F\big|_{\Omega}$ in ${\mathcal{D}}'(\Omega)$, and \eqref{HGaYga.3}, one estimates  
\begin{align}\label{ut543i-Ai}
\big\|\widetilde\gamma_N(f,F)\big\|_{L^2(\partial\Omega)} &\leq C
\big(\|\nabla f\|_{[H^{1/2}(\Omega)]^n}
+\|\Delta(\nabla f)\|_{[H^{-(3/2)+\varepsilon}(\Omega)]^n}\big)
\nonumber\\[2pt]
&\leq C\big(\|f\|_{H^{3/2}(\Omega)}+\|\Delta f\|_{H^{-(1/2)+\varepsilon}(\Omega)}\big)
\nonumber\\[2pt]
&=C\big(\|f\|_{H^{3/2}(\Omega)}+\|F\|_{H^{-(1/2)+\varepsilon}(\bbR^n)}\big)
\end{align}
for some constant $C\in(0,\infty)$, independent of $(f,F)$. 

The operator $\Upsilon_N$ in \eqref{2.88X-NN-ii} corresponding to $s=\tfrac32$ is defined
as in \eqref{2.88X-NN-44}, in which the boundary value problem is now understood as 
\begin{equation}\label{eNN2.iia}
\begin{cases}
u\in C^\infty(\Omega),\quad (-\Delta+1)u=0\,\text{ in }\,\Omega,
\\[4pt] 
{\mathcal{N}}_\kappa u,\,{\mathcal{N}}_\kappa(\nabla u)\in L^2(\partial\Omega),
\\[4pt]  
\nu\cdot\big(\nabla u\big|^{\kappa-{\rm n.t.}}_{\partial\Omega}\big)=\psi
\,\text{ $\sigma$-a.e.~on }\,\partial\Omega,\quad\,\psi\in L^2(\partial\Omega).
\end{cases}    
\end{equation}
Work in \cite{Mi94}, \cite[Theorem~6.1]{MT99}, prove that the latter problem is well posed.
Moreover, since this boundary value problem as well as the one intervening in \eqref{2.88X-NN-44} are 
solved using the same formalism based on boundary layer potentials, it follows that the corresponding
solution operators $\Upsilon_N$ act in a coherent manner. By \eqref{eq:MM3BIS}, \eqref{eq:Nnan7yg-P}, 
and \eqref{eqn:gammaDs.2auxBBB}, one deduces that
\begin{equation}\label{2.88X-NN2-ii-ZZ}
\widetilde\gamma_N\big(\Upsilon_N\psi,\widetilde{\Delta(\Upsilon_N\psi)}\,\big)=\psi,
\quad\forall\,\psi\in L^2(\partial\Omega),
\end{equation}
justifying \eqref{2.88X-NN2-ii} in the case when $s=\tfrac32$. Of course, this also proves
the surjectivity of the weak Neumann trace operator in the current case.
\\[1mm] 
\noindent{\bf Case~3:} {\it Assume} $s=\tfrac{1}{2}$. In this scenario, we 
begin by assigning a meaning to the weak Neumann boundary trace $\widetilde\gamma_N(f,F)$ 
when, for some $\varepsilon\in(0,1)$,  
\begin{equation}\label{u4338jh}
\text{$f\in H^{1/2}(\Omega)$ and $F\in H^{-(3/2)+\varepsilon}_0(\Omega)$ satisfy 
$\Delta f=F\big|_{\Omega}$ in ${\mathcal{D}}'(\Omega)$}. 
\end{equation}
Specifically, in a first stage we extend $f$ by zero outside $\Omega$, 
to a function $\widetilde{f}\in L^2({\mathbb{R}}^n)$, and consider 
\begin{align}\label{eNjyg-u51}
\begin{split}
& \eta:=\big(E_1\ast(-F+\widetilde{f})\big)\big|_{\Omega}\,\text{ so that }\,\eta\in H^{(1/2)+\varepsilon}(\Omega)
\\[2pt]
& \quad\text{with }\,\|\eta\|_{H^{(1/2)+\varepsilon}(\Omega)}
\leq C\big(\|f\|_{L^2(\Omega)}+\|F\|_{H^{-(3/2)+\varepsilon}({\mathbb{R}}^n)}\big),
\end{split}
\end{align}
for some $C\in(0,\infty)$ independent of $(f,F)$. We also note that
\begin{align}\label{eNjyg-u53}
(-\Delta+1)\eta &=(-\Delta+1)\Big[\big(E_1\ast(-F+\widetilde{f})\big)\big|_{\Omega}\Big]
\nonumber\\[2pt]
&=\big[(-\Delta+1)\big(E_1\ast(-F+\widetilde{f})\big)\big]\big|_{\Omega}
\nonumber\\[2pt]
&=\big[\big((-\Delta+1)E_1\big)\ast(-F+\widetilde{f})\big]\big|_{\Omega}
\nonumber\\[2pt]
&=(-F+\widetilde{f})\big|_{\Omega}\,\text{ in }\,{\mathcal{D}}'(\Omega).
\end{align}
In particular, if $\widetilde\eta\in L^2({\mathbb{R}}^n)$ is the extension by zero of $\eta$ 
to the entire Euclidean space, one has $F-\widetilde{f}+\widetilde\eta\in H^{-(3/2)+\varepsilon}_0(\Omega)$ 
and $\Delta\eta=(F-\widetilde{f}+\widetilde\eta\,)\big|_{\Omega}$. 
Given these facts, Theorem~\ref{GNT} applies and gives that
\begin{equation}\label{eNjyg-u54}
\widetilde\gamma_N(\eta,F-\widetilde{f}+\widetilde\eta\,)\in H^{-1+\varepsilon}(\partial\Omega)
\end{equation}
and, for some constant $C\in(0,\infty)$ independent of $(f,F)$, we have
\begin{align}\label{y65rd-PPT}
& \big\|\widetilde\gamma_N(\eta, F-\widetilde{f}+\widetilde\eta\,)\big\|_{H^{-1+\varepsilon}(\partial\Omega)}
\nonumber\\[2pt]
& \quad\leq C\big(\|\eta\|_{H^{(1/2)+\varepsilon}(\Omega)}
+\big\|F-\widetilde{f}+\widetilde\eta\,\big\|_{H^{-(3/2)+\varepsilon}({\mathbb{R}}^n)}\big)
\nonumber\\[2pt]
& \quad\leq C\big(\|\eta\|_{H^{(1/2)+\varepsilon}(\Omega)}
+\|F\|_{H^{-(3/2)+\varepsilon}({\mathbb{R}}^n)}
\nonumber\\
&\qquad +\big\|\widetilde{f}\,\big\|_{H^{-(3/2)+\varepsilon}({\mathbb{R}}^n)}
+\|\widetilde\eta\,\|_{H^{-(3/2)+\varepsilon}({\mathbb{R}}^n)}\big)
\nonumber\\[2pt]
& \quad\leq C\big(\|\eta\|_{H^{(1/2)+\varepsilon}(\Omega)}+\|F\|_{H^{-(3/2)+\varepsilon}({\mathbb{R}}^n)}
+\big\|\widetilde{f}\,\big\|_{L^2({\mathbb{R}}^n)}+\|\widetilde\eta\,\|_{L^2({\mathbb{R}}^n)}\big)
\nonumber\\[2pt]
& \quad\leq C\big(\|\eta\|_{H^{(1/2)+\varepsilon}(\Omega)}+\|F\|_{H^{-(3/2)+\varepsilon}({\mathbb{R}}^n)}
+\|f\|_{L^2(\Omega)}+\|\eta\|_{L^2(\Omega)}\big)
\nonumber\\[2pt]
& \quad\leq C\big(\|f\|_{L^2(\Omega)}+\|F\|_{H^{-(3/2)+\varepsilon}({\mathbb{R}}^n)}\big),
\end{align}
where the last inequality uses \eqref{eNjyg-u51}.
In a second stage, we consider the Neumann boundary problem 
\begin{equation}\label{eNN2-Q}
\begin{cases}
(-\Delta+1)\vartheta=0\,\text{ in }\Omega,\quad\vartheta\in H^{(1/2)+\varepsilon}(\Omega),     
\\[2pt]  
\widetilde\gamma_N(\vartheta,\widetilde{\vartheta})
=\widetilde\gamma_N(\eta,F-\widetilde{f}+\widetilde\eta\,)\in H^{-1+\varepsilon}(\partial\Omega),
\end{cases}    
\end{equation}
where $\widetilde{\vartheta}\in L^2({\mathbb{R}}^n)$ is the extension of $\vartheta$ by zero to 
${\mathbb{R}}^n$. From the work in \cite{FMM98}, \cite{Mi96}, \cite{MT00}, it follows that this 
has a unique solution which, by \eqref{y65rd-PPT}, satisfies
\begin{align}\label{eNjyg-u57}
\|\vartheta\|_{H^{(1/2)+\varepsilon}(\Omega)}
&\leq C\big\|\widetilde\gamma_N\big(\eta,F-\widetilde{f}
+\widetilde\eta\,\big)\big\|_{H^{-1+\varepsilon}(\partial\Omega)}
\nonumber\\[2pt]
&\leq C\big(\|f\|_{L^2(\Omega)}+\|F\|_{H^{-(3/2)+\varepsilon}({\mathbb{R}}^n)}\big),
\end{align}
for some constant $C\in(0,\infty)$, independent of $(f,F)$. In a third stage, define 
\begin{equation}\label{eNjyg-u58}
v:=(f-\eta+\vartheta)\in H^{1/2}(\Omega)
\end{equation}
and note that, in the sense of distributions in $\Omega$, 
\begin{align}\label{eNN2-Qe}
(-\Delta+1)v &=(-\Delta+1)f-(-\Delta+1)\eta
\nonumber\\[2pt]
&=(-\Delta+1)f-(-F+\widetilde{f})\big|_{\Omega}
\nonumber\\[2pt]
&=(-\Delta+1)f+\Delta f-f=0,
\end{align}
by \eqref{eNjyg-u58}, \eqref{eNN2-Q}, \eqref{eNjyg-u53}, and the last 
condition in \eqref{u4338jh}. In particular, $v\in{\mathscr{V}}(\Omega)$, the 
space introduced in \eqref{eq:Dcaav.1}. Given this, it then makes sense to finally define 
\begin{align}\label{eNN2-Qe2}
\widetilde\gamma_N(f,F):=\gamma^{\mathscr{V}}_N v\in H^{-1}(\partial\Omega),
\end{align}
with $\gamma^{\mathscr{V}}_N v$ defined in the sense of \eqref{eq:Gva444rt}--\eqref{eq:Jba-hr4ee}.
As a consequence of this definition, one confirms that the assignment 
$(f,F)\mapsto\widetilde\gamma_N(f,F)$ is linear and 
\begin{equation}\label{eNjyg-u599}
\big\|\widetilde\gamma_N(f,F)\big\|_{H^{-1}(\partial\Omega)}
\leq C\big(\|f\|_{H^{1/2}(\Omega)}+\|F\|_{H^{-(3/2)+\varepsilon}({\mathbb{R}}^n)}\big),
\end{equation}
for some constant $C\in(0,\infty)$, independent of $(f,F)$. Indeed, on the one hand, \eqref{eNN2-Qe2}, 
the boundedness of \eqref{eq:Gva444rt}, and \eqref{eNjyg-u58} permit us to estimate 
\begin{align}\label{eNjyg-u599-bfe}
\big\|\widetilde\gamma_N(f,F)\big\|_{H^{-1}(\partial\Omega)}
&=\big\|\gamma^{\mathscr{V}}_N v\big\|_{H^{-1}(\partial\Omega)}
\leq C\|v\big\|_{H^{1/2}(\Omega)}
\nonumber\\[2pt]
&\leq C\big(\|f\|_{H^{1/2}(\Omega)}+\|\eta\|_{H^{1/2}(\Omega)}
+\|\vartheta\|_{H^{1/2}(\Omega)}\big),
\end{align}
while, on the other hand, \eqref{eq:DDEEn.2}, \eqref{eNjyg-u51}, and \eqref{eNjyg-u57} give
\begin{align}\label{eNjyg-u51-TT.1}
& \|\eta\big\|_{H^{1/2}(\Omega)}\leq C\|\eta\|_{H^{(1/2)+\varepsilon}(\Omega)}
\leq C\big(\|f\|_{L^2(\Omega)}+\|F\|_{H^{-(3/2)+\varepsilon}({\mathbb{R}}^n)}\big),
\\[2pt]
& \|\vartheta\big\|_{H^{1/2}(\Omega)}\leq C\|\vartheta\|_{H^{(1/2)+\varepsilon}(\Omega)}
\leq C\big(\|f\|_{L^2(\Omega)}+\|F\|_{H^{-(3/2)+\varepsilon}({\mathbb{R}}^n)}\big).
\label{eNjyg-u51-TT.2}
\end{align}
Collectively, \eqref{eNjyg-u599-bfe}--\eqref{eNjyg-u51-TT.2} prove \eqref{eNjyg-u599}.

For future references, it is useful to observe that 
\begin{align}\label{8yhGVV}
\begin{split}
& \text{for each }\,v\in{\mathscr{V}}(\Omega)\,\text{ one has }\, 
\widetilde\gamma_N(v,\widetilde{v}\,)=\gamma^{\mathscr{V}}_N v\,\text{ where}
\\[2pt]
& \quad\widetilde{v}\in L^2({\mathbb{R}}^n)\,\text{ is the extension of $v$ by zero outside $\Omega$.}
\end{split}
\end{align}
Indeed, if $v\in{\mathscr{V}}(\Omega)$ then formula \eqref{eNjyg-u51} written for $f:=v$ and $F:=\widetilde{v}$ 
implies that $\eta=0$. In turn, the unique solution of the Neumann problem \eqref{eNN2-Q} for $\eta=0$, 
$f:=v$, and $F:=\widetilde{v}$ is $\vartheta=0$. Having established that $\eta=\vartheta=0$ in this case, the 
conclusion in \eqref{8yhGVV} is seen by appropriately translating \eqref{eNjyg-u58} and \eqref{eNN2-Qe2}.

Next, we shall show that the Neumann trace defined in \eqref{eNN2-Qe2} is compatible with the 
Neumann traces from Case~1. To this end, assume that one is given a function $f\in H^{s}(\Omega)$ 
with $s\in\big(\tfrac{1}{2},\tfrac{3}{2}\big)$ along with some $F\in H^{s-2}_0(\Omega)$ 
satisfying $\Delta f=F\big|_{\Omega}$ in ${\mathcal{D}}'(\Omega)$. Then all conditions 
in \eqref{u4338jh} hold if one chooses 
\begin{equation}\label{eq:EEE.WACO}
\varepsilon:=s-(1/2)\in(0,1). 
\end{equation}
Next, given any function $\phi\in H^{1}(\partial\Omega)\subset H^{(3/2)-s}(\partial\Omega)$, let 
$v$ be as in \eqref{eNjyg-u58}, and take $w$ to be the unique function in 
${\mathscr{W}}(\Omega)\subset H^{3/2}(\Omega)\subset H^{2-s}(\Omega)$ such that 
$\gamma_D w=\phi$. Bear in mind that the mere membership of $w$ to ${\mathscr{W}}(\Omega)$ 
entails $\Delta w=w=\widetilde{w}\big|_{\Omega}$ (where tilde denotes the extension by 
zero outside $\Omega$) and $w\in H^{3/2}(\Omega)\subset H^{2-s}(\Omega)$ 
(in particular, $\widetilde{w}\in L^2({\mathbb{R}}^n)$).
Then \eqref{eNjyg-u58} forces
\begin{align}\label{2.NNa-P}
\big(\gamma_D v,\,\nu\cdot\gamma_D(\nabla w)\big)_{L^2(\partial\Omega)}=I-II,
\end{align}
where, by \eqref{eq:Nnan7yg-P} and Green's formula \eqref{GGGRRR-prim}, 
\begin{align}\label{2.NNa.1-P}
I &:=\big(\gamma_D f,\nu\cdot\gamma_D(\nabla w)\big)_{L^2(\partial\Omega)}
=\big(\gamma_D f,\widetilde\gamma_N(w,\widetilde{w})\big)_{L^2(\partial\Omega)}
\nonumber\\[2pt]
&\,\,={}_{H^{s-(1/2)}(\partial\Omega)}\big\langle\gamma_D f,\widetilde\gamma_N(w,\widetilde{w})
\big\rangle_{(H^{s-(1/2)}(\partial\Omega))^*}
\nonumber\\[2pt]
&\,\,={}_{(H^{(3/2)-s}(\partial\Omega))^*}\big\langle\widetilde\gamma_N(f,F),
\phi\big\rangle_{H^{(3/2)-s}(\partial\Omega)}
\nonumber\\[2pt]
&\quad\,\,+(f,w)_{L^2(\Omega)}
-{}_{(H^{2-s}(\Omega))^*}\langle F,w\rangle_{H^{2-s}(\Omega)},
\end{align}
and where, with 
\begin{align}\label{lubu8uh}
u:=(\eta-\vartheta)\in H^{(1/2)+\varepsilon}(\Omega)=H^{s}(\Omega)
\end{align}
(thanks to the choice of $\varepsilon$ in \eqref{eq:EEE.WACO}), we abbreviated
\begin{align}\label{2.NNa.2-P}
II & :=\big(\gamma_D u,\nu\cdot\gamma_D(\nabla w)\big)_{L^2(\partial\Omega)}
=\big(\gamma_D u,\widetilde\gamma_N(w,\widetilde{w})\big)_{L^2(\partial\Omega)}
\nonumber\\[2pt] 
&\,\,={}_{H^{s-(1/2)}(\partial\Omega)}
\big\langle\gamma_D u,\widetilde\gamma_N(w,\widetilde{w})\big\rangle_{(H^{s-(1/2)}(\partial\Omega))^*}
\nonumber\\[2pt]
&\,\,=(u,w)_{L^2(\Omega)}
-{}_{(H^{2-s}(\Omega))^*}\big\langle F-\widetilde{f}+\widetilde{u},w\big\rangle_{H^{2-s}(\Omega)}
\nonumber\\[2pt]
&\,\,=(f,w)_{L^2(\Omega)}
-{}_{(H^{2-s}(\Omega))^*}\langle F,w\rangle_{H^{2-s}(\Omega)}. 
\end{align}
Here \eqref{eq:Nnan7yg-P} and Green's formula \eqref{GGGRRR-prim}, 
keeping in mind that \eqref{eNjyg-u53} and \eqref{eNN2-Q} yield  
$(-\Delta+1)u=(-\Delta+1)\eta=(-F+\widetilde{f})\big|_{\Omega}$, 
one obtains 
\begin{align}\label{2.NNa.2-Prre-i}
\text{$\Delta u=(F-\widetilde{f}+\widetilde{u})\big|_{\Omega}$,
with $\big(F-\widetilde{f}+\widetilde{u}\big)  
\in H^{-(3/2)+\varepsilon}_0(\Omega)$,} 
\end{align}
and 
\begin{align}\label{2.NNa.2-Prre}
\widetilde\gamma_N(u,F-\widetilde{f}+\widetilde{u})
&=\widetilde\gamma_N(u+\vartheta,F-\widetilde{f}+\widetilde{u}+\widetilde{\vartheta})
-\widetilde\gamma_N(\vartheta,\widetilde{\vartheta})
\nonumber\\[2pt]
&=\widetilde\gamma_N(\eta,F-\widetilde{f}+\widetilde{\eta}\,)
-\widetilde\gamma_N(\eta,F-\widetilde{f}+\widetilde\eta\,)
\nonumber\\[2pt]
&=0,
\end{align}
by \eqref{lubu8uh} and \eqref{eNN2-Q}. Collectively, \eqref{2.NNa-P}, 
\eqref{2.NNa.1-P}, \eqref{2.NNa.2-P}, and \eqref{eq:Jba-hr4ee} prove that, 
with $\widetilde\gamma_N(f,F)$ interpreted in the sense discussed in Case~1,
\begin{align}\label{2.9bab-P}
{}_{(H^{(3/2)-s}(\partial\Omega))^*}\big\langle
\widetilde\gamma_N(f,F),\phi\big\rangle_{H^{(3/2)-s}(\partial\Omega)}
&=\big(\gamma_D v,\nu\cdot\gamma_D(\nabla w)\big)_{L^2(\partial\Omega)}
\nonumber\\[2pt]
&={}_{H^{-1}(\partial\Omega)}\big\langle\gamma^{\mathscr{V}}_N v,\phi\big\rangle_{H^{1}(\partial\Omega)},
\end{align}
which, after unraveling definitions (cf.~\eqref{eNN2-Qe2}), shows the desired compatibility result for the 
two weak Neumann trace operators. Moreover, that the weak Neumann trace operator in the current context 
is surjective is a direct consequence of the fact that $\gamma^{\mathscr{V}}_N$ in \eqref{eq:Gva444rt}
is an isomorphism.

Corresponding to the case $s=\tfrac12$, we shall let the operator $\Upsilon_N$ in \eqref{2.88X-NN-ii}  
act on a given $\psi\in H^{-1}(\partial\Omega)$ according to $\Upsilon_N\psi:=f$, where 
$f\in{\mathscr{V}}(\Omega)$ is the unique function with the property that 
$\gamma^{\mathscr{V}}_N f=\psi$ (cf.\  Lemma~\ref{YTfBBVa}). Then  
\begin{equation}\label{2.88X-NN2-ii-ZZ2}
\widetilde\gamma_N\big(\Upsilon_N\psi,\widetilde{\Delta(\Upsilon_N\psi)}\,\big)
=\widetilde\gamma_N\big(f,\widetilde{\Delta f}\,\big)
=\gamma^{\mathscr{V}}_N f=\psi,
\end{equation}
due to the manner in which we defined the weak Neumann trace operator 
$\widetilde\gamma_N\big(f,F\big)$ with $f$ as above and $F:=\widetilde{\Delta f}$ 
in the present case. Indeed, this is seen from \eqref{eNN2-Qe2} since  
both, $\eta$ in \eqref{eNjyg-u51} and $\vartheta$ in \eqref{eNN2-Q}, now vanish
(given the choice of $F$), hence $v$ in \eqref{eNjyg-u58} is now equal to $f$.
In turn, \eqref{2.88X-NN2-ii-ZZ2} justifies \eqref{2.88X-NN2-ii} in the case when $s=\tfrac12$
(and also proves the surjectivity of the weak Neumann trace operator in the current case).
Since, as seen from the proof of Lemma~\ref{YTfBBVa}, solving  
\begin{equation}\label{j6f4e-y22}
f\in{\mathscr{V}}(\Omega),\quad\gamma^{\mathscr{V}}_N f=\psi\in H^{-1}(\partial\Omega),
\end{equation}
uses the same formalism based on boundary layer potentials employed in the  
treatment of the boundary value problem intervening in \eqref{2.88X-NN-44}, it follows 
that the corresponding solution operators $\Upsilon_N$ are compatible. 
\\[1mm]
{\it Proof of $(ii)$.}
In a first stage we will show that, whenever $s\in\big[\tfrac{1}{2},\tfrac{3}{2}\big]$, 
then for any two functions $f\in H^s(\Omega)$ with $\Delta f\in L^2(\Omega)$ and 
$g\in H^{2-s}(\Omega)$ with $\Delta g\in L^2(\Omega)$ the following Green's formula holds:
\begin{align}\label{GGGRRR-n}
& {}_{H^{(3/2)-s}(\partial\Omega)}\big\langle\gamma_D g,\widetilde\gamma_N(f,\widetilde{\Delta f})
\big\rangle_{(H^{(3/2)-s}(\partial\Omega))^*}
\nonumber\\[2pt]
&\qquad
-{}_{(H^{s-(1/2)}(\partial\Omega))^*}\big\langle\widetilde\gamma_N(g,\widetilde{\Delta g}),\gamma_D f
\big\rangle_{H^{s-(1/2)}(\partial\Omega)}     
\nonumber\\[2pt]
&\quad=(g,\Delta f)_{L^2(\Omega)}-(\Delta g,f)_{L^2(\Omega)},
\end{align}
where $\widetilde{\Delta f},\widetilde{\Delta g}\in L^2({\mathbb{R}}^n)$ denote the extensions 
of $\Delta f,\Delta g\in L^2(\Omega)$ by zero to ${\mathbb{R}}^n$. 

To justify this particular case of formula \eqref{GGGRRR-prim}, 
one invokes Lemma~\ref{Dense-LLLe} in order to find two sequences 
$\{f_j\}_{j\in\bbN},\{g_j\}_{j\in\bbN}\subset C^\infty(\overline{\Omega})$ with the 
property that, as $j\to\infty$, 
\begin{align}\label{eq:Nba55}
\begin{split}
& f_j\to f\,\text{ in }\,H^{s}(\Omega),\quad
\Delta f_j\to\Delta f\,\text{ in }\,L^{2}(\Omega),   
\\[2pt]
& g_j\to g\,\text{ in }\,H^{2-s}(\Omega),\quad
\Delta g_j\to\Delta g\,\text{ in }\,L^{2}(\Omega).
\end{split}
\end{align}
As a consequence of \eqref{eq:Nba55}, the continuity of the boundary 
traces already proved, and \eqref{eq:Nnan7yg-P} one infers that 
\begin{align}\label{eq:JBakyt}
\begin{split}
& \gamma_Df_j\to\gamma_Df\,\text{ in }\,H^{s-(1/2)}(\partial\Omega),
\\[2pt]
& \nu\cdot\gamma_D(\nabla f_j)
=\widetilde\gamma_N(f_j,\widetilde{\Delta f_j})\to\widetilde\gamma_N(f,\widetilde{\Delta f})
\,\text{ in }\,H^{s-(3/2)}(\partial\Omega),    
\\[2pt] 
& \gamma_Dg_j\to\gamma_Dg\,\text{ in }\,H^{(3/2)-s}(\partial\Omega),
\\[2pt]
&\nu\cdot\gamma_D(\nabla g_j)
=\widetilde\gamma_N(g_j,\widetilde{\Delta g_j})\to\widetilde\gamma_N(g,\widetilde{\Delta g})
\,\text{ in }\,H^{(3/2)-s}(\partial\Omega), 
\end{split}
\end{align}
as $j\to\infty$. Now \eqref{GGGRRR-n} written for $f,g$ as above follows from 
\eqref{eq:Nba55}, \eqref{eq:JBakyt}, and the ordinary Green's formula for functions in 
$C^\infty(\overline{\Omega})$ (itself, a consequence of Theorem~\ref{banff-3}), 
via a limiting argument.

Going forward, having fixed some $\varepsilon\in(0,1)$ along with $s\in\big[\tfrac{1}{2},\tfrac{3}{2}\big]$, pick two pairs, 
$(f,F)\in H^s(\Omega)\times H^{s-2+\varepsilon}_0(\Omega)$ such that $\Delta f=F|_{\Omega}$ in $\mathcal{D}'(\Omega)$, 
and $(g,G)\in H^{2-s}(\Omega)\times H^{-s+\varepsilon}_0(\Omega)$ such that $\Delta g=G|_{\Omega}$ 
in $\mathcal{D}'(\Omega)$. The validity of Green's formula \eqref{GGGRRR-prim} for the aforementioned 
pairs when $s\in\big(\tfrac12,\tfrac32\big)$ has been already established in Theorem~\ref{GNT} 
(even in the limiting case $\varepsilon=0$). As such, there remains to treat the situation when
$\varepsilon\in(0,1)$ and $s\in\big\{\tfrac12,\tfrac32\big\}$. Moreover, simple symmetry considerations
actually reduce matters to considering just one of these two extreme values of $s$, say $s=\tfrac12$. 

Corresponding to this choice of the parameter $s$, assume that $\varepsilon\in(0,1)$ and that 
two pairs, $(f,F)\in H^{1/2}(\Omega)\times H^{-(3/2)+\varepsilon}_0(\Omega)$ 
such that $\Delta f=F|_{\Omega}$ in $\mathcal{D}'(\Omega)$, and 
$(g,G)\in H^{3/2}(\Omega)\times H^{-(1/2)+\varepsilon}_0(\Omega)$ such that $\Delta g=G|_{\Omega}$ 
in $\mathcal{D}'(\Omega)$ have been given. Then Lemma~\ref{Dense-LLLe}
ensures the existence of a sequence $\{g_j\}_{j\in\bbN}\subset C^\infty(\overline{\Omega})$ with the 
property that, as $j\to\infty$, 
\begin{align}\label{eq:Nba55-pp}
g_j\to g\,\text{ in }\,H^{3/2}(\Omega),\quad
\Delta g_j\to\Delta g\,\text{ in }\,H^{-(1/2)+\varepsilon}(\Omega).
\end{align}
In particular, the continuity of $\gamma_D$ in \eqref{eqn:gammaDs.2aux} gives
\begin{align}\label{eq:Nba55-pp.i}
\gamma_D g_j\to\gamma_D g\,\text{ in }\,H^{1}(\partial\Omega)\,\text{ as }\,j\to\infty.
\end{align}
In addition, 
\begin{align}\label{eq:Nba55-ww}
\widetilde{\Delta g_j}\to\widetilde{\Delta g}=G\,\text{ in }\,H^{-(1/2)+\varepsilon}_0(\Omega)
\end{align}
by \eqref{eq:12d-yw5} which, by virtue of the continuity of the weak Neumann trace operator,
further implies that
\begin{align}\label{eq:Nba55-ww.2}
\widetilde\gamma_N(g_j,\widetilde{\Delta g_j})\to\widetilde\gamma_N(g,G)
\,\text{ in }\,L^2(\partial\Omega)\,\text{ as }\,j\to\infty.
\end{align}
Next, if $v$ is an in \eqref{eNjyg-u58}, based on \eqref{eNN2-Qe2}, \eqref{8yhGVV}, \eqref{eq:Nba55-pp.i}, 
and \eqref{GGGRRR-n} with ($s=\tfrac12$), one computes  
\begin{align}\label{eq:Nba55-99i.a}
& {}_{H^1(\partial\Omega)}\big\langle\gamma_Dg, \widetilde\gamma_N(f,F)\big\rangle_{H^{-1}(\partial\Omega)}
\nonumber\\[2pt]
& \quad={}_{H^1(\partial\Omega)}\big\langle\gamma_Dg,\gamma^{\mathscr{V}}_N v\big\rangle_{H^{-1}(\partial\Omega)}
\nonumber\\[2pt]
& \quad={}_{H^1(\partial\Omega)}\big\langle\gamma_Dg,\widetilde\gamma_N(v,\widetilde{v}\,)\big\rangle_{H^{-1}(\partial\Omega)}
\\[2pt]
& \quad=\lim_{j\to\infty}
{}_{H^1(\partial\Omega)}\big\langle\gamma_Dg_j,\widetilde\gamma_N(v,\widetilde{\Delta v}\,)\big\rangle_{H^{-1}(\partial\Omega)}
\nonumber\\[2pt]
& \quad=\lim_{j\to\infty}\Big\{\big(\widetilde\gamma_N(g_j,\widetilde{\Delta g_j}),\gamma_D v\big)_{L^2(\partial\Omega)}     
+(g_j,\Delta v)_{L^2(\Omega)}-(\Delta g_j,v)_{L^2(\Omega)}\Big\}.
\nonumber\\[2pt]
& \quad=\lim_{j\to\infty}\Big\{\big(\widetilde\gamma_N(g_j,\widetilde{\Delta g_j}),\gamma_D v\big)_{L^2(\partial\Omega)}     
+(g_j,v)_{L^2(\Omega)}-(\Delta g_j,v)_{L^2(\Omega)}\Big\},
\nonumber
\end{align}
where the third equality and the last equality use the fact that $\Delta v=v$
(given that $v\in{\mathscr{V}}(\Omega)$). On the other hand, from \eqref{eNjyg-u58} and \eqref{lubu8uh} one concludes $v=f-u$, 
hence for each $j\in{\mathbb{N}}$ we have
\begin{align}\label{eq:Nba55-99i.b}
\big(\widetilde\gamma_N(g_j,\widetilde{\Delta g_j}),\gamma_D v\big)_{L^2(\partial\Omega)}
&=\big(\widetilde\gamma_N(g_j,\widetilde{\Delta g_j}),\gamma_D f\big)_{L^2(\partial\Omega)}
\nonumber\\[2pt]
&\quad-\big(\widetilde\gamma_N(g_j,\widetilde{\Delta g_j}),\gamma_D u\big)_{L^2(\partial\Omega)},
\end{align}
since we currently have $\gamma_D f\in L^2(\partial\Omega)$ and $\gamma_D u\in L^2(\partial\Omega)$ by 
\eqref{eqn:gammaDs.2aux}. In addition, \eqref{eq:Nnan7yg-P} and \eqref{GGGRRR-prim}, used here with 
$s=\tfrac32-\varepsilon\in\big(\tfrac12,\tfrac32\big)$, give, on account 
of \eqref{2.NNa.2-Prre-i} and \eqref{2.NNa.2-Prre},
\begin{align}\label{eq:Nba55-99i.c}
& \big(\widetilde\gamma_N(g_j,\widetilde{\Delta g_j}),\gamma_D u\big)_{L^2(\partial\Omega)}
={}_{H^{-\varepsilon}(\partial\Omega)}\Big\langle\widetilde 
\gamma_N(g_j,\widetilde{\Delta g_j}),
\gamma_D u\Big\rangle_{H^{\varepsilon}(\partial\Omega)}
\nonumber\\[2pt]
& \quad={}_{(H^{(1/2)+\varepsilon}(\Omega))^*}\big\langle\Delta g_j,u\big\rangle_{H^{(1/2)+\varepsilon}(\Omega)}
-{}_{H^{(3/2)-\varepsilon}(\Omega)}\big\langle g_j,F-\widetilde{f}+\widetilde{u}\big\rangle_{H^{-(3/2)+\varepsilon}_0(\Omega)}
\nonumber\\[2pt]
& \quad=\big(\Delta g_j,u\big)_{L^2(\Omega)}-{}_{H^{(3/2)-\varepsilon}(\Omega)}\big\langle g_j,
F\big\rangle_{H^{-(3/2)+\varepsilon}_0(\Omega)}+\big(g_j,v\big)_{L^2(\Omega)}.
\end{align}
From \eqref{eq:Nba55-99i.a}--\eqref{eq:Nba55-99i.c} and \eqref{eq:Nnan7yg-P} one then concludes 
(recalling $u+v=f$) that
\begin{align}\label{eq:Nba55-99i.e}
& {}_{H^1(\partial\Omega)}\big\langle\gamma_Dg,\widetilde\gamma_N(f,F)\big\rangle_{H^{-1}(\partial\Omega)}
\nonumber\\[2pt]
& \quad=\lim_{j\to\infty}\Big\{\big(\widetilde\gamma_N(g_j,\widetilde{\Delta g_j}),\gamma_D f\big)_{L^2(\partial\Omega)}     
-\big(\widetilde{\Delta g_j},f\big)_{L^2(\Omega)}
\nonumber\\[2pt]
& \qquad +{}_{H^{(3/2)-\varepsilon}(\Omega)}\big\langle g_j,
F\big\rangle_{H^{-(3/2)+\varepsilon}_0(\Omega)}\Big\}
\nonumber\\[2pt]
& \quad=\lim_{j\to\infty}\Big\{\big(\widetilde\gamma_N(g_j,\widetilde{\Delta g_j}),\gamma_D f\big)_{L^2(\partial\Omega)}     
-{}_{H^{-(1/2)}_0(\Omega)}\big\langle\widetilde{\Delta g_j},f\big\rangle_{H^{1/2}(\Omega)}
\nonumber\\[2pt]
& \qquad +{}_{H^{(3/2)-\varepsilon}(\Omega)}\big\langle g_j,
F\big\rangle_{H^{-(3/2)+\varepsilon}_0(\Omega)}\Big\}
\nonumber\\[2pt]
& \quad=\big(\widetilde\gamma_N(g,G),\gamma_D f\big)_{L^2(\partial\Omega)}
-{}_{H^{-(1/2)}_0(\Omega)}\big\langle G,f\big\rangle_{H^{1/2}(\Omega)}
\nonumber\\[2pt]
& \qquad+{}_{H^{(3/2)-\varepsilon}(\Omega)}\big\langle g,
F\big\rangle_{H^{-(3/2)+\varepsilon}_0(\Omega)},
\end{align}
by \eqref{eq:Nba55-ww.2}, \eqref{eq:Nba55-ww}, and \eqref{eq:Nba55-pp}.
This finishes the proof of the desired version of Green's formula. \\[1mm] 
{\it Proof of $(iii)$.} 
To treat the claim in \eqref{gafvv.6588-P}, we assume that some 
$f\in H^{1/2}(\Omega)$ and $F\in H^{-(3/2)+\varepsilon}_0(\Omega)$ with 
$0<\varepsilon\leq 1$ satisfy $\Delta f=F\big|_{\Omega}$ in ${\mathcal{D}}'(\Omega)$ and 
$\widetilde\gamma_N(f,F)=0$. One recalls from \eqref{eNN2-Qe2} that the latter condition means 
$\gamma^{\mathscr{V}}_N v=0$ in $H^{-1}(\partial\Omega)$, where $v\in{\mathscr{V}}(\Omega)$ 
is given in \eqref{eNjyg-u58}. The fact that the operator \eqref{eq:Gva444rt} is an isomorphism 
then forces $v=0$ which, in light of \eqref{eNjyg-u58}, entails 
\begin{equation}\label{eNjyg-u58-UU}
f=(\eta-\vartheta)\in H^{(1/2)+\varepsilon}(\Omega), 
\end{equation}
given that both memberships, $\eta\in H^{(1/2)+\varepsilon}(\Omega)$ in \eqref{eNjyg-u51}
and $\vartheta\in H^{(1/2)+\varepsilon}(\Omega)$ in \eqref{eNN2-Q}, are valid in the range 
$0<\varepsilon\leq 1$. Finally, the estimate in \eqref{gafvv.6588-P} is a consequence of 
the estimate in \eqref{eNjyg-u51} and \eqref{eNjyg-u57}, both of which continue to hold 
for $0<\varepsilon\leq 1$. This finishes the proof of the claim made in \eqref{gafvv.6588-P}. 
\\[1mm] 
{\it Proof of $(iv)$.} 
As noted earlier, \eqref{2.9MMabn-P} implies \eqref{eq:Nnan7yg-P}.
Finally, the estimate claimed in \eqref{ut543i-Ai.WACO} has been justified in \eqref{ut543i-Ai}. 
\\[1mm] 
{\it Proof of $(v)$.} Working under the assumption that \eqref{ihg65f+AAA} holds, 
consider $f\in H^s(\Omega)$ with $\Delta f\in H^{s-2+\varepsilon}_{z}(\Omega)$. In view of 
\eqref{restr-0.EMB}, there exists $F\in H^{s-2+\varepsilon}_0(\Omega)$ satisfying 
$\Delta f=F\big|_{\Omega}$ in ${\mathcal{D}}'(\Omega)$. This implies that ${\mathcal{I}}(f,F)=f$ 
which, in turn, proves that the mapping ${\mathcal{I}}$ is surjective in the context of \eqref{eqn:gammaN-III.a}.
Obviously, ${\mathcal{I}}$ is linear. To show that ${\mathcal{I}}$ is also injective, assume 
$(f,F)\in H^s(\Omega)\times H^{s-2+\varepsilon}_0(\Omega)$ satisfy $\Delta f=F\big|_{\Omega}$ in 
${\mathcal{D}}'(\Omega)$ and ${\mathcal{I}}(f,F)=0$. The latter implies $f=0$, hence $F\big|_{\Omega}=0$
in ${\mathcal{D}}'(\Omega)$. Since, by design (cf.~\eqref{u54ruy}), one has  
${\rm supp}\,F\subseteq\overline{\Omega}$, and one concludes that 
$F\in H^{s-2+\varepsilon}({\mathbb{R}}^n)$ has 
${\rm supp}\,F\subseteq\partial\Omega$. In view of \eqref{SDErr99-1} and the fact that $s-2+\varepsilon>-\tfrac{1}{2}$ 
(cf.~\eqref{ihg65f+AAA}), one deduces that $F=0$. Ultimately, this proves that ${\mathcal{I}}$ is injective 
in the context of \eqref{eqn:gammaN-III.a}. Since by design ${\mathcal{I}}$ is also bounded, one finally 
concludes that  ${\mathcal{I}}$ is, in fact, a continuous linear isomorphism. 
All other claims readily follow from these facts. 
\end{proof}
%%%%%%%

The next two remarks are designed to clarify the scope of Theorem~\ref{YTfdf.NNN.2-Main}, by further shedding light on the relationship between the weak Neumann trace operator defined in \eqref{2.9NEW} and its ``classical" version. 

%%%%%%%
\begin{remark}\label{r4r-WACO.a}
As in Theorem~\ref{YTfdf.NNN.2-Main}, assume $\Omega\subset\mathbb R^n$ is a bounded Lipschitz domain and denote by $\nu$ the outward unit normal vector to $\Omega$. In this context, suppose some function 
\begin{equation}\label{8qagg-TEXAS.aaa}
f\in H^{s_o}(\Omega)\,\text{ with }\,s_o>3/2 
\end{equation}
has been given. Pick $s\in\big(\tfrac{3}{2},\tfrac{5}{2}\big)$ with $s<s_o$ 
and note that $f\in H^{s_o}(\Omega)\hookrightarrow H^s(\Omega)$, while 
\eqref{eq:DDEEn.3}, \eqref{eq:DDEEn.2}, \eqref{Rdac.WACO}, and \eqref{restr-0.EMB} imply that 
\begin{equation}\label{8qagg-TEXAS}
\Delta f\in H^{s_o-2}(\Omega)\hookrightarrow H^{s-2}(\Omega)=H^{s-2}_{z}(\Omega)
=\big\{u|_\Omega\,\big|\,u\in H^{s-2}_0(\Omega)\big\}.
\end{equation}
In particular, there exists $F\in H^{s-2}_0(\Omega)$ such that $\Delta f=F\big|_{\Omega}$ in ${\mathcal{D}}'(\Omega)$.
Granted these facts, we may invoke \eqref{eq:Nnan7yg-P} to conclude that 
\begin{equation}\label{eq:Nnan7yg-P.TEXAS}
\widetilde\gamma_N(f,F)=\nu\cdot\gamma_D(\nabla f)\in L^2(\partial\Omega) 
\,\text{ with the Dirichlet trace taken as in \eqref{eqn:gammaDs.1}.}
\end{equation}
More directly, one can invoke \eqref{eqn:gammaN-III.e}, with the same effect. 
This discussion may be interpreted as saying that the weak Neumann trace operator $(f,F)\mapsto\widetilde\gamma_N(f,F)$ 
defined in \eqref{2.9NEW} is in fact compatible with the ``classical" Neumann boundary trace operator acting on arbitrary 
functions $f$ as in \eqref{8qagg-TEXAS.aaa} according to $f\mapsto\nu\cdot\gamma_D(\nabla f)$ (with the Dirichlet trace 
understood in the sense of \eqref{eqn:gammaDs.1}). \hfill $\diamond$
\end{remark}
%%%%%%%

%%%%%%%
\begin{remark}\label{r4r-WACO.b}
We wish to emphasize that the weak Neumann trace operator $(f,F)\mapsto\widetilde\gamma_N(f,F)$ 
defined in \eqref{2.9NEW} is a renormalization of the ``classical" Neumann boundary 
trace operator $f\mapsto\nu\cdot\gamma_D(\nabla f)$, which requires $f$ to be more regular 
(say $f\in H^{(3/2)+\varepsilon}(\Omega)$ for some $\varepsilon>0$) than assumed
in Theorem~\ref{YTfdf.NNN.2-Main}, relative to the extension of 
$\Delta f\in H^{s-2}(\Omega)$ to a functional $F$ in the space 
\begin{equation}\label{yree}
\big(H^{2-s}(\Omega)\big)^*=H^{s-2}_0(\Omega)
=\big\{F\in H^{s-2}({\mathbb{R}}^n)\,\big|\,\,{\rm supp}\,F\subseteq\overline{\Omega}\,\big\}. 
\end{equation}
More specifically, suppose that $f\in H^{s}(\Omega)$ with $s\in\big[\tfrac{1}{2},\tfrac{3}{2}\big]$
is such that there exists some $F\in\big(H^{2-s}(\Omega)\big)^*=H^{s-2}_0(\Omega)$ with the 
property that for each $\varphi\in C^\infty_0(\Omega)$ one has  
${}_{{\mathcal{D}}'({\mathbb{R}}^n)}\langle F,\widetilde\varphi\,\rangle_{{\mathcal{D}}({\mathbb{R}}^n)}
={}_{{\mathcal{D}}'(\Om)}\langle\Delta f,\varphi\rangle_{{\mathcal{D}}(\Omega)}$, where 
$\widetilde\varphi$ is the extension of $\varphi$ by zero to $\bbR^n$. Then $F$ is not uniquely determined 
by these qualities (since altering $F$ additively by any distribution in $H^{s-2}_0(\Omega)$ supported
on $\partial\Omega$ also does the job), and the specific choice of such an extension $F$ of $\Delta f$ 
strongly affects the manner in which $\gamma_N(f,F)$ is defined in \eqref{2.9NEW}. \hfill $\diamond$
\end{remark}
%%%%%%%

In applications, the following special case of Theorem~\ref{YTfdf.NNN.2-Main} will play a major role. 

%%%%%%%
\begin{corollary}\label{YTfdf.NNN.2}
Assume that $\Omega\subset\bbR^n$ is a bounded Lipschitz domain with outward unit normal $\nu$.
Then the Neumann trace map, originally defined as $u\mapsto\nu\cdot(\nabla u)|_{\partial\Omega}$ 
for $u\in C^\infty(\overline{\Omega})$, extends uniquely to linear continuous operators  
\begin{equation}\label{eqn:gammaN-pp}
\gamma_N:\big\{u\in H^s(\Omega)\,\big|\,\Delta u\in L^2(\Omega)\big\} 
\rightarrow H^{s-(3/2)}(\partial\Omega),\quad s\in\big[\tfrac{1}{2},\tfrac{3}{2}\big],
\end{equation}
{\rm (}throughout, the space on the left-hand side of \eqref{eqn:gammaN-pp} equipped with the 
natural graph norm $u\mapsto\|u\|_{H^{s}(\Omega)}+\|\Delta u\|_{L^2(\Omega)}${\rm )} that are 
compatible with one another, as well as surjective. In fact, there exist linear and bounded operators
\begin{equation}\label{2.88X-NN-ii-RR}
\Upsilon_N:H^{s-(3/2)}(\partial\Omega)\rightarrow 
\big\{u\in H^s(\Omega)\,\big|\,\Delta u\in L^2(\Omega)\big\},\quad s\in\big[\tfrac12,\tfrac32\big],
\end{equation}
which are compatible with one another and are right-inverses for the Neumann trace, that is, 
\begin{equation}\label{2.88X-NN2-ii-RR}
\gamma_N(\Upsilon_N\psi)=\psi,
\quad\forall\,\psi\in H^{s-(3/2)}(\partial\Omega)\,\text{ with }\,s\in\big[\tfrac12,\tfrac32\big].
\end{equation}

In addition, the following properties are valid:

\begin{enumerate}
\item[(i)] If $s\in\big[\tfrac{1}{2},\tfrac{3}{2}\big]$, 
then for any functions $f\in H^s(\Omega)$ with $\Delta f\in L^2(\Omega)$ and 
$g\in H^{2-s}(\Omega)$ with $\Delta g\in L^2(\Omega)$ the following Green's formula holds:
\begin{align}\label{GGGRRR}
& {}_{H^{(3/2)-s}(\partial\Omega)}\big\langle\gamma_D g,\gamma_N f
\big\rangle_{(H^{(3/2)-s}(\partial\Omega))^*}
\nonumber\\[2pt]
& \qquad -{}_{(H^{s-(1/2)}(\partial\Omega))^*}\big\langle\gamma_N g,\gamma_D f
\big\rangle_{H^{s-(1/2)}(\partial\Omega)}     
\nonumber\\[2pt]
& \quad=(g,\Delta f)_{L^2(\Omega)}-(\Delta g,f)_{L^2(\Omega)}.
\end{align}

\item[(ii)] For each $s\in\big[\tfrac{1}{2},\tfrac{3}{2}\big]$, the null space of 
the Neumann boundary trace operator \eqref{eqn:gammaN-pp} satisfies
\begin{equation}\label{eq:EFFa.111}
{\rm ker}(\gamma_N)\subseteq H^{3/2}(\Omega).
\end{equation}
In fact, the inclusion in \eqref{eq:EFFa.111} is quantitative in the sense that
there exists a constant $C\in(0,\infty)$ with the property that
\begin{align}\label{gafvv.6588}
\begin{split}
& \text{whenever $u\in H^{1/2}(\Omega)$ has $\Delta u\in L^2(\Omega)$ and $\gamma_N u=0$ then}
\\[2pt]
& \quad\text{$u\in H^{3/2}(\Omega)$ and 
$\|u\|_{H^{3/2}(\Omega)}\leq C\big(\|u\|_{L^2(\Omega)}+\|\Delta u\|_{L^2(\Omega)}\big)$.}
\end{split}
\end{align}

\item[(iii)] The following property holds:
\begin{align}\label{eq:Nnan7yg-P.CCC}
\begin{split}
& \text{if $u\in H^{3/2}(\Omega)$ has $\Delta u\in L^2(\Omega)$ then $\gamma_Nu=\nu\cdot\gamma_D(\nabla u)$} 
\\[2pt]
& \quad\text{with the Dirichlet trace taken as in \eqref{eqn:gammaDs.2aux}.}
\end{split}
\end{align}
\end{enumerate}
\end{corollary}
%%%%%%%
\begin{proof}
The key is establishing a relationship between the weak Neumann trace operator 
from Theorem~\ref{YTfdf.NNN.2-Main} and the present Neumann trace operator. 
To accomplish this, assume some $s\in\big[\tfrac{1}{2},\tfrac{3}{2}\big]$ has been given 
and choose $0<\varepsilon<\min\{1,2-s\}$. If one denotes by 
\begin{align}\label{2.10X}
&\iota:\big\{u\in H^s(\Omega)\,\big|\,\Delta u\in L^2(\Omega)\big\}\rightarrow
\nonumber\\[2pt]
& \quad\;\big\{(f,F)\in H^s(\Omega)\times H^{s-2+\varepsilon}_0(\Omega)\,\big|\,\Delta f=F\big|_{\Omega}
\,\text{ in }\,{\mathcal{D}}'(\Omega)\big\}
\end{align}
the continuous injection given by 
\begin{equation}\label{2.11X}
\iota(u):=(u,\widetilde{\Delta u}),\quad\forall\,u\in H^s(\Omega)\,\text{ with }\,\Delta u\in L^2(\Omega),
\end{equation}
(as usual, tilde denotes the extension by zero outside $\Omega$), then 
\begin{equation}\label{2.12X}
\gamma_N:=\widetilde\gamma_N\circ\iota
\end{equation}
yields a well defined, linear, and bounded mapping in the context of \eqref{eqn:gammaN-pp}.
To illustrate the manner in which $\gamma_N$ operates, consider the case where  
$s\in\big(\tfrac{1}{2},\tfrac{3}{2}\big)$. Then, given $u\in H^s(\Omega)$ 
with $\Delta u\in L^2(\Omega)$, along with $\phi\in H^{(3/2)-s}(\partial\Omega)$
and $\Phi\in H^{2-s}(\Omega)$ such that $\gamma_D\Phi=\phi$, then the action of 
$\gamma_N u\in H^{s-(3/2)}(\partial\Omega)=\big(H^{(3/2)-s}(\partial\Omega)\big)^*$ on 
$\phi\in H^{(3/2)-s}(\partial\Omega)$ is concretely given by 
\begin{align}\label{2.9}
& {}_{H^{(3/2)-s}(\partial\Omega)}\big\langle\phi,\gamma_N u\big\rangle_{(H^{(3/2)-s}(\partial\Omega))^*}
\nonumber\\[2pt]
& \quad={}_{H^{(3/2)-s}(\partial\Omega)}\big\langle\phi,
\widetilde\gamma_N(u,\widetilde{\Delta u})\big\rangle_{(H^{(3/2)-s}(\partial\Omega))^*}
\nonumber\\[2pt]
& \quad=\sum_{j=1}^n {}_{H^{1-s}(\Omega)}\big\langle\partial_j\Phi,\partial_j u\big\rangle_{(H^{1-s}(\Omega))^*}
+{}_{H^{2-s}(\Omega)}\big\langle\Phi,\widetilde{\Delta u}\big\rangle_{(H^{2-s}(\Omega))^*}
\nonumber\\[2pt]
& \quad=\sum_{j=1}^n {}_{H^{1-s}(\Omega)}\big\langle\partial_j\Phi,\partial_j u\big\rangle_{(H^{1-s}(\Omega))^*}
+(\Phi,\Delta u)_{L^2(\Omega)}.
\end{align}

Next, we remark that retaining the operators $\Upsilon_N$ as in \eqref{2.88X-NN-ii} implies,
in light of \eqref{2.12X}, \eqref{2.11X}, and \eqref{2.88X-NN2-ii},
\begin{align}\label{2.88X-NN2-ii-RR-2}
\begin{split}
& \gamma_N(\Upsilon_N\psi)=\widetilde\gamma_N
\big(\Upsilon_N\psi,\widetilde{\Delta(\Upsilon_N\psi)}\big)=\psi,
\\[2pt]
& \quad\forall\,\psi\in H^{s-(3/2)}(\partial\Omega)
\,\text{ with }\,s\in\big[\tfrac12,\tfrac32\big].
\end{split}
\end{align}
This justifies \eqref{2.88X-NN2-ii-RR} (which also proves the surjectivity of $\gamma_N$ 
in \eqref{eqn:gammaN-pp}). Moreover, from \eqref{2.12X}, \eqref{eq:Nnan7yg-P}, and the 
discussion pertaining to the nature of \eqref{eqn:gammaDs.1}, one  concludes that
\begin{equation}\label{2.12X.uu}
\gamma_N u=\widetilde\gamma_N(u,\widetilde{\Delta u})=\nu\cdot(\nabla u)\big|_{\partial\Omega},
\quad\forall\,u\in C^\infty(\overline{\Omega}),
\end{equation}
proving that, indeed, our $\gamma_N$ is a genuine extension of the classical (strong) Neumann 
trace operator acting on $C^\infty(\overline{\Omega})$. Since by Lemma~\ref{Dense-LLLe}  
the latter space is dense in $\big\{u\in H^s(\Omega)\,\big|\,\Delta u\in L^2(\Omega)\big\}$, it follows that
the said extension is unique. 

Next, \eqref{GGGRRR} is a particular case of the more general Green's formula in \eqref{GGGRRR-prim}. 
In turn, \eqref{eq:EFFa.111} and \eqref{gafvv.6588} are a direct consequence of \eqref{gafvv.6588-P} 
(used here with $\varepsilon=1$), keeping in mind that since 
$L^2({\mathbb{R}}^n)\hookrightarrow H^{-1/2}({\mathbb{R}}^n)$ continuously,  one has 
\begin{equation}\label{u6gygv}
\|\widetilde{\Delta u}\|_{H^{-1/2}({\mathbb{R}}^n)}\leq C\|\widetilde{\Delta u}\|_{L^2({\mathbb{R}}^n)}
=C\|\Delta u\|_{L^2(\Omega)},
\end{equation}
for some constant $C\in(0,\infty)$, independent of $u$. Finally, \eqref{eq:Nnan7yg-P.CCC} is implied by 
\eqref{eq:Nnan7yg-P}.
\end{proof}
%%%%%%%%

%%%%%%%
\begin{remark}\label{WACO.1}
For higher-order Sobolev spaces, characterizations in the spirit of \eqref{incl-Yb.EE}
have been proved in \cite{MMS10} and \cite{MM13}. For us it is useful to know that
\begin{equation}\label{Tan-C3}
\accentset{\circ}{H}^2(\Omega)=\big\{f\in H^2(\Omega)\,\big|\,\gamma_D f=\gamma_N f=0\big\} 
\end{equation}
for any bounded Lipschitz domain $\Omega\subset\bbR^n$. \hfill $\diamond$ 
\end{remark}
%%%%%%%

The result discussed in the remark below answers a question posed to us by Selim Sukhtaiev.

%%%%%%%
\begin{remark}\label{WACO.2}
Given an arbitrary bounded Lipschitz domain $\Omega\subset\bbR^n$, abbreviate $H^1_{\Delta}(\Omega):=H^{1,0}_{\Delta}(\Omega)$
(where the latter space is as in \eqref{gag-8gb.222} with $s_1:=1$ and $s_2:=0$), that is, define 
\begin{equation}\label{gag-8gb.222.WACO}
H^1_{\Delta}(\Omega):=\big\{u\in H^1(\Omega)\,\big|\,\Delta u\in L^2(\Omega)\big\}
\end{equation}
equipped with the natural graph norm $u\mapsto\|u\|_{H^1(\Omega)}+\|\Delta u\|_{L^2(\Omega)}$.
Since Corollary~\ref{YTfdf-T.NNN} and Corollary~\ref{YTfdf.NNN.2} guarantee that the trace maps 
\begin{align}\label{eqn:gammaDs.2.WACO}
&\gamma_D:H^1_\Delta(\Omega)\rightarrow H^{1/2}(\partial\Omega),
\\[4pt]
&\gamma_N:H^1_\Delta(\Omega)\rightarrow H^{-1/2}(\partial\Omega), 
\label{eqn:gammaN-pp.WACO}
\end{align}
are well defined, linear, and continuous, it follows that the joint trace map 
\begin{align}\label{eqn:gamma.WACO}
\begin{array}{c}
\gamma_{(D,N)}:H^1_\Delta(\Omega)\rightarrow H^{1/2}(\partial\Omega)\times H^{-1/2}(\partial\Omega),
\\[4pt]
\gamma_{(D,N)}u:=\big(\gamma_Du,\gamma_Nu)\,\text{ for each }\,u\in H^1_\Delta(\Omega),
\end{array}
\end{align}
is also well defined, linear, and continuous. However, while Corollary~\ref{YTfdf-T.NNN} and 
Corollary~\ref{YTfdf.NNN.2} imply that the individual Dirichlet and Neumann trace maps from 
\eqref{eqn:gammaDs.2.WACO}--\eqref{eqn:gammaN-pp.WACO} are surjective, we claim that the joint 
trace map \eqref{eqn:gamma.WACO} fails to be surjective. 

To justify this claim, observe that any function $u\in H^1_\Delta(\Omega)$ is uniquely determined by 
$f:=(-\Delta+1)u\in L^2(\Omega)$ and $\phi:=\gamma_Du\in H^{1/2}(\partial\Omega)$. Indeed, from \cite{MT00}
we know that for each given $f\in L^2(\Omega)$ and $\phi\in H^{1/2}(\partial\Omega)$
the inhomogeneous Dirichlet problem 
\begin{equation}\label{eqn:bvp-REG.222-Mi.WACO}
\begin{cases}
(-\Delta+1)u=f\,\text{ in }\,\Omega,\quad u\in H^1(\Omega),     
\\[2pt]  
\gamma_D u=\phi\,\text{ on }\,\partial\Omega,
\end{cases}    
\end{equation}
has a unique solution, which is actually given by 
\begin{align}\label{eqn:gamma.WACO.2}
u=\Pi f+{\mathscr{S}}\Big(S^{-1}\big(\phi-\gamma_D(\Pi f)\big)\Big)\,\text{ in }\,\Omega.
\end{align}
Above, with the fundamental solution $E_1$ as in \eqref{trree.P}, 
\begin{equation}\label{eqn.bweyq-MMM.WACO}
\Pi:\begin{cases} 
L^2(\Omega)\rightarrow H^1(\Omega), 
\\
L^2(\Omega) \ni h \mapsto (\Pi h)(x):=\int_\Omega E_1(x-y)h(y)\,d^ny,\quad x\in\Omega, 
\end{cases} 
\end{equation}
is the volume (Newtonian) potential operator in $\Omega$, while 
\begin{align}\label{eq:GCa.3.WACO.aaa}
& {\mathscr{S}}:H^{-1/2}(\partial\Omega)\rightarrow H^1(\Omega),
\\[4pt]
& S:H^{-1/2}(\partial\Omega)\rightarrow H^{1/2}(\partial\Omega),
\label{eq:GCa.3.WACO.bbb}
\end{align}
are, respectively, the boundary-to-domain single layer potential operator and the boundary-to-boundary 
single layer potential operator associated with the Helmholtz operator $-\Delta+1$ in $\Omega$ 
(cf.~\eqref{eq:GCa.1}--\eqref{eq:GCa.2}). As a consequence of work in \cite{MT00}, these operators 
are well defined, linear, and continuous in each of the indicated contexts. Moreover, $\Pi$ in 
\eqref{eqn.bweyq-MMM.WACO} is actually compact, as 
\begin{equation}\label{eq:eqrqrq.WACO.mmm}
\begin{array}{c}
\text{$\Pi$ maps $L^2(\Omega)$ continuously into $H^2(\Omega)$},
\\[4pt]
\text{which further embeds compactly into $H^1(\Omega)$},
\end{array}
\end{equation}
and $S$ in \eqref{eq:GCa.3.WACO.bbb} is actually an isomorphism. Hence, $u$ in \eqref{eqn:gamma.WACO.2} is well defined 
and, given that 
\begin{equation}\label{eq:ss-WACO}
\gamma_D{\mathscr{S}}=S\,\,\text{ in the setting of \eqref{eq:GCa.3.WACO.aaa}--\eqref{eq:GCa.3.WACO.bbb}},
\end{equation}
it can be checked without difficulty that the function $u$ satisfies \eqref{eqn:bvp-REG.222-Mi.WACO}. 

In light of this discussion, the issue whether the joint trace $\gamma_{(D,N)}$ in \eqref{eqn:gamma.WACO} is surjective boils 
down to the following question: Given an arbitrary $\phi\in H^{1/2}(\partial\Omega)$ along with an arbitrary 
$\psi\in H^{-1/2}(\partial\Omega)$, is it possible to find some $f\in L^2(\Omega)$ with the property that 
$u$ defined as in \eqref{eqn:gamma.WACO.2} satisfies $\gamma_N u=\psi$\,? 

To better understand the latter property we bring in the double layer potential operator, 
originally introduced in \eqref{eq:GCa.7}, presently considered in the context 
\begin{align}\label{eq:GCa.3.TEXAS}
K:H^{1/2}(\partial\Omega)\rightarrow H^{1/2}(\partial\Omega).
\end{align}
Work in \cite{MT00} guarantees that this is well defined, linear, bounded, 
and (with $I$ denoting the identity) satisfies
\begin{equation}\label{eq:WACO.texas}
\gamma_N{\mathscr{S}}=-\tfrac{1}{2}I+K^*\,\text{ as operators on }\,H^{-1/2}(\partial\Omega).
\end{equation}
Bearing these properties in mind, having $\gamma_N u=\psi$ then comes down to solving 
\begin{align}\label{eqn:gamma.WACO.3}
\gamma_N(\Pi f)+\big(-\tfrac{1}{2}I+K^*\big)\Big(S^{-1}\big(\phi-\gamma_D(\Pi f)\big)\Big)=\psi
\end{align}
or, equivalently, 
\begin{align}\label{eqn:gamma.WACO.4}
Tf=\eta, 
\end{align}
where
\begin{align}\label{eqn:gamma.WACO.5}
Tf:=\gamma_N(\Pi f)-\big(-\tfrac{1}{2}I+K^*\big)\Big(S^{-1}\big(\gamma_D(\Pi f)\big)\Big)
\end{align}
and
\begin{align}\label{eqn:gamma.WACO.6}
\eta:=\psi-\big(-\tfrac{1}{2}I+K^*\big)\big(S^{-1}\phi\big).
\end{align}
In view of the compactness of \eqref{eqn.bweyq-MMM.WACO} and the mapping 
properties of $\gamma_N$, $\gamma_D$, $K^\ast$, $S^{-1}$, it follows that 
\begin{equation}\label{eq:waco.tx}
T:L^2(\Omega)\rightarrow H^{-1/2}(\partial\Omega)
\end{equation}
is a linear compact operator. We also note that as $\phi$ and $\psi$ range freely in $H^{1/2}(\partial\Omega)$ 
and $H^{-1/2}(\partial\Omega)$, respectively, $\eta$ can become any function in $H^{-1/2}(\partial\Omega)$. 
Granted this observation, the ability of solving \eqref{eqn:gamma.WACO.4} hinges on whether the operator 
\eqref{eq:waco.tx} is also surjective, which would contradict its compactness. Specifically, if $T$ were 
surjective, the Open Mapping Theorem would imply that $T$ is open. Hence, if $B_{L^2(\Omega)}$ and 
$B_{H^{-1/2}(\partial\Omega)}$ denote the unit balls in $L^2(\Omega)$ and $H^{-1/2}(\partial\Omega)$, 
respectively, we would conclude that there exists $c\in(0,\infty)$ such that 
\begin{equation}\label{eq:agf-WACO}
c\,B_{H^{-1/2}(\partial\Omega)}\subseteq T\big(B_{L^2(\Omega)}\big).
\end{equation}
Given that $T\big(B_{L^2(\Omega)}\big)$ is relatively compact in $H^{-1/2}(\partial\Omega)$, we would then be able 
to conclude that $B_{H^{-1/2}(\partial\Omega)}$ is a relatively compact set. However, according to Riesz's Theorem 
this would further force $H^{-1/2}(\partial\Omega)$ to be a finite-dimensional space, which is certainly not the case. 
The contradiction just reached ultimately proves that the joint trace map \eqref{eqn:gamma.WACO} is not surjective. \hfill $\diamond$
\end{remark}
%%%%%%%

\begin{corollary}\label{MMM.WACO.c}
Let $\Omega\subset\bbR^n$ be an arbitrary bounded Lipschitz domain, and recall the space $H^1_{\Delta}(\Omega)$ 
defined in \eqref{gag-8gb.222.WACO}. Then $H^1_{\Delta}(\Omega)\cap\accentset{\circ}{H}^1(\Omega)$ becomes a 
Banach space when equipped with the norm
\begin{equation}\label{eq:eaMab.WACO}
H^1_{\Delta}(\Omega)\cap\accentset{\circ}{H}^1(\Omega)\ni u\mapsto\|u\|_{H^1(\Omega)}+\|\Delta u\|_{L^2(\Omega)},
\end{equation}
and the Neumann trace map \eqref{eqn:gammaN-pp} induces a well defined, linear, compact operator, in the context 
\begin{equation}\label{eqn:gammaN-pp.WACO.yall}
\gamma_N:H^1_{\Delta}(\Omega)\cap\accentset{\circ}{H}^1(\Omega)\rightarrow L^2(\partial\Omega),
\end{equation}
when the space in the left-hand side is equipped with the norm \eqref{eq:eaMab.WACO}. As a corollary, 
\begin{equation}\label{eqn:gammaN-pp.WACO.yall.222}
\gamma_N:H^1_{\Delta}(\Omega)\cap\accentset{\circ}{H}^1(\Omega)\rightarrow L^2(\partial\Omega)
\,\,\text{ is not surjective}.
\end{equation}
\end{corollary}

\begin{proof}
To justify the first claim, suppose $\{u_j\}_{j\in\bbN}$ is a Cauchy sequence in the space 
$H^1_{\Delta}(\Omega)\cap\accentset{\circ}{H}^1(\Omega)$, equipped with the norm \eqref{eq:eaMab.WACO}.
Then $\{u_j\}_{j\in\bbN}$ is Cauchy in $H^1(\Omega)$ and $\{\Delta u_j\}_{j\in\bbN}$ is Cauchy in $L^2(\Omega)$.
Given that the latter spaces are complete, we conclude that there exist $u\in H^1(\Omega)$ along with $v\in L^2(\Omega)$ 
such that, as $j\to\infty$, 
\begin{align}\label{eq:Nba55.WACO.m}
u_j\to u\,\text{ in }\,H^1(\Omega)\,\,\text{ and }\,\,\Delta u_j\to v\,\text{ in }\,L^{2}(\Omega).
\end{align}
Then, as a consequence of \eqref{eq:Nba55.WACO.m} and the continuity of the Dirichlet trace map \eqref{eqn:gammaDs.2},
$0=\gamma_D u_j\to\gamma_D u\,\text{ in }\,H^{1/2}(\partial\Omega)$ as $j\to\infty$. Hence, $\gamma_D u=0$ which places 
$u$ in $\accentset{\circ}{H}^1(\Omega)$ (cf. \eqref{incl-Yb.EE}). In addition, \eqref{eq:Nba55.WACO.m} implies that
$u_j\to u$ in ${\mathcal{D}}'(\Omega)$ as $j\to\infty$, hence also $\Delta u_j\to\Delta u$ in ${\mathcal{D}}'(\Omega)$ as $j\to\infty$, 
and $\Delta u_j\to v$ in ${\mathcal{D}}'(\Omega)$ as $j\to\infty$. In view of the fact that ${\mathcal{D}}'(\Omega)$ is a 
Hausdorff topological space, these properties force $\Delta u=v\in L^2(\Omega)$, hence $u$ belongs to $H^1_{\Delta}(\Omega)$ as well. 
As such, $u\in H^1_{\Delta}(\Omega)\cap\accentset{\circ}{H}^1(\Omega)$ and, as seen from \eqref{eq:Nba55.WACO.m}, the 
sequence $\{u_j\}_{j\in\bbN}$ converges to $u$ (with respect to the norm \eqref{eq:eaMab.WACO}). This finishes the proof 
of the fact that $H^1_{\Delta}(\Omega)\cap\accentset{\circ}{H}^1(\Omega)$ is a Banach space when endowed with 
the norm \eqref{eq:eaMab.WACO}.

Let us now deal with the second claim, pertaining to the well definiteness, linearity, and compactness of 
\eqref{eqn:gammaN-pp.WACO.yall}. To establish that this Neumann trace is a well defined linear map we first observe 
from \eqref{eq:EFFa.NNN} and \eqref{incl-Yb.EE} that
\begin{equation}\label{eqn:gammaN-pp.WACO.yall.X}
H^1_{\Delta}(\Omega)\cap\accentset{\circ}{H}^1(\Omega)\subseteq H^{3/2}(\Omega).
\end{equation}
Granted this, \eqref{eqn:gammaN-pp} with $s:=\tfrac{3}{2}$ gives that $\gamma_N$ is indeed a well defined linear 
map in the context of \eqref{eqn:gammaN-pp.WACO.yall}. Next we shall prove that said map is also compact. 
To justify this, we shall freely borrow results from, and notation employed in, Remark~\ref{WACO.2}. 
To get started, define the map
\begin{equation}\label{eqn.bweyq-MMM.WACO.123}
\Theta: 
\begin{cases} 
L^2(\Omega)\rightarrow H^1_{\Delta}(\Omega)\cap\accentset{\circ}{H}^1(\Omega), 
\\[6pt]
L^2(\Omega)\ni f\mapsto\Theta f:=\Pi f-{\mathscr{S}}\Big(S^{-1}\big(\gamma_D(\Pi f)\big)\Big).
\end{cases} 
\end{equation}
That this is well defined, linear, and bounded, follows from \eqref{eq:eqrqrq.WACO.mmm} and the 
discussion in the proof of Lemma~\ref{YTfBBVa} where, among other things, it was pointed out that
\begin{align}\label{eq:GCa.3.WACO.aaa.M}
& {\mathscr{S}}:L^2(\partial\Omega)\rightarrow H^{3/2}(\Omega)\,\,\text{ boundedly},
\\[4pt]
& S:L^2(\partial\Omega)\rightarrow H^1(\partial\Omega)\,\,\text{ isomorphically}.
\label{eq:GCa.3.WACO.bbb.M}
\end{align}

We claim that $\Theta$ is actually an isomorphism in the context of \eqref{eqn.bweyq-MMM.WACO.123}.
To justify that $\Theta$ is injective, let $f\in L^2(\Omega)$ be such that $\Theta f=0$. Then 
$0=(-\Delta+1)\Theta f=f$, as wanted. The surjectivity of $\Theta$ follows from the observation that, for each given 
$f\in L^2(\Omega)$, the boundary value problem \eqref{eqn:bvp-REG.222-Mi.WACO} written with $\phi:=0$ has 
a unique solution, which is actually given by \eqref{eqn:gamma.WACO.2} with $\phi:=0$, which is precisely $\Theta f$. 
Hence, $\Theta$ is an isomorphism and, according to the Open Mapping Theorem (whose applicability is ensured by 
the completeness result established in the first part of the proof), $\Theta^{-1}$ is linear and bounded.

Consequently, proving the compactness of $\gamma_N$ in the context of \eqref{eqn:gammaN-pp.WACO.yall} is equivalent to 
showing that 
\begin{equation}\label{eqn:gammaN-pp.WACO.yall.XYZ}
Q:=\gamma_N\circ\Theta:L^2(\Omega)\rightarrow L^2(\partial\Omega)\,\,\text{ is compact}.
\end{equation}
Denote by $\nu$ the outward unit normal vector to $\Omega$.
From \eqref{eqn.bweyq-MMM.WACO.123}, \eqref{eq:eqrqrq.WACO.mmm}, \eqref{eq:Nnan7yg-P}, and \eqref{eq:WACO.texas}
we then see that for each $f\in L^2(\Omega)$ we have
\begin{equation}\label{eqn.bweyq-MMM.WACO.123.QQQ.1}
Qf=\nu\cdot\gamma_D(\nabla\Pi f)-\big(-\tfrac{1}{2}I+K^\ast\big)\Big(S^{-1}\big(\gamma_D(\Pi f)\big)\Big).
\end{equation}
Since the assignment $L^2(\Omega)\ni f\mapsto\gamma_D(\nabla\Pi f)\in H^{1/2}(\partial\Omega)$ is bounded, 
and the embedding $H^{1/2}(\partial\Omega)\hookrightarrow L^2(\partial\Omega)$ is compact, it follows that
\begin{equation}\label{eq:fafafa.qqE.1}
L^2(\Omega)\ni f\mapsto\gamma_D(\nabla\Pi f)\in L^2(\partial\Omega)\,\,\text{ is compact}.
\end{equation}
Also, bearing in mind that the Newtonian potential operator $\Pi$ maps $L^2(\Omega)$ continuously into 
$H^2(\Omega)$ which, for each fixed $\varepsilon\in\big(0,\tfrac{1}{2}\big)$, further embeds compactly 
into the space $\big\{u\in H^{3/2}(\Omega):\,\Delta u\in H^{-(1/2)+\varepsilon}(\Omega)\big\}$ 
(equipped with the natural graph norm), we conclude from \eqref{eqn:gammaDs.2aux}, used with $s:=\tfrac{3}{2}$,  
that the assignment 
\begin{equation}\label{eq:fafafa.qqE.2}
L^2(\Omega)\ni f\mapsto\gamma_D(\Pi f)\in H^1(\partial\Omega)\,\,\text{ is compact}.
\end{equation}
Collectively, \eqref{eqn.bweyq-MMM.WACO.123.QQQ.1}, \eqref{eq:fafafa.qqE.1}, \eqref{eq:fafafa.qqE.2}, \eqref{eq:GCa.3.WACO.bbb.M}, 
and the fact that $K$ is a well defined and bounded operator on $L^2(\partial\Omega)$ then prove that the operator 
\eqref{eqn:gammaN-pp.WACO.yall.XYZ} is indeed compact.

At this stage, there remains to justify the claim made in \eqref{eqn:gammaN-pp.WACO.yall.222}
For this, we reason by contradiction, as in the last part of Remark~\ref{WACO.2} with natural alterations. 
Specifically, if $Q$ were surjective, the Open Mapping Theorem would imply that $Q$ is open. As such, 
if $B_{L^2(\Omega)}$ and $B_{L^2(\partial\Omega)}$ denote, respectively, the unit balls in $L^2(\Omega)$ and $L^2(\partial\Omega)$, 
we would conclude that there exists some constant $c\in(0,\infty)$ with the property that 
\begin{equation}\label{eq:agf-WACO.XXX}
c\,B_{L^2(\partial\Omega)}\subseteq Q\big(B_{L^2(\Omega)}\big).
\end{equation}
Since $Q\big(B_{L^2(\Omega)}\big)$ is relatively compact in $L^2(\partial\Omega)$, we would then be able 
to conclude that $B_{L^2(\partial\Omega)}$ is a relatively compact set in $L^2(\partial\Omega)$. However, according to 
Riesz's Theorem this would further force $L^2(\partial\Omega)$ to be a finite-dimensional space, which is clearly not the case. 
This contradiction ultimately establishes \eqref{eqn:gammaN-pp.WACO.yall.222}.
\end{proof}

We conclude this section by establishing the counterpart of Corollary~\ref{MMM.WACO.c} for the Dirichlet trace map.

\begin{corollary}\label{MMM.WACO.c.D}
Let $\Omega\subset\bbR^n$ be an arbitrary bounded Lipschitz domain, and recall the space $H^1_{\Delta}(\Omega)$ defined in 
\eqref{gag-8gb.222.WACO}. Then $\big\{u\in H^1_{\Delta}(\Omega)|\,\gamma_Nu=0\big\}$ becomes a Banach space when equipped 
with the norm inherited from $H^1_{\Delta}(\Omega)$, and the Dirichlet trace map \eqref{eqn:gammaDs.2} induces 
a well defined, linear, compact operator, in the context 
\begin{equation}\label{eqn:gammaN-pp.WACO.yall.D}
\gamma_D:\big\{u\in H^1_{\Delta}(\Omega)|\,\gamma_Nu=0\big\}\rightarrow H^1(\partial\Omega).
\end{equation}
As a consequence, 
\begin{equation}\label{eqn:gammaN-pp.WACO.yall.222.D}
\gamma_D:\big\{u\in H^1_{\Delta}(\Omega)|\,\gamma_Nu=0\big\}\rightarrow H^1(\partial\Omega)
\,\,\text{ is not surjective}.
\end{equation}
\end{corollary}

\begin{proof}
Lemma~\ref{Dense-LLLe} tells us that $H^1_{\Delta}(\Omega)$ is a Banach space, while from Corollary~\ref{YTfdf.NNN.2}
we know that 
\begin{equation}\label{eqn:gammaN-pp.MmMm}
\begin{array}{c}
\gamma_N:H^1_\Delta(\Omega)\rightarrow H^{-1/2}(\partial\Omega)\,\,\text{ is well defined, linear, bounded,}
\\[4pt]
\text{and }\,\,{\rm ker}(\gamma_N)=\big\{u\in H^1_{\Delta}(\Omega)|\,\gamma_Nu=0\big\}\subseteq H^{3/2}(\Omega).
\end{array}
\end{equation}
Together, these properties allow us to conclude that $\big\{u\in H^1_{\Delta}(\Omega)|\,\gamma_Nu=0\big\}$ is a closed 
subspace of $H^1_{\Delta}(\Omega)$, hence a Banach space itself when equipped with the norm inherited from $H^1_{\Delta}(\Omega)$.

Consider next the claim regarding the well definiteness, linearity, and compactness of \eqref{eqn:gammaN-pp.WACO.yall.D}. 
Granted the inclusion in \eqref{eqn:gammaN-pp.MmMm}, from \eqref{eqn:gammaDs.2} with $s:=3/2$ 
we conclude that $\gamma_D$ is a well defined linear map in the context of \eqref{eqn:gammaN-pp.WACO.yall.D}.
Let us now show that this map is also compact. To justify this, we shall freely borrow results and notation 
from Remark~\ref{WACO.2} and Corollary~\ref{MMM.WACO.c}. We begin by defining 
\begin{equation}\label{eqn.bweyq-MMM.WACO.123.D}
\Psi: 
\begin{cases} 
L^2(\Omega)\rightarrow\big\{u\in H^1_{\Delta}(\Omega)|\,\gamma_Nu=0\big\}, 
\\[6pt]
L^2(\Omega)\ni f\mapsto\Psi f:=\Pi f-{\mathscr{S}}\Big(\big(-\tfrac{1}{2}I+K^\ast\big)^{-1}\big(\gamma_N(\Pi f)\big)\Big).
\end{cases} 
\end{equation}
That this is well defined, linear, and bounded, follows from \eqref{eq:eqrqrq.WACO.mmm} and the 
discussion in the proof of Lemma~\ref{YTfBBVa} where it was noted that
\begin{align}\label{eq:GCa.3.WACO.aaa.M.D}
& {\mathscr{S}}:L^2(\partial\Omega)\rightarrow H^{3/2}(\Omega)\,\,\text{ boundedly},
\\[4pt]
& -\tfrac{1}{2}I+K^\ast:L^2(\partial\Omega)\rightarrow L^2(\partial\Omega)\,\,\text{ isomorphically}.
\label{eq:GCa.3.WACO.bbb.M.D}
\end{align}

We claim that $\Psi$ is actually an isomorphism in the context of \eqref{eqn.bweyq-MMM.WACO.123.D}.
To see that $\Psi$ is injective, suppose $f\in L^2(\Omega)$ satisfies $\Psi f=0$. Then 
$0=(-\Delta+1)\Psi f=f$, as desired. To show that $\Psi$ is surjective, pick an arbitrary $f\in L^2(\Omega)$.
From \cite{MT00} we know that the inhomogeneous Neumann problem 
\begin{equation}\label{eqn:bvp-REG.222-Mi.WACO.D}
\begin{cases}
(-\Delta+1)u=f\,\text{ in }\,\Omega,\quad u\in H^1(\Omega),     
\\[2pt]  
\gamma_N u=0\,\text{ on }\,\partial\Omega,
\end{cases}    
\end{equation}
has a unique solution, which is actually given by 
\begin{align}\label{eqn:gamma.WACO.2.D}
u=\Pi f-{\mathscr{S}}\Big(\big(-\tfrac{1}{2}I+K^\ast\big)^{-1}\big(\gamma_N(\Pi f)\big)\Big)=\Psi f.
\end{align}
Hence, $\Psi$ is an isomorphism and, according to the Open Mapping Theorem (whose applicability is ensured by 
the completeness result established in the first part of the proof), $\Psi^{-1}$ is linear and bounded.

As a result, proving the compactness of $\gamma_D$ in the context of \eqref{eqn:gammaN-pp.WACO.yall.D} 
becomes equivalent to showing that 
\begin{equation}\label{eqn:gammaN-pp.WACO.yall.XYZ.D}
R:=\gamma_D\circ\Psi:L^2(\Omega)\rightarrow H^1(\partial\Omega)\,\,\text{ is compact}.
\end{equation}
To proceed, denote by $\nu$ the outward unit normal vector to $\Omega$.
From \eqref{eqn.bweyq-MMM.WACO.123.D}, \eqref{eq:eqrqrq.WACO.mmm}, \eqref{eq:Nnan7yg-P}, \eqref{eq:GCa.3}, 
\eqref{eKJa77}, and Lemma~\ref{Hgg-uh} we then see that for each $f\in L^2(\Omega)$ we have
\begin{equation}\label{eqn.bweyq-MMM.WACO.123.QQQ.1.D}
Rf=\gamma_D(\Pi f)-S\Big(\big(-\tfrac{1}{2}I+K^\ast\big)^{-1}\big(\nu\cdot\gamma_D(\nabla\Pi f)\big)\Big).
\end{equation}
As before (cf. \eqref{eq:fafafa.qqE.1}, \eqref{eq:fafafa.qqE.2}), 
\begin{equation}\label{eq:fafafa.qqE.1.D}
L^2(\Omega)\ni f\mapsto\gamma_D(\nabla\Pi f)\in L^2(\partial\Omega)\,\,\text{ is compact},
\end{equation}
and
\begin{equation}\label{eq:fafafa.qqE.2.D}
L^2(\Omega)\ni f\mapsto\gamma_D(\Pi f)\in H^1(\partial\Omega)\,\,\text{ is compact}.
\end{equation}
Gathering \eqref{eqn.bweyq-MMM.WACO.123.QQQ.1.D}, \eqref{eq:fafafa.qqE.1.D}, \eqref{eq:fafafa.qqE.2.D}, 
\eqref{eq:GCa.3.WACO.bbb.M.D}, and \eqref{eq:GCa.6} then establishes \eqref{eqn:gammaN-pp.WACO.yall.XYZ.D}. 
Finally, \eqref{eqn:gammaN-pp.WACO.yall.222.D} is justified by reasoning as in the last part in the 
proof of Corollary~\ref{MMM.WACO.c}.
\end{proof}

%%%%%%%%%%%%%%%%%%%%%%%%%%%%%%
%%%%%%%%%%%%%%%%%%%%%%%%%%%%%%
\section{Schr\"{o}dinger Operators on Open Sets and Bounded Lipschitz Domains}  
\label{s6}
%%%%%%%%%%%%%%%%%%%%%%%%%%%%%%
%%%%%%%%%%%%%%%%%%%%%%%%%%%%%%

This section is devoted to a study of minimal and maximal Schr\"{o}dinger operators on nonempty 
open sets and bounded Lipschitz domains $\Omega\subseteq\bbR^n$. Furthermore, the self-adjoint 
Friedrichs extension and the self-adjoint Dirichlet and Neumann realizations are discussed. 

In the beginning of this section we make the following general assumption.

%%%%%%%%%
\begin{hypothesis}\label{h3.1} 
Let $n\in\bbN\backslash\{1\}$, assume that $\Omega\subseteq\bbR^n$ is a nonempty open set, 
and suppose that $V\in L^\infty(\Omega)$ is real-valued. 
\end{hypothesis}
%%%%%%%%%
In the following we denote the essential infimum of $V\in L^\infty(\Omega)$ by $v_-$, i.e., 
\begin{equation}\label{essinfv}
v_-:=\essinf_{x\in\Omega} V(x).
\end{equation}

We are interested in operator realizations of the differential expression $-\Delta+V$ in 
the Hilbert space $L^2(\Omega)$. We define the \textit{preminimal} realization $A_{p,\Omega}$ 
of $-\Delta+V$ by
\begin{equation}\label{We-Q.1}
A_{p,\Omega}:=-\Delta+V,\quad\dom(A_{p,\Omega}):=C^\infty_0(\Omega).
\end{equation}
Thus, $A_{p,\Omega}$ is a densely defined, symmetric operator in $L^2(\Omega)$, and hence 
closable. Next, the \textit{minimal} realization $A_{min,\Omega}$ of $-\Delta+V$ is defined 
as the closure of $A_{p,\Omega}$ in $L^2(\Omega)$, 
\begin{equation}\label{eqn:Amin1}
A_{min,\Omega}:=\overline{A_{p,\Omega}}.
\end{equation}
It follows that $A_{min,\Omega}$ is a densely defined, closed, symmetric operator in $L^2(\Omega)$. 
Finally, the \textit{maximal} realization $A_{max,\Omega}$ of $-\Delta+V$ is given by
\begin{equation}\label{We-Q.2}
A_{max,\Omega}:=-\Delta+V,\quad
\dom(A_{max,\Omega}):=\big\{f\in L^2(\Omega)\,\big|\,\Delta f\in L^2(\Omega)\big\},    
\end{equation}
where the expression $\Delta f$, $f\in L^2(\Omega)$, is understood in the sense of distributions. 
We mention that the assumption $V\in L^\infty(\Omega)$ in Hypothesis~\ref{h3.1} 
yields that for $f\in L^2(\Omega)$ one has $\Delta f\in L^2(\Omega)$ if and only if 
$-\Delta f+Vf\in L^2(\Omega)$.

Next ,we collect some well-known properties of the operators $A_{p,\Omega}$, $A_{min,\Omega}$, 
and $A_{max,\Omega}$ which follow from a standard distribution-type argument, see, for instance, 
\cite[Section~6.2]{Tr92}. 

%%%%%%
\begin{lemma}\label{l3.2}
Assume Hypothesis~\ref{h3.1}. Let $A_{p,\Omega}$, $A_{min,\Omega}$, and $A_{max,\Omega}$ be as 
introduced above. Then the  operators $A_{min,\Omega}$ and $A_{max,\Omega}$ are adjoints of each 
other, that is, 
\begin{equation}\label{eqn:AminAmax}
A_{min,\Omega}^*=A_{p,\Omega}^*=A_{max,\Omega}\,\text{ and }\, 
A_{min,\Omega}=\overline{A_{p,\Omega}}=A_{max,\Omega}^*,
\end{equation}
and the closed symmetric operator $A_{min,\Omega}$ is semibounded from below by $v_-$, that is, 
\begin{equation}\label{aminsemi}
(A_{min,\Omega}f,f)_{L^2(\Omega)}\geq v_{-}\|f\|^2_{L^2(\Omega)},\quad\forall\,f\in\dom(A_{min,\Omega}).
\end{equation}
\end{lemma}
%%%%%%
\begin{proof} 
The assumption $V\in L^\infty(\Omega)$ implies that $V$ is a bounded operator in $L^2(\Omega)$. 
Thus, the domains and adjoints of $A_{p,\Omega}$, $A_{min,\Omega}$, and $A_{max,\Omega}$ 
do not depend on $V$ and hence one can assume without loss of generality that $V\equiv 0$ in the following.
Since $A_{min,\Omega}$ is the closure of $A_{p,\Omega}$ in $L^2(\Omega)$ their adjoints 
$A_{min,\Omega}^*$ and $A_{p,\Omega}^*$ coincide. We first establish the inclusion 
$A_{p,\Omega}^*\subseteq A_{max,\Omega}$. For this purpose, let $f\in\dom(A_{p,\Omega}^*)$ 
be arbitrary. Then one has $f\in L^2(\Omega)$ and $A_{p,\Omega}^* f\in L^2(\Omega)$, 
hence for each function $\varphi\in C^\infty_0(\Omega)$ one can write 
\begin{equation}\label{ua4ftyoki}
\begin{split}
{}_{\cD'(\Omega)}\big\langle\,\overline{A_{p,\Omega}^* f},\varphi\big\rangle_{\cD(\Omega)} 
&=\big(A_{p,\Omega}^* f,\varphi\big)_{L^2(\Omega)}=\big(f,A_{p,\Omega}\varphi\big)_{L^2(\Omega)}
\\[2pt] 
&=(f,-\Delta\varphi)_{L^2(\Omega)}
={}_{\cD'(\Omega)}\langle\,\overline{f},-\Delta\varphi\rangle_{\cD(\Omega)}
\\[2pt] 
&={}_{\cD'(\Omega)}\langle\,\overline{-\Delta f},\varphi\rangle_{\cD(\Omega)},
\end{split}
\end{equation}
by definition of the adjoint and \eqref{We-Q.1} with $V\equiv 0$. Hence, in the sense of distributions, 
one obtains $-\Delta f=A_{p,\Omega}^* f\in L^2(\Omega)$, thus $f\in\dom(A_{max,\Omega})$
and $A_{p,\Omega}^*f=A_{max,\Omega}f$, implying $A_{p,\Omega}^*\subseteq A_{max,\Omega}$. 
Next, we verify the inclusion $A_{max,\Omega}\subseteq A_{p,\Omega}^*$. Pick some 
$f\in\dom(A_{max,\Omega})$. Then $-\Delta f$, considered in the sense of distributions,
belongs to $L^2(\Omega)$, and one may write  
\begin{equation}\label{iu7gt6}
(-\Delta f,\varphi)_{L^2(\Omega)}=(f,-\Delta\varphi)_{L^2(\Omega)} 
=(f,A_{p,\Omega}\varphi)_{L^2(\Omega)}
\end{equation}
for each $\varphi\in\dom(A_{p,\Omega})=C^\infty_0(\Omega)$. In turn, this implies 
$f\in\dom(A_{p,\Omega}^*)$ and $A_{max,\Omega}f=A_{p,\Omega}^*f$, and hence $A_{max,\Omega}\subseteq A_{p,\Omega}^*$. 
The reasoning so far proves the first equality in \eqref{eqn:AminAmax}. The second equality in 
\eqref{eqn:AminAmax} follows by taking adjoints. 

It remains to show that $A_{min,\Omega}$ is semibounded from below by $v_{-}$.
Since $V\geq v_{-}$, for each $f\in C^\infty_0(\Omega)$, repeated integrations by parts yields
\begin{equation}\label{u633}
((A_{p,\Omega}-v_{-})f,f)_{L^2(\Omega)}=(-\Delta f+(V-v_{-})f,f)_{L^2(\Omega)}
\geq\sum_{j=1}^n\|{\partial_j f}\|^2_{L^2(\Omega)}\geq 0.
\end{equation}
This proves that $A_{p,\Omega}-v_{-}$ is nonnegative, and the same 
holds for the closure $A_{min,\Omega}-v_{-}$, that is, \eqref{aminsemi} holds.
\end{proof}
%%%%%%%%

In the next lemma we consider the minimal realization $A_{min,\Omega}$ in the case that $\Omega$ 
is a bounded open set. For the definition of the Sobolev space $\accentset{\circ}{W}^2(\Omega)$ see \eqref{Rdac.WWW}.

%%%%%%%%
\begin{lemma}\label{l3.3}
Assume Hypothesis~\ref{h3.1} and suppose, in addition, that $\Omega$ is bounded. 
Then the closed symmetric operator $A_{min,\Omega}$ is given by
\begin{equation}\label{eqn:Amin2-BBB}
A_{min,\Omega}=-\Delta+V,\quad\dom(A_{min,\Omega})=\accentset{\circ}{W}^2(\Omega).
\end{equation}
Furthermore, $A_{min,\Omega}-v_-$ is strictly positive and $A_{min,\Omega}$  has 
infinite deficiency indices,
\begin{equation}\label{3.12}
\dim\big(\ker(A_{max,\Omega}-zI)\big)=\dim\big(\ker(A_{max,\Omega}-v_{-})\big)=\infty,
\end{equation}
for all $z\in\mathbb{C}\backslash[v_{-},\infty)$.
\end{lemma}
%%%%%%%%
\begin{proof}
The assumption that $\Omega$ is a bounded nonempty open subset of ${\mathbb{R}}^n$ guarantees 
the classical Poincar\'{e} inequality holds. This readily implies that the norm  
\begin{equation}\label{eqn:Deltanorm}
f\mapsto\biggl(\|f\|_{L^2(\Omega)}^2+\sum_{j,k=1}^n
\|\partial_j\partial_k f\|^2_{L^2(\Omega)}\biggr)^{1/2},\quad\forall\,f\in\accentset{\circ}{W}^2(\Omega),
\end{equation}
is equivalent with the norm $\accentset{\circ}{W}^2(\Omega)$ inherits from $W^2(\Omega)$
(cf., e.g., \cite[Theorem~7.6]{Wl87}). For any fixed $f\in C^\infty_0(\Omega)$, 
successive integrations by parts yield
\begin{align}\label{enemenemu}
\sum_{j,k=1}^n\|\partial_j\partial_k f\|^2_{L^2(\Omega)}
&=\sum_{j,k=1}^n\big(\partial_j\partial_k f,\partial_j\partial_k f\big)_{L^2(\Omega)} 
\nonumber\\[2pt] 
&=\sum_{j,k=1}^n\big(\partial_j^2 f,\partial_k^2 f\big)_{L^2(\Omega)}
=\|\Delta f\|_{L^2(\Omega)}^2,
\end{align}
and, as $ C^\infty_0(\Omega)$ is dense in $\accentset{\circ}{W}^2(\Omega)$, the equality of the most 
extreme terms in \eqref{enemenemu} remains to hold for all $f\in\accentset{\circ}{W}^2(\Omega)$. 
Together with the earlier observation pertaining the nature of \eqref{eqn:Deltanorm}, 
this implies that the graph norm 
\begin{equation}\label{ihgffd}
f\mapsto\big(\|f\|_{L^2(\Omega)}^2+\|\Delta f\|^2_{L^2(\Omega)}\big)^{1/2}, 
\quad\forall\,f\in\accentset{\circ}{W}^2(\Omega),
\end{equation}
is equivalent with the norm $\accentset{\circ}{W}^2(\Omega)$ inherits from $W^2(\Omega)$. 
As such, the closure of $-\Delta|_{ C^\infty_0(\Omega)}$ in $L^2(\Omega)$ is the 
operator $-\Delta$ with domain 
\begin{equation}
\overline{ C^\infty_0(\Omega)}^{W^2(\Omega)}=\accentset{\circ}{W}^2(\Omega).
\end{equation}
As the potential $V$ is bounded, this fact remains valid for $-\Delta+V$, and hence 
\eqref{eqn:Amin2-BBB} follows.

In order to see that $A_{min,\Omega}-v_-$ is strictly positive, one again makes use of 
the classical Poincar\'{e} inequality. This permits one to estimate as in \eqref{u633}, 
\begin{equation}\label{u6333}
((A_{p,\Omega}-v_-)f,f)_{L^2(\Omega)}\geq\sum_{j=1}^n  
\|{\partial_j f}\|^2_{L^2(\Omega)}
\geq c\,\|f\|^2_{L^2(\Omega)}
\end{equation}
for some constant $c>0$ independent of $f$. This proves that $A_{p,\Omega}-v_-$ is strictly positive and hence the same holds for the closure $A_{min,\Omega}-v_-$ of $A_{p,\Omega}-v_-$. 

To show that the deficiency numbers of $A_{min,\Omega}$ equal $\infty$, one can argue as follows: 
First, since relatively bounded perturbations with relative bound strictly less than $1$ leave 
deficiency indices invariant as shown in \cite{BF77}, one can again assume $V\equiv 0$ (and, hence, $v_{-}=0$). 
Next, since the set $\Omega\subset\bbR^n$ is bounded, one can contain $\Omega$ in the Euclidean ball 
$B(0,R)\subset\bbR^n$ centered at $0\in{\mathbb{R}}^n$ and having a sufficiently large radius $R>0$. 
Using spherical coordinates and decomposing $-\Delta$ as well as $L^2\big(B(0,R)\big)$ with respect to angular 
momenta (cf., e.g., \cite[Appendix to Section~X.1]{RS75}), employing $n$-dimensional spherical harmonics, 
proves that $A_{max,B(0,R)}$ has infinite deficiency indices. Restricting the elements of $\ker(A_{max,B(0,R)})$ 
to $\Omega\subset B(0,R)$, and using the fact that by the unique continuation property for harmonic functions on 
an open set (see, e.g., \cite[Theorems~6.25, 6.26]{Mi13}), arbitrary finite linear combinations of linearly 
independent harmonic functions on $B(0,R)$ remain linearly independent when restricted to $\Omega$, one obtains 
$\dim\big(\ker(A_{max,\Omega})\big)=\infty$. 

Finally, \eqref{3.12} follows from the fact that $A_{min,\Omega}-v_{-}$ is strictly positive and the 
defect indices are constant on the (connected) set of points of regular type of the closed symmetric operator 
$A_{min,\Omega}-v_{-}$; in particular, the set of regular points of $A_{min,\Omega}-v_{-}$ contains the set 
$\bbC\backslash[v_-,\infty)$ (cf. \cite[Propositions~2.4 and 3.2, and p.~39]{Sc12}). 
\end{proof}
%%%%%%%%

Before taking a closer look at quadratic forms associated to Schr\"odinger-type operators, we briefly introduce the 
basic facts underlying sesquilinear forms drawing primarily from \cite[Ch.~VI]{Ka80}: Let $\cD$ be a linear subspace 
of a complex, separable Hilbert space $\cH$, then 
\begin{equation}\label{RDFG}
\ga:\begin{cases} 
\cD\times\cD\to\bbC, 
\\
(u,v)\mapsto\ga(u,v), 
\end{cases}
\end{equation} 
is called a {\it sesquilinear form} (in short, a {\it form}) in $\cH$ if $\ga(\dott,\dott)$ is linear in the second 
argument and antilinear in the first; $\cD$ then equals the domain of $\ga$ (i.e., $\dom(\ga)=\cD$). 
The underlying {\it quadratic form} is given by $\ga(u,u)$, $u\in\dom(\ga)$. One calls $\ga$ {\it symmetric} if 
$\ga(u,v)=\ga(v,u)^*$, $u,v\in\dom(\ga)$ (with $*$ denoting complex conjugation to distinguish it from the operation of closure). 
A symmetric form $\gs$ is called {\it bounded from below} if there exists $c\in\bbR$ such that 
$\gs(u,u)\geq c\|u\|_{\cH}^2$ for every $u\in\dom(\gs)$. The sesquilinear form $\gt$ is called {\it closed} 
\begin{align}\label{RDcf}
& \text{if }\,\{u_j\}_{j\in\bbN}\subset\dom(\gt)\; u\in\cH\,\text{ satisfying }\,\|u_j-u\|_{\cH}\underset{j\to\infty}\longrightarrow 0 
\\ 
& \quad\text{and }\,\gt(u_j-u_k,u_j-u_k)\underset{j,k\to\infty}\longrightarrow 0\,\text{ implies }\,u\in\dom(\gt)
\\ 
& \quad\text{and }\,\gt(u_j-u,u_j-u)\underset{j\to\infty}\longrightarrow 0. 
\end{align}
A sesquilinear form $\gt$ is called {\it closable} if it has a closed extension; the smallest closed 
extension of a sequilinear form $\ga$ is called its {\it closure} and denoted by $\ol\ga$. Finally, a 
linear subspace $\cD_0$ of $\cH$ is called a {\it core} of the closed sesquilinear form $\gt$ if 
$\ol{\gt|_{\cD_0}}=\gt$. 

The celebrated second representation theorem (combined with a special case of the first representation theorem) 
for forms then reads as follows.

%%%%%%%%
\begin{theorem}\lb{t6.formrep} 
Let $\gt$ be a densely defined, closed sesquilinear form bounded from below by some $c\in\bbR$ in $\cH$.
Then there exists a self-adjoint operator $T$ in $\cH$ such that $T\geq cI_{\cH}$ and the following properties hold: 
\\[1mm]
$(i)$ One has $\dom(T)\subseteq\dom(\gt)$ and 
\begin{equation}\label{WACO.TX.123}
\gt(u,v)=(u,Tv)_{\cH}\quad\forall\,u\in\dom(\gt),\,\forall\,v\in\dom(T).
\end{equation}  
$(ii)$ The linear subspace $\dom(T)$ is a core of $\gt$. 
\\[1mm] 
$(iii)$ If $v\in\dom(\gt)$, $w\in\cH$ and 
\begin{equation}
\gt(u,v)=(u,w)_{\cH} 
\end{equation} 
holds for all $u$ belonging to a core of $\gt$, then $v\in\dom(T)$ and $Tv=w$. 
The self-adjoint operator $T$ is uniquely determined by condition $(i)$. 
\\[1mm] 
$(iv)$ One has $\dom\big(|T|^{1/2}\big)=\dom\big((T-cI_{\cH})^{1/2}\big)=\dom(\gt)$ and 
\begin{equation}\label{WACO.yall}
\gt(u,v)=\big((T-cI_{\cH})^{1/2}u,(T-cI_{\cH})^{1/2}v\big)_{\cH}+c(u,v)_{\cH}\quad\forall\,u,v\in\dom(\gt). 
\end{equation}
Moreover, $\cD_0\subseteq\dom(\gt)$ is a core of $\gt$ if and only if it is a core of $(T-cI_{\cH})^{1/2}$. 
\end{theorem}
%%%%%%%%

Another particularly useful special case of Theorem~\ref{t6.formrep} is the following result:

%%%%%%%%
\begin{theorem}\lb{t6.form}
Let $\cH_j$, $j=1,2$, be complex separable Hilbert spaces, assume that the linear operator $S$ maps 
$\dom(S)\subseteq\cH_1$ into $\cH_2$, and introduce the nonnegative sesquilinear form $\gt_S$ via 
\begin{equation}\label{WACO.TX.TX}
\gt_S(u,v)=(Su,Sv)_{\cH_2},\quad u,v\in\dom(\gt_S)=\dom(S). 
\end{equation} 
Then the following properties hold: 
\\[1mm] 
$(i)$ The form $\gt_S$ is closable $($resp., closed\,$)$ if and only $S$ is closable $($resp., closed\,$)$. 
\\[1mm] 
$(ii)$ If $\gt_S$ is closed, then $\cD_0\subseteq\dom(\gt_S)=\dom(S)$ is a core of $\gt_S$ if and only if it is a core of $S$. 
\\[1mm]
$(iii)$ Suppose $S$ is densely defined and closed. Then the self-adjoint, nonnegative operator $T_S$ in 
$\cH_1$, uniquely associated to $\gt_S$ via Theorem~\ref{t6.formrep}, is given by $T_S=S^*S\geq 0$. 
Moreover, $\dom(S^* S)$ is a core of $\gt_S$ and hence of $S$. 
\end{theorem}
%%%%%%%% 

Item $(iii)$ of Theorem~\ref{t6.form} independently proves a well-known theorem of von Neumann 
\cite[Satz~3]{Ne32} (see also \cite{GS19} and the references therein).

We continue with a brief outline of the connection between the Friedrichs extension of closed, symmetric 
operators bounded from below and the theory of sequilinear forms. Let $S$ be a densely defined, closed, 
symmetric, linear operator in $\cH$ satisfying $S\geq cI_{\cH}$ for some $c\in\bbR$. Then Freudenthal's 
intrinsic description (cf. \cite{Fr36}) of the self-adjoint Friedrichs extension $S_F$ of $S$ 
(satisfying $S_F\geq cI_{\cH}$) is given by 
\begin{align}\label{Fr-2} 
& S_Fu:=S^*u\,\,\text{ for each }\,\,u\in\dom(S_F),\,\,\text{ where}
\nonumber\\[4pt]
& \dom(S_F):=\Big\{v\in\dom(S^*)\,\big|\,\text{ there exists }\,\{v_j\}_{j\in\bbN}\subset\dom(S)\,\text{ with} 
\\[4pt]
&\quad\lim_{j\to\infty}\|v_j-v\|_{\cH}=0\,\text{ and }\,\big((v_j-v_k),S(v_j-v_k)\big)_\cH\to 0 
\,\text{ as }\,j,k\to\infty\Big\}.    
\nonumber 
\end{align}

%%%%%%%%
\begin{theorem}\lb{t6.Fr} 
Suppose that $S$ is a densely defined, symmetric, linear operator in $\cH$ bounded from below, 
and introduce the sesquilinear form $\gs$ in $\cH$ by
\begin{equation}\label{WACO.TX.TX.TX}
\gs(u,v):=(u,Sv)_{\cH},\quad u,v\in\dom(\gs)=\dom(S). 
\end{equation}
Then the following properties hold: 
\\[1mm] 
$(i)$ The form $\gs$ is densely defined, symmetric, and closable. Denoting its closure by $\ol\gs$, 
the self-adjoint operator uniquely associated to $\ol\gs$ via Theorem~\ref{t6.formrep} is precisely 
the Friedrichs extension $S_F$ of $S$. 
\\[1mm]
$(ii)$ Among all self-adjoint extensions $\wti S$ of $S$ bounded from below, $S_F$ has the smallest 
form domain $\Big($i.e., the form domain $\dom\big(|S_F|^{1/2}\big)$ of the sesquilinear form of $S_F$ 
is contained in the form domain $\dom\Big(\big|\wti S\big|^{1/2}\Big)$ of any $\wti S$$\Big)$. 
\\[1mm]
$(iii)$ The Friedrichs extension $S_F$ of $S$ is the only self-adjoint extension bounded from below 
whose domain is contained in $\dom({\ol\gs})$. 
\end{theorem}
%%%%%%%%

Next, retaining Hypothesis~\ref{h3.1}, introduce the sesquilinear form
\begin{equation}\label{afform}
\mathfrak a_{F,\Omega}(f,g):=\big(\nabla f,\nabla g\big)_{[L^2(\Omega)]^n}
+(f,Vg)_{L^2(\Omega)},\quad\dom(\mathfrak a_{F,\Omega})=\accentset{\circ}{W}^1(\Omega),
\end{equation}
which is densely defined, closed, symmetric, and semibounded from below 
(by $v_{-}$) in $L^2(\Omega)$. Hence, it follows from the first representation theorem 
\cite[Theorem VI.2.1]{Ka80}, and here recorded in Theorem~\ref{t6.formrep}\,$(i)$, that there 
is a unique self-adjoint operator $A_{F,\Omega}$ in $L^2(\Omega)$ such that the identity
\begin{equation}\label{formopaf}
\mathfrak a_{F,\Omega}(f,g)=\big(f,A_{F,\Omega}g\big)_{L^2(\Omega)}
\end{equation}
holds for all $f\in\dom(\mathfrak a_{F,\Omega})=\accentset{\circ}{W}^1(\Omega)$ and all 
$g\in\dom(A_{F,\Omega})\subset\dom(\mathfrak a_{F,\Omega})$. Making use of \eqref{Rdac.WWW} 
and Green's formula it follows that
\begin{equation}\label{OUYTf}
A_{F,\Omega}=-\Delta+V,\quad\dom(A_{F,\Omega})=\big\{f\in\accentset{\circ}{W}^1(\Omega)
\,\big|\,\Delta f\in L^2(\Omega)\big\}, 
\end{equation}
and hence $A_{F,\Omega}$ is a self-adjoint extension of the minimal realization $A_{min,\Omega}$ 
of $-\Delta+V$ defined in \eqref{eqn:Amin1}. By \cite[Subsection~VI.2.3]{Ka80}, as recalled in 
Theorem~\ref{t6.Fr}\,$(iii)$, $A_{F,\Omega}$ represents the Friedrichs extension of $A_{min,\Omega}$.

The next well-known theorem collects some properties of the Friedrichs extension $A_{F,\Omega}$ 
in the present setting (see, for instance, \cite[Section~6.1]{EE18}). 

%%%%%%%%
\begin{theorem}\label{tfried}
Assume Hypothesis~\ref{h3.1}. Then the Friedrichs extension $A_{F,\Omega}$ of $A_{min,\Omega}$ 
is a self-adjoint operator in $L^2(\Omega)$ with spectrum contained in $[v_-,\infty)$. 
If, in addition, $\Omega$ is a bounded domain then the resolvent of $A_{F,\Omega}$ is 
compact, the spectrum is purely discrete and contained in $(v_-,\infty)$. In particular, 
$\sigma_{ess}(A_{F,\Om})=\varnothing$. 
\end{theorem}
%%%%%%%

Next, we study the Dirichlet and Neumann realizations of $-\Delta+V$ on a bounded Lipschitz 
domain $\Omega$ in $\bbR^n$. In this context we now strengthen Hypothesis~\ref{h3.1} and use 
the following set of assumptions until and including Section~\ref{s10}:

%%%%%%%%%
\begin{hypothesis}\label{h4.2} 
Let $n\in\bbN\backslash\{1\}$, assume that $\Omega\subset\bbR^n$ is a bounded Lipschitz domain, 
and suppose that $V\in L^\infty(\Omega)$ is real-valued.
\end{hypothesis}
%%%%%%%%%

In the setting of bounded Lipschitz domains it follows from \eqref{eq:WH} and \eqref{incl-Yb.EE} 
with $s=1$ that $\dom(\mathfrak a_{F,\Omega})=\accentset{\circ}{H}^1(\Omega)$ and 
the Friedrichs extension $A_{F,\Omega}$ coincides with the self-adjoint Dirichlet 
operator 
\begin{equation}\label{We-Q.10EE}
\begin{split}
& A_{D,\Omega}=-\Delta+V, 
\\[2pt]
& \dom(A_{D,\Omega})=\big\{f\in H^{1}(\Omega)\cap\dom(A_{max,\Omega})\,\big|\,\gamma_D f=0\big\}
\\[2pt]
& \hskip 0.72in
=\big\{f\in\accentset{\circ}{H}^1(\Omega)\,\big|\,\Delta f\in L^2(\Omega)\big\}.
\end{split}
\end{equation}
Next, we collect further properties of the self-adjoint Dirichlet operator.

%%%%%%%%%
\begin{theorem}\label{t4.4}
Assume Hypothesis~\ref{h4.2} and let $A_{D,\Omega}$ be the Dirichlet realization of $-\Delta+V$ 
in \eqref{We-Q.10EE}. Then the functions in $\dom(A_{D,\Omega})$ possess $H^{3/2}$-regularity, 
that is, $\dom(A_{D,\Omega})\subset H^{3/2}(\Omega)$, 
\begin{equation}\label{We-Q.10EE-jussi}
\begin{split}
& A_{D,\Omega}=-\Delta+V, 
\\[2pt]
&\dom(A_{D,\Omega})=\big\{f\in H^{3/2}(\Omega)\cap\dom(A_{max,\Omega})\,\big|\,\gamma_D f=0\big\}
\\[2pt]
&\hskip 0.72in
=\big\{f\in H^{3/2}(\Omega)\cap\accentset{\circ}{H}^1(\Omega)\,\big|\,\Delta f\in L^2(\Omega)\big\},
\end{split}
\end{equation}
and on $\dom(A_{D,\Omega})$ the norms 
\begin{equation}\label{eq:UUjh-jussi}
f\mapsto\|f\|_{H^s(\Omega)}+\|\Delta f\|_{L^2(\Omega)},\quad s\in\big[0,\tfrac{3}{2}\big],
\end{equation}
are equivalent. In addition, $A_{D,\Omega}$ is self-adjoint in $L^2(\Omega)$, with compact 
resolvent, and purely discrete spectrum, contained in $(v_-,\infty)$. In particular, 
$\sigma_{ess}(A_{D,\Om})=\varnothing$. Moreover,
\begin{equation}\label{2.4Hba}
\dom\big(|A_{D,\Om}|^{1/2}\big)=\accentset{\circ}{H}^1(\Om).   
\end{equation}
\end{theorem}
%%%%%%%%%
\begin{proof} 
The additional $H^{3/2}$-regularity of the function in $\dom(A_{D,\Omega})$ follows from
\eqref{eq:EFFa.NNN} with $s=1$, which together with \eqref{We-Q.10EE} also yields 
\eqref{We-Q.10EE-jussi}. For $s\in [1,\tfrac{3}{2}]$ the claim in \eqref{eq:UUjh-jussi}
is a consequence of \eqref{gafvv.6577}, and for $s\in [0,1]$ the reasoning is as follows. 
For $f\in\dom(A_{D,\Omega})$ and $s=1$ one obtains from \eqref{2.9}
\begin{equation}\label{greendirichlet}
\begin{split}
0 &=(\gamma_D f,\gamma_N f)_{L^2(\partial\Omega)}
={}_{H^{1/2}(\partial\Omega)}\big\langle\gamma_D f,\gamma_N f\big\rangle_{H^{-1/2}(\partial\Omega)}
\\[2pt]
&=\big(\nabla f,\nabla f)_{[L^2(\Omega)]^n}+(f,\Delta f)_{L^2(\Omega)},
\end{split}
\end{equation}
which leads to 
\begin{equation}\label{jaffkg3s-WACO}
\|\nabla f\|^2_{[L^2(\Omega)]^n}\leq\|f\|_{L^2(\Omega)}\,\|\Delta f\|_{L^2(\Omega)}
\leq\big(\|f\|_{L^2(\Omega)}+\|\Delta f\|_{L^2(\Omega)}\big)^2
\end{equation}
for $f\in\dom(A_{D,\Omega})$. Therefore, $\|f\|_{H^1(\Omega)}\leq C(\|f\|_{L^2(\Omega)}+\|\Delta f\|_{L^2(\Omega)})$ 
on $\dom(A_{D,\Omega})$ which in turn implies \eqref{eq:UUjh-jussi} for $s\in [0,1]$. 
The remaining statements follow from Theorem~\ref{tfried} and the second representation theorem 
\cite[Theorem~VI.2.23]{Ka80} gives \eqref{2.4Hba}, see also \cite[Theorem~2.10]{GM08} 
and \cite[Theorem~4.6]{GM09b} for the case $V=0$. 
\end{proof}
%%%%%%%%%

Next, we introduce the sesquilinear form
\begin{equation}\label{formN}
\mathfrak a_{N,\Omega}(f,g):=\big(\nabla f,\nabla g\big)_{[L^2(\Omega)]^n}
+(Vf,g)_{L^2(\Omega)},\quad\dom(\mathfrak a_{N,\Omega})=H^1(\Omega),
\end{equation}
which is densely defined, closed, symmetric, and semibounded from below (by $v_-$) 
in $L^2(\Omega)$. One observes that $\mathfrak a_{N,\Omega}$ is an extension of the form 
$\mathfrak a_{F,\Omega}$ in \eqref{afform} since
\begin{equation}\label{formN-dom}
\dom(\mathfrak a_{F,\Omega})=\accentset{\circ}{H}^1(\Omega)\subset H^1(\Omega)=\dom(\mathfrak a_{N,\Omega}). 
\end{equation}
As above, it follows from the First Representation Theorem (cf., e.g., \cite[Theorem~VI.2.1]{Ka80}; see also Theorem~\ref{t6.formrep}) 
that there is a unique self-adjoint operator $A_{N,\Omega}$ in $L^2(\Omega)$ such that the identity
\begin{equation}\label{formopan}
\mathfrak a_{N,\Omega}(f,g)=\big(f,A_{N,\Omega} g\big)_{L^2(\Omega)}
\end{equation}
holds for all $f\in\dom(\mathfrak a_{N,\Omega})=H^1(\Omega)$ and all 
$g\in\dom(A_{N,\Omega})\subset\dom(\mathfrak a_{N,\Omega})$. Making use of 
\eqref{formN}, \eqref{formopan}, and \eqref{2.9} for $s=1$ one obtains 
\begin{equation}\label{neumannjussi}
(f,A_{N,\Omega} g)_{L^2(\Omega)}=\big(f,(-\Delta+V)g\big)_{L^2(\Omega)}+(\gamma_D f,\gamma_N g)_{L^2(\partial\Omega)}
\end{equation}
for $g\in\dom(A_{N,\Omega})$ and all $f\in H^1(\Omega)$. 
By considering $f\in\accentset{\circ}{H}^1(\Omega)$ only it follows in a first step from 
\eqref{neumannjussi} that $A_{N,\Omega}=-\Delta+V$. In a second step, taking into 
account that the range of $\gamma_D$ restricted to $\dom(\mathfrak a_{N,\Omega})=H^1(\Omega)$ 
is the dense subspace $H^{1/2}(\partial\Omega)$ of $L^2(\partial\Omega)$ (see \eqref{eqn:gammaDs.2aux} with $s=1$), 
one finds $\gamma_N g=0$ for all functions $g\in\dom(A_{N,\Omega})$. Thus, one obtains 
\begin{align}\label{We-Q.10EENN}
\begin{split}
& A_{N,\Omega}=-\Delta+V, 
\\[2pt]
&\dom(A_{N,\Omega}) 
=\big\{f\in H^{1}(\Omega)\cap\dom(A_{max,\Omega})\,\big|\,\gamma_N f=0\big\}
\\[2pt]
&\hskip 0.72in
=\big\{f\in H^{1}(\Omega)\,\big|\,\Delta f\in L^2(\Omega)\,\text{ and }\,\gamma_N f=0\big\},
\end{split}
\end{align}
and hence $A_{N,\Omega}$ is a self-adjoint extension of the minimal realization 
$A_{min,\Omega}$ of $-\Delta+V$ defined in \eqref{eqn:Amin1}. In the following we 
shall refer to $A_{N,\Omega}$ as the Neumann extension of $A_{min,\Omega}$.

Next, we list some useful properties of the Neumann realization.  

%%%%%%%%%%
\begin{theorem}\label{t4.5} 
Assume Hypothesis~\ref{h4.2} and let $A_{N,\Omega}$ be the Neumann realization of 
$-\Delta+V$ in \eqref{We-Q.10EENN}. Then the functions in $\dom(A_{N,\Omega})$ possess 
$H^{3/2}$-regularity, that is, $\dom(A_{N,\Omega})\subset H^{3/2}(\Omega)$, 
\begin{equation}\label{We-Q.10EE-jussiNN}
\begin{split}
& A_{N,\Omega}=-\Delta+V, 
\\[2pt]
&\dom(A_{N,\Omega})=\big\{f\in H^{3/2}(\Omega)\cap\dom(A_{max,\Omega})\,\big|\,\gamma_N f=0\big\}
\\[2pt]
&\hskip 0.72in
=\big\{f\in H^{3/2}(\Omega)\,\big|\,\Delta f\in L^2(\Omega)\,\text{ and }\,\gamma_N f=0\big\},
\end{split}
\end{equation}
and on $\dom(A_{N,\Omega})$ the norms 
\begin{equation}\label{eq:UUjh-jussiNN}
f\mapsto\|f\|_{H^s(\Omega)}+\|\Delta f\|_{L^2(\Omega)},\quad s\in\big[0,\tfrac{3}{2}\big],
\end{equation}
are equivalent. In addition, $A_{N,\Omega}$ is self-adjoint in $L^2(\Omega)$, with compact 
resolvent, and purely discrete spectrum, contained in $[v_-,\infty)$. In particular, 
$\sigma_{ess}(A_{N,\Om})=\varnothing$. Moreover,
\begin{equation}\label{2.4HyN}
\dom\big(|A_{N,\Om}|^{1/2}\big)=H^1(\Omega).   
\end{equation}
\end{theorem}
%%%%%%%%%%
\begin{proof}
The $H^{3/2}$-regularity of the functions in $\dom(A_{N,\Omega})$ is a consequence of \eqref{eq:EFFa.111} 
(used with $s=1$), while the claim in \eqref{eq:UUjh-jussiNN} follows immediately from \eqref{gafvv.6588}. 
The remaining statements can be found in \cite[Theorem~2.6]{GM08} and \cite[Theorem 4.5]{GM09b} for the 
case $V=0$. The proof in the case $V\not=0$ is analogous. We note that the spectrum of  $A_{N,\Omega}$ 
is bounded from below by $v_-$ since the corresponding form $\mathfrak a_{N,\Omega}$ in \eqref{formN} 
is bounded from below by $v_-$.
\end{proof}
%%%%%%%%%%

Next, as an immediate consequence of Lemma~\ref{l3.3} and 
\eqref{eq:WH},  we state a lemma describing the domain of the minimal operator $A_{min,\Omega}$.

%%%%%%%%%%
\begin{lemma}\label{l4.3}
Assume Hypothesis~\ref{h4.2}. Then the closed symmetric operator $A_{min,\Omega}$ is given by
\begin{equation}\label{eqn:Amin2}
A_{min,\Omega}=-\Delta+V,\quad\dom(A_{min,\Omega})=\accentset{\circ}{H}^2(\Omega).
\end{equation}
\end{lemma}
%%%%%%%%%

Finally we show that $A_{D,\Omega}$ and $A_{N,\Omega}$ are relatively prime (or disjoint), 
a fact that will play a prominent role later on.

%%%%%%%%%%
\begin{theorem}\label{t4.6}
Assume Hypothesis~\ref{h4.2}. Then the operators $A_{D,\Omega}$ and $A_{N,\Omega}$ 
are relatively prime, that is,
\begin{equation}\label{3.30}
\dom(A_{D,\Omega})\cap\dom(A_{N,\Omega})=\dom(A_{min,\Omega})=\accentset{\circ}{H}^2(\Omega).
\end{equation}
\end{theorem}
%%%%%%%%%%
\begin{proof} 
Let $f\in\dom(A_{D,\Omega})\cap\dom(A_{N,\Omega})$. Then from \eqref{We-Q.10EE-jussi} 
and \eqref{We-Q.10EE-jussiNN} one deduces $f\in H^{3/2}(\Omega)$ and $\gamma_D f=\gamma_N f=0$.
Together with \eqref{GGGRRR}, these conditions ensure that for every 
$\psi\in C^\infty(\overline{\Omega})$ one can write
\begin{equation}\label{uartUTR}
(\overline{f},\Delta\psi)_{L^2(\Omega)}=(\overline{\Delta f},\psi)_{L^2(\Omega)}.
\end{equation}
As in the past, using tilde to denote the extension of a function, originally defined in $\Omega$, 
to the entire space $\bbR^n$ by taking said extension to be zero outside $\Omega$, the fact that 
$\widetilde{f}\in L^2(\bbR^n)$ and \eqref{uartUTR} imply
\begin{align}\label{jh7g8}
{}_{\mathcal D'({\mathbb{R}}^n)}\big\langle\Delta\widetilde f,\varphi\big\rangle_{\mathcal D({\mathbb{R}}^n)}
&={}_{\mathcal D'({\mathbb{R}}^n)}\big\langle\widetilde f,\Delta\varphi\big\rangle_{\mathcal D({\mathbb{R}}^n)}
=\big(\overline{f},\Delta\varphi|_\Omega\big)_{L^2(\Omega)}
\nonumber\\[2pt]
&=\big(\overline{\Delta f},\varphi|_\Omega\big)_{L^2(\Omega)}=\big(\overline{\widetilde{\Delta f}},\varphi\big)_{L^2(\bbR^n)}
\nonumber\\[2pt]
&={}_{\mathcal D'({\mathbb{R}}^n)}\big\langle\widetilde{\Delta f},\varphi\big\rangle_{\mathcal D({\mathbb{R}}^n)}
\end{align}
for all $\varphi\in C_0^\infty(\bbR^n)$. As such, $\Delta\widetilde{f}=\widetilde{\Delta f}$ in $\cD'(\bbR^n)$. 
Since $\widetilde{\Delta f}\in L^2(\bbR^n)$, invoking standard elliptic regularity one concludes that 
$\widetilde{f}\in H^2_{\rm loc}(\bbR^n)$, which further implies $f\in H^2(\Omega)$. With this in hand, 
we are in a position to invoke Lemma~\ref{l4.3} and \eqref{Tan-C3} to conclude that 
\begin{equation}\label{eq:WACO.MMM}
\dom(A_{D,\Omega})\cap\dom(A_{N,\Omega})\subset\accentset{\circ}{H}^2(\Omega)=\dom(A_{min,\Omega}). 
\end{equation}
This establishes the left-to-right inclusion in \eqref{3.30}. The opposite inclusion follows from 
Lemma~\ref{l4.3} and the fact that $A_{D,\Omega}$ and $A_{N,\Omega}$ are both extensions of $A_{min,\Omega}$.
\end{proof}
%%%%%%%%%%

%%%%%%%%%%%%%%%%%%%%%%%%%%%%%%
%%%%%%%%%%%%%%%%%%%%%%%%%%%%%%
\section{Weyl--Titchmarsh Operators for Schr\"{o}dinger Operators on Bounded Lipschitz Domains} 
\label{s7}
%%%%%%%%%%%%%%%%%%%%%%%%%%%%%%
%%%%%%%%%%%%%%%%%%%%%%%%%%%%%%

In this section we study $z$-dependent Dirichlet-to-Neumann maps, that is, 
Weyl--Titchmarsh operators, for Schr\"odinger operators on bounded Lipschitz domains, assuming 
Hypothesis~\ref{h4.2} throughout this section. 

For each complex number $z$ not in the spectrum of the self-adjoint Dirichlet operator 
$A_{D,\Omega}$, that is, for $z\in\rho(A_{D,\Omega})=\bbC\backslash\sigma(A_{D,\Omega})$, 
and for each $s\in[0,\frac{3}{2}]$, the characterization in \eqref{We-Q.10EE-jussi} implies the 
following direct sum decompositions of $\dom(A_{max,\Omega})\cap H^{s}(\Omega)$:
\begin{equation}\label{eqn:decAmax}
\begin{split}
\dom(A_{max,\Omega})\cap H^{s}(\Omega) 
&=\big[\dom(A_{D,\Omega})\,\dot{+}\,\ker(A_{max,\Omega}-z I)\big]\cap H^{s}(\Omega) 
\\[2pt] 
&=\dom(A_{D,\Omega})\,\dot{+}\,\big\{f\in H^s(\Omega)\,\big|\,-\Delta f+Vf=z f\big\}.
\end{split}
\end{equation}
In a similar manner, \eqref{We-Q.10EE-jussiNN} ensures that the following direct sum decomposition holds 
for the self-adjoint Neumann operator $A_{N,\Omega}$, $z\in\rho(A_{N,\Omega})$, $s\in[0,\frac{3}{2}]$:
\begin{equation}\label{eqn:decAmax2}
\begin{split}
\dom(A_{max,\Omega})\cap H^{s}(\Omega) 
&=\big[\dom(A_{N,\Omega})\,\dotplus\,\ker(A_{max,\Omega}-z I)\big]\cap H^{s}(\Omega)
\\[2pt] 
&=\dom(A_{N,\Omega})\,\dotplus\,\big\{f\in H^s(\Omega)\,\big|\,-\Delta f+Vf=zf\big\}.
\end{split}
\end{equation}

For further reference, we also note that if $z\in\rho(A_{D,\Om})$ then 
\begin{equation}\label{3.2SD1}
\gamma_N(A_{D,\Om}-zI)^{-1}\in\cB\big(L^2(\Omega),L^2(\partial\Omega)\big),
\end{equation}
by \eqref{We-Q.10EE-jussi}, \eqref{eq:UUjh-jussi} with $s=0$ and $s=\tfrac{3}{2}$, 
and \eqref{eqn:gammaN-pp} with $s=\tfrac{3}{2}$. In particular, \eqref{3.2SD1} entails
\begin{equation}\label{3.2SD2}
\big[\gamma_N(A_{D,\Om}-zI)^{-1}\big]^*\in\cB\big(L^2(\partial\Omega),L^2(\Omega)\big).  
\end{equation}
Similarly, if $z\in\rho(A_{N,\Om})$ then \eqref{We-Q.10EE-jussiNN}, \eqref{eq:UUjh-jussiNN} 
with $s=0$ and $s=\tfrac{3}{2}$, and \eqref{eqn:gammaDs.2} with $s=\tfrac{3}{2}$, imply that
\begin{equation}\label{3.2SD1N}
\gamma_D(A_{N,\Om}-zI)^{-1}\in\cB\big(L^2(\Omega),H^1(\partial\Omega)\big),
\end{equation}
hence
\begin{equation}\label{3.2SD2N}
\big[\gamma_D(A_{N,\Om}-zI)^{-1}\big]^*\in\cB\big(H^{-1}(\partial\Omega),L^2(\Omega)\big).  
\end{equation}

To be able to proceed, we also need the following useful results contained in the next two lemmas: 

%%%%%%%%%%
\begin{lemma}\label{l5.1} 
Assume Hypothesis~\ref{h4.2} and fix an arbitrary $z\in\rho(A_{D,\Omega})\cup\rho(A_{N,\Omega})$. 
Then $\ker(A_{max,\Omega}-zI)\cap H^{3/2}(\Omega)$ is dense in $\ker(A_{max,\Omega}-zI)$ 
when the latter space is equipped with the $L^2(\Omega)$-norm.
\end{lemma}
%%%%%%%%%%
\begin{proof}
Fix $z\in\rho(A_{D,\Omega})$ (the case when $z\in\rho(A_{N,\Omega})$ is similar).
Employing the density result \eqref{gag-8gb} (with $s_1=s_2=0$) shows that 
given any $f\in\ker(A_{max,\Omega}-zI)$ there exists a sequence 
$\{g_j\}_{j\in\bbN}\subset C^\infty(\overline{\Omega})$ with the property 
that $g_j\to f$ and $\Delta g_j\to\Delta f$ in $L^2(\Omega)$ as $j\to\infty$.
Then 
\begin{equation}\label{74r4}
f_j:=\big[g_j-(A_{D,\Omega}-zI)^{-1}(-\Delta+V-zI)g_j\big]
\in\ker(A_{max,\Omega}-z I)\cap H^{3/2}(\Omega)
\end{equation}
for every $j\in\bbN$, and since and since $V\in L^\infty(\Omega)$ it follows that
\begin{equation}\label{oitre}
(-\Delta+V-zI)g_j\underset{j\to\infty}{\longrightarrow}(-\Delta+V-zI)f=0\,\text{ in $L^2(\Omega)$}.
\end{equation}  
Therefore, one concludes that $f_j\rightarrow f$ in $L^2(\Omega)$ as $j\to\infty$. 
\end{proof}
%%%%%%%%%

Here is the second density result alluded to above. 

%%%%%%%%%%
\begin{lemma}\label{l5.1MM} 
Assume Hypothesis~\ref{h4.2}. Then $\dom(A_{max,\Omega})\cap H^{3/2}(\Omega)$
is a dense subspace of $\dom(A_{max,\Omega})$, when the latter space is equipped with the 
natural graph norm $f\mapsto\|f\|_{L^2(\Omega)}+\|\Delta f\|_{L^2(\Omega)}$.
\end{lemma}
%%%%%%%%%%
\begin{proof}
Fix some $z\in\rho(A_{D,\Omega})$ and select an arbitrary $f\in\dom(A_{max,\Omega})$.
Use \eqref{eqn:decAmax} (with $s=0$) to decompose $f=g+h$ with $g\in\dom(A_{D,\Omega})$ 
and $h\in\ker(A_{max,\Omega}-zI)$. By \eqref{We-Q.10EE-jussi} this entails 
\begin{equation}\label{Mi-tr-MM.1}
g\in\dom(A_{max,\Omega})\cap H^{3/2}(\Omega).
\end{equation}  
Then invoke Lemma~\ref{l5.1} to produce a sequence 
\begin{equation}\label{eq:r4rr}
\{h_j\}_{j\in{\mathbb{N}}}\subset\ker(A_{max,\Omega}-zI)\cap H^{3/2}(\Omega)
\subset\dom(A_{max,\Omega})\cap H^{3/2}(\Omega),
\end{equation}
such that $h_j\to h$ in $L^2(\Omega)$ as $j\to\infty$. 
Since $V\in L^\infty(\Omega)$, one also has  
\begin{equation}\label{Mi-tr-MM.2}
\Delta h_j=(V-zI)h_j\underset{j\to\infty}{\longrightarrow}(V-zI)h=\Delta h\,\text{ in $L^2(\Omega)$}.
\end{equation}  
Hence $g+h_j\in\dom(A_{max,\Omega})\cap H^{3/2}(\Omega)$ for each $j\in{\mathbb{N}}$, and
\begin{equation}\label{Mi-tr-MM.3}
g+h_j\underset{j\to\infty}{\longrightarrow} f\,\text{ in $L^2(\Omega)$}\,\text{ and }\,
\Delta(g+h_j)\underset{j\to\infty}{\longrightarrow}\Delta f\,\text{ in $L^2(\Omega)$},
\end{equation}  
from which the desired conclusion follows. 
\end{proof}
%%%%%%%%%

Our next result extends \cite[Theorem~3.6, Corollary~3.3]{GM08} and \cite[Theorem~5.3]{GM11}.

%%%%%%%%%
\begin{lemma}\label{LamAA.2}
Assume Hypothesis~\ref{h4.2}. Then for each $z\in\rho(A_{D,\Omega})$ and $s\in [0,1]$ 
the boundary value problem
\begin{equation}\label{eqn:bvp}
\begin{cases}
(-\Delta+V-z)f=0\,\text{ in $\Omega,\quad f\in H^{s+(1/2)}(\Omega)\cap\dom(A_{max,\Omega})$,}     
\\[2pt]  
\gamma_D f=\varphi\,\text{ on $\partial\Omega,\quad\varphi\in H^s(\partial\Omega)$,}
\end{cases}    
\end{equation}
is well posed, with unique solution $f=f_D(z,\varphi)$ given by 
\begin{equation}\label{4.7}
f_D(z,\varphi)=-\big[\gamma_N(A_{D,\Omega}-{\ol z}I)^{-1}\big]^*\varphi, 
\end{equation} 
with the adjoint understood in the sense of \eqref{3.2SD2}.
\end{lemma}
%%%%%%%%%
\begin{proof}
That \eqref{eqn:bvp} is uniquely solvable is a consequence of the surjectivity of the boundary 
trace map $\gamma_D$ in \eqref{eqn:gammaDs.2} and the decomposition in \eqref{eqn:decAmax}. 
Regarding \eqref{4.7}, we denote by $f_D$ the unique solution of \eqref{eqn:bvp}. 
Based on \eqref{3.2SD1}--\eqref{3.2SD2} and Green's formula \eqref{GGGRRR}, for each 
$v\in L^2(\Omega)$ one computes 
\begin{align}\label{yrer}
(f_D,v)_{L^2(\Omega)} &=
\big(f_D,(-\Delta+V-\overline{z})(A_{D,\Omega}-\overline{z}I)^{-1}v\big)_{L^2(\Omega)}
\nonumber\\[2pt] 
&=\big((-\Delta+V-z)f_D,(A_{D,\Omega}-\overline{z}I)^{-1}v\big)_{L^2(\Omega)}
\nonumber\\[2pt] 
&\quad +{}_{H^{-1}(\partial\Omega)}
\big\langle\gamma_N f_D,\gamma_D(A_{D,\Omega}-\overline{z}I)^{-1}v\big\rangle
_{H^1(\partial\Omega)}  
\nonumber\\[2pt] 
&\quad -\big(\gamma_D f_D,\gamma_N(A_{D,\Omega}-\overline{z}I)^{-1}v\big)_{L^2(\partial\Omega)}
\nonumber\\[2pt] 
&=-\big(\varphi,\gamma_N(A_{D,\Omega}-\overline{z}I)^{-1}v\big)_{L^2(\partial\Omega)} 
\nonumber\\[2pt] 
&=-\big(\big[\gamma_N(A_{D,\Omega}-\overline{z}I)^{-1}\big]^\ast\varphi,v\big)_{L^2(\Omega)}.
\end{align}
In light of the arbitrariness of $v$ in $L^2(\Omega)$, this proves \eqref{4.7}.
\end{proof}
%%%%%%%%%

We continue by discussing an extension of  
\cite[Theorems~3.2, 4.3, Corollaries~3.3, 4.4]{GM08}, \cite[Theorem~5.4]{GM11}.

%%%%%%%%
\begin{lemma}\label{LamAA.2N}
Assume Hypothesis~\ref{h4.2}. Then for each $z\in\rho(A_{N,\Omega})$ and $s\in[0,1]$ 
the boundary value problem
\begin{equation}\label{eqn:bvp2}
\begin{cases}
(-\Delta+V-z)f=0\,\text{ in $\Omega,\quad f\in H^{s+(1/2)}(\Omega)\cap\dom(A_{max,\Omega})$,}     
\\[2pt] 
-\gamma_N f=\varphi\,\text{ in $H^{s-1}(\partial\Omega),\quad\varphi\in H^{s-1}(\partial\Omega)$,}
\end{cases}    
\end{equation}
is well posed, with unique solution $f=f_N(z,\varphi)$ given by 
\begin{equation}\label{4.9}
f_N(z,\varphi)=-\big[\gamma_D(A_{N,\Omega}-{\ol z}I)^{-1}\big]^*\varphi,   
\end{equation}
with the adjoint understood in the sense of \eqref{3.2SD2N}.
\end{lemma}
%%%%%%%%
\begin{proof}
Together, the fact that the boundary trace map $\gamma_N$ in \eqref{eqn:gammaN-pp} is surjective
and the decomposition in \eqref{eqn:decAmax2} imply that the boundary value problem \eqref{eqn:bvp2} 
is uniquely solvable. To justify \eqref{4.9}, denote by $f_N$ the unique solution of \eqref{eqn:bvp2}. 
Relying on \eqref{3.2SD1N}--\eqref{3.2SD2N} and Green's formula \eqref{GGGRRR}, for each $v\in L^2(\Omega)$ 
one may write 
\begin{align}\label{yrer.N}
(f_N,v)_{L^2(\Omega)} &=
\big(f_N,(-\Delta+V-\overline{z})(A_{N,\Omega}-\overline{z}I)^{-1}v\big)_{L^2(\Omega)}
\nonumber\\[2pt] 
&=\big((-\Delta+V-z)f_N,(A_{N,\Omega}-\overline{z}I)^{-1}v\big)_{L^2(\Omega)}
\nonumber\\[2pt] 
&\quad +{}_{H^{-1}(\partial\Omega)}
\big\langle\gamma_N f_N,\gamma_D(A_{N,\Omega}-\overline{z}I)^{-1}v\big\rangle
_{H^1(\partial\Omega)}  
\nonumber\\[2pt] 
&\quad -\big(\gamma_D f_N,\gamma_N(A_{N,\Omega}-\overline{z}I)^{-1}v\big)_{L^2(\partial\Omega)}
\nonumber\\[2pt] 
&=-{}_{H^{-1}(\partial\Omega)}
\big\langle\varphi,\gamma_D(A_{N,\Omega}-\overline{z}I)^{-1}v\big\rangle
_{H^1(\partial\Omega)}  
\nonumber\\[2pt] 
&=-\big(\big[\gamma_D(A_{N,\Omega}-\overline{z}I)^{-1}\big]^\ast\varphi,v\big)_{L^2(\Omega)}.
\end{align}
Given that $v\in L^2(\Omega)$ is arbitrary, this proves \eqref{4.9}.
\end{proof}
%%%%%%%%%

Next, we bring into play the solution operator corresponding to the boundary value problems 
\eqref{eqn:bvp} and \eqref{eqn:bvp2}.

%%%%%%%%%
\begin{theorem}\label{t5.2}
Assume Hypothesis~\ref{h4.2}. Then the following assertions hold:
\\
$(i)$ For $z\in\rho(A_{D,\Omega})$ and $s\in[0,1]$, define  
\begin{equation}\label{We-Q.13}
P_{s,D,\Omega}(z): 
\begin{cases} 
H^s(\partial\Omega)\rightarrow H^{s+(1/2)}(\Omega)\cap\dom(A_{max,\Omega}), 
\\[2pt]  
\varphi\mapsto P_{s,D,\Omega}(z)\varphi:=f_D(z,\varphi),
\end{cases}    
\end{equation}
where $f_D(z,\varphi)$ is the unique solution of the boundary value problem \eqref{eqn:bvp}. Then for each 
$z\in\rho(A_{D,\Omega})$ and $s\in[0,1]$ the operator $\big[\gamma_N(A_{D,\Omega}-{\ol z}I)^{-1}\big]^*$, 
originally considered as in \eqref{3.2SD2}, induces a mapping
\begin{equation}\label{GCC-Jna.1}
\big[\gamma_N(A_{D,\Omega}-{\ol z}I)^{-1}\big]^* 
\in\cB\big(H^s(\partial\Omega),H^{s+(1/2)}(\Omega)\cap\dom(A_{max,\Omega})\big)
\end{equation}
{\rm (}where the space $H^{s+(1/2)}(\Omega)\cap\dom(A_{max,\Omega})$ is equipped with the natural 
norm $f\mapsto\|f\|_{H^{s+1/2}(\Omega)}+\|\Delta f\|_{L^2(\Omega)}${\rm )}, 
and 
\begin{equation}\label{eqn:P_D*}
P_{s,D,\Omega}(z)=-\big[\gamma_N(A_{D,\Omega}-{\ol z}I)^{-1}\big]^*\,\text{ on }\, 
H^s(\partial\Omega).
\end{equation}
Moreover, $P_{s,D,\Omega}(z)$ is injective with 
\begin{equation}\label{eqn:ranP_D} 
\ran(P_{s,D,\Omega}(z))=\ker(A_{max,\Omega}-zI)\cap H^{s+(1/2)}(\Omega).
\end{equation}
In particular, $\ran(P_{s,D,\Omega}(z))$ is dense in $\ker(A_{max,\Omega}-zI)$ 
with respect to the $L^2(\Omega)$-norm. \\[1mm] 
\noindent $(ii)$ For $z\in\rho(A_{N,\Omega})$ and $s\in[0,1]$, define  
\begin{equation}\label{eqn:P_N}
P_{s,N,\Omega}(z): 
\begin{cases} 
H^{s-1}(\partial\Omega)\rightarrow H^{s+(1/2)}(\Omega)\cap\dom(A_{max,\Omega}), 
\\[2pt]
\varphi\mapsto P_{s,N,\Omega}(z)\varphi:=f_N(z,\varphi),
\end{cases}     
\end{equation}
where $f_N(z,\varphi)$ is the unique solution of the boundary value problem \eqref{eqn:bvp2}. 
Then for each $z\in\rho(A_{N,\Omega})$ and $s\in[0,1]$ the operator 
$\big[\gamma_D(A_{N,\Omega}-{\ol z}I)^{-1}\big]^*$, initially viewed as in \eqref{3.2SD2N}, 
induces a mapping
\begin{equation}\label{GCC-Jna.1N}
\big[\gamma_D(A_{N,\Omega}-{\ol z}I)^{-1}\big]^* 
\in\cB\big(H^{s-1}(\partial\Omega),H^{s+(1/2)}(\Omega)\cap\dom(A_{max,\Omega})\big)
\end{equation}
{\rm (}where the space $H^{s+(1/2)}(\Omega)\cap\dom(A_{max,\Omega})$ is equipped with 
the natural norm $f\mapsto\|f\|_{H^{s+1/2}(\Omega)}+\|\Delta f\|_{L^2(\Omega)}${\rm )}, 
and 
\begin{equation}\label{eqn:P_N*}
P_{s,N,\Omega}(z)=-\big[\gamma_D(A_{N,\Omega}-{\ol z}I)^{-1}\big]^*\,\text{ on }\, 
H^{s-1}(\partial\Omega).    
\end{equation}
In addition, $P_{s,N,\Omega}(z)$ is injective with  
\begin{equation}\label{eqn:ranP_N} 
\ran(P_{s,N,\Omega}(z))=\ker(A_{max,\Omega}-z I)\cap H^{s+(1/2)}(\Omega).
\end{equation} 
In particular, $\ran(P_{s,N,\Omega}(z))$ is dense in $\ker(A_{max,\Omega}-zI)$ 
with respect to the $L^2(\Omega)$-norm. \\[1mm] 
$(iii)$ For $z\in\rho(A_{D,\Omega})$ and $s\in[0,1]$, the Dirichlet-to-Neumann 
operator defined by
\begin{equation}\label{We-Q.14}
M_{s,\Omega}(z):
\begin{cases} 
H^s(\partial\Omega)\rightarrow H^{s-1}(\partial\Omega),   
\\[2pt] 
\varphi\mapsto M_{s,\Omega}(z)\varphi:=-\gamma_N P_{s,D,\Omega}(z)\varphi,
\end{cases}    
\end{equation}
satisfies
\begin{equation}\label{4.17}
M_{s,\Omega}(z)=\gamma_N\big[\gamma_N(A_{D,\Omega}-{\ol z}I)^{-1}\big]^* 
\in\cB\big(H^s(\partial\Omega),H^{s-1}(\partial\Omega)\big).     
\end{equation}
Moreover, for each $z\in\rho(A_{D,\Omega})$ and each $s\in[0,1]$,
\begin{align}\label{4.17AD}
\begin{split}
& \text{the adjoint of $M_{s,\Omega}(z)\in\cB\big(H^s(\partial\Omega),H^{s-1}(\partial\Omega)\big)$ is}
\\[2pt]
& \quad\text{the operator 
$M_{1-s,\Omega}(\overline{z})\in\cB\big(H^{1-s}(\partial\Omega),H^{-s}(\partial\Omega)\big)$}.     
\end{split}
\end{align}
$(iv)$ For $z\in\rho(A_{N,\Omega})$ and $s\in[0,1]$, the Neumann-to-Dirichlet 
operator defined by
\begin{equation}\label{We-Q.15}
N_{s,\Omega}(z): 
\begin{cases} 
H^{s-1}(\partial\Omega)\rightarrow H^s(\partial\Omega),   
\\[2pt] 
\varphi\mapsto N_{s,\Omega}(z)\varphi:=-\gamma_D P_{s,N,\Omega}(z)\varphi,
\end{cases}    
\end{equation}
satisfies
\begin{equation}\label{4.19}
N_{s,\Omega}(z)=\gamma_D\big[\gamma_D(A_{N,\Omega}-{\ol z} I)^{-1}\big]^*  
\in\cB\big(H^{s-1}(\partial\Omega),H^s(\partial\Omega)\big).     
\end{equation} 
In addition, for each $z\in\rho(A_{N,\Omega})$ and each $s\in[0,1]$,
\begin{align}\label{4.17AN}
\begin{split}
& \text{the adjoint of $N_{s,\Omega}(z)\in\cB\big(H^{s-1}(\partial\Omega),H^s(\partial\Omega)\big)$ is}
\\[2pt]
& \quad\text{the operator $N_{1-s,\Omega}(\overline{z})\in\cB\big(H^{-s}(\partial\Omega),H^{1-s}(\partial\Omega)\big)$}.     
\end{split}
\end{align}
$(v)$ If $z\in\rho(A_{D,\Omega})\cap\rho(A_{N,\Omega})$, then for each $s\in[0,1]$ 
the Dirichlet-to-Neumann operator $M_{s,\Omega}(z)$ maps $H^{s}(\partial\Omega)$ bijectively 
onto $H^{s-1}(\partial\Omega)$, the Neumann-to-Dirichlet operator $N_{s,\Omega}(z)$ maps 
$H^{s-1}(\partial\Omega)$ bijectively onto $H^{s}(\partial\Omega)$, and their inverses satisfy 
\begin{align}\label{5.36JJJ.1}
& M_{s,\Omega}(z)^{-1}=-N_{s,\Omega}(z)\in\cB\big(H^{s-1}(\partial\Omega),H^{s}(\partial\Omega)\big),
\\[2pt] 
& N_{s,\Omega}(z)^{-1}=-M_{s,\Omega}(z)\in\cB\big(H^{s}(\partial\Omega),H^{s-1}(\partial\Omega)\big).
\label{5.36JJJ.2}
\end{align} 
\end{theorem}
%%%%%%%%%
\begin{proof}
Most of the claims in $(i)$--$(ii)$ follow from Lemmas~\ref{LamAA.2}--\ref{LamAA.2N} in a 
straightforward manner. For the membership in \eqref{GCC-Jna.1} one first observe that 
$[\gamma_N(A_{D,\Omega}-{\ol z}I)^{-1}]^*$ regarded as mapping from $H^s(\partial\Omega)$
to $H^{s+(1/2)}(\Omega)\cap\dom(A_{max,\Omega})$ (where the latter space is equipped with the 
norm $f\mapsto\|f\|_{H^{s+1/2}(\Omega)}+\|\Delta f\|_{L^2(\Omega)}$) is closed.
In fact, if $\{\varphi_j\}_{j\in{\mathbb{N}}}\subset H^s(\partial\Omega)$ is sequence which converges to 
$\varphi\in H^s(\partial\Omega)$ in the norm of $H^s(\partial\Omega)$ and 
\begin{equation}\label{convss}
\lim_{j\rightarrow\infty}\big[\gamma_N(A_{D,\Omega}-{\ol z}I)^{-1}\big]^*\varphi_j 
=\psi\in H^{s+(1/2)}(\Omega)\cap\dom(A_{max,\Omega})
\end{equation}
with respect to the graph norm on $H^{s+(1/2)}(\Omega)\cap\dom(A_{max,\Omega})$ 
then it follows that also $\varphi_j\rightarrow\varphi$ in $L^2(\partial\Omega)$
as $j\to\infty$ and the limit in \eqref{convss} exists also in $L^2(\Omega)$. 
Hence it follows from the boundedness of $[\gamma_N(A_{D,\Omega}-{\ol z}I)^{-1}]^*$ 
when regarded as a mapping from $L^2(\partial\Omega)$ to $L^2(\Omega)$ 
(see \eqref{3.2SD2}) that 
\begin{equation}
\big[\gamma_N(A_{D,\Omega}-{\ol z}I)^{-1}\big]^*\varphi
=\psi\in H^{s+(1/2)}(\Omega)\cap\dom(A_{max,\Omega}).
\end{equation}
Therefore,  
\begin{equation}\label{ssbound}
\big[\gamma_N(A_{D,\Omega}-{\ol z}I)^{-1}\big]^*:H^s(\partial\Omega)
\rightarrow H^{s+(1/2)}(\Omega)\cap\dom(A_{max,\Omega})
\end{equation}
is closed and defined on the whole space $H^s(\partial\Omega)$ by the well-posedness of the boundary 
value problem \eqref{eqn:bvp} and the representation \eqref{4.7}. This yields the boundedness of 
\eqref{ssbound} and hence shows \eqref{GCC-Jna.1}. The membership in \eqref{GCC-Jna.1N} follows 
from a similar reasoning, employing the well-posedness of \eqref{eqn:bvp2} and \eqref{4.9}. 
In addition, the fact that $\ran(P_{s,D,\Omega}(z))$ and $\ran(P_{s,N,\Omega}(z))$, $s\in [0,1]$, 
are dense in $\ker(A_{max,\Omega}-z I)$ with respect to the $L^2(\Omega)$-norm follows from 
combining Lemma~\ref{l5.1} with \eqref{eqn:ranP_D} and \eqref{eqn:ranP_N}. 

Next, the first claim in $(iii)$, that is,  \eqref{4.17}, follows from combining 
\eqref{eqn:gammaN-pp}, \eqref{GCC-Jna.1}--\eqref{eqn:P_D*}, and \eqref{We-Q.14}. 
To verify \eqref{4.17AD}, fix $z,z'\in\rho(A_{D,\Omega})$, $s\in[0,1]$, and pick
$\varphi_1\in H^{1-s}(\partial\Omega)$, $\varphi_2\in H^s(\partial\Omega)$, arbitrary. 
Then, noticing 
\begin{align}\label{eq:FFGG.gh}
\begin{split}
& P_{s,D,\Omega}(z)\varphi_2\in H^{s+(1/2)}(\Omega)\cap\dom(A_{max,\Omega}), 
\\[2pt] 
& P_{1-s,D,\Omega}(z')\varphi_1\in H^{(3/2)-s}(\Omega)\cap\dom(A_{max,\Omega}),
\end{split}
\end{align}
one observes that by design, 
\begin{align}\label{8uhg}
\begin{split}
& \gamma_D P_{s,D,\Omega}(z)\varphi_2=\varphi_2,\quad
\gamma_N P_{s,D,\Omega}(z)\varphi_2=-M_{s,\Omega}(z)\varphi_2,    
\\[2pt]
& \gamma_D P_{1-s,D,\Omega}(z')\varphi_1=\varphi_1,\quad
\gamma_N P_{1-s,D,\Omega}(z')\varphi_1=-M_{1-s,\Omega}(z')\varphi_1.  
\end{split}
\end{align}
As such, Green's identity \eqref{GGGRRR} applied to the functions from 
\eqref{eq:FFGG.gh} implies 
\begin{align}\label{5.40}
& {}_{H^{1-s}(\partial\Omega)}\big\langle\varphi_1, 
M_{s,\Omega}(z)\varphi_2\big\rangle_{H^{s-1}(\partial\Omega)}
-{}_{H^{-s}(\partial\Omega)}\big\langle M_{1-s,\Omega}(z')\varphi_1,\varphi_2
\big\rangle_{H^{s}(\partial\Omega)}
\nonumber\\[2pt] 
& \quad={}_{H^{-s}(\partial\Omega)}\big\langle\gamma_N P_{1-s,D,\Omega}(z')\varphi_1, 
\gamma_D P_{s,D,\Omega}(z)\varphi_2\big\rangle_{H^{s}(\partial\Omega)}
\nonumber\\[2pt] 
& \qquad 
-{}_{H^{1-s}(\partial\Omega)}\big\langle\gamma_D P_{1-s,D,\Omega}(z')\varphi_1, 
\gamma_N P_{s,D,\Omega}(z)\varphi_2\big\rangle_{H^{s-1}(\partial\Omega)}
\nonumber\\[2pt] 
& \quad=\big(P_{1-s,D,\Omega}(z')\varphi_1, 
A_{max,\Omega}P_{s,D,\Omega}(z)\varphi_2\big)_{L^2(\Omega)}  
\nonumber\\[2pt] 
& \qquad -\big(A_{max,\Omega}P_{1-s,D,\Omega}(z')\varphi_1, 
P_{s,D,\Omega}(z)\varphi_2\big)_{L^2(\Omega)}  
\nonumber\\[2pt]  
& \quad=\big(z-\ol{z'}\,\big)\big(P_{1-s,D,\Omega}(z')\varphi_1, 
P_{s,D,\Omega}(z)\varphi_2\big)_{L^2(\Omega)}.    
\end{align}
Specializing the above formula to the case when $z'={\ol z}$ then proves that 
for every $z\in\rho(A_{D,\Omega})$, every $s\in[0,1]$, and each 
$\varphi_1\in H^{1-s}(\partial\Omega)$, $\varphi_2\in H^s(\partial\Omega)$, one has  
\begin{equation}\label{eq:RGbb}
{}_{H^{1-s}(\partial\Omega)}\big\langle\varphi_1, M_{s,\Omega}(z)\varphi_2\big\rangle_{H^{s-1}(\partial\Omega)}
={}_{H^{-s}(\partial\Omega)}\big\langle M_{1-s,\Omega}({\ol z})\varphi_1,\varphi_2
\big\rangle_{H^{s}(\partial\Omega)}.
\end{equation}
In turn, this identity justifies the claim in \eqref{4.17AD}. The treatment of $(iii)$ 
is therefore complete and the claims in part $(iv)$ are handled in a similar fashion. 

Finally, we consider the claims made in part $(v)$. 
To this end, select $z\in\rho(A_{D,\Omega})\cap\rho(A_{N,\Omega})$ and fix $s\in[0,1]$. 
In addition, let $\psi\in H^{1-s}(\partial\Omega)$ and $\varphi\in H^{s-1}(\partial\Omega)$ be arbitrary.
Then, observing 
\begin{align}\label{eq:FFGG.gh.2}
\begin{split}
& P_{1-s,D,\Omega}({\ol z})\psi\in H^{(3/2)-s}(\Omega)\cap\dom(A_{max,\Omega}),    
\\[2pt] 
& P_{s,N,\Omega}(z)\varphi\in H^{s+(1/2)}(\Omega)\cap\dom(A_{max,\Omega}),
\end{split}
\end{align}
our definitions ensure that 
\begin{align}\label{8uhg.t4}
\begin{split} 
& \gamma_D P_{1-s,D,\Omega}({\ol z})\psi=\psi,\quad 
\gamma_N P_{1-s,D,\Omega}({\ol z})\psi=- M_{1-s,\Omega}({\ol z})\psi,   
\\[2pt] 
& \gamma_D P_{s,N,\Omega}(z)\varphi=-N_{s,\Omega}(z)\varphi,\quad
\gamma_N P_{s,N,\Omega}(z)\varphi=-\varphi.
\end{split} 
\end{align}
Keeping these facts in mind and employing Green's identity \eqref{GGGRRR} for 
the functions in \eqref{eq:FFGG.gh.2}, one concludes that 
\begin{align}\label{iugta}
& {}_{H^{-s}(\partial\Omega)}\big\langle M_{1-s,\Omega}(\ol z)\psi,  
N_{s,\Omega}(z)\varphi\big\rangle_{H^{s}(\partial\Omega)} 
\nonumber\\[2pt]
& \quad={}_{H^{-s}(\partial\Omega)}\big\langle 
\gamma_N P_{1-s,D,\Omega}(\ol z)\psi,\gamma_D P_{s,N,\Omega}(z)\varphi
\big\rangle_{H^{s}(\partial\Omega)}
\nonumber\\[2pt]
& \quad={}_{H^{1-s}(\partial\Omega)}\big\langle\gamma_D P_{1-s,D,\Omega}(\ol z)\psi, 
\gamma_N P_{s,N,\Omega}(z)\varphi\big\rangle_{H^{s-1}(\partial\Omega)}    
\nonumber\\[2pt] 
& \qquad +\big(P_{1-s,D,\Omega}(\ol z)\psi, 
A_{max,\Omega}P_{s,N,\Omega}(z)\varphi\big)_{L^2(\Omega)}    
\nonumber\\[2pt]
& \qquad -\big(A_{max,\Omega}P_{1-s,D,\Omega}(\ol z)\psi, 
P_{s,N,\Omega}(z)\varphi\big)_{L^2(\Omega)}      
\nonumber\\[2pt]
&\quad={}_{H^{1-s}(\partial\Omega)}\big\langle\psi,(-\varphi)\big\rangle_{H^{s-1}(\partial\Omega)}
\nonumber\\[2pt]
& \qquad +\,\big(P_{1-s,D,\Omega}(\ol z)\psi,zP_{s,N,\Omega}(z)\varphi\big)_{L^2(\Omega)}    
\nonumber\\[2pt]
& \qquad -\big(\ol zP_{1-s,D,\Omega}(\ol z)\psi,P_{s,N,\Omega}(z)\varphi\big)_{L^2(\Omega)}    
\nonumber\\[2pt]
& \quad={}_{H^{1-s}(\partial\Omega)}\big\langle\psi,(-\varphi)\big\rangle_{H^{s-1}(\partial\Omega)}.    
\end{align}
In view of \eqref{4.17AD}, \eqref{4.17AN}, and the arbitrariness of 
$\psi\in H^{1-s}(\partial\Omega)$ and $\varphi\in H^{s-1}(\partial\Omega)$, this further implies
\begin{align}\label{ytr4322.AA}
& M_{s,\Omega}(z)N_{s,\Omega}(z)=-I\in\cB\big(H^{s-1}(\partial\Omega)\big),    
\\[2pt] 
& N_{1-s,\Omega}({\ol z})M_{1-s,\Omega}({\ol z})=-I\in\cB\big(H^{1-s}(\partial\Omega)\big).
\label{ytr4322.BB}
\end{align}
Since $s\in[0,1]$ and $z\in\rho(A_{D,\Omega})\cap\rho(A_{N,\Omega})$ have been arbitrarily selected, 
all desired claims follow from \eqref{ytr4322.AA}--\eqref{ytr4322.BB}.
\end{proof}
%%%%%%%%%%%

In the next lemma we collect some important properties of the Dirichlet-to-Neumann map in 
the case $s=1$. In this case, for each $z\in\rho(A_{D,\Omega})$ we define
\begin{align}\label{5.34}
\begin{split}
& M_{\Omega}(z):=M_{1,\Omega}(z)\,\text{ as an unbounded operator on $L^2(\partial\Omega)$}
\\[2pt]
& \quad\text{with dense domain }\dom\!\big(M_{\Omega}(z)\big):=H^1(\partial\Omega).
\end{split}
\end{align} 

%%%%%%%%
\begin{lemma}\label{t5.4part1}
Assume Hypothesis~\ref{h4.2} and let $z\in\rho(A_{D,\Omega})\cap\rho(A_{N,\Omega})$. 
Then the operator $M_{\Omega}(z)$ maps $H^1(\partial\Omega)$ bijectively onto 
$L^2(\partial\Omega)$. One has
\begin{equation}\label{5.35}
M_{\Omega}(\ol z)=M_{\Omega}(z)^*
\end{equation}
where the adjoint is understood in $L^2(\partial\Omega)$, and 
\begin{equation}\label{5.36}
M_{\Omega}(z)^{-1}=-N_{1,\Omega}(z)\in\cB_{\infty}\big(L^2(\partial\Omega)\big).
\end{equation}
\end{lemma}
%%%%%%%%%
\begin{proof}
First, according to \eqref{4.19}, one has $N_{1,\Omega}(z)\in\cB(L^2(\partial\Omega),H^1(\partial\Omega))$ 
and hence $N_{1,\Omega}(z)\in\cB(L^2(\partial\Omega))$ for all $z\in\rho(A_{N,\Omega})$. Moreover, since $H^1(\dOm)$ 
embeds compactly into $L^2(\partial\Omega)$ it follows that $N_{1,\Omega}(z)\in\cB_{\infty}(L^2(\partial\Omega))$ 
for $z\in\rho(A_{N,\Omega})$. From \eqref{5.36JJJ.2} one obtains  
\begin{equation}\label{halali}
N_{1,\Omega}(z)=-M_{1,\Omega}(z)^{-1}
=-M_{\Omega}(z)^{-1},\quad z\in\rho(A_{D,\Omega})\cap\rho(A_{N,\Omega}),
\end{equation}
and hence one concludes assertion \eqref{5.36}. 

In order to prove \eqref{5.35} we verify the identity
\begin{equation}\label{nidjussi}
N_{1,\Omega}(\ol z)=N_{1,\Omega}(z)^*
\end{equation}
for $z\in\rho(A_{N,\Omega})$, where the adjoint is understood in 
$L^2(\partial\Omega)$. Pick $\varphi,\psi\in L^2(\partial\Omega)$
and notice that
\begin{equation}\label{eq:FFGG.gh-jussi}
P_{1,N,\Omega}(\ol z)\varphi,\,P_{1,N,\Omega}( z)\psi\in H^{3/2}(\Omega)\cap\dom(A_{max,\Omega}), 
\end{equation}
by \eqref{eqn:P_N} and
\begin{align}\label{8uhg-jussi2}
\begin{split}
& \gamma_N P_{1,N,\Omega}(\ol z)\varphi=-\varphi,\quad
\gamma_D P_{1,N,\Omega}(\ol z)\varphi=-N_{1,\Omega}(\ol z)\varphi\in H^1(\partial\Omega),    
\\[2pt]
& \gamma_N P_{1,N,\Omega}( z)\psi=-\psi,\quad
\gamma_D P_{1,N,\Omega}( z)\psi=-N_{1,\Omega}( z)\psi\in H^1(\partial\Omega);  
\end{split}
\end{align}
cf.~\eqref{We-Q.15}. From Green's identity \eqref{GGGRRR} one obtains  
\begin{align}\label{5.40-jussi}
& (\varphi,N_{1,\Omega}(z)\psi)_{L^2(\partial\Omega)}-( N_{1,\Omega}(\ol z)\varphi,\psi)_{L^2(\partial\Omega)}
\nonumber\\[2pt] 
& \quad={}_{H^{-1}(\partial\Omega)}\big\langle\varphi, 
N_{1,\Omega}(z)\psi\big\rangle_{H^1(\partial\Omega)}
-{}_{H^1(\partial\Omega)}\big\langle N_{1,\Omega}(\ol z)\varphi,\psi
\big\rangle_{H^{-1}(\partial\Omega)}
\nonumber\\[2pt]
& \quad={}_{H^{-1}(\partial\Omega)}\big\langle
\gamma_N P_{1,N,\Omega}(\ol z)\varphi,\gamma_D P_{1,N,\Omega}(z)\psi\big\rangle_{H^1(\partial\Omega)}
\nonumber\\[2pt]
& \qquad-{}_{H^1(\partial\Omega)}\big\langle\gamma_D P_{1,D,\Omega}(\ol z)\varphi,
\gamma_N P_{1,N,\Omega}(z)\psi\big\rangle_{H^{-1}(\partial\Omega)}
\nonumber\\[2pt]
& \quad=\big(P_{1,N,\Omega}(\ol z)\varphi,A_{max,\Omega}P_{1,N,\Omega}(z)\psi\big)_{L^2(\Omega)}  
\nonumber\\[2pt]
& \qquad-\big(A_{max,\Omega}P_{1,N,\Omega}(\ol z)\varphi,P_{1,N,\Omega}(z)\varphi\big)_{L^2(\Omega)}  
\nonumber\\[2pt] 
& \quad=\big(z-z\,\big)\big(P_{1,N,\Omega}(\ol z)\varphi, 
P_{1,N,\Omega}(z)\varphi\big)_{L^2(\Omega)}=0,   
\end{align}
which implies \eqref{nidjussi}. For $z\in\rho(A_{D,\Omega})\cap\rho(A_{N,\Omega})$ 
one then finally concludes with the help of \eqref{halali} and \eqref{nidjussi} that
\begin{equation}\label{nidjussi-2}
M_{\Omega}(\ol z)=-N_{1,\Omega}(\ol z)^{-1}=\big(-N_{1,\Omega}(z)^{*}\big)^{-1}
=\big(-N_{1,\Omega}(z)^{-1}\big)^{*}=M_{\Omega}(z)^*,
\end{equation}
where the adjoint is understood in $L^2(\partial\Omega)$.
\end{proof}
%%%%%%%%%

Next, from \eqref{eqn:P_D*}, \eqref{eqn:P_N*}, the resolvent identity,
and the self-adjointness of $A_{D,\Omega},A_{N,\Omega}$, the following useful relations 
on $L^2(\partial\Omega)$ may be deduced:
\begin{equation}
\begin{split}
P_{0,D,\Omega}(z) &=\big(I+(z-z')(A_{D,\Omega}-z I)^{-1}\big)P_{0,D,\Omega}(z'),
\quad\forall\,z,z'\in\rho(A_{D,\Omega}),    
\\[2pt] 
P_{1,N,\Omega}(z) &=\big(I+(z-z')(A_{N,\Omega}-z I)^{-1}\big)P_{1,N,\Omega}(z'),
\quad\forall\,z,z'\in\rho(A_{N,\Omega}).
\end{split}
\end{equation}
By \eqref{We-Q.13} (with $s=0,1$), \eqref{eqn:P_D*}, and \eqref{3.2SD2}, 
one infers that for each $z\in\rho(A_{D,\Omega})$ the operator 
$P_{1,D,\Omega}(z)$, originally defined on $H^1(\partial\Omega)$ and presently viewed   
as a densely defined operator on $L^2(\partial\Omega)$, has the bounded 
$L^2(\partial\Omega)$--$L^2(\Omega)$-closure
\begin{equation}\label{domjussi}
\ol{P_{1,D,\Omega}(z)}=P_{0,D,\Omega}(z)\in\cB\big(L^2(\partial\Omega),H^{1/2}(\Omega)\big) 
\subset\cB\big(L^2(\partial\Omega),L^2(\Omega)\big).
\end{equation}
As such, 
\begin{equation}\label{5.22}
P_{1,D,\Omega}(z)^*=P_{0,D,\Omega}(z)^*\in\cB\big(L^2(\Omega),L^2(\partial\Omega)\big), 
\quad\forall\,z\in\rho(A_{D,\Omega}). 
\end{equation}
In particular, we emphasize that 
\begin{align}\label{WACO.WTF.1}
& P_{0,D,\Omega}(z):L^2(\partial\Omega)\rightarrow H^{1/2}(\Omega)\hookrightarrow L^2(\Omega), 
\quad\forall\,z\in\rho(A_{D,\Omega}),    
\\[2pt] 
& P_{1,N,\Omega}(z):L^2(\partial\Omega)\rightarrow H^{3/2}(\Omega)\hookrightarrow L^2(\Omega), 
\quad\forall\,z\in\rho(A_{N,\Omega}),
\label{WACO.WTF.2}
\end{align} 
and it is in this sense that the adjoint symbol ${}^*$ is understood for 
$L^2(\Omega)$--$L^2(\partial\Omega)$ operators in \eqref{5.22}, as well as 
in the remainder of this and the following section. 

Next, we note that collectively Lemma~\ref{l5.1}, \eqref{eqn:ranP_D}, and \eqref{eqn:ranP_N}, imply 
\begin{align}\label{tfc6g4d4.A}
& \begin{split}
\big(\ker(P_{0,D,\Omega}(z)^*)\big)^\bot &=\overline{\ran(P_{0,D,\Omega}(z))}    
\\[2pt]
&=\ker(A_{max,\Omega}-z I),\quad\forall\,z\in\rho(A_{D,\Omega}),
\end{split}   
\\[2pt]
& \begin{split}
\big(\ker(P_{1,N,\Omega}(z)^*)\big)^\bot &=\overline{\ran(P_{1,N,\Omega}(z))}    
\\[2pt] 
&=\ker(A_{max,\Omega}-z I),\quad\forall\,z\in\rho(A_{N,\Omega}),
\end{split}
\label{tfc6g4d4.B}
\end{align}
which further yield the orthogonal decompositions 
\begin{equation}\label{eqn:odec}
\begin{split}
L^2(\Omega) &=\ker(P_{0,D,\Omega}(z)^*)\oplus\ker(A_{max,\Omega}-zI),\quad\forall\,z\in\rho(A_{D,\Omega}),    
\\[2pt] 
L^2(\Omega) &=\ker(P_{1,N,\Omega}(z)^*)\oplus\ker(A_{max,\Omega}-zI),\quad\forall\,z\in\rho(A_{N,\Omega}).
\end{split}
\end{equation}

%%%%%%%%%
%%%%%%%%%
\section{Maximal Extensions of the Dirichlet and Neumann Trace on Bounded Lipschitz Domains} 
\lb{s8} 
%%%%%%%%%
%%%%%%%%%

The main objective of this section is to extend the Dirichlet and Neumann trace operator by 
continuity onto the domain of the maximal operator $A_{max,\Omega}$, with $\Omega\subset\bbR^n$ a 
bounded Lipschitz domain. Again it will be assumed throughout this section that Hypothesis~\ref{h4.2} holds. 

The following trace spaces equipped with a suitable topology will play 
the key role in the extension procedure discussed below (cf. \cite{BM13}). 

%%%%%%%%
\begin{definition}\label{Tggg.55}
Assuming Hypothesis~\ref{h4.2}, consider the spaces
\begin{equation}\label{eqn:G0G1}
\mathscr{G}_D(\partial\Omega):=\ran\big(\gamma_D\big|_{\dom(A_{N,\Omega})}\big)
\,\text{ and }\,\mathscr{G}_N(\partial\Omega):=\ran\big(\gamma_N\big|_{\dom(A_{D,\Omega})}\big).
\end{equation} 
\end{definition}
%%%%%%%%

To get a better insight into the nature of the spaces just introduced we observe that in the case 
when $\Omega$ is smooth (e.g., $\Omega$ of class $ C^{1,r}$ for some $r>1/2$ will do; see \cite{GM11}) one has 
\begin{equation}\label{mentionalso}
\mathscr{G}_D(\partial\Omega)=H^{3/2}(\partial\Omega)\,\text{ and }\,
\mathscr{G}_N(\partial\Omega)=H^{1/2}(\partial\Omega).
\end{equation}
We also point out that, in the case when $\Omega$ is a bounded quasi-convex domain in the 
sense of \cite{GM11} (hence, in particular, if  $\Omega$ is a bounded convex open set, or 
a bounded domain of class $ C^{1,r}$ for some $r>1/2$) then the spaces in \eqref{eqn:G0G1} 
may be explicitly described as
\begin{align}\label{Tfccc.6gH}
\begin{split} 
\mathscr{G}_N(\partial\Omega) &=\big\{g\in L^2(\partial\Omega)\,\big|\, 
g\nu_j\in H^{1/2}(\partial\Omega),\, 1\leq j\leq n\big\},
\\[2pt]  
\mathscr{G}_D(\partial\Omega) &=\big\{g\in H^1(\partial\Omega)\,\big|\, 
\nabla_{tan}g\in\big[H^{1/2}(\partial\Omega)\big]^n\big\},
\end{split} 
\end{align} 
where the $\nu_j$'s are the components of the outward unit normal $\nu$, and $\nabla_{tan}$ 
is the tangential gradient on $\partial\Omega$ (see \cite{GM11} for a proof and further comments). 

Here we emphasize that in the more general class of arbitrary bounded Lipschitz domains in 
${\mathbb{R}}^n$ the descriptions in \eqref{mentionalso} and \eqref{Tfccc.6gH} are no longer valid 
(the root of the problem being the failure of the inclusions in \eqref{Df-H1}), though, the following inclusions 
continue to hold:
\begin{align}\label{Tfccc.6gH.WACO}
\begin{split} 
&\big\{g\in L^2(\partial\Omega)\,\big|\,g\nu_j\in H^{1/2}(\partial\Omega),\, 1\leq j\leq n\big\}
\subseteq\mathscr{G}_N(\partial\Omega), 
\\[2pt]  
&\big\{g\in H^1(\partial\Omega)\,\big|\,\nabla_{tan}g\in\big[H^{1/2}(\partial\Omega)\big]^n\big\}
\subseteq\mathscr{G}_D(\partial\Omega).
\end{split} 
\end{align}  

Returning to the mainstream discussion (in the setting of Hypothesis~\ref{h4.2}), 
from \eqref{eqn:gammaDs.2}, \eqref{We-Q.10EE-jussi}, \eqref{eqn:gammaN-pp}, and 
\eqref{We-Q.10EE-jussiNN} we remark that
\begin{equation}\label{We-Q.17}
\begin{split}
\mathscr{G}_D(\partial\Omega) &=\big\{\gamma_D f\,\big|\,
f\in H^{3/2}(\Omega)\cap\dom(A_{max,\Omega}),\,\gamma_N f=0\big\}\subset H^1(\partial\Omega),
\\[2pt]
\mathscr{G}_N(\partial\Omega) &=\big\{\gamma_N f\,\big|\, 
f\in H^{3/2}(\Omega)\cap\dom(A_{max,\Omega}),\,\gamma_D f=0\big\}\subset L^2(\partial\Omega).
\end{split}
\end{equation}
One also observes that \eqref{eqn:P_D*}, \eqref{eqn:P_N*}, and \eqref{eqn:G0G1} entail 
\begin{align}\label{eqn:ranP*}
\begin{split} 
& \ran\!\big(P_{0,D,\Omega}(z)^*\big)=\mathscr{G}_N(\partial\Omega),
\quad\forall\,z\in\rho(A_{D,\Omega}),
\\[2pt]
& \ran\!\big(P_{1,N,\Omega}(z)^*\big)=\mathscr{G}_D(\partial\Omega),
\quad\forall\,z\in\rho(A_{N,\Omega}).
\end{split}
\end{align}

%%%%%%%%%%%
\begin{lemma}\label{l5.3}
Assume Hypothesis~\ref{h4.2}. Then $\mathscr{G}_N(\partial\Omega)$ is a dense proper linear subspace of $L^2(\partial\Omega)$, 
while $\mathscr{G}_D(\partial\Omega)$ is a dense proper linear subspace of $H^1(\partial\Omega)$ 
{\rm (}hence also dense in $L^2(\partial\Omega)${\rm )}.
\end{lemma}
%%%%%%%%%%%
\begin{proof}
That $\mathscr{G}_N(\partial\Omega)$ is a proper linear subspace of $L^2(\partial\Omega)$ is seen from \eqref{We-Q.17}, 
\eqref{Tggg.55}, and \eqref{eqn:gammaN-pp.WACO.yall.222}, bearing in mind that (cf. \eqref{We-Q.10EE-jussi})
\begin{equation}\label{eq:414542.M.WACO.1}
\dom(A_{D,\Omega})=H^1_{\Delta}(\Omega)\cap\accentset{\circ}{H}^1(\Omega).
\end{equation}
Likewise, that $\mathscr{G}_D(\partial\Omega)$ is a proper linear subspace of $H^1(\partial\Omega)$ is seen from \eqref{We-Q.17}, 
\eqref{Tggg.55}, and \eqref{eqn:gammaN-pp.WACO.yall.222.D}, bearing in mind that (cf. \eqref{We-Q.10EE-jussiNN})
\begin{equation}\label{eq:414542.M.WACO.2}
\dom(A_{N,\Omega})=\big\{u\in H^1_{\Delta}(\Omega)|\,\gamma_Nu=0\big\}.
\end{equation}

There remains to deal with the density claimed in the statement. To this end, suppose that the function 
$\phi\in L^2(\partial\Omega)$ is orthogonal to the subspace $\mathscr{G}_N(\partial\Omega)$ of $L^2(\partial\Omega)$. 
In view of \eqref{We-Q.17} this implies 
\begin{equation}\label{eq:tr3dd}
\big(\phi,\gamma_N f\big)_{L^2(\partial\Omega)}=0\,\text{ for all }\,
f\in H^{3/2}(\Omega)\cap\dom(A_{max,\Omega})\,\text{ with }\,\gamma_D f=0.
\end{equation}
Using the fact that $\gamma_D$ in \eqref{eqn:gammaDs.2} with $s=1/2$ is surjective, 
it follows that there exists 
\begin{equation}\label{eq:tr3dd.2}
g\in H^{1/2}(\Omega)\cap\dom(A_{max,\Omega})\,\text{ with }\,\gamma_D g=\phi.
\end{equation}
Hence, for each $f\in H^{3/2}(\Omega)\cap\dom(A_{max,\Omega})$ with $\gamma_D f=0$,
Green's formula \eqref{GGGRRR} yields 
\begin{align}\label{GGGRann}
0 &=\big(\phi,\gamma_N f\big)_{L^2(\partial\Omega)}
=\big(\gamma_D g,\gamma_N f\big)_{L^2(\partial\Omega)}
\nonumber\\[2pt] 
&={}_{(H^{1}(\partial\Omega))^*}\big\langle\gamma_N g,\gamma_D f\big\rangle_{H^{1}(\partial\Omega)}
+(g,\Delta f)_{L^2(\Omega)}-(\Delta g,f)_{L^2(\Omega)}
\nonumber\\[2pt] 
&=(g,\Delta f)_{L^2(\Omega)}-(\Delta g,f)_{L^2(\Omega)}.
\end{align}
By \eqref{We-Q.10EE-jussi}, one can rephrase the above condition as 
\begin{equation}\label{eq:4ffr88}
\big(g,A_{D,\Omega}f\big)_{L^2(\Omega)}
=\big((-\Delta+V)g,f\big)_{L^2(\Omega)},\quad\forall\,f\in\dom(A_{D,\Omega}),
\end{equation}
which, in turn, forces $g\in\dom(A_{D,\Omega}^*)$ and hence $g\in\dom(A_{D,\Omega})$
by the self-adjointness of $A_{D,\Omega}$ (cf.\  Theorem~\ref{t4.4}). As a consequence 
of this membership, \eqref{eq:tr3dd.2}, and \eqref{We-Q.10EE-jussi}, one obtains $\phi=\gamma_Dg=0$. 
This ultimately proves that the space $\mathscr{G}_N(\partial\Omega)$ is dense in $L^2(\partial\Omega)$.

Next, assume that the functional $\psi\in H^{-1}(\partial\Omega)=\big(H^1(\partial\Omega)\big)^*$ 
annihilates the subspace $\mathscr{G}_D(\partial\Omega)$ of $H^1(\partial\Omega)$. 
By\eqref{We-Q.17}, this translates into  
\begin{align}\label{eq:tr3dd.CC} 
\begin{split} 
& {}_{(H^1(\partial\Omega))^*}\big\langle\psi,\gamma_D f\big\rangle_{H^1(\partial\Omega)}=0
\,\text{ for all functions}    
\\[2pt] 
& \quad f\in H^{3/2}(\Omega)\cap\dom(A_{max,\Omega})\,\text{ with }\,\gamma_N f=0.
\end{split} 
\end{align}
Given that $\gamma_N$ in \eqref{eqn:gammaN-pp} with $s=1/2$ is surjective, one concludes that there exists 
\begin{equation}\label{eq:tr3dd.2CC}
g\in H^{1/2}(\Omega)\cap\dom(A_{max,\Omega})\,\text{ with }\,\gamma_N g=\psi.
\end{equation}
As such, for each $f\in H^{3/2}(\Omega)\cap\dom(A_{max,\Omega})$ with $\gamma_N f=0$,
Green's formula \eqref{GGGRRR} allows us to write 
\begin{align}\label{GGGRann.CC}
0 &={}_{(H^1(\partial\Omega))^*}\big\langle\psi,\gamma_D f\big\rangle_{H^1(\partial\Omega)}
={}_{(H^1(\partial\Omega))^*}\big\langle\gamma_N g,\gamma_D f\big\rangle_{H^1(\partial\Omega)}
\nonumber\\[2pt] 
&=\big(\gamma_D g,\gamma_N f\big)_{L^2(\partial\Omega)}
-(g,\Delta f)_{L^2(\Omega)}+(\Delta g,f)_{L^2(\Omega)}
\nonumber\\[2pt] 
&=-(g,\Delta f)_{L^2(\Omega)}+(\Delta g,f)_{L^2(\Omega)}.
\end{align}
By virtue of \eqref{We-Q.10EE-jussiNN}, this may be rephrased as 
\begin{equation}\label{eq:4ffr88.CC}
\big(g,A_{N,\Omega}f\big)_{L^2(\Omega)}
=\big((-\Delta+V)g,f\big)_{L^2(\Omega)}, 
\quad\forall\,f\in\dom(A_{N,\Omega}),
\end{equation}
which further  entails $g\in\dom(A_{N,\Omega}^*)$. Thus, $g\in\dom(A_{N,\Omega})$
by the self-adjointness of $A_{N,\Omega}$ (cf.\  Theorem~\ref{t4.5}). 
This fact, \eqref{eq:tr3dd.2CC}, and \eqref{We-Q.10EE-jussiNN} imply $\psi=\gamma_N g=0$. 
By the Hahn--Banach theorem, this proves that the space $\mathscr{G}_D(\partial\Omega)$ 
is dense in $H^1(\partial\Omega)$ (hence also dense in $L^2(\partial\Omega)$).
\end{proof}
%%%%%%%%%%%

In the next theorem we list some important properties of the imaginary 
part of the 
Dirichlet-to-Neumann map and its inverse in the case $s=1$. For this purpose, we recall (cf.~\eqref{5.34})  
that we employ the notation $M_{\Omega}(z):=M_{1,\Omega}(z)$ for $z\in\rho(A_{D,\Omega})$.

%%%%%%%%%%%
\begin{theorem}\label{t5.4}
Assume Hypothesis~\ref{h4.2}. Then the following assertions hold: \\[1mm] 
\noindent $(i)$ If  $z\in\bbC_+$ $($resp., $z\in\bbC_-$$)$ then 
\begin{align}\label{We-Q.19}
\begin{split} 
& \,\Im(M_{\Omega}(z)):=\frac{1}{2i}\big(M_{\Omega}(z)-M_{\Omega}(\bar{z})\big) 
=\Im(z)\,P_{1,D,\Omega}(z)^*P_{1,D,\Omega}(z),   
\\[2pt]   
& \dom\big(\Im(M_{\Omega}(z))\big):=H^1(\partial\Omega),
\end{split} 
\end{align}
is a densely defined bounded operator in $L^2(\partial\Omega)$ with bounded closure 
\begin{equation}\label{We-EE.1}
\overline{\Im(M_{\Omega}(z))}=\Im(z)\,P_{0,D,\Omega}(z)^* P_{0,D,\Omega}(z) 
\in\cB\big(L^2(\partial\Omega)\big). 
\end{equation} 
In addition, $\overline{\Im(M_{\Omega}(z))}$ is a 
nonnegative $($resp., nonpositive\,$)$ self-adjoint operator in 
$L^2(\partial\Omega)$ which is invertible with an unbounded inverse. \\[1mm] 
$(ii)$ If $z\in\bbC_+$ $($resp., $z\in\bbC_-$$)$ then 
\begin{equation}\label{eqn:imMinv}
\Im\big(-M_{\Omega}(z)^{-1}\big)=\Im(z)\,P_{1,N,\Omega}(z)^* P_{1,N,\Omega}(z) 
\in\cB\big(L^2(\partial\Omega)\big), 
\end{equation} 
is a nonnegative $($resp., nonpositive\,$)$, bounded, self-adjoint operator 
in $L^2(\partial\Omega)$ which is invertible with an unbounded inverse. 
\end{theorem}
%%%%%%%%%%
\begin{proof}
Concerning $(i)$, one observes that the same argument as in equation \eqref{5.40} implies 
\begin{equation}\label{i654rfd4e} 
M_{\Omega}(z)-M_{\Omega}(z')^*=\big(z-\ol{z'}\,\big)P_{0,D,\Omega}(z')^*P_{1,D,\Omega}(z) 
\end{equation} 
for every $z,z'\in\rho(A_{D,\Omega})$. Setting $z=z'$ and taking into account \eqref{5.22} 
and \eqref{5.35} yields \eqref{We-Q.19}. Next, fix $z\in\rho(A_{D,\Omega})$,  then 
\begin{equation}\label{olala}
\Im(M_{\Omega}(z))=\Im(z)\,P_{1,D,\Omega}(z)^*P_{1,D,\Omega}(z)=\Im(z)\,P_{0,D,\Omega}(z)^*P_{1,D,\Omega}(z)
\end{equation}
(see \eqref{5.22}) together with \eqref{domjussi} yields 
\begin{equation}\label{nYYH}
\overline{\Im(M_{\Omega}(z))}=\Im(z)\,P_{0,D,\Omega}(z)^* P_{0,D,\Omega}(z) 
\in\cB\big(L^2(\partial\Omega)\big),
\end{equation} 
which goes to show that for each $z\in\bbC_{+}$ (resp., each $z\in\bbC_{-}$) the bounded operator 
$\overline{\Im(M_{\Omega}(z))}$ is nonnegative (resp., nonpositive) and self-adjoint in 
$L^2(\partial\Omega)$.

Next, fix $z\in\bbC_{-}\cup\bbC_{+}$. According to Lemma~\ref{l5.3}, the space 
$\mathscr{G}_N(\partial\Omega)$ is dense in $L^2(\partial\Omega)$ hence one obtains from 
\begin{equation}\label{yr544}
\ker(P_{0,D,\Omega}(z))=\big(\ran(P_{0,D,\Omega}(z)^*)\big)^\bot, 
\quad\ran(P_{0,D,\Omega}(z)^*)=\mathscr{G}_N(\partial\Omega),
\end{equation}
(cf.~\eqref{eqn:ranP*}), that 
\begin{align}\label{766GFF}
\ker\big(\overline{\Im(M_{\Omega}(z))}\big)
&=\ker\big(P_{0,D,\Omega}(z)^* P_{0,D,\Omega}(z)\big)=\ker\big(P_{0,D,\Omega}(z)\big)  
\nonumber\\[2pt] 
&=\big(\ran(P_{0,D,\Omega}(z)^*)\big)^\bot=\mathscr{G}_N(\partial\Omega)^\bot=\{0\}.
\end{align}
Thus, $\overline{\Im(M_{\Omega}(z))}$ is injective. From the representation \eqref{nYYH} 
and the second identity in \eqref{yr544} it follows that the inclusion
\begin{equation}\label{ranm1}
\ran\big(\overline{\Im(M_{\Omega}(z))}\big)\subset\mathscr{G}_N(\partial\Omega)
\end{equation}
holds. As the operator $\overline{\Im(M_{\Omega}(z))}$ is self-adjoint, one  concludes that 
its range is a dense subspace of $L^2(\partial\Omega)$ and from \eqref{ranm1} and Lemma~\ref{l5.3} it is clear 
that the range of $\overline{\Im(M_{\Omega}(z))}$ is a proper subspace of $L^2(\partial\Omega)$. 
Hence, the inverse is an unbounded operator in $L^2(\partial\Omega)$.

Finally, item $(ii)$ follows in the same way as item $(i)$ by interchanging the roles of 
$M_{\Omega}(z)$ and $-M_{\Omega}(z)^{-1}$, $P_{0,D,\Omega}(z)$ and $P_{1,N,\Omega}(z)$, $\gamma_D$ 
and  $-\gamma_N$, $A_{D,\Omega}$ and $A_{N,\Omega}$, and $\mathscr{G}_N(\partial\Omega)$ 
and $\mathscr{G}_D(\partial\Omega)$.   
\end{proof}
%%%%%%%%%%%

The following theorem builds on \cite{BM13}, \cite{GM08}, \cite{GM11} under various assumptions 
on the underlying domain and the regularity of functions involved. Here we now present the most 
general PDE result in this spirit. The notion of equivalence of norms in different Banach spaces 
used in item (vi) of Theorem~\ref{t5.5} is explained in Lemma~\ref{uyrre} below.

%%%%%%%%%%%%
\begin{theorem}\label{t5.5}
Assume Hypothesis~\ref{h4.2} and consider
\begin{equation}\label{eqn:SigmaLambda}
\Sigma:=\Im\big(-M_{\Omega}(i)^{-1}\big),\quad\Lambda:=\overline{\Im(M_{\Omega}(i))},
\end{equation}
which, according to Theorem~\ref{t5.4}, are bounded, nonnegative, self-adjoint 
operators in $L^2(\partial\Omega)$, that are invertible, with unbounded inverses. 
Then the following statements hold: \\[1mm] 
$(i)$ One has 
\begin{equation}\label{eqn:domLambda}
\begin{split}
\mathscr{G}_D(\partial\Omega)&=\dom\big(\Sigma^{-1/2}\big)=\ran\big(\Sigma^{1/2}\big), 
\\[2pt]
\mathscr{G}_N(\partial\Omega)&=\dom\big(\Lambda^{-1/2}\big)=\ran\big(\Lambda^{1/2}\big),
\end{split}
\end{equation} 
and when equipped with the scalar products 
\begin{equation}\label{scalarjussi}
\begin{split}
(\varphi,\psi)_{\mathscr{G}_D(\partial\Omega)} 
&:=\big(\Sigma^{-1/2}\varphi,\Sigma^{-1/2}\psi\big)_{L^2(\partial\Omega)},
\quad\forall\,\varphi,\psi\in\mathscr{G}_D(\partial\Omega),    
\\[2pt] 
(\varphi,\psi)_{\mathscr{G}_N(\partial\Omega)} 
&:=\big(\Lambda^{-1/2}\varphi,\Lambda^{-1/2}\psi\big)_{L^2(\partial\Omega)},
\quad\forall\,\varphi,\psi\in\mathscr{G}_N(\partial\Omega),
\end{split}
\end{equation}
the spaces $\mathscr{G}_D(\partial\Omega),\mathscr{G}_N(\partial\Omega)$ become Hilbert spaces. 
\\[1mm] 
$(ii)$ The Dirichlet trace operator $\gamma_D$ {\rm (}as defined in \eqref{eqn:gammaDs.2}{\rm )} 
and the Neumann trace operator $\gamma_N$ {\rm (}as defined in \eqref{eqn:gammaN-pp}{\rm )} 
extend by continuity {\rm (}hence in a compatible manner{\rm )} to continuous surjective mappings 
\begin{equation}\label{4.55}
\begin{split}
\widetilde{\gamma}_D:\dom(A_{max,\Omega})&\to\mathscr{G}_N(\partial\Omega)^*,
\\[2pt]  
\widetilde{\gamma}_N:\dom(A_{max,\Omega})&\to\mathscr{G}_D(\partial\Omega)^*,  
\end{split}
\end{equation}
where $\dom(A_{max,\Omega})$ is endowed with the graph norm of $A_{max,\Omega}$, 
and $\mathscr{G}_D(\partial\Omega)^*$, $\mathscr{G}_N(\partial\Omega)^*$ are, respectively, 
the adjoint {\rm (}conjugate dual{\rm )} spaces of $\mathscr{G}_D(\partial\Omega)$,
$\mathscr{G}_N(\partial\Omega)$ carrying the natural topology induced by 
\eqref{scalarjussi} on $\mathscr{G}_D(\partial\Omega)$, 
$\mathscr{G}_N(\partial\Omega)$, respectively, such that
\begin{equation}\label{4.30}
\ker(\widetilde{\gamma}_D)=\dom(A_{D,\Omega})\,\text{ and }\,   
\ker(\widetilde{\gamma}_N)=\dom(A_{N,\Omega}).
\end{equation} 
Furthermore, for each $s\in[0,1]$ there exists a constant $C\in(0,\infty)$ with the property that
\begin{align}\label{4.55.REE.D}
\begin{split}
& f\in\dom(A_{max,\Omega})\,\text{ and }\,\widetilde{\gamma}_D f\in H^s(\partial\Omega)
\,\text{ imply }\,f\in H^{s+(1/2)}(\Omega)
\\[2pt]
& \quad\text{and }\,\|f\|_{H^{s+(1/2)}(\Omega)}\leq C\big(\|\Delta f\|_{L^2(\Omega)}
+\|\widetilde{\gamma}_D f\|_{H^s(\partial\Omega)}\big),
\end{split}
\end{align}
and
\begin{align}\label{4.55.REE.N}
\begin{split}
& f\in\dom(A_{max,\Omega})\,\text{ and }\,\widetilde{\gamma}_N f\in H^{-s}(\partial\Omega)
\,\text{ imply }\,f\in H^{-s+(3/2)}(\Omega)
\\[2pt]
& \quad\text{and }\,\|f\|_{H^{-s+(3/2)}(\Omega)}\leq C\big(\|f\|_{L^2(\Omega)}+\|\Delta f\|_{L^2(\Omega)}
+\|\widetilde{\gamma}_N f\|_{H^{-s}(\partial\Omega)}\big).
\end{split}
\end{align} 
$(iii)$ With $\widetilde{\gamma}_D,\widetilde{\gamma}_N$ as in \eqref{4.55}, one has
{\rm (}compare to \eqref{Tan-C3}{\rm )}
\begin{align}\label{ut444}
\accentset{\circ}{H}^2(\Omega)=\big\{f\in\dom(A_{max,\Omega})\,\big|\,
& \wti{\gamma}_D f=0\,\text{ in }\,\mathscr{G}_N(\partial\Omega)^*
\nonumber\\[2pt]
&\text{and }\,\wti{\gamma}_N f=0\,\text{ in }\,\mathscr{G}_D(\partial\Omega)^*\big\}.  
\end{align}
$(iv)$ The manner in which the mapping $\widetilde{\gamma}_D$ in \eqref{4.55} 
operates is as follows: Given $f\in\dom(A_{max,\Omega})$, the action of the functional 
$\widetilde{\gamma}_D f\in\mathscr{G}_N(\partial\Omega)^*$ on some arbitrary 
$\phi\in\mathscr{G}_N(\partial\Omega)$ is given by 
\begin{equation}\label{eq:33f4iU}
{}_{\mathscr{G}_N(\partial\Omega)^\ast}\big\langle\widetilde{\gamma}_D f,
\phi\big\rangle_{\mathscr{G}_N(\partial\Omega)}
=(f,\Delta g)_{L^2(\Omega)}-(\Delta f,g)_{L^2(\Omega)},
\end{equation}
for any $g\in H^{3/2}(\Omega)\cap\dom(A_{max,\Om})$ such that $\gamma_D g=0$ 
and $\gamma_N g=\phi$ {\rm (}the existence of such $g$ being ensured by \eqref{We-Q.17}{\rm )}.
As a consequence, the following Green's formula holds: 
\begin{equation}\label{eq:33VaV}
{}_{\mathscr{G}_N(\partial\Omega)^\ast}\big\langle\widetilde{\gamma}_D f,
\gamma_N g\big\rangle_{\mathscr{G}_N(\partial\Omega)}
=(f,\Delta g)_{L^2(\Omega)}-(\Delta f,g)_{L^2(\Omega)},
\end{equation}
for each $f\in\dom(A_{max,\Omega})$ and each $g\in\dom(A_{D,\Om})$. \\[1mm] 
$(v)$ The mapping $\widetilde{\gamma}_N$ in \eqref{4.55} 
operates in the following fashion: Given a function $f\in\dom(A_{max,\Omega})$, the 
action of the functional $\widetilde{\gamma}_N f\in\mathscr{G}_D(\partial\Omega)^*$ 
on some arbitrary $\psi\in\mathscr{G}_D(\partial\Omega)$ is given by 
\begin{equation}\label{eq:33f4iU.NN}
{}_{\mathscr{G}_D(\partial\Omega)^\ast}\big\langle\widetilde{\gamma}_N f,
\psi\big\rangle_{\mathscr{G}_D(\partial\Omega)}
=-(f,\Delta g)_{L^2(\Omega)}+(\Delta f,g)_{L^2(\Omega)},
\end{equation}
for any $g\in H^{3/2}(\Omega)\cap\dom(A_{max,\Om})$ such that $\gamma_N g=0$ 
and $\gamma_D g=\psi$ {\rm (}the existence of such $g$ being ensured by \eqref{We-Q.17}{\rm )}.
In particular, the following Green's formula holds: 
\begin{equation}\label{eq:33VaV.NN}
{}_{\mathscr{G}_D(\partial\Omega)^\ast}\big\langle\widetilde{\gamma}_N f,
\gamma_D g\big\rangle_{\mathscr{G}_D(\partial\Omega)}
=-(f,\Delta g)_{L^2(\Omega)}+(\Delta f,g)_{L^2(\Omega)},
\end{equation}
for each $f\in\dom(A_{max,\Omega})$ and each $g\in\dom(A_{N,\Om})$. \\[1mm] 
$(vi)$ The operators 
\begin{align}\label{4.5ASdf.1}
&\gamma_D:\dom(A_{N,\Om})=H^{3/2}(\Omega)\cap\dom(A_{max,\Omega})\cap\ker(\gamma_N)
\to\mathscr{G}_D(\partial\Omega),
\\[2pt] 
&\gamma_N:\dom(A_{D,\Om})=H^{3/2}(\Omega)\cap\dom(A_{max,\Omega})\cap\ker(\gamma_D)
\to\mathscr{G}_N(\partial\Omega),    
\label{4.5ASdf.2}
\end{align}
are well defined, linear, surjective, and continuous if for some $s\in[0,\tfrac{3}{2}]$ both 
spaces on the left-hand sides of \eqref{4.5ASdf.1}, \eqref{4.5ASdf.2} are equipped with the 
norm $f\mapsto\|f\|_{H^s(\Omega)}+\|\Delta f\|_{L^2(\Omega)}$ {\rm (}which are all equivalent; 
cf.~\eqref{eq:UUjh-jussi} and \eqref{eq:UUjh-jussiNN}{\rm )}. In addition, 
\begin{equation}\label{eq:Ajf}
\text{the kernel of $\gamma_D$ and $\gamma_N$ in \eqref{4.5ASdf.1}--\eqref{4.5ASdf.2} 
is $\accentset{\circ}{H}^2(\Omega)$}.
\end{equation}
Moreover, 
\begin{align}\label{i7g5r}
\|\phi\|_{\mathscr{G}_D(\partial\Omega)}& \approx\inf_{\substack{f\in
H^{3/2}(\Omega)\cap\dom(A_{max,\Omega})\\ \gamma_N f=0,\,\,\gamma_D f=\phi}}
\big(\|f\|_{H^{3/2}(\Omega)}+\|\Delta f\|_{L^2(\Omega)}\big)
\nonumber\\[2pt]
& \approx\inf_{\substack{f\in
H^{3/2}(\Omega)\cap\dom(A_{max,\Omega})\\ \gamma_N f=0,\,\,\gamma_D f=\phi}}
\big(\|f\|_{L^2(\Omega)}+\|\Delta f\|_{L^2(\Omega)}\big)
\nonumber\\[2pt]
& \approx\inf_{\substack{f\in\dom(A_{max,\Omega})\\ \widetilde{\gamma}_N f=0,\,\,\widetilde{\gamma}_D f=\phi}}
\big(\|f\|_{L^2(\Omega)}+\|\Delta f\|_{L^2(\Omega)}\big),
\end{align}
uniformly for $\phi\in\mathscr{G}_D(\partial\Omega)$, and 
\begin{align}\label{i7g5r.2N}
\|\psi\|_{\mathscr{G}_N(\partial\Omega)}& \approx\inf_{\substack{g\in
H^{3/2}(\Omega)\cap\dom(A_{max,\Omega})\\ \gamma_D g=0,\,\,\gamma_N g=\psi}}
\big(\|g\|_{H^{3/2}(\Omega)}+\|\Delta g\|_{L^2(\Omega)}\big)
\nonumber\\[2pt]
& \approx\inf_{\substack{g\in
H^{3/2}(\Omega)\cap\dom(A_{max,\Omega})\\ \gamma_D g=0,\,\,\gamma_N g=\psi}}
\big(\|g\|_{L^2(\Omega)}+\|\Delta g\|_{L^2(\Omega)}\big)
\nonumber\\[2pt]
& \approx\inf_{\substack{g\in\dom(A_{max,\Omega})\\ \widetilde\gamma_D g=0,\,\,\widetilde\gamma_N g=\psi}}
\big(\|g\|_{L^2(\Omega)}+\|\Delta g\|_{L^2(\Omega)}\big)
\nonumber\\[2pt]
& \approx\inf_{\substack{g\in\dom(A_{max,\Omega})\\ \widetilde{\gamma}_D g=0,\,\,\widetilde{\gamma}_N g=\psi}}
\|\Delta g\|_{L^2(\Omega)},
\end{align}
uniformly for $\psi\in\mathscr{G}_N(\partial\Omega)$.

As a consequence,
\begin{align}\label{ytrr555e}
\begin{split}
& \mathscr{G}_D(\partial\Omega)\hookrightarrow H^1(\partial\Omega)
\hookrightarrow L^2(\partial\Omega)\hookrightarrow H^{-1}(\partial\Omega)
\hookrightarrow\mathscr{G}_D(\partial\Omega)^*,  
\\[2pt] 
& \mathscr{G}_N(\partial\Omega)\hookrightarrow L^2(\partial\Omega)
\hookrightarrow\mathscr{G}_N(\partial\Omega)^*, 
\end{split}
\end{align} 
with all embeddings linear, continuous, and with dense range. Moreover, the duality 
pairings between $\mathscr{G}_D(\partial\Omega)$ and $\mathscr{G}_D(\partial\Omega)^\ast$,
as well as between $\mathscr{G}_N(\partial\Omega)$ and $\mathscr{G}_N(\partial\Omega)^\ast$, 
are both compatible with the inner product in $L^2(\partial\Omega)$. \\[1mm] 
\noindent $(vii)$ For each $z\in\rho(A_{D,\Omega})$, the boundary value problem 
\begin{equation}\label{eqn:bvp.2}
\begin{cases} 
(-\Delta+V-z)f=0\,\text{ in $\Omega,\quad f\in\dom(A_{max,\Omega})$,}   
\\[2pt]  
\widetilde{\gamma}_D f=\varphi
\,\text{ in $\mathscr{G}_N(\partial\Omega)^*,\quad\varphi\in\mathscr{G}_N(\partial\Omega)^*$,}
\end{cases} 
\end{equation}
is well posed. In particular, for each $z\in\rho(A_{D,\Omega})$ there exists a constant $C\in(0,\infty)$, which depends 
only on $\Omega$, $n$, $z$, and $V$, with the property that
\begin{align}\label{eqn:bvp.2WACO}
\begin{split}
& \|f\|_{L^2(\Omega)}\leq C\|\widetilde{\gamma}_D f\|_{\mathscr{G}_N(\partial\Omega)^*}\,\text{ for each }\, f\in\dom(A_{max,\Omega})
\\[2pt]
& \quad\text{with }\,(-\Delta+V-z)f=0\,\text{ in }\,\Omega.
\end{split} 
\end{align}
Moreover, if 
\begin{equation}\label{6h8u}
\widetilde P_{D,\Omega}(z): 
\begin{cases} 
\mathscr{G}_N(\partial\Omega)^*\rightarrow\dom(A_{max,\Omega}),   
\\[2pt] 
\varphi\mapsto\widetilde P_{D,\Omega}(z)\varphi:=\wti{f}_{D,\Omega}(z,\varphi),
\end{cases} 
\end{equation}
where $\wti{f}_{D,\Omega}(z,\varphi)$ is the unique solution of \eqref{eqn:bvp.2}, 
then the solution operator $\widetilde{P}_{D,\Omega}(z)$ is an extension of 
$P_{0,D,\Omega}(z)$ in \eqref{We-Q.13}, and  $\widetilde{P}_{D,\Omega}(z)$ is continuous, when 
the adjoint space $\mathscr{G}_N(\partial\Omega)^*$ and $\dom(A_{max,\Omega})$ are endowed 
with the norms in item $(ii)$. \\[1mm] 
\noindent $(viii)$ For each $z\in\rho(A_{N,\Omega})$, the boundary value problem 
\begin{equation}\label{eqn:bvpNeumann}
\begin{cases} 
(-\Delta+V-z)f=0\,\text{ in $\Omega,\quad f\in\dom(A_{max,\Omega})$,}   
\\[2pt]  
-\widetilde{\gamma}_N f=\varphi 
\,\text{ in $\mathscr{G}_D(\partial\Omega)^*,\quad\varphi\in\mathscr{G}_D(\partial\Omega)^*$,}
\end{cases} 
\end{equation}
is well posed. In particular, for each $z\in\rho(A_{N,\Omega})$ there exists a constant $C\in(0,\infty)$, which depends 
only on $\Omega$, $n$, $z$, and $V$, with the property that
\begin{align}\label{eqn:bvpNeumann.WACO}
\begin{split}
& \|f\|_{L^2(\Omega)}\leq C\|\widetilde{\gamma}_N f\|_{\mathscr{G}_D(\partial\Omega)^*}\,\text{ for each }\,f\in\dom(A_{max,\Omega})     
\\[2pt]
& \quad\text{with }\,(-\Delta+V-z)f=0\,\text{ in }\,\Omega.
\end{split} 
\end{align}
Moreover, if 
\begin{equation}\label{755fv}
\widetilde{P}_{N,\Omega}(z): 
\begin{cases} 
\mathscr{G}_D(\partial\Omega)^*\rightarrow\dom(A_{max,\Omega}),  
\\[2pt] 
\varphi\mapsto\widetilde{P}_{N,\Omega}(z)\varphi:=\wti{f}_{N,\Omega}(z,\varphi), 
\end{cases} 
\end{equation}
where $\wti{f}_{N,\Omega}(z,\varphi)$ is the unique solution of \eqref{eqn:bvpNeumann}, 
then the solution operator $\widetilde{P}_{N,\Omega}(z)$ is an extension of 
$P_{1,N,\Omega}(z)$ in \eqref{eqn:P_N}, and $\widetilde{P}_{N,\Omega}(z)$ is continuous, 
when the adjoint space $\mathscr{G}_D(\partial\Omega)^*$ and $\dom(A_{max,\Omega})$ 
are endowed with the norms in item $(ii)$. \\[1mm] 
\noindent $(ix)$ For all $z\in\rho(A_{D,\Omega})$ the Dirichlet-to-Neumann map $M_{\Omega}(z)$ in  
\eqref{5.34} admits an extension 
\begin{equation}\label{iutee}
\widetilde M_{\Omega}(z): 
\begin{cases} 
\mathscr{G}_N(\partial\Omega)^*\rightarrow\mathscr{G}_D(\partial\Omega)^*,  
\\[2pt]
\varphi\mapsto\widetilde{M}_{\Omega}(z)\varphi
:=-\widetilde{\gamma}_N\widetilde{P}_{D,\Omega}(z)\varphi,
\end{cases} 
\end{equation}
and $\widetilde{M}_{\Omega}(z)$ is continuous, when the adjoint spaces 
$\mathscr{G}_D(\partial\Omega)^*$, $\mathscr{G}_N(\partial\Omega)^*$ carry the 
natural topology induced by \eqref{scalarjussi} on $\mathscr{G}_D(\partial\Omega)$, 
$\mathscr{G}_N(\partial\Omega)$, respectively.
\end{theorem}
%%%%%%%%%%

As a preamble to the proof of this theorem, we first deal with a couple of useful elementary results.

%%%%%%%%
\begin{lemma}\label{uyrre}
Let $X,Y$ be two Banach spaces and assume that $T\in{\mathcal{B}}(X,Y)$ is surjective. Then 
\begin{equation}\label{eq:Rfav}
\|y\|_Y\approx\inf_{x\in X,\,\,Tx=y}\|x\|_X\,\text{ uniformly in }\,y\in Y,
\end{equation}
that is, there exists a constant $C\in(1,\infty)$, independent of $y\in Y$, such that
\begin{equation}\label{explain}
C^{-1 }\|y\|_Y\leq\inf_{x\in X,\,\,Tx=y}\|x\|_X\leq C\|y\|_Y.
\end{equation}
Moreover, if the space $X$ is reflexive then $Y$ is also reflexive.
\end{lemma}
%%%%%%%%
\begin{proof}
The fact that $T:X\rightarrow Y$ is linear and continuous implies that $\ker T$ is a 
closed subspace of $X$. Moreover, given that $T$ is surjective, $T$ induces a  
continuous isomorphism 
\begin{equation}\label{eq:Taagb}
\widehat{T}:X/\ker T\rightarrow Y,\quad\widehat{T}(x+\ker T):=Tx,\quad\forall\,x\in X,
\end{equation}
where the space on the left-hand side of \eqref{eq:Taagb} is equipped with the quotient norm 
\begin{equation}\label{eq:eBB}
\|x+\ker T\|_{X/\ker T}:=\inf_{z\in\ker T}\|x+z\|_X,\quad\forall\,x\in X.
\end{equation}
Then \eqref{eq:Rfav} becomes a consequence of \eqref{eq:Taagb}--\eqref{eq:eBB}
and the Open Mapping Theorem. Next, we recall that in general, 
\begin{equation}\label{eq:Nvav}
\parbox{9.80cm}
{every closed subspace of a reflexive Banach space is reflexive, 
every quotient of a reflexive Banach space by a closed subspace 
is reflexive, and every Banach space continuously isomorphic with a 
reflexive Banach space is itself reflexive.}
\end{equation}
Granted these facts and assuming that $X$ is a reflexive Banach space, it follows from 
\eqref{eq:Taagb} that $Y$ is also reflexive. 
\end{proof}
%%%%%%%%

%%%%%%%%
\begin{lemma}\label{uyrre.LD}
Let $X,Y$ be two Banach spaces with the property that $X\subset Y$ densely, and the inclusion 
$\iota:X\hookrightarrow Y$ is continuous. Then the following hold. \\[1mm] 
$(i)$ The operator $\iota^\ast:Y^\ast\rightarrow X^\ast$ is linear, continuous, and injective. 
In particular, identifying $Y^\ast$ with $\ran(\iota^\ast)$ yields the continuous embedding 
$Y^\ast\hookrightarrow X^\ast$. \\[1mm] 
$(ii)$ In the special case when $Y$ is a Hilbert space,  one has
\begin{equation}\label{eq:tre3e}
X\overset{\iota}{\hookrightarrow} Y\equiv Y^\ast\overset{\,\,\iota^\ast}{\hookrightarrow} X^\ast,
\end{equation}
where $Y\equiv Y^\ast$ is the canonical identification between the Hilbert space $Y$ and its dual,
and the duality pairing between $X$ and $X^\ast$ is compatible with the inner product in $Y$. 
\\[1mm] 
$(iii)$ If $X$ is reflexive, then the embedding $Y^\ast\hookrightarrow X^\ast$ 
has dense range. 
\end{lemma}
%%%%%%%%%
\begin{proof}
The main claim in part $(i)$ is a particular case of the well-known general result 
to the effect that if $X,Y$ are Banach spaces and $T\in{\mathcal{B}}(X,Y)$ then $\ran(T)$ 
is dense in $Y$ if and only if $T^\ast$ is injective. Regarding $(ii)$, assume that
$Y$ is a Hilbert space with inner product $(\cdot,\cdot)_Y$. Then the identification 
$Y\equiv Y^\ast$ manifests itself in the following manner: $Y\ni y\mapsto\Lambda_y:=(y,\cdot)_Y\in Y^\ast$.
Consequently, if $x\in X$ and $y\in Y$, then $\iota(x)\in Y$ and  
\begin{equation}\label{eq:jugf97h}
{}_{X^\ast}\langle\iota^\ast(\Lambda_y),x\rangle_X
={}_{Y^\ast}\langle\Lambda_y,\iota(x)\rangle_Y=(y,\iota(x))_Y.
\end{equation}
This proves that the duality pairing between $X$ and $X^\ast$ is compatible with the inner product in $Y$. 
Finally, to deal with the claim in item $(iii)$,
assume that a functional in $(X^\ast)^\ast=X$ has been fixed with the property that
its restriction to $Y^\ast$ (identified with $\ran(\iota^\ast)$, as a subspace of $X^\ast$) 
vanishes identically. This comes down to having some $x\in X$ such that $\Lambda(\iota(x))=0$ for 
each $\Lambda\in Y^\ast$, and the density of the embedding $Y^\ast\hookrightarrow X^\ast$
follows as soon as one shows that $x=0$. The latter conclusion is, however, implied by the Hahn--Banach 
theorem. 
\end{proof}
%%%%%%%%%

We are now ready to present the proof of Theorem~\ref{t5.5}. 

\smallskip

%%%%%%%%%
\begin{proof}[Proof of Theorem~\ref{t5.5}]
Regarding $(i)$, one verifies that $\mathscr{G}_N(\partial\Omega)=\dom\big(\Lambda^{-1/2}\big)$. 
The assertion on $\mathscr{G}_D(\partial\Omega)$ in \eqref{eqn:domLambda} follows in a similar form by interchanging the 
roles of $\Lambda$ with $\Sigma$, $P_{0,D,\Omega}$ with $P_{1,N,\Omega}$ and 
$\mathscr{G}_N(\partial\Omega)$ with $\mathscr{G}_D(\partial\Omega)$ (see \cite[Section~2]{BM13}).   

According to Theorem~\ref{t5.4}, the operator 
\begin{equation}\label{y566}
\Lambda=P_{0,D,\Omega}(i)^* P_{0,D,\Omega}(i)\in\cB\big(L^2(\partial\Omega)\big)
\end{equation}  
is self-adjoint, injective, and non-negative. Hence $\ran\big(\Lambda\big)$ and 
$\ran\big(\Lambda^{1/2}\big)$ are both dense in $L^2(\partial\Omega)$. The space 
\begin{equation}\label{eq:Nbfd66}
\mathscr{G}:=\ran\big(\Lambda^{1/2}\big)=\dom\big(\Lambda^{-1/2}\big) 
\end{equation}
is now equipped with the inner product
\begin{align}\label{tr54edd}
\begin{split}
& (\varphi,\psi)_{\mathscr{G}}:=\big(\Lambda^{-1/2}\varphi,\Lambda^{-1/2}\psi\big)_{L^2(\partial\Omega)}, 
\\[2pt]
& \forall\,\varphi,\psi\in\mathscr{G}=\ran\big(\Lambda^{1/2}\big)=\dom\big(\Lambda^{-1/2}\big).
\end{split}
\end{align}
Then $\mathscr{G}$ is a Hilbert space which is densely embedded in $L^2(\partial\Omega)$ 
and hence gives rise to a Gelfand triple 
$\mathscr{G}\hookrightarrow L^2(\partial\Omega)\hookrightarrow\mathscr{G}^*$,
where the adjoint (antidual) space $\mathscr{G}^*$ coincides with the completion of 
$L^2(\partial\Omega)$ equipped with the inner product 
\begin{equation}\label{freitagder13.}
\big(\Lambda^{1/2}u,\Lambda^{1/2} v\big)_{L^2(\partial\Omega)},\quad\forall\,u,v\in L^2(\partial\Omega).
\end{equation}

For $\varphi\in L^2(\partial\Omega)$ one computes
\begin{align}\label{We-EE.10} 
\big\|P_{0,D,\Omega}(i)\varphi\big\|_{L^2(\Omega)}^2
&=\big(P_{0,D,\Omega}(i)\varphi, P_{0,D,\Omega}(i)\varphi\big)_{L^2(\Omega)}   
\nonumber\\[2pt] 
&=\big(P_{0,D,\Omega}(i)^* P_{0,D,\Omega}(i)\varphi,\varphi\big)_{L^2(\partial\Omega)} 
\nonumber\\[2pt] 
&=(\Lambda\varphi,\varphi)_{L^2(\partial\Omega)} 
=\big(\Lambda^{1/2}\varphi,\Lambda^{1/2}\varphi\big)_{L^2(\partial\Omega)} 
\nonumber\\[2pt] 
&=\big\|\Lambda^{1/2}\varphi\big\|^2_{L^2(\partial\Omega)}=\norm{\varphi}^2_{\mathscr{G}^*}.
\end{align}
As the range of $P_{0,D,\Omega}(i)$ is dense in the space $\ker(A_{max,\Omega}-iI)$ with respect to 
the $L^2(\Omega)$-norm (see Theorem~\ref{t5.2}\,$(i)$), it follows from \eqref{We-EE.10} that
\begin{align}\label{eMNa}
\begin{split}
& \text{$P_{0,D,\Omega}(i)$ admits a continuation to an isometry}  
\\[2pt] 
& \quad\text{$\widetilde{P}_{D,\Omega}(i)$ acting from $\mathscr{G}^*$ onto $\ker(A_{max,\Omega}-iI)$,}
\end{split}
\end{align}
where the latter space is equipped with the $L^2(\Omega)$-norm. Furthermore, as 
$P_{1,D,\Omega}(i)$ is a restriction of $P_{0,D,\Omega}(i)$ and 
$\dom(P_{1,D,\Omega}(i))=H^1(\partial\Omega)$ is dense in $\mathscr{G}^*$ 
(a consequence of $H^1(\partial\Omega)$ being dense in $L^2(\partial\Omega)$
and the definition of the norm in $\mathscr{G}^*$), it follows that $P_{1,D,\Omega}(i)$ also 
admits a continuation to an isometry from $\mathscr{G}^*$ onto $\ker(A_{max,\Omega}-iI)$
which coincides with $\widetilde P_{D,\Omega}(i)$. Furthermore, for 
$\varphi\in L^2(\partial\Omega)\subset\mathscr{G}^*$ and $f\in L^2(\Omega)$ one concludes from
\begin{equation}\label{yav-WACO}
\begin{split}
\big(P_{0,D,\Omega}(i)^*f,\varphi\big)_{L^2(\partial\Omega)}
&=\big(f,P_{0,D,\Omega}(i)\varphi\big)_{L^2(\Omega)}
=\big(f,\widetilde P_{D,\Omega}(i)\varphi\big)_{L^2(\Omega)}
\\[2pt]
&={}_{\mathscr{G}}\big\langle\widetilde P_{D,\Omega}(i)^*f,\varphi\big\rangle_{\mathscr{G}^*}
=\big(\widetilde P_{D,\Omega}(i)^*f,\varphi\big)_{L^2(\Omega)}
\end{split}
\end{equation}
(here \eqref{WACO.WTF.1} and the subsequent discussion is relevant) 
that the adjoint of the operator $\widetilde P_{D,\Omega}(i):\mathscr{G}^*\rightarrow L^2(\Omega)$ 
coincides with $P_{0,D,\Omega}(i)^*$. Together with \eqref{eMNa} this shows that $P_{0,D,\Omega}(i)^*$ 
is a continuous map from $L^2(\Omega)$ onto $\mathscr{G}$.

In a similar way as in \eqref{We-EE.10} the fact that for each $\varphi\in L^2(\partial\Omega)$ one has 
\begin{equation}\label{iute66j}
\|\Lambda\varphi\|^2_{\mathscr{G}}=(\Lambda\varphi,\Lambda\varphi)_{\mathscr{G}}
=\big(\Lambda^{1/2}\varphi,\Lambda^{1/2}\varphi\big)_{L^2(\partial\Omega)} 
=\|\varphi\|^2_{\mathscr{G}^*}
\end{equation}
shows that the operators $\Lambda=\overline{\Im(M_{\Omega}(i))}$ and $\Im(M_{\Omega}(i))$ 
admit continuations to an isometry $\widetilde\Lambda$ from $\mathscr{G}^*$ onto $\mathscr{G}$ 
(one observes that $\ran(\Lambda)$ is a dense subspace in the Hilbert space 
$(\mathscr{G},(\cdot,\cdot)_{\mathscr{G}})$) with 
\begin{equation}\label{5.65}
\Im(M_{\Omega}(i))\subset\Lambda\subset\widetilde\Lambda
=P_{0,D,\Omega}(i)^*\widetilde P_{D,\Omega}(i);
\end{equation}  
in the last equality in \eqref{5.65} we have used \eqref{y566} and the fact 
that both operators $\widetilde P_{D,\Omega}(i):\mathscr{G}^*\rightarrow L^2(\Omega)$ 
and $P_{0,D,\Omega}(i)^*:L^2(\Omega)\rightarrow\mathscr{G}$ are continuous.

From \eqref{5.65} and the fact that $P_{0,D,\Omega}(i)^*|_{\ker(A_{max,\Omega}-iI)}$ 
is a bijection onto $\mathscr{G}_N(\partial\Omega)$ (as seen from \eqref{eqn:odec} and 
\eqref{eqn:ranP*}) one deduces that
\begin{equation}\label{tafff965}
\mathscr{G}=\ran\big(\widetilde\Lambda\big)=\ran\big(P_{0,D,\Omega}(i)^* 
\widetilde{P}_{D,\Omega}(i)\big)=\ran(P_{0,D,\Omega}(i)^*)=\mathscr{G}_N(\partial\Omega).
\end{equation}
This completes the treatment of $(i)$.

Next, we proceed to verify the claims made in relation to $\gamma_D$ in item $(ii)$. First, we define 
\begin{equation}\label{We-EE.12}
\widetilde{\gamma}_D:\dom(A_{max,\Omega})\rightarrow\mathscr{G}_N(\partial\Omega)^*
\end{equation}
as follows: Given any $f\in\dom(A_{max,\Omega})=\dom(A_{D,\Omega})
\dotplus\ker(A_{max,\Omega}-i I)$, write $f=f_D+f_i$ with 
\begin{align}\label{eq:trVva.1}
\begin{split}
f_D &:=(A_{D,\Omega}-iI)^{-1}(A_{max,\Omega}-iI)f\in\dom(A_{D,\Omega}), 
\\[2pt]
f_i &:=f-(A_{D,\Omega}-iI)^{-1}(A_{max,\Omega}-iI)f\in\ker(A_{max,\Omega}-iI),
\end{split}
\end{align}
then set
\begin{equation}\label{We-EE.12pp}
\widetilde{\gamma}_D f:=\widetilde P_{D,\Omega}(i)^{-1}f_i\in\mathscr{G}_N(\partial\Omega)^*, 
\end{equation}
where the membership in \eqref{We-EE.12pp} follows from \eqref{eMNa} and \eqref{tafff965}.
Upon noting that
\begin{align}\label{yarfJvv}
\|f_D\|_{L^2(\Omega)}&=\big\|(A_{D,\Omega}-iI)^{-1}(A_{max,\Omega}-iI)f\big\|_{L^2(\Omega)}
\nonumber\\[2pt]
&\leq C\big\|(A_{max,\Omega}-iI)f\big\|_{L^2(\Omega)}
\nonumber\\[2pt]
&\leq C\big\{\|f\|_{L^2(\Omega)}+\|A_{max,\Omega}f\|_{L^2(\Omega)}\big\},
\end{align}
for some constant $C\in(0,\infty)$, independent of $f$, one estimates  
\begin{align}\label{yarfrRR}
\big\|\widetilde\gamma_D f\big\|_{\mathscr{G}_N(\partial\Omega)^*} 
&=\big\|\widetilde{P}_{D,\Omega}(i)^{-1}f_i\big\|_{\mathscr{G}_N(\partial\Omega)^*}
=\big\|\widetilde{P}_{D,\Omega}(i)^{-1}(f-f_D)\big\|_{\mathscr{G}_N(\partial\Omega)^*}    
\nonumber\\[2pt]
&\leq\big\|\widetilde{P}_{D,\Omega}(i)^{-1}\big\|_{\cB(L^2(\Omega),\mathscr{G}_N(\partial\Omega)^*)}
\Big\{\|f\|_{L^2(\Omega)}+\|f_D\|_{L^2(\Omega)}\Big\}   
\nonumber\\[2pt]
&\leq C\big\{\|f\|_{L^2(\Omega)}+\|A_{max,\Omega}f\|_{L^2(\Omega)}\big\},
\end{align}
proving that the operator $\widetilde{\gamma}_D$ in \eqref{We-EE.12} 
is continuous with respect to the graph norm of $A_{max,\Omega}$ in $L^2(\Omega)$ 
and the norm on $\mathscr{G}_N(\partial\Omega)^*$ induced by \eqref{freitagder13.}.
To see that $\widetilde{\gamma}_D$ is compatible with ${\gamma}_D$, 
consider the case when $f\in\dom(A_{max,\Omega})\cap H^{1/2}(\Omega)$ which forces 
$f_i\in\ker(A_{max,\Omega}-i I)\cap H^{1/2}(\Omega)$ (cf.~\eqref{eqn:decAmax}).
In particular, $\gamma_D f_i\in L^2(\partial\Omega)$ by \eqref{eqn:gammaDs.2} with $s=1/2$.
In this scenario, 
\begin{equation}\label{5.63}
\begin{split}
\widetilde{\gamma}_D f &=\widetilde P_{D,\Omega}(i)^{-1}f_i 
=\widetilde{P}_{D,\Omega}(i)^{-1}P_{0,D,\Omega}(i)\gamma_D f_i   
\\[2pt] 
&=\widetilde{P}_{D,\Omega}(i)^{-1}\widetilde{P}_{D,\Omega}(i)\gamma_D f_i 
=\gamma_D f_i=\gamma_D f.
\end{split}
\end{equation}
The first equality in \eqref{5.63} follows from \eqref{We-EE.12pp}. The second 
equality in \eqref{5.63} employs the fact that $f_i=P_{0,D,\Omega}(i)\gamma_D f_i$, 
which in turn is a consequence of the fact that both $f_i$ and 
$P_{0,D,\Omega}(i)\gamma_D f_i$ solve the boundary value problem
\begin{equation}\label{tree-75}
\begin{cases}
(-\Delta+V-i)f=0\,\text{ in $\Omega,\quad f\in H^{1/2}(\Omega)\cap\dom(A_{max,\Omega})$,}     
\\[2pt]
\gamma_D f=\gamma_D f_i\in L^2(\partial\Omega), 
\end{cases}
\end{equation} 
which is well posed, by Lemma~\ref{LamAA.2} with $s=0$ and $z=i$. The third equality in \eqref{5.63} is clear 
from the fact that $\widetilde{P}_{D,\Omega}(i)$ is an extension of $P_{0,D,\Omega}(i)$ (cf.~\eqref{eMNa}). 
Hence, $\widetilde{\gamma}_D$ is an extension of $\gamma_D$, implying the first assertion in item $(ii)$. 
Next, the claim that $\ker(\widetilde{\gamma}_D)=\dom(A_{D,\Omega})$ is an immediate consequence of the 
definition of $\widetilde{\gamma}_D$ in \eqref{We-EE.12pp} since $\widetilde{P}_{D,\Omega}(i)^{-1}$ acts 
isometrically from $\ker(A_{max,\Omega}-iI)$ onto $\mathscr{G}_N(\partial\Omega)^*$.

Concerning \eqref{4.55.REE.D}, in a first stage we shall prove that there exists 
a constant $C\in(0,\infty)$ with the property that
\begin{align}\label{4.55.REE.Dbbb}
\begin{split}
& f\in\dom(A_{max,\Omega})\,\text{ and }\,\widetilde{\gamma}_D f=0\,\text{ imply}
\\[2pt]
& \quad f\in H^{3/2}(\Omega)\,\text{ and }\,\|f\|_{H^{3/2}(\Omega)}\leq C\|\Delta f\|_{L^2(\Omega)}.
\end{split}
\end{align}
To this end, assume that $f\in\dom(A_{max,\Omega})$ satisfies $\widetilde{\gamma}_D f=0$ 
in $\mathscr{G}_N(\partial\Omega)^*$. Then \eqref{We-EE.12pp} forces $f_i=0$, hence $f=f_D$. 
Introduce $g:=-\Delta f\in L^2(\Omega)$. Since $f_D$ from \eqref{eq:trVva.1} belongs to 
$H^{3/2}(\Omega)$ (cf.\  Theorem~\ref{t4.4}), it follows that $f_D$ solves the boundary value 
problem
\begin{equation}\label{trKM75CX}
\begin{cases}
-\Delta u=g\,\text{ in }\,\Omega,\quad u\in H^{3/2}(\Omega),     
\\[2pt]
\gamma_D u=0\,\text{ on }\,\partial\Omega,
\end{cases}
\end{equation} 
and satisfies the naturally accompanying estimate 
$\|f_D\|_{H^{3/2}(\Omega)}\leq C\|\Delta f_D\|_{L^2(\Omega)}$ (cf. \cite{JK95}).
In turn, this implies (since $f=f_D$)
\begin{equation}\label{eq:Hcac.3}
\|f\|_{H^{3/2}(\Omega)}=\|f_D\|_{H^{3/2}(\Omega)}\leq C\|\Delta f\|_{L^2(\Omega)},
\end{equation}
and \eqref{4.55.REE.Dbbb} follows. Having established \eqref{4.55.REE.Dbbb}, we now prove 
\eqref{4.55.REE.D} by reasoning as follows. Given $s\in[0,1]$ and any $f\in\dom(A_{max,\Omega})$ 
with $\varphi:=\widetilde{\gamma}_D f\in H^s(\partial\Omega)$, use the surjectivity of the map 
\eqref{eqn:gammaDs.2} in order to find $g\in H^{s+(1/2)}(\Omega)\cap\dom(A_{max,\Omega})$
with $\gamma_D g=\varphi$. Moreover, by the Open Mapping Theorem and the surjectivity
of \eqref{eqn:gammaDs.2}, matters may be arranged so that 
\begin{equation}\label{eq:vgtg6t}
\|g\|_{H^{s+(1/2)}(\Omega)}+\|\Delta g\|_{L^2(\Omega)}\leq C\|\varphi\|_{H^s(\partial\Omega)}
=C\|\widetilde{\gamma}_D f\|_{H^s(\partial\Omega)}
\end{equation}
for some constant $C\in(0,\infty)$ independent of $\varphi$. Then $h:=(f-g)\in\dom(A_{max,\Omega})$
has $\widetilde{\gamma}_D h=0$, so \eqref{4.55.REE.Dbbb} implies $h\in H^{3/2}(\Omega)$ 
and $\|h\|_{H^{3/2}(\Omega)}\leq C\|\Delta h\|_{L^2(\Omega)}$. Consequently, 
$f=g+h\in H^{s+(1/2)}(\Omega)$ and 
\begin{align}\label{ut4esws}
\|f\|_{H^{s+(1/2)}(\Omega)}& \leq\|g\|_{H^{s+(1/2)}(\Omega)}+\|h\|_{H^{s+(1/2)}(\Omega)}
\nonumber\\[2pt]
& \leq\|g\|_{H^{s+(1/2)}(\Omega)}+\|h\|_{H^{3/2}(\Omega)}
\nonumber\\[2pt] 
& \leq C\|\widetilde{\gamma}_D f\|_{H^s(\partial\Omega)}+C\|\Delta h\|_{L^2(\Omega)}
\nonumber\\[2pt] 
& \leq C\|\widetilde{\gamma}_D f\|_{H^s(\partial\Omega)}
+C\big(\|\Delta f\|_{L^2(\Omega)}+\|\Delta g\|_{L^2(\Omega)}\big) 
\nonumber\\[2pt] 
& \leq C\|\widetilde{\gamma}_D f\|_{H^s(\partial\Omega)}+C\|\Delta f\|_{L^2(\Omega)},
\end{align}
finishing the proof of \eqref{4.55.REE.D}.

Next, we verify the assertions for $\gamma_N$ in item $(ii)$. 
Denote by $\widetilde P_{N,\Omega}(i)$ the extension of $P_{1,N,\Omega}(i)$ 
to an isometry from $\mathscr{G}_D(\partial\Omega)^*$ onto $\ker(A_{max,\Omega}-iI)$ 
(which is constructed in a similar way as $\widetilde P_{D,\Omega}(i)$ above) and define
\begin{equation}\label{We-EE.12n}
\widetilde{\gamma}_N:\dom(A_{max,\Omega})\rightarrow\mathscr{G}_D(\partial\Omega)^*
\end{equation}
as follows: Given any $f\in\dom(A_{max,\Omega})=\dom(A_{N,\Omega})\dotplus\ker(A_{max,\Omega}-i I)$, 
write $f=f_N+f_i$ with $f_N\in\dom(A_{N,\Omega})$ and $f_i\in\ker(A_{max,\Omega}-iI)$, then set
\begin{equation}\label{We-EE.12ppn}
\widetilde{\gamma}_N f:=\widetilde P_{N,\Omega}(i)^{-1}f_i\in\mathscr{G}_D(\partial\Omega)^*.
\end{equation}
The same arguments as in \eqref{yarfrRR} (with $\widetilde\gamma_D$, 
$\mathscr{G}_N(\partial\Omega)^*$, $\widetilde P_{D,\Omega}(i)$ and $A_{D,\Omega}$ 
replaced by $\widetilde\gamma_N$, $\mathscr{G}_D(\partial\Omega)^*$, 
$\widetilde P_{N,\Omega}(i)$ and $A_{N,\Omega}$, respectively) show that $\widetilde\gamma_N$ 
in \eqref{We-EE.12n} is continuous with respect to the natural graph norm in $\dom(A_{max,\Omega})$ 
and the norm on $\mathscr{G}_D(\partial\Omega)^*$.

To see that $\widetilde{\gamma}_N$ is compatible with ${\gamma}_N$ in \eqref{eqn:gammaN-pp}, 
we first consider the case when $f\in\dom(A_{max,\Omega})\cap H^{3/2}(\Omega)$ which forces 
$f_i\in\ker(A_{max,\Omega}-i I)\cap H^{3/2}(\Omega)$ (cf.~\eqref{eqn:decAmax2}).
In particular, $\gamma_N f_i\in L^2(\partial\Omega)$ by \eqref{eqn:gammaN-pp}. One can then write  
\begin{equation}\label{5.63n}
\begin{split}
\widetilde{\gamma}_N f &=\widetilde P_{N,\Omega}(i)^{-1}f_i 
=\widetilde{P}_{N,\Omega}(i)^{-1}P_{1,N,\Omega}(i)\gamma_N f_i   
\\[2pt] 
&=\widetilde{P}_{N,\Omega}(i)^{-1}\widetilde{P}_{N,\Omega}(i)\gamma_N f_i 
=\gamma_N f_i=\gamma_N f.
\end{split}
\end{equation}
In \eqref{5.63n}, the first equality follows from \eqref{We-EE.12ppn}, while the second 
equality employs the fact that $f_i=P_{1,N,\Omega}(i)\gamma_N f_i$, 
which in turn is a consequence of the fact that both $f_i$ and 
$P_{1,N,\Omega}(i)\gamma_N f_i$ solve the boundary value problem
\begin{equation}
\begin{cases}
(-\Delta+V-i)f=0\,\text{ in $\Omega,\quad f\in H^{3/2}(\Omega)\cap\dom(A_{max,\Omega})$,}     
\\[2pt]
-\gamma_N f=-\gamma_N f_i\in L^2(\partial\Omega), 
\end{cases}
\end{equation} 
which is well posed, by Lemma~\ref{LamAA.2N} with $s=1$ and $z=i$.
Finally, the third equality in \eqref{5.63n} is clear from the fact that 
$\widetilde{P}_{N,\Omega}(i)$ is an extension of $P_{N,D,\Omega}(i)$. 

Having established \eqref{5.63n}, we conclude that $\widetilde\gamma_N$ 
in \eqref{We-EE.12n} is an extension of the Neumann trace operator
$\gamma_N:H^{3/2}(\Omega)\cap\dom(A_{max,\Omega})\rightarrow L^2(\partial\Omega)$,
that is, 
\begin{equation}\label{eq:RTTYH}
\text{$\widetilde\gamma_N$ is compatible with $\gamma_N$ in
\eqref{eqn:gammaN-pp} when $s=\tfrac{3}{2}$}. 
\end{equation}

It turns out that the compatibility property established in \eqref{eq:RTTYH} suffices 
to prove $(v)$, a task to which we now turn. Specifically, fix two arbitrary functions 
$f\in\dom(A_{max,\Omega})$ and $\phi\in\mathscr{G}_D(\partial\Omega)$. Then 
$\phi\in H^1(\partial\Omega)$ and \eqref{We-Q.17} ensures the existence of a function 
\begin{equation}\label{eq:gDDc-NnN}
\text{$g\in H^{3/2}(\Omega)\cap\dom(A_{max,\Om})$ such that $\gamma_N g=0$ 
and $\gamma_D g=\phi$.} 
\end{equation}
Making use of Lemma~\ref{l5.1MM}, it is possible to find 
$\{f_j\}_{j\in{\mathbb{N}}}\subset\dom(A_{max,\Omega})\cap H^{3/2}(\Omega)$ 
with the property that
\begin{equation}\label{Mi-tr-MM.4-NnN}
f_j\underset{j\to\infty}{\longrightarrow}f\,\text{ in $L^2(\Omega)$} 
\,\text{ and }\, 
\Delta f_j\underset{j\to\infty}{\longrightarrow}\Delta f 
\,\text{ in $L^2(\Omega)$}.
\end{equation}  
Then $f_j\to f$ in the natural graph norm of $\dom(A_{max,\Omega})$ as $j\to\infty$, and one concludes that 
$\widetilde{\gamma}_N f_j\to\widetilde{\gamma}_N f$ in 
$\mathscr{G}_D(\partial\Omega)^*$ as $j\to\infty$ by the continuity of the second map in 
\eqref{4.55}. Furthermore, $\widetilde{\gamma}_N f_j=\gamma_N f_j\in L^2(\partial\Omega)$ 
for each $j\in{\mathbb{N}}$ by \eqref{eq:RTTYH} and the fact that $f_j\in H^{3/2}(\Omega)$. 
With the help of these remarks, \eqref{eq:gDDc-NnN}, \eqref{Mi-tr-MM.4-NnN}, and Green's 
formula \eqref{GGGRRR}, one computes  
\begin{align}\label{eq:33Nnab-NnN}
{}_{\mathscr{G}_D(\partial\Omega)^\ast}\big\langle\widetilde{\gamma}_N f,
\phi\big\rangle_{\mathscr{G}_D(\partial\Omega)}
&=\lim_{j\to\infty}{}_{\mathscr{G}_D(\partial\Omega)^\ast}\big\langle\widetilde{\gamma}_N f_j,
\phi\big\rangle_{\mathscr{G}_D(\partial\Omega)}
\nonumber\\[2pt] 
&=\lim_{j\to\infty}\big(\gamma_N f_j,\gamma_D g\big)_{L^2(\partial\Omega)}
\nonumber\\[2pt] 
&=\lim_{j\to\infty}\big\{-(f_j,\Delta g)_{L^2(\Omega)}+(\Delta f_j,g)_{L^2(\Omega)}\big\}
\nonumber\\[2pt] 
&=-(f,\Delta g)_{L^2(\Omega)}+(\Delta f,g)_{L^2(\Omega)},
\end{align}
and \eqref{eq:33f4iU.NN} follows. 

Next, we shall employ \eqref{eq:33f4iU.NN} in order to show that $\widetilde\gamma_N$ 
is also compatible with $\gamma_N$ in \eqref{eqn:gammaN-pp} when $s\in[\tfrac{1}{2},\tfrac{3}{2})$.
In this regard, it suffices to treat the case $s=1/2$. With this goal in mind, fix 
$f\in\dom(A_{max,\Omega})\cap H^{1/2}(\Omega)$ and let $\psi\in\mathscr{G}_D(\partial\Omega)$ be arbitrary. 
Then by \eqref{eq:33f4iU.NN} one has 
\begin{equation}\label{okgut}
{}_{\mathscr{G}_D(\partial\Omega)^\ast}\big\langle\widetilde{\gamma}_N f,
\psi\big\rangle_{\mathscr{G}_D(\partial\Omega)}
=-(f,\Delta g)_{L^2(\Omega)}+(\Delta f,g)_{L^2(\Omega)}
\end{equation}
for any $g\in H^{3/2}(\Omega)\cap\dom(A_{max,\Om})$ such that $\gamma_N g=0$ 
and $\gamma_D g=\psi$. On the other hand by Green's identity \eqref{GGGRRR} one also has  
\begin{equation}\label{okgut2}
{}_{(H^{1}(\partial\Omega))^\ast}\big\langle\gamma_N f,
\psi\big\rangle_{H^1(\partial\Omega)}
=-(f,\Delta g)_{L^2(\Omega)}+(\Delta f,g)_{L^2(\Omega)},
\end{equation}
and it follows from \eqref{okgut}--\eqref{okgut2} that the functionals $\widetilde\gamma_N f$ 
and $\gamma_N f$ coincide on $\mathscr{G}_D(\partial\Omega)\subset H^1(\partial\Omega)$
whenever $f\in\dom(A_{max,\Omega})\cap H^{1/2}(\Omega)$. This finishes the proof of the claim 
that $\widetilde\gamma_N$ is compatible with $\gamma_N$ in \eqref{eqn:gammaN-pp}. 

Regarding the second formula in \eqref{4.30}, the statement 
$\ker(\widetilde{\gamma}_N)=\dom(A_{N,\Omega})$ is an immediate consequence 
of the definition of $\widetilde{\gamma}_N$ in \eqref{We-EE.12ppn} since 
$\widetilde{P}_{N,\Omega}(i)^{-1}$ acts isometrically from $\ker(A_{max,\Omega}-iI)$ onto 
$\mathscr{G}_D(\partial\Omega)^*$. Finally, \eqref{4.55.REE.N} is proved in a similar manner 
to \eqref{4.55.REE.D}, where instead of \eqref{trKM75CX} one has to make use of the well-posedness of 
the boundary value problem
\begin{equation}\label{trKM75CX-2}
\begin{cases}
(-\Delta+1)u=g\in L^2(\Omega),\quad u\in H^{3/2}(\Omega),     
\\[2pt]
\gamma_N u=0\,\text{ on }\,\partial\Omega,
\end{cases}
\end{equation} 
and the naturally accompanying estimate $\|u\|_{H^{3/2}(\Omega)}\leq C\|g\|_{L^{2}(\Omega)}$; see \cite{FMM98}.

At this point we note that the surjectivity of the maps in \eqref{4.55} 
can be used to show that
\begin{equation}\label{4.55EEdc}
\text{the Banach spaces $\mathscr{G}_N(\partial\Omega)$, $\mathscr{G}_D(\partial\Omega)$ 
are reflexive}  
\end{equation}
(which also follows directly from  part $(i)$). Specifically, one first observes that when 
$\dom(A_{max,\Omega})$ is equipped with the natural graph norm, the mapping 
\begin{equation}\label{mapping}
\dom(A_{max,\Omega})\ni f\mapsto(f,\Delta f)\in L^2(\Omega)\oplus L^2(\Omega)
\end{equation} 
is a continuous isomorphism onto its range and this yields (cf.\  the discussion in 
\eqref{eq:Nvav}) that $\dom(A_{max,\Omega})$ is a reflexive Banach space. With this 
in hand, \eqref{4.55EEdc} follows from the surjectivity of the maps in \eqref{4.55}, 
Lemma~\ref{uyrre}, and the well-known fact that 
\begin{equation}\label{eq:kig}
\text{a Banach space is reflexive if and only if its dual is reflexive}.
\end{equation}

Turning to $(iii)$, identity \eqref{ut444} is a direct consequence of \eqref{3.30} 
and \eqref{4.30}. Regarding the first claim in part $(iv)$, we start by fixing some 
arbitrary functions $f\in\dom(A_{max,\Omega})$ and $\phi\in\mathscr{G}_N(\partial\Omega)$. 
Then $\phi\in L^2(\partial\Omega)$ and \eqref{We-Q.17} ensures the existence of a function 
\begin{equation}\label{eq:gDDc}
\text{$g\in H^{3/2}(\Omega)\cap\dom(A_{max,\Om})$ such that $\gamma_D g=0$ and $\gamma_N g=\phi$.} 
\end{equation}
Making use of Lemma~\ref{l5.1MM}, it is possible to find 
$\{f_j\}_{j\in{\mathbb{N}}}\subset\dom(A_{max,\Omega})\cap H^{3/2}(\Omega)$ 
with the property that
\begin{equation}\label{Mi-tr-MM.4}
f_j\underset{j\to\infty}{\longrightarrow} f\,\text{ in $L^2(\Omega)$}\,\text{ and }\, 
\Delta f_j\underset{j\to\infty}{\longrightarrow}\Delta f\,\text{ in $L^2(\Omega)$}.
\end{equation}  
Then $f_j\to f$ in the natural graph norm of $\dom(A_{max,\Omega})$ as $j\to\infty$, from which one 
deduces that $\widetilde{\gamma}_D f_j\to\widetilde{\gamma}_D f$ in $\mathscr{G}_N(\partial\Omega)^*$ 
as $j\to\infty$ due to the continuity of the first map in \eqref{4.55}. Moreover, for each 
$j\in{\mathbb{N}}$ one has $\widetilde{\gamma}_D f_j=\gamma_D f_j\in L^2(\partial\Omega)$ 
since $f_j\in H^{3/2}(\Omega)$ and $\widetilde{\gamma}_D$ is compatible with $\gamma_D$.
In turn, these observations and Green's formula \eqref{GGGRRR} permit us to write 
(keeping in mind that $\gamma_D g=0$) 
\begin{align}\label{eq:33Nnab}
{}_{\mathscr{G}_N(\partial\Omega)^\ast}\big\langle\widetilde{\gamma}_D f,
\phi\big\rangle_{\mathscr{G}_N(\partial\Omega)}
&=\lim_{j\to\infty}{}_{\mathscr{G}_N(\partial\Omega)^\ast}\big\langle\widetilde{\gamma}_D f_j,
\phi\big\rangle_{\mathscr{G}_N(\partial\Omega)}
\nonumber\\[2pt] 
&=\lim_{j\to\infty}\big(\gamma_D f_j,\gamma_N g\big)_{L^2(\partial\Omega)}
\nonumber\\[2pt] 
&=\lim_{j\to\infty}\big\{(f_j,\Delta g)_{L^2(\Omega)}-(\Delta f_j,g)_{L^2(\Omega)}\big\}
\nonumber\\[2pt] 
&=(f,\Delta g)_{L^2(\Omega)}-(\Delta f,g)_{L^2(\Omega)},
\end{align}
finishing the proof of \eqref{eq:33f4iU}. 

Next, we deal with the claims in item $(vi)$. Pick an arbitrary  
$f\in H^{3/2}(\Omega)\cap\dom(A_{max,\Omega})$ such that $\gamma_N f=0$ and note that, by 
\eqref{We-Q.17}, the function $\gamma_D f$ is well defined and belongs to $\mathscr{G}_D(\partial\Omega)$, 
which is a reflexive Banach space (cf.~\eqref{4.55EEdc}). As such the norm of 
$\gamma_D f\in\mathscr{G}_D(\partial\Omega)=\big(\mathscr{G}_D(\partial\Omega)^\ast\big)^\ast$ 
may be computed as
\begin{equation}\label{eq:Dfav76}
\|\gamma_D f\|_{\mathscr{G}_D(\partial\Omega)}=
\sup_{\substack{\xi\in\mathscr{G}_D(\partial\Omega)^\ast\\ \|\xi\|_{\mathscr{G}_D(\partial\Omega)^\ast}\leq 1}}
\Big|{}_{\mathscr{G}_D(\partial\Omega)^\ast}\big\langle\xi,\gamma_D f
\big\rangle_{\mathscr{G}_D(\partial\Omega)}\Big|.
\end{equation}
One recalls from part $(ii)$ that the operator 
$\widetilde{\gamma}_N:\dom(A_{max,\Omega})\rightarrow\mathscr{G}_D(\partial\Omega)^*$ 
is linear, surjective, and continuous when $\dom(A_{max,\Omega})$ is equipped with the natural graph 
norm $f\mapsto\|f\|_{L^2(\Omega)}+\|\Delta f\|_{L^2(\Omega)}$ (or with any of the other equivalent norms 
$f\mapsto\|f\|_{H^s(\Omega)}+\|\Delta f\|_{L^2(\Omega)}$, $s\in[0,\tfrac{3}{2}]$; cf.~\eqref{We-Q.10EE-jussi}). 
As a consequence of this and the Open Mapping Theorem it then follows that there exists a constant 
$C\in(0,\infty)$ with the property that 
\begin{align}\label{We-EHba.1w}
\begin{split}
& \text{for each $\xi\in\mathscr{G}_D(\partial\Omega)^\ast$ satisfying 
$\|\xi\|_{\mathscr{G}_D(\partial\Omega)^\ast}\leq 1$ there exists}
\\[2pt]
& \quad\text{$g\in\dom(A_{max,\Omega})$ with $\widetilde{\gamma}_N g=\xi$ and
$\|g\|_{L^2(\Omega)}+\|\Delta g\|_{L^2(\Omega)}\leq C$.}
\end{split}
\end{align}
Given now an arbitrary $\xi\in\mathscr{G}_D(\partial\Omega)^\ast$ with 
$\|\xi\|_{\mathscr{G}_D(\partial\Omega)^\ast}\leq 1$, let $g$ be as in 
\eqref{We-EHba.1w} and compute 
\begin{align}\label{We-EHba.1w22}
\big|{}_{\mathscr{G}_D(\partial\Omega)^\ast}\big\langle\xi,\gamma_D f
\big\rangle_{\mathscr{G}_D(\partial\Omega)}\big|
&=\big|{}_{\mathscr{G}_D(\partial\Omega)^\ast}\big\langle\widetilde{\gamma}_N g,\gamma_D f
\big\rangle_{\mathscr{G}_D(\partial\Omega)}\big|
\nonumber\\[2pt] 
&=\big|(g,\Delta f)_{L^2(\Omega)}-(\Delta g,f)_{L^2(\Omega)}\big|
\nonumber\\[2pt] 
&\leq\big(\|f\|_{L^2(\Omega)}+\|\Delta f\|_{L^2(\Omega)}\big)
\big(\|g\|_{L^2(\Omega)}+\|\Delta g\|_{L^2(\Omega)}\big)
\nonumber\\[2pt] 
&\leq C\big(\|f\|_{L^2(\Omega)}+\|\Delta f\|_{L^2(\Omega)}\big)
\nonumber\\[2pt] 
&\leq C\big(\|f\|_{H^{s}(\Omega)}+\|\Delta f\|_{L^2(\Omega)}\big)
\end{align}
for $s\in [0,\tfrac{3}{2}]$, where the second equality above is a consequence of 
\eqref{eq:33f4iU.NN}--\eqref{eq:33VaV.NN}. Together, \eqref{eq:Dfav76} and \eqref{We-EHba.1w22} 
yield 
\begin{equation}\label{eq:Dfav76.BB}
\|\gamma_D f\|_{\mathscr{G}_D(\partial\Omega)}
\leq C\big(\|f\|_{H^s(\Omega)}+\|\Delta f\|_{L^2(\Omega)}\big),\quad s\in\big[0,\tfrac{3}{2}\big],
\end{equation}
proving the continuity of the operator $\gamma_D$ in \eqref{4.5ASdf.1}. 
The case of the operator $\gamma_N$ in \eqref{4.5ASdf.2} is handled similarly. 
Continuing the treatment of $(vi)$, one observes that the claim in \eqref{eq:Ajf} is 
a direct consequence of Theorem~\ref{t4.6}, while the equivalences in 
\eqref{i7g5r}--\eqref{i7g5r.2N} are seen from the surjectivity of the operators 
in \eqref{4.5ASdf.1}--\eqref{4.5ASdf.2}, Lemma~\ref{uyrre}, and 
\eqref{4.55.REE.D}--\eqref{4.55.REE.N}; the last equivalence in 
\eqref{i7g5r.2N} is due to the fact that the Dirichlet Laplacian is strictly positive. 
Next, \eqref{We-Q.17} yields $\mathscr{G}_D(\partial\Omega)\subset H^1(\partial\Omega)$ 
and $\mathscr{G}_N(\partial\Omega)\subset L^2(\partial\Omega)$. Given any 
$\phi\in\mathscr{G}_D(\partial\Omega)$, making use of \eqref{i7g5r} 
and the boundedness of $\gamma_D$ in \eqref{eqn:gammaDs.2} with $s=\tfrac{3}{2}$, one obtains 
\begin{align}\label{i7g5r-BBB}
\|\phi\|_{\mathscr{G}_D(\partial\Omega)}& \geq C\inf_{\substack{f\in
H^{3/2}(\Omega)\cap\dom(A_{max,\Omega})\\ \gamma_N f=0,\,\,\gamma_D f=\phi}}
\big(\|f\|_{H^{3/2}(\Omega)}+\|\Delta f\|_{L^2(\Omega)}\big)
\nonumber\\[2pt]
& \geq C\inf_{\substack{f\in
H^{3/2}(\Omega)\cap\dom(A_{max,\Omega})\\ \gamma_N f=0,\,\,\gamma_D f=\phi}}
\big(\|\gamma_D f\|_{L^2(\partial\Omega)}\big)
\nonumber\\[2pt]
&=C\|\phi\|_{H^1(\partial\Omega)},
\end{align}
for some constant $C\in(0,\infty)$ independent of $\phi$. This proves that the inclusion
$\mathscr{G}_D(\partial\Omega)\hookrightarrow H^1(\partial\Omega)$ is continuous, 
and a similar argument shows that the inclusion 
$\mathscr{G}_N(\partial\Omega)\hookrightarrow L^2(\partial\Omega)$ 
continuously as well. Since these inclusions also have dense ranges (cf.\  Lemma~\ref{l5.3}), 
the claims pertaining \eqref{ytrr555e} follow with the help of Lemma~\ref{uyrre.LD}
(also keeping \eqref{4.55EEdc} in mind). 

Next, the claims $(vii)$ and $(viii)$ follow from item $(ii)$ and the direct sum decompositions 
\begin{align}\label{yarf7h5f}
\dom(A_{max,\Omega}) &=\dom(A_{D,\Omega})\dotplus\ker(A_{max,\Omega}-zI)  
\nonumber\\[2pt]
&=\ker(\widetilde{\gamma}_D)\dotplus\ker(A_{max,\Omega}-zI),\quad\forall\,z\in\rho(A_{D,\Omega}), 
\\[2pt] 
\dom(A_{max,\Omega}) &=\dom(A_{N,\Omega})\dotplus\ker(A_{max,\Omega}-zI)  
\nonumber\\[2pt] 
&=\ker(\widetilde{\gamma}_N)\dotplus\ker(A_{max,\Omega}-zI),\quad\forall\,z\in\rho(A_{N,\Omega}).
\label{yarf7h5f.BB}
\end{align}
Finally, statement $(ix)$ is a consequence of items $(ii)$, $(vii)$ and $(viii)$; 
see also \cite[Corollary~4.2]{BM13}.
\end{proof}
%%%%%%%%%%

In the following remarks we will elaborate on the links to abstract boundary triples and their 
$\gamma$-fields and Weyl functions from extension theory of symmetric operators.
%%%%%%%%%
\begin{remark}\label{jussirem}
Consider the operator 
\begin{equation}\label{t32jussi}
T_{3/2,\Omega}=-\Delta+V,\quad\dom(T_{3/2,\Omega})=H^{3/2}(\Omega)\cap\dom(A_{max,\Omega}),
\end{equation}
and note that Lemma~\ref{l5.1MM} and Lemma~\ref{l3.2} imply $\overline{T_{3/2,\Omega}}=A_{max,\Omega}=A_{min,\Omega}^*$. 
It is immediate from Corollary~\ref{YTfdf-T.NNN} and Corollary~\ref{YTfdf.NNN.2} for $s=3/2$ 
that $f,g\in\dom (T_{3/2,\Omega})$ satisfy $\gamma_D f,\gamma_D g\in H^1(\partial\Omega)$ and 
$\gamma_N f,\gamma_N g\in L^2(\partial\Omega)$. Furthermore, the following Green's formula is a consequence 
of Corollary~\ref{YTfdf.NNN.2}\,$(i)$ with $s=3/2$, bearing in mind that $\gamma_N f\in L^2(\partial\Omega)$:
\begin{equation}\label{green32jussi}
\begin{split}
&(T_{3/2,\Omega}f,g)_{L^2(\Omega)}-(f,T_{3/2,\Omega} g)_{L^2(\Omega)}
\\
&\qquad\qquad=(\gamma_D f,\gamma_N g)_{L^2(\partial\Omega)}-{}_{(H^{1}(\partial\Omega))^*}\big\langle\gamma_N f,\gamma_D g
\big\rangle_{H^{1}(\partial\Omega)}
\\
&\qquad\qquad=(\gamma_D f,\gamma_N g)_{L^2(\partial\Omega)}-(\gamma_N f,\gamma_D g)_{L^2(\partial\Omega)}.
\end{split}
\end{equation}
Observe that $\mathscr{G}_D(\partial\Omega)\times \{0\}$ and $\{0\}\times\mathscr{G}_N(\partial\Omega)$ are both contained 
in the range of the map
\begin{equation}\label{rangegamma}
(\gamma_D,-\gamma_N):\dom(T_{3/2,\Omega})\rightarrow L^2(\partial\Omega)\times L^2(\partial\Omega)
\end{equation}
by Lemma~\ref{l5.3}, hence the range of \eqref{rangegamma} is dense. Furthermore, for $s=3/2$ 
Corollary~\ref{YTfdf.NNN.2} shows that $\gamma_N:\dom(T_{3/2,\Omega})\rightarrow L^2(\partial\Omega)$ is surjective. 
It is also clear from Theorem~\ref{t4.4} and Theorem~\ref{t4.5} that
\begin{equation}\label{xcv}
\begin{split}
A_{D,\Omega}&=T_{3/2,\Omega}\upharpoonright\big\{f\in \dom(T_{3/2,\Omega})\,\big|\,\gamma_D f=0\big\},
\\
A_{N,\Omega}&=T_{3/2,\Omega}\upharpoonright\big\{f\in \dom(T_{3/2,\Omega})\,\big|\,\gamma_N f=0\big\},
\end{split}
\end{equation}
are both self-adjoint restrictions of the operator $T_{3/2,\Omega}$ in $L^2(\Omega)$. 

From the above observations it follows that $\{L^2(\partial\Omega),\gamma_D,-\gamma_N\}$ is a so-called 
quasi boundary triple for $T_{3/2,\Omega}\subset A_{max,\Omega}$ with corresponding $\gamma$-field $P_{1,D,\Omega}$ 
and Weyl function $M_\Omega=M_{1,\Omega}$ from Theorem~\ref{t5.2}\,$(i)$ and $(iii)$ (see \cite{BL07,BL12}). The transposed triple 
$\{L^2(\partial\Omega),\gamma_N,\gamma_D\}$ is even a $B$-generalized boundary triple for $T_{3/2,\Omega}\subset A_{max,\Omega}$ 
with corresponding $\gamma$-field $P_{1,N,\Omega}$ and Weyl function $N_{1,\Omega}$ from Theorem~\ref{t5.2}\,$(ii)$ and $(iv)$
(see \cite{DHMS06,DM95}). The abstract theory of quasi boundary triples and $B$-generalized boundary triples yields the 
continuity of the $\gamma$-fields as mappings from $L^2(\partial\Omega)$ to $L^2(\Omega)$ and the representations
of the adjoints in Theorem~\ref{t5.2}\,$(i)$ and $(ii)$. Similarly, the formulas \eqref{We-Q.19} and \eqref{eqn:imMinv} 
in Theorem~\ref{t5.4} for the imaginary parts of $M_{1,\Omega}$ and $N_{1,\Omega}=-M_{1,\Omega}^{-1}$ (see Lemma~\ref{t5.4part1}) 
reflect the connection between the $\gamma$-field and Weyl function of a quasi boundary triple or $B$-generalized boundary triple. 

In this context we mention that the extension of the Dirichlet trace operator $\gamma_D$ and Neumann trace operator 
$\gamma_N$ onto $\dom(A_{max,\Omega})$ in Theorem~\ref{t5.5}\,$(ii)$ is based on an abstract technique developed for 
quasi boundary triples in \cite{BM13}. In the case of Schr\"odinger operators on bounded Lipschitz domains this 
method gives rise to a certain regularization of the Neumann trace operator such that a modified second Green's identity 
holds on $\dom(A_{max,\Omega})$. Using $\widetilde\gamma_D$ in Theorem~\ref{t5.5}\,$(ii)$ and replacing the Neumann trace 
operator $\widetilde\gamma_N$ by such a regularized version leads to an ordinary boundary triplet; cf. \cite{BM13} for details.
For domains with smooth boundary the corresponding construction of a boundary triple
(including regularization) and parametrization of all proper extensions was proposed in different
manners in \cite{Vi63} and \cite{Gr68} (see also \cite{Ma10}). Besides, the corresponding $\gamma$-field and the Weyl
function $M$ corresponding to this ordinary boundary triple were computed in \cite{Ma10}.
\end{remark}
%%%%%%%%%

%%%%%%%%%
\begin{remark}\label{jussirem2}
Consider the operator 
\begin{equation}\label{t1jussi}
T_{1,\Omega}=-\Delta+V, \quad \dom(T_{1,\Omega})= H^{1}(\Omega)\cap\dom(A_{max,\Omega}).
\end{equation}
As in Remark~\ref{jussirem} we have $\overline{T_{1,\Omega}}=A_{max,\Omega}=A_{min,\Omega}^*$ and it follows from
Corollary~\ref{YTfdf-T.NNN} and Corollary~\ref{YTfdf.NNN.2} for $s=1$ that $f,g\in\dom (T_{1,\Omega})$ satisfy 
$\gamma_D f,\gamma_D g\in H^{1/2}(\partial\Omega)$ and $\gamma_N f,\gamma_N g\in H^{-1/2}(\partial\Omega)$. 
Furthermore, Corollary~\ref{YTfdf.NNN.2}\,$(i)$ with $s:=1$ shows that
\begin{equation}\label{green1jussi}
\begin{split}
&(T_{1,\Omega} f,g)_{L^2(\Omega)}-(f,T_{1,\Omega} g)_{L^2(\Omega)}\\
&\qquad= {}_{H^{1/2}(\partial\Omega)}\big\langle\gamma_D f,\gamma_N g
\big\rangle_{(H^{1/2}(\partial\Omega))^*}
-{}_{(H^{1/2}(\partial\Omega))^*}\big\langle\gamma_N f,\gamma_D g\big\rangle_{H^{1/2}(\partial\Omega)}
\end{split}
\end{equation}
for all $f,g\in\dom(T_{1,\Omega})$. Since 
$H^{1/2}(\partial\Omega)\hookrightarrow L^2(\partial\Omega)\hookrightarrow(H^{1/2}(\partial\Omega))^*$ 
we can fix a uniformly positive self-adjoint operator $\jmath$ in $L^2(\partial\Omega)$ with 
$\dom(\jmath)=H^{1/2}(\partial\Omega)$ such that $\jmath:H^{1/2}(\partial\Omega)\rightarrow L^2(\partial\Omega)$ 
is an isomorphism and $\jmath^{-1}$ admits an extension to an isomorphism
$\widetilde{\jmath^{-1}}:(H^{1/2}(\partial\Omega))^*\rightarrow L^2(\partial\Omega)$ and the duality pairing 
between $H^{1/2}(\partial\Omega)$ and $(H^{1/2}(\partial\Omega))^*$ is compatible with the scalar product in 
$L^2(\partial\Omega)$. Hence \eqref{green1jussi} can be written in the form
\begin{equation}\label{green1jussi2}
\begin{split}
&(T_{1,\Omega}f,g)_{L^2(\Omega)}-(f,T_{1,\Omega}g)_{L^2(\Omega)}
\\
&\qquad=\big(\jmath \gamma_D f,\widetilde{\jmath^{-1}}\gamma_N g\big)_{L^2(\partial\Omega)}
-\big(\widetilde{\jmath^{-1}}\gamma_N f,\jmath\gamma_D g\big)_{L^2(\partial\Omega)}
\end{split}
\end{equation}
for all $f,g\in\dom(T_{1,\Omega})$. Furthermore, the mappings 
\begin{equation*}
\jmath\gamma_D:\dom(T_{1,\Omega})\rightarrow L^2(\partial\Omega)\quad\text{and}\quad 
\widetilde{\jmath^{-1}}\gamma_N:\dom(T_{1,\Omega})\rightarrow L^2(\partial\Omega)
\end{equation*}
are both surjective by Corollary~\ref{YTfdf-T.NNN}, Corollary~\ref{YTfdf.NNN.2}, and the properties of $\jmath$ and 
$\widetilde{\jmath^{-1}}$. As in \eqref{xcv} one sees that
\begin{equation}\label{xcv2}
\begin{split}
A_{D,\Omega}&=T_{1,\Omega}\upharpoonright\big\{f\in \dom(T_{1,\Omega})\,\big|\,\jmath\gamma_D f=0\big\},
\\
A_{N,\Omega}&=T_{1,\Omega}\upharpoonright\big\{f\in \dom(T_{1,\Omega})\,\big|\, \widetilde{\jmath^{-1}}\gamma_N f=0\big\},
\end{split}
\end{equation}
are both self-adjoint restrictions of the operator $T_{1,\Omega}$. Therefore, it follows that 
$\{L^2(\partial\Omega),\jmath \gamma_D,-\widetilde{\jmath^{-1}}\gamma_N\}$ is a so-called double 
$B$-generalized boundary triple for $T_{1,\Omega}\subset A_{max,\Omega}$ in the sense of \cite[Definition 2.1]{BMN17}.
The corresponding $\gamma$-field is given by $P_{1/2,D,\Omega}(\cdot)\jmath^{-1}$ and the corresponding Weyl function 
is given by $\widetilde{\jmath^{-1}}M_{1/2,\Omega}(\cdot)\jmath^{-1}$.
In the case of a smooth boundary such a double $B$-generalized boundary triple was
constructed in \cite{BMN17} and the corresponding $\gamma$-field and Weyl function were also
provided there.
\end{remark}
%%%%%%%%%%

%%%%%%%%%%%%%%%%%%%%%%%%%%%%%%
%%%%%%%%%%%%%%%%%%%%%%%%%%%%%%
\section{The Krein--von Neumann Extension on Bounded Lipschitz Domains} \label{s9}
%%%%%%%%%%%%%%%%%%%%%%%%%%%%%%
%%%%%%%%%%%%%%%%%%%%%%%%%%%%%%

The principal purpose of this section is to describe the Krein--von Neumann extension for perturbed Laplacians on bounded Lipschitz domains. Special emphasis is given to its spectral properties, the corresponding boundary conditions in terms of extended Dirichlet and Neumann traces and the Dirichlet-to-Neumann map at $z=0$, Krein-type resolvent formulas connecting the Krein--von Neumann and Dirichlet resolvent, and finally to 
the Weyl asymptotics of perturbed Krein Laplacians. 

In this section we now strengthen Hypothesis~\ref{h4.2} by assuming, in addition, that $V\in L^\infty(\Omega)$ is nonnegative~a.e. 

%%%%%%%%%
\begin{hypothesis}\label{h6.1} 
Let $n\in\bbN\backslash\{1\}$, assume that $\Omega\subset\bbR^n$ is a bounded Lipschitz domain, 
and suppose that $V\in L^\infty(\Omega)$ is nonnegative~a.e. 
\end{hypothesis}
%%%%%%%%%

It then follows from Lemma~\ref{l3.3} that the minimal operator 
\begin{equation}\label{aminnon}
A_{min,\Omega}=-\Delta+V,\quad\dom(A_{min,\Omega})=\accentset{\circ}{H}^2(\Omega),
\end{equation}
is strictly positive, and the same holds for the Friedrichs extension $A_{F,\Omega}$ of 
$A_{min,\Omega}$ by Theorem~\ref{tfried}. One recalls from the paragraph 
preceding Theorem~\ref{t4.4} that $A_{F,\Omega}$ coincides with the 
Dirichlet realization $A_{D,\Omega}$ of $-\Delta+V$. 
Next, we recall that the {\it Krein--von Neumann} extension
$A_{K,\Omega}$ of $A_{min,\Omega}$ is given by
\begin{equation}\label{We-Q.9jussi}
A_{K,\Omega}=-\Delta+V,\quad\dom(A_{K,\Omega})=\dom(A_{min,\Omega})\,\dotplus\,\ker(A_{max,\Omega}). 
\end{equation}
We remark that, collectively, the functions in $\dom(A_{K,\Omega})$ do not possess any additional 
Sobolev regularity, that is, $\dom(A_{K,\Omega})\not\subset H^s(\Omega)$ for any $s>0$.

In the following theorem we briefly collect some well-known properties of the Krein--von Neumann extension 
$A_{K,\Omega}$ which were shown by M.G. Krein in \cite{Kr47} (see also \cite{AS80,AGMT10,AGMST10,AGMST13}, 
and \cite[Section~2]{GLMS14}).

%%%%%%%%%%
\begin{theorem}\label{t6.3}
Assume Hypothesis~\ref{h6.1} and let $A_{K,\Omega}$ be the Krein--von Neumann extension of $A_{min,\Omega}$. 
Then the following assertions hold: \\[1mm]
\noindent $(i)$ $A_{K,\Omega}$ is a nonnegative self-adjoint operator in $L^2(\Omega)$ and $\sigma(A_{K,\Omega})$ 
consists of eigenvalues only. In addition, the eigenvalue $\lambda=0$ has infinite multiplicity, 
$\dim(\ker(A_{K,\Omega}))=\infty$, and the restriction $A_{K,\Omega}|_{(\ker(A_{K,\Omega}))^\bot}$ 
is a strictly positive self-adjoint operator in the Hilbert space $(\ker(A_{K,\Omega}))^\bot$ 
with compact resolvent. \\[1mm] 
\noindent $(ii)$ A nonnegative self-adjoint operator $A_{\Omega}$ in $L^2(\Omega)$ is a 
self-adjoint extension of $A_{min,\Omega}$ if and only if 
\begin{equation}\label{afakjussi}
(A_{D,\Omega}-\mu)^{-1}\leq(A_{\Omega}-\mu)^{-1}\leq(A_{K,\Omega}-\mu)^{-1}
\end{equation}
holds for some $($and, hence for all\,$)$ $\mu<0$.
\end{theorem}
%%%%%%%%%%

We note that \eqref{afakjussi} is equivalent to the inequality  
$A_{K,\Omega}\leq A_{\Omega}\leq A_{F,\Omega}$, when interpreted in the sense of 
quadratic forms (see \cite[Section~I.6]{Fa75} and \cite[Theorem~VI.2.21]{Ka80}).
In the next lemma we explicitly verify that the Dirichlet and the Krein--von Neumann extension 
are relatively prime (or disjoint), see, for instance, \cite[Lemma 2.8]{AGMT10}.

%%%%%%%%%
\begin{lemma}
Assume Hypothesis~\ref{h6.1} and let $A_{D,\Omega}$ be the Dirichlet extension and let 
$A_{K,\Omega}$ be  the Krein--von Neumann extension of $A_{min,\Omega}$ in \eqref{We-Q.9jussi}. 
Then 
\begin{equation}\label{ytrre}
\dom(A_{D,\Omega})\cap\dom(A_{K,\Omega})=\dom(A_{min,\Omega})=\accentset{\circ}{H}^2(\Omega).
\end{equation}
\end{lemma}
%%%%%%%%%
\begin{proof}
Suppose that $f\in\dom(A_{D,\Omega})\cap\dom(A_{K,\Omega})$ and decompose $f$ according to 
\eqref{We-Q.9jussi} in the form $f=f_{min}+f_0$ with $f_{min}\in\dom(A_{min,\Omega})$ and 
$f_0\in\ker(A_{max,\Omega})$. It follows that $f_0\in\dom(A_{D,\Omega})\cap\ker(A_{max,\Omega})$ 
and since $A_{D,\Omega}$ is strictly positive one concludes that $f_0=0$. 
Thus $f=f_{min}\in\dom(A_{min,\Omega})$. The inclusion
\begin{equation}\label{clear}
\dom(A_{min,\Omega})\subset\big(\dom(A_{D,\Omega})\cap\dom(A_{K,\Omega})\big)
\end{equation}
is clear as both $A_{D,\Omega}$ and $A_{K,\Omega}$ are extensions of $A_{min,\Omega}$. 
The last equality in \eqref{ytrre} was shown in Lemma~\ref{l4.3}. 
\end{proof} 
%%%%%%%%

Alternatively, this result follows abstractly from \cite[Lemma 2.8]{AGMT10} upon noting that the Dirichlet, $A_{D,\Omega}$, and the Friedrichs realization, $A_{F,\Omega}$, of $A_{min,\Omega}$, coincide (cf.\ \eqref{We-Q.10EE}).

Our next goal is to obtain an explicit description of the domain of the Krein--von Neumann 
extension $A_{K,\Omega}$ in terms of the extended Dirichlet and Neumann trace operators 
in Theorem~\ref{t5.5}. The Dirichlet-to-Neumann map at $z=0$ will enter as regularization 
parameter here. One observes that $M_{\Omega}(0)$ and its extension $\widetilde{M}_{\Omega}(0)$ 
in the context of Theorem~\ref{t5.5} are well defined as $A_{D,\Omega}$ is strictly positive 
by Theorem~\ref{t4.4}. We mention that for smooth domains and elliptic differential 
operators with smooth coefficients, this description of the Krein--von Neumann extension $A_{K,\Omega}$ goes back to a remarkably early 1952 paper (translated into English in 1963) by Vi{\u s}ik \cite{Vi63}, followed by work of Grubb \cite{Gr68} in 1968. For quasi-convex domains, Theorem~\ref{t5.6} below coincides with \cite[Theorem 5.5]{AGMST10} and \cite[Theorem 13.1]{GM11}; for the abstract setting we refer to \cite[Example 3.9]{BM13}. For Lipschitz domains this result was recently obtained in \cite[Theorem 3.3]{BGMM14}.

%%%%%%%%%%
\begin{theorem}\label{t5.6}
Assume Hypothesis~\ref{h6.1} and let $\widetilde\gamma_D$, $\widetilde\gamma_N$ and 
$\widetilde M_{\Omega}$ be as in Theorem~\ref{t5.5}. Then the  Krein--von Neumann 
extension $A_{K,\Omega}$ of $A_{min,\Omega}$ is given by 
\begin{align} 
\begin{split} 
& A_{K,\Omega}=-\Delta+V,    
\\[2pt] 
& \dom(A_{K,\Omega})=\big\{f\in\dom(A_{max,\Omega})\,\big|\, 
\widetilde{\gamma}_N f+\widetilde M_{\Omega}(0)\widetilde{\gamma}_D f=0\big\}.
\label{eqn:A_K}
\end{split} 
\end{align}
\end{theorem}
%%%%%%%%%%
\begin{proof}   
Let $A_{K,\Omega}$ be the Krein--von Neumann extension of $A_{min,\Omega}$ and note that
\begin{equation}\label{akakdom}
\dom(A_{K,\Omega})=\dom(A_{min,\Omega})\,
\dotplus\,\ker(A_{max,\Omega})=\accentset{\circ}{H}^2(\Omega)\,
\dotplus\,\ker(A_{max,\Omega})
\end{equation}
by \eqref{We-Q.9jussi} and \eqref{aminnon}. Consider $f\in\dom(A_{K,\Omega})$. 
Then $f\in\dom(A_{max,\Omega})$ and by \eqref{akakdom} $f$ can be decomposed in the 
form $f=f_{min}+f_0$, where $f_{min}\in\accentset{\circ}{H}^2(\Omega)$ and $f_0\in\ker(A_{max,\Omega})$. 
Thus, $\gamma_D f_{min}=\widetilde{\gamma}_D f_{min}=0$ and 
$\gamma_N f_{min}=\widetilde{\gamma}_N f_{min}=0$, and hence it follows from 
Theorem~\ref{t5.5}\,$(vii)$,\,$(ix)$ that
\begin{align}\label{7yhBG}
\widetilde{M}_{\Omega}(0)\widetilde{\gamma}_D f
&=\widetilde{M}_{\Omega}(0)\widetilde{\gamma}_D(f_{min}+f_0)
=\widetilde{M}_{\Omega}(0)\widetilde{\gamma}_D f_0     
\\[2pt]
&=-\widetilde{\gamma}_N f_0=-\widetilde{\gamma}_N(f_{min}+f_0)
=-\widetilde{\gamma}_N f.
\end{align}
Hence, 
\begin{equation}\label{o7554}
\dom(A_{K,\Omega})\subseteq\big\{f\in\dom(A_{max,\Omega})\,\big|\,\widetilde{\gamma}_N f
+\widetilde{M}_{\Omega}(0)\widetilde{\gamma}_D f=0\big\}.
\end{equation}

Next, we verify the opposite inclusion of the domains in \eqref{eqn:A_K}. 
To this end, pick $f\in\dom(A_{max,\Omega})$ satisfying the boundary condition
$\widetilde{M}_{\Omega}(0)\widetilde{\gamma}_D f+\widetilde{\gamma}_N f=0$. 
According to the decomposition \eqref{eqn:decAmax} one can write $f$ 
in the form 
$f=f_D+f_0$, where $f_D\in\dom(A_{D,\Omega})$ and $f_0\in\ker(A_{max,\Omega})$. 
Then $\gamma_D f_D=\widetilde{\gamma}_D f_D=0$ and with the help of 
Theorem~\ref{t5.5}\,$(vii)$,\,$(ix)$ one computes  
\begin{equation}\label{u654}
\widetilde{M}_{\Omega}(0)\widetilde{\gamma}_D f
=\widetilde{M}_{\Omega}(0)\widetilde{\gamma}_D(f_D+f_0)
=\widetilde{M}_{\Omega}(0)\widetilde\gamma_D f_0=-\widetilde{\gamma}_N f_0.
\end{equation} 
Taking into account the boundary condition 
$\widetilde{M}_{\Omega}(0)\widetilde{\gamma}_D f=-\widetilde{\gamma}_N f$ one obtains  
\begin{equation}\label{Lo7g44}
0=\widetilde{\gamma}_N(f-f_0)=\widetilde{\gamma}_N f_D,
\end{equation}
and hence $f_D\in\ker(\widetilde{\gamma}_N)=\ker(\gamma_N)=\dom(A_{N,\Omega})$ 
(cf.\ Theorem~\ref{t5.5}\,$(ii)$). Thus, making use of Theorem~\ref{t4.6} and 
\eqref{aminnon} one obtains 
\begin{equation}\label{o765f4}
f_D\in\dom(A_{D,\Omega})\cap\dom(A_{N,\Omega})=\dom(A_{min,\Omega})
=\accentset{\circ}{H}^2(\Omega), 
\end{equation}
implying $f=f_D+f_0\in\accentset{\circ}{H}^2(\Omega)\dotplus\ker(A_{max,\Omega})$, 
that is, $f\in\dom(A_{K,\Omega})$. 
\end{proof}
%%%%%%%%%%

Next, we prove a variant of Krein's resolvent formula relating the resolvent of the 
Krein--von Neumann extension $A_{K,\Omega}$ to the resolvent of the Dirichlet (and hence, Friedrichs) realization 
$A_{D,\Omega}$. For variants of Krein's formula discussed here see \cite{BL07}, \cite{BL12}, 
\cite{BMN17}, \cite{BM13}, \cite{BGW09}, \cite{GM08}, \cite{Ma10}, \cite{PR09}, and Section~\ref{s10}.

%%%%%%%%%%
\begin{theorem}\label{YrarTfc}
Assume Hypothesis~\ref{h6.1}, and let $A_{K,\Omega}$ be the Krein--von Neumann extension 
of $A_{min,\Omega}$. Let $\widetilde{P}_{D,\Omega}(z)$ be the solution operator of the 
boundary value problem \eqref{eqn:bvp.2} and let $\widetilde{M}_\Omega(z)$ be the extended 
Dirichlet-to-Neumann map in \eqref{iutee}. Then, for each $z\in\rho(A_{K,\Omega})\cap\rho(A_{D,\Omega})$,
\begin{equation}\label{mtildebij}
\widetilde{M}_\Omega(z)-\widetilde{M}_\Omega(0):\mathscr{G}_N(\partial\Omega)^*
\rightarrow\mathscr{G}_D(\partial\Omega)^*
\end{equation}
is a linear, continuous, injective mapping, with range 
\begin{equation}\label{range}
\ran\big(\widetilde{M}_\Omega(z)-\widetilde{M}_\Omega(0)\big)=\mathscr{G}_N(\partial\Omega).
\end{equation}

Moreover, for each $z\in\rho(A_{K,\Omega})\cap\rho(A_{D,\Omega})$, the operator 
\begin{align}\label{eq:RVav-9jj}
\begin{split}
& \widetilde{M}_\Omega(z)-\widetilde{M}_\Omega(0):\mathscr{G}_N(\partial\Omega)^*
\rightarrow\mathscr{G}_N(\partial\Omega)
\\[2pt]
& \quad\text{is a continuous linear isomorphism}
\end{split}
\end{align}
and, with $\big(\widetilde{M}_\Omega(z)-\widetilde{M}_\Omega(0)\big)^{-1}
\in{\mathcal{B}}\big(\mathscr{G}_N(\partial\Omega),\mathscr{G}_N(\partial\Omega)^*\big)$, 
the following Krein-type resolvent formula holds in ${\mathcal{B}}\big(L^2(\Omega)\big)$:
\begin{align}\label{resform}
& (A_{K,\Omega}-zI)^{-1}-(A_{D,\Omega}-zI)^{-1}
\nonumber\\[2pt]
& \quad=-\widetilde{P}_{D,\Omega}(z)\big(\widetilde{M}_\Omega(z)-\widetilde{M}_\Omega(0)\big)^{-1}
\big(\widetilde{P}_{D,\Omega}(\bar{z})\big)^*,
\end{align}
where $\big(\widetilde{P}_{D,\Omega}(z)\big)^*\in{\mathcal{B}}\big(L^2(\Omega),\mathscr{G}_N(\partial\Omega)\big)$ 
is the adjoint of the operator $\widetilde{P}_{D,\Omega}(z)$ in \eqref{6h8u}
{\rm (}viewed here as a linear and continuous mapping from $\mathscr{G}_N(\partial\Omega)^*$ 
into $L^2(\Omega)${\rm )}.  
\end{theorem}
%%%%%%%%%%
\begin{proof}
Fix $z\in\rho(A_{K,\Omega})\cap\rho(A_{D,\Omega})$. We start by noting  
that \eqref{iutee} and the fact that $0\in\rho(A_{D,\Omega})$ guarantee that the operator 
$\widetilde{M}_\Omega(z)-\widetilde{M}_\Omega(0)$ in \eqref{mtildebij} is well defined,
linear, and continuous. To see that $\widetilde{M}_\Omega(z)-\widetilde{M}_\Omega(0)$ 
is also injective, assume that $\varphi\in\mathscr{G}_N(\partial\Omega)^*$ is such that
\begin{equation}\label{mtildebij2}
\big(\widetilde{M}_\Omega(z)-\widetilde{M}_\Omega(0)\big)\varphi=0
\,\text{ in }\,\mathscr{G}_D(\partial\Omega)^*.
\end{equation}
By design,
\begin{equation}\label{okgut22}
\widetilde{f}_{D,\Omega}(z,\varphi):=\widetilde{P}_{D,\Omega}(z)\varphi\in\dom(A_{max,\Omega})
\end{equation}
is the unique solution of the boundary value problem \eqref{eqn:bvp.2}, hence 
\begin{align}\label{bcbc}
\begin{split}
& \widetilde{\gamma}_D\widetilde{f}_{D,\Omega}(z,\varphi)
=\widetilde{\gamma}_D\widetilde{P}_{D,\Omega}(z)\varphi=\varphi,
\\[2pt]
& \quad\text{and $\widetilde{f}_{D,\Omega}(z,\varphi)\in\ker(A_{max,\Omega}-zI)$}.
\end{split}
\end{align}
It follows from \eqref{iutee}, \eqref{mtildebij2}, and \eqref{bcbc}, that
\begin{align}\label{uredvhu}
\widetilde{\gamma}_N\widetilde{f}_{D,\Omega}(z,\varphi)
&=\widetilde{\gamma}_N\widetilde{P}_{D,\Omega}(z)\varphi
=-\widetilde{M}_{\Omega}(z)\varphi=-\widetilde{M}_{\Omega}(0)\varphi
\nonumber\\[2pt]
&=-\widetilde{M}_{\Omega}(0)\widetilde{\gamma}_D\widetilde{f}_{D,\Omega}(z,\varphi).
\end{align}
Consequently, $\widetilde{f}_{D,\Omega}(z,\varphi)\in\dom(A_{K,\Omega})$ by \eqref{uredvhu} 
and \eqref{eqn:A_K}. Given this fact and keeping in mind \eqref{bcbc} one  deduces that 
$\widetilde{f}_{D,\Omega}(z,\varphi)\in\ker(A_{K,\Omega}-z I)$. 
In turn, this forces $\widetilde{f}_{D,\Omega}(z,\varphi)=0$, given that we are 
presently assuming $z\in\rho(A_{K,\Omega})$. With this in hand, by once again appealing to 
\eqref{bcbc} one finally concludes that $\varphi=0$. Therefore, have 
shown that the operator 
$\widetilde{M}_\Omega(z)-\widetilde{M}_\Omega(0)$ in \eqref{mtildebij} is injective.

In order to prove the range condition in \eqref{range} one first notes that for 
$\varphi\in\mathscr{G}_N(\partial\Omega)^*$ one has, by definition (cf.~\eqref{iutee}), 
\begin{equation}\label{abcdef}
\big(\widetilde{M}_\Omega(z)-\widetilde{M}_\Omega(0)\big)\varphi
=-\widetilde{\gamma}_N\big(\widetilde{P}_{D,\Omega}(z)-\widetilde{P}_{D,\Omega}(0)\big)\varphi.
\end{equation}
On the other hand, $\widetilde{\gamma}_D\big(\widetilde{P}_{D,\Omega}(z)
-\widetilde{P}_{D,\Omega}(0)\big)\varphi=\varphi-\varphi=0$ which goes to show that 
\begin{equation}\label{eq:tgfd}
\big(\widetilde{P}_{D,\Omega}(z)-\widetilde{P}_{D,\Omega}(0)\big)\varphi
\in\ker(\widetilde{\gamma}_D)=\dom(A_{D,\Omega})\subset H^{3/2}(\Omega)\cap\dom(A_{max,\Omega})
\end{equation}
by the first relation in \eqref{4.30} and \eqref{We-Q.10EE-jussi}. This fact, \eqref{eq:RTTYH}, 
and the definition of $\mathscr{G}_N(\partial\Omega)$ in \eqref{We-Q.17}, imply that the function 
in \eqref{abcdef} belongs to $\mathscr{G}_N(\partial\Omega)$. This yields the left-to-right inclusion 
in \eqref{range}. In order to verify the right-to-left inclusion in \eqref{range}, consider some arbitrary 
$\psi\in\mathscr{G}_N(\partial\Omega)$. Then there exists a function $f\in H^{3/2}(\Omega)\cap\dom(A_{max,\Omega})$ 
such that $\gamma_D f=0$ and $\gamma_N f=\psi$ (cf.~\eqref{We-Q.17}). In particular,
\begin{equation}\label{y545f4}
\widetilde{\gamma}_N f+\widetilde{M}_{\Omega}(0)\widetilde{\gamma}_D f=\gamma_N f=\psi.
\end{equation}
Since $z\in\rho(A_{K,\Omega})$, this ensures the direct sum decomposition 
\begin{equation}\label{domdecoak}
\dom(A_{max,\Omega})=\dom(A_{K,\Omega})\dotplus\ker(A_{max,\Omega}-zI). 
\end{equation}
Using this in relation to the function $f\in\dom(A_{max,\Omega})$ and observing that we have
$\widetilde{\gamma}_N g+\widetilde{M}_{\Omega}(0)\widetilde{\gamma}_D g=0$ for each 
$g\in\dom(A_{K,\Omega})$ by Theorem~\ref{t5.6}, it follows from \eqref{y545f4} that
there exists 
\begin{equation}\label{eqn:bvpNeumannkn} 
\eta\in\ker(A_{max,\Omega}-zI)\,\text{ such that }\,
\widetilde{\gamma}_N\eta+\widetilde{M}(0)\widetilde{\gamma}_D\eta=\psi. 
\end{equation}
Setting $\varphi:=-\widetilde{\gamma}_D\eta\in\mathscr{G}_N(\partial\Omega)^*$, one concludes from  
\eqref{iutee} and \eqref{eqn:bvpNeumannkn} that
\begin{align}\label{oi6g4we}
\big(\widetilde{M}_\Omega(z)-\widetilde{M}_\Omega(0)\big)\varphi
=\widetilde{\gamma}_N\eta+\widetilde{M}_\Omega(0)\widetilde{\gamma}_D\eta=\psi.
\end{align}
The conclusion is that $\widetilde{M}_\Omega(z)-\widetilde{M}_\Omega(0)$ maps onto 
$\mathscr{G}_N(\partial\Omega)$, finishing the proof of \eqref{range}.

Regarding \eqref{eq:RVav-9jj}, one only needs to establish the continuity of the operator 
in question. By \eqref{abcdef}--\eqref{eq:tgfd}, \eqref{eq:RTTYH}, the fact that 
the operator $\gamma_N$ in \eqref{4.5ASdf.2} is continuous, and by the significance of  
$\widetilde{P}_{D,\Omega}(z),\widetilde{P}_{D,\Omega}(0)$ in the context of \eqref{eqn:bvp.2}
and their memberships to ${\mathcal{B}}\big(\mathscr{G}_N(\partial\Omega)^*,L^2(\Omega)\big)$, one estimates for each 
$\varphi\in\mathscr{G}_N(\partial\Omega)^*$,   
\begin{align}\label{abcdef.22}
&\big\|\big(\widetilde{M}_\Omega(z)-\widetilde{M}_\Omega(0)\big)\varphi\big\|_{\mathscr{G}_N(\partial\Omega)}
\nonumber\\[2pt]
& \quad=\big\|\widetilde{\gamma}_N\big(\widetilde{P}_{D,\Omega}(z)-\widetilde{P}_{D,\Omega}(0)\big)\varphi
\big\|_{\mathscr{G}_N(\partial\Omega)}
\nonumber\\[2pt]
& \quad\leq C\big(\big\|\big(\widetilde{P}_{D,\Omega}(z)-\widetilde{P}_{D,\Omega}(0)\big)\varphi\big\|_{L^2(\Omega)}
+\big\|\Delta\big(\widetilde{P}_{D,\Omega}(z)-\widetilde{P}_{D,\Omega}(0)\big)\varphi\big\|_{L^2(\Omega)}\big)
\nonumber\\[2pt]
& \quad\leq C\big(\big\|\widetilde{P}_{D,\Omega}(z)\varphi\big\|_{L^2(\Omega)}+\big\|\widetilde{P}_{D,\Omega}(0)\varphi
\big\|_{L^2(\Omega)} 
\nonumber\\[2pt]
&\qquad 
+\big\|(V-z)\widetilde{P}_{D,\Omega}(z)\varphi\big\|_{L^2(\Omega)}+\big\|V\widetilde{P}_{D,\Omega}(0)\varphi
\big\|_{L^2(\Omega)}\big)
\nonumber\\[2pt]
& \quad\leq C\|\varphi\|_{\mathscr{G}_N(\partial\Omega)^*},
\end{align}
for some finite constants $C$ independent of $\varphi$. This justifies the claim in \eqref{eq:RVav-9jj}.

It remains to prove the resolvent formula \eqref{resform}. For this purpose, pick $h\in L^2(\Omega)$ and consider
\begin{equation}\label{fff}
f:=(A_{D,\Omega}-zI)^{-1}h-\widetilde{P}_{D,\Omega}(z)\big(\widetilde{M}_\Omega(z)
-\widetilde{M}_\Omega(0)\big)^{-1}\big(\widetilde{P}_{D,\Omega}(\bar{z})\big)^*h.
\end{equation}
From what has been proved up to this point, and from the fact that the operator $(\widetilde{P}_{D,\Omega}(z))^*$ 
maps $L^2(\Omega)$ into $\mathscr{G}_N(\partial\Omega)$, one concludes that the function $f$ is well defined and 
belongs to $\dom(A_{max,\Omega})$. Moreover, as the solution operator $\widetilde{P}_{D,\Omega}(z)$ of the boundary 
value problem \eqref{eqn:bvp.2} maps into $\ker(A_{max,\Omega}-zI)$, one has 
\begin{equation}\label{amaxf}
(A_{max,\Omega}-zI)f=(A_{max,\Omega}-zI)(A_{D,\Omega}-zI)^{-1}h=h.
\end{equation}
At this stage we claim that $f\in\dom(A_{max,\Omega})$ in \eqref{fff} satisfies 
the boundary condition
\begin{equation}\label{bckn}
\widetilde{\gamma}_N f+\widetilde{M}_{\Omega}(0)\widetilde{\gamma}_D f=0.
\end{equation}
Indeed, from \eqref{eqn:P_D*}, Theorem~\ref{t5.5}~$(vii)$, and \eqref{iutee}, one obtains 
\begin{align}\label{nummer1}
\widetilde{\gamma}_N f &=\widetilde{\gamma}_N(A_{D,\Omega}-zI)^{-1}h
-\widetilde{\gamma}_N\widetilde{P}_{D,\Omega}(z)
\big(\widetilde{M}_\Omega(z)-\widetilde{M}_\Omega(0)\big)^{-1}\big(\widetilde{P}_{D,\Omega}(\bar{z})\big)^*h
\nonumber\\[2pt] 
&=-\big(\widetilde{P}_{D,\Omega}(\bar{z})\big)^*h+\widetilde{M}_{\Omega}(z)
\big(\widetilde{M}_\Omega(z)-\widetilde{M}_\Omega(0)\big)^{-1}\big(\widetilde{P}_{D,\Omega}(\bar{z})\big)^*h
\nonumber\\[2pt] 
&=\widetilde{M}_{\Omega}(0)\big(\widetilde{M}_\Omega(z)-\widetilde{M}_\Omega(0)\big)^{-1}
\big(\widetilde{P}_{D,\Omega}(\bar{z})\big)^*h.
\end{align}
On the other hand, since $\widetilde\gamma_D(A_{D,\Omega}-zI)^{-1}h=0$, 
relying on Theorem~\ref{t5.5}~$(vii)$ one computes  
\begin{align}\label{nummer2}
\widetilde{M}_{\Omega}(0)\widetilde{\gamma}_D f
&=-\widetilde{M}_{\Omega}(0)\widetilde{\gamma}_D\widetilde{P}_{D,\Omega}(z)
\big(\widetilde{M}_{\Omega}(z)-\widetilde{M}_{\Omega}(0)\big)^{-1}
\big(\widetilde{P}_{D,\Omega}(\bar{z})\big)^*h
\nonumber\\[2pt] 
&=-\widetilde{M}_{\Omega}(0)\big(\widetilde{M}_{\Omega}(z)-\widetilde{M}_{\Omega}(0)\big)^{-1}
\big(\widetilde{P}_{D,\Omega}(\bar{z})\big)^*h.
\end{align}
Now the claim in \eqref{bckn} is seen from \eqref{nummer1}--\eqref{nummer2}.
To proceed, from \eqref{bckn} and Theorem~\ref{t5.6} we conclude that $f\in\dom(A_{K,\Omega})$.
As such, \eqref{amaxf} gives
\begin{equation}\label{bitteschoen}
(A_{K,\Omega}-zI)f=(A_{max,\Omega}-zI)f=h,
\end{equation}
and since $z\in\rho(A_{K,\Omega})$ one finally infers from \eqref{bitteschoen} and \eqref{fff} that
\begin{align}\label{it4yh}
(A_{K,\Omega}-zI)^{-1}h=f= &\,\,(A_{D,\Omega}-zI)^{-1}h
\nonumber\\[2pt] 
&-\widetilde{P}_{D,\Omega}(z)
\big(\widetilde{M}_\Omega(z)-\widetilde{M}_\Omega(0)\big)^{-1}\big(\widetilde{P}_{D,\Omega}(z)\big)^*h.
\end{align}
This readily implies \eqref{resform}, finishing the proof of Theorem~\ref{YrarTfc}.
\end{proof}
%%%%%%%%%%

As a final result in this section we derive the Weyl spectral asymptotics of $A_{K,\Omega}$ in 
Theorem~\ref{t6.2} below. Here we follow the lines of \cite{AGMT10,AGMST10}, where the case of 
so-called quasi-convex domains was investigated. We first recall a basic result due to Kozlov 
\cite{Ko83}. Let $W_{\Omega}$ be a closed subspace in $H^2(\Omega)$ containing $\accentset{\circ}{H}^2(\Omega)$,  
\begin{equation}\label{6.1}
\accentset{\circ}{H}^2(\Omega)\subseteq W_{\Omega}\subseteq H^2(\Omega), 
\end{equation}
in particular, 
\begin{equation}\label{6.2}
W_{\Omega}\hookrightarrow L^2(\Omega)\,\text{ compactly.}     
\end{equation}  
In addition, consider the following forms in $L^2(\Omega)$:
\begin{align}\label{kko-1}
\mathfrak a_{\Omega}(f,g) &:=\sum_{0\leq|\alpha|,|\beta|\leq 2}
\int_{\Omega}a_{\alpha,\beta}(x)\ol{(\partial^\beta f)(x)}(\partial^\alpha g)(x)\,d^n x,  
\quad\dom(\mathfrak a_{\Omega})=W_{\Omega},      
\\[2pt]
\mathfrak b_{\Omega}(f,g) &:=\sum_{0\leq|\alpha|,|\beta|\leq 1}
\int_{\Omega}b_{\alpha,\beta}(x)\ol{(\partial^\beta f)(x)}(\partial^\alpha g)(x)\,d^n x,   
\quad\dom(\mathfrak b_{\Omega})=W_{\Omega}.
\label{kko-2}
\end{align} 
Suppose that they are both symmetric, that the leading coefficients of $\mathfrak a_{\Omega}$ and 
$\mathfrak b_{\Omega}$ are Lipschitz functions, while the coefficients of all lower-order 
terms are bounded, measurable functions in $\Omega$. Furthermore, assume that the following 
coercivity, nondegeneracy, and nonnegativity conditions hold for some $c\in(0,\infty)$, 
\begin{align}\label{kko-3}
&\,\mathfrak a_{\Omega}(f,f)\geq c\,\|f\|^2_{H^2(\Omega)},\quad\forall\,f\in\dom(\mathfrak a_{\Omega}),
\\[2pt]
&\sum_{|\alpha|=|\beta|=1}b_{\alpha,\beta}(x)\,\xi^{\alpha+\beta}\not=0, 
\quad\forall\,x\in\ol{\Om},\,\forall\,\xi\in{\mathbb{R}}^n\backslash\{0\},
\label{kko-4}
\\[2pt]
&\,\mathfrak b_{\Omega}(f,f)\geq 0,\quad\forall\,f\in\dom(\mathfrak b_{\Omega}). 
\label{kko-5}
\end{align}
Recall that for each multi-index $\gamma=(\gamma_1,\dots,\gamma_n)$ and each vector $\xi=(\xi_1,\dots,\xi_n)$, 
the symbol $\xi^\gamma$ stands for $\xi_1^{\gamma_1}\cdots\xi_n^{\gamma_n}$ (this is relevant in \eqref{kko-4}).

Under the above assumptions, $W_{\Omega}$ can be regarded as a Hilbert space when
equipped with the inner product $\mathfrak a_{\Omega}(\dott,\dott)$. Next, consider the operator 
$T_{\Omega}\in\cB(W_{\Omega})$, uniquely defined by the requirement that
\begin{equation}\label{kko-6}
\mathfrak a_{\Omega}(f,T_{\Omega}\,g)=\mathfrak b_{\Omega}(f,g),\quad\forall\,f,g\in W_{\Omega}.
\end{equation} 
It follows from \eqref{6.2}, \eqref{kko-3}, and \eqref{kko-5}, that the operator $T_{\Omega}$ is compact, 
nonnegative, and self-adjoint in the Hilbert space $(W_{\Omega},\mathfrak a_{\Omega}(\dott,\dott))$.
Denoting by
\begin{equation}\label{kko-7}
0\leq\cdots\leq\mu_{j+1}(T_{\Omega})\leq\mu_j(T_{\Omega})\leq\cdots\leq\mu_1(T_{\Omega}),
\end{equation} 
the eigenvalues of $T_{\Omega}$ listed according to their multiplicity, we set
\begin{equation}\label{kko-8}
\cN(\lambda,T_{\Omega}):=\#\,\big\{j\in\bbN\,\big|\,\mu_j(T_{\Omega})\geq\lambda^{-1}\big\},  
\quad\lambda>0.
\end{equation} 
The following Weyl asymptotic formula is a particular case of a general result due to 
Kozlov \cite{Ko83}. We also note that various related results can be found in \cite{Ko79}, \cite{Ko84}.

%%%%%%%%% 
\begin{theorem}\label{t6.1}
Let $\Omega\subset{\mathbb{R}}^n$, $n\geq 2$, be a bounded Lipschitz domain and retain the above 
notation and assumptions on $\mathfrak a_{\Omega}$, $\mathfrak b_{\Omega}$, and $T_{\Omega}$. 
Then the distribution function of the spectrum of $T_{\Omega}$ introduced in \eqref{kko-8} 
satisfies the following asymptotics 
\begin{equation}\label{kko-9}
\cN(\lambda,T_{\Omega})\underset{\lambda\to\infty}{=}
\omega_{\mathfrak a,\mathfrak b,\Omega}\,\lambda^{n/2}+O\big(\lambda^{(n-(1/2))/2}\big), 
\end{equation} 
where, 
\begin{equation}\label{kko-10}
\omega_{\mathfrak a,\mathfrak b,\Omega}:=\frac{1}{n(2\pi)^n}\int_{\Omega}\left(\int_{\bbS^{n-1}}
\left[\frac{\sum\limits_{|\alpha|=|\beta|=1}b_{\alpha,\beta}(x)\xi^{\alpha+\beta}}
{\sum\limits_{|\alpha|=|\beta|=2}a_{\alpha,\beta}(x)\xi^{\alpha+\beta}}
\right]^{\frac{n}{2}}d\omega_{n-1}(\xi)\right)d^n x, 
\end{equation} 
with $d\omega_{n-1}$ denoting the surface measure on the unit sphere $\bbS^{n-1}$ in $\bbR^n$. 
\end{theorem}
%%%%%%%%% 

The Weyl asymptotics for perturbed Krein Laplacians on a bounded Lipschitz domain 
now follow from Kozlov's Theorem~\ref{t6.1} in a similar way as in \cite{AGMST10,AGMST13}; cf.
\cite[Theorem 4.1]{BGMM14}.

%%%%%%%%%%
\begin{theorem}\label{t6.2} 
Assume Hypothesis~\ref{h6.1}. Let $\{\lambda_j\}_{j\in\bbN}\subset(0,\infty)$ be the strictly 
positive eigenvalues of the  Krein--von Neumann extension $A_{K,\Omega}$ enumerated in nondecreasing 
order counting multiplicity, and let
\begin{equation}\label{4455}
N(\lambda,A_{K,\Omega}):=\#\{j\in\bbN\,|\,0<\lambda_j\leq\lambda\},\quad\forall\,\lambda>0,
\end{equation}
be the eigenvalue distribution function for $A_{K,\Omega}$. 
Then the following Weyl asymptotic formula holds, 
\begin{equation}
N(\lambda,A_{K,\Omega})\underset{\lambda\to\infty}{=}(2\pi)^{-n}v_n\,|\Omega|\,\lambda^{n/2}
+O\big(\lambda^{(n-(1/2))/2}\big),
\end{equation} 
where $v_n$ denotes the volume of the unit ball in $\bbR^n$ and $|\Omega|$ is the volume of $\Omega$.
\end{theorem}
%%%%%%%%%%
\begin{proof}
Consider the densely defined symmetric forms $\mathfrak{a}_{K,\Omega}$ and 
$\mathfrak{b}_{K,\Omega}$ in $L^2(\Omega)$,  
\begin{align}
\mathfrak{a}_{K,\Omega}(f,g) & :=\big(A_{min,\Omega}f,A_{min,\Omega}g\big)_{L^2(\Omega)},
\quad\dom(\mathfrak{a}_{K,\Omega})=\accentset{\circ}{H}^2(\Omega),
\\[2pt]
\begin{split} 
\mathfrak{b}_{K,\Omega}(f,g) & :=\big(f,A_{min,\Omega}g\big)_{L^2(\Omega)},\quad  
\dom(\mathfrak{b}_{K,\Omega})=\accentset{\circ}{H}^2(\Omega).
\end{split} 
\end{align}
We note that  $\dom(\mathfrak{a}_{K,\Omega})=\dom(\mathfrak{b}_{K,\Omega})=\dom A_{min,\Omega}$ 
holds by Lemma~\ref{l4.3}. One can then verify that conditions \eqref{kko-3}--\eqref{kko-5} are 
satisfied by $\mathfrak{a}_{K,\Omega}$ and $\mathfrak{b}_{K,\Omega}$ with $W_\Omega=\accentset{\circ}{H}^2(\Omega)$. 
In this context one observes that the graph norm of $-\Delta+V$ is equivalent to the $H^2$-norm on 
$\accentset{\circ}{H}^2(\Omega)$, that is, there exists $C\in(1,\infty)$ such that  
\begin{equation}\label{6.19}
C^{-1}\|f\|^2_{H^2(\Omega)}\leq\mathfrak{a}_{K,\Omega}(f,f)\leq C\|f\|^2_{H^2(\Omega)}, 
\quad\forall\,f\in\accentset{\circ}{H}^2(\Omega),     
\end{equation}
(cf.\ the proof of Lemma~\ref{l3.3} and \eqref{eq:WH}). One observes that the self-adjoint 
operator in $L^2(\Omega)$ uniquely associated with the form $\mathfrak{a}_{K,\Omega}$ is given 
by $A_{max,\Omega} A_{min,\Omega}$ (cf. \cite[Example~VI.2.13]{Ka80}). In particular,
\begin{equation}\label{formaa}
\mathfrak{a}_{K,\Omega}(f,g)=\big(f,A_{max,\Omega}A_{min,\Omega}g\big)_{L^2(\Omega)}
\end{equation}
holds for all $f\in\dom\mathfrak{a}_{K,\Omega}$ and 
$g\in\dom(A_{max,\Omega}A_{min,\Omega})\subset\dom(\mathfrak{a}_{K,\Omega})$.

We introduce the operator $T_{K,\Omega}$ via the demand that 
\begin{equation}\label{tktktk}
\mathfrak{a}_{K,\Omega}(f,T_{K,\Omega}\,g) 
=\mathfrak{b}_{K,\Omega}(f,g),\quad\forall\,f,g\in\accentset{\circ}{H}^2(\Omega). 
\end{equation}
As discussed at the beginning of this section, $T_{K,\Omega}$ is compact, nonnegative, 
and self-adjoint on $W_{K,\Omega}$, the Hilbert space $\accentset{\circ}{H}^2(\Omega)$ 
equipped with the scalar product $\mathfrak{a}_{K,\Omega}(\dott,\dott)$. 
Moreover, one has  
\begin{equation}\label{6.20}
\lambda\in\sigma(A_{K,\Omega})\backslash\{0\}\,\text{ if and only if }\, 
\lambda^{-1}\in\sigma(T_{K,\Omega})    
\end{equation}
counting multiplicity, that is, the eigenvalues of $T_{K,\Omega}$ are precisely the reciprocals of 
the nonzero eigenvalues of $A_{K,\Omega}$, counting multiplicity. In fact, in order to verify \eqref{6.20},  
assume first that $\lambda>0$ is an eigenvalue of $A_{K,\Omega}$ corresponding to the eigenfunction 
$h\in\dom(A_{K,\Omega})$, that is, 
\begin{equation}\label{jodeldi}
A_{K,\Omega}h=\lambda h,
\end{equation}
and according to \eqref{We-Q.9jussi} the function $h$ admits a decomposition in $h=h_{min}+h_0$, where 
$h_{min}\in\dom(A_{min,\Omega})$ and $h_0\in\ker(A_{max,\Omega})$. One observes that $\lambda>0$ 
and \eqref{jodeldi} imply $h_{min}\not=0$ and $A_{min,\Omega}h_{min}=A_{K,\Omega}h$. Therefore,
\begin{equation}\label{jodeldi2}
A_{min,\Omega}h_{min}-\lambda h_{min}
=A_{K,\Omega}h-\lambda h_{min}
=\lambda h-\lambda h_{min}
=\lambda h_0
\end{equation}
belongs to $\ker(A_{max,\Omega})$ and it follows that
\begin{equation}\label{jodeldi3}
A_{max,\Omega}A_{min,\Omega}h_{min}=\lambda A_{max,\Omega} h_{min} 
=\lambda A_{min,\Omega} h_{min}.
\end{equation}
Together with \eqref{formaa} and \eqref{tktktk} this yields
\begin{equation}
\begin{split}
\mathfrak{a}_{K,\Omega}(f,\lambda^{-1}h_{min})
&=\big(f,\lambda^{-1}A_{max,\Omega}A_{min,\Omega}h_{min}\big)_{L^2(\Omega)} 
\\[2pt]
&=\big(f,A_{min,\Omega}h_{min}\big)_{L^2(\Omega)} 
\\[2pt]
&=\mathfrak{b}_{K,\Omega}(f, h_{min})
\\[2pt]
&=\mathfrak{a}_{K,\Omega}(f,T_{K,\Omega} h_{min})
\end{split}
\end{equation}
for all $f\in\dom(\mathfrak{a}_{K,\Omega})$ and hence
\begin{equation}\label{tkg}
T_{K,\Omega}h_{min}=\frac{1}{\lambda}h_{min}.
\end{equation}
Conversely, assume that $h_{min}\in\dom(T_{K,\Omega})=\accentset{\circ}{H}^2(\Omega)$ and 
$\lambda\not=0$ are such that \eqref{tkg} holds. Then
\begin{equation}
\begin{split}
\mathfrak{a}_{K,\Omega}(f,h_{min}) &=\mathfrak{a}_{K,\Omega}(f,\lambda T_{K,\Omega}h_{min})
\\[2pt]
&=\mathfrak b_{K,\Omega}(f,\lambda h_{min})
\\[2pt]
&=\big(f,\lambda A_{min,\Omega}h_{min}\big)_{L^2(\Omega)}
\end{split}
\end{equation}
for all $f\in\dom(\mathfrak{a}_{K,\Omega})$. The fact that $A_{max,\Omega}A_{min,\Omega}$ 
is the representing operator for $\mathfrak{a}_{K,\Omega}$ and the first representation 
theorem for quadratic forms \cite[Theorem~VI.2.1\,$(iii)$]{Ka80} imply 
\begin{equation}\label{jodeldi4}
h_{min}\in\dom(A_{max,\Omega}A_{min,\Omega})\,\text{ and }\,A_{max,\Omega}A_{min,\Omega}h_{min}
=\lambda A_{min,\Omega}h_{min}.
\end{equation}
Next, we consider $h:=\lambda^{-1}A_{min,\Omega}h_{min}$. It then follows from \eqref{jodeldi4} that
\begin{equation}
A_{max,\Omega}(h-h_{min})=\lambda^{-1}A_{max,\Omega}A_{min,\Omega}h_{min}-A_{min,\Omega}h_{min}=0
\end{equation}
and hence one has  
\begin{equation}
h=h_{min}+(h-h_{min}),\quad h_{min}\in\dom(A_{min,\Omega}),\,\,h-h_{min}\in\ker(A_{max,\Omega}).
\end{equation}
From \eqref{We-Q.9jussi} one concludes $h\in\dom(A_{K,\Omega})$ 
and the definition of $h$ and \eqref{jodeldi4} yield
\begin{equation}
A_{K,\Omega}h=A_{max,\Omega}h=\lambda^{-1}A_{max,\Omega}A_{min,\Omega}h_{min}=A_{min,\Omega}h_{min}=\lambda  h,
\end{equation}
that is, $h$ is an eigenfunction of $A_{K,\Omega}$ corresponding to the eigenvalue $\lambda$. 
This completes the proof of the equivalence \eqref{6.20}.

Next, introducing 
\begin{equation}\label{55tii}
\cN(\lambda,T_{K,\Omega}):=\#\big\{j\in\bbN\,\big|\,\mu_j(T_{K,\Omega})\geq\lambda^{-1}\big\}, 
\quad\forall\,\lambda>0, 
\end{equation}
where $\{\mu_j(T_{K,\Omega})\}_{j\in\bbN}$ is the ascending sequence of eigenvalues of 
$T_{K,\Omega}$ counting multiplicity, then $\cN(\lambda,T_{K,\Omega})=N(\lambda,A_{K,\Omega})$ 
for all $\lambda>0$, and Theorem~\ref{t6.1} yields the asymptotic formula, 
\begin{equation}
N(\lambda,A_{K,\Omega})=\cN(\lambda,T_{K,\Omega})\underset{\lambda\to\infty}{=} 
\omega_{K,\Omega}\,\lambda^{n/2}+O\big(\lambda^{(n-(1/2))/2}\big),
\end{equation}
with 
\begin{align}\label{trray45y}
\omega_{K,\Omega} &:=\frac{1}{n(2\pi)^n}\int_\Omega\Bigg(\int_{\bbS^{n-1}} 
\Bigg[\frac{\sum_{j=1}^n\xi_j^2}{\sum_{j,k=1}^n\xi^2_j\xi^2_k}\Bigg]^{\frac{n}{2}}\,d\omega_{n-1}(\xi)\Bigg)d^n x
\nonumber\\[2pt]
&\,\,=(2\pi)^{-n}\,v_n\,|\Omega|,
\end{align}
since the surface area of $\bbS^{n-1}$ is $nv_n$.
\end{proof}
%%%%%%%%%%%

In closing, we note that for the special case of the so-called quasi-convex domains, 
Theorem~\ref{t6.2} coincides with \cite[Theorem~8.2]{AGMT10}.

%%%%%%%%%%%%%%%%%%%%%%%%%%%%%%
%%%%%%%%%%%%%%%%%%%%%%%%%%%%%%
\section{A Description of All Self-Adjoint Extensions and Krein-Type Resolvent Formulas for Schr\"odinger Operators on Bounded Lipschitz Domains} 
\label{s10}
%%%%%%%%%%%%%%%%%%%%%%%%%%%%%%
%%%%%%%%%%%%%%%%%%%%%%%%%%%%%%

In this section we describe all self-adjoint realizations of the Schr\"{o}dinger 
differential expression $-\Delta+V$ on a bounded Lipschitz domain via explicit boundary 
conditions, and we express their resolvents in a Krein-type resolvent formula. 
Throughout this section it is assumed that Hypothesis~\ref{h4.2} holds.

First of all, we fix some real point $\mu$ which is not in the spectrum of the Dirichlet realization 
$A_{D,\Omega}$, that is, $\mu\in\rho(A_{D,\Omega})\cap\mathbb R$ and remark that such a point $\mu$ 
exists since $A_{D,\Omega}$ is semibounded from below. Moreover, by  \eqref{eqn:decAmax} one obtains the decomposition
\begin{equation}\label{fdeco0}
\begin{split}
\dom(A_{max,\Omega})
&=\dom(A_{D,\Omega})\,\dot{+}\,\ker(A_{max,\Omega}-\mu)
\\[2pt]
&=\dom(A_{D,\Omega})\,\dot{+}\,\big\{f_\mu\in\dom(A_{max,\Omega})\,\big|\,-\Delta f_\mu+ Vf_\mu=\mu f_\mu\big\}, 
\end{split}
\end{equation}
to be used in the following. We agree to decompose functions $f$ in the domain 
of $A_{max,\Omega}$ accordingly, that is, for $f\in\dom(A_{max,\Omega})$ we write
\begin{equation}\label{fdeco}
f=f_D+f_\mu,\quad f\in\dom(A_{D,\Omega}),\quad f_\mu\in\ker(A_{max,\Omega}-\mu). 
\end{equation}

In the following we will make use of the extended Dirichlet trace 
\begin{equation}\label{k7YTG}
\widetilde\gamma_D:\dom(A_{max,\Omega})\rightarrow\mathscr{G}_N(\partial\Omega)^*
\end{equation}
in Theorem~\ref{t5.5}, where $\mathscr{G}_N(\partial\Omega)^*$ is the dual of the space 
\begin{equation}\label{jagff-WACO}
\mathscr{G}_N(\partial\Omega)=\ran(\gamma_N|_{\dom(A_{D,\Omega})})
\end{equation}
introduced in Definition~\ref{Tggg.55}. Since
\begin{equation}
\mathscr{G}_N(\partial\Omega)\hookrightarrow L^2(\partial\Omega)\hookrightarrow\mathscr{G}_N(\partial\Omega)^*
\end{equation}
forms a Gelfand triple (see, e.g., \cite{Wl87}) there exist two isometric isomorphisms 
$\iota_{+}:\mathscr{G}_N(\partial\Omega)\rightarrow L^2(\partial\Omega)$ and 
$\iota_{-}:\mathscr{G}_N(\partial\Omega)^*\rightarrow L^2(\partial\Omega)$ such that
\begin{equation}
(\iota_{+}\varphi,\iota_{-}\psi)_{L^2(\partial\Omega)}
={}_{\mathscr{G}_N(\partial\Omega)}\big\langle\varphi,\psi\big\rangle_{\mathscr{G}_N(\partial\Omega)^*}
\end{equation}
holds for all $\varphi\in\mathscr{G}_N(\partial\Omega)$ and $\psi\in\mathscr{G}_N(\partial\Omega)^*$. 
For a closed subspace $\mathscr{X}\subset\mathscr{G}_N(\partial\Omega)^*$ set 
\begin{equation}\label{raded-TEXAS}
\mathscr{X}^*:=\iota_{+}^{-1}\iota_{-}(\mathscr{X})\subset\mathscr{G}_N(\partial\Omega),
\end{equation}
so that 
\begin{equation}\label{jussixx}
\iota_{+}(\mathscr{X}^*)=\iota_{-}(\mathscr{X})\subset L^2(\partial\Omega).
\end{equation}
If $Q_{\iota_{+}(\mathscr{X}^*)}$ denotes the orthogonal projection in $L^2(\partial\Omega)$ 
onto the closed subspace $\iota_{+}(\mathscr{X}^*)\subset L^2(\partial\Omega)$ then we say that 
\begin{equation}\label{jussipq}
P_{\mathscr{X}^*}:=\iota_{+}^{-1}Q_{\iota_{+}(\mathscr{X}^*)}\iota_{+}
\end{equation}
is the orthogonal projection in $\mathscr{G}_N(\partial\Omega)$ onto the closed subspace 
$\mathscr{X}^*\subset\mathscr{G}_N(\partial\Omega)$. We note that for all $\varphi\in\mathscr{G}_N(\partial\Omega)$ 
and all $\psi\in\mathscr{X}$ one has $P_{\mathscr{X}^*}\varphi\in\mathscr{X}^*$ and
\begin{equation}\label{jussipxx}
\begin{split}
{}_{\mathscr{X}^*}\big\langle P_{\mathscr{X}^*}\varphi,\psi\big\rangle_{\mathscr{X}}
&={}_{\mathscr{G}_N(\partial\Omega)}\big\langle P_{\mathscr{X}^*}\varphi,\psi\big\rangle_{\mathscr{G}_N(\partial\Omega)^*} 
=\big(\iota_{+}P_{\mathscr{X}^*}\varphi,\iota_{-}\psi\big)_{L^2(\partial\Omega)}
\\[2pt]
&=\big(Q_{\iota_{+}(\mathscr{X}^*)}\iota_{+}\varphi,\iota_{-}\psi\big)_{L^2(\partial\Omega)}
=\big(\iota_+\varphi,Q_{\iota_+(\mathscr X^*)}\iota_-\psi\big)_{L^2(\partial\Omega)}
\\[2pt]
&=\big(\iota_{+}\varphi,\iota_{-}\psi\big)_{L^2(\partial\Omega)}
={}_{\mathscr{G}_N(\partial\Omega)}\big\langle\varphi,\psi\big\rangle_{\mathscr{G}_N(\partial\Omega)^*},
\end{split}
\end{equation}
where \eqref{jussipq} and $\iota_{-}\psi\in\iota_{-}(\mathscr{X})=\iota_{+}(\mathscr{X}^*)$ were used. 
We denote by $(\mathscr{X}^*)^\bot$ the corresponding orthogonal complement of $\mathscr{X}^*$, 
that is, $(\mathscr{X}^*)^\bot=\iota_+^{-1}(\iota_+(\mathscr{X}^*)^{\bot_{L^2}})$,
and the corresponding orthogonal projection in $\mathscr{G}_N(\partial\Omega)$ 
is denoted by $P_{(\mathscr{X}^*)^\bot}$. In  the same style we write $P_{\mathscr{X}}$ 
and $P_{\mathscr{X}^\bot}$ for the orthogonal projections onto $\mathscr{X}$ and 
$\mathscr{X}^\bot$, respectively. The canonical embedding of $\mathscr{X}$ into 
$\mathscr{G}_N(\partial\Omega)^\ast$ will be denoted by $\iota_{\mathscr{X}}$.

Let again $\mathscr{X}\subset\mathscr{G}_N(\partial\Omega)^\ast$ be a closed subspace and let 
$\mathscr{X}^*=\iota_{+}^{-1}\iota_{-}(\mathscr{X})\subset\mathscr{G}_N(\partial\Omega)$.
We shall say that a densely defined operator $T:\mathscr{X}\supset\dom(T)\rightarrow\mathscr{X}^*$ 
is symmetric if 
\begin{equation}\label{jussitsym}
{}_{\mathscr{X}^*}\langle T\varphi,\psi\rangle_{\mathscr{X}}
={}_{\mathscr{X}}\langle\varphi,T\psi\rangle_{\mathscr{X}^*}\,\text{ for all }\,\varphi,\psi\in\dom(T)
\end{equation}
and $T:\mathscr{X}\supset\dom(T)\rightarrow\mathscr{X}^*$ is said to be self-adjoint if 
\begin{align}\label{jussitsa}
\begin{split}
& {}_{\mathscr{X}^*}\langle T\varphi,\psi\rangle_{\mathscr{X}}
={}_{\mathscr{X}}\langle\varphi,\psi^\prime\rangle_{\mathscr{X}^*}
\,\text{ for all }\,\varphi\in\dom(T)
\\[2pt] 
& \quad\text{ implies }\,\psi\in\dom(T)\,\text{ and }\,T\psi=\psi^\prime.
\end{split}
\end{align}
We note that $T:\mathscr{X}\supset\dom(T)\rightarrow\mathscr{X}^*$ is symmetric (self-adjoint) in 
this sense if and only if the operator
\begin{equation}\label{yta5r-WACO}
\iota_{+}T\iota_{-}^{-1},\quad\dom(\iota_{+}T\iota_{-}^{-1})
=\iota_{-}(\dom(T))\subset\iota_{-}(\mathscr{X})
\end{equation}
is symmetric (self-adjoint, respectively) in the Hilbert space 
$\iota_{-}(\mathscr{X})=\iota_{+}(\mathscr{X}^*)\subset L^2(\partial\Omega)$.

In the following theorem all self-adjoint realizations of $-\Delta+V$ are characterized via explicit 
boundary conditions in terms of closed subspaces $\mathscr{X}\subset\mathscr{G}_N(\partial\Omega)^*$ and self-adjoint operators $T$. In this context we note that the first description of all self-adjoint realizations of second-order proper elliptic operators with smooth coefficients on smooth domains in terms of boundary conditions was obtained by Vi{\u s}ik in his celebrated 1952 memoir, see \cite[Section~6]{Vi63}. The result below, is along the lines of the classical parametrization 
due to Grubb in \cite{Gr68}, is given here a complete and self-contained proof. For earlier work, 
see also \cite[Corollary~4.4]{BM13}, \cite[Theorem~14.3]{GM11}, and \cite[Propositions~3.5, 3.6]{Ma10}. 

%%%%%%%%
\begin{theorem}\label{jussisaext}
Assume Hypothesis~\ref{h4.2}, let $\widetilde\gamma_D$ be the extension of 
the Dirichlet trace operator onto $\dom(A_{max,\Omega})$, fix some point 
$\mu\in\rho(A_{D,\Omega})\cap\mathbb R$ and decompose $f\in\dom(A_{max,\Omega})$ in the form \eqref{fdeco}. 

Then there is a one-to-one correspondence between the self-adjoint extensions of $A_{min,\Omega}$ in $L^2(\Omega)$
and the family of pairs $\{\mathscr{X},T\}$, consisting of a closed subspace $\mathscr{X}$ of $\mathscr{G}_N(\partial\Omega)^*$ 
and a self-adjoint operator $T:\mathscr{X}\supset\dom(T)\rightarrow\mathscr{X}^*$ as follows: For every closed subspace 
$\mathscr{X}\subset\mathscr{G}_N(\partial\Omega)^*$ and every self-adjoint operator 
$T:\mathscr{X}\supset\dom(T)\rightarrow\mathscr{X}^*$ the operator
\begin{equation}\label{atjussi}
\begin{split}
& A_{T,\Omega}=-\Delta+V,
\\[2pt]
& \dom(A_{T,\Omega})=\big\{f\in\dom(A_{max,\Omega})\,\big|\,T\,\widetilde\gamma_D f=P_{\mathscr{X}^*}\gamma_N f_D\big\}
\end{split}
\end{equation}
is a self-adjoint extension of $A_{min,\Omega}$ in $L^2(\Omega)$, where $P_{\mathscr{X}^*}$ denotes 
the orthogonal projection in $\mathscr{G}_N(\partial\Omega)$ onto $\mathscr{X}^*$. Conversely, for every 
self-adjoint extension $A$ of $A_{min,\Omega}$ in $L^2(\Omega)$ there exists a closed subspace 
$\mathscr{X}\subset\mathscr{G}_N(\partial\Omega)^*$ and a self-adjoint operator 
$T:\mathscr{X}\supset\dom(T)\rightarrow\mathscr{X}^*$ such that $A=A_{T,\Omega}$, that is,
\begin{equation}\label{gafac.WACO}
\begin{split}
& A=-\Delta+V,
\\[2pt]
& \dom(A)=\big\{f\in\dom(A_{max,\Omega})\,\big|\,T\,\widetilde\gamma_D f=P_{\mathscr{X}^*}\gamma_N f_D\big\}.
\end{split}
\end{equation}
\end{theorem}
%%%%%%%%%
\begin{proof}
Let $f,g\in\dom(A_{max,\Omega})$ and decompose $f,g$ in the form 
\begin{equation}\label{decojussi}
f=f_D+f_\mu\,\text{ and }\,g=g_D+g_\mu
\end{equation}
as in \eqref{fdeco0}--\eqref{fdeco}. It then follows from the self-adjointness of $A_D$, the properties 
$A_{max,\Omega}f_\mu=\mu f_\mu$ and $A_{max,\Omega}g_\mu=\mu g_\mu$, and the extended Green's formula 
\eqref{eq:33VaV} that
\begin{align}\label{greenjussi}
&(A_{max,\Omega}f,g)_{L^2(\Omega)}-(f,A_{max,\Omega}g)_{L^2(\Omega)}
\nonumber\\[2pt]
&\quad=(A_{D,\Omega}f_D+A_{max,\Omega}f_\mu,g_D+g_\mu)_{L^2(\Omega)}
-(f_D+f_\mu,A_{D,\Omega}g_D+A_{max,\Omega} g_\mu)_{L^2(\Omega)}
\nonumber\\[2pt]
&\quad=(A_{D,\Omega}f_D,g_\mu)_{L^2(\Omega)}-(f_D,A_{max,\Omega}g_\mu)_{L^2(\Omega)} 
\nonumber\\[2pt]
&\qquad+(A_{max,\Omega}f_\mu,g_D)_{L^2(\Omega)}-( f_\mu,A_{D,\Omega}g_D)_{L^2(\Omega)}
\nonumber\\[2pt]
&\quad={}_{\mathscr{G}_N(\partial\Omega)}\big\langle\gamma_N f_D,\widetilde\gamma_D g_\mu
\big\rangle_{\mathscr{G}_N(\partial\Omega)^*} 
-{}_{\mathscr{G}_N(\partial\Omega)^*}\big\langle\widetilde\gamma_D f_\mu,\gamma_N g_D\big\rangle_{\mathscr{G}_N(\partial\Omega)}
\nonumber\\[2pt]
&\quad={}_{\mathscr{G}_N(\partial\Omega)}\big\langle\gamma_N f_D,\widetilde\gamma_D g\big\rangle_{\mathscr{G}_N(\partial\Omega)^*} 
-{}_{\mathscr{G}_N(\partial\Omega)^*}\big\langle\widetilde\gamma_D f,\gamma_N g_D\big\rangle_{\mathscr{G}_N(\partial\Omega)},
\end{align}
where $\dom(A_{D,\Omega})=\ker\widetilde\gamma_D$ was used in the last step (cf.~\eqref{4.30}). \\[1mm] 

Next, assume that $\mathscr{X}\subset\mathscr{G}_N(\partial\Omega)^*$ is a closed subspace 
of $\mathscr{G}_N(\partial\Omega)^*$ and let $T$ be a self-adjoint operator which is defined on 
the dense subspace $\dom(T)\subset\mathscr{X}$ and maps into $\mathscr{X}^*$ (cf.~\eqref{jussitsa}). 
We consider the operator $A_{T,\Omega}=-\Delta+V$ defined on the linear subspace 
\begin{equation}\label{domatjussi}
\dom(A_{T,\Omega})=\big\{f\in\dom(A_{max,\Omega})\,\big|\,T\,\widetilde\gamma_D f=P_{\mathscr{X}^*}\gamma_N f_D\big\}.
\end{equation}
As $\dom(A_{min,\Omega})$ is contained in $\ker\widetilde\gamma_D\cap\ker\gamma_N$, it follows that  
\begin{equation}\label{gar-MM.TEXAS}
\dom(A_{min,\Omega})\subset\dom(A_{T,\Omega})
\end{equation}
and the inclusion $\dom(A_{T,\Omega})\subset\dom A_{max,\Omega}$ is clear from \eqref{domatjussi}. Hence, 
\begin{equation}\label{incl}
\dom(A_{min,\Omega})\subset\dom(A_{T,\Omega})\subset\dom(A_{max,\Omega}),
\end{equation}
and therefore, the operator $A_{T,\Omega}$ is an extension of $A_{min,\Omega}$, and a restriction of 
$A_{max,\Omega}$. Next we verify that the operator $A_{T,\Omega}$ is symmetric in $L^2(\Omega)$. 
For this purpose, let $f,g\in\dom(A_{T,\Omega})$. By \eqref{incl} the functions $f,g$ belong to 
$\dom(A_{max,\Omega})$ and hence they can be decomposed as in \eqref{decojussi}. Then one has 
\begin{equation}\label{einmal}
\widetilde\gamma_D f\in\dom(T)\subset\mathscr{X},\quad 
T\,\widetilde\gamma_D f=P_{\mathscr{X}^*}\gamma_N f_D\subset\mathscr{X}^*,
\end{equation}
and 
\begin{equation}\label{zweimal}
\widetilde\gamma_D g\in\dom(T)\subset\mathscr{X},\quad 
T\,\widetilde\gamma_D g=P_{\mathscr{X}^*}\gamma_N g_D\subset\mathscr{X}^*.
\end{equation}
Thus, one concludes from \eqref{greenjussi} together with \eqref{einmal}, 
\eqref{zweimal}, and \eqref{jussipxx} that
\begin{align}\label{gafnn-WACO}
&(A_{T,\Omega}f,g)_{L^2(\Omega)}-(f,A_{T,\Omega}g)_{L^2(\Omega)}
\nonumber\\[2pt]
&\quad=(A_{max,\Omega}f,g)_{L^2(\Omega)}-(f,A_{max,\Omega}g)_{L^2(\Omega)}
\nonumber\\[2pt]
&\quad={}_{\mathscr{G}_N(\partial\Omega)}
\big\langle\gamma_N f_D,\widetilde\gamma_D g\big\rangle_{\mathscr{G}_N(\partial\Omega)^*}
-{}_{\mathscr{G}_N(\partial\Omega)^*}\big\langle\widetilde\gamma_D f,\gamma_N g_D\big\rangle_{\mathscr{G}_N(\partial\Omega)}
\nonumber\\[2pt]
&\quad={}_{\mathscr{X}^*}\big\langle P_{\mathscr{X}^*}\gamma_N f_D,\widetilde\gamma_D g\big\rangle_{\mathscr{X}} 
-{}_{\mathscr{X}}\big\langle\widetilde\gamma_D f,P_{\mathscr{X}^*}\gamma_N g_D\big\rangle_{\mathscr{X}^*}
\nonumber\\[2pt]
&\quad={}_{\mathscr{X}^*}\big\langle T\,\widetilde\gamma_D f,\widetilde\gamma_D g\big\rangle_{\mathscr{X}} 
-{}_{\mathscr{X}}\big\langle\widetilde\gamma_D f,T\,\widetilde\gamma_D g\big\rangle_{\mathscr{X}^*}=0,
\end{align}
using that $T$ is symmetric in the last step (cf.~\eqref{jussitsym}). 
This proves that the operator $A_{T,\Omega}$ is symmetric in $L^2(\Omega)$.

Next, it will be verified that the inclusion $\dom(A_{T,\Omega}^*)\subset\dom(A_{T,\Omega})$ holds. 
To accomplish this goal, pick some $g\in\dom(A_{T,\Omega}^*)$. We will then show that
\begin{equation}\label{toshow}
\widetilde\gamma_D g\in\dom(T)\,\text{ and }\,T\,\widetilde\gamma_D g=P_{\mathscr{X}^*}\gamma_N g_D.
\end{equation}
In fact, note first that the mapping 
\begin{equation}\label{jussibm}
\dom(A_{max,\Omega})\ni f=f_D+f_\mu\mapsto\big\{\widetilde\gamma_D f,\gamma_N f_D\big\} 
\in\mathscr{G}_N(\partial\Omega)^*\times\mathscr{G}_N(\partial\Omega)
\end{equation}
is surjective; this is an immediate consequence of Theorem~\ref{t5.5}\,$(i)$, \eqref{4.30}, and 
Definition~\ref{Tggg.55}. Next, we check that $\widetilde\gamma_D g\in\mathscr{X}$; in fact, we will show that
\begin{equation}\label{goodtohave}
P_{\mathscr{X}^\bot}\widetilde\gamma_D g=0,
\end{equation}
where $P_{\mathscr{X}^\bot}$ is the orthogonal projection in $\mathscr{G}_N(\partial\Omega)^*$ 
onto $\mathscr{X}^\bot$. For $\varphi\in(\mathscr{X}^*)^\bot$ choose $f\in\dom(A_{\max,\Omega})$ 
such that $\widetilde\gamma_D f=0$ and $\gamma_N f_D=\varphi$; this is possible since \eqref{jussibm} 
is surjective. In that case one has $P_{\mathscr{X}^*}\gamma_N f_D=P_{\mathscr{X}^*}\varphi=0$
and hence $f\in\dom(A_{T,\Omega})$ by \eqref{domatjussi}. It now follows from $g\in\dom(A_{T,\Omega}^*)$ 
and \eqref{greenjussi} that
\begin{equation}\label{yarf-WACO}
\begin{split}
0 &=(A_{T,\Omega}f,g)_{L^2(\Omega)}-(f,A_{T,\Omega}^*g)_{L^2(\Omega)}
\\[2pt]
&=(A_{max,\Omega}f,g)_{L^2(\Omega)}-(f,A_{max,\Omega}g)_{L^2(\Omega)}
\\[2pt]
&={}_{\mathscr{G}_N(\partial\Omega)}\big\langle\gamma_N f_D,\widetilde\gamma_D g\big\rangle_{\mathscr{G}_N(\partial\Omega)^*}
-{}_{\mathscr{G}_N(\partial\Omega)^*}\big\langle\widetilde\gamma_D f,\gamma_N g_D\big\rangle_{\mathscr{G}_N(\partial\Omega)}
\\[2pt]
&={}_{(\mathscr{X}^*)^\bot}\big\langle\varphi,P_{\mathscr{X}^\bot}\widetilde\gamma_D g\big\rangle_{\mathscr{X}^\bot}.
\end{split}
\end{equation}
Since this identity holds for all $\varphi\in(\mathscr{X}^*)^\bot$ one concludes \eqref{goodtohave}, hence
\begin{equation}\label{goodtohave2}
\widetilde\gamma_D g\in\mathscr{X}.
\end{equation}
In the sequel we again make use of the surjectivity of the map \eqref{jussibm}. In particular, the space 
$\mathscr{X}\times\mathscr{X}^*$ as a subspace of $\mathscr{G}_N(\partial\Omega)^*\times\mathscr{G}_N(\partial\Omega)$ 
is contained in the range of the map in \eqref{jussibm}. Hence for $\varphi\in\dom(T)\subset\mathscr{X}$ 
there exists $f\in\dom(A_{max,\Omega})$ such that
\begin{equation}\label{varphijussi}
\varphi=\widetilde\gamma_D f\in\dom(T)\subset\mathscr{X}\,\text{ and }\, 
T\varphi=\gamma_N f_D=P_{\mathscr{X}^*}\gamma_N f_D\subset\mathscr{X}^*,
\end{equation}
and from \eqref{domatjussi} one concludes that $f\in\dom(A_{T,\Omega})$.
Making use of \eqref{varphijussi}, \eqref{greenjussi}, 
\begin{equation}\label{yarffv-gr4.TEXAS}
f\in\dom(A_{T,\Omega})\,\text{ and }\,g\in\dom(A_{T,\Omega}^*),
\end{equation}
one computes together with \eqref{goodtohave2}, 
\begin{align}\label{jaffnbv-65f.WACO}
&{}_{\mathscr{X}^*}\big\langle T\varphi,\widetilde\gamma_D g\big\rangle_{\mathscr{X}}
={}_{\mathscr{X}^*}\big\langle P_{\mathscr{X}^*}\gamma_N f_D,\widetilde\gamma_D g\big\rangle_{\mathscr{X}} 
={}_{\mathscr{G}_N(\partial\Omega)}\big\langle\gamma_N f_D,\widetilde\gamma_D g\big\rangle_{\mathscr{G}_N(\partial\Omega)^*}
\nonumber\\[2pt]
&\quad=(A_{max,\Omega}f,g)_{L^2(\Omega)}-(f,A_{max,\Omega}g)_{L^2(\Omega)}
+{}_{\mathscr{G}_N(\partial\Omega)^*}\big\langle\widetilde\gamma_D f,\gamma_N g_D\big\rangle_{\mathscr{G}_N(\partial\Omega)}
\nonumber\\[2pt]
&\quad=(A_{T,\Omega}f,g)_{L^2(\Omega)}-(f,A_{T,\Omega}^*g)_{L^2(\Omega)}
+{}_{\mathscr{X}}\big\langle\varphi,P_{\mathscr{X}^*}\gamma_N g_D\big\rangle_{\mathscr{X}^*}
\nonumber\\[2pt]
&\quad={}_{\mathscr{X}}\big\langle\varphi,P_{\mathscr{X}^*}\gamma_N g_D\big\rangle_{\mathscr{X}^*}.
\end{align}
This relation holds for all $\varphi\in\dom(T)$ and as $T$ is assumed to be self-adjoint 
(cf.~\eqref{jussitsa}) this implies $\widetilde\gamma_D g\in\dom(T)$ and 
$T\,\widetilde\gamma_D g=P_{\mathscr{X}^*}\gamma_N g_D$, that is, \eqref{toshow} holds.
But then \eqref{domatjussi} immediately implies $g\in\dom(A_{T,\Omega})$. This establishes  
the inclusion $\dom(A_{T,\Omega}^*)\subset\dom(A_{T,\Omega})$. All together, it follows that 
for a self-adjoint operator $T:\mathscr{X}\supset\dom(T)\rightarrow\mathscr{X}^*$ mapping 
from some closed subspace $\mathscr{X}\subset\mathscr{G}_N(\partial\Omega)^*$ into 
$\mathscr{X}^*\subset\mathscr{G}_N(\partial\Omega)$ the operator $A_{T,\Omega}$ in 
\eqref{atjussi} is self-adjoint in $L^2(\Omega)$. \\[1mm] 
Next, we prove the converse statement. Suppose in this context that $A=A^*$ is some 
self-adjoint extension of $A_{min,\Omega}$ in $L^2(\Omega)$, that is,
\begin{equation}\label{yafcwq-WACO}
A_{min,\Omega}\subset A=A^*\subset A_{max,\Omega}.
\end{equation}
In particular, $A$ acts as $-\Delta +V$ on $\dom(A)\subset\dom(A_{max,\Omega})$.
We now define a closed subspace $\mathscr{X}$ in $\mathscr{G}_N(\partial\Omega)^*$ by
\begin{equation}\label{xjussi}
\mathscr{X}:=\overline{\big\{\varphi\in\mathscr{G}_N(\partial\Omega)^*\,\big|\, 
\varphi=\widetilde\gamma_D f\,\text{ for some }\,f\in\dom A\big\}}.
\end{equation}
At this point we introduce the linear operator $T$ mapping from 
$\mathscr{X}\supset\dom(T)\rightarrow\mathscr{X}^*$ by
\begin{equation}\label{tjussi}
\begin{split}
& T\,\widetilde\gamma_D f:=P_{\mathscr{X}^*}\gamma_N f_D,
\\[2pt]
& \dom(T)=\big\{\varphi\in\mathscr{X}\,\big|\,\varphi=\widetilde\gamma_D f\,\text{ for some }\,f\in\dom(A)\big\}.
\end{split}
\end{equation}
One observes that $T$ is a well defined linear operator. In fact, if for some function $f\in\text{dom}(A)$ 
one has $\widetilde\gamma_D f=0$, then for every $g\in\dom(A)$ one may write 
\begin{align}\label{98yyg-WACO}
&{}_{\mathscr{X}^*}\big\langle P_{\mathscr{X}^*}\gamma_N f_D,\widetilde\gamma_D g\big\rangle_{\mathscr{X}}
\nonumber\\[2pt]
&\quad={}_{\mathscr{G}_N(\partial\Omega)}\big\langle\gamma_N f_D,\widetilde\gamma_D g\big\rangle_{\mathscr{G}_N(\partial\Omega)^*}   
-{}_{\mathscr{G}_N(\partial\Omega)^*}\big\langle\widetilde\gamma_D f,\gamma_N g_D\big\rangle_{\mathscr{G}_N(\partial\Omega)}
\nonumber\\[2pt]
&\quad=(A_{max,\Omega}f,g)_{L^2(\Omega)}-(f,A_{max,\Omega}g)_{L^2(\Omega)}
\nonumber\\[2pt]
&\quad=(Af,g)_{L^2(\Omega)}-(f,Ag)_{L^2(\Omega)}=0,
\end{align}
where also \eqref{greenjussi} and the symmetry of $A$ was used. By the definition of the space 
$\mathscr{X}$ in \eqref{xjussi}, the elements $\widetilde\gamma_D g$ with $g\in\dom(A)$ form a dense 
set in $\mathscr{X}$. This implies $P_{\mathscr{X}^*}\gamma_N f_D=0$ and hence the operator $T$ 
in \eqref{tjussi} is well defined.

Next, it will be shown that $T:\mathscr{X}\supset\dom(T)\rightarrow\mathscr{X}^*$ is self-adjoint 
in the sense of \eqref{jussitsa}. Assume in this context that $\psi\in\mathscr{X}$ and $\psi^\prime\in\mathscr{X}^*$
are such that
\begin{equation}\label{jussipsipsi}
{}_{\mathscr{X}^*}\langle T\varphi,\psi\rangle_{\mathscr{X}}
={}_{\mathscr{X}}\langle\varphi,\psi^\prime\rangle_{\mathscr{X}^*}
\end{equation}
holds for all $\varphi\in\dom(T)$. Next, choose $g\in\dom(A_{max,\Omega})$ such that 
\begin{equation}\label{jaf-jvc-WACO}
\psi=\widetilde\gamma_D g\,\text{ and }\,\psi^\prime=\gamma_N g_D=P_{\mathscr{X}^*}\gamma_N g_D,
\end{equation}
which is possible due to the surjectivity of the map \eqref{jussibm}. Clearly, for $\varphi\in\dom(T)$ 
there exists $f\in\dom(A)$ such that $\varphi=\widetilde\gamma_D f$, hence $T\varphi=P_{\mathscr{X}^*}\gamma_N f_D$. 
Then one concludes from \eqref{jussipsipsi} that 
\begin{equation}\label{8yfda99.WACO}
\begin{split}
0 &={}_{\mathscr{X}^*}\langle T\varphi,\psi\rangle_{\mathscr{X}}
-{}_{\mathscr{X}}\langle\varphi,\psi^\prime\rangle_{\mathscr{X}^*}
\\[2pt]
&={}_{\mathscr{X}^*}\langle P_{\mathscr{X}^*}\gamma_N f_D,\widetilde\gamma_D g\rangle_{\mathscr{X}}
-{}_{\mathscr{X}}\langle\widetilde\gamma_D f,\gamma_N g_D\rangle_{\mathscr{X}^*}
\\[2pt]
&={}_{\mathscr{G}_N(\partial\Omega)}
\big\langle\gamma_N f_D,\widetilde\gamma_D g\big\rangle_{\mathscr{G}_N(\partial\Omega)^*}
-{}_{\mathscr{G}_N(\partial\Omega)^*}\big\langle\widetilde\gamma_D f,\gamma_N g_D\big\rangle_{\mathscr{G}_N(\partial\Omega)}
\\[2pt]
&=(A_{max,\Omega}f,g)_{L^2(\Omega)}-(f,A_{max,\Omega}g)_{L^2(\Omega)}
\\[2pt]
&=(Af,g)_{L^2(\Omega)}-(f,A_{max,\Omega}g)_{L^2(\Omega)}.
\end{split}
\end{equation}
The above equality holds for all $\varphi=\widetilde\gamma_D f\in\dom(T)$ or, equivalently, for all 
$f\in\dom(A)$. As $A$ is assumed to be self-adjoint in $L^2(\Omega)$ one infers that $g\in\dom(A)$ and 
$Ag=A_{max,\Omega}g$. In particular, 
\begin{equation}\label{ut4ed.WACO}
\psi=\widetilde\gamma_D g\in\dom(T)\,\text{ and }\,  
T\psi=T\,\widetilde\gamma_D g=P_{\mathscr{X}^*}\gamma_N g_D=\psi^\prime.
\end{equation}
Therefore, by \eqref{jussipsipsi} and \eqref{jussitsa} the operator 
$T:\mathscr{X}\supset\dom(T)\rightarrow\mathscr{X}^*$ is self-adjoint. 
This completes the proof of Theorem~\ref{jussisaext}.
\end{proof}
%%%%%%%%%

Given Theorem \ref{jussisaext}, one can now attempt a spectral analysis of self-adjoint extensions other than those 
discussed in this monograph. Interesting candidates can be found, for instance, in \cite{Ag13}, \cite[Chs.~11,\,12]{Ag15}.  

\smallskip 

It is worth noting that for $\mathscr{X}:=\mathscr{G}_N(\partial\Omega)^*$ and $T:=0$ the self-adjoint 
realization in \eqref{atjussi} coincides with the Krein--von Neumann extension $A_{K,\Omega}$. 
From this point of view, the following theorem may be viewed as a generalization of Theorem~\ref{YrarTfc}, 
where the resolvents of $A_{K,\Omega}$ and $A_{D,\Omega}$ have been related via a Krein-type resolvent formula. 
In fact, setting $\mathscr{X}:=\mathscr{G}_N(\partial\Omega)^*$, $T:=0$, and choosing $\mu:=0$,  
the resolvent formula in the next theorem reduces to the one in Theorem~\ref{YrarTfc}. 
Let us now turn to the general situation. 

%%%%%%%%
\begin{theorem}\label{taff-WACO}
Assume Hypothesis~\ref{h4.2} and let $\widetilde\gamma_D$ be the extension of the Dirichlet 
trace operator onto $\dom(A_{max,\Omega})$. Let $\mathscr{X}\subset\mathscr{G}_N(\partial\Omega)^*$ 
be a closed subspace, let $T:\mathscr{X}\supset\dom(T)\rightarrow\mathscr{X}^*$ be a self-adjoint 
operator and let 
\begin{equation}\label{atjussi2}
\begin{split}
& A_{T,\Omega}=-\Delta+V,
\\[2pt]
& \dom(A_{T,\Omega})=\big\{f\in\dom A_{max,\Omega}\,\big|\,T\,\widetilde\gamma_D f=P_{\mathscr{X}^*}\gamma_N f_D\big\}
\end{split}
\end{equation}
be the corresponding self-adjoint realization of $-\Delta+V$ in $L^2(\Omega)$ in \eqref{atjussi}. Then the operator
\begin{equation}\label{ytfddd-WACO}
T+P_{\mathscr{X}^*}\big(\widetilde M_\Omega(z)-\widetilde M_\Omega(\mu)\big)
\iota_{\mathscr{X}}:\mathscr{X}\supset\dom(T)\rightarrow\mathscr{X}^*
\end{equation}
is bijective and with inverse in $\mathcal{B}(\mathscr{X}^*,\mathscr{X})$ whenever 
$z\in\rho(A_{T,\Omega})\cap\rho(A_{D,\Omega})$, and the following Krein-type resolvent 
formula holds in ${\mathcal{B}}(L^2(\Omega))$:
\begin{align}\label{6543ww.TEXAS}
\begin{split}
&(A_{T,\Omega}-zI)^{-1}-(A_{D,\Omega}-zI)^{-1}
\\[2pt]
&\quad=-\widetilde P_{D,\Omega}(z)\iota_{\mathscr{X}}\big(T+P_{\mathscr{X}^*}
\big(\widetilde M_\Omega(z)-\widetilde M_\Omega(\mu)\big)\iota_{\mathscr{X}}\big)^{-1}
P_{\mathscr{X}^*}\big(\widetilde P_{D,\Omega}(\bar z)\big)^*.
\end{split}
\end{align}
\end{theorem}
%%%%%%%%
\begin{proof}
For $z\in\rho(A_{D,\Omega})$ define the operator $H(z):\mathscr{X}\rightarrow\mathscr{X}^*$ by setting
\begin{equation}\label{6d4d6g-WACO}
H(z):=P_{\mathscr{X}^*}\big(\widetilde M_\Omega(z)-\widetilde M_\Omega(\mu)\big)\iota_{\mathscr{X}}.
\end{equation}
Note that $H(z)$ is well defined, as the range of $\widetilde M_\Omega(z)-\widetilde M_\Omega(\mu)$ 
is contained in $\mathscr{G}_N(\partial\Omega)$ (this can be verified in the same way as in the proof 
of Theorem~\ref{YrarTfc}). Furthermore, $H(z)$ is bounded (cf. the proof of Theorem~\ref{YrarTfc}). 
Let $T:\mathscr{X}\supset\dom(T)\rightarrow\mathscr{X}^*$ be a self-adjoint operator. We shall 
show that the operator
\begin{equation}\label{jussitm}
T+H(z)=T+P_{\mathscr{X}^*}\big(\widetilde M_\Omega(z)-\widetilde M_\Omega(\mu)\big)\iota_{\mathscr{X}}
:\mathscr{X}\supset\dom(T)\rightarrow\mathscr{X}^*
\end{equation}
is bijective for all $z\in\rho(A_{D,\Omega})\cap\rho(A_{T,\Omega})$. To this end, first suppose that for 
some $\varphi\in\dom(T)$ we have
\begin{equation}
\big(T+H(z)\big)\varphi=T\varphi+P_{\mathscr{X}^*}\big(\widetilde M_\Omega(z)-\widetilde M_\Omega(\mu)\big)
\iota_{\mathscr{X}}\varphi=0.
\end{equation}
There exists $f_z\in\ker(A_{max,\Omega}-zI)$ such that $\widetilde\gamma_D f_z=\varphi$. 
As $\varphi\in\mathscr{X}$, one has $\iota_{\mathscr{X}}\widetilde\gamma_D f_z=\widetilde\gamma_D f_z$.
Decompose $f_z$ as in \eqref{fdeco} in the form $f_z=f_{D,z}+f_{\mu,z}$, where 
$f_{D,z}\in\dom(A_{D,\Omega})$ and $f_{\mu,z}\in\ker(A_{max,\Omega}-\mu I)$. One then computes 
\begin{equation}\label{43ddc-WACO}
\begin{split}
T\,\widetilde\gamma_D f_z=T\varphi
&=-P_{\mathscr{X}^*}\big(\widetilde M_\Omega(z)-\widetilde M_\Omega(\mu)\big)\iota_{\mathscr{X}}\varphi
\\[2pt]
&=-P_{\mathscr{X}^*}\big(\widetilde M_\Omega(z)-\widetilde M_\Omega(\mu)\big)\widetilde\gamma_D f_z
\\[2pt]
&=-P_{\mathscr{X}^*}\big(\widetilde M_\Omega(z)\widetilde\gamma_D f_z
-\widetilde M_\Omega(\mu)\widetilde\gamma_D(f_{D,z}+f_{\mu,z})\big)
\\[2pt]
&=P_{\mathscr{X}^*}\big(\widetilde\gamma_N f_z-\widetilde\gamma_N f_{\mu,z}\big)=P_{\mathscr{X}^*}\gamma_N f_{D,z}.
\end{split}
\end{equation}
Hence $f_z\in\dom(A_{T,\Omega})\cap\ker(A_{max,\Omega}-zI)$, which implies $f_z\in\ker(A_{T,\Omega}-zI)$.
This yields $f_z=0$ as $z\in\rho(A_{T,\Omega})$ by assumption. Consequently, $\varphi=\widetilde\gamma_D f_z=0$ 
which ultimately implies that the operator $T+H(z)$ in \eqref{jussitm} is invertible for all 
$z\in\rho(A_{D,\Omega})\cap\rho(A_{T,\Omega})$.

Next, we shall show that $T+H(z)$ maps onto $\mathscr{X}^*$ whenever $z\in\rho(A_{D,\Omega})\cap\rho(A_{T,\Omega})$. 
For this purpose, let $\psi\in\mathscr{X}^*$ and choose $f\in\dom(A_{max,\Omega})$ such that
\begin{equation}\label{condis1}
\widetilde\gamma_D f=0\,\text{ and }\,P_{\mathscr{X}^*}\gamma_N f_D=\psi
\end{equation}
(here we once again use that the mapping \eqref{jussibm} is surjective).
Note that thanks to the first condition we have $f=f_D$. Since, by assumption, $z\in\rho(A_{T,\Omega})$ 
we also have the direct sum decomposition
\begin{equation}\label{yrdda-WACO}
\dom(A_{max,\Omega})=\dom(A_{T,\Omega})\,\dot{+}\,\ker(A_{max,\Omega}-zI).
\end{equation}
As such, $f$ may also be written in the form
\begin{equation}\label{deco7}
f=f_T+f_z,\,\text{ where }\,f_T\in\dom(A_{T,\Omega})\,\text{ and }\,f_z\in\ker(A_{max,\Omega}-zI).
\end{equation}
Next will make use of the decomposition of $f_T\in\dom(A_{T,\Omega})$ with respect 
to \eqref{fdeco}, that is, write $f_T$ in the form
\begin{equation}\label{deco8}
f_T=f_{D,T}+f_{\mu,T},\,\text{ where }\,f_{D,T}\in\dom(A_{D,\Omega})
\,\text{ and }\,f_{\mu,T}\in\ker(A_{max,\Omega}-\mu).
\end{equation}
One notes that $f_T\in\dom(A_{T,\Omega})$ implies $T\,\widetilde\gamma_D f_T=P_{\mathscr{X}^*}\gamma_N f_{D,T}$.
In particular, $\widetilde\gamma_D f_T\in\dom(T)\subset\mathscr{X}$ and therefore,  
$\iota_{\mathscr{X}}\widetilde\gamma_D f_T=\widetilde\gamma_D f_T$. 
It then follows from the first condition in \eqref{condis1} and \eqref{deco7} that
\begin{equation}\label{kgewdv.W.A.TEXAS}
\widetilde\gamma_D f_T=-\widetilde\gamma_D f_z.
\end{equation}
One computes 
\begin{equation}\label{jagvvv-ge3.WACO}
\begin{split}
\big(T+H(z)\big)\widetilde\gamma_D f_T &=\big(T+P_{\mathscr{X}^*}\big(\widetilde M_\Omega(z)-\widetilde M_\Omega(\mu)\big)
\iota_{\mathscr{X}}\big)\widetilde\gamma_D f_T
\\[2pt]
&=T\,\widetilde\gamma_D f_T+P_{\mathscr{X}^*}\big(\widetilde M_\Omega(z)\widetilde\gamma_D f_T
-\widetilde M_\Omega(\mu)\widetilde\gamma_D f_T\big)
\\[2pt]
&=P_{\mathscr{X}^*}\gamma_N f_{D,T}+P_{\mathscr{X}^*}\big(-\widetilde M_\Omega(z)\widetilde\gamma_D f_z 
-\widetilde M_\Omega(\mu)\widetilde\gamma_D(f_{D,T}+f_{\mu,T})\big)
\\[2pt]
&=P_{\mathscr{X}^*}\big(\gamma_N f_{D,T}+\widetilde\gamma_N f_z-\widetilde M_\Omega(\mu)\widetilde\gamma_D f_{\mu,T}\big)
\\[2pt]
&=P_{\mathscr{X}^*}\big(\gamma_N f_{D,T}+\widetilde\gamma_N f_z+\widetilde\gamma_N f_{\mu,T}\big) 
\\[2pt]
&=P_{\mathscr{X}^*}\widetilde\gamma_N\big(f_{D,T}+f_z+f_{\mu,T}\big) 
\\[2pt]
&=P_{\mathscr{X}^*}\widetilde\gamma_N f=P_{\mathscr{X}^*}\gamma_N f_D=\psi, 
\end{split}
\end{equation}
and hence it follows that the operator $T+H(z)$ in \eqref{jussitm} maps onto $\mathscr{X}^*$. 
We have shown that $T+H(z)$ in \eqref{jussitm} is bijective for all $z\in\rho(A_{D,\Omega})\cap\rho(A_{T,\Omega})$.

As $H(z)$ is a bounded operator from $\mathscr{X}$ to $\mathscr{X}^*$ and $T$ is self-adjoint it 
follows that $T+H(z)$ is closed as an operator from $\mathscr{X}$ onto $\mathscr{X}^*$. This implies 
that the inverse is closed as well, and hence bounded by the Closed Graph Theorem. \\[1mm] 
Next, it will be shown that the resolvent formula in the theorem holds. To get started, 
pick $f\in L^2(\Omega)$ and define
\begin{equation}\label{jussisg}
g:=(A_{D,\Omega}-zI)^{-1}f-\widetilde P_{D,\Omega}(z)\iota_{\mathscr{X}}\big(T+H(z)\big)^{-1}
P_{\mathscr{X}^*}\big(\widetilde P_{D,\Omega}(\bar z)\big)^*f,
\end{equation}
where, as above,  
\begin{equation}\label{katfrdC-WACO}
T+H(z)=T+P_{\mathscr{X}^*}\big(\widetilde M_\Omega(z)-\widetilde M_\Omega(\mu)\big)
\iota_{\mathscr{X}}:\mathscr{X}\supset\dom(T)\rightarrow\mathscr{X}^*.
\end{equation}
First, observe that $g\in\dom(A_{max,\Omega}-zI)$ and that 
$\ran\widetilde P_{D,\Omega}(z)\subset\ker(A_{max,\Omega}-zI)$ yields
\begin{equation}\label{amaxgf}
(A_{max,\Omega}-zI)g=f.
\end{equation}
We claim that $g$ belongs to $\dom(A_{T,\Omega})$. To justify this, it suffices to verify that the boundary condition
\begin{equation}\label{jussibc1}
T\,\widetilde\gamma_D g=P_{\mathscr{X}^*}\gamma_N g_D
\end{equation}
is satisfied. Making use of the decomposition $g=g_D+g_\mu$ one rewrites 
\begin{equation}\label{jgdX-WACO}
\begin{split}
\gamma_N g_D &=\gamma_N(g-g_\mu)=\widetilde\gamma_N g-\widetilde\gamma_N g_\mu
=\widetilde\gamma_N g+\widetilde M_\Omega(\mu)\widetilde\gamma_D g_\mu
\\[2pt]
&=\widetilde\gamma_N g+\widetilde M_\Omega(\mu)\widetilde\gamma_D g.
\end{split}
\end{equation}
Thus, the boundary condition \eqref{jussibc1} is equivalent to 
\begin{equation}\label{jussibc2}
T\,\widetilde\gamma_D g=P_{\mathscr{X}^*}\big(\widetilde\gamma_N g+\widetilde M_\Omega(\mu)\widetilde\gamma_D g\big).
\end{equation}
Next we verify that $g$ in \eqref{jussisg} satisfies \eqref{jussibc2}. First, we note that
\begin{equation}\label{rewws-WACO}
\begin{split}
\widetilde\gamma_D g &=-\iota_{\mathscr{X}}
\big(T+H(z)\big)^{-1}P_{\mathscr{X}^*}\big(\widetilde P_{D,\Omega}(\bar z)\big)^*f,
\\[2pt]
\widetilde\gamma_N g &=\widetilde\gamma_N(A_{D,\Omega}-zI)^{-1}f
-\widetilde\gamma_N\widetilde P_{D,\Omega}(z)\iota_{\mathscr{X}}\big(T+H(z)\big)^{-1}
P_{\mathscr{X}^*}\big(\widetilde P_{D,\Omega}(\bar z)\big)^*f
\\[2pt]
&=-\big(\widetilde P_{D,\Omega}(\bar z)\big)^*f 
+\widetilde M_\Omega(z)\iota_{\mathscr{X}}\big(T+H(z)\big)^{-1}
P_{\mathscr{X}^*}\big(\widetilde P_{D,\Omega}(\bar z)\big)^*f.
\end{split}
\end{equation}
This implies
\begin{equation}\label{jussibc3}
\begin{split}
T\,\widetilde\gamma_D g &=-T\big(T+H(z)\big)^{-1}
P_{\mathscr{X}^*}\big(\widetilde P_{D,\Omega}(\bar z)\big)^*f
\\[2pt]
&=-P_{\mathscr{X}^*}\big(\widetilde P_{D,\Omega}(\bar z)\big)^*f 
+H(z)\big(T+H(z)\big)^{-1}P_{\mathscr{X}^*}\big(\widetilde P_{D,\Omega}(\bar z)\big)^*f
\end{split}
\end{equation}
and
\begin{equation}\label{tr6fy-WACO}
\begin{split}
& \widetilde\gamma_N g+\widetilde M_\Omega(\mu)\widetilde\gamma_D g
\\[2pt]
&\quad=-\big(\widetilde P_{D,\Omega}(\bar z)\big)^*f 
+\big(\widetilde M_\Omega(z)-\widetilde M_\Omega(\mu)\big)\iota_{\mathscr{X}}\big(T+H(z)\big)^{-1}
P_{\mathscr{X}^*}\big(\widetilde P_{D,\Omega}(\bar z)\big)^*f,
\end{split} 
\end{equation}
hence
\begin{equation}\label{jussibc4}
\begin{split}
& P_{\mathscr{X}^*}\big(\widetilde\gamma_N g+\widetilde M_\Omega(\mu)\widetilde\gamma_D g\big)
\\[2pt]
&\quad=-P_{\mathscr{X}^*}\big(\widetilde P_{D,\Omega}(\bar z)\big)^*f+ H(z)\big(T+H(z)\big)^{-1}
P_{\mathscr{X}^*}\big(\widetilde P_{D,\Omega}(\bar z)\big)^*f.
\end{split}
\end{equation}
It now follows from \eqref{jussibc3} and \eqref{jussibc4} that \eqref{jussibc2} holds. 
Thus, $g\in\dom(A_{T,\Omega})$, and from \eqref{amaxgf} one concludes that 
\begin{equation}\label{trfcc-WACO}
(A_{T,\Omega}-zI)g=f,
\end{equation}
or equivalently, as $z\in\rho(A_{T,\Omega})$,
\begin{equation}\label{hafvav-WACO}
g=(A_{T,\Omega}-zI)^{-1}f.
\end{equation}
Thus, \eqref{jussisg} completes the proof.
\end{proof}
%%%%%%%

%%%%%%%%%%%%%%%%%%%%%%%%%%%%%%
%%%%%%%%%%%%%%%%%%%%%%%%%%%%%%
\section{The Case of Variable Coefficient Operators}  
\label{s11} 
%%%%%%%%%%%%%%%%%%%%%%%%%%%%%%
%%%%%%%%%%%%%%%%%%%%%%%%%%%%%%

The principal purpose of this section is to initiate a treatment of Laplace--Beltrami operators 
$-\Delta_g$ (and hence the case of variable coefficients induced by a metric $g$), perturbed by 
a scalar potential $V$. While this circle of ideas is worth pursuing further, we will at this point 
provide the basic results to demonstrate how the bulk of the material in Sections~\ref{s2}--\ref{s10} 
extends to perturbed Laplace--Beltrami operators on Lipschitz subdomains of compact boundaryless Riemannian manifolds. 

Throughout this final section we let $(M,g)$ be a compact, smooth ($C^\infty$), boundaryless 
manifold of (real) dimension $n\in\mathbb{N}$, $n\geq 2$, equipped with a $C^{1,1}$ Riemannian metric $g$. 
That is, in local coordinates the metric tensor $g$ may be expressed by
\begin{equation}\label{eqn.aaou}
g=\sum_{j,k=1}^n g_{jk}\,dx_j\otimes dx_k, 
\end{equation}
where the coefficients $g_{jk}$ are functions of class $C^{1,1}$. 
Hereafter, we shall often invoke Einstein's summation convention over repeated indices and suppress the
sigma symbol. The letter $g$ is also used to abbreviate 
\begin{equation}\label{eqn.aaoik}
g:=\det\big[(g_{jk})_{1\leq j,k\leq n}\big],
\end{equation}
and we shall use $(g^{jk})_{1\leq j,k\leq n}$ to denote the inverse of 
the matrix $(g_{jk})_{1\leq j,k\leq n}$, that is,  
\begin{equation}\label{eqn.altou}
(g^{jk})_{1\leq j,k\leq n}:=\big[(g_{jk})_{1\leq j,k\leq n}\big]^{-1}.
\end{equation}
The volume element $d{\mathcal V}_{\!g}$ on $M$ (with respect to the Riemannian metric $g$ from 
\eqref{eqn.aaou}) then can be written in local coordinates as
\begin{equation}\label{eqn.aa-pp}
d{\mathcal V}_{\!g}(x)=\sqrt{g(x)}\,d^n x.
\end{equation}
Consequently, given any relatively compact subset ${\mathcal{O}}$ of a coordinate patch 
(which we canonically identify with an open subset of the Euclidean space) it follows from \eqref{eqn.aa-pp}
that for any absolutely integrable function $f$ on ${\mathcal{O}}$ we have
\begin{equation}\label{VCnnd.M.WACO.EEE}
\int_{\mathcal{O}}f\,d{\mathcal V}_{\!g}=\int_{\mathcal{O}}f(x)\sqrt{g(x)}\,d^n x.
\end{equation}

As is customary, we use $\{\partial_j\}_{1\leq j\leq n}$ to denote a local basis in the tangent 
bundle $TM$. This implies that if $X,Y\in TM$ are locally expressed as $X=X_j\partial_j$, $Y=Y_j\partial_j$, 
then 
\begin{equation}\label{eqn.aa-phREED}
\langle X,Y\rangle_{TM}=X_jY_k g_{jk},
\end{equation}
where $\langle\,\cdot\,,\,\cdot\,\rangle_{TM}$ stands for the pointwise inner product in $TM$. 

Given an open set $\Omega\subset M$, for any scalar function $f\in C^1(\Omega)$,  
and any vector field $X\in C^1(\Omega,TM)$ locally written as $X=X_j\partial_j$, 
we may locally write (with the summation convention over repeated indices understood throughout)
\begin{equation}\label{rdf96-2DC.1RR}
{\rm grad}_gf:=(\partial_jf)g^{jk}\partial_k,
\quad X(f)=X_j(\partial_j f)=\langle{\rm grad}_gf,X\rangle_{TM},
\end{equation}
and
\begin{equation}\label{gMM-24GBN}
{\rm div}_gX:=g^{-1/2}\partial_j(g^{1/2}X_j)=\partial_jX_j+\Gamma^j_{jk}X_k,
\end{equation}
where $\Gamma^i_{jk}$ are the Christoffel symbols associated with the metric \eqref{eqn.aaou}. 
Moreover, for any scalar functions $f,h\in C^{\,1}(\Omega)$ and 
any vector field $X\in C^{\,1}(\Omega,TM)$, one has the product formulas
\begin{align}\label{UInabb}
\begin{split}
& {\rm grad}_g(fh)=f\,{\rm grad}_gh+h\,{\rm grad}_gf,\quad X(fh)=X(f)h+fX(h),
\\[2pt]
& {\rm div}_g(fX)=X(f)+f\,{\rm div}_g X.
\end{split}
\end{align}
Also, if ${\mathcal{O}}$ is a relatively compact subset of a coordinate patch 
(canonically identified with an open subset of the Euclidean ambient) then for any two scalar functions 
$\phi,\psi\in C^{\,1}(\overline{\mathcal{O}})$ we have
\begin{equation}\label{VCnnd.M.WACO}
\int_{\mathcal{O}}\langle{\rm grad}_g\phi, {\rm grad}_g\psi\rangle_{TM}\,d{\mathcal V}_{\!g}
=\int_{\mathcal{O}}\sum_{j,k=1}^n(\partial_j\phi)(x)(\partial_k\psi)(x)g^{jk}(x)\sqrt{g(x)}\,d^n x,
\end{equation}
thanks to \eqref{rdf96-2DC.1RR}, \eqref{eqn.aa-phREED}, and \eqref{eqn.aa-pp}.

The Laplace--Beltrami operator  
\begin{equation}\label{eqn.thu-TELL}
\Delta_g:={\rm div}_g\,{\rm grad}_g
\end{equation}
is expressed locally as
\begin{equation}\label{eqn.thuaou}
\Delta_gu=g^{-1/2}\partial_j\big(g^{jk}g^{1/2}\partial_ku\big).
\end{equation}
It satisfies the product formula
\begin{equation}\label{rd-jtR7h}
\Delta_g(uv)=v\Delta_gu+u\Delta_gv+2\langle{\rm grad}_g u,{\rm grad}_g v\rangle_{TM}.
\end{equation}
In the first part of this section we are interested in working with the 
formally symmetric Schr\"odinger differential expression 
\begin{equation}\label{eqn.lmnae}
L:=-\Delta_g+V,
\end{equation}
where the potential $V$ is a real-valued, essentially bounded, scalar-valued function on $M$. 

Given a nonempty open (necessarily bounded) set $\Omega\subset M$, for each integer $k\in{\mathbb{N}}$ we let 
$W^k(\Omega)$ stand for the $L^2$-based Sobolev space of order $k$ in $\Omega$.  
For each $k\in{\mathbb{N}}$ we also define
\begin{equation}\label{tr6ytG-MMM} 
\accentset{\circ}{W}^k(\Omega):=\overline{C^\infty_0(\Omega)}^{W^k(\Omega)}, 
\end{equation}
and equip the latter space with the norm inherited from $W^k(\Omega)$.
Corresponding to $\Omega=M$, for each $k\in{\mathbb{N}}$, we also set $W^{-k}(M):=(W^{k}(M))^\ast$. 

%%%%%%%%%
\begin{lemma}\label{h3.LLa} 
Assume $\Omega\subset M$ is a nonempty open set, pick $V\in L^\infty(M)$, and define $L$ as in \eqref{eqn.lmnae}. 
Then the graph norm 
\begin{equation}\label{ihgffd-YTG-MMM}
f\mapsto\|f\|_{L^2(\Omega)}+\|Lf\|_{L^2(\Omega)}, 
\quad\forall\,f\in\accentset{\circ}{W}^2(\Omega),
\end{equation}
is equivalent with the norm $\accentset{\circ}{W}^2(\Omega)$ inherits from $W^2(\Omega)$. 
\end{lemma}
%%%%%%%%%
\begin{proof}
From the work in \cite{MT99} one knows that if $\lambda>0$ is a sufficiently large real number then
the linear and bounded operator
\begin{equation}\label{h6f5fH-2-MMM}
L_\lambda:=L+\lambda:W^1(M)\rightarrow W^{-1}(M)
\end{equation}
is invertible, with bounded inverse
\begin{equation}\label{h6f5fH-3-MMM}
L_\lambda^{-1}:W^{-1}(M)\rightarrow W^1(M).
\end{equation}
In such a scenario, one can consider $E_\lambda\in{\mathcal{D}}'(M\times M)$, 
the Schwartz kernel of $L_{\lambda}^{-1}$, which is a distribution regular on 
the complement of the diagonal in $M\times M$. From \cite[Proposition~6.1]{MT00} one knows  
that the volume (Newtonian) potential operator
\begin{equation}\label{eqn.bweyq-MMM}
\Pi_\lambda f(x):=\int_M E_\lambda(x,y)f(y)\,d{\mathcal V}_{\!g}(y),\quad x\in M,
\end{equation}
is a linear and bounded mapping in the context 
\begin{equation}\label{eqn.ngqw-MMM}
\Pi_\lambda:L^2(M)\rightarrow W^{2}(M),
\end{equation}
which satisfies
\begin{equation}\label{eqn.ol-7767-MMM}
\Pi_\lambda(L_\lambda f)=f\,\text{ on }\,M,\quad\forall\,f\in W^{2}(M).
\end{equation}
Thus, for every $f\in C^\infty_0(\Omega)$ one estimates  
(recalling that tilde denotes the extension by zero to the entire ambient 
manifold $M$) 
\begin{align}\label{eqn.ol-7767-MMM-2}
\|f\|_{W^2(\Omega)} &=\big\|\widetilde{f}\,\big\|_{W^2(M)}
=\big\|\Pi_\lambda(L_\lambda\widetilde{f}\,)\big\|_{W^2(M)}
\nonumber\\[2pt]
&\leq C\big\|L_\lambda\widetilde{f}\,\big\|_{L^2(M)} 
\nonumber\\[2pt] 
&\leq C\big(\big\|L\widetilde{f}\,\big\|_{L^2(M)}+\big\|\widetilde{f}\,\big\|_{L^2(M)}\big)
\nonumber\\[2pt]
&=C\big(\|Lf\|_{L^2(\Omega)}+\|f\|_{L^2(\Omega)}\big),
\end{align}
for some constant $C\in(0,\infty)$, independent of $f$. In view of \eqref{tr6ytG-MMM} this implies
\begin{equation}\label{ihgffd-YTG-MMM-2}
\|f\|_{W^2(\Omega)}\leq C\big(\|f\|_{L^2(\Omega)}+\|Lf\|_{L^2(\Omega)}\big), 
\quad\forall\,f\in\accentset{\circ}{W}^2(\Omega).
\end{equation}
Since the opposite inequality is clear, the desired conclusion follows.
\end{proof}
%%%%%%%

Given an open nonempty set $\Omega\subset M$ and a real-valued potential $V\in L^\infty(M)$, 
we consider operator realizations of the differential expression $-\Delta_g+V$ in the Hilbert 
space $L^2(\Omega)$. We first define the \textit{preminimal} realization $L_{p,\Omega}$ 
of $-\Delta_g+V$ by
\begin{equation}\label{We-Q.1-MMM}
L_{p,\Omega}:=-\Delta_g+V,\quad\dom(L_{p,\Omega}):=C^\infty_0(\Omega).
\end{equation}
As such, $L_{p,\Omega}$ is a densely defined, symmetric operator in $L^2(\Omega)$, hence 
closable. Next, the \textit{minimal} realization $L_{min,\Omega}$ of $-\Delta_g+V$ is defined 
as the closure of $L_{p,\Omega}$ in $L^2(\Omega)$, that is, 
\begin{equation}\label{eqn:Amin1-MMM}
L_{min,\Omega}:=\overline{L_{p,\Omega}}.
\end{equation}
It follows that $L_{min,\Omega}$ is a densely defined, closed, symmetric operator in $L^2(\Omega)$. 
The \textit{maximal} realization $L_{max,\Omega}$ of $-\Delta_g+V$ is given by
\begin{equation}\label{We-Q.2-MMM}
L_{max,\Omega}:=-\Delta_g+V,\quad
\dom(L_{max,\Omega}):=\big\{f\in L^2(\Omega)\,\big|\,\Delta_g f\in L^2(\Omega)\big\},    
\end{equation}
where, much as in the Euclidean case, the expression $\Delta_g f$, $f\in L^2(\Omega)$, is understood in the sense of distributions. 
The assumption $V\in L^\infty(M)$ ensures that for $f\in L^2(\Omega)$ implies $\Delta_g f\in L^2(\Omega)$ 
if and only if $(-\Delta_g+V)f\in L^2(\Omega)$. 

Some of the most basic properties of the operators $L_{p,\Omega}$, $L_{min,\Omega}$, 
$L_{max,\Omega}$ are discussed below. 

%%%%%%%%%%
\begin{lemma}\label{l3.2-MMM}
Suppose $\Omega\subset M$ is an open nonempty set, and pick a real-valued potential $V\in L^\infty(M)$.
In this setting, let $L_{p,\Omega}$, $L_{min,\Omega}$, and $L_{max,\Omega}$ be as above. 
Then the operators $L_{min,\Omega}$ and $L_{max,\Omega}$ are adjoints of each other, that is, 
\begin{equation}\label{eqn:AminAmax-MMM}
L_{min,\Omega}^*=L_{p,\Omega}^*=L_{max,\Omega}\,\text{ and }\, 
L_{min,\Omega}=\overline{L_{p,\Omega}}=L_{max,\Omega}^*,
\end{equation}
and the closed symmetric operator $L_{min,\Omega}$ is semibounded from below by 
\begin{equation}\label{essinfv-MMM}
v_-:=\essinf_{x\in\Omega} V(x),
\end{equation}
that is, 
\begin{equation}\label{aminsemi-MMM}
(L_{min,\Omega}f,f)_{L^2(\Omega)}\geq v_{-}\|f\|^2_{L^2(\Omega)},\quad\forall\,f\in\dom(L_{min,\Omega}).
\end{equation}

In fact, the closed symmetric operator $L_{min,\Omega}$ is given by
\begin{equation}\label{eqn:Amin2-BBB-MMM}
L_{min,\Omega}=-\Delta_g+V,\quad\dom(L_{min,\Omega})=\accentset{\circ}{W}^2(\Omega),
\end{equation}
and $L_{min,\Omega}-v_{-}$ is strictly positive.
\end{lemma}
%%%%%%%
\begin{proof} 
Once Lemma~\ref{h3.LLa} has been established, all conclusions follow along the lines of the Euclidean case treated in Lemmas~\ref{l3.2}--\ref{l3.3}.
\end{proof}
%%%%%%%

%%%%%%%%%%%%%%%%
\subsection{Sobolev spaces on Lipschitz subdomains of a Riemannian manifold}\label{ss.MM.1}
%%%%%%%%%%%%%%%%

The reader is reminded that Sobolev spaces of fractional smoothness on $M$ are 
defined in a natural fashion, via localization (using a smooth partition of unity subordinate to a 
finite cover of $M$ with local coordinate charts) and pullback to the Euclidean model. This scale 
of spaces is then adapted to open subsets of $M$ via restriction, in analogy to the case $M={\mathbb{R}}^n$ 
considered earlier in Subsections~\ref{ss2.2}--\ref{ss2.3}, by setting
\begin{equation}\label{u64rLL-1}
H^s(\Omega):=\big\{u\big|_{\Omega}\,\big|\,u\in H^s(M)\big\},\quad s\in{\mathbb{R}}.
\end{equation}
In particular, $H^0(\Omega)$ coincides with $L^2(\Omega)$, the space of square integrable functions 
with respect to volume element $d{\mathcal V}_{\!g}$ in $\Omega$. 

Since bounded Lipschitz domains in the Euclidean setting are invariant under $C^1$-diffeomorphisms (cf. \cite{HMT07}), this class may be canonically defined on the 
manifold $M$, using local coordinate charts. If $\Omega\subset M$ is a Lipschitz domain 
then, as in the Euclidean setting, $H^k(\Omega)=W^k(\Omega)$ for every $k\in{\mathbb{N}}$. 
Given a Lipschitz domain $\Omega\subset M$, it is also possible to define (again, in a canonical manner, 
via localization and pullback) fractional Sobolev spaces on its boundary, $H^s(\partial\Omega)$, 
for $s\in[-1,1]$. In such a scenario one has $\big(H^s(\partial\Omega)\big)^\ast=H^{-s}(\partial\Omega)$ 
for each $s\in[-1,1]$, and $H^0(\partial\Omega)$ coincides with $L^2(\partial\Omega)$, the space of 
square integrable functions with respect to the surface measure $\sigma_g$ induced by the ambient 
Riemannian metric on $\partial\Omega$. Moreover, 
\begin{equation}\label{sop-uGG}
\big\{\phi\big|_{\partial\Omega}\,\big|\,\phi\in C^\infty(M)\big\}\,\text{ is dense in each }\,
H^s(\partial\Omega),\quad s\in[-1,1],
\end{equation}
and
\begin{equation}\label{sop-uGG.2}
H^{s_1}(\partial\Omega)\hookrightarrow H^{s_0}(\partial\Omega)\,\text{ continuously, whenever }\,
-1\leq s_0\leq s_1\leq 1.
\end{equation}

Next, if $\Omega\subset M$ is a given Lipschitz domain, the (Euclidean) nontangential 
approach region defined in \eqref{eq:MM1} has a natural version on $M$, simply replacing the standard 
Euclidean distance in ${\mathbb{R}}^n$ by the geodesic distance on $M$. With this interpretation, 
the nontangential maximal operator and nontangential boundary trace are then defined on Lipschitz 
subdomains of the manifold $M$ as in \eqref{eq:MM2} and \eqref{eq:MM2-BIS}, respectively. 
Then, virtually by design, it follows that all these objects satisfy similar properties to those of their Euclidean 
counterparts. See, for instance, \cite{MT99}, \cite{MT00}, \cite{Ta96}, \cite{Wl87}, and the references therein. 

Next, we record a regularity result which is a particular case of \cite[Proposition~4.9]{MT00}.
The reader is alerted to the fact that in Theorems~\ref{thm.Regularity-1} and \ref{thm.Regularity} 
we shall deviate from our typical condition $V\in L^{\infty}(M)$ and assumed $V\in L^p(M)$, with $p>n$, instead. 
This has its origins in the Calder\'on--Zygmund-type results in \cite{MT99}, \cite{MT00}, 
culminating in the mapping properties \eqref{eqn.plSSS}--\eqref{eqn.plSSS-2}.

%%%%%%%
\begin{theorem}\label{thm.Regularity-1}
Suppose $\Omega\subset M$ is a Lipschitz domain and pick a real-valued potential $V\in L^p(M)$ 
with $p>n$, where $n$ is the dimension of $M$. Then for any function $u\in C^1(\Omega)$ solving 
\begin{equation}\label{eqn.sjul-A}
Lu=0\,\text{  in }\,{\mathcal{D}}'(\Omega)
\end{equation}
one has
\begin{equation}\label{eqn.jslrh-111}
\mathcal{N}_{\kappa}u\in L^2(\partial\Omega)\Longleftrightarrow u\in H^{1/2}(\Omega)
\end{equation}
and, in fact, 
\begin{equation}\label{eqn.jslrh-22-aaa}
\big\|\mathcal{N}_{\kappa}u\big\|_{L^2(\partial\Omega)}
\approx\|u\|_{H^{1/2}(\Omega)},
\end{equation}
uniformly for $u\in C^1(\Omega)$ satisfying \eqref{eqn.sjul-A}.

Moreover, 
\begin{align}\label{eiBBB}
\begin{split}
& \text{if $\mathcal{N}_{\kappa}u\in L^2(\partial\Omega)$, then 
$u\big|^{\kappa-{\rm n.t.}}_{\partial\Omega}$ exists $\sigma$-a.e.~on $\partial\Omega$,
belongs to $L^2(\partial\Omega)$,} 
\\[4pt] 
& \quad\text{and satisfies $\|u\big|^{\kappa-{\rm n.t.}}_{\partial\Omega}\big\|_{L^2(\partial\Omega)}
\leq\big\|\mathcal{N}_{\kappa}u\big\|_{L^2(\partial\Omega)}$.}
\end{split} 
\end{align}
\end{theorem}
%%%%%%%%

The goal here is to establish a result similar in spirit to Theorem~\ref{thm.Regularity-1},
at a higher regularity level. Specifically, we shall prove the following theorem.
%%%%%%
\begin{theorem}\label{thm.Regularity}
Assume $\Omega\subseteq M$ is a Lipschitz domain, and pick a real-valued potential $V\in L^p(M)$ 
with $p>n$, where $n$ is the dimension of $M$. If the function $u\in C^1(\Omega)$ is such that
\begin{equation}\label{eqn.sjul}
Lu=0\,\text{ in }\,{\mathcal{D}}'(\Omega),
\end{equation}
then 
\begin{equation}\label{eqn.jslrh}
\mathcal{N}_{\kappa}({\rm grad}_g u)\in L^2(\partial\Omega)\Longleftrightarrow u\in H^{3/2}(\Omega)
\end{equation}
and, in fact, 
\begin{equation}\label{eqn.jslrh-22}
\big\|\mathcal{N}_{\kappa}u\big\|_{L^2(\partial\Omega)}
+\big\|\mathcal{N}_{\kappa}({\rm grad}_g u)\big\|_{L^2(\partial\Omega)}\approx\|u\|_{H^{3/2}(\Omega)},
\end{equation}
uniformly for $u\in C^1(\Omega)$ satisfying \eqref{eqn.sjul}.
Moreover, 
\begin{align}\label{eiBBB-NBVC}
& \text{if $\mathcal{N}_{\kappa}({\rm grad}_gu)\in L^2(\partial\Omega)$, then 
$u\big|^{\kappa-{\rm n.t.}}_{\partial\Omega}$ exists $\sigma$-a.e.~on $\partial\Omega$},  
\nonumber\\[2pt] 
& \quad\text{belongs to the Sobolev space $H^1(\partial\Omega)$, and satisfies}
\\[2pt] 
&\quad\big\|u\big|^{\kappa-{\rm n.t.}}_{\partial\Omega}\big\|_{H^1(\partial\Omega)}
\leq C\big\|\mathcal{N}_{\kappa}({\rm grad}_gu)\big\|_{L^2(\partial\Omega)}
+C\big\|\mathcal{N}_{\kappa}u\big\|_{L^2(\partial\Omega)},    
\nonumber
\end{align}
for some constant $C\in(0,\infty)$, independent of $u$. 
\end{theorem}
%%%%%%

As a preamble to the proof of Theorem~\ref{thm.Regularity}, we record a 
regularity result pertaining to the membership to fractional order Sobolev spaces
in Lipschitz domains, which is a slight variant of \cite[Lemma~2.34, p.~59]{MM13}. 
See \cite[Theorem~9.45, p.~444]{MMMT16} for a proof.

%%%%%%
\begin{lemma}\label{Bpq:L2}
Let $\Omega\subset M$ be a Lipschitz domain and suppose $u\in C^0(\Omega)\cap H^1_{\rm loc}(\Omega)$
is a function satisfying
\begin{align}\label{in-Pt:a-1bb}
{\mathcal{N}}_\kappa u\in L^2(\partial\Omega)\,\text{ and }\, 
\int_{\Omega}|({\rm grad}_g u)(x)|^2\,{\rm dist}_g(x,\partial\Omega)\,d{\mathcal V}_{\!g}(x)<\infty,
\end{align} 
where ${\rm dist}_g(x,\partial\Omega)$ denotes the geodesic distance from $x$ to $\partial\Omega$. 

Then $u\in H^{1/2}(\Omega)$ and there exists a constant $C\in(0,\infty)$, 
independent of $u$, with the property that 
\begin{align}\label{in-Pt:a-1}
\|u\|_{H^{1/2}(\Omega)}
&\leq C\|{\mathcal{N}}_\kappa u\|_{L^2(\partial\Omega)}
\nonumber\\[2pt]
&\quad
+C\bigg(\int_{\Omega}|({\rm grad}_g u)(x)|^2\,{\rm dist}_g(x,\partial\Omega)\,d{\mathcal V}_{\!g}(x)\bigg)^{1/2}.
\end{align} 
\end{lemma}
%%%%%%

We are now ready to present the proof of Theorem~\ref{thm.Regularity}.

%%%%%%%
\begin{proof}[Proof of Theorem~\ref{thm.Regularity}]
First, we note that given the nature of the conclusion we presently seek to confirm, 
there is no loss of generality in assuming that the differential operator $L$ 
satisfies the non-singularity hypothesis:
\begin{equation}\label{u6r44} 
\begin{rcases}
\text{for every Lipschitz domain }\,D\subseteq M
\\[2pt]
\text{(including the case $D=M$), and every}
\\[2pt]
\text{function }\,u\in\accentset{\circ}{H}^1(D),\,\text{ with }\,Lu=0\,\text{ in }\,D 
\end{rcases} 
\Longrightarrow u=0\,\text{ in }\,D.
\end{equation}
Indeed, since $L$ is elliptic and formally symmetric, 
by arguing as in the proof of \cite[Proposition~4.9]{MT00}, it is possible 
(after first arranging to work in a domain $\Omega$ which is 
very small relative to $M$, as in proof of \cite[Proposition~4.9]{MT00})
to suitably alter $L$ away from $\overline{\Omega}$ so that it becomes strictly 
positive definite on $M$, in the sense that there exists some $\varkappa>0$ such that
\begin{equation}\label{uy63333}
{}_{H^{-1}(M)}\big\langle Lw,w\big\rangle_{H^1(M)}\geq\varkappa\,\|w\|^2_{H^1(M)},
\quad\forall\,w\in H^1(M).
\end{equation}
Assume that this is the case, and pick a Lipschitz domain
$D\subseteq M$ along with some $u\in\accentset{\circ}{H}^1(D)$ satisfying $Lu=0$ in $D$.
Then, with tilde denoting extension by zero outside $D$ to the entire manifold $M$, 
it follows that $\widetilde{u}\in H^1_0(D)\subset H^1(M)$ satisfies 
${\rm supp}\,(L\widetilde{u})\subseteq\partial D$. In particular, this entails 
\begin{equation}\label{uy6g55RR}
{}_{H^{-1}(M)}\big\langle L\widetilde{u},\widetilde{u}\big\rangle_{H^1(M)}=0 
\end{equation}
as seen by approximating $\widetilde{u}\in H^1_0(D)$ in $H^1(M)$ with test 
functions on $M$ which are compactly supported in $D$ (cf.~\eqref{eq:12ddC} for the Euclidean setting). 
Since we are assuming that $L$ is strictly positive on $M$ (in the sense of \eqref{uy63333}), 
this forces $\widetilde{u}=0$ on $M$, hence ultimately $u=0$ in $D$. This 
concludes the justification of the fact that, for the current purposes, 
we may assume that the non-singularity hypothesis \eqref{u6r44} holds. 

The usefulness of the non-singularity hypothesis mentioned above  
is already apparent from choosing $D=M$ in \eqref{u6r44}, which  
implies that the linear and bounded operator
\begin{equation}\label{h6f5fH-2}
L:H^1(M)\rightarrow H^{-1}(M)
\end{equation}
is invertible, with bounded inverse
\begin{equation}\label{h6f5fH-3}
L^{-1}:H^{-1}(M)\rightarrow H^1(M).
\end{equation}
In particular, it makes sense to consider the Schwartz kernel of $L^{-1}$, a distribution on $M\times M$ 
which we denote by $E_L\in{\mathcal{D}}'(M\times M)$. From \cite{MT99} one knows the behavior of $E_L$ off 
the diagonal ${\rm diag}(M)$ of the Cartesian product $M\times M$, specifically, 
\begin{equation}\label{h6f5776t}
E_L\in C^1\big(M\times M\backslash {\rm diag}(M)\big).
\end{equation}

In turn, these considerations permit us to introduce the (boundary-to-boundary) single 
layer operator $S_L$ associated associated with $L$, by defining its action on any 
$\psi\in H^s(\partial\Omega)$ with $s\in[-1,0]$ according to the formula 
\begin{align}\label{eqn.plSSS} 
(S_L\psi)(x):={}_{H^{-s}(\partial\Omega)}\big\langle E_L(x,\cdot),\psi\big\rangle_{H^s(\partial\Omega)},
\quad\forall\,x\in\partial\Omega.
\end{align}
Then work in \cite{MT00} (involving the more general scale of Besov spaces) implies that
\begin{align}\label{eqn.plSSS-2}
S_L:H^{s}(\partial\Omega)\rightarrow H^{s+1}(\partial\Omega),\quad s\in[-1,0],
\end{align}
are invertible operators, with bounded, compatible inverses
\begin{align}\label{eqn.plSSS-3}
S_L^{-1}:H^{s+1}(\partial\Omega)\rightarrow H^{s}(\partial\Omega),\quad s\in[-1,0].
\end{align}
We also define the action of the boundary-to-domain version of the single 
layer operator ${\mathscr{S}}_L$ associated associated with $L$ on any 
$\psi\in H^s(\partial\Omega)$ with $s\in[-1,0]$ to be (compare with \eqref{eqn.plSSS})
\begin{align}\label{eqn.plSSS-4} 
({\mathscr{S}}_L\psi)(x):={}_{H^{-s}(\partial\Omega)}\big\langle E_L(x,\cdot),
\psi\big\rangle_{H^s(\partial\Omega)},\quad\forall\,x\in\Omega.
\end{align}
This operator satisfies the nontangential maximal function estimates 
(cf. \cite{MT99}, \cite{MT00})
\begin{align}\label{eqn.NTE}
\|\mathcal{N}_\kappa({\rm grad}_g\mathscr{S}_L\psi)\|_{L^2(\partial\Omega)}
&\leq C\|\psi\|_{L^2(\partial\Omega)},\quad\forall\,\psi\in L^2(\partial\Omega),
\\[2pt]
\|\mathcal{N}_\kappa(\mathscr{S}_L\psi)\|_{L^2(\partial\Omega)}
&\leq C\|\psi\|_{H^{-1}(\partial\Omega)},\quad\forall\,\psi\in H^{-1}(\partial\Omega),
\label{eqn.NTE.2}
\end{align}
as well as the square function estimates (cf. \cite{HMMM22}, \cite{MMMS15}, \cite{MMMT16})
\begin{align}\label{eqn.SFE}
\int_\Omega\big|\nabla^2(\mathscr{S}_L\psi)(x)\big|^2{\rm dist}_g(x,\partial\Omega)\,d{\mathcal V}_{\!g}(x)
&\leq C\|\psi\|^2_{L^2(\partial\Omega)},\quad\forall\,\psi\in L^2(\partial\Omega), 
\\[2pt]
\int_\Omega\big|\nabla(\mathscr{S}_L\psi)(x)\big|^2{\rm dist}_g(x,\partial\Omega)\,d{\mathcal V}_{\!g}(x)
&\leq C\|\psi\|^2_{H^{-1}(\partial\Omega)},\quad\forall\,\psi\in H^{-1}(\partial\Omega), 
\label{eqn.SFE.2}
\end{align}
for some constant $C\in(0,\infty)$ independent of $\psi$ (here and elsewhere $\nabla^2$ denotes the Hessian operator). 
These properties are going to be of basic importance for us later on.

After this preamble, we begin by considering the left-pointing implication 
in \eqref{eqn.jslrh}. To this end, assume a function $u\in C^1(\Omega)\cap H^{3/2}(\Omega)$, 
solving \eqref{eqn.sjul} (i.e., $Lu=0$ in ${\mathcal{D}}'(\Omega)$), has been given. 
Fix a smooth tangent field $X\in C^\infty(M,TM)$ and, with $\nabla_X$ denoting the covariant derivative 
along $X$, define 
\begin{equation}\label{eqn.dfws}
v:=\nabla_Xu\in C^0(\Omega)\cap H^{1/2}(\Omega). 
\end{equation}
Then there exists $C=C(\Omega,X)\in(0,\infty)$ such that
\begin{equation}\label{eqn.b-1ytrr}
\|v\|_{H^{1/2}(\Omega)}\leq C\|u\|_{H^{3/2}(\Omega)}.
\end{equation}
Moreover, since $Lu=0$ in $\Omega$, one can write 
\begin{equation}\label{eqn.bedou}
Lv=L(\nabla_Xu)=[L,\nabla_X]u\,\text{ in }\,\Omega,
\end{equation}
where, $[A,B]:=AB-BA$ abbreviates the commutator of the differential expressions $A$ and $B$. 
Locally, if $X=\sum_{\ell=1}^n a_\ell\partial_\ell$, where 
the coefficients $a_\ell$ are ($C^\infty$-) smooth functions, then a direct 
computation gives
\begin{align}\label{eqn.jlljs-11}
[L,\nabla_X]u &=\sum_{\ell=1}^n [L,a_\ell\partial_\ell]u
=-{\rm I}_1+{\rm I}_2-{\rm I}_3+{\rm I}_4-{\rm I}_5,
\end{align}
where 
\begin{align}\label{eqn.jlljs-12}
{\rm I}_1 &:=\sum_{j,k,\ell=1}^n g^{-1/2}\partial_j\big(g^{jk}g^{1/2}
(\partial_ka_\ell)(\partial_\ell u)\big)\in H^{-1/2}(\Omega),
\nonumber\\[2pt]
{\rm I}_2 &:=\sum_{j,k,\ell=1}^n a_\ell(\partial_\ell g^{-1/2})
\partial_j(g^{jk}g^{1/2}\partial_ku)\in H^{-1/2}(\Omega),
\nonumber\\[2pt]
{\rm I}_3 &:=\sum_{j,k,\ell=1}^n\Big\{(\partial_ja_\ell)g^{jk}\partial_\ell\partial_ku 
-g^{-1/2}a_\ell\partial_j\big(\partial_\ell(g^{jk}g^{1/2})\partial_ku\big)\Big\}\in H^{-1/2}(\Omega),
\nonumber\\[2pt]
{\rm I}_4 &:=\sum_{\ell=1}^n Va_\ell\partial_\ell u\in H^{-1}(\Omega), 
\nonumber\\[2pt]
{\rm I}_5 &:=\sum_{\ell=1}^n a_\ell\partial_\ell(Vu)\in H^{-1}(\Omega).
\end{align}
The memberships of ${\rm I}_1,{\rm I}_2,{\rm I}_3$ to $H^{-1/2}(\Omega)$ are readily justified by the 
fact that multiplication with functions from $\operatorname{Lip}(\Omega)$ preserves $H^s(\Omega)$ 
whenever $s\in[-1,1]$ (this follows in the same way as in the proof of Lemma~\ref{utfr}). 
To place ${\rm I}_4$ in $H^{-1}(\Omega)$ one observes that
\begin{align}\label{eqn.jl44}
{\rm I}_4=\sum_{\ell=1}^n Va_\ell\partial_\ell u\,\in\,\, & L^p(\Omega)\cdot H^{1/2}(\Omega)
\hookrightarrow L^{2n/3}(\Omega)\cdot L^{2n/(n-1)}(\Omega)
\nonumber\\[2pt]
&\hookrightarrow L^{2n/(n+2)}(\Omega)\hookrightarrow H^{-1}(\Omega)
\end{align}
(with continuous embeddings), by standard embedding results. Finally, to place ${\rm I}_5$ in $H^{-1}(\Omega)$ it suffices to note that 
\begin{equation}\label{eqn.hnerw.WACO}
Vu\in L^p(\Omega)\cdot H^{3/2}(\Omega)\hookrightarrow L^{2n/3}(\Omega)\cdot L^{2n/(n-3)}(\Omega)
\hookrightarrow L^2(\Omega). 
\end{equation}
The bottom line is that, as seen from \eqref{eqn.jlljs-11}-\eqref{eqn.jlljs-12},
\begin{equation}\label{eqn.hnerw}
f:=[L,\nabla_X]u\in H^{-1}(\Omega)\,\text{  and }\,\|f\|_{H^{-1}(\Omega)}\leq C\|u\|_{H^{3/2}(\Omega)},
\end{equation}
for some constant $C\in(0,\infty)$ which depends only on $\Omega,L,V,X$. In particular \eqref{eqn.bedou} becomes
\begin{equation}\label{eqn.mjbt}
Lv=f\in H^{-1}(\Omega).
\end{equation}

To proceed, we recall that $E_L(x,y)$ denotes the Schwartz kernel of $L^{-1}$ in 
\eqref{h6f5fH-3}. In \cite[Proposition~6.1]{MT00} it is shown that the volume 
(Newtonian) potential operator
\begin{equation}\label{eqn.bweyq}
\Pi_Lf(x):=\int_M E_L(x,y)f(y)\,d{\mathcal V}_{\!g}(y),\quad x\in M,
\end{equation}
originally acting on functions $f\in L^2(M)$, extends to a linear and bounded mapping
\begin{equation}\label{eqn.ngqw}
\Pi_L:\big(H^{1-s}(M)\big)^*=H^{s-1}(M)\rightarrow H^{s+1}(M),\quad
\forall\,s\in[-1,1],
\end{equation}
which satisfies
\begin{equation}\label{eqn.ol-7767}
L(\Pi_LF)=F\,\text{ in }\,{\mathcal{D}}'(M),\quad\forall\,F\in H^{s-1}(M),\,\,s\in[-1,1].
\end{equation}
Thus, we consider $F\in H^{-1}(M)$ such that $F|_{\Omega}=f$ as
distributions in $\Omega$, and $\|F\|_{H^{-1}(M)}\leq 2\|f\|_{H^{-1}(\Omega)}$. 
Then $w_X:=(\Pi_LF)\big|_{\Omega}\in H^1(\Omega)$ satisfies
\begin{equation}\label{eqn.olqwr}
Lw_X=(L\Pi_LF)\big|_{\Omega}=F\big|_{\Omega}=f\,\text{ in }\,\Omega,
\end{equation}
and
\begin{align}\label{eqn.hyTF}
\|w_X\|_{H^1(\Omega)} &\leq\|\Pi_L F\|_{H^{1}(M)}\leq C\|F\|_{H^{-1}(M)}
\nonumber\\[2pt]
&\leq C\|f\|_{H^{-1}(\Omega)}\leq C\|u\|_{H^{3/2}(\Omega)},
\end{align}
for some constant $C=C(\Omega,L,V,X)\in(0,\infty)$. 
In particular, if we now introduce $\vartheta_X:=v-w_X\in H^{1/2}(\Omega)$, then  
\begin{align}\label{eqn.mncv}
\begin{split}
& L\vartheta_X=Lv-Lw_X=f-f=0\,\text{ in }\,\Omega,\,\text{ and}
\\[2pt]
& \quad\|\vartheta_X\|_{H^{1/2}(\Omega)}\leq\|v\|_{H^{1/2}(\Omega)}+\|w_X\|_{H^{1/2}(\Omega)}
\leq C\|u\|_{H^{3/2}(\Omega)},
\end{split}
\end{align}
for some constant $C=C(\Omega,L,V,X)\in(0,\infty)$. Moreover, by the local 
elliptic regularity result established in \cite[Proposition~3.1]{MT00} one has
\begin{equation}\label{eqn.cvxwe}
\vartheta_X\in\bigcap_{1<p<\infty}W^{2,p}_{\rm loc}(\Omega),
\end{equation}
where $W^{2,p}_{\rm loc}(\Omega)$ is the subspace of $L^1_{\rm loc}(\Omega)$ consisting 
of functions with distributional derivatives of order $\leq 2$ belonging to $L^p_{\rm loc}(\Omega)$.
Since standard embedding results yield $W^{2,p}_{\rm loc}(\Omega)\subseteq C^1(\Omega)$ if $p>n$, one concludes that $\vartheta_X\in C^1(\Omega)$. In addition, Theorem~\ref{thm.Regularity-1} 
implies that $\mathcal{N}_\kappa(\vartheta_X)\in L^2(\partial\Omega)$ and
\begin{equation}\label{eqn.mncv-65r}
\big\|\mathcal{N}_\kappa(\vartheta_X)\big\|_{L^2(\partial\Omega)}
\leq C\|\vartheta_X\|_{H^{1/2}(\Omega)}\leq C\|u\|_{H^{3/2}(\Omega)},
\end{equation}
for some constant $C=C(\Omega,L,V,X)\in(0,\infty)$.

In summary, for every smooth vector field $X$ on $M$ we proved the decomposition 
\begin{align}\label{eqn.ikqe-ytt}
\begin{split}
& \nabla_X u=\vartheta_X+w_X\,\text{ in }\,\Omega,\,\text{ for some function}
\\[2pt]
& \quad\vartheta_X\in H^{1/2}(\Omega)\cap C^1(\Omega)\,\text{ satisfying }\,
\mathcal{N}_\kappa(\vartheta_X)\in L^2(\partial\Omega)
\\[2pt]
& \quad\text{as well as }\,\big\|\mathcal{N}_\kappa(\vartheta_X)\big\|_{L^2(\partial\Omega)}
\leq C\|u\|_{H^{3/2}(\Omega)},\,\text{ and some}
\\[2pt]
& \quad\text{function }\,w_X\in H^1(\Omega)\,\text{  with }\,\|w_X\|_{H^{1/2}(\Omega)}\leq C\|u\|_{H^{3/2}(\Omega)},
\end{split}
\end{align}
for some constant $C=C(\Omega,L,V,X)\in(0,\infty)$.

Next, we claim that the function $u\in H^{3/2}(\Omega)$ has the property 
\begin{equation}\label{eqn.plq-gtt}
\gamma_Du\in H^1(\partial\Omega)\,\,\text{ and }\,\,\|\gamma_Du\|_{H^1(\partial\Omega)}\leq C\|u\|_{H^{3/2}(\partial\Omega)},
\end{equation}
for some constant $C\in(0,\infty)$ independent of $u$. 
Since membership to $H^1(\partial\Omega)$ is a local property, we may work in local 
coordinates. For this portion of our proof one can assume that $M=\mathbb{R}^n$. 
Granted this fact, we adjust the notation in \eqref{eqn.ikqe-ytt}, namely, 
\begin{align}\label{eqn.ikqe-ytt-222}
\begin{split}
& \text{for }\,i\in\{1,\dots,n\}\,\text{ we write }\,\partial_i u=\vartheta_i+w_i\,\text{ in }\,\Omega,
\\[2pt]
& \quad\text{where }\,\vartheta_i\in H^{1/2}(\Omega)\cap C^1(\Omega)\,\text{ satisfies }\,  
\mathcal{N}_\kappa(\vartheta_i)\in L^2(\partial\Omega)
\\[2pt]
& \quad\text{as well as }\,\big\|\mathcal{N}_\kappa(\vartheta_i)\big\|_{L^2(\partial\Omega)}
\leq C\|u\|_{H^{3/2}(\Omega)},\,\text{ and where}
\\[2pt]
& \quad w_i\in H^1(\Omega)\,\text{ satisfies }\,\|w_i\|_{H^{1/2}(\Omega)}\leq C\|u\|_{H^{3/2}(\Omega)},
\end{split}
\end{align}
for some constant $C=C(\Omega,L,V,X)\in(0,\infty)$.

The strategy for proving the claim made in \eqref{eqn.plq-gtt}
is to fix an arbitrary test function $\psi\in C^\infty_0(\mathbb{R}^n)$ 
along with two arbitrary indices $j,k\in\{1,\dots,n\}$, with the intent of applying 
the divergence theorem to the vector field
\begin{equation}\label{eqn.plqasx}
\vec{F}:=u(\partial_k\psi)e_j-u(\partial_j\psi)e_k\,\text{ in }\,\Omega.
\end{equation}
With this goal in mind, one first observes that
\begin{equation}\label{eqn.p-RED}
\vec{F}\in\big[H^{3/2}(\Omega)\big]^n
\end{equation}
and, in the sense of distributions, 
\begin{align}\label{eqn.lasup}
\operatorname{div}\vec{F}
&=\partial_j(u\,\partial_k\psi)-\partial_k(u\,\partial_j\psi)
\nonumber
\\[2pt]
&=(\partial_ju)(\partial_k\psi)-(\partial_ku)(\partial_j\psi)\,\text{ in }\,\Omega,
\end{align}
hence 
\begin{align}\label{eqn.qwsao}
\operatorname{div}\vec{F}\in H^{1/2}(\Omega)\subset L^2(\Omega)\subset L^1(\Omega).
\end{align}
In addition, with $\nu=(\nu_1,\dots,\nu_n)$ denoting the outward unit normal to $\Omega$, one has 
\begin{align}\label{eqn.lpsr}
\nu\cdot\gamma_D\vec{F} 
&=(\gamma_Du)\big(\nu_j\big(\partial_k\psi\big)\big|_{\partial\Omega}
-\nu_k\big(\partial_j\psi\big)\big|_{\partial\Omega}\big)
\nonumber\\[2pt]
&=(\gamma_Du)(\partial_{\tau_{jk}}\psi)\,\text{ on }\,\partial\Omega.
\end{align}
Moreover, \eqref{eqn.p-RED} implies that for every $\varepsilon\in(0,1)$,
\begin{equation}\label{eqn.uaore}
\Delta\vec{F}\in\big[H^{-1/2}(\Omega)\big]^n
\subset\big[H^{-(3/2)+\varepsilon}(\Omega)\big]^n.
\end{equation}
Hence, Theorem~\ref{Ygav-75} applies (we recall that we are currently working in the Euclidean setting) 
and, if $\sigma$ denotes the canonical surface measure on $\partial\Omega$, one computes 
\begin{align}\label{eqn.ljamm}
\int_{\partial\Omega}(\gamma_Du)(\partial_{\tau_{jk}}\psi)\,d^{n-1}\sigma 
&=\int_{\partial\Omega}\nu\cdot\gamma_D\vec{F}\,d^{n-1}\sigma 
\nonumber\\[2pt]
&=\int_\Omega\operatorname{div}\vec{F}\,d^n x 
\nonumber\\[2pt]
&=\int_\Omega\big\{(\partial_ju)(\partial_k\psi)
-(\partial_ku)(\partial_j\psi)\big\}\,d^n x, 
\end{align}
by \eqref{eqn.lpsr} and \eqref{eqn.lasup}. At this point one introduces an approximating sequence, 
$\Omega_\ell\nearrow\Omega$ as $\ell\to\infty$, in the sense of Lemma~\ref{OM-OM}. From the local 
elliptic regularity result proved in \cite[Proposition~3.1]{MT00} one infers that 
\begin{equation}\label{eqn.jmner}
u\in\bigcap_{1<p<\infty}W^{2,p}_{\rm loc}(\Omega).
\end{equation}
In particular, 
\begin{equation}\label{eqn.bbvae}
u\in H^2(\Omega_\ell),\,\text{ for each }\,\ell\in\mathbb{N}.
\end{equation}
In turn, this implies that the vector field 
\begin{equation}\label{eqn.plqasx-G}
\vec{G}:=\psi(\partial_ju)e_k-\psi(\partial_ku)e_j\in\big[H^{1/2}(\Omega)\big]^n
\end{equation}
satisfies
\begin{align}\label{eqn.lasup-G}
{\rm div}\,\vec{G}
&=\partial_k(\psi\,\partial_j u)-\partial_j(\psi\,\partial_ku)
\nonumber
\\[2pt]
&=(\partial_k\psi)(\partial_ju)-(\partial_j\psi)(\partial_ku)
\end{align}
in the sense of distributions in $\Omega$.
In light of \eqref{eqn.bbvae}, this implies that for each $\ell\in{\mathbb{N}}$, 
\begin{equation}\label{eqn.plqasx-G2}
\vec{G}\big|_{\Omega_\ell}\in\big[H^1(\Omega_\ell)\big]^n
\end{equation}
and (cf.~\eqref{eqn:gammaDs.1-TR})
\begin{align}\label{eqn.lpsr-G}
\gamma_{\ell,D}\big(\vec{G}\big|_{\Omega_\ell}\big)
&=\big(\psi\big|_{\partial\Omega_\ell}\big)
\gamma_{\ell,D}\big(\partial_ju\big|_{\Omega_\ell}\big)e_k
\nonumber\\[2pt]
&\quad
-\big(\psi\big|_{\partial\Omega_\ell}\big)
\gamma_{\ell,D}\big(\partial_ku\big|_{\Omega_\ell}\big)e_j
\,\text{ on }\,\partial\Omega_\ell.
\end{align}
Invoking the last part of Theorem~\ref{Ygav-75} for the 
the vector field \eqref{eqn.plqasx-G2} in the Lipschitz domain $\Omega_\ell$, 
as well as employing the decomposition in \eqref{eqn.ikqe-ytt-222} 
(again, we recall that we are currently working in the Euclidean setting) this permits us to write: 
\begin{align}\label{eqn.syayy}
\int_\Omega\big\{ &(\partial_ju)(\partial_k\psi)
-(\partial_ku)(\partial_j\psi)\big\}\,d^n x
\nonumber\\[2pt]
&=\lim_{\ell\to\infty}\int_{\Omega_\ell}\big\{(\partial_ju)(\partial_k\psi)
-(\partial_ku)(\partial_j\psi)\big\}\,d^n x
\nonumber\\[2pt]
&=\lim_{\ell\to\infty}\int_{\Omega_\ell}{\rm div}\,\big(\vec{G}\big|_{\Omega_\ell}\big)\,d^n x
=\lim_{\ell\to\infty}\int_{\Omega_\ell}
\nu^\ell\cdot\gamma_{\ell,D}\big(\vec{G}\big|_{\Omega_\ell}\big)\,d^{n-1}\sigma_\ell
\nonumber\\[2pt]
&=\lim_{\ell\to\infty}
\int_{\partial\Omega_\ell}\big\{\nu_k^\ell\gamma_{\ell,D}\big(\partial_ju\big|_{\Omega_\ell}\big)
-\nu_j^\ell\gamma_{\ell,D}\big(\partial_ku\big|_{\Omega_\ell}\big)\big\}
\big(\psi\big|_{\partial\Omega_\ell}\big)\,d^{n-1}\sigma_\ell
\nonumber\\[2pt]
&=\lim_{\ell\to\infty}\int_{\partial\Omega_\ell}
\Big\{\nu_k^\ell\big[\vartheta_j\big|_{\partial\Omega_\ell} 
+\gamma_{\ell,D}\big(w_j\big|_{\Omega_\ell}\big)\big]
\nonumber\\[2pt]
&\hspace{2.1cm}
-\nu_j^\ell\big[\vartheta_k\big|_{\partial\Omega_\ell} 
+\gamma_{\ell,D}\big(w_k\big|_{\Omega_\ell}\big)\big]\Big\}
\big(\psi\big|_{\partial\Omega_\ell}\big)\,d^{n-1}\sigma_\ell
\nonumber\\[2pt]
&=\lim_{\ell\to\infty}\int_{\partial\Omega}
\Big\{\nu_k^\ell\circ\Lambda_\ell\big[\big(\vartheta_j\big|_{\partial\Omega_\ell}\big)\circ\Lambda_\ell 
+\gamma_{\ell,D}\big(w_j\big|_{\Omega_\ell}\big)\circ\Lambda_\ell\big]
\\[2pt]
&\hspace*{1.5cm}
-\nu_j^\ell\circ\Lambda_\ell\big[\big(\vartheta_k\big|_{\partial\Omega_\ell}\big)\circ\Lambda_\ell
+\gamma_{\ell,D}\big(w_k\big|_{\Omega_\ell}\big)\circ\Lambda_\ell\big]\Big\}\,
\big(\psi\big|_{\partial\Omega_\ell}\big)\circ\Lambda_\ell\,\omega_\ell\,d^{n-1}\sigma.
\nonumber
\end{align}
Keeping in mind that for every $j\in\{1,\dotsc,n\}$ and every $\ell\in\mathbb{N}$ one has 
\begin{equation}\label{eqn.ddac}
\big|\big(\vartheta_j\big|_{\partial\Omega_\ell}\big)\circ\Lambda_\ell\big|
\leq\mathcal{N}_\kappa(\vartheta_j)\,\text{ pointwise on }\,\partial\Omega,
\end{equation}
one then deduces from \eqref{eqn.ljamm}, \eqref{eqn.syayy}, and \eqref{eqn.ddac} that
\begin{align}\label{eqn.aljuu}
& \biggl\lvert\int_{\partial\Omega}(\gamma_Du) (\partial_{\tau_{jk}}\psi)\,d^{n-1}\sigma\biggr\rvert
\\[2pt]
& \quad\leq C\limsup_{\ell\to\infty}\int\limits_{\partial\Omega}
\sum_{m=1}^n\Big\{\mathcal{N}_\kappa(\vartheta_m) 
+\bigl\lvert\gamma_{\ell,D}\big(w_m\big|_{\Omega_\ell}\big)\circ\Lambda_\ell\big\rvert\Big\}
\Big|\big(\psi\big|_{\partial\Omega_\ell}\big)\circ\Lambda_\ell\Big|\,d^{n-1}\sigma.
\nonumber
\end{align}
Now, for each $m\in\{1,\dotsc,n\}$ and $\ell\in\mathbb{N}$, one estimates  
\begin{align}\label{eqn.kahhs}
& \int_{\partial\Omega}\big|\gamma_{\ell,D}\big(w_m\big|_{\Omega_\ell}\big)\circ\Lambda_\ell\big|^2\,d^{n-1}\sigma 
\nonumber\\[2pt]
& \quad\leq C\int_{\partial\Omega}\big|\gamma_{\ell,D}\big(w_m\big|_{\Omega_\ell}\big)
\circ\Lambda_\ell\big|^2\omega_\ell\,d^{n-1}\sigma
=C\int_{\partial\Omega_\ell}\big|\gamma_{\ell,D}\big(w_m\big|_{\Omega_\ell}\big)\big|^2\,d^{n-1}\sigma_\ell
\nonumber\\[2pt]
& \quad=C\big\|\gamma_{\ell,D}\big(w_m\big|_{\Omega_\ell}\big)\big\|^2_{L^2(\partial\Omega_\ell)}
\leq C\big\|\gamma_{\ell,D}\big(w_m\big|_{\Omega_\ell}\big)\big\|^2_{H^{1/2}(\partial\Omega_\ell)},
\end{align}
for some constant $C\in(0,\infty)$, independent of $\ell\in{\mathbb{N}}$. However, 
\begin{equation}\label{eqn.uwqs}
\gamma_{\ell,D}:H^1(\Omega_\ell)\rightarrow H^{1/2}(\partial\Omega_\ell)\,\text{ is bounded,}
\end{equation}
with operator norm controlled in terms of the Lipschitz character of $\Omega_\ell$. 
Hence, there exists $C\in(0,\infty)$ independent of $\ell\in\mathbb{N}$ such that 
\begin{equation}\label{eqn.uitw}
\big\|\gamma_{\ell,D}w\big\|_{H^{1/2}(\partial\Omega_\ell)}
\leq C\|w\|_{H^1(\Omega_\ell)},\quad\forall\,w\in H^1(\Omega_\ell).
\end{equation}
Based on this and \eqref{eqn.ikqe-ytt-222} one concludes that
\begin{align}\label{eqn.ujslaj}
\big\|\gamma_{\ell,D}\big(w_m\big|_{\Omega_\ell}\big)\big\|^2_{H^{1/2}(\partial\Omega_\ell)}
&\leq C\big\|w_m\big|_{\Omega_\ell}\big\|^2_{H^1(\Omega_\ell)}
\nonumber\\[2pt]
&\leq C\|w_m\|^2_{H^1(\Omega)}\leq C\|u\|^2_{H^{3/2}(\Omega)}, 
\end{align}
where the constant $C\in(0,\infty)$ remains independent of the index $\ell\in\mathbb{N}$. 
Thus, for each $m\in\{1,\dotsc,n\}$, from \eqref{eqn.kahhs} and \eqref{eqn.ujslaj} one obtains 
\begin{equation}\label{eqn.qyuna}
\int_{\partial\Omega}\big|\gamma_{\ell,D}\big(w_m\big|_{\Omega_\ell}\big)
\circ\Lambda_\ell\big|^2\,d^{n-1}\sigma\leq C\|u\|^2_{H^{3/2}(\Omega)}
\end{equation}
for some constant $C\in(0,\infty)$, independent of $\ell\in\mathbb{N}$. 
This ultimately proves that 
\begin{equation}\label{eqn.ngverq}
\sup_{\ell\in\mathbb{N}}\sum_{m=1}^n
\big\|\gamma_{\ell,D}\big(w_m\big|_{\Omega_\ell}\big)\circ\Lambda_\ell\big\|_{L^2(\partial\Omega)}
\leq C\|u\|_{H^{3/2}(\Omega)},
\end{equation}
for some constant $C\in(0,\infty)$, independent of $u$. 
Returning to \eqref{eqn.aljuu}, with the help of \eqref{eqn.ikqe-ytt-222} and 
\eqref{eqn.ngverq} one estimates 
\begin{align}\label{eqn.uiewq}
& \bigg|\int_{\partial\Omega} (\gamma_Du)(\partial_{\tau_{jk}}\psi)\,d^{n-1}\sigma\bigg|
\nonumber\\[2pt]
&\quad 
\leq C\limsup_{\ell\to\infty}\sum_{m=1}^n
\Big\{\big\|\mathcal{N}_\kappa(\vartheta_m)\big\|_{L^2(\partial\Omega)}  
\nonumber\\[2pt]
& \qquad\qquad\qquad
+\big\|\gamma_{\ell,D}\big(w_m\big|_{\Omega_\ell}\big)\circ\Lambda_\ell\big\|_{L^2(\partial\Omega)}\Big\}
\big\|\big(\psi\big|_{\partial\Omega_\ell}\big)\circ\Lambda_\ell\big\|_{L^2(\partial\Omega)}
\nonumber\\[2pt]
& \quad\leq C\|u\|_{H^{3/2}(\Omega)}\big\|\psi\big|_{\partial\Omega}\big\|_{L^2(\partial\Omega)},
\end{align}
where $C=C(\Omega)\in(0,\infty)$ is independent of $u$ and $\psi$. 
Since also
\begin{align}\label{eqn.uiTTRRR}
\|\gamma_Du\|_{L^2(\partial\Omega)}\leq
\|\gamma_Du\|_{H^{1/2}(\partial\Omega)}\leq
C\|u\|_{H^1(\Omega)}\leq C\|u\|_{H^{3/2}(\Omega)},
\end{align}
by reasoning much as before, based on \eqref{y6rf8HB.ttr}, \eqref{eqn:gammaDs.1} (with $s=1$), 
\eqref{eq:DDEEn.2}, and the characterization of $H^1(\partial\Omega)$ proved in Lemma~\ref{Dida}, it follows 
that $\gamma_Du\in H^1(\partial\Omega)$, as claimed. Moreover, from \eqref{eqn.uiewq}, 
\eqref{eqn.uiTTRRR}, and \eqref{trRFj6t}, one concludes the existence 
of a constant $C\in(0,\infty)$, independent of $u$, with the property that 
\begin{align}\label{eqn.uiTT-hr44}
\|\gamma_Du\|_{H^1(\partial\Omega)}\leq C\|u\|_{H^{3/2}(\partial\Omega)}.
\end{align}

Next, note that since $u\in H^{3/2}(\Omega)\cap C^1(\Omega)
\subset H^{1/2}(\Omega)\cap C^1(\Omega)$ and $Lu=0$ in $\Omega$, 
Theorem~\ref{thm.Regularity-1} applies and yields that
\begin{align}\label{eqn.tyqxa}
\mathcal{N}_\kappa u\in L^2(\partial\Omega)\,\text{ and }\, 
\big\|\mathcal{N}_\kappa u\big\|_{L^2(\partial\Omega)}\leq C\|u\|_{H^{1/2}(\Omega)},
\end{align}
for some constant $C\in(0,\infty)$, independent of $u$. 
On the other hand, given that
\begin{equation}\label{eqn.unare}
\begin{cases}
Lu=0\,\,\text{ in }\,\,\Omega,\quad u\in C^1(\Omega),
\\[2pt]
\mathcal{N}_\kappa u\in L^2(\partial\Omega),
\end{cases}
\end{equation}
we know from \cite[Proposition~3.1]{MT01} that the pointwise non-tangential trace
$u\big|^{\kappa-{\rm n.t.}}_{\partial\Omega}$ exists $\sigma$-a.e. on $\partial\Omega$. Hence we may invoke the 
(manifold version of) Lemma~\ref{Hgg-uh} to conclude that
\begin{equation}\label{eqn.waxyu}
\psi:=u\big|^{\kappa-{\rm n.t.}}_{\partial\Omega}=\gamma_Du\in H^1(\partial\Omega).
\end{equation}
Regarding $\psi$ as a function in $L^2(\partial\Omega)$, this means that $u$ solves the 
Dirichlet boundary value problem
\begin{equation}\label{eqn.eryqs}
\begin{cases}
Lu=0\,\text{ in }\,\Omega,\quad u\in C^1(\Omega),
\\[2pt]
\mathcal{N}_\kappa u\in L^2(\partial\Omega),
\\[4pt]
u\big|^{\kappa-{\rm n.t.}}_{\partial\Omega}=\psi\,\text{ on }\,\partial\Omega.
\end{cases}
\end{equation}
Granted the non-singularity condition \eqref{u6r44} we are currently assuming, 
it follows from \cite[Proposition~9.1]{MT99} that the solution of \eqref{eqn.eryqs} is unique 
and may be represented as 
\begin{equation}\label{wqybd}
u=\mathscr{S}_L\big(S_L^{-1}\psi\big)\,\text{ in }\,\Omega,
\end{equation}
where ${\mathscr{S}}_L$ is given by \eqref{eqn.plSSS-4} and $S_L^{-1}$ by  
\eqref{eqn.plSSS-3}. Since actually $\psi\in H^1(\partial\Omega)$, 
it follows that $S_L^{-1}\psi\in L^2(\partial\Omega)$. Consequently, for some constant $C\in(0,\infty)$, independent of $u$, one estimates  
\begin{align}\label{eqn.suary}
\big\|\mathcal{N}_\kappa({\rm grad}_g u)\big\|_{L^2(\partial\Omega)}
&=\big\|\mathcal{N}_\kappa\big({\rm grad}_g\mathscr{S}_L(S_L^{-1}\psi)\big)\big\|_{L^2(\partial\Omega)}
\nonumber\\[2pt]
&\leq C\|S_L^{-1}\psi\|_{L^2(\partial\Omega)}
\leq C\|\psi\|_{H^1(\partial\Omega)}
\nonumber\\[2pt]
&=\|\gamma_Du\|_{H^1(\partial\Omega)}\leq C\|u\|_{H^{3/2}(\partial\Omega)}, 
\end{align}
where the first inequality in \eqref{eqn.suary} is a consequence of \eqref{eqn.NTE}, 
while the last inequality has been proved in \eqref{eqn.uiTT-hr44}.
This completes the proof of the left-pointing implication in \eqref{eqn.jslrh}.

Turning our attention to the proof of the right-pointing implication 
in \eqref{eqn.jslrh}, we assume that $u\in C^1(\Omega)$ is such that
$\mathcal{N}_\kappa(\nabla u)\in L^2(\partial\Omega)$ and $Lu=0$ 
in $\Omega$. In light of the current assumptions on $u$, it follows from the manifold version 
of \eqref{eqn.hAAbvvb-1} that one also has $\mathcal{N}_\kappa u\in L^2(\partial\Omega)$. 
Having established this fact, \cite[Proposition~2.7]{MM13} implies that
\begin{align}\label{eqn.hnwte}
& \text{the nontangential trace $u\big|^{\kappa-{\rm n.t.}}_{\partial\Omega}$ exists  at $\sigma$-a.e.~point 
on $\partial\Omega$,}
\nonumber\\[2pt]
& \quad\text{the function $\phi:=u\big|^{\kappa-{\rm n.t.}}_{\partial\Omega}$ belongs 
to the space $H^1(\partial\Omega)$,} 
\\[2pt]
& \quad\text{and }\,
\|\phi\|_{H^1(\partial\Omega)}\leq C\big(\big\|\mathcal{N}_\kappa u\big\|_{L^2(\partial\Omega)}
+\|\mathcal{N}_\kappa({\rm grad}_g u)\|_{L^2(\partial\Omega)}\big),    
\nonumber
\end{align}
for some constant $C\in(0,\infty)$ independent of $u$. As such, it follows that the function 
$u$ solves the so-called regularity boundary value problem 
\begin{equation}\label{eqn.qsxct}
\begin{cases}
Lu=0\,\text{ in }\,\Omega,\quad u\in C^1(\Omega),
\\[2pt]
\mathcal{N}_\kappa u,\,\mathcal{N}_\kappa({\rm grad}_g u)\in L^2(\partial\Omega),
\\[4pt]
u\big|^{\kappa-{\rm n.t.}}_{\partial\Omega}=\phi\,\text{ on }\,\partial\Omega.
\end{cases} 
\end{equation}
Since we are presently assuming the non-singularity condition \eqref{u6r44}, it follows 
from Proposition~9.2 in \cite{MT99} (and its proof) that
\begin{equation}\label{eqn.uiqsd}
u=\mathscr{S}_L\big(S_L^{-1}\phi\big)\,\text{ in }\,\Omega.
\end{equation}
In addition, it follows from the local elliptic regularity 
result established in \cite[Proposition~3.1]{MT00} that
\begin{equation}\label{eqn.jmnMM}
u\in H^2_{\rm loc}(\Omega).
\end{equation}
If $\nabla^2$, as before, denotes the Hessian operator, then a combination of \eqref{eqn.uiqsd}, 
\eqref{eqn.SFE}, \eqref{eqn.plSSS-3}, and \eqref{eqn.hnwte} yields 
\begin{align}\label{eqn.vcqwr}
&\int_\Omega\big|(\nabla^2u)(x)\big|^2{\rm dist}_g(x,\partial\Omega)\,d{\mathcal V}_{\!g}(x)
\nonumber\\[2pt]
&\quad
=\int_\Omega\big|\big(\nabla^2\mathscr{S}_L\big(S_L^{-1}\phi\big)\big)(x)\big|^2
{\rm dist}_g(x,\partial\Omega)\,d{\mathcal V}_{\!g}(x)
\nonumber\\[2pt]
&\quad
\leq C\big\|S_L^{-1}\phi\big\|^2_{L^2(\partial\Omega)}
\leq C\|\phi\|^2_{H^1(\partial\Omega)}
\nonumber\\[2pt]
&\quad
\leq C\big(\big\|\mathcal{N}_\kappa({\rm grad}_g u)\big\|^2_{L^2(\partial\Omega)}
+\|\mathcal{N}_\kappa u\|^2_{L^2(\partial\Omega)}\big).
\end{align}
With these in hand, Lemma~\ref{Bpq:L2} implies that for each smooth vector field $X$ on $M$,  
\begin{align}\label{eqn.vcqw-11}
\nabla_X u\in H^{1/2}(\Omega), 
\end{align}
and for some constant $C=C(\Omega,X)\in(0,\infty)$, independent of $u$, 
\begin{align}\label{eqn.vcqw-22}
\|\nabla_X u\|_{H^{1/2}(\Omega)}
\leq C\big(\big\|\mathcal{N}_\kappa({\rm grad}_g u)\big\|_{L^2(\partial\Omega)}
+\|\mathcal{N}_\kappa u\|_{L^2(\partial\Omega)}\big).
\end{align}
In addition, from the manifold version of \eqref{eq:HBab.2} one deduces that
\begin{align}\label{eqn.lsuur1}
u\in L^{\frac{2n}{n-1}}(\Omega)\subset L^2(\Omega)\,\text{ and }\, 
\|u\|_{L^2(\Omega)}\leq C\|\mathcal{N}_\kappa u\|_{L^2(\partial\Omega)},
\end{align}
for some constant $C\in(0,\infty)$, independent of $u$.  
Having proved \eqref{eqn.vcqw-11}--\eqref{eqn.vcqw-22} and \eqref{eqn.lsuur1}, a quantitative 
lifting result (much as the one recorded in \eqref{eq:TRfav}) applies and yields 
\begin{equation}\label{eqn.wubf}
u\in H^{3/2}(\Omega)\,\text{ and }\,\|u\|_{H^{3/2}(\Omega)}
\leq C\big(\big\|\mathcal{N}_\kappa({\rm grad}_g u)\big\|_{L^2(\partial\Omega)}
+\|\mathcal{N}_\kappa u\|_{L^2(\partial\Omega)}\big), 
\end{equation}
for some constant $C\in(0,\infty)$, independent of $u$. This completes the justification 
of the right-pointing implication in \eqref{eqn.jslrh}. Since \eqref{eqn.wubf} also takes care of 
\eqref{eiBBB-NBVC}, the proof of Theorem~\ref{thm.Regularity} is complete.
\end{proof}
%%%%%%%

%%%%%%%%%%%%%%%%
\subsection{Sharp Dirichlet and Neumann traces on Lipschitz subdomains of Riemannian manifolds}\label{ss.MM.2}
%%%%%%%%%%%%%%%%
Much as in the Euclidean setting, if $\Omega\subset M$ is a Lipschitz domain, then 
the Dirichlet boundary trace map $C^\infty(\overline\Omega)\ni f\mapsto f\big|_{\partial\Omega}$
extends to operators (compatible with one another) 
\begin{align}\label{eqn:gammaDs.1-MMM}
\gamma_D:H^s(\Omega)\rightarrow H^{s-(1/2)}(\partial\Omega),\quad\forall\, 
s\in\big(\tfrac{1}{2},\tfrac{3}{2}\big),
\end{align}
that are linear, continuous, and surjective. We aim to further refine and extend this trace result 
in the theorem below, which the manifold counterpart of Theorem~\ref{YTfdf-T}, by also considering the 
end-point cases $s\in\big\{\tfrac{1}{2},\tfrac{3}{2}\big\}$ in the class of functions mapped by the 
Laplace--Beltrami operator in a better-than-expected Sobolev space.

%%%%%%%%%%
\begin{theorem}\label{YTfdf-T-MMM}
Assume that $\Omega\subset M$ is a Lipschitz domain and fix an arbitrary 
$\varepsilon>0$. Then the restriction of the boundary trace operator \eqref{eqn:gammaDs.1-MMM} 
to the space $\big\{u\in H^s(\Omega)\,\big|\,\Delta_g u\in H^{s-2+\varepsilon}(\Omega)\big\}$, 
originally considered for $s\in\big(\tfrac{1}{2},\tfrac{3}{2}\big)$, induces 
a well defined, linear, continuous operator 
\begin{equation}\label{eqn:gammaDs.2aux-MMM}
\gamma_D:\big\{u\in H^s(\Omega)\,\big|\,\Delta_g u\in H^{s-2+\varepsilon}(\Omega)\big\}
\rightarrow H^{s-(1/2)}(\partial\Omega),\quad\forall\,s\in\big[\tfrac{1}{2},\tfrac{3}{2}\big]     
\end{equation}
{\rm (}throughout, the space on the left-hand side of \eqref{eqn:gammaDs.2aux-MMM} equipped with the natural graph 
norm $u\mapsto\|u\|_{H^{s}(\Omega)}+\|\Delta_g u\|_{H^{s-2+\varepsilon}(\Omega)}${\rm )}, which continues 
to be compatible with \eqref{eqn:gammaDs.1-MMM} when $s\in\big(\tfrac{1}{2},\tfrac{3}{2}\big)$.
Thus defined, the Dirichlet trace operator possesses the following additional properties: \\[1mm] 
$(i)$ The Dirichlet boundary trace operator in \eqref{eqn:gammaDs.2aux-MMM} is surjective.
In fact, there exist linear and bounded operators
\begin{equation}\label{2.88X-NN-ii-RRDD-MMM}
\Upsilon_D:H^{s-(1/2)}(\partial\Omega)\rightarrow
\big\{u\in H^s(\Omega)\,\big|\,\Delta_g u\in L^2(\Omega)\big\},\quad s\in\big[\tfrac12,\tfrac32\big],
\end{equation}
which are compatible with one another and serve as right-inverses for the Dirichlet trace, that is, 
\begin{equation}\label{2.88X-NN2-ii-RRDD-MMM}
\gamma_D(\Upsilon_D\psi)=\psi,
\quad\forall\,\psi\in H^{s-(1/2)}(\partial\Omega)\,\text{ with }\,s\in\big[\tfrac12,\tfrac32\big].
\end{equation}
In fact, matters may be arranged so that each function in the range of $\Upsilon_D$ is harmonic, that is, 
\begin{equation}\label{2.88X-NN2-ii-RRDD.bis.3}
\Delta_g(\Upsilon_D\psi)=0,\quad\forall\,\psi\in H^{s-(1/2)}(\partial\Omega)
\,\text{ with }\,s\in\big[\tfrac{1}{2},\tfrac{3}{2}\big].
\end{equation}
$(ii)$ The Dirichlet boundary trace operator \eqref{eqn:gammaDs.2aux-MMM} is compatible
with the pointwise nontangential trace in the sense that:
\begin{align}\label{eqn:gammaDs.2auxBBB-MMM}
\begin{split}
& \text{if $u\in H^s(\Omega)$ has $\Delta_g u\in H^{s-2+\varepsilon}(\Omega)$ for some 
$s\in\big[\tfrac{1}{2},\tfrac{3}{2}\big]$,}    
\\[2pt]
& \quad\text{and if $u\big|^{\kappa-{\rm n.t.}}_{\partial\Omega}$ exists $\sigma$-a.e.~on $\partial\Omega$,  
then $u\big|^{\kappa-{\rm n.t.}}_{\partial\Omega}=\gamma_D u\in H^{s-(1/2)}(\partial\Omega)$.}
\end{split}
\end{align}
$(iii)$ The Dirichlet boundary trace operator $\gamma_D$ in 
\eqref{eqn:gammaDs.2aux-MMM} is the unique extension by continuity and density of the mapping 
$C^\infty(\overline{\Omega})\ni f\mapsto f\big|_{\partial\Omega}$. \\[1mm] 
$(iv)$ For each $s\in\big[\tfrac{1}{2},\tfrac{3}{2}\big]$ the Dirichlet 
boundary trace operator satisfies
\begin{align}\label{eqn:gammaDs.1-TR.2-MMM}
\begin{split}
& \gamma_D(\Phi u)=\big(\Phi\big|_{\partial\Omega}\big)\gamma_D u
\,\text{ at $\sigma$-a.e.~point on $\partial\Omega$, for all}
\\[2pt]
& \quad u\in H^s(\Omega)\,\text{ with }\,\Delta_g u\in H^{s-2+\varepsilon}(\Omega)
\,\text{ and }\,\Phi\in C^\infty(\overline\Omega).
\end{split}
\end{align}
$(v)$ For each $s\in\big[\tfrac{1}{2},\tfrac{3}{2}\big]$ such that $\varepsilon\not=\tfrac{3}{2}-s$, 
the null space of the Dirichlet boundary trace operator \eqref{eqn:gammaDs.2aux-MMM} satisfies
\begin{equation}\label{eq:EFFa-MMM}
{\rm ker}(\gamma_D)\subseteq H^{\,\min\!\{s+\varepsilon,3/2\}}(\Omega).
\end{equation}
In fact, the inclusion recorded in \eqref{eq:EFFa-MMM} is quantitative in the sense that,
whenever $s\in\big[\tfrac{1}{2},\tfrac{3}{2}\big]$ is such that $\varepsilon\not=\tfrac{3}{2}-s$, 
there exists a constant $C\in(0,\infty)$ with the property that
\begin{align}\label{gafvv.655-MMM} 
& \text{if $u\in H^s(\Omega)$ satisfies $\Delta_g u\in H^{s-2+\varepsilon}(\Omega)$ 
and $\gamma_D u=0$}   
\nonumber\\[2pt]
& \quad\text{then the function $u$ belongs to $H^{\,\min\!\{s+\varepsilon,3/2\}}(\Omega)$ and}    
\\[2pt]
& \quad\text{$\|u\|_{H^{\,\min\!\{s+\varepsilon,3/2\}}(\Omega)}\leq C\big(\|u\|_{H^s(\Omega)}
+\|\Delta_g u\|_{H^{s-2+\varepsilon}(\Omega)}\big)$.}   
\nonumber 
\end{align} 
\end{theorem}
%%%%%%%%%%
\begin{proof}
This may be established using the proof of Theorem~\ref{YTfdf-T} as a blue-print, 
substituting Theorems~\ref{thm.Regularity-1}--\ref{thm.Regularity} to the regularity 
and Fatou-type results in the Euclidean setting from Subsection~\ref{ss2.5}.
In addition, all relevant well-posedness results for the Dirichlet problem for the Laplace--Beltrami 
operator on Lipschitz subdomains of Riemannian manifolds may be found in \cite{MT99} and \cite{MT00}.
\end{proof}
%%%%%%%

As in the past, we will use the same symbol $\gamma_D$ in connection with 
either \eqref{eqn:gammaDs.1-MMM} or \eqref{eqn:gammaDs.2aux-MMM}, as the setting in 
which this is used will always be clear from the context. A particular case of 
Theorem~\ref{YTfdf-T-MMM}, which is particularly useful in applications, 
is singled out next.

%%%%%%%%%%
\begin{corollary}\label{YTfdf-T.NNN-MMM}
Suppose $\Omega\subset M$ is a given Lipschitz domain. Then the restriction of the operator 
\eqref{eqn:gammaDs.1-MMM} to 
$\big\{u\in H^s(\Omega)\,\big|\,\Delta_g u\in L^2(\Omega)\big\}$, originally 
considered for $s\in\big(\tfrac{1}{2},\tfrac{3}{2}\big)$, induces 
a well defined, linear, continuous operator 
\begin{align}\label{eqn:gammaDs.2-MMM}
\gamma_D:\big\{u\in H^s(\Omega)\,\big|\,\Delta_g u\in L^2(\Omega)\big\}
\rightarrow H^{s-(1/2)}(\partial\Omega),\quad\forall\, 
s\in\big[\tfrac{1}{2},\tfrac{3}{2}\big]    
\end{align}
{\rm (}throughout, the space on the left-hand side of \eqref{eqn:gammaDs.2-MMM} being equipped with the 
natural graph norm $u\mapsto\|u\|_{H^{s}(\Omega)}+\|\Delta_g u\|_{L^{2}(\Omega)}${\rm )}, which continues 
to be compatible with \eqref{eqn:gammaDs.1-MMM} when $s\in\big(\tfrac{1}{2},\tfrac{3}{2}\big)$,
and also with the pointwise nontangential trace, whenever the latter exists.

In addition, the following properties are true:

\begin{enumerate}
\item[(i)] The Dirichlet boundary trace operator in \eqref{eqn:gammaDs.2-MMM} is surjective
and, in fact, there exist linear and bounded operators
\begin{equation}\label{2.88X-NN-ii-RRDD-2-MMM}
\Upsilon_D:H^{s-(1/2)}(\partial\Omega)\to
\big\{u\in H^s(\Omega)\,\big|\,\Delta_g u\in L^2(\Omega)\big\},\quad s\in\big[\tfrac12,\tfrac32\big],
\end{equation}
which are compatible with one another and serve as right-inverses for the Dirichlet trace, that is, 
\begin{equation}\label{2.88X-NN2-ii-RRDD-2-MMM}
\gamma_D(\Upsilon_D\psi)=\psi,
\quad\forall\,\psi\in H^{s-(1/2)}(\partial\Omega)\,\text{ with }\,s\in\big[\tfrac12,\tfrac32\big].
\end{equation}
Actually, matters may be arranged so that each function in the range of $\Upsilon_D$ is harmonic, that is, 
\begin{equation}\label{2.88X-NN2-ii-RRDD.bis.4}
\Delta_g(\Upsilon_D\psi)=0,\quad\forall\,\psi\in H^{s-(1/2)}(\partial\Omega)
\,\text{ with }\,s\in\big[\tfrac12,\tfrac32\big].
\end{equation}

\item[(ii)] For each $s\in\big[\tfrac{1}{2},\tfrac{3}{2}\big]$, the null space   
of the Dirichlet boundary trace operator \eqref{eqn:gammaDs.2-MMM} satisfies
\begin{equation}\label{eq:EFFa.NNN-MMM}
{\rm ker}(\gamma_D)\subseteq H^{3/2}(\Omega).
\end{equation}
In fact, the inclusion in \eqref{eq:EFFa.NNN-MMM} is quantitative in the sense that
there exists a constant $C\in(0,\infty)$ with the property that
\begin{align}\label{gafvv.6577-MMM}
\begin{split}
&\text{whenever $u\in H^{1/2}(\Omega)$ with $\Delta_g u\in L^2(\Omega)$ satisfies $\gamma_D u=0$, then}
\\[2pt]
& \quad\text{$u\in H^{3/2}(\Omega)$ and 
$\|u\|_{H^{3/2}(\Omega)}\leq C\big(\|u\|_{L^2(\Omega)}+\|\Delta_g u\|_{L^2(\Omega)}\big)$.}
\end{split}
\end{align}
\end{enumerate}
\end{corollary}
%%%%%%%%%%%
\begin{proof}
All claims up to, and including, \eqref{eq:EFFa.NNN-MMM} are particular cases of 
the corresponding statement in Theorem~\ref{YTfdf-T-MMM}, choosing $\varepsilon=2-s$.
Finally, the proof of \eqref{gafvv.6577-MMM} follows the same pattern as that of its 
Euclidean counterpart in \eqref{gafvv.6577} (granted the well-posedness results in
\cite{MT99} and \cite{MT00}).
\end{proof}
%%%%%%%

To proceed, we make the following definition: 

%%%%%%%
\begin{definition}\label{DDD-1}
Given a nonempty open set $\Omega\subset M$ along with two numbers $s_0,s_1\in{\mathbb{R}}$ 
satisfying $s_0-1\geq s_1$, define $H^{s_0,s_1}_{\Delta_g}(\Omega,TM)$ as the collection 
of all vector fields $\vec{F}\in H^{s_0}(\Omega,TM)$ with the property that for every 
$x\in\overline{\Omega}$ there exists a local coordinate patch $U$ on $M$ which contains 
$x$ and such that if $\vec{F}=F_j\partial_j$ is the local writing of $\vec{F}$ in 
$U\cap\Omega$, then $\Delta_g F_j\in H^{s_1}(U\cap\Omega)$ for each $j\in\{1,\dots,n\}$. 
\end{definition}
%%%%%%%

In the context of Definition~\ref{DDD-1}, it is clear that $H^{s_0,s_1}_{\Delta_g}(\Omega,TM)$ is a vector space. 
The condition that $s_0-1\geq s_1$ ensures that this space is actually a module over $C^\infty(\overline{\Omega})$, that is, 
\begin{align}\label{uy54FF}
\begin{array}{c}
\psi\vec{F}\in H^{s_0,s_1}_{\Delta_g}(\Omega,TM)\,\text{ whenever} 
\\[4pt]
\psi\in C^\infty(\overline{\Omega})\,\text{ and }\,\vec{F}\in H^{s_0,s_1}_{\Delta_g}(\Omega,TM). 
\end{array}
\end{align}

%%%%%%%
\begin{definition}\label{DDD-2}
Given a Lipschitz domain $\Omega\subset M$, along with some real number $s\in\big[\tfrac{1}{2},\tfrac{3}{2}\big]$
and a vector field $\vec{F}\in H^{s,s-2+\varepsilon}_{\Delta_g}(\Omega,TM)$ with $\varepsilon\in(0,1)$, define 
\begin{equation}\label{75tfF}
\gamma_D\vec{F}\in H^{s-(1/2)}(\partial\Omega,TM) 
\end{equation}
as follows. First, one covers $\partial\Omega$ with finitely many coordinate patches $\{U_j\}_{1\leq j\leq N}$ 
and considers a smooth partition of unity associated to this cover. That is, one picks $\psi_j\in C^\infty_0(U_j)$, $1\leq j\leq N$, such that $\sum_{j=1}^N\psi_j=1$ near $\partial\Omega$. Then one sets 
\begin{equation}\label{hyr.a1-MMM}
\gamma_D\vec{F}:=\sum_{j=1}^N\gamma_D(\psi_j\vec{F})
\end{equation}
where, for each $j\in\{1,\dots,N\}$, if $\vec{F}=F^{(j)}_k\partial_k$ 
is the local writing of $\vec{F}$ in $U_j\cap\Omega$, and we have set
\begin{equation}\label{hyr.a2-MMM}
\gamma_D(\psi_j\vec{F}):=\gamma_D\big(\psi_j F^{(j)}_k\big)\partial_k\in H^{s-(1/2)}(\partial\Omega,TM).
\end{equation}
\end{definition}
%%%%%%%

That the Dirichlet traces in the right-hand side of \eqref{hyr.a2-MMM} make sense as functions
in $H^{s-(1/2)}(\partial\Omega)$ is a consequence of the membership 
$\vec{F}\in H^{s,s-2+\varepsilon}_{\Delta_g}(\Omega,TM)$ and \eqref{eqn:gammaDs.1-MMM}.   

The goal now is to state and prove a version of the divergence theorem which extends 
Theorem~\ref{Ygav-75} from the Euclidean setting to the context of Riemannian manifolds. 
As a preamble, we recall a few basic facts from differential geometry. 
Suppose that $\Omega\subset M$ is a Lipschitz domain. In local coordinates, if
\begin{align}\label{gMM-2}
\begin{split}
& \big(\nu^{\rm E}_j\big)_{1\leq j\leq n}\,\text{ is the outward unit normal on $\partial\Omega$}
\\[2pt]
& \quad\text{with respect to the Euclidean metric in ${\mathbb{R}}^n$},
\end{split}
\end{align}
and  
\begin{equation}\label{gMM-3}
{\mathfrak{G}}:=g^{rs}\nu^{\rm E}_r\nu^{\rm E}_s,
\end{equation}
then the unit outward normal to $\partial\Omega$ with respect to the Riemannian metric 
\begin{equation}\label{gMM-1}
g:=g_{jk}\,dx_j\otimes dx_k
\end{equation}
is given by (compare with \cite[Section~5.1 p.~2763,  Section~5.3, p.~2773]{HMT10})
\begin{equation}\label{gMM-4}
\nu=\nu_j\,\partial_j\in TM,\,\text{ where }\,\nu_j:=g^{jk}{\mathfrak{G}}^{-1/2}\nu^{\rm E}_k.
\end{equation}
In particular, 
\begin{equation}\label{gMM-4-BBa}
\nu^{\rm E}_j=g_{jk}{\mathfrak{G}}^{1/2}\nu_k.
\end{equation}
In addition, if locally we denote by $\sigma^{\rm E}$ the Euclidean surface measure on $\partial\Omega$,  then the surface measure $\sigma_g$ induced by the Riemannian metric \eqref{gMM-1} on $\partial\Omega$ is given by
\begin{equation}\label{gMM-7SD}
\sigma_g=\sqrt{g}\,{\mathfrak{G}}^{1/2}\,\sigma^{\rm E}.
\end{equation}
We are now ready to present the divergence theorem alluded to earlier. 

%%%%%%%%%%
\begin{theorem}\label{Ygav-75-MMM}
Consider a Lipschitz domain $\Omega\subset M$, with surface measure $\sigma_g$ and outward 
unit normal $\nu\in L^\infty(\partial\Omega,TM)$. Then for every given vector field 
$\vec{F}\in H^{1/2,-(3/2)+\varepsilon}_{\Delta_g}(\Omega,TM)$  
with $\varepsilon\in(0,1)$, satisfying ${\rm div}_g\vec{F}\in L^1(\Omega)$, 
one has
\begin{equation}\label{eq:GCxa-MMM}
\int_{\Omega}{\rm div}_g\vec{F}\,d{\mathcal{V}}_g
=\int_{\partial\Omega}\langle\nu,\gamma_D\vec{F}\rangle_{TM}\,d\sigma_g,
\end{equation}
where $\gamma_D\vec{F}$ is considered in the sense of Definition~\ref{DDD-2}
with $s=1/2$ {\rm (}implying $\gamma_D\vec{F}\in L^2(\partial\Omega,TM)${\rm )}.

As a corollary, \eqref{eq:GCxa-MMM} holds for every vector field 
$\vec{F}\in H^{(1/2)+\varepsilon}(\Omega,TM)$ for some $\varepsilon>0$
with the property that ${\rm div}_g\vec{F}\in L^1(\Omega)$ {\rm (}hence, in particular, 
for every vector field $\vec{F}\in H^1(\Omega,TM)${\rm )}.
\end{theorem}
%%%%%%%%%%%
\begin{proof}
We shall first prove \eqref{eq:GCxa-MMM} under the additional assumption that
there exists a local coordinate patch $U$ on $M$ such that
\begin{equation}\label{jy65tGV}
{\rm supp}\,\big(\vec{F}\big)\subset U\cap\overline{\Omega},
\end{equation}
and if $\vec{F}=F_j\partial_j$ is the local writing of $\vec{F}$ in $U\cap\Omega$, then 
\begin{equation}\label{jy65tGV-22222}
\Delta_g F_j\in H^{-(3/2)+\varepsilon}(U\cap\Omega)\,\text{ for each }\,j\in\{1,\dots,n\}.
\end{equation}
Assuming this to be the case, we identify $U$ with an Euclidean 
open set (via the corresponding local chart), and consider 
a Euclidean Lipschitz domain $\Omega'$ satisfying
\begin{align}\label{jy65tGV.222}
\begin{split}
& \Omega'\subset\Omega\cap U,\quad
{\rm supp}\,\big(\vec{F}\big)\cap\partial\Omega'\subset\partial\Omega,\quad
{\rm supp}\,\big(\vec{F}\big)\subset\overline{\Omega'},\quad
\\[2pt]
& \quad\text{and }\, 
\sigma'_g\lfloor(\partial\Omega'\cap\partial\Omega)=\sigma_g\lfloor(\partial\Omega'\cap\partial\Omega),
\end{split}
\end{align}
where $\sigma_g'$ is the surface measure induced by the Riemannian metric $g$ on $\partial\Omega'$. 

To proceed, for each $j\in\{1,\dots,n\}$ we invoke \cite{MT00} in order to solve the boundary value problem 
\begin{equation}\label{eqn:bNN.1-MMM}
\begin{cases}
\Delta_g G_j=\Delta_g F_j\,\text{ in $\Omega',\quad 
G_j\in H^{(1/2)+\varepsilon}(\Omega')$,}     
\\[2pt]  
\gamma_D G_j=0\,\text{ at $\sigma_g'$-a.e.~point on $\partial\Omega'$}.
\end{cases}   
\end{equation}
Then consider the vector field $\vec{G}:=G_j\partial_j$ in $\Omega'$ and set 
\begin{equation}\label{jh6gf5F}
\vec{h}:=\vec{F}-\vec{G}\,\text{ in }\,\Omega'. 
\end{equation}
It follows that $\vec{h}=h_j\partial_j$ with each component $h_j$ satisfying 
$\Delta_g h_j=0$ in $\Omega'$. Thus, $\vec{h}\in C^\infty(\Omega',TM)$ which, 
in particular, implies
\begin{equation}\label{e-MMa.1}
{\rm div}_g\,\vec{G}={\rm div}_g\vec{F}-{\rm div}_g\vec{h}\in L^1_{\rm loc}(\Omega').
\end{equation}
Moreover, $\vec{h}\in H^{1/2}(\Omega',TM)$ hence, if ${\mathcal{N}}'_\kappa$ denotes the nontangential 
maximal operator associated with the Lipschitz domain $\Omega'$, one concludes via Theorem~\ref{thm.Regularity-1} 
that ${\mathcal{N}}'_\kappa\vec{h}\in L^2(\partial\Omega')$ and $\vec{h}\big|^{\kappa-{\rm n.t.}}_{\partial\Omega'}$ exists 
$\sigma'_g$-a.e.~on $\partial\Omega'$, and belongs to $L^2(\partial\Omega',TM)$. If $\gamma'_D$ denotes 
the Dirichlet trace operator associated with the Lipschitz domain $\Omega'$, together with the last condition in 
\eqref{eqn:bNN.1-MMM} this forces 
\begin{equation}\label{e-MMaBBB}
\gamma'_D F_j=\gamma'_D h_j=h_j\big|^{\kappa-{\rm n.t.}}_{\partial\Omega'}\,\text{  on $\partial\Omega'$ for each $j$},
\end{equation}
where the last equality is a consequence of item $(ii)$ in Theorem~\ref{YTfdf-T} 
(cf.~\eqref{eqn:gammaDs.2auxBBB} for the Euclidean setting).

To proceed, we consider an approximating family 
$\Omega_\ell\nearrow\Omega'$ as $\ell\to\infty$ 
of the sort described in Lemma~\ref{OM-OM}, and recall that 
$\nu_{\ell}\circ\Lambda_\ell\to{\nu\,'}^{\,{\rm E}}$ as $\ell\to\infty$ both pointwise 
${\sigma\,'}^{\,{\rm E}}$-a.e.~on $\partial\Omega'$ and in $\big[L^2(\partial\Omega',{\sigma\,'}^{\,{\rm E}})\big]^n$. Moreover, the properties of the homeomorphisms $\Lambda_\ell$ allow one to conclude that  
for each $j\in\{1,\dots,n\}$,
\begin{align}\label{u554rfF}
\begin{array}{c}
\big(h_j\big|_{\partial\Omega_\ell}\big)\circ\Lambda_\ell\to h_j\big|^{\kappa-{\rm n.t.}}_{\partial\Omega'} 
\,\,\text{ as }\,\,\ell\to\infty
\\[4pt]
\text{both pointwise and in $\big[L^2(\partial\Omega',{\sigma\,'}^{\,{\rm E}})\big]^n$}, 
\end{array}
\end{align}
by Lebesgue's dominated convergence theorem (with uniform domination provided 
by a multiple of ${\mathcal{N}}'_\kappa\vec{h}\in L^2(\partial\Omega')$). 
Finally, one notes that the $\omega_\ell$'s appearing in the change of variable formula 
\eqref{eQQ-14} are uniformly bounded, and converge to $1$ as $\ell\to\infty$ pointwise 
${\sigma\,'}^{\,{\rm E}}$-a.e.~on $\partial\Omega$. Given these facts and keeping in mind that 
$\vec{h}\in C^\infty(\Omega',TM)$, one computes  
\begin{align}\label{Utfv-MMM}
& \lim_{\ell\to\infty}\int_{\partial\Omega_\ell}g^{1/2}\nu_{\ell,j}
\big(h_j\big|_{\partial\Omega_\ell}\big)\,d\sigma_\ell
\nonumber\\[2pt]
& \quad=\lim_{\ell\to\infty}\int_{\partial\Omega'}g^{1/2}(\nu_{\ell,j}\circ\Lambda_\ell)\cdot
\big(h_j\big|_{\partial\Omega_\ell}\big)\circ\Lambda_\ell\,\omega_\ell\,d{\sigma\,'}^{\,{\rm E}}
\nonumber\\[2pt]
& \quad=\int_{\partial\Omega'}g^{1/2}{\nu\,'}_{\!\!\!\!j}^{\,{\rm E}}
\Big(h_j\big|^{\kappa-{\rm n.t.}}_{\partial\Omega'}\Big)\,d{\sigma\,'}^{\,{\rm E}}
=\int_{\partial\Omega'}g_{jk}\nu\,'_{\!\!k}\gamma'_D F_j\,d\sigma'_g
\nonumber\\[2pt]
& \quad=\int_{\partial\Omega}g_{jk}\nu_k\gamma_D F_j\,d\sigma_g
=\int_{\partial\Omega}\big\langle\nu,\gamma_D\vec{F}\big\rangle_{TM}\,d\sigma_g.
\end{align}
Above, we used that (cf.~\eqref{gMM-4-BBa}--\eqref{gMM-7SD})
\begin{equation}\label{gMM-4-BBa-2}
{\nu\,'}_{\!\!\!\!j}^{\,{\rm E}}=g_{jk}{\mathfrak{G}}^{1/2}\nu\,'_{\!\!k}\,\,\text{ and }\,\,
{\sigma\,'}^{\,{\rm E}}=g^{-1/2}\,{\mathfrak{G}}^{-1/2}\,\sigma'_g,
\end{equation}
as well as \eqref{jy65tGV.222} and \eqref{eqn.aa-phREED}.

On the other hand, applying the (Euclidean) divergence theorem in each Euclidean Lipschitz 
domain $\Omega_\ell$ for the Euclidean vector field 
\begin{align}\label{Utfv-MMaaG}
\big(g^{1/2}h_j\big|_{\Omega_\ell}\big)_{1\leq j\leq n}\in\big[C^\infty(\overline{\Omega_\ell})\big]^n
\end{align}
(cf.\  Theorem~\ref{banff-3}), relying on Lebesgue's dominated convergence theorem, 
and invoking Lemma~\ref{YgLLam}, yields (cf.~\eqref{eqn.aa-pp}, \eqref{gMM-24GBN}, \eqref{jh6gf5F}),  
\begin{align}\label{Utfv.2-MMM}
& \lim_{\ell\to\infty}\int_{\partial\Omega_\ell}g^{1/2}\nu_{\ell,j}  
\big(h_j\big|_{\partial\Omega_\ell}\big)\,d\sigma_\ell
\nonumber\\[2pt]
& \quad=\lim_{\ell\to\infty}\int_{\Omega_\ell}\partial_j\big(g^{1/2}h_j\big)\,d^n x
=\lim_{\ell\to\infty}\int_{\Omega_\ell}g^{-1/2}\partial_j\big(g^{1/2}h_j\big)\,\sqrt{g}\,d^n x
\nonumber\\[2pt]
& \quad=\lim_{\ell\to\infty}\int_{\Omega_\ell}{\rm div}_g\vec{h}\,d{\mathcal{V}}_g
=\lim_{\ell\to\infty}\int_{\Omega_\ell}{\rm div}_g\vec{F}\,d{\mathcal{V}}_g
-\lim_{\ell\to\infty}\int_{\Omega_\ell}{\rm div}_g\vec{G}\,d{\mathcal{V}}_g
\nonumber\\[2pt]
& \quad=\lim_{\ell\to\infty}\int_{\Omega_\ell}{\rm div}_g\vec{F}\,d{\mathcal{V}}_g
-\lim_{\ell\to\infty}\int_{\Omega_\ell}\partial_j\big(g^{1/2}G_j\big)\,d^n x
\nonumber\\[2pt]
& \quad=\int_{\Omega}{\rm div}_g\vec{F}\,d{\mathcal{V}}_g
-\lim_{\ell\to\infty}\int_{\partial\Omega_\ell}
\nu_{\ell,j}\gamma_{\ell,D}\big(G_j\big|_{\Omega_\ell}\big)\,d\sigma_\ell, 
\end{align}
where, for each $\ell\in{\mathbb{N}}$, we denoted by $\gamma_{\ell,D}$ the Dirichlet boundary 
trace operator associated with the Lipschitz domain $\Omega_\ell$. The next step is to 
pick a small number $\delta\in\big(0,\min\{\tfrac{1}{2},\varepsilon\}\big)$ and then estimate 
\begin{align}\label{jab0uhb-MMM}
\begin{split} 
\bigg|\int_{\partial\Omega_\ell}\nu_{\ell,j}\gamma_{\ell,D}
\big(G_j\big|_{\Omega_\ell}\big)\,d\sigma_\ell\bigg|
&\leq\sum_{j=1}^n\big\|\gamma_{\ell,D}\big(G_j\big|_{\Omega_\ell}\big)\big\|
_{L^1(\partial\Omega_\ell,\sigma_\ell)}    
\\[2pt] 
&\leq C\sum_{j=1}^n
\big\|\gamma_{\ell,D}\big(G_j\big|_{\Omega_\ell}\big)\big\|_{H^\delta(\partial\Omega_\ell)}
\end{split} 
\end{align}
for some constant $C\in(0,\infty)$, independent of $\ell\in{\mathbb{N}}$. 
Since by \eqref{incl-Yb.EE} and \eqref{eqn:bNN.1-MMM}, $G_j\in\accentset{\circ}{H}^{(1/2)+\delta}(\Omega)$, 
it follows from Lemma~\ref{Tgav9jy} (used with $s=\tfrac{1}{2}+\delta\in(\tfrac{1}{2},1)$) that 
\begin{align}\label{jab0uhb.2-MMM}
\lim_{\ell\to\infty}\sum_{j=1}^n\big\|\gamma_{\ell,D}
\big(G_j\big|_{\Omega_\ell}\big)\big\|_{H^\delta(\partial\Omega_\ell)}=0.
\end{align}
At this stage, \eqref{eq:GCxa-MMM} follows from \eqref{Utfv-MMM}--\eqref{jab0uhb.2-MMM}.

Finally, it remains to dispense with the additional assumption \eqref{jy65tGV}. 
To this end, one covers $\overline{\Omega}$ with finitely many coordinate patches $\{U_k\}_{1\leq k\leq N}$ 
and consider a family of functions $\psi_k\in C^\infty_0(U_k)$, $1\leq k\leq N$, 
such that $\sum_{k=1}^N\psi_k=1$ near $\overline{\Omega}$. Then, by \eqref{uy54FF}, each vector 
field $\psi_k\vec{F}$ satisfies the hypotheses that permits one to conclude that
\begin{equation}\label{eq:GCxa-MMM-12ew}
\int_{\Omega}{\rm div}_g\big(\psi_k\vec{F}\big)\,d{\mathcal{V}}_g
=\int_{\partial\Omega}\big\langle\nu,\gamma_D\big(\psi_k\vec{F}\big)\big\rangle_{TM}\,d\sigma_g,\quad\forall\,k\in\{1,\dots,N\}.
\end{equation}
Since the sub-collection of $\{\psi_k\}_{1\leq k\leq N}$ consisting 
of those functions whose support intersects $\partial\Omega$ does constitute a smooth 
partition of unity near $\partial\Omega$, \eqref{eq:GCxa-MMM-12ew} and \eqref{hyr.a1-MMM} imply that
\begin{align}\label{eq:GCxa-MMM-12333}
\int_{\partial\Omega}\langle\nu,\gamma_D\vec{F}\rangle_{TM}\,d\sigma_g
&=\sum_{k=1}^N\int_{\partial\Omega}\langle\nu,\gamma_D(\psi_k\vec{F})\rangle_{TM}\,d\sigma_g
\nonumber\\[2pt]
&=\sum_{k=1}^N\int_{\Omega}{\rm div}_g(\psi_k\vec{F})\,d{\mathcal{V}}_g
=\int_{\Omega}{\rm div}_g\vec{F}\,d{\mathcal{V}}_g,
\end{align}
as wanted. 
\end{proof}
%%%%%%%%%%%

We shall find it useful to have a version of the divergence theorem, complementing Theorem~\ref{Ygav-75-MMM}, 
for vector fields whose divergence is not necessarily an absolutely integrable function. This task is accomplished below. 

%%%%%%%%%%
\begin{theorem}\label{Ygav-75-BIS-MMM}
Suppose $\Omega\subset M$ is a Lipschitz domain with surface measure $\sigma_g$ 
and outward unit normal $\nu$. Consider a vector field
$\vec{F}\in H^{1/2,-(3/2)+\varepsilon}_{\Delta_g}(\Omega,TM)$ for some $\varepsilon\in(0,1)$
with the property that ${\rm div}_g\vec{F}\in H^{-(1/2)+\varepsilon}(\Omega)$. 
Then 
\begin{equation}\label{eq:GCxa-BIS-MMM}
{}_{H^{(1/2)-\varepsilon}(\Omega)}\big\langle{\bf 1},{\rm div}_g\vec{F}
\big\rangle_{H^{-(1/2)+\varepsilon}(\Omega)}
=\int_{\partial\Omega}\langle\nu,\gamma_D\vec{F}\rangle_{TM}\,d\sigma_g,
\end{equation}
where ${\bf 1}$ denotes the constant function identically to $1$ in $\Omega$, 
and the action of $\gamma_D$ on $\vec{F}$ is considered in the sense
of Definition~\ref{DDD-2} with $s=1/2$ {\rm (}implying 
$\gamma_D\vec{F}\in L^2(\partial\Omega,TM)${\rm )}.
\end{theorem}
%%%%%%%%%%%
\begin{proof}
As in the proof of Theorem~\ref{Ygav-75-MMM}, making use of a smooth partition of unity, 
there is no loss of generality in assuming that there exists a local coordinate patch $U$ on $M$ such that
\begin{equation}\label{jy65tGV-MMM}
{\rm supp}\,\big(\vec{F}\big)\subset U\cap\overline{\Omega},
\end{equation}
and if $\vec{F}=F_j\partial_j$ is the local writing of $\vec{F}$ in $U\cap\Omega$, then 
\begin{equation}\label{jy65tGV-MMM-2222}
\Delta_g F_j\in H^{-(3/2)+\varepsilon}(U\cap\Omega)\,\text{ for each }\,j\in\{1,\dots,n\}.
\end{equation}
Assuming this to be the case, we identify $U$ with an Euclidean 
open set (via the corresponding local chart), and consider 
a Euclidean Lipschitz domain $\Omega'$ satisfying
\begin{align}\label{jy65tGV.222-MMM}
\begin{split}
& \Omega'\subset\Omega\cap U,\quad
{\rm supp}\,\big(\vec{F}\big)\cap\partial\Omega'\subset\partial\Omega,\quad
{\rm supp}\,\big(\vec{F}\big)\subset\overline{\Omega'}, 
\\[2pt] 
& \quad\text{and }\, 
\sigma'_g\lfloor(\partial\Omega'\cap\partial\Omega)=\sigma_g\lfloor(\partial\Omega'\cap\partial\Omega),
\end{split}
\end{align}
where $\sigma_g'$ is the surface measure induced by the Riemannian metric $g$ on $\partial\Omega'$. 
In particular, if we let $\vec{G}=G_j\partial_j$ solve \eqref{eqn:bNN.1-MMM} and set 
\begin{equation}\label{u6f434edD}
\vec{h}:=\vec{F}-\vec{G}=h_j\partial_j\,\text{ in }\,\Omega'
\end{equation}
then, as before, 
\begin{align}\label{NEWDVT.1-MMM}
& h_j\in C^\infty(\Omega')\cap H^{1/2}(\Omega'),
\\[2pt]
& \Delta_g h_j=0\,\text{ in }\,\Omega',\quad{\mathcal{N}}'_\kappa h_j\in L^2(\partial\Omega'),
\label{NEWDVT.2-MMM}
\\[2pt]
& \gamma_D\vec{F}=\Big(h_j\big|^{\kappa-{\rm n.t.}}_{\partial\Omega}\Big)\partial_j\in L^2(\partial\Omega',TM).
\label{NEWDVT.3-MMM}
\end{align}
Granted the current hypotheses, one also has  
\begin{align}\label{NEWDVT.4-MMM}
{\rm div}_g\vec{h}={\rm div}_g\vec{F}-{\rm div}_g\vec{G}\in L^1_{\rm loc}(\Omega')
\cap H^{-(1/2)+\varepsilon}(\Omega').
\end{align}
Since each $G_j\in\accentset{\circ}{H}^{(1/2)+\varepsilon}(\Omega')$, by  
\eqref{eqn:bNN.1-MMM} and \eqref{incl-Yb.EE}, it follows that there exists a sequence 
$\{G_j^{k}\}_{k\in{\mathbb{N}}}\subset C^\infty_0(\Omega')$ with the property that
\begin{align}\label{NEWDVT.5-MMM}
G_j^{k}\to G_j\,\text{ in }\,H^{(1/2)+\varepsilon}(\Omega')\,\text{ as }\,k\to\infty.
\end{align}
As a consequence, if for each $k\in{\mathbb{N}}$ one sets $\vec{G}^k :=G_j^{k}\partial_j$, then
\begin{align}\label{NEWDVT.6-MMM}
{\rm div}_g\vec{G}^k\to{\rm div}_g\vec{G}\,\text{ in }\,H^{-(1/2)+\varepsilon}(\Omega')\,\text{ as }\,k\to\infty,
\end{align}
and hence, 
\begin{align}\label{NEWDVT.7-MMM}
{}_{H^{(1/2)-\varepsilon}(\Omega')}\big\langle{\bf 1},{\rm div}_g\vec{G}
\big\rangle_{H^{-(1/2)+\varepsilon}(\Omega')}
&=\lim_{k\to\infty}{}_{H^{(1/2)-\varepsilon}(\Omega')}\big\langle{\bf 1},{\rm div}_g\vec{G}^{k}
\big\rangle_{H^{-(1/2)+\varepsilon}(\Omega')}
\nonumber\\[2pt]
&=\lim_{k\to\infty}\int_{\Omega'}g^{-1/2}\partial_j\big(g^{1/2}G^k_j\big)\sqrt{g}\,d^nx
\nonumber\\[2pt]
&=\lim_{k\to\infty}\int_{\Omega'}\partial_j\big(g^{1/2}G^k_j\big)\,d^nx
\nonumber\\[2pt]
&=\lim_{k\to\infty}\int_{\partial\Omega'}\nu'_j
\big(G^k_j\big|_{\partial\Omega'}\big)\,d\sigma^{'{\rm E}}=0,
\end{align}
given that $\vec{G}^{k}\in\big[C^\infty_0(\Omega')\big]^n$ for every $k\in{\mathbb{N}}$. 
We then proceed to write
\begin{align}\label{NEWDVT.8-MMM}
{}_{H^{(1/2)-\varepsilon}(\Omega)}\big\langle{\bf 1},{\rm div}_g\vec{F}
\big\rangle_{H^{-(1/2)+\varepsilon}(\Omega)}
&={}_{H^{(1/2)-\varepsilon}(\Omega')}\big\langle{\bf 1},{\rm div}_g\vec{F}
\big\rangle_{H^{-(1/2)+\varepsilon}(\Omega')}
\nonumber\\[2pt]
&={}_{H^{(1/2)-\varepsilon}(\Omega')}\big\langle{\bf 1},{\rm div}_g\vec{h}
\big\rangle_{H^{-(1/2)+\varepsilon}(\Omega')}.
\end{align}
The first equality above is implied by \eqref{trans-U}, \eqref{jy65tGV-MMM}, and the first line 
of \eqref{jy65tGV.222-MMM}, while the second equality is a consequence of \eqref{NEWDVT.7-MMM} 
and \eqref{NEWDVT.4-MMM}.

As in the past, we introduce an approximating family 
$\Omega_\ell\nearrow\Omega'$ as $\ell\to\infty$ 
(described in Lemma~\ref{OM-OM}). Then one can write 
\begin{align}\label{Utfv.2BIS-MMM}
{}_{H^{(1/2)-\varepsilon}(\Omega')}\big\langle{\bf 1},{\rm div}_g\vec{h}
\big\rangle_{H^{-(1/2)+\varepsilon}(\Omega')}
&=\lim_{\ell\to\infty}\int_{\Omega_\ell}{\rm div}_g\vec{h}\,\sqrt{g}\,d^n x
\nonumber\\[2pt]
&=\lim_{\ell\to\infty}\int_{\Omega_\ell}\partial_j\big(g^{1/2}h\big)\,d^n x
\nonumber\\[2pt]
&=\lim_{\ell\to\infty}\int_{\partial\Omega_\ell}g^{1/2}\nu_{\ell,j} 
\big(h_j\big|_{\partial\Omega_\ell}\big)\,d\sigma_\ell
\nonumber\\[2pt]
&=\int_{\partial\Omega}\big\langle\nu,\gamma_{D}\vec{F}\big\rangle_{TM}\,d\sigma_g, 
\end{align}
where the first equality is implied by Lemma~\ref{Lapp-L1} and \eqref{NEWDVT.4-MMM}, 
the second equality relies on \eqref{gMM-24GBN}, the third equality is a consequence of 
\eqref{NEWDVT.1-MMM} and the divergence theorem in the Lipschitz domain $\Omega_\ell$ for the vector field 
$\big(h_j\big|_{\Omega_\ell}\big)_{1\leq j\leq n}\in\big[C^\infty(\overline{\Omega_\ell})\big]^n$ 
(Theorem~\ref{banff-3} is more than adequate in this context), while the fourth 
equality is seen from \eqref{Utfv-MMM}. Formula \eqref{eq:GCxa-BIS-MMM} now follows 
by combining \eqref{NEWDVT.8-MMM} and \eqref{Utfv.2BIS-MMM}.
\end{proof}
%%%%%%%

Having dealt with the Dirichlet trace $\gamma_D$ earlier in this section, 
we now turn 
our attention to the task of defining the Neumann boundary trace operator $\gamma_N$
in the class of Lipschitz subdomains of Riemannian manifolds. As in the Euclidean setting, in a first stage
we shall introduce a weak version $\widetilde\gamma_N$ of the aforementioned Neumann boundary trace operator,
whose definition is inspired by the ``half" Green's formula for the Laplace--Beltrami operator.
Specifically, we make the following definition. 

%%%%%%%
\begin{definition}\label{h6r4d5UU-MMM}
Let $\Omega\subset M$ be a Lipschitz domain. For some fixed 
$s\in\big(\tfrac12,\tfrac32\big)$, the {\tt weak Neumann trace operator} 
is considered acting in the context  
\begin{equation}\label{2.88X-MMM}
\widetilde\gamma_N:\big\{(f,F)\in H^s(\Omega)\times H^{s-2}_0(\Omega) 
\,\big|\,\Delta_g f=F|_{\Omega}\text{ in }\mathcal{D}'(\Omega)\big\} 
\rightarrow H^{s-(3/2)}(\partial\Omega).
\end{equation}
Specifically, suppose that some function $f\in H^s(\Omega)$ along with some distribution 
$F\in H^{s-2}_0(\Omega)\subset H^{s-2}(M)$ 
satisfying $\Delta_g f=F|_{\Omega}$ in $\mathcal{D}'(\Omega)$ have been given. 
In particular, 
\begin{equation}\label{u543d-MMM}
{\rm grad}_g f\in H^{s-1}(\Omega,TM)=\big(H^{1-s}(\Omega,TM)\big)^*. 
\end{equation}
Then the manner in which $\widetilde\gamma_N(f,F)$ is now defined as a functional in the space
$H^{s-(3/2)}(\partial\Omega)=\big(H^{(3/2)-s}(\partial\Omega)\big)^*$ is as follows: 
Given $\phi\in H^{(3/2)-s}(\partial\Omega)$, then for any $\Phi\in H^{2-s}(\Omega)$ 
such that $\gamma_D\Phi=\phi$ {\rm (}whose existence is ensured by the surjectivity of 
\eqref{eqn:gammaDs.1-MMM}{\rm )}, set
\begin{align}\label{2.9NEW-MMM}
&\hskip -0.20in
{}_{H^{(3/2)-s}(\partial\Omega)}\big\langle\phi,
\widetilde\gamma_N(f,F)\big\rangle_{(H^{(3/2)-s}(\partial\Omega))^*}
\nonumber\\[2pt]
& \quad:={}_{H^{1-s}(\Omega,TM)}\big\langle{\rm grad}_g\Phi,
{\rm grad}_g f\big\rangle_{(H^{1-s}(\Omega,TM))^*}
\nonumber\\[2pt]
& \qquad+{}_{H^{2-s}(\Omega)}\big\langle\Phi,F\big\rangle_{(H^{2-s}(\Omega))^*}.
\end{align}
\end{definition}
%%%%%%%

Concerning Definition~\ref{h6r4d5UU-MMM} one observes that in the 
context described there, ${\rm grad}_g\Phi$ belongs to $H^{1-s}(\Omega,TM)$.
Utilizing \eqref{u543d-MMM}, this membership shows that the first pairing in the right-hand side 
of \eqref{2.9NEW-MMM} is meaningful. In addition, here we canonically identify the distribution $F$, 
originally belonging to $H^{s-2}_0(\Omega)$, with a functional in $(H^{2-s}(\Omega))^*$ 
(compare with the discussion pertaining to \eqref{eq:Redxax.2} in the Euclidean setting), 
so the last pairing in \eqref{2.9NEW-MMM} is also meaningfully defined as
\begin{align}\label{uyree-y5-MMM}
\begin{split}
& {}_{H^{2-s}(\Omega)}\big\langle\Phi,F\big\rangle_{(H^{2-s}(\Omega))^*}
={}_{H^{2-s}(M)}\big\langle\Theta,F\big\rangle_{H^{s-2}(M)}
\\[2pt]
& \quad\text{for any $\Theta\in H^{2-s}(M)$ satisfying $\Theta\big|_{\Omega}=\Phi$ 
in ${\mathcal{D}}'(\Omega)$}.
\end{split}
\end{align}

Here is a theorem which elaborates on the main properties of the weak Neumann trace operator defined above. 

%%%%%%%%%%%%%%%%%%%%%
\begin{theorem}\label{GNT-MMM}
Let $\Omega\subset M$ be a Lipschitz domain, and fix $s\in\big(\tfrac12,\tfrac32\big)$. 
Then the weak Neumann trace map $\widetilde\gamma_N$ 
from Definition~\ref{h6r4d5UU-MMM} yields an operator which is unambiguously defined, 
linear, and bounded {\rm (}assuming the space on the left-hand side of \eqref{2.88X-MMM} 
is equipped with the natural norm $(f,F)\mapsto\|f\|_{H^{s}(\Omega)}+\|F\|_{H^{s-2}(M)}${\rm )}.  
In addition, the following properties are true:  
\\[1mm] 
$(i)$ The weak Neumann trace operators corresponding to various values of the 
parameter $s\in\big(\tfrac12,\tfrac32\big)$ are compatible with one another and each 
of them is surjective. In fact, there exist linear and bounded operators
\begin{equation}\label{2.88X-NN-MMM}
\Upsilon_N:H^{s-(3/2)}(\partial\Omega)\to
\big\{u\in H^s(\Omega)\,\big|\,\Delta_g u\in L^2(\Omega)\big\},\quad s\in\big(\tfrac12,\tfrac32\big),
\end{equation}
which are compatible with one another and satisfy {\rm (}with tilde denoting the extension by zero 
outside $\Omega${\rm )} 
\begin{equation}\label{2.88X-NN2-MMM}
\widetilde\gamma_N\big(\Upsilon_N\psi,\widetilde{\Delta_g(\Upsilon_N\psi)}\,\big)=\psi,
\quad\forall\,\psi\in H^{s-(3/2)}(\partial\Omega)\,\text{ with }\,s\in\big(\tfrac12,\tfrac32\big).
\end{equation}
$(ii)$ Given any two pairs, 
\begin{align}\label{ju7653es-MMM}
\begin{split}
& \text{$(f_1,F_1)\in H^s(\Omega)\times H^{s-2}_0(\Omega)$ 
such that $\Delta_g f_1=F_1|_{\Omega}$ in $\mathcal{D}'(\Omega)$, and}
\\[2pt] 
& \quad\text{$(f_2,F_2)\in H^{2-s}(\Omega)\times H^{-s}_0(\Omega)$ such that $\Delta_g f_2=F_2|_{\Omega}$ 
in $\mathcal{D}'(\Omega)$},
\end{split}
\end{align}
the following Green's formula holds:
\begin{align}\label{GGGRRR-prim-MMM}
& {}_{H^{(3/2)-s}(\partial\Omega)}\big\langle\gamma_D f_2,\widetilde\gamma_N(f_1,F_1)
\big\rangle_{(H^{(3/2)-s}(\partial\Omega))^*}
\nonumber\\[2pt]
& \qquad-{}_{(H^{s-(1/2)}(\partial\Omega))^*}\big\langle\widetilde\gamma_N(f_2,F_2),\gamma_D f_1
\big\rangle_{H^{s-(1/2)}(\partial\Omega)}     
\nonumber\\[2pt]
&\quad={}_{H^{2-s}(\Omega)}\big\langle f_2,F_1\big\rangle_{(H^{2-s}(\Omega))^*}
-{}_{(H^s(\Omega))^*}\big\langle F_2,f_1\big\rangle_{H^s(\Omega)}.
\end{align} 
\end{theorem}
%%%%%%%
\begin{proof}
The proof follows along the lines of the proof of Theorem~\ref{GNT}, making use of the 
well-posedness results for the Neumann problem for the Laplace--Beltrami operator on 
Lipschitz subdomains of Riemannian manifolds from \cite{MT00}.
\end{proof}
%%%%%%%

We are prepared to state our main result concerning the Neumann boundary trace operator on Lipschitz 
subdomains of Riemannian manifolds in the theorem below. As in the case of the Dirichlet trace, 
by restricting ourselves to functions mapped by the Laplace--Beltrami operator into a 
better-than-expected Sobolev space, we are able to include the end-point cases $s=\tfrac12$ 
and $s=\tfrac32$ in \eqref{2.88X-MMM}. 

%%%%%%%%%%
\begin{theorem}\label{YTfdf.NNN.2-Main-MMM}
Assume that $\Omega\subset M$ is a Lipschitz domain. Then for each $\varepsilon>0$ the
weak Neumann boundary trace map, originally introduced in Definition~\ref{h6r4d5UU-MMM}, 
induces linear and continuous operators in the context  
\begin{equation}\label{eqn:gammaN-MMM}
\begin{array}{c}
\widetilde\gamma_N:\big\{(f,F)\in H^s(\Omega)\times H^{s-2+\varepsilon}_0(\Omega)\,|\,
\Delta_g f=F\big|_{\Omega}\,\text{ in }\,{\mathcal{D}}'(\Omega)\big\}\rightarrow H^{s-(3/2)}(\partial\Omega)
\\[4pt]
\text{with }\,\,s\in\big[\tfrac{1}{2},\tfrac{3}{2}\big]
\end{array}
\end{equation}
{\rm (}where the space on the left-hand side of \eqref{eqn:gammaN-MMM} is equipped with the natural norm 
$(f,F)\mapsto\|f\|_{H^{s}(\Omega)}+\|F\|_{H^{s-2+\varepsilon}(M)}${\rm )} which are compatible with 
those in Definition~\ref{h6r4d5UU-MMM} when $s\in\big(\tfrac{1}{2},\tfrac{3}{2}\big)$. 
Thus defined, the weak Neumann boundary trace map possesses the following properties: 
\\[1mm] 
$(i)$ The weak Neumann boundary trace map in \eqref{eqn:gammaN-MMM} is surjective. 
In fact, there exist linear and bounded operators
\begin{equation}\label{2.88X-NN-ii-MMM}
\Upsilon_N:H^{s-(3/2)}(\partial\Omega)\to
\big\{u\in H^s(\Omega)\,\big|\,\Delta_g u\in L^2(\Omega)\big\},\quad s\in\big[\tfrac12,\tfrac32\big],
\end{equation}
which are compatible with one another and satisfy {\rm (}with tilde denoting the extension by zero 
outside $\Omega${\rm )} 
\begin{equation}\label{2.88X-NN2-ii-MMM}
\widetilde\gamma_N\big(\Upsilon_N\psi,\widetilde{\Delta_g(\Upsilon_N\psi)}\,\big)=\psi,
\quad\forall\,\psi\in H^{s-(3/2)}(\partial\Omega)\,\text{ with }\,s\in\big[\tfrac12,\tfrac32\big].
\end{equation}
$(ii)$ If $\varepsilon\in(0,1)$ and $s\in\big[\tfrac{1}{2},\tfrac{3}{2}\big]$ then for any two pairs
\begin{align}\label{n7b66-MMM}
\begin{split}
& \text{$(f_1,F_1)\in H^s(\Omega)\times H^{s-2+\varepsilon}_0(\Omega)$ 
such that $\Delta_g f_1=F_1|_{\Omega}$ in $\mathcal{D}'(\Omega)$, and}
\\[2pt]
& \quad\text{$(f_2,F_2)\in H^{2-s}(\Omega)\times H^{-s+\varepsilon}_0(\Omega)$ 
such that $\Delta_g f_2=F_2|_{\Omega}$ in $\mathcal{D}'(\Omega)$}, 
\end{split}
\end{align}
the following Green's formula holds:
\begin{align}\label{GGGRRR-prim222-MMM}
& {}_{H^{(3/2)-s}(\partial\Omega)}\big\langle\gamma_D f_2,\widetilde\gamma_N(f_1,F_1)
\big\rangle_{(H^{(3/2)-s}(\partial\Omega))^*}
\nonumber\\[2pt]
& \qquad -{}_{(H^{s-(1/2)}(\partial\Omega))^*}\big\langle\widetilde\gamma_N(f_2,F_2),\gamma_D f_1
\big\rangle_{H^{s-(1/2)}(\partial\Omega)}     
\nonumber\\[2pt]
&\quad={}_{H^{2-s}(\Omega)}\big\langle f_2,F_1\big\rangle_{(H^{2-s}(\Omega))^*}
-{}_{(H^s(\Omega))^*}\big\langle F_2,f_1\big\rangle_{H^s(\Omega)}.
\end{align}
$(iii)$ There exists a constant $C\in(0,\infty)$ with the property that
\begin{align}\label{gafvv.6588-P-MMM}
\begin{split}
& \text{if $f\in H^{1/2}(\Omega)$ and $F\in H^{-(3/2)+\varepsilon}_0(\Omega)$ with 
$0<\varepsilon\leq 1$ satisfy} 
\\[2pt]
& \quad\text{$\Delta_g f=F\big|_{\Omega}$ in ${\mathcal{D}}'(\Omega)$ and 
$\widetilde\gamma_N(f,F)=0$, then $f\in H^{(1/2)+\varepsilon}(\Omega)$ }
\\[2pt]
& \quad\text{and $\|f\|_{H^{(1/2)+\varepsilon}(\Omega)}\leq 
C\big(\|f\|_{L^2(\Omega)}+\|F\|_{H^{-(3/2)+\varepsilon}(M)}\big)$.}
\end{split}
\end{align} 
\end{theorem}
%%%%%%%%%%
\begin{proof}
In the case when $s\in\big(\tfrac{1}{2},\tfrac{3}{2}\big)$, all desired conclusions follow from 
Theorem~\ref{GNT-MMM} simply by observing that  
\begin{equation}\label{5rfc-WACO}
\big\{(f,F)\in H^s(\Omega)\times H^{s-2+\varepsilon}_0(\Omega)\,\big|\,
\Delta_g f=F\big|_{\Omega}\,\text{ in }\,{\mathcal{D}}'(\Omega)\big\}, 
\end{equation} 
the domain of the weak Neumann trace operator $\widetilde\gamma_N$ in \eqref{eqn:gammaN-MMM}, is a subspace of 
\begin{equation}\label{hafvv-WACO}
\big\{(f,F)\in H^s(\Omega)\times H^{s-2}_0(\Omega)\,\big|\,
\Delta_g f=F\big|_{\Omega}\,\text{ in }\,{\mathcal{D}}'(\Omega)\big\}, 
\end{equation}
the domain of $\widetilde\gamma_N$ in \eqref{2.88X-MMM}. In this context one can employ the operators 
$\Upsilon_N$ in \eqref{2.88X-NN-MMM}. 

Next, we consider the case when $s=\tfrac{3}{2}$. 
For the goals we have in mind, there is no loss of generality in assuming that $\varepsilon\in(0,1)$.
Suppose some $f\in H^{3/2}(\Omega)$ along with some $F\in H^{-(1/2)+\varepsilon}_0(\Omega)$ satisfying 
$\Delta_g f=F\big|_{\Omega}$ in ${\mathcal{D}}'(\Omega)$ have been given. In particular, 
\begin{equation}\label{j6f432-MMM}
{\rm grad}_gf\in H^{1/2}(\Omega,TM)\,\text{ and }\,\Delta_g f\in H^{-(1/2)+\varepsilon}(\Omega).
\end{equation}
In addition, for each $X\in C^\infty(M,TM)$, the function $\nabla_X f\in H^{1/2}(\Omega)$ satisfies 
\begin{align}\label{j6f433-MMM}
\Delta_g(\nabla_X f) &=\big[\Delta_g,\nabla_X\big]f+\nabla_X(\Delta_g f)
\nonumber\\[2pt]
&=\big[\Delta_g,\nabla_X\big]f+\nabla_X\big(F\big|_{\Omega}\big)
\nonumber\\[2pt]
&=\big[\Delta_g,\nabla_X\big]f+(\nabla_X F)\big|_{\Omega}\in H^{-(3/2)+\varepsilon}(\Omega),
\end{align}
since the commutator $\big[\Delta_g,\nabla_X\big]$ is a second-order differential expression. 
Moreover, there is a naturally accompanying estimate to the effect that for each vector field $X\in C^\infty(M,TM)$
there exists $C\in(0,\infty)$ independent of $f$ and $F$ such that  
\begin{align}\label{j6f433-hREE-MMM}
\|\Delta_g(\nabla_X f)\|_{H^{-(3/2)+\varepsilon}(\Omega)}\leq C\left(
\|f\|_{H^{3/2}(\Omega)}+\|F\|_{H^{-(1/2)+\varepsilon}(\Omega)}\right).
\end{align}
From \eqref{j6f432-MMM} and \eqref{j6f433-MMM} one concludes that
\begin{equation}\label{lk8ABba-MMM}
{\rm grad}_gf\in H^{1/2,-(3/2)+\varepsilon}_{\Delta_g}(\Omega,TM).
\end{equation}
In turn, from \eqref{lk8ABba-MMM} and Definition~\ref{DDD-2} (used with $s=1/2$) one infers that
\begin{equation}\label{eq:Gav7gt5-MMM}
\text{$\gamma_D({\rm grad}_gf)$ exists in $L^2(\partial\Omega,TM)$}.
\end{equation}
Moreover, \eqref{j6f433-hREE-MMM} implies that 
\begin{align}\label{j6f433-hREE.222-MMM}
\|\gamma_D({\rm grad}_gf)\|_{L^2(\partial\Omega,TM)}\leq C\left(
\|f\|_{H^{3/2}(\Omega)}+\|F\|_{H^{-(1/2)+\varepsilon}(\Omega)}\right).
\end{align}

To proceed further, pick an arbitrary $\Phi\in C^\infty(\overline{\Omega})$, set 
$\phi:=\Phi\big|_{\partial\Omega}$, and consider the vector field 
\begin{equation}\label{elhgtf-MMM}
\vec{F}:=\overline{\Phi}\,{\rm grad}_g f\,\text{ in }\,\Omega.
\end{equation}
In light of the manifold counterpart of \eqref{eq:DDEj6g5}, the above definition implies  
\begin{equation}\label{lk8AAA-MMM}
\vec{F}\in H^{1/2}(\Omega,TM).
\end{equation}
Moreover, from \eqref{elhgtf-MMM}, \eqref{UInabb}, \eqref{rdf96-2DC.1RR}, \eqref{eqn.thu-TELL}, 
and \eqref{j6f432-MMM} one infers that
\begin{equation}\label{lk8ht5-MMM}
{\rm div}_g\vec{F}=\big\langle\overline{{\rm grad}_g\Phi},{\rm grad}_gf\big\rangle_{TM}
+\overline{\Phi}\Delta_g f\in H^{-(1/2)+\varepsilon}(\Omega).
\end{equation}
Moreover, locally, 
\begin{equation}\label{elhgtf-MMM-2}
\vec{F}=F_j\partial_j,\,\text{ with }\,F_j=\overline{\Phi}g^{jk}\partial_k f,
\end{equation}
and for each $j$ one has locally, 
\begin{align}\label{eq:NBb65a-MMM}
\Delta_g F_j &=(\partial_k f)\Delta_g(\overline{\Phi}g^{jk})
+\overline{\Phi}g^{jk}\Delta_g(\partial_k f)
\nonumber\\[2pt]
&\quad +2\big\langle\,{\rm grad}_g(\overline{\Phi}g^{jk}),{\rm grad}_g(\partial_k f)\big\rangle_{TM}.
\end{align}
From \eqref{lk8AAA-MMM}, \eqref{eq:NBb65a-MMM}, and \eqref{j6f433-MMM} one concludes that
\begin{equation}\label{lk8AAA-MaGG-MMM}
\vec{F}\in H^{1/2,-(3/2)+\varepsilon}_{\Delta_g}(\Omega,TM).
\end{equation}

Given \eqref{lk8AAA-MaGG-MMM}, Theorem~\ref{Ygav-75-BIS-MMM} applies to the vector field $\vec{F}$. 
Specifically, let $\nu$ and $\sigma_g$ denote, respectively, the outward unit normal and 
surface measure on $\partial\Omega$. Then, with the Dirichlet trace $\gamma_D({\rm grad}_gf)$ 
understood in the sense of \eqref{eq:Gav7gt5-MMM}, one has 
\begin{align}\label{2.9MMa-ii-MMM}
& \big(\phi\,,\,\langle\nu,\gamma_D({\rm grad}_gf)\rangle_{TM}\big)_{L^2(\partial\Omega)}
=\int_{\partial\Omega}\langle\nu,\gamma_D\vec{F}\rangle_{TM}\,d\sigma_g
\nonumber\\[2pt]
& \quad={}_{H^{(1/2)-\varepsilon}(\Omega)}\big\langle{\bf 1},{\rm div}_g\vec{F}
\big\rangle_{H^{-(1/2)+\varepsilon}(\Omega)}
\nonumber\\[2pt]
& \quad={}_{H^{(1/2)-\varepsilon}(\Omega)}\Big\langle{\bf 1},
\big\langle\overline{{\rm grad}_g\Phi},{\rm grad}_gf\big\rangle_{TM}
\Big\rangle_{H^{-(1/2)+\varepsilon}(\Omega)}
\nonumber\\[2pt]
& \qquad+{}_{H^{(1/2)-\varepsilon}(\Omega)}\big\langle{\bf 1},\overline{\Phi}\Delta_g f
\big\rangle_{H^{-(1/2)+\varepsilon}(\Omega)}
\nonumber\\[2pt]
& \quad={}_{H^{(1/2)-\varepsilon}(\Omega)}\big\langle
{\rm grad}_g\Phi,{\rm grad}_gf\big\rangle_{H^{-(1/2)+\varepsilon}(\Omega)}
\nonumber\\[2pt]
& \qquad+{}_{H^{(1/2)-\varepsilon}(\Omega)}\big\langle\Phi,\Delta_g f
\big\rangle_{H^{-(1/2)+\varepsilon}(\Omega)}
\nonumber\\[2pt]
& \quad=\big({\rm grad}_g\Phi,{\rm grad}_gf\big)_{L^2(\Omega)}
+{}_{H^{(1/2)-\varepsilon}(\Omega)}\big\langle\Phi,F
\big\rangle_{(H^{(1/2)-\varepsilon}(\Omega))^*}, 
\end{align}
where the last step relies on the manner in which 
$(H^{(1/2)-\varepsilon}(\Omega))^*$ is identified with $H^{-(1/2)+\varepsilon}(\Omega)$
(see \eqref{u5iiui}--\eqref{u5iiui-bbb} for the Euclidean setting). 

The fact that $f\in H^{3/2}(\Omega)$ entails $f\in H^{s}(\Omega)$ for each 
$s\in\big(\tfrac{1}{2},\tfrac{3}{2}\big)$ and, as such, a direct comparison of 
\eqref{2.9MMa-ii-MMM} and \eqref{2.9NEW-MMM} reveals that 
\begin{align}\label{2.9MMabn-P-MMM}
\begin{split}
& {}_{H^{(3/2)-s}(\partial\Omega)}\big\langle\phi,
\widetilde\gamma_N(f,F)\big\rangle_{(H^{(3/2)-s}(\partial\Omega))^*}
=\big(\phi\,,\,\langle\nu,\gamma_D({\rm grad}_gf)\rangle_{TM}\big)_{L^2(\partial\Omega)} 
\\[2pt]
& \quad\text{for every $s\in\big(\tfrac{1}{2},\tfrac{3}{2}\big)$ and every function
$\phi\in\big\{\Phi\big|_{\partial\Omega}\,\big|\,\,\Phi\in C^\infty(\overline{\Omega})\big\}$.}
\end{split}
\end{align}
Since the latter space is dense in $L^2(\partial\Omega)$, this ultimately proves that 
\begin{align}\label{eq:Nnan7yg-P-MMM}
& \text{if $f\in H^{3/2}(\Omega)$ and $F\in H^{-(1/2)+\varepsilon}_0(\Omega)$ 
for some $\varepsilon\in(0,1)$ satisfy} 
\nonumber\\[2pt]
& \quad\text{$\Delta_g f=F\big|_{\Omega}$ in ${\mathcal{D}}'(\Omega)$, 
then actually $\widetilde\gamma_N(f,F)\in L^2(\partial\Omega)$ and,} 
\\[2pt]
& \quad\text{$\widetilde\gamma_N(f,F)=\langle\nu,\gamma_D({\rm grad}_gf)\rangle_{TM}$
with the Dirichlet trace as in \eqref{eq:Gav7gt5-MMM}.}  
\nonumber
\end{align}
Moreover, from \eqref{j6f433-hREE.222-MMM} one infers that
\begin{align}\label{ut543i-Ai-MMM}
\big\|\widetilde\gamma_N(f,F)\big\|_{L^2(\partial\Omega)} 
\leq C\big(\|f\|_{H^{3/2}(\Omega)}+\|F\|_{H^{-(1/2)+\varepsilon}(\Omega)}\big)
\end{align}
for some constant $C\in(0,\infty)$, independent of $(f,F)$. 

At this stage, all remaining claims in the statement of the current theorem may be justified 
based on what we have proved already by reasoning along the lines of the proof of Theorem~\ref{YTfdf.NNN.2-Main},
with natural alterations. The well-posedness results for boundary value problems for the 
Laplace--Beltrami operator on Lipschitz subdomains of Riemannian manifolds which are relevant 
for us here are available from the work in \cite{MT99} and \cite{MT00}.
\end{proof}
%%%%%%%

The following special case of Theorem~\ref{YTfdf.NNN.2-Main-MMM} plays a significant role
in applications.

%%%%%%%%%%
\begin{corollary}\label{YTfdf.NNN.2-MMM}
Assume that $\Omega\subset M$ is a Lipschitz domain, and denote by $\nu$ its outward unit normal.
Then the Neumann trace map, originally defined for each for $u\in C^\infty(\overline{\Omega})$
as $u\mapsto\langle\nu,{\rm grad}_gu\rangle_{TM}$ on $\partial\Omega$, extends uniquely to linear 
continuous operators  
\begin{equation}\label{eqn:gammaN-pp-MMM}
\gamma_N:\big\{u\in H^s(\Omega)\,\big|\,\Delta_g u\in L^2(\Omega)\big\}\to H^{s-(3/2)}(\partial\Omega), 
\quad s\in\big[\tfrac{1}{2},\tfrac{3}{2}\big] 
\end{equation}
{\rm (}throughout, the space on the left-hand side of \eqref{eqn:gammaN-pp-MMM} equipped with the natural graph 
norm $u\mapsto\|u\|_{H^{s}(\Omega)}+\|\Delta_g u\|_{L^2(\Omega)}${\rm )}, that are compatible with one another.
In addition, the following properties are true:

\begin{enumerate}
\item[(i)] The Neumann trace map \eqref{eqn:gammaN-pp-MMM} is surjective. In fact, there exist linear and bounded operators
\begin{equation}\label{2.88X-NN-ii-RR-MMM}
\Upsilon_N:H^{s-(3/2)}(\partial\Omega)\to
\big\{u\in H^s(\Omega)\,\big|\,\Delta_g u\in L^2(\Omega)\big\},\quad s\in\big[\tfrac12,\tfrac32\big],
\end{equation}
which are compatible with one another and are right-inverses for the Neumann trace, that is, 
\begin{equation}\label{2.88X-NN2-ii-RR-MMM}
\gamma_N(\Upsilon_N\psi)=\psi,
\quad\forall\,\psi\in H^{s-(3/2)}(\partial\Omega)\,\text{ with }\,s\in\big[\tfrac12,\tfrac32\big].
\end{equation}

\item[(ii)] If $s\in\big[\tfrac{1}{2},\tfrac{3}{2}\big]$, then for any functions $f\in H^s(\Omega)$ with 
$\Delta_g f\in L^2(\Omega)$ and $h\in H^{2-s}(\Omega)$ with $\Delta_g h\in L^2(\Omega)$ the following 
Green's formula holds:
\begin{align}\label{GGGRRR-MMM}
& {}_{H^{(3/2)-s}(\partial\Omega)}\big\langle\gamma_D h,\gamma_N f
\big\rangle_{(H^{(3/2)-s}(\partial\Omega))^*}
\nonumber\\[2pt]
& \qquad-{}_{(H^{s-(1/2)}(\partial\Omega))^*}\big\langle\gamma_N h,\gamma_D f
\big\rangle_{H^{s-(1/2)}(\partial\Omega)}     
\nonumber\\[2pt]
& \quad=(h,\Delta_g f)_{L^2(\Omega)}-(\Delta_g h,f)_{L^2(\Omega)}.
\end{align}

\item[(iii)] For each $s\in\big[\tfrac{1}{2},\tfrac{3}{2}\big]$, the null space of 
the Neumann boundary trace operator \eqref{eqn:gammaN-pp-MMM} satisfies
\begin{equation}\label{eq:EFFa.111-MMM}
{\rm ker}(\gamma_N)\subseteq H^{3/2}(\Omega).
\end{equation}
In fact, the inclusion in \eqref{eq:EFFa.111-MMM} is quantitative in the sense that
there exists a constant $C\in(0,\infty)$ with the property that
\begin{align}\label{gafvv.6588-MMM}
\begin{split}
& \text{whenever $u\in H^{1/2}(\Omega)$ satisfies $\Delta_g u\in L^2(\Omega)$ and $\gamma_N u=0$, then}
\\[2pt]
& \quad\text{$u\in H^{3/2}(\Omega)$ and 
$\|u\|_{H^{3/2}(\Omega)}\leq C\big(\|u\|_{L^2(\Omega)}+\|\Delta_g u\|_{L^2(\Omega)}\big)$.}
\end{split}
\end{align}
\end{enumerate}
\end{corollary}
%%%%%%%%%%
\begin{proof}
The idea is to produce a formula restricting the weak Neumann trace operator from Theorem~\ref{YTfdf.NNN.2-Main-MMM} to the present 
setting. With this goal 
in mind, we assume an $s\in\big[\tfrac{1}{2},\tfrac{3}{2}\big]$ has been fixed 
and choose $0<\varepsilon<\min\{1,2-\varepsilon\}$. Next, we denote by 
\begin{align}\label{2.10X-MMM}
\begin{split} 
&\iota:\big\{u\in H^s(\Omega)\,\big|\,\Delta_g u\in L^2(\Omega)\big\}       
\\[2pt]
&\quad\;\rightarrow
\big\{(f,F)\in H^s(\Omega)\times H^{s-2+\varepsilon}_0(\Omega)\,\big|\,\Delta_g f=F\big|_{\Omega}
\,\text{ in }\,{\mathcal{D}}'(\Omega)\big\}, 
\end{split} 
\end{align}
the continuous injection given by 
\begin{equation}\label{2.11X-MMM}
\iota(u):=\big(u,\widetilde{\Delta_g u}\big),\quad\forall\,u\in H^s(\Omega)\,\text{ with }\,\Delta_g u\in L^2(\Omega),
\end{equation}
where, as usual, tilde denotes the extension by zero outside $\Omega$. We then define 
\begin{equation}\label{2.12X-MMM}
\gamma_N:=\widetilde\gamma_N\circ\iota
\end{equation}
and note that this is a well defined, linear, and bounded mapping in the context of \eqref{eqn:gammaN-pp-MMM}.
With this in hand, all other claims in the statement are established as in the proof of 
Corollary~\ref{YTfdf.NNN.2}. 
\end{proof}
%%%%%%%

To exemplify the manner in which the mapping $\gamma_N$ introduced in \eqref{2.12X-MMM} operates, we 
consider the case where $s\in\big(\tfrac{1}{2},\tfrac{3}{2}\big)$. Given $u\in H^s(\Omega)$ 
with $\Delta_g u\in L^2(\Omega)$, along with $\phi\in H^{(3/2)-s}(\partial\Omega)$
and $\Phi\in H^{2-s}(\Omega)$ such that $\gamma_D\Phi=\phi$, then the action of 
$\gamma_N u\in H^{s-(3/2)}(\partial\Omega)=\big(H^{(3/2)-s}(\partial\Omega)\big)^*$ on 
$\phi\in H^{(3/2)-s}(\partial\Omega)$ is concretely given by 
\begin{align}\label{2.9-Mi}
& {}_{H^{(3/2)-s}(\partial\Omega)}\big\langle\phi,\gamma_N u\big\rangle_{(H^{(3/2)-s}(\partial\Omega))^*}
\nonumber\\[2pt]
& \quad={}_{H^{(3/2)-s}(\partial\Omega)}\big\langle\phi,
\widetilde\gamma_N(u,\widetilde{\Delta_g u})\big\rangle_{(H^{(3/2)-s}(\partial\Omega))^*}
\nonumber\\[2pt]
& \quad={}_{H^{1-s}(\Omega,TM)}\big\langle{\rm grad}_g\Phi,{\rm grad}_g f\big\rangle_{(H^{1-s}(\Omega,TM))^*}
\nonumber\\[2pt]
& \qquad+{}_{H^{2-s}(\Omega)}\big\langle\Phi,\widetilde{\Delta_g u}\big\rangle_{(H^{2-s}(\Omega))^*}
\nonumber\\[2pt]
& \quad={}_{H^{1-s}(\Omega,TM)}\big\langle{\rm grad}_g\Phi,{\rm grad}_g f\big\rangle_{(H^{1-s}(\Omega,TM))^*}
\nonumber\\[2pt]
& \qquad+(\Phi,\Delta_g u)_{L^2(\Omega)}.
\end{align}

%%%%%%%%%%%%%%%%
\subsection{Schr\"odinger operators on Lipschitz subdomains of a Riemannian manifold}\label{ss.MM.3}
%%%%%%%%%%%%%%%%
The goal here is to study Schr\"odinger operators $L$ on Lipschitz subdomains of the 
compact Riemannian manifold $M$. To set the stage, given a Lipschitz domain $\Omega\subset M$ 
and an essentially bounded real-valued potential $V$, we first introduce the sesquilinear form
\begin{equation}\label{afform-Mi}
\mathfrak l_{F,\Omega}(f,h):=\big({\rm grad}_g f,{\rm grad}_g h\big)_{L^2(\Omega,TM)}+(f,Vh)_{L^2(\Omega)}, 
\quad\dom(\mathfrak l_{F,\Omega}):=\accentset{\circ}{H}^1(\Omega), 
\end{equation}
which is densely defined, closed, symmetric, and semibounded from below in $L^2(\Omega)$. Hence, it 
follows from the First Representation Theorem (cf. \cite[Theorem~VI.2.1]{Ka80}) that there is a unique 
self-adjoint operator $L_{F,\Omega}$ in $L^2(\Omega)$ such that the identity
\begin{equation}\label{formopaf-Mi}
\mathfrak l_{F,\Omega}(f,h)=\big(f,L_{F,\Omega}h\big)_{L^2(\Omega)}
\end{equation}
holds for all $f\in\dom(\mathfrak l_{F,\Omega})=\accentset{\circ}{H}^1(\Omega)$ and all 
$h\in\dom(L_{F,\Omega})\subset\dom(\mathfrak l_{F,\Omega})$. Making use of \eqref{tr6ytG-MMM} 
and Green's formula it follows that
\begin{equation}\label{OUYTf-Mi}
L_{F,\Omega}=-\Delta_g+V,\quad\dom(L_{F,\Omega})=\big\{f\in\accentset{\circ}{H}^1(\Omega)
\,\big|\,\Delta_g f\in L^2(\Omega)\big\},
\end{equation}
and hence $L_{F,\Omega}$ is a self-adjoint extension of the minimal realization 
$L_{min,\Omega}$ of $-\Delta_g+V$ defined in \eqref{eqn:Amin1-MMM}. 
Again, by \cite[Subsection~VI.2.3]{Ka80}, $L_{F,\Omega}$ represents the Friedrichs extension of 
$L_{min,\Omega}$. Abstract functional theoretic results (cf., e.g., \cite[Section~6.1]{EE18}) 
then yield the following theorem. 

%%%%%%%%
\begin{theorem}\label{tfried-Mi}
For a Lipschitz domain $\Omega\subset M$, the Friedrichs extension $L_{F,\Omega}$ of $L_{min,\Omega}$ 
is a self-adjoint operator in $L^2(\Omega)$, whose resolvent is compact, and whose spectrum is 
purely discrete and contained in $(v_-,\infty)$ {\rm (}where $v_{-}$ is as in \eqref{essinfv}{\rm )}.
In particular, $\sigma_{ess}(L_{F,\Om})=\varnothing$. 
\end{theorem}
%%%%%%%

Our next goal is to study the Dirichlet and Neumann realizations of $-\Delta_g+V$ 
on a Lipschitz subdomain $\Omega$ of the compact manifold $M$. Assuming, as before, that  
$V$ is an essentially bounded real-valued potential, it follows from \eqref{eq:WH} and \eqref{incl-Yb.EE} with 
$s=1$ that $\dom(\mathfrak l_{F,\Omega})=\accentset{\circ}{H}^1(\Omega)$ and the Friedrichs extension 
$L_{F,\Omega}$ coincides with the self-adjoint Dirichlet operator 
\begin{equation}\label{We-Q.10EE-Mi}
\begin{split}
& \,L_{D,\Omega}=-\Delta_g+V, 
\\[2pt]
& \dom(L_{D,\Omega})=\big\{f\in H^{1}(\Omega)\cap\dom(L_{max,\Omega})
\,\big|\,\gamma_D f=0\big\}.
\end{split}
\end{equation}
Our next theorem collects further useful properties of this operator. 

%%%%%%%%%
\begin{theorem}\label{t4.4-Mi}
Assume $\Omega\subset M$ is a bounded Lipschitz domain, and pick some $V\in L^\infty(M)$.
In this setting, let $L_{D,\Omega}$ be the Dirichlet realization of $-\Delta_g+V$ introduced 
in \eqref{We-Q.10EE-Mi}. Then $\dom(L_{D,\Omega})\subset H^{3/2}(\Omega)$, hence
\begin{equation}\label{We-Q.10EE-jussi-Mi}
\begin{split}
& \,L_{D,\Omega}=-\Delta_g+V, 
\\[2pt]
&\dom(L_{D,\Omega})=\big\{f\in H^{3/2}(\Omega)\cap\dom(L_{max,\Omega})\,\big|\,\gamma_D f=0\big\}.
\end{split}
\end{equation}
In addition, on $\dom(L_{D,\Omega})$ the norms 
\begin{equation}\label{eq:UUjh-jussi-Mi}
f\mapsto\|f\|_{H^s(\Omega)}+\|\Delta_g f\|_{L^2(\Omega)},\quad s\in\big[0,\tfrac{3}{2}\big],
\end{equation}
are equivalent. Furthermore, $L_{D,\Omega}$ is self-adjoint in $L^2(\Omega)$, with compact 
resolvent, and purely discrete spectrum contained in $(v_-,\infty)$. In particular, $\sigma_{ess}(L_{D,\Om})=\varnothing$. 
Moreover,
\begin{equation}\label{2.4Hba-Mi}
\dom\big(|L_{D,\Om}|^{1/2}\big)=\accentset{\circ}{H}^1(\Om).   
\end{equation}
\end{theorem}
%%%%%%%%%
\begin{proof} 
That functions in $\dom(L_{D,\Omega})$ exhibit $H^{3/2}$-regularity is a consequence  
of \eqref{eq:EFFa.NNN-MMM} (used with $s=1$). Together with \eqref{We-Q.10EE-Mi} this also 
proves \eqref{We-Q.10EE-jussi-Mi}. When $s\in [1,\tfrac{3}{2}]$ the claim in \eqref{eq:UUjh-jussi-Mi}
is implied by \eqref{gafvv.6577-MMM}, while for $s\in [0,1]$ one reasons as follows. 
Given $f\in\dom(L_{D,\Omega})$, from \eqref{2.9-Mi} written for $\Phi:=f$, $F:=\Delta_g f$, 
and $s=1$, one obtains  
\begin{equation}\label{greendirichlet-Mi}
\begin{split}
0 &={}_{H^{1/2}(\partial\Omega)}\big\langle\gamma_D f,\gamma_N f\big\rangle_{H^{-1/2}(\partial\Omega)}
\\[2pt]
&=\big({\rm grad}_g f,{\rm grad}_g f)_{L^2(\Omega,TM)}+(f,\Delta_g f)_{L^2(\Omega)},
\end{split}
\end{equation}
which further implies that for all $f\in\dom(L_{D,\Omega})$, 
\begin{align}\label{u76gVV-Mi}
\|{\rm grad}_g f\|^2_{L^2(\Omega,TM)} &\leq\|f\|_{L^2(\Omega)}\,\|\Delta_g f\|_{L^2(\Omega)}
\nonumber\\[2pt]
&\leq\big(\|f\|_{L^2(\Omega)}+\|\Delta_g f\|_{L^2(\Omega)}\big)^2.
\end{align}
Thus, $\|f\|_{H^1(\Omega)}\leq C(\|f\|_{L^2(\Omega)}+\|\Delta_g f\|_{L^2(\Omega)})$ for all
$f\in\dom(L_{D,\Omega})$ which establishes \eqref{eq:UUjh-jussi-Mi} for $s\in [0,1]$. 
The Second Representation Theorem (see \cite[Theorem~VI.2.23]{Ka80}) gives \eqref{2.4Hba-Mi}, and 
the remaining claims in the statement of the theorem are consequences of Theorem~\ref{tfried-Mi}. 
\end{proof}
%%%%%%%%%

Next, introduce the sesquilinear form
\begin{align}\label{formN-Mi} 
\begin{split} 
& \mathfrak l_{N,\Omega}(f,h):=\big({\rm grad}_g f,{\rm grad}_g h\big)_{L^2(\Omega,TM)}+(f,Vh)_{L^2(\Omega)},
\\[2pt] 
& \dom(\mathfrak l_{N,\Omega}):=H^1(\Omega), 
\end{split}
\end{align}
which is densely defined, closed, symmetric, and semibounded from below in $L^2(\Omega)$. 
One notes that $\mathfrak l_{N,\Omega}$ is an extension of the form $\mathfrak l_{F,\Omega}$ 
in \eqref{afform-Mi} since
\begin{equation}\label{formN-dom-Mi}
\dom(\mathfrak l_{F,\Omega})=\accentset{\circ}{H}^1(\Omega)\subset H^1(\Omega)
=\dom(\mathfrak l_{N,\Omega}). 
\end{equation}
Once again, the first representation theorem \cite[Theorem~VI.2.1]{Ka80} implies that 
there is a unique self-adjoint operator $L_{N,\Omega}$ in $L^2(\Omega)$ such that the identity
\begin{equation}\label{formopan-Mi}
\mathfrak l_{N,\Omega}(f,h)=\big(f,L_{N,\Omega}h\big)_{L^2(\Omega)}
\end{equation}
is valid for all $f\in\dom(\mathfrak l_{N,\Omega})=H^1(\Omega)$ and all 
$h\in\dom(L_{N,\Omega})\subset\dom(\mathfrak l_{N,\Omega})$. Having fixed such $f,h$, 
one makes use of \eqref{formN-Mi}, \eqref{formopan-Mi}, and \eqref{2.9-Mi} 
(written for $\Phi:=f$, $f:=h$, $F:=\Delta_g h$, and $s=1$) in order to obtain
\begin{align}\label{neumannjussi-Mi}
(f,L_{N,\Omega}h)_{L^2(\Omega)} &=\big(f,(-\Delta_g+V)h\big)_{L^2(\Omega)}
\nonumber\\[2pt]
&\quad+{}_{H^{1/2}(\partial\Omega)}\big\langle\gamma_D f,\gamma_N h\big\rangle_{H^{-1/2}(\partial\Omega)}
\end{align}
for all $h\in\dom(L_{N,\Omega})$ and all $f\in H^1(\Omega)$. First restricting $f\in\accentset{\circ}{H}^1(\Omega)$ 
in \eqref{neumannjussi-Mi} then implies that $L_{N,\Omega}=-\Delta_g+V$. Next, taking into 
account that the range of $\gamma_D$ acting from $\dom(\mathfrak l_{N,\Omega})=H^1(\Omega)$ 
equals $H^{1/2}(\partial\Omega)$, which in turn is a dense subspace of $L^2(\partial\Omega)$
(cf. \eqref{eqn:gammaDs.2aux-MMM} with $s=\varepsilon=1$), one infers that \eqref{neumannjussi-Mi} 
forces $\gamma_N h=0$ for each $h\in\dom(L_{N,\Omega})$. Altogether, this proves that
\begin{equation}\label{We-Q.10EENN-Mi}
\begin{split}
& L_{N,\Omega}=-\Delta_g+V, 
\\[2pt]
&\dom(L_{N,\Omega})=\big\{f\in H^{1}(\Omega)\cap\dom(L_{max,\Omega})\,\big|\,\gamma_N f=0\big\}.
\end{split}
\end{equation}
Hence, $L_{N,\Omega}$ is a self-adjoint extension of the minimal realization 
$L_{min,\Omega}$ of $-\Delta_g+V$ defined in \eqref{eqn:Amin1-MMM}. Henceforth 
we shall refer to $L_{N,\Omega}$ as the Neumann extension (or Neumann realization) 
of $L_{min,\Omega}$. Our next theorem contains further properties of this Neumann realization. 

%%%%%%%%%%
\begin{theorem}\label{t4.5-Mi} 
Assume $\Omega\subset M$ is a bounded Lipschitz domain, and pick a potential $V\in L^\infty(M)$.
In this context, let $L_{N,\Omega}$ be the Neumann realization of $-\Delta_g+V$ defined as in 
\eqref{We-Q.10EENN-Mi}. Then $\dom(L_{N,\Omega})\subset H^{3/2}(\Omega)$, hence
\begin{equation}\label{We-Q.10EE-jussiNN-Mi}
\begin{split}
& L_{N,\Omega}=-\Delta_g+V, 
\\[2pt]
&\dom(L_{N,\Omega})=\big\{f\in H^{3/2}(\Omega)\cap\dom(L_{max,\Omega})\,\big|\,\gamma_N f=0\big\}.
\end{split}
\end{equation}
Moreover, on $\dom(L_{N,\Omega})$ the norms 
\begin{equation}\label{eq:UUjh-jussiNN-Mi}
f\mapsto\|f\|_{H^s(\Omega)}+\|\Delta_g f\|_{L^2(\Omega)},\quad s\in\big[0,\tfrac{3}{2}\big],
\end{equation}
are equivalent. In addition, $L_{N,\Omega}$ is self-adjoint in $L^2(\Omega)$, with compact 
resolvent, and purely discrete spectrum, contained in $[v_-,\infty)$. In particular, 
$\sigma_{ess}(L_{N,\Om})=\varnothing$. Moreover,
\begin{equation}\label{2.4HyN-Mi}
\dom\big(|L_{N,\Om}|^{1/2}\big)=H^1(\Omega).   
\end{equation}
\end{theorem}
%%%%%%%%%%
\begin{proof}
That $\dom(L_{N,\Omega})$ is contained in $H^{3/2}(\Omega)$ is seen from \eqref{eq:EFFa.111-MMM} 
(used with $s=1$), while the claim in \eqref{eq:UUjh-jussiNN-Mi} a direct consequence of \eqref{gafvv.6588-MMM}. 
All other claims may be justified by reasoning as in the proofs of \cite[Theorem~2.6]{GM08} and 
\cite[Theorem 4.5]{GM09b}. Here we just remark that the spectrum of $L_{N,\Omega}$ is bounded from 
below by $v_-$ since the corresponding form $\mathfrak l_{N,\Omega}$ in \eqref{formN-Mi} is bounded from below by $v_-$.
\end{proof}
%%%%%%%%%%

We continue by describing the domain of the minimal operator $L_{min,\Omega}$.

%%%%%%%%%%
\begin{lemma}\label{l4.3-Mi}
Assume that $\Omega\subset M$ is a bounded Lipschitz domain, and suppose that $V\in L^\infty(M)$.
Then the closed symmetric operator $L_{min,\Omega}$ is given by
\begin{equation}\label{eqn:Amin2-Mi}
L_{min,\Omega}=-\Delta+V,\quad\dom(L_{min,\Omega})=\accentset{\circ}{H}^2(\Omega).
\end{equation}
\end{lemma}
%%%%%%%%%
\begin{proof}
This is an immediate consequence of Lemma~\ref{l3.2-MMM} and \eqref{eq:WH}.
\end{proof}
%%%%%%%

Our last result shows that, as in the Euclidean setting, the operators 
$L_{D,\Omega}$ and $L_{N,\Omega}$ are relatively prime.

%%%%%%%%%%
\begin{theorem}\label{t4.6-Mi}
Assume that $\Omega\subset M$ is a bounded Lipschitz domain, and suppose that $V\in L^\infty(M)$.
Then the operators $L_{D,\Omega}$ and $L_{N,\Omega}$ are relatively prime, that is,
\begin{equation}\label{3.30-Mi}
\dom(L_{D,\Omega})\cap\dom(L_{N,\Omega})=\dom(L_{min,\Omega})=\accentset{\circ}{H}^2(\Omega).
\end{equation}
\end{theorem}
%%%%%%%%%%
\begin{proof} 
Given any $f\in\dom(L_{D,\Omega})\cap\dom(L_{N,\Omega})$, \eqref{We-Q.10EE-jussi-Mi} 
and \eqref{We-Q.10EE-jussiNN-Mi} imply that $f\in H^{3/2}(\Omega)$ and $\gamma_D f=\gamma_N f=0$. 
Together with \eqref{GGGRRR-MMM}, these conditions ensure that for every 
$\psi\in C^\infty(\overline{\Omega})$ one may write
\begin{equation}\label{uartUTR-Mi}
(f,\Delta\psi)_{L^2(\Omega)}=(\Delta f,\psi)_{L^2(\Omega)}.
\end{equation}
As in analogous contexts before, we denote by tilde the zero extension of a function, originally defined
in $\Omega$, to the entire manifold $M$. Then $\widetilde{f}\in L^2(M)$ and \eqref{uartUTR-Mi} implies that for 
each $\varphi\in C_0^\infty(M)$ we may write 
\begin{equation}\label{jh7g8-Mi}
\begin{split}
(\Delta\widetilde f,\varphi)_{L^2(M)} &=(\widetilde f,\Delta\varphi)_{L^2(M)}
=\big(f,\Delta\varphi|_\Omega\big)_{L^2(\Omega)}
\\[2pt]
&=\big(\Delta f,\varphi|_\Omega\big)_{L^2(\Omega)}=(\widetilde{\Delta f},\varphi)_{L^2(M)}.
\end{split}
\end{equation}
Hence, $\Delta\widetilde{f}=\widetilde{\Delta f}$ in $\cD'(M)$. Since $\widetilde{\Delta_g f}\in L^2(M)$, 
invoking standard elliptic regularity implies that $\widetilde{f}\in H^2(M)$, 
which further implies $f\in H^2(\Omega)$. With this in hand, one invokes  Lemma~\ref{l4.3-Mi} and 
\eqref{Tan-C3} in order to conclude that $\dom(L_{D,\Omega})\cap\dom(L_{N,\Omega})
\subset\accentset{\circ}{H}^2(\Omega)=\dom(L_{min,\Omega})$. This establishes the left-to-right 
inclusion in \eqref{3.30-Mi}. The opposite inclusion follows from Lemma~\ref{l4.3-Mi}
and the fact that $L_{D,\Omega}$ and $L_{N,\Omega}$ are both extensions of $L_{min,\Omega}$.
\end{proof}
%%%%%%%%%%

The machinery developed up to this point in this section makes it possible 
to study $z$-dependent Dirichlet-to-Neumann maps, that is, Weyl--Titchmarsh operators, 
for Schr\"odinger operators in Lipschitz subdomains of the compact Riemannian manifold $M$, 
in a very similar manner to the treatment in Section~\ref{s7} of the Euclidean setting. Deferring a 
detailed treatment of this circle of ideas to future work, a typical sample result in this connection 
reads as follows. 

%%%%%%%%%
\begin{theorem}\label{t5.2-Mi}
Assume that $\Omega\subset M$ is a Lipschitz domain, and suppose that $V\in L^\infty(M)$.
Then the following assertions hold: \\[1mm] 
$(i)$ For each $z\in\rho(L_{D,\Omega})$ and $s\in [0,1]$ the boundary value problem
\begin{equation}\label{eqn:bvp-Mi}
\begin{cases}
(-\Delta_g+V-z)f=0\,\text{ in $\Omega,\quad f\in H^{s+(1/2)}(\Omega)\cap\dom(L_{max,\Omega})$,}     
\\[2pt]  
\gamma_D f=\varphi\,\text{ on $\partial\Omega,\quad\varphi\in H^s(\partial\Omega)$,}
\end{cases}    
\end{equation}
is well posed, with unique solution $f=f_D(z,\varphi)$ given by 
\begin{equation}\label{4.7-Mi}
f_D(z,\varphi)=-\big[\gamma_N(L_{D,\Omega}-{\ol z}I)^{-1}\big]^*\varphi, 
\end{equation} 
where the star indicates the adjoint of 
\begin{equation}\label{3.2SD1-Mi}
\gamma_N(L_{D,\Om}-zI)^{-1}\in\cB\big(L^2(\Omega),L^2(\partial\Omega)\big).
\end{equation}
Moreover, if for each $z\in\rho(L_{D,\Omega})$ and $s\in[0,1]$ one defines  
\begin{equation}\label{We-Q.13-Mi}
P_{s,D,\Omega}(z): 
\begin{cases} 
H^s(\partial\Omega)\rightarrow H^{s+(1/2)}(\Omega)\cap\dom(L_{max,\Omega}), 
\\[2pt]  
\varphi\mapsto P_{s,D,\Omega}(z)\varphi:=f_D(z,\varphi),
\end{cases}    
\end{equation}
then the operator $\big[\gamma_N(L_{D,\Omega}-{\ol z}I)^{-1}\big]^*$, 
originally understood as the adjoint of \eqref{3.2SD1-Mi}, induces a mapping
\begin{equation}\label{GCC-Jna.1-Mi}
\big[\gamma_N(L_{D,\Omega}-{\ol z}I)^{-1}\big]^* 
\in\cB\big(H^s(\partial\Omega),H^{s+(1/2)}(\Omega)\cap\dom(L_{max,\Omega})\big)
\end{equation}
{\rm (}where the space $H^{s+(1/2)}(\Omega)\cap\dom(L_{max,\Omega})$ is equipped with the natural 
norm $f\mapsto\|f\|_{H^{s+1/2}(\Omega)}+\|\Delta_g f\|_{L^2(\Omega)}${\rm )}, and 
\begin{equation}\label{eqn:P_D*-Mi}
P_{s,D,\Omega}(z)=-\big[\gamma_N(L_{D,\Omega}-{\ol z}I)^{-1}\big]^*\,\text{ on }\, 
H^s(\partial\Omega).
\end{equation}
In addition, $P_{s,D,\Omega}(z)$ is injective with 
\begin{equation}\label{eqn:ranP_D-Mi} 
\ran(P_{s,D,\Omega}(z))=\ker(L_{max,\Omega}-zI)\cap H^{s+(1/2)}(\Omega).
\end{equation}
In particular, $\ran(P_{s,D,\Omega}(z))$ is dense in $\ker(L_{max,\Omega}-zI)$ 
with respect to the $L^2(\Omega)$-norm. \\[1mm] 
$(ii)$ For each $z\in\rho(L_{N,\Omega})$ and $s\in[0,1]$ the boundary value problem
\begin{equation}\label{eqn:bvp2-Mi}
\begin{cases}
(-\Delta_g+V-z)f=0\,\text{ in $\Omega,\quad f\in H^{s+(1/2)}(\Omega)\cap\dom(L_{max,\Omega})$,}     
\\[2pt] 
-\gamma_N f=\varphi\,\text{ in $H^{s-1}(\partial\Omega),\quad\varphi\in H^{s-1}(\partial\Omega)$,}
\end{cases}    
\end{equation}
is well posed, with unique solution $f=f_N(z,\varphi)$ given by 
\begin{equation}\label{4.9-Mi}
f_N(z,\varphi)=-\big[\gamma_D(L_{N,\Omega}-{\ol z}I)^{-1}\big]^*\varphi,   
\end{equation}
where the star indicates the adjoint of
\begin{equation}\label{3.2SD1N-Mi}
\gamma_D(L_{N,\Om}-zI)^{-1}\in\cB\big(L^2(\Omega),H^1(\partial\Omega)\big).
\end{equation}
Moreover, if for each $z\in\rho(L_{N,\Omega})$ and $s\in[0,1]$ one defines  
\begin{equation}\label{eqn:P_N-Mi}
P_{s,N,\Omega}(z): 
\begin{cases} 
H^{s-1}(\partial\Omega)\rightarrow H^{s+(1/2)}(\Omega)\cap\dom(L_{max,\Omega}), 
\\[2pt]
\varphi\mapsto P_{s,N,\Omega}(z)\varphi:=f_N(z,\varphi),
\end{cases}     
\end{equation}
then for each $z\in\rho(L_{N,\Omega})$ and $s\in[0,1]$ the operator 
$\big[\gamma_D(L_{N,\Omega}-{\ol z}I)^{-1}\big]^*$, initially regarded as the adjoint of \eqref{3.2SD1N-Mi}, 
induces a mapping
\begin{equation}\label{GCC-Jna.1N-Mi}
\big[\gamma_D(L_{N,\Omega}-{\ol z}I)^{-1}\big]^* 
\in\cB\big(H^{s-1}(\partial\Omega),H^{s+(1/2)}(\Omega)\cap\dom(L_{max,\Omega})\big)
\end{equation}
{\rm (}where the space $H^{s+(1/2)}(\Omega)\cap\dom(L_{max,\Omega})$ is equipped with 
the natural norm $f\mapsto\|f\|_{H^{s+1/2}(\Omega)}+\|\Delta_g f\|_{L^2(\Omega)}${\rm )}, 
and 
\begin{equation}\label{eqn:P_N*-Mi}
P_{s,N,\Omega}(z)=-\big[\gamma_D(L_{N,\Omega}-{\ol z}I)^{-1}\big]^*\,\text{ on }\, 
H^{s-1}(\partial\Omega).    
\end{equation}
In addition, $P_{s,N,\Omega}(z)$ is injective with  
\begin{equation}\label{eqn:ranP_N-Mi} 
\ran(P_{s,N,\Omega}(z))=\ker(L_{max,\Omega}-z I)\cap H^{s+(1/2)}(\Omega).
\end{equation} 
In particular, $\ran(P_{s,N,\Omega}(z))$ is dense in $\ker(L_{max,\Omega}-zI)$ 
with respect to the $L^2(\Omega)$-norm. \\[1mm] 
$(iii)$ For $z\in\rho(L_{D,\Omega})$ and $s\in[0,1]$, the Dirichlet-to-Neumann 
operator defined by
\begin{equation}\label{We-Q.14-Mi}
M_{s,\Omega}(z):
\begin{cases} 
H^s(\partial\Omega)\rightarrow H^{s-1}(\partial\Omega),   
\\[2pt] 
\varphi\mapsto M_{s,\Omega}(z)\varphi:=-\gamma_N P_{s,D,\Omega}(z)\varphi,
\end{cases}    
\end{equation}
satisfies
\begin{equation}\label{4.17-Mi}
M_{s,\Omega}(z)=\gamma_N\big[\gamma_N(L_{D,\Omega}-{\ol z}I)^{-1}\big]^* 
\in\cB\big(H^s(\partial\Omega),H^{s-1}(\partial\Omega)\big).     
\end{equation}
Moreover, for each $z\in\rho(L_{D,\Omega})$ and each $s\in[0,1]$,
\begin{align}\label{4.17AD-Mi}
\begin{split}
& \text{the adjoint of $M_{s,\Omega}(z)\in\cB\big(H^s(\partial\Omega),H^{s-1}(\partial\Omega)\big)$}
\\[2pt]
& \quad\text{is the operator 
$M_{1-s,\Omega}(\overline{z})\in\cB\big(H^{1-s}(\partial\Omega),H^{-s}(\partial\Omega)\big)$}.     
\end{split}
\end{align}
$(iv)$ For $z\in\rho(L_{N,\Omega})$ and $s\in[0,1]$, the Neumann-to-Dirichlet 
operator defined by
\begin{equation}\label{We-Q.15-Mi}
N_{s,\Omega}(z): 
\begin{cases} 
H^{s-1}(\partial\Omega)\rightarrow H^s(\partial\Omega),   
\\[2pt] 
\varphi\mapsto N_{s,\Omega}(z)\varphi:=-\gamma_D P_{s,N,\Omega}(z)\varphi,
\end{cases}    
\end{equation}
satisfies
\begin{equation}\label{4.19-Mi}
N_{s,\Omega}(z)=\gamma_D\big[\gamma_D(L_{N,\Omega}-{\ol z} I)^{-1}\big]^*  
\in\cB\big(H^{s-1}(\partial\Omega),H^s(\partial\Omega)\big).     
\end{equation} 
In addition, for each $z\in\rho(L_{N,\Omega})$ and each $s\in[0,1]$,
\begin{align}\label{4.17AN-Mi}
\begin{split}
& \text{the adjoint of $N_{s,\Omega}(z)\in\cB\big(H^{s-1}(\partial\Omega),H^s(\partial\Omega)\big)$}
\\[2pt]
& \quad\text{is the operator $N_{1-s,\Omega}(\overline{z})\in\cB\big(H^{-s}(\partial\Omega),H^{1-s}(\partial\Omega)\big)$}.     
\end{split}
\end{align}
$(v)$ If $z\in\rho(L_{D,\Omega})\cap\rho(L_{N,\Omega})$, then for each $s\in[0,1]$ 
the Dirichlet-to-Neumann operator $M_{s,\Omega}(z)$ maps $H^{s}(\partial\Omega)$ bijectively 
onto $H^{s-1}(\partial\Omega)$, the Neumann-to-Dirichlet operator $N_{s,\Omega}(z)$ maps 
$H^{s-1}(\partial\Omega)$ bijectively onto $H^{s}(\partial\Omega)$, and their inverses satisfy 
\begin{align}\label{5.36JJJ.1-Mi}
& M_{s,\Omega}(z)^{-1}=-N_{s,\Omega}(z)\in\cB\big(H^{s-1}(\partial\Omega),H^{s}(\partial\Omega)\big),
\\[2pt] 
& N_{s,\Omega}(z)^{-1}=-M_{s,\Omega}(z)\in\cB\big(H^{s}(\partial\Omega),H^{s-1}(\partial\Omega)\big).
\label{5.36JJJ.2-Mi}
\end{align} 
\end{theorem}
%%%%%%%%%
\begin{proof}
All claims may be justified in a similar fashion to their Euclidean counterparts proved 
in Theorem~\ref{t5.2}, by relying on the trace theory in Corollary~\ref{YTfdf-T.NNN-MMM} 
and Corollary~\ref{YTfdf.NNN.2-MMM}.
\end{proof}
%%%%%%%

In turn, having established Theorem~\ref{t5.2-Mi}, makes it possible to prove the following extension of Theorem~\ref{t5.5} to the 
setting of Lipschitz subdomains of Riemannian manifolds and with the Laplace--Beltrami operator replacing the ordinary 
flat-space Laplacian.

%%%%%%%%%%%%
\begin{theorem}\label{t5.5.MAN}
Assume that $\Omega\subset M$ is a Lipschitz domain, and suppose that $V\in L^\infty(M)$.
Consider the spaces
\begin{equation}\label{eqn:G0G1.MAN}
\mathscr{G}_D(\partial\Omega):=\ran\big(\gamma_D\big|_{\dom(L_{N,\Omega})}\big),\quad
\mathscr{G}_N(\partial\Omega):=\ran\big(\gamma_N\big|_{\dom(L_{D,\Omega})}\big),
\end{equation} 
and, the Dirichlet-to-Neumann map $M_{\Omega}(z):=M_{1,\Omega}(z)$ as in \eqref{We-Q.14-Mi}, define
\begin{equation}\label{eqn:SigmaLambda.MAN}
\Sigma:=\Im\big(-M_{\Omega}(i)^{-1}\big),\quad\Lambda:=\overline{\Im(M_{\Omega}(i))}.
\end{equation}

Then the following statements hold: \\[1mm]
$(i)$ Both $\Sigma$ and $\Lambda$ are bounded, nonnegative, self-adjoint operators in $L^2(\partial\Omega)$, 
that are invertible and have unbounded inverses. 
\\[1mm]
$(ii)$ One has
\begin{equation}\label{We-Q.17.MAN}
\begin{split}
\mathscr{G}_D(\partial\Omega) &=\big\{\gamma_D f\,\big|\,
f\in H^{3/2}(\Omega)\cap\dom(L_{max,\Omega}),\,\gamma_N f=0\big\}\subset H^1(\partial\Omega),
\\[2pt]
\mathscr{G}_N(\partial\Omega) &=\big\{\gamma_N f\,\big|\, 
f\in H^{3/2}(\Omega)\cap\dom(L_{max,\Omega}),\,\gamma_D f=0\big\}\subset L^2(\partial\Omega).
\end{split}
\end{equation}

\noindent $(iii)$ One has 
\begin{equation}\label{eqn:domLambda.MAN}
\begin{split}
\mathscr{G}_D(\partial\Omega)&=\dom\big(\Sigma^{-1/2}\big)=\ran\big(\Sigma^{1/2}\big), 
\\[2pt]
\mathscr{G}_N(\partial\Omega)&=\dom\big(\Lambda^{-1/2}\big)=\ran\big(\Lambda^{1/2}\big),
\end{split}
\end{equation} 
and when equipped with the scalar products 
\begin{equation}\label{scalarjussi.MAN}
\begin{split}
(\varphi,\psi)_{\mathscr{G}_D(\partial\Omega)} 
&:=\big(\Sigma^{-1/2}\varphi,\Sigma^{-1/2}\psi\big)_{L^2(\partial\Omega)},
\quad\forall\,\varphi,\psi\in\mathscr{G}_D(\partial\Omega),    
\\[2pt] 
(\varphi,\psi)_{\mathscr{G}_N(\partial\Omega)} 
&:=\big(\Lambda^{-1/2}\varphi,\Lambda^{-1/2}\psi\big)_{L^2(\partial\Omega)},
\quad\forall\,\varphi,\psi\in\mathscr{G}_N(\partial\Omega),
\end{split}
\end{equation}
the spaces $\mathscr{G}_D(\partial\Omega),\mathscr{G}_N(\partial\Omega)$ become Hilbert spaces. 
\\[1mm] 
$(iv)$ The Dirichlet trace operator $\gamma_D$ {\rm (}as defined in \eqref{eqn:gammaDs.2-MMM}{\rm )} 
and the Neumann trace operator $\gamma_N$ {\rm (}as defined in \eqref{eqn:gammaN-pp-MMM}{\rm )} 
extend by continuity {\rm (}hence in a compatible manner{\rm )} to continuous surjective mappings 
\begin{equation}\label{4.55.MAN}
\begin{split}
\widetilde{\gamma}_D:\dom(L_{max,\Omega})&\to\mathscr{G}_N(\partial\Omega)^*,
\\[2pt]  
\widetilde{\gamma}_N:\dom(L_{max,\Omega})&\to\mathscr{G}_D(\partial\Omega)^*,  
\end{split}
\end{equation}
where $\dom(L_{max,\Omega})$ is endowed with the graph norm of $L_{max,\Omega}$, 
and $\mathscr{G}_D(\partial\Omega)^*$, $\mathscr{G}_N(\partial\Omega)^*$ are, respectively, 
the adjoint {\rm (}conjugate dual{\rm )} spaces of $\mathscr{G}_D(\partial\Omega)$,
$\mathscr{G}_N(\partial\Omega)$ carrying the natural topology induced by 
\eqref{scalarjussi.MAN} on $\mathscr{G}_D(\partial\Omega)$, 
$\mathscr{G}_N(\partial\Omega)$, respectively, such that
\begin{equation}\label{4.30.MAN}
\ker(\widetilde{\gamma}_D)=\dom(L_{D,\Omega})\,\text{ and }\,   
\ker(\widetilde{\gamma}_N)=\dom(L_{N,\Omega}).
\end{equation} 
Furthermore, for each $s\in[0,1]$ there exists a constant $C\in(0,\infty)$ with the property that
\begin{align}\label{4.55.REE.D.MAN}
\begin{split}
& f\in\dom(L_{max,\Omega})\,\text{ and }\,\widetilde{\gamma}_D f\in H^s(\partial\Omega)
\,\text{ imply }\,f\in H^{s+(1/2)}(\Omega)
\\[2pt]
& \quad\text{and }\,\|f\|_{H^{s+(1/2)}(\Omega)}\leq C\big(\|\Delta_g f\|_{L^2(\Omega)}
+\|\widetilde{\gamma}_D f\|_{H^s(\partial\Omega)}\big),
\end{split}
\end{align}
and
\begin{align}\label{4.55.REE.N.MAN}
\begin{split}
& f\in\dom(L_{max,\Omega})\,\text{ and }\,\widetilde{\gamma}_N f\in H^{-s}(\partial\Omega)
\,\text{ imply }\,f\in H^{-s+(3/2)}(\Omega)
\\[2pt]
& \quad\text{and }\,\|f\|_{H^{-s+(3/2)}(\Omega)}\leq C\big(\|f\|_{L^2(\Omega)}+\|\Delta_g f\|_{L^2(\Omega)}
+\|\widetilde{\gamma}_N f\|_{H^{-s}(\partial\Omega)}\big).
\end{split}
\end{align} 

$(v)$ With $\widetilde{\gamma}_D,\widetilde{\gamma}_N$ as in \eqref{4.55.MAN}, one has
\begin{align}\label{ut444.MAN}
\accentset{\circ}{H}^2(\Omega)=\big\{f\in\dom(L_{max,\Omega})\,\big|\,
& \wti{\gamma}_D f=0\,\text{ in }\,\mathscr{G}_N(\partial\Omega)^*
\nonumber\\[2pt]
&\text{and }\,\wti{\gamma}_N f=0\,\text{ in }\,\mathscr{G}_D(\partial\Omega)^*\big\}.  
\end{align}

$(vi)$ The manner in which the mapping $\widetilde{\gamma}_D$ in \eqref{4.55.MAN} 
operates is as follows: Given $f\in\dom(L_{max,\Omega})$, the action of the functional 
$\widetilde{\gamma}_D f\in\mathscr{G}_N(\partial\Omega)^*$ on some arbitrary 
$\phi\in\mathscr{G}_N(\partial\Omega)$ is given by 
\begin{equation}\label{eq:33f4iU.MAN}
{}_{\mathscr{G}_N(\partial\Omega)^\ast}\big\langle\widetilde{\gamma}_D f,
\phi\big\rangle_{\mathscr{G}_N(\partial\Omega)}
=(f,\Delta_g h)_{L^2(\Omega)}-(\Delta_g f,h)_{L^2(\Omega)},
\end{equation}
for any $h\in H^{3/2}(\Omega)\cap\dom(L_{max,\Om})$ such that $\gamma_D h=0$ 
and $\gamma_N h=\phi$ {\rm (}the existence of such $h$ being ensured by \eqref{We-Q.17.MAN}{\rm )}.
As a consequence, the following Green's formula holds: 
\begin{equation}\label{eq:33VaV.MAN}
{}_{\mathscr{G}_N(\partial\Omega)^\ast}\big\langle\widetilde{\gamma}_D f,
\gamma_N h\big\rangle_{\mathscr{G}_N(\partial\Omega)}
=(f,\Delta_g h)_{L^2(\Omega)}-(\Delta_g f,h)_{L^2(\Omega)},
\end{equation}
for each $f\in\dom(L_{max,\Omega})$ and each $h\in\dom(L_{D,\Om})$. \\[1mm] 
$(vii)$ The mapping $\widetilde{\gamma}_N$ in \eqref{4.55.MAN} 
operates in the following fashion: Given a function $f\in\dom(L_{max,\Omega})$, the 
action of the functional $\widetilde{\gamma}_N f\in\mathscr{G}_D(\partial\Omega)^*$ 
on some arbitrary $\psi\in\mathscr{G}_D(\partial\Omega)$ is given by 
\begin{equation}\label{eq:33f4iU.NN.MAN}
{}_{\mathscr{G}_D(\partial\Omega)^\ast}\big\langle\widetilde{\gamma}_N f,
\psi\big\rangle_{\mathscr{G}_D(\partial\Omega)}
=-(f,\Delta_g h)_{L^2(\Omega)}+(\Delta_g f,h)_{L^2(\Omega)},
\end{equation}
for any $h\in H^{3/2}(\Omega)\cap\dom(L_{max,\Om})$ such that $\gamma_N h=0$ 
and $\gamma_D h=\psi$ {\rm (}the existence of such $h$ being ensured by \eqref{We-Q.17.MAN}{\rm )}.
In particular, the following Green's formula holds: 
\begin{equation}\label{eq:33VaV.NN.MAN}
{}_{\mathscr{G}_D(\partial\Omega)^\ast}\big\langle\widetilde{\gamma}_N f,
\gamma_D h\big\rangle_{\mathscr{G}_D(\partial\Omega)}
=-(f,\Delta_g h)_{L^2(\Omega)}+(\Delta_g f,h)_{L^2(\Omega)},
\end{equation}
for each $f\in\dom(L_{max,\Omega})$ and each $h\in\dom(L_{N,\Om})$. \\[1mm] 
$(viii)$ The operators 
\begin{align}\label{4.5ASdf.1.MAN}
&\gamma_D:\dom(L_{N,\Om})=H^{3/2}(\Omega)\cap\dom(L_{max,\Omega})\cap\ker(\gamma_N)
\to\mathscr{G}_D(\partial\Omega),
\\[2pt] 
&\gamma_N:\dom(L_{D,\Om})=H^{3/2}(\Omega)\cap\dom(L_{max,\Omega})\cap\ker(\gamma_D)
\to\mathscr{G}_N(\partial\Omega),    
\label{4.5ASdf.2.MAN}
\end{align}
are well defined, linear, surjective, and continuous if for some $s\in[0,\tfrac{3}{2}]$ both 
spaces on the left-hand sides of \eqref{4.5ASdf.1.MAN}, \eqref{4.5ASdf.2.MAN} are equipped with the 
norm $f\mapsto\|f\|_{H^s(\Omega)}+\|\Delta_g f\|_{L^2(\Omega)}$ {\rm (}which are all equivalent{\rm )}. In addition, 
\begin{equation}\label{eq:Ajf.MAN}
\text{the kernel of $\gamma_D$ and $\gamma_N$ in \eqref{4.5ASdf.1.MAN}--\eqref{4.5ASdf.2.MAN} 
is $\accentset{\circ}{H}^2(\Omega)$}.
\end{equation}
Moreover, 
\begin{align}\label{i7g5r.MAN}
\|\phi\|_{\mathscr{G}_D(\partial\Omega)}& \approx\inf_{\substack{f\in
H^{3/2}(\Omega)\cap\dom(L_{max,\Omega})\\ \gamma_N f=0,\,\,\gamma_D f=\phi}}
\big(\|f\|_{H^{3/2}(\Omega)}+\|\Delta_g f\|_{L^2(\Omega)}\big)
\nonumber\\[2pt]
& \approx\inf_{\substack{f\in H^{3/2}(\Omega)\cap\dom(L_{max,\Omega})\\ \gamma_N f=0,\,\,\gamma_D f=\phi}}
\big(\|f\|_{L^2(\Omega)}+\|\Delta_g f\|_{L^2(\Omega)}\big)
\nonumber\\[2pt]
& \approx\inf_{\substack{f\in\dom(L_{max,\Omega})\\ \widetilde{\gamma}_N f=0,\,\,\widetilde{\gamma}_D f=\phi}}
\big(\|f\|_{L^2(\Omega)}+\|\Delta_g f\|_{L^2(\Omega)}\big),
\end{align}
uniformly for $\phi\in\mathscr{G}_D(\partial\Omega)$, and 
\begin{align}\label{i7g5r.2N.MAN}
\|\psi\|_{\mathscr{G}_N(\partial\Omega)}& 
\approx\inf_{\substack{h\in H^{3/2}(\Omega)\cap\dom(L_{max,\Omega})\\ \gamma_D h=0,\,\,\gamma_N h=\psi}}
\big(\|h\|_{H^{3/2}(\Omega)}+\|\Delta_g h\|_{L^2(\Omega)}\big)
\nonumber\\[2pt]
& \approx\inf_{\substack{h\in H^{3/2}(\Omega)\cap\dom(L_{max,\Omega})\\ \gamma_D h=0,\,\,\gamma_N h=\psi}}
\big(\|h\|_{L^2(\Omega)}+\|\Delta_g h\|_{L^2(\Omega)}\big)
\nonumber\\[2pt]
& \approx\inf_{\substack{h\in\dom(L_{max,\Omega})\\ \widetilde\gamma_D h=0,\,\,\widetilde\gamma_N h=\psi}}
\big(\|h\|_{L^2(\Omega)}+\|\Delta_g h\|_{L^2(\Omega)}\big)
\nonumber\\[2pt]
& \approx\inf_{\substack{h\in\dom(L_{max,\Omega})\\ \widetilde{\gamma}_D h=0,\,\,\widetilde{\gamma}_N h=\psi}}
\|\Delta_g h\|_{L^2(\Omega)},
\end{align}
uniformly for $\psi\in\mathscr{G}_N(\partial\Omega)$.

As a consequence,
\begin{align}\label{ytrr555e.MAN}
\begin{split}
& \mathscr{G}_D(\partial\Omega)\hookrightarrow H^1(\partial\Omega)
\hookrightarrow L^2(\partial\Omega)\hookrightarrow H^{-1}(\partial\Omega)
\hookrightarrow\mathscr{G}_D(\partial\Omega)^*,  
\\[2pt] 
& \mathscr{G}_N(\partial\Omega)\hookrightarrow L^2(\partial\Omega)
\hookrightarrow\mathscr{G}_N(\partial\Omega)^*, 
\end{split}
\end{align} 
with all embeddings linear, continuous, and with dense range. Moreover, the duality 
pairings between $\mathscr{G}_D(\partial\Omega)$ and $\mathscr{G}_D(\partial\Omega)^\ast$,
as well as between $\mathscr{G}_N(\partial\Omega)$ and $\mathscr{G}_N(\partial\Omega)^\ast$, 
are both compatible with the inner product in $L^2(\partial\Omega)$. \\[1mm] 
\noindent $(ix)$ For each $z\in\rho(L_{D,\Omega})$, the boundary value problem 
\begin{equation}\label{eqn:bvp.2.MAN}
\begin{cases} 
(-\Delta_g+V-z)f=0\,\text{ in $\Omega,\quad f\in\dom(L_{max,\Omega})$,}   
\\[2pt]  
\widetilde{\gamma}_D f=\varphi
\,\text{ in $\mathscr{G}_N(\partial\Omega)^*,\quad\varphi\in\mathscr{G}_N(\partial\Omega)^*$,}
\end{cases} 
\end{equation}
is well posed. In particular, for each $z\in\rho(L_{D,\Omega})$ there exists a constant $C\in(0,\infty)$, which depends 
only on $\Omega$, $n$, $z$, and $V$, with the property that
\begin{align}\label{eqn:bvp.2WACO.MAN}
\begin{split}
& \|f\|_{L^2(\Omega)}\leq C\|\widetilde{\gamma}_D f\|_{\mathscr{G}_N(\partial\Omega)^*}\,\text{ for each }\,f\in\dom(L_{max,\Omega})     
\\[2pt]
& \quad\text{with }\,(-\Delta_g+V-z)f=0\,\text{ in }\,\Omega.
\end{split} 
\end{align}
Moreover, if 
\begin{equation}\label{6h8u.MAN}
\widetilde P_{D,\Omega}(z): 
\begin{cases} 
\mathscr{G}_N(\partial\Omega)^*\rightarrow\dom(L_{max,\Omega}),   
\\[2pt] 
\varphi\mapsto\widetilde P_{D,\Omega}(z)\varphi:=\wti{f}_{D,\Omega}(z,\varphi),
\end{cases} 
\end{equation}
where $\wti{f}_{D,\Omega}(z,\varphi)$ is the unique solution of \eqref{eqn:bvp.2.MAN}, 
then the solution operator $\widetilde{P}_{D,\Omega}(z)$ is an extension of 
$P_{0,D,\Omega}(z)$ in \eqref{We-Q.13-Mi}, and $\widetilde{P}_{D,\Omega}(z)$ is continuous, when 
the adjoint space $\mathscr{G}_N(\partial\Omega)^*$ and $\dom(L_{max,\Omega})$ are endowed 
with the norms in the current item $(iv)$. \\[1mm] 
\noindent $(x)$ For each $z\in\rho(L_{N,\Omega})$, the boundary value problem 
\begin{equation}\label{eqn:bvpNeumann.MAN}
\begin{cases} 
(-\Delta_g+V-z)f=0\,\text{ in $\Omega,\quad f\in\dom(L_{max,\Omega})$,}   
\\[2pt]  
-\widetilde{\gamma}_N f=\varphi 
\,\text{ in $\mathscr{G}_D(\partial\Omega)^*,\quad\varphi\in\mathscr{G}_D(\partial\Omega)^*$,}
\end{cases} 
\end{equation}
is well posed. In particular, for each $z\in\rho(L_{N,\Omega})$ there exists a constant $C\in(0,\infty)$, which depends 
only on $\Omega$, $n$, $z$, and $V$, with the property that
\begin{align}\label{eqn:bvpNeumann.WACO.MAN}
\begin{split}
& \|f\|_{L^2(\Omega)}\leq C\|\widetilde{\gamma}_N f\|_{\mathscr{G}_D(\partial\Omega)^*}\,\text{ for each }
\, f\in\dom(L_{max,\Omega})     
\\[2pt]
& \quad\text{with }\,(-\Delta_g+V-z)f=0\,\text{ in }\,\Omega.
\end{split} 
\end{align}
Moreover, if 
\begin{equation}\label{755fv.MAN}
\widetilde{P}_{N,\Omega}(z): 
\begin{cases} 
\mathscr{G}_D(\partial\Omega)^*\rightarrow\dom(L_{max,\Omega}),  
\\[2pt] 
\varphi\mapsto\widetilde{P}_{N,\Omega}(z)\varphi:=\wti{f}_{N,\Omega}(z,\varphi), 
\end{cases} 
\end{equation}
where $\wti{f}_{N,\Omega}(z,\varphi)$ is the unique solution of \eqref{eqn:bvpNeumann.MAN}, 
then the solution operator $\widetilde{P}_{N,\Omega}(z)$ is an extension of 
$P_{1,N,\Omega}(z)$ in \eqref{eqn:P_N-Mi}, and $\widetilde{P}_{N,\Omega}(z)$ is continuous, 
when the adjoint space $\mathscr{G}_D(\partial\Omega)^*$ and $\dom(L_{max,\Omega})$ 
are endowed with the norms in the current item $(iv)$. \\[1mm] 
\noindent $(xi)$ For all $z\in\rho(L_{D,\Omega})$ the Dirichlet-to-Neumann map $M_{\Omega}(z):=M_{1,\Omega}(z)$ in  
\eqref{We-Q.14-Mi} admits an extension 
\begin{equation}\label{iutee.MAN}
\widetilde M_{\Omega}(z): 
\begin{cases} 
\mathscr{G}_N(\partial\Omega)^*\rightarrow\mathscr{G}_D(\partial\Omega)^*,  
\\[2pt]
\varphi\mapsto\widetilde{M}_{\Omega}(z)\varphi
:=-\widetilde{\gamma}_N\widetilde{P}_{D,\Omega}(z)\varphi,
\end{cases} 
\end{equation}
and $\widetilde{M}_{\Omega}(z)$ is continuous, when the adjoint spaces 
$\mathscr{G}_D(\partial\Omega)^*$, $\mathscr{G}_N(\partial\Omega)^*$ carry the 
natural topology induced by \eqref{scalarjussi.MAN} on $\mathscr{G}_D(\partial\Omega)$, 
$\mathscr{G}_N(\partial\Omega)$, respectively.
\end{theorem}
%%%%%%%%%%
\begin{proof}
We may establish all claims reasoning analogously to the proof of the Euclidean result in Theorem~\ref{t5.5}, 
now relying on the trace theory in Corollary~\ref{YTfdf-T.NNN-MMM} and Corollary~\ref{YTfdf.NNN.2-MMM}, as well as 
the theory of Weyl--Titchmarsh operators for Schr\"odinger operators in Lipschitz subdomains of the compact 
Riemannian manifold $M$ developed in Theorem~\ref{t5.2}. 
\end{proof}
%%%%%%%

%%%%%%%%%%%%%%%%
\subsection{Variable coefficient elliptic operators in Euclidean Lipschitz domains}\label{ss.MM.EEE}
%%%%%%%%%%%%%%%%

Virtually everything we have established so far in this chapter for the perturbed Laplace--Beltrami operator $\Delta_g+V$ 
on Lipschitz subdomains of Riemannian manifolds yields corresponding results for variable coefficient Schr\"odinger 
operators in Euclidean Lipschitz domains, in a natural way. The goal in this section is to briefly elaborate on this aspect.
For example, having proved Theorem~\ref{thm.Regularity}, we can now establish 
regularity results in the spirit of \eqref{eq:MM3}--\eqref{eq:MM3-AC}, and \eqref{eq:MM3BIS}--\eqref{eq:MM3BIS-AC} (with $k=1$),  
for variable coefficient elliptic operators in place of the standard Laplacian in ${\mathbb{R}}^n$.

To set the stage, given a nonempty, bounded open set 
$\Omega\subset{\mathbb{R}}^n$, we agree to introduce
\begin{align}
\begin{split} 
& C^{1,1}(\overline{\Omega}):=\big\{\varphi:\Omega\to\bbC\,\big|\,\text{there exists an open neighborhood ${\mathcal{O}}$ 
of $\overline{\Omega}$}  
\\ 
& \hspace*{3.8cm} 
\text{and $\Phi\in C^{1,1}({\mathcal{O}})$ such that $\Phi\big|_{\Omega}=\varphi$}\big\}. 
\end{split} 
\end{align} 

%%%%%%%
\begin{theorem}\label{thm.Regularity-CC}
Let $\Omega\subset{\mathbb{R}}^n$ be a bounded Lipschitz domain, $u\in C^1(\Omega)$, and consider a second-order 
divergence-form differential expression ${\mathcal{L}}$, acting according to 
\begin{equation}\label{hyrdfr4}
{\mathcal{L}}u:=\sum\limits_{j,k=1}^n
\partial_j\big(a_{jk}(x)\partial_ku\big)\,\text{ in }\,\Omega,
\end{equation}
in the sense of distributions, where $A(x)=\big(a_{jk}(x)\big)_{1\leq j,k\leq n}$, with $x\in\Omega$, is a symmetric, 
positive definite matrix, with real-valued entries $a_{jk}\in C^{1,1}(\overline{\Omega})$. Moreover, pick a 
real-valued potential $V\in L^p(\Omega)$, with $p>n$, and introduce 
\begin{equation}\label{eqn.lmnae-CC}
L:=-{\mathcal{L}}+V\,\text{ in }\,\Omega.
\end{equation}

Then for any function $u\in C^1(\Omega)$ solving 
\begin{equation}\label{eqn.sjul-A-CC}
Lu=0\,\text{ in }\,{\mathcal{D}}'(\Omega)
\end{equation}
one has
\begin{align}\label{eqn.jslrh-111-CC}
\begin{split}
& \mathcal{N}_{\kappa}u\in L^2(\partial\Omega)\Longleftrightarrow u\in H^{1/2}(\Omega), 
\\[2pt]
& \big\|\mathcal{N}_{\kappa}u\big\|_{L^2(\partial\Omega)}
\approx\|u\|_{H^{1/2}(\Omega)},
\end{split}
\end{align}
as well as
\begin{align}\label{eqn.jslrh-CC}
\begin{split} 
& \mathcal{N}_{\kappa}(\nabla u)\in L^2(\partial\Omega)\Longleftrightarrow u\in H^{3/2}(\Omega), 
\\[2pt]
& \big\|\mathcal{N}_{\kappa}u\big\|_{L^2(\partial\Omega)}
+\big\|\mathcal{N}_{\kappa}(\nabla u)\big\|_{L^2(\partial\Omega)}\approx\|u\|_{H^{3/2}(\Omega)},
\end{split}
\end{align}
uniformly for $u\in C^1(\Omega)$ satisfying \eqref{eqn.sjul-A-CC}.
\end{theorem}
%%%%%%%
\begin{proof}
The key observation is that any divergence-form operator ${\mathcal{L}}$ as in \eqref{hyrdfr4}
coincides, up to left multiplication by a power of ${\rm det}\,A(x)$, with the Laplace--Beltrami operator 
$\Delta_g$ of the manifold $\Omega$ equipped with the Riemannian metric tensor
\begin{equation}\label{hyrdfr5}
g:=\big(\text{det}(A)\big)^{1/(n-2)}\sum_{j,k=1}^n a^{jk}\,dx_j\otimes dx_k,
\end{equation}
where the $a^{jk}$'s are the entries in the matrix $A^{-1}$.
Specifically, if $\Delta_g$ is the Laplace--Beltrami operator associated as in \eqref{eqn.thuaou}
with the metric tensor $g$ given in \eqref{hyrdfr5}, then 
\begin{equation}\label{eqn.thu-TTG}
{\mathcal{L}}=\big(\text{det} (A)\big)^{1/(n-2)}\Delta_g.
\end{equation}
In particular, for any function $u\in C^1(\Omega)$ one has
\begin{align}\label{eqn.sjul-A-CC2}
\begin{split}
& Lu=0\Longleftrightarrow(-\Delta_g+V_A)u=0, 
\\[2pt]
& \quad\text{where }\,V_A:=\big(\text{det}(A)\big)^{-1/(n-2)}V\in L^p(\Omega).
\end{split}
\end{align}
Then all desired conclusions will follow from Theorem~\ref{thm.Regularity} as soon as one succeeds in 
viewing $\Omega$ as a subset of a local coordinate patch of a smooth, compact, boundaryless, Riemannian 
manifold $M$, whose metric tensor agrees with \eqref{hyrdfr5} near $\overline{\Omega}$. 

With this aim in mind, let ${\mathcal{O}}$ be an open neighborhood of $\overline{\Omega}$
with the property that the entries of the matrix $A$ extend to real-valued functions in 
$C^{1,1}({\mathcal{O}})$. We retain the same notation $a_{jk}$ for these entries
and observe that there is no loss of generality in assuming that the matrix
$\big(a_{jk}(x)\big)_{1\leq j,k\leq n}$ continues to be symmetric and positive definite for each 
$x\in\mathcal{O}$. To proceed, pick a function $\eta\in C^\infty_0({\mathcal{O}})$ satisfying $0\leq\eta\leq 1$ 
as well as $\eta=1$ near $\overline{\Omega}$, and consider the Riemannian metric in ${\mathbb{R}}^n$ 
given by
\begin{align}\label{hyrdfr5-RR}
\begin{split}
& \displaystyle
g:=\sum_{j,k=1}^n g_{jk}\,dx_j\otimes dx_k\,\text{ where, for }\,1\leq j,k\leq n,
\\[2pt]
& \quad\text{we have set }\,
g_{jk}:=(1-\eta)\delta_{jk}+\eta\big(\text{det}(A)\big)^{1/(n-2)}a^{jk}.
\end{split}
\end{align}
It is apparent from \eqref{hyrdfr5-RR} that 
near $\overline{\Omega}$ we have 
\begin{equation}\label{eq:52d2cv.1}
\sqrt{g}=\big(\text{det}(A)\big)^{1/(n-2)}
\end{equation}
and
\begin{equation}\label{eq:52d2cv.2}
g^{jk}=\big(\text{det}(A)\big)^{-1/(n-2)}a_{jk}\,\text{ for }\,1\leq j,k\leq n.
\end{equation}

In addition, select a sufficiently large number $R>0$ such that $\overline{\mathcal{O}}\subset(0,R)^n$, 
and define the torus
\begin{equation}\label{hyrdfr5-Ra2}
M:={\mathbb{R}}^n\big/\sim
\end{equation}
where $\sim$ is the equivalence relation in ${\mathbb{R}}^n$ given by 
\begin{equation}\label{hyrdfr5-Ra3}
x\sim y\Longleftrightarrow x-y\in\{0,\pm Re_1,\dots,\pm Re_n\}
\end{equation}
for every $x,y\in{\mathbb{R}}^n$. Then $M$ is a ($C^\infty$) smooth, compact, boundaryless, 
manifold, of real dimension $n$, which contains $\overline{\Omega}$ in a single coordinate chart. 
Moreover, since for $1\leq j,k\leq n$ one has $g_{jk}=\delta_{jk}$ near the boundary of 
the cube $(0,R)^n$, it follows that \eqref{hyrdfr5-RR} induces a Riemannian metric on $M$ 
which has $C^{1,1}$-coefficients and which coincides with the metric \eqref{hyrdfr5} near 
$\overline{\Omega}$. 
\end{proof}
%%%%%%%

By the same token, we may painlessly reformulate results proved earlier in Subsections~\ref{ss.MM.1}--\ref{ss.MM.3}
in the language of elliptic differential operators with variable coefficients, of class $C^{1,1}$, on the closure of 
a bounded Lipschitz domain $\Omega\subset{\mathbb{R}}^n$. Given their intrinsic importance, we shall elaborate the 
variable-coefficient versions of the Euclidean trace results from Theorem~\ref{YTfdf-T} and Theorem~\ref{GNT}, 
starting with the former.  

%%%%%%%%%%
\begin{theorem}\label{YTfdf-T-MMM.Euclid}
Fix an arbitrary $\varepsilon>0$, let $\Omega\subset{\mathbb{R}}^n$ be a bounded Lipschitz domain,  and consider a second-order divergence-form differential expression ${\mathcal{L}}$, acting on each distribution $u\in H^{-1}_{\rm loc}(\Omega)$ according to 
\begin{equation}\label{hyrdfr4.Euclid}
{\mathcal{L}}u:=\sum\limits_{j,k=1}^n
\partial_j\big(a_{jk}(x)\partial_ku\big)\,\text{ in }\,\Omega,
\end{equation}
in the sense of distributions, where $A(x)=\big(a_{jk}(x)\big)_{1\leq j,k\leq n}$, with $x\in\Omega$, is a symmetric, positive definite matrix, with real-valued entries $a_{jk}\in C^{1,1}(\overline{\Omega})$ (see \eqref{eqn.thu-TTG.Euclid.4} below, and the subsequent comment).

Then the restriction of the boundary trace operator $\gamma_D$ from \eqref{eqn:gammaDs.1} 
to the space $\big\{u\in H^s(\Omega)\,\big|\,{\mathcal{L}}u\in H^{s-2+\varepsilon}(\Omega)\big\}$, 
originally considered for $s\in\big(\tfrac{1}{2},\tfrac{3}{2}\big)$, induces 
a well defined, linear, continuous operator 
\begin{equation}\label{eqn:gammaDs.2aux-MMM.Euclid}
\gamma_D:\big\{u\in H^s(\Omega)\,\big|\,{\mathcal{L}}u\in H^{s-2+\varepsilon}(\Omega)\big\}
\rightarrow H^{s-(1/2)}(\partial\Omega),\quad\forall\,s\in\big[\tfrac{1}{2},\tfrac{3}{2}\big]     
\end{equation}
{\rm (}throughout, the space on the left-hand side of \eqref{eqn:gammaDs.2aux-MMM.Euclid} is equipped with the natural graph 
norm $u\mapsto\|u\|_{H^{s}(\Omega)}+\|{\mathcal{L}}u\|_{H^{s-2+\varepsilon}(\Omega)}${\rm )}, which continues 
to be compatible with \eqref{eqn:gammaDs.1} when $s\in\big(\tfrac{1}{2},\tfrac{3}{2}\big)$.
Thus defined, the Dirichlet trace operator possesses the following additional properties: \\[1mm] 
$(i)$ The Dirichlet boundary trace operator in \eqref{eqn:gammaDs.2aux-MMM.Euclid} is surjective.
In fact, there exist linear and bounded operators
\begin{equation}\label{2.88X-NN-ii-RRDD-MMM.Euclid}
\Upsilon_D:H^{s-(1/2)}(\partial\Omega)\rightarrow
\big\{u\in H^s(\Omega)\,\big|\,{\mathcal{L}}u\in L^2(\Omega)\big\},\quad s\in\big[\tfrac12,\tfrac32\big],
\end{equation}
which are compatible with one another and serve as right-inverses for the Dirichlet trace, that is, 
\begin{equation}\label{2.88X-NN2-ii-RRDD-MMM.Euclid}
\gamma_D(\Upsilon_D\psi)=\psi,
\quad\forall\,\psi\in H^{s-(1/2)}(\partial\Omega)\,\text{ with }\,s\in\big[\tfrac12,\tfrac32\big].
\end{equation}
In fact, matters may be arranged so that each function in the range of $\Upsilon_D$ is a null-solution of ${\mathcal{L}}$, that is, 
\begin{equation}\label{2.88X-NN2-ii-RRDD.bis.3.Euclid}
{\mathcal{L}}(\Upsilon_D\psi)=0,\quad\forall\,\psi\in H^{s-(1/2)}(\partial\Omega)
\,\text{ with }\,s\in\big[\tfrac{1}{2},\tfrac{3}{2}\big].
\end{equation}
$(ii)$ The Dirichlet boundary trace operator \eqref{eqn:gammaDs.2aux-MMM.Euclid} is compatible
with the pointwise nontangential trace in the sense that:
\begin{align}\label{eqn:gammaDs.2auxBBB-MMM.Euclid}
\begin{split}
& \text{if $u\in H^s(\Omega)$ has ${\mathcal{L}}u\in H^{s-2+\varepsilon}(\Omega)$ for some 
$s\in\big[\tfrac{1}{2},\tfrac{3}{2}\big]$,}    
\\[2pt]
& \quad\text{and if $u\big|^{\kappa-{\rm n.t.}}_{\partial\Omega}$ exists $\sigma$-a.e.~on $\partial\Omega$,  
then $u\big|^{\kappa-{\rm n.t.}}_{\partial\Omega}=\gamma_D u\in H^{s-(1/2)}(\partial\Omega)$.}
\end{split}
\end{align}
$(iii)$ The Dirichlet boundary trace operator $\gamma_D$ in 
\eqref{eqn:gammaDs.2aux-MMM.Euclid} is the unique extension by continuity and density of the mapping 
$C^\infty(\overline{\Omega})\ni f\mapsto f\big|_{\partial\Omega}$. \\[1mm] 
$(iv)$ For each $s\in\big[\tfrac{1}{2},\tfrac{3}{2}\big]$ the Dirichlet 
boundary trace operator satisfies
\begin{align}\label{eqn:gammaDs.1-TR.2-MMM.Euclid}
\begin{split}
& \gamma_D(\Phi u)=\big(\Phi\big|_{\partial\Omega}\big)\gamma_D u
\,\text{ at $\sigma$-a.e.~point on $\partial\Omega$, for all}
\\[2pt]
& \quad u\in H^s(\Omega)\,\text{ with }\,{\mathcal{L}}u\in H^{s-2+\varepsilon}(\Omega)
\,\text{ and }\,\Phi\in C^\infty(\overline\Omega).
\end{split}
\end{align}
$(v)$ For each $s\in\big[\tfrac{1}{2},\tfrac{3}{2}\big]$ such that $\varepsilon\not=\tfrac{3}{2}-s$, 
the null space of the Dirichlet boundary trace operator \eqref{eqn:gammaDs.2aux-MMM.Euclid} satisfies
\begin{equation}\label{eq:EFFa-MMM.Euclid}
{\rm ker}(\gamma_D)\subseteq H^{\,\min\!\{s+\varepsilon,3/2\}}(\Omega).
\end{equation}
In fact, the inclusion recorded in \eqref{eq:EFFa-MMM.Euclid} is quantitative in the sense that,
whenever $s\in\big[\tfrac{1}{2},\tfrac{3}{2}\big]$ is such that $\varepsilon\not=\tfrac{3}{2}-s$, 
there exists a constant $C\in(0,\infty)$ with the property that
\begin{align}\label{gafvv.655-MMM.Euclid} 
& \text{if $u\in H^s(\Omega)$ satisfies ${\mathcal{L}}u\in H^{s-2+\varepsilon}(\Omega)$ 
and $\gamma_D u=0$}   
\nonumber\\[2pt]
& \quad\text{then the function $u$ belongs to $H^{\,\min\!\{s+\varepsilon,3/2\}}(\Omega)$ and}    
\\[2pt]
& \quad\text{$\|u\|_{H^{\,\min\!\{s+\varepsilon,3/2\}}(\Omega)}\leq C\big(\|u\|_{H^s(\Omega)}
+\|{\mathcal{L}}u\|_{H^{s-2+\varepsilon}(\Omega)}\big)$.}   
\nonumber 
\end{align} 
\end{theorem}
%%%%%%%%%%
\begin{proof}
To set the stage, we claim that if $M_\psi$ denotes the operator of pointwise multiplication by a given 
function $\psi\in C^{1,1}(\overline{\Omega})$ then 
\begin{equation}\label{eqn.thu-TTG.Euclid.3}
M_\psi:H^s(\Omega)\rightarrow H^{s}(\Omega),\quad\forall\,s\in[-2,2], 
\end{equation}
is a linear and bounded mapping (compare with \eqref{eq:DDEj6g5}). 
Indeed, the case when $s\in[0,2]$ is seen via interpolation between $s=0$ and $s=2$. 
Moreover, since pointwise multiplication with a function does not 
increase the support, pointwise multiplication by $\psi\in C^{1,1}(\overline{\Omega})$
induces a well defined, linear, and bounded operator from $H^s_0(\Omega)$ into itself for each $s\in[0,2]$. 
Based on this and duality (cf. \eqref{eq:Redxax.1REF}) we then conclude that $M_\psi$ maps $\big(H^s_0(\Omega)\big)^*=H^{-s}(\Omega)$ 
linearly and boundedly into itself for every $s\in[0,2]$. As such, \eqref{eqn.thu-TTG.Euclid.3} is established.

As an immediate consequence of \eqref{eqn.thu-TTG.Euclid.3} and \eqref{tgBBn.344f} we see that,  
given any function $\psi\in C^{1,1}(\overline{\Omega})$, it follows that the operator 
\begin{equation}\label{eqn.thu-TTG.Euclid.4}
\text{$M_\psi$ maps $H^s_{\rm loc}(\Omega)$ into itself, for each $s\in[-2,2]$.}
\end{equation}
In particular, from \eqref{eqn.thu-TTG.Euclid.4} \big(considered with $s=-2$ and $\psi$ any of the entries 
$a_{ij}\in C^{1,1}(\overline{\Omega})$ of the coefficient matrix $A=\big(a_{jk}\big)_{1\leq j,k\leq n}$\big) 
we conclude that the differential expression ${\mathcal{L}}$ acts in a meaningful manner (as in indicated in \eqref{hyrdfr4.Euclid}) 
on any given distribution $u\in H^{-1}_{\rm loc}(\Omega)$ and, in fact, ${\mathcal{L}}u\in H^{-2}_{\rm loc}(\Omega)$.
Let us also note here that, as seen from \eqref{eqn.thu-TTG.Euclid.3} and the duality formula recorded in \eqref{eq:Redxax.1}, 
for each function $\psi\in C^{1,1}(\overline{\Omega})$ it follows that 
\begin{equation}\label{eqn.thu-TTG.Euclid.3.xxx}
M_\psi:H^s_0(\Omega)\rightarrow H^{s}_0(\Omega),\quad\forall\,s\in[-2,2], 
\end{equation}
is a well defined, linear, and bounded mapping.

Next, from the proof of Theorem~\ref{thm.Regularity-CC} we know that there exists a ($C^\infty$) smooth, compact, boundaryless, 
manifold $M$, of real dimension $n$, which contains $\overline{\Omega}$ in a single coordinate chart and which may be equipped  
with a Riemannian metric tensor $g$ possessing $C^{1,1}$-coefficients such that 
\begin{equation}\label{eqn.thu-TTG.Euclid}
{\mathcal{L}}=\big(\text{det}(A)\big)^{1/(n-2)}\Delta_g\,\text{ near }\,\overline{\Omega},
\end{equation}
where $\Delta_g$ denotes the Laplace--Beltrami operator on the Riemannian manifold $M$, 
associated (as in \eqref{eqn.thuaou}) with the metric tensor $g$. One also observes that 
\begin{equation}\label{eqn.thu-TTG.Euclid.2}
\big(\text{det}(A)\big)^{1/(n-2)}\in C^{1,1}(\overline{\Omega}),\quad
\big(\text{det}(A)\big)^{-1/(n-2)}\in C^{1,1}(\overline{\Omega}).
\end{equation}
Collectively, \eqref{eqn.thu-TTG.Euclid.3}, \eqref{eqn.thu-TTG.Euclid}, and \eqref{eqn.thu-TTG.Euclid.2} prove that 
for any given distribution $u\in H^{-1}_{\rm loc}(\Omega)$ and any given index $s\in[-2,2]$ we have 
\begin{equation}\label{eqn.thu-TTG.Euclid.4ii}
{\mathcal{L}}u\in H^{s}(\Omega)\,\text{ if and only if }\,\Delta_g u\in H^{s}(\Omega)
\end{equation}
in a quantitative fashion (i.e., with naturally accompanying estimates), as well as 
\begin{equation}\label{eqn.thu-TTG.Euclid.5}
{\mathcal{L}}u=0\,\text{ in }\,\Omega\,\text{ if and only if }\,\Delta_g u=0\,\text{ in }\,\Omega. 
\end{equation}
Given \eqref{eqn.thu-TTG.Euclid.4ii}--\eqref{eqn.thu-TTG.Euclid.5}, all conclusions in 
Theorem~\ref{YTfdf-T} (formulated in relation to the Laplace--Beltrami operator $\Delta_g$) translate into the 
properties claimed in the current statement. 
\end{proof}
%%%%%%%

Following past conventions, we will use the same symbol $\gamma_D$ in connection with 
either \eqref{eqn:gammaDs.1}, or \eqref{eqn:gammaDs.2aux-MMM.Euclid}. A special case of 
Theorem~\ref{YTfdf-T-MMM.Euclid}, which is especially useful in applications, is recorded below.

%%%%%%%%%%
\begin{corollary}\label{YTfdf-T.NNN-MMM.Euclid}
Fix an arbitrary $\varepsilon>0$, suppose $\Omega\subset{\mathbb{R}}^n$ is a bounded Lipschitz domain, and consider a second-order 
divergence-form differential expression ${\mathcal{L}}$, acting on each distribution 
$u\in H^{-1}_{\rm loc}(\Omega)$ according to 
\begin{equation}\label{hyrdfr4.Euclid.CCC}
{\mathcal{L}}u:=\sum\limits_{j,k=1}^n\partial_j\big(a_{jk}(x)\partial_ku\big)\,\text{ in }\,\Omega,
\end{equation}
in the sense of distributions, where $A(x)=\big(a_{jk}(x)\big)_{1\leq j,k\leq n}$, with $x\in\Omega$, is a symmetric, 
positive definite matrix, with real-valued entries $a_{jk}\in C^{1,1}(\overline{\Omega})$. 

Then the restriction of the operator \eqref{eqn:gammaDs.1} to $\big\{u\in H^s(\Omega)\,\big|\,{\mathcal{L}}u\in L^2(\Omega)\big\}$, 
originally considered for $s\in\big(\tfrac{1}{2},\tfrac{3}{2}\big)$, induces 
a well defined, linear, continuous operator 
\begin{align}\label{eqn:gammaDs.2-MMM.Euclid}
\gamma_D:\big\{u\in H^s(\Omega)\,\big|\,{\mathcal{L}}u\in L^2(\Omega)\big\}
\rightarrow H^{s-(1/2)}(\partial\Omega),\quad\forall\, 
s\in\big[\tfrac{1}{2},\tfrac{3}{2}\big]    
\end{align}
{\rm (}throughout, the space on the left-hand side of \eqref{eqn:gammaDs.2-MMM.Euclid} being equipped with the 
natural graph norm $u\mapsto\|u\|_{H^{s}(\Omega)}+\|{\mathcal{L}}u\|_{L^{2}(\Omega)}${\rm )}, which continues 
to be compatible with \eqref{eqn:gammaDs.1} when $s\in\big(\tfrac{1}{2},\tfrac{3}{2}\big)$,
and also with the pointwise nontangential trace, whenever the latter exists.

In addition, the following properties are true:

\begin{enumerate}
\item[$(i)$] The Dirichlet boundary trace operator in \eqref{eqn:gammaDs.2-MMM.Euclid} is surjective
and, in fact, there exist linear and bounded operators
\begin{equation}\label{2.88X-NN-ii-RRDD-2-MMM.Euclid}
\Upsilon_D:H^{s-(1/2)}(\partial\Omega)\to
\big\{u\in H^s(\Omega)\,\big|\,{\mathcal{L}}u\in L^2(\Omega)\big\},\quad s\in\big[\tfrac12,\tfrac32\big],
\end{equation}
which are compatible with one another and serve as right-inverses for the Dirichlet trace, that is, 
\begin{equation}\label{2.88X-NN2-ii-RRDD-2-MMM.Euclid}
\gamma_D(\Upsilon_D\psi)=\psi,
\quad\forall\,\psi\in H^{s-(1/2)}(\partial\Omega)\,\text{ with }\,s\in\big[\tfrac12,\tfrac32\big].
\end{equation}
Actually, matters may be arranged so that each function in the range of $\Upsilon_D$ is a null-solution of ${\mathcal{L}}$, that is, 
\begin{equation}\label{2.88X-NN2-ii-RRDD.bis.4.Euclid}
{\mathcal{L}}(\Upsilon_D\psi)=0,\quad\forall\,\psi\in H^{s-(1/2)}(\partial\Omega)
\,\text{ with }\,s\in\big[\tfrac12,\tfrac32\big].
\end{equation}

\item[$(ii)$] For each $s\in\big[\tfrac{1}{2},\tfrac{3}{2}\big]$, the null space   
of the Dirichlet boundary trace operator \eqref{eqn:gammaDs.2-MMM.Euclid} satisfies
\begin{equation}\label{eq:EFFa.NNN-MMM.Euclid}
{\rm ker}(\gamma_D)\subseteq H^{3/2}(\Omega).
\end{equation}
In fact, the inclusion in \eqref{eq:EFFa.NNN-MMM.Euclid} is quantitative in the sense that
there exists a constant $C\in(0,\infty)$ with the property that
\begin{align}\label{gafvv.6577-MMM.Euclid}
\begin{split}
&\text{whenever $u\in H^{1/2}(\Omega)$ with ${\mathcal{L}}u\in L^2(\Omega)$ satisfies $\gamma_D u=0$, then}
\\[2pt]
& \quad\text{$u\in H^{3/2}(\Omega)$ and 
$\|u\|_{H^{3/2}(\Omega)}\leq C\big(\|u\|_{L^2(\Omega)}+\|{\mathcal{L}}u\|_{L^2(\Omega)}\big)$.}
\end{split}
\end{align}
\end{enumerate}
\end{corollary}
%%%%%%%%%%%
\begin{proof}
All claims are obtained from their counterparts in the statement of Theorem~\ref{YTfdf-T-MMM.Euclid}, specialized to 
the case when $\varepsilon:=2-s$.
\end{proof}
%%%%%%%

After introducing the weak Neumann trace operator in the present setting, we continue by presenting a 
variable-coefficient version of the Euclidean weak Neumann trace result from Theorems~\ref{GNT} and \ref{YTfdf.NNN.2-Main}.
%%%%%%%
\begin{definition}\label{gycMM}
Let $\Omega\subset{\mathbb{R}}^n$ be a bounded Lipschitz domain. Then for some fixed smoothness exponent 
$s\in\big(\tfrac12,\tfrac32\big)$, the {\tt weak Neumann trace operator} 
\begin{equation}\label{2.88X-MMM.Euclid}
\widetilde\gamma_{N,\mathcal{L}}:\big\{(f,F)\in H^s(\Omega)\times H^{s-2}_0(\Omega) 
\,\big|\,{\mathcal{L}}f=F|_{\Omega}\text{ in }\mathcal{D}'(\Omega)\big\} 
\rightarrow H^{s-(3/2)}(\partial\Omega)
\end{equation}
is defined as follows: Suppose that some function $f\in H^s(\Omega)$ along with some distribution 
$F\in H^{s-2}_0(\Omega)\subset H^{s-2}({\mathbb{R}}^n)$ satisfying ${\mathcal{L}}f=F|_{\Omega}$ 
in $\mathcal{D}'(\Omega)$ have been given. In particular, 
\begin{equation}\label{u543d-MMM.Euclid}
\nabla f\in[H^{s-1}(\Omega)]^n=\big([H^{1-s}(\Omega)]^n\big)^*. 
\end{equation}
Then the manner in which $\widetilde\gamma_{N,\mathcal{L}}(f,F)$ is now defined as a functional 
in the space $H^{s-(3/2)}(\partial\Omega)=\big(H^{(3/2)-s}(\partial\Omega)\big)^*$ is as follows: 
Given $\phi\in H^{(3/2)-s}(\partial\Omega)$, then for any $\Phi\in H^{2-s}(\Omega)$ such that $\gamma_D\Phi=\phi$ 
{\rm (}whose existence is ensured by the surjectivity of \eqref{eqn:gammaDs.1}{\rm )}, set
\begin{align}\label{2.9NEW-MMM.Euclid}
\begin{split} 
& {}_{H^{(3/2)-s}(\partial\Omega)}\big\langle\phi,
\widetilde\gamma_{N,\mathcal{L}}(f,F)\big\rangle_{(H^{(3/2)-s}(\partial\Omega))^*}  \\[2pt]
& \quad:={}_{[H^{1-s}(\Omega)]^n}\big\langle A\nabla\Phi,\nabla f\big\rangle_{([H^{1-s}(\Omega)]^n)^*}
+{}_{H^{2-s}(\Omega)}\big\langle\Phi,F\big\rangle_{(H^{2-s}(\Omega))^*}.
\end{split} 
\end{align}
Then the weak Neumann trace mapping \eqref{2.88X-MMM.Euclid} is an operator which is unambiguously defined, 
linear, and bounded {\rm (}assuming the space on the left-hand side of \eqref{2.88X-MMM.Euclid} 
is equipped with the natural norm 
$(f,F)\mapsto\|f\|_{H^{s}(\Omega)}+\|F\|_{H^{s-2}({\mathbb{R}}^n)}${\rm )}. 
\end{definition}
%%%%%%%

The above definition plays a basic role in the following theorem, which may be regarded as a variable-coefficient version 
of the Neumann trace result established (for the ordinary Laplacian) earlier in Theorems~\ref{GNT} and \ref{YTfdf.NNN.2-Main}.

%%%%%%%
\begin{theorem}\label{GNT.Euclid}
Assume $\Omega\subset{\mathbb{R}}^n$ is a bounded Lipschitz domain and consider a second-order 
divergence-form differential expression ${\mathcal{L}}$, acting on each distribution 
$u\in H^{-1}_{\rm loc}(\Omega)$ according to 
\begin{equation}\label{hyrdfr4.Euclid.iii}
{\mathcal{L}}u:=\sum\limits_{j,k=1}^n\partial_j\big(a_{jk}(x)\partial_ku\big)\,\text{ in }\,\Omega,
\end{equation}
in the sense of distributions, where $A(x)=\big(a_{jk}(x)\big)_{1\leq j,k\leq n}$, with $x\in\Omega$, is a symmetric, positive definite matrix, with real-valued entries $a_{jk}\in C^{1,1}(\overline{\Omega})$.

Then for each $\varepsilon>0$, the weak Neumann boundary trace map, originally introduced as in \eqref{2.88X-MMM.Euclid}, induces linear and continuous operators in the context  
\begin{align}\label{eqn:gammaN-MMM.Euclid.XXX}
& \widetilde\gamma_{N,\mathcal{L}}:\big\{(f,F)\in H^s(\Omega)\times H^{s-2+\varepsilon}_0(\Omega)\,|\,
{\mathcal{L}}f=F\big|_{\Omega}\,\text{ in }\,{\mathcal{D}}'(\Omega)\big\}\rightarrow H^{s-(3/2)}(\partial\Omega)
\no \\
& \quad \text{with }\,\,s\in\big[\tfrac{1}{2},\tfrac{3}{2}\big]
\end{align}
{\rm (}where the space on the left-hand side of \eqref{eqn:gammaN-MMM.Euclid.XXX} is equipped with the natural norm 
$(f,F)\mapsto\|f\|_{H^{s}(\Omega)}+\|F\|_{H^{s-2+\varepsilon}({\mathbb{R}}^n)}${\rm )} which are compatible with 
those in \eqref{2.88X-MMM.Euclid} when $s\in\big(\tfrac{1}{2},\tfrac{3}{2}\big)$. 
Thus defined, the weak Neumann boundary trace map possesses the following properties: \\[1mm] 
$(i)$ The weak Neumann trace operators corresponding to various values of the 
parameter $s\in\big[\tfrac12,\tfrac32\big]$ are compatible with one another and each 
of them is surjective. In fact, there exist linear and bounded operators
\begin{equation}\label{2.88X-NN.Euclid}
\Upsilon_{N,\mathcal{L}}:H^{s-(3/2)}(\partial\Omega)\rightarrow 
\big\{u\in H^s(\Omega)\,\big|\,{\mathcal{L}}u\in L^2(\Omega)\big\},\quad s\in\big[\tfrac12,\tfrac32\big],
\end{equation}
which are compatible with one another and satisfy {\rm (}with tilde denoting the extension by zero 
outside $\Omega${\rm )} 
\begin{equation}\label{2.88X-NN2.Euclid}
\widetilde\gamma_{N,\mathcal{L}}\big(\Upsilon_{N,\mathcal{L}}\psi,\widetilde{{\mathcal{L}}(\Upsilon_{N,\mathcal{L}}\psi)}\,\big)=\psi,
\quad\forall\,\psi\in H^{s-(3/2)}(\partial\Omega)\,\text{ with }\,s\in\big[\tfrac12,\tfrac32\big].
\end{equation}
$(ii)$ If $\varepsilon\in(0,1)$ and $s\in\big[\tfrac{1}{2},\tfrac{3}{2}\big]$ then for any two pairs
\begin{align}\label{n7b66-MMM.Euclid}
\begin{split}
& \text{$(f_1,F_1)\in H^s(\Omega)\times H^{s-2+\varepsilon}_0(\Omega)$ 
such that ${\mathcal{L}}f_1=F_1|_{\Omega}$ in $\mathcal{D}'(\Omega)$, and}
\\[2pt]
& \quad\text{$(f_2,F_2)\in H^{2-s}(\Omega)\times H^{-s+\varepsilon}_0(\Omega)$ 
such that ${\mathcal{L}}f_2=F_2|_{\Omega}$ in $\mathcal{D}'(\Omega)$}, 
\end{split}
\end{align}
the following Green's formula holds:
\begin{align}\label{GGGRRR-prim222-MMM.Euclid}
& {}_{H^{(3/2)-s}(\partial\Omega)}\big\langle\gamma_D f_2,\widetilde\gamma_{N,\mathcal{L}}(f_1,F_1)
\big\rangle_{(H^{(3/2)-s}(\partial\Omega))^*}
\nonumber\\[2pt]
& \qquad -{}_{(H^{s-(1/2)}(\partial\Omega))^*}\big\langle\widetilde\gamma_{N,\mathcal{L}}(f_2,F_2),\gamma_D f_1
\big\rangle_{H^{s-(1/2)}(\partial\Omega)}     
\nonumber\\[2pt]
&\quad={}_{H^{2-s}(\Omega)}\big\langle f_2,F_1\big\rangle_{(H^{2-s}(\Omega))^*}
-{}_{(H^s(\Omega))^*}\big\langle F_2,f_1\big\rangle_{H^s(\Omega)}.
\end{align}
$(iii)$ There exists a constant $C\in(0,\infty)$ with the property that
\begin{align}\label{gafvv.6588-P-MMM.Euclid}
\begin{split}
& \text{if $f\in H^{1/2}(\Omega)$ and $F\in H^{-(3/2)+\varepsilon}_0(\Omega)$ with 
$0<\varepsilon\leq 1$ satisfy} 
\\[2pt]
& \quad\text{${\mathcal{L}}f=F\big|_{\Omega}$ in ${\mathcal{D}}'(\Omega)$ and 
$\widetilde\gamma_{N,\mathcal{L}}(f,F)=0$, then $f\in H^{(1/2)+\varepsilon}(\Omega)$ }
\\[2pt]
& \quad\text{and $\|f\|_{H^{(1/2)+\varepsilon}(\Omega)}
\leq C\big(\|f\|_{L^2(\Omega)}+\|F\|_{H^{-(3/2)+\varepsilon}({\mathbb{R}}^n)}\big)$.}
\end{split}
\end{align} 
\end{theorem}
%%%%%%%
\begin{proof}
Bringing back the ($C^\infty$) smooth, compact, boundaryless, manifold $M$, of real dimension $n$, from the proof of 
Theorem~\ref{thm.Regularity-CC}, this has the property that $\overline{\Omega}$ is contained in a single coordinate 
chart of $M$. Moreover, if we equip $M$ with the $C^{1,1}$ Riemannian metric tensor $g$ defined as in \eqref{hyrdfr5-RR}, then 
\begin{equation}\label{eqn.thu-TTG.Euclid.xxx}
{\mathcal{L}}=\big(\text{det}(A)\big)^{1/(n-2)}\Delta_g\,\text{ near }\,\overline{\Omega},
\end{equation}
where $\Delta_g$ denotes the Laplace--Beltrami operator on the Riemannian manifold $M$, associated (as in \eqref{eqn.thuaou}) with the metric tensor $g$.

Based on \eqref{eq:52d2cv.1}, \eqref{eq:52d2cv.2}, and \eqref{VCnnd.M.WACO.EEE} we conclude that 
for any $\phi\in H^{2-s}(\Omega)$ and $\psi\in H^{s}(\Omega)$ with $s\in\big(\tfrac12,\tfrac32\big)$ we have 
\begin{align}\label{VCnnd.M.WACO.EEE.ii}
& {}_{H^{2-s}(\Omega)}\big\langle\phi,\psi\big\rangle_{(H^{2-s}(\Omega))^*}, 
\,\text{ with $\Omega$ viewed as a set in $M$,}
\no \\[4pt]
& \quad \text{coincides with }\,{}_{H^{2-s}(\Omega)}\Big\langle\phi,\big(\text{det}(A)\big)^{1/(n-2)}\psi\Big\rangle_{(H^{2-s}(\Omega))^*}, 
\\[4pt]
& \quad \text{where $\Omega$ is now viewed as an open set in ${\mathbb{R}}^n$}.   \no 
\end{align}
On account of \eqref{eq:52d2cv.1}, \eqref{eq:52d2cv.2}, and \eqref{VCnnd.M.WACO} we also see that if $\phi\in H^{2-s}(\Omega)$
and $\psi\in H^{s}(\Omega)$ for some $s\in\big(\tfrac12,\tfrac32\big)$ then 
\begin{align}\label{VCnnd.M.WACO.ii}
\begin{split} 
& {}_{H^{2-s}(\Omega,TM)}\big\langle{\rm grad}_g\phi,  {\rm grad}_g\psi\big\rangle_{(H^{1-s}(\Omega,TM))^*}
\\[2pt]
& \quad ={}_{[H^{1-s}(\Omega)]^n}\big\langle A\nabla\phi,\nabla\psi\big\rangle_{([H^{1-s}(\Omega)]^n)^*}.
\end{split} 
\end{align}

In addition, define the operator ${\mathcal{M}}$ mapping the space 
\begin{equation}\label{eq:Space.1}
\big\{(f,F)\in H^{-1}_{\rm loc}(\Omega)\times H^{-2}_0(\Omega) 
\,\big|\,{\mathcal{L}}f=F|_{\Omega}\text{ in }\mathcal{D}'(\Omega)\big\} 
\end{equation}
(where $\Omega$ is viewed as an open set in ${\mathbb{R}}^n$) 
into the space 
\begin{equation}\label{eq:Space.2}
\big\{(f,F)\in H^{-1}_{\rm loc}(\Omega)\times H^{-2}_0(\Omega) 
\,\big|\,\Delta_g f=F|_{\Omega}\text{ in }\mathcal{D}'(\Omega)\big\} 
\end{equation}
(where $\Omega$ is now regarded as an open subset of the Riemannian manifold $(M,g)$) according to 
\begin{equation}\label{eq:Operator}
{\mathcal{M}}(f,F):=\Big(f\,,\,\big(\text{det}(A)\big)^{1/(n-2)}F\Big).
\end{equation}
Thanks to \eqref{eqn.thu-TTG.Euclid.xxx}, \eqref{eqn.thu-TTG.Euclid.2}, and \eqref{eqn.thu-TTG.Euclid.3.xxx}, 
this is a well defined linear operator. In this regard, the key observation is that for each $s\in\big[\tfrac12,\tfrac32\big]$ 
we have
\begin{align}\label{2.88X-MMM.Euclid.tV}
& \widetilde\gamma_{N,\mathcal{L}}=\widetilde\gamma_N\circ{\mathcal{M}}\,\text{ as operators acting from the space}      \no 
\\[2pt]
& \quad \big\{(f,F)\in H^s(\Omega)\times H^{s-2}_0(\Omega)\,\big|\,{\mathcal{L}}f=F|_{\Omega}\text{ in }\mathcal{D}'(\Omega)\big\}
\\[2pt]
& \quad \text{and taking values into the space }\,H^{s-(3/2)}(\partial\Omega),     \no
\end{align}
where $\widetilde\gamma_N$ is the weak Neumann trace operator associated as in Theorem~\ref{YTfdf.NNN.2-Main-MMM} 
when $\Omega$ is regarded as a subdomain of the Riemannian manifold $(M,g)$ (see also  \eqref{2.9NEW-MMM}). Indeed, if $s\in\big(\tfrac12,\tfrac32\big)$ 
then \eqref{2.88X-MMM.Euclid.tV} is seen directly from \eqref{2.9NEW-MMM.Euclid}, \eqref{eq:Operator}, 
\eqref{eqn.thu-TTG.Euclid.xxx}, \eqref{VCnnd.M.WACO.EEE.ii}, and \eqref{VCnnd.M.WACO.ii}. Since the scale of Sobolev spaces 
is nested, this also covers (a posteriori) the end-point case $s=\tfrac32$. Finally, in the case $s=\tfrac12$ we take 
\eqref{2.88X-MMM.Euclid.tV} as a definition of the weak Neumann trace operator $\widetilde\gamma_{N,\mathcal{L}}$. 

Given that for each $s\in\big[\tfrac12,\tfrac32\big]$ and $\varepsilon\in(0,1)$ the operator ${\mathcal{M}}$ becomes an 
isomorphism of the space 
\begin{equation}\label{eq:Space.1.A}
\big\{(f,F)\in H^{s}(\Omega)\times H^{s-2+\varepsilon}_0(\Omega) 
\,\big|\,{\mathcal{L}}f=F|_{\Omega}\text{ in }\mathcal{D}'(\Omega)\big\} 
\end{equation}
(where $\Omega$ is viewed as an open set in ${\mathbb{R}}^n$) onto the space 
\begin{equation}\label{eq:Space.2.B}
\big\{(f,F)\in H^{2}(\Omega)\times H^{s-2+\varepsilon}_0(\Omega) 
\,\big|\,\Delta_g f=F|_{\Omega}\text{ in }\mathcal{D}'(\Omega)\big\} 
\end{equation}
(where $\Omega$ is now regarded as an open subset of the Riemannian manifold $(M,g)$), all claims in the statement of the current theorem become relatively straightforward consequences of \eqref{2.88X-MMM.Euclid.tV} and the 
corresponding properties of the weak Neumann trace operator $\widetilde\gamma_N$ from Theorem~\ref{YTfdf.NNN.2-Main-MMM}
(while also bearing in mind \eqref{eqn.thu-TTG.Euclid.xxx}, \eqref{eqn.thu-TTG.Euclid.2}, \eqref{eqn.thu-TTG.Euclid.3}, 
and \eqref{eqn.thu-TTG.Euclid.3.xxx}). 
\end{proof}
%%%%%%%

We conclude by presenting the following special case of Theorem~\ref{GNT.Euclid}, which plays a significant role
in applications.

%%%%%%%%%%
\begin{corollary}\label{YTfdf.NNN.2-MMM.CCC}
Suppose $\Omega\subset{\mathbb{R}}^n$ is a bounded Lipschitz domain, and denote by $\nu$ its outward unit normal. In addition, consider a second-order divergence-form differential expression ${\mathcal{L}}$, acting on each distribution 
$u\in H^{-1}_{\rm loc}(\Omega)$ according to 
\begin{equation}\label{hyrdfr4.Euclid.jjj}
{\mathcal{L}}u:=\sum\limits_{j,k=1}^n\partial_j\big(a_{jk}(x)\partial_ku\big)\,\text{ in }\,\Omega,
\end{equation}
in the sense of distributions, where $A(x)=\big(a_{jk}(x)\big)_{1\leq j,k\leq n}$, with $x\in\Omega$, is a symmetric, 
positive definite matrix, with real-valued entries $a_{jk}\in C^{1,1}(\overline{\Omega})$.

Then the Neumann trace map, originally defined for each for $u\in C^\infty(\overline{\Omega})$
as $u\mapsto\langle\nu,A\nabla u\rangle$ on $\partial\Omega$, extends uniquely to linear 
continuous operators  
\begin{equation}\label{eqn:gammaN-pp-MMM.CCC}
\gamma_N:\big\{u\in H^s(\Omega)\,\big|\,{\mathcal{L}}u\in L^2(\Omega)\big\}\to H^{s-(3/2)}(\partial\Omega), 
\quad s\in\big[\tfrac{1}{2},\tfrac{3}{2}\big] 
\end{equation}
{\rm (}throughout, the space on the left-hand side of \eqref{eqn:gammaN-pp-MMM.CCC} is equipped with the natural graph 
norm $u\mapsto\|u\|_{H^{s}(\Omega)}+\|{\mathcal{L}}u\|_{L^2(\Omega)}${\rm )}, that are compatible with one another.
In addition, the following properties are true:

\begin{enumerate}
\item[$(i)$] The Neumann trace map \eqref{eqn:gammaN-pp-MMM.CCC} is surjective. In fact, there exist linear and bounded operators
\begin{equation}\label{2.88X-NN-ii-RR-MMM.CCC}
\Upsilon_N:H^{s-(3/2)}(\partial\Omega)\to
\big\{u\in H^s(\Omega)\,\big|\,{\mathcal{L}}u\in L^2(\Omega)\big\},\quad s\in\big[\tfrac12,\tfrac32\big],
\end{equation}
which are compatible with one another and are right-inverses for the Neumann trace, that is, 
\begin{equation}\label{2.88X-NN2-ii-RR-MMM.CCC}
\gamma_N(\Upsilon_N\psi)=\psi,\quad\forall\,\psi\in H^{s-(3/2)}(\partial\Omega)\,\text{ with }\,s\in\big[\tfrac12,\tfrac32\big].
\end{equation}

\item[$(ii)$] If $s\in\big[\tfrac{1}{2},\tfrac{3}{2}\big]$, then for any functions $f\in H^s(\Omega)$ with 
${\mathcal{L}}f\in L^2(\Omega)$ and $h\in H^{2-s}(\Omega)$ with ${\mathcal{L}}h\in L^2(\Omega)$ the following 
Green's formula holds:
\begin{align}\label{GGGRRR-MMM.CCC}
& {}_{H^{(3/2)-s}(\partial\Omega)}\big\langle\gamma_D h,\gamma_N f
\big\rangle_{(H^{(3/2)-s}(\partial\Omega))^*}
\nonumber\\[2pt]
& \qquad-{}_{(H^{s-(1/2)}(\partial\Omega))^*}\big\langle\gamma_N h,\gamma_D f
\big\rangle_{H^{s-(1/2)}(\partial\Omega)}     
\nonumber\\[2pt]
& \quad=(h,{\mathcal{L}}f)_{L^2(\Omega)}-({\mathcal{L}}h,f)_{L^2(\Omega)}.
\end{align}

\item[$(iii)$] For each $s\in\big[\tfrac{1}{2},\tfrac{3}{2}\big]$, the null space of 
the Neumann boundary trace operator \eqref{eqn:gammaN-pp-MMM.CCC} satisfies
\begin{equation}\label{eq:EFFa.111-MMM.CCC}
{\rm ker}(\gamma_N)\subseteq H^{3/2}(\Omega).
\end{equation}
In fact, the inclusion in \eqref{eq:EFFa.111-MMM.CCC} is quantitative in the sense that
there exists a constant $C\in(0,\infty)$ with the property that
\begin{align}\label{gafvv.6588-MMM.CCC}
\begin{split}
& \text{whenever $u\in H^{1/2}(\Omega)$ satisfies ${\mathcal{L}}u\in L^2(\Omega)$ and $\gamma_N u=0$, then}
\\[2pt]
& \quad\text{$u\in H^{3/2}(\Omega)$ and 
$\|u\|_{H^{3/2}(\Omega)}\leq C\big(\|u\|_{L^2(\Omega)}+\|{\mathcal{L}}u\|_{L^2(\Omega)}\big)$.}
\end{split}
\end{align}
\end{enumerate}
\end{corollary}
%%%%%%%%%%
\begin{proof}
Having fixed $s\in\big[\tfrac{1}{2},\tfrac{3}{2}\big]$, pick $0<\varepsilon<\min\{1,2-\varepsilon\}$ and define 
\begin{align}\label{2.10X-MMM.CCC}
\begin{split} 
&\iota:\big\{u\in H^s(\Omega)\,\big|\,{\mathcal{L}}u\in L^2(\Omega)\big\}       
\\[2pt]
&\quad\;\rightarrow
\big\{(f,F)\in H^s(\Omega)\times H^{s-2+\varepsilon}_0(\Omega)\,\big|\,{\mathcal{L}}f=F\big|_{\Omega}
\,\text{ in }\,{\mathcal{D}}'(\Omega)\big\}, 
\end{split} 
\end{align}
as being the continuous injection given by 
\begin{equation}\label{2.11X-MMM.CCC}
\iota(u):=\Big(u,\widetilde{{\mathcal{L}}u}\Big),\quad\forall\,u\in H^s(\Omega)\,\text{ with }\,{\mathcal{L}}u\in L^2(\Omega),
\end{equation}
where tilde denotes the extension by zero outside $\Omega$. With the weak Neumann trace operator 
$\widetilde\gamma_{N,\mathcal{L}}$ associated with ${\mathcal{L}}$ as in Theorem~\ref{GNT.Euclid}, 
we then set 
\begin{equation}\label{2.12X-MMM.CCC}
\gamma_N:=\widetilde\gamma_{N,\mathcal{L}}\circ\iota.
\end{equation}
Thanks to the continuity of $\iota$ in \eqref{2.10X-MMM.CCC} and $\widetilde\gamma_{N,\mathcal{L}}$ in 
\eqref{eqn:gammaN-MMM.Euclid.XXX}, this is a well defined, linear, and bounded mapping in the context of \eqref{eqn:gammaN-pp-MMM.CCC}.
In fact, all other claims in the statement are clear from \eqref{2.12X-MMM.CCC} and Theorem~\ref{GNT.Euclid}.
\end{proof}
%%%%%%

%%%%%%%%%%%%%%%%%%%%%%%%%%%%%%%%
\noindent 
{\bf Acknowledgments.} We are indebted to the referees for a very careful reading of our manuscript. 
The level of detail recorded in these reports was extraordinarily helpful and constructive. 
%%%%%%%%%%%%%%%%%%%%%%%%%%%%%%%%

%%%%%%%%%%%%%%%%%%%%%%%%%%%%%%%%
%%%%%%%%%%%%%%%%%%%%%%%%%%%%%%%%
 
\end{document}